\documentclass[11pt]{article}
\usepackage{amsthm,amsmath,amssymb,bbm}
\usepackage{color}
\def\bbe{\mathbb E}

\def\bbr{\mathbb R}

\def\bbone{{\mathbbm 1}}
\theoremstyle{plain}
\newtheorem{thm}{Theorem}[section]

\newtheorem{cor}{Corollary}[section]
\newtheorem{prop}{Proposition}[section]
\newtheorem{lem}{Lemma}[section]
\newtheorem{rem}{Remark}[section]
\newtheorem{defi}{Definition}[section]
\allowdisplaybreaks[1]
\numberwithin{equation}{section}

\begin{document}
\title{Around the Sobolev Inequalities for the Stable Heat Semigroups}

\author{Benjamin Arras\thanks{Univ. Lille, CNRS, UMR 8524 - Laboratoire Paul Painlevé, F-59000 Lille, France; benjamin.arras@univ-lille.fr}\; and Christian Houdr\'e\thanks{Georgia Institute of Technology, 
School of Mathematics, Atlanta, GA 30332-0160, USA; houdre@math.gatech.edu. 
Research supported in part by the grant \# $524678$ from the Simons Foundation. CH would like to thank the Laboratoire Paul Painlev\'e of the Universit\'e de Lille for its hospitality while part of this research was carried out.
}}

\maketitle

\vspace{\fill}
\begin{abstract}
\noindent
We develop a general distributional theory of fractional (an)isotropic Sobolev spaces associated with the non-degenerate symmetric $\alpha$-stable, $\alpha \in (1,2)$, probability measures on $\mathbb{R}^d$. 
\end{abstract}
\vspace{\fill}

\tableofcontents

\section{Introduction}
\subsection{Overview and main results}
The classical Sobolev and isoperimetric inequalities are fundamental objects in mathematical and geometric analysis. They are ubiquitous in the modern theory of Sobolev spaces and partial differential equations (see, e.g., \cite{Stein_1,Lieb_Loss_book01,AF_book03,Mazya_book,B_FA11,Ponce_book16}). 
Moreover, these inequalities share tight connections with probability theory. Indeed, they naturally appear in the analysis of Markov diffusion operators as put forward in \cite{BGL14}. Starting in the nineteen-eighties with the work \cite{Varo_85}, a close connection has been systematically investigated between Sobolev-type inequalities and the mapping properties of certain semigroups of linear operators on Lebesgue spaces. Since then, this connection has been thoroughly studied as in \cite{CKS_aihp87,VSC_book92,CM_tams93,Coulhon_jfa96,AB_15,Kim_potal16}. In all those works, the central property of the semigroups under consideration is their ultracontractivity. Let $(X,\mathcal{X},m)$ be a sigma-finite measure space and, for every $p \in [1,+\infty]$, let $L^p(X,m)$ be the classical Lebesgue space over $(X,\mathcal{X},m)$ of exponent $p$. A semigroup of linear bounded operators $(P_t)_{t\geq 0}$ on $L^p(X,m)$, $p \in [1,+\infty]$, is ultracontractive with ``dimension" $n>0$ if, for all $f \in L^1(X,m)$ and all $t>0$, $P_t(f) \in L^{\infty}(X,m)$ and  
\begin{align}\label{eq:ultracontractivity_general_sg}
\|P_t(f)\|_{L^{\infty}(X,m)} \leq C t^{-\frac{n}{2}} \|f\|_{L^1(X,m)}, 
\end{align}
for some $C>0$ and where $\|\cdot\|_{L^1(X,m)}$ and $\|\cdot\|_{L^{\infty}(X,m)}$ are the classical norms on the respective Lebesgue spaces. A prototypical example of such a semigroup is the classical heat semigroup on the Euclidean Lebesgue spaces endowed with the Lebesgue measure. Recall that the heat semigroup $(P^H_t)_{t \geq 0}$ can be defined through an integral representation formula: for all $f \in L^1(\mathbb{R}^d,dx)$, all $t \geq 0$ and all $x \in \mathbb{R}^d$, 
\begin{align}\label{eq:classical_heat_sg}
P^H_t(f)(x) = \int_{\mathbb{R}^d} f\left(x + \sqrt{2 t} y\right) \gamma(dy),
\end{align} 
where $\gamma$ is the standard Gaussian measure on $\mathbb{R}^d$. From \eqref{eq:classical_heat_sg}, it is clear that $(P_t^H)_{t \geq 0}$ is a semigroup of linear Markovian contractions on $L^p(\mathbb{R}^d,dx)$, for $p \in [1,+\infty]$, symmetric on $L^2(\mathbb{R}^d,dx)$ and ultracontractive in the sense of \eqref{eq:ultracontractivity_general_sg} with dimension $n=d$. Here and in the sequel, Markovian means that the semigroup is positivity preserving and mass conservative.~Applying for example \cite[Theorem $1$]{Varo_85}, for all $d>2$ integer and all $f \in \mathcal{C}_c^{\infty}(\mathbb{R}^d)$,
\begin{align}\label{ineq:standard_Sobolev-type_inequality_p2}
\|f\|_{L^{2^*}(\mathbb{R}^d,dx)} \leq C \|(-\Delta)^{\frac{1}{2}}(f)\|_{L^2(\mathbb{R}^d,dx)}, 
\end{align}  
where $C>0$, $2^* = 2d/(d-2)$, $(-\Delta)^{1/2}$ stands for the square root of the Laplacian operator with Fourier symbol $\|\xi\|$, for all $\xi \in \mathbb{R}^d$, and $\mathcal{C}_c^{\infty}(\mathbb{R}^d)$ is the set of infinitely differentiable functions on $\mathbb{R}^d$ with compact support. (Here and in the sequel, $\|\cdot\|$ is the Euclidean norm on $\mathbb{R}^d$ and $\langle \cdot ; \cdot \rangle$ the Euclidean inner product.) Combined with the well-known equivalence of norms induced by the mapping property of the vectorial Riesz transform on $L^2(\mathbb{R}^d,dx)$, \eqref{ineq:standard_Sobolev-type_inequality_p2} implies the classical Sobolev inequality with the exponent $2$: for all $f \in \mathcal{C}_c^{\infty}(\mathbb{R}^d)$,
\begin{align}\label{ineq:standard_Sobolev_inequality_p2}
\|f\|_{L^{2^*}(\mathbb{R}^d,dx)} \leq C \| \nabla(f)\|_{L^2(\mathbb{R}^d,dx)},  
\end{align}
where $C>0$ and $\nabla$ is the gradient operator. Below, we are interested in the Sobolev-type and isoperimetric inequalities naturally associated with another family of ultracontractive semigroups: the stable heat semigroups. Stable probability measures on $\mathbb{R}^d$ encompass the standard Gaussian probability measure as a limiting case and are defined as follows. Let $\alpha \in (0,2)$ and let $\nu_\alpha$ be a positive Borel measure on $\mathbb{R}^d$ such that $\nu_\alpha(\{0\}) = 0$, $\int_{\mathbb{R}^d} (1 \wedge \|u\|^2) \nu_\alpha(du) < +\infty$ and, for all $c>0$,
\begin{align}\label{eq:scale}
c^{-\alpha}T_c(\nu_\alpha)(du)=\nu_\alpha(du),
\end{align}
where $T_c(\nu_\alpha)(B):=\nu_{\alpha}(B/c)$, for all $B$ Borel subset of $\mathbb{R}^d$. The measure $\nu_\alpha$ is called an $\alpha$-stable L\'evy measure, and (see \cite[Theorem $14.3.$]{S}) it admits the following polar decomposition 
\begin{align}\label{eq:polar}
\nu_\alpha(du) = \bbone_{(0,+\infty)}(r) \bbone_{\mathbb{S}^{d-1}}(y)\dfrac{dr}{r^{\alpha+1}}\sigma(dy),
\end{align}
where $\mathbb{S}^{d-1}$ is the Euclidean unit sphere of $\mathbb{R}^d$ and $\sigma$ is a finite positive measure on $\mathbb{S}^{d-1}$. In the sequel, $\alpha \in (1,2)$ and $\sigma$ is assumed to be \textit{symmetric}, i.e., $\sigma(B) = \sigma(-B)$, for all $B$ Borel subset of $\mathbb{S}^{d-1}$. Moreover, $\nu_\alpha$ is assumed to be \textit{non-degenerate} in the sense that 
\begin{align}\label{eq:non_deg}
\underset{y \in \mathbb{S}^{d-1}}{\inf} \int_{\mathbb{S}^{d-1}} |\langle y;x \rangle|^{\alpha} \lambda_1(dx) > 0,
\end{align}
where $\lambda_1$ is a finite positive measure on $\mathbb{S}^{d-1}$ called the \textit{spectral measure} and defined by 
\begin{align}\label{eq:spectral_measure}
\lambda_1(dx)= - \cos \left(\alpha \frac{\pi}{2}\right)\dfrac{\Gamma(2-\alpha)}{\alpha(\alpha-1)}\sigma(dx),
\end{align}
with $\Gamma$ the Euler Gamma function. Next, a probability measure $\mu_\alpha$ on $\mathbb{R}^d$ is called a non-degenerate symmetric $\alpha$-stable probability measure on $\mathbb{R}^d$ if its Fourier transform is given, for all $\xi \in \mathbb{R}^d$, by 
\begin{align}\label{stable:characteristic}
\widehat{\mu_\alpha}(\xi):= \int_{\mathbb{R}^d} e^{i \langle y ; \xi\rangle} \mu_\alpha(dy) = \exp\left(\int_{\mathbb{R}^d} (e^{i \langle u;\xi \rangle}-1-i\langle \xi;u\rangle) \nu_\alpha(du)\right),
\end{align} 
with $\nu_\alpha$ as above. Since $\sigma$ is symmetric, \cite[Theorem 14.13.]{S} provides a useful alternative representation (the L\'evy representation) for $\widehat{\mu_\alpha}$, namely, for all $\xi \in \mathbb{R}^d$,
\begin{align}\label{eq:rep_spectral_measure}
\widehat{\mu_\alpha}(\xi) = \exp\left(- \sigma_\alpha(\xi)^\alpha\right), \quad \sigma_\alpha(\xi)= \left(\int_{\mathbb{S}^{d-1}} |\langle y;\xi \rangle|^\alpha \lambda_1(dy)\right)^{\frac{1}{\alpha}}.
\end{align}
The stable heat semigroup $(P_t^\alpha)_{t \geq 0}$ associated with $\mu_\alpha$ is defined, for all $f\in L^1(\mathbb{R}^d,dx)$, all $t \geq 0$ and all $x \in \mathbb{R}^d$, by
\begin{align}\label{eq:StheatSM}
P^{\alpha}_t(f)(x) = \int_{\mathbb{R}^d} f\left(x+ t^{\frac{1}{\alpha}} y\right) \mu_\alpha(dy).
\end{align}
This semigroup is a special case of convolution semigroups for which a full theory has been well-developed (see \cite{NJ02_1,NJ02_2,NJ02_3}). In particular, $(P_t^{\alpha})_{t \geq 0}$ is a semigroup of linear Markovian contractions on $L^p(\mathbb{R}^d,dx)$, for $p \in [1,+\infty]$, symmetric on $L^2(\mathbb{R}^d,dx)$, strongly continuous on $L^p(\mathbb{R}^d,dx)$, for $p \in [1,+\infty)$, and ultracontractive with dimension $n=2d/\alpha$ in the sense of \eqref{eq:ultracontractivity_general_sg}. Thanks to \cite[Theorem $1$]{Varo_85}, for all integer $d \geq 2$ and all $f \in \mathcal{C}_c^{\infty}(\mathbb{R}^d)$, 
\begin{align}\label{ineq:stable_Sobolev-type_inequality_p2}
\|f\|_{L^{\frac{2n}{n-2}}(\mathbb{R}^d,dx)} \leq C \left\|\left(- \mathcal{A}_{\alpha}\right)^{\frac{1}{2}}(f)\right\|_{L^2(\mathbb{R}^d,dx)},
\end{align}
where $C>0$ and $\mathcal{A}_\alpha$ is the generator of $(P^\alpha_t)_{t \geq 0}$ given, for all $f \in \mathcal{C}_c^{\infty}(\mathbb{R}^d)$ and all $x \in \mathbb{R}^d$, by 
\begin{align}\label{eq:Stheatgen}
\mathcal{A}_\alpha(f)(x) = \int_{\mathbb{R}^d} \left(f(x+u)-f(x)-\langle u; \nabla(f)(x) \rangle\right) \nu_\alpha(du). 
\end{align}
This integro-differential operator is a non-local, fractional and (an)isotropic extension of the Laplace operator $\Delta$ and admits an alternative form thanks to an integration by parts in the radial coordinate: for all $f \in \mathcal{C}_c^{\infty}(\mathbb{R}^d)$ and all $x \in \mathbb{R}^d$, 
\begin{align}\label{eq:Stheatgen_nice_decomposition}
\mathcal{A}_\alpha(f)(x) = \frac{1}{\alpha} \int_{\mathbb{R}^d} \langle \nabla(f)(x+u) - \nabla(f)(x)  ; u\rangle \nu_\alpha(du) = \frac{1}{\alpha} \sum_{k=1}^d D^{\alpha-1}_k \partial_k (f)(x),
\end{align}
where $\partial_k$ is the first order partial derivative in the direction $e_k$ and $D^{\alpha-1}$ is the fractional gradient operator defined, for all $f \in \mathcal{C}_c^{\infty}(\mathbb{R}^d)$ and all $x \in \mathbb{R}^d$, by 
\begin{align}\label{eq:perim_fracGrad}
D^{\alpha-1}(f)(x)=\int_{\mathbb{R}^d} (f(x+u)-f(x)) u \nu_\alpha(du). 
\end{align}
From the standpoint of the classical Sobolev inequality \eqref{ineq:standard_Sobolev_inequality_p2}, it is rather natural to expect a fractional Sobolev inequality which puts into play the fractional gradient operator $D^{\alpha-1}$ instead of $(-\mathcal{A}_\alpha)^{1/2}$ (up to a change of Lebesgue exponent). This is the purpose of our first main result. 

\begin{thm}\label{thm:FSI_NDS_full}
Let $d \geq 2$ be an integer, let $\alpha \in (1,2)$ and let $\nu_\alpha$ be a non-degenerate symmetric $\alpha$-stable L\'evy measure on $\mathbb{R}^d$. Then, for all $p \in (1,d/(\alpha-1))$ and all $f \in \mathcal{C}^{\infty}_c(\mathbb{R}^d)$, 
\begin{align}\label{ineq:FSI_NDS_full}
\|f\|_{L^{p_\alpha^*}(\mathbb{R}^d,dx)} \leq C_{\alpha,p,d} \|D^{\alpha-1}(f)\|_{L^p(\mathbb{R}^d,dx)}, 
\end{align}
where $C_{\alpha,p,d}>0$ depends on $\alpha$, $p$ and $d$ and $p_\alpha^* = pd/(d-p(\alpha-1))$. 
\end{thm}
\noindent
The proof of Theorem \ref{thm:FSI_NDS_full} is deferred to Section \ref{sec:Sobolev_embeddings}.~It is based on ultracontractive estimates, interpolation arguments and the functional calculus of symmetric contractive Markovian semigroups in its Fourier analytic form stemming from the convolution structure of $(P^\alpha_t)_{t \geq 0}$. This method of proof is taken from \cite{Auscher_book07} (see also \cite[Theorem $1.2.$]{CM_tams93}) and adapted to our fractional setting. Moreover, an important step in the proof of Theorem \ref{thm:FSI_NDS_full} is based on a result (see Lemma \ref{lem:Lp_bound_FractionalPower} below) which ensures that the following inequality holds true: for all $p \in (1,+\infty)$ and all $f \in \mathcal{C}^{\infty}_c(\mathbb{R}^d)$, 
\begin{align}\label{ineq:upper_bound_Frac_Riesz}
\left\| \left(- \mathcal{A}_\alpha\right)^{\frac{\alpha-1}{\alpha}}(f) \right\|_{L^p(\mathbb{R}^d,dx)} \leq C_{p,\alpha,d} \|D^{\alpha-1}(f)\|_{L^p(\mathbb{R}^d,dx)},
\end{align}
where $C_{p,\alpha,d}>0$ depends on $p$, $\alpha$ and $d$ and is made explicit in Lemma \ref{lem:Lp_bound_FractionalPower}. Indeed, Inequality \eqref{ineq:upper_bound_Frac_Riesz} is the fractional analogue of one side of the equivalence of norms used to infer \eqref{ineq:standard_Sobolev_inequality_p2} from \eqref{ineq:standard_Sobolev-type_inequality_p2}.~The continuity properties of the associated Riesz transforms were investigated in \cite{AH20_4} using Bismut-type representation formulas. 

Regarding the geometric regime ($p = 1$), it is classical since the work \cite{VSC_book92} to prove isoperimetric inequalities using semigroup tools. Indeed, this is the approach followed in \cite{MLedoux_94} to prove the Euclidean isoperimetric inequality using the ultracontractivity of $(P^H_t)_{t \geq 0}$ together with the Bismut-type representation formula: for all $t>0$, all $f \in \mathcal{C}_c^{\infty}(\mathbb{R}^d)$ and all $x \in \mathbb{R}^d$, 
\begin{align}\label{eq:Gaussian_Bismut_formula}
\nabla P^H_t(f)(x) = \frac{1}{\sqrt{2t}} \int_{\mathbb{R}^d} y f\left(x + \sqrt{2t}y\right) \gamma(dy). 
\end{align}
In the non-degenerate symmetric $\alpha$-stable case, Bismut-type formulas have been obtained in \cite[Proposition $2.1$ and Lemma $3.7$]{AH20_4} for the actions of $\nabla$ and $D^{\alpha-1}$ on $(P_t^{\alpha})_{t \geq 0}$. As expected, they lead to the following versions of the pseudo-Poincar\'e inequality. Below, for $p \in [1,+\infty)$, $L^p(\mu_\alpha)$ denotes the Lebesgue space over the probability space $(\mathbb{R}^d, \mathcal{B}(\mathbb{R}^d),\mu_\alpha)$ with integrability exponent $p$ and endowed with the usual norm $\|\cdot\|_{L^p(\mu_\alpha)}$.

\begin{prop}\label{prop:pseudo_poincare}
Let $\alpha \in (1, 2)$, let $\nu_\alpha$ be a non-degenerate symmetric $\alpha$-stable L\'evy measure on $\mathbb{R}^d$, $d \geq 1$, with spectral measure $\lambda_1$. Let $\mu_\alpha$ be the associated $\alpha$-stable probability measure on $\mathbb{R}^d$, let $p_\alpha$ be its positive Lebesgue density and let us assume that, for all $p \in [1 , +\infty)$, 
\begin{align}\label{ineq:condint_logarithmic_derivative}
\left\| \frac{\nabla(p_\alpha)}{p_\alpha} \right\|_{L^p(\mu_\alpha)} < +\infty. 
\end{align}
Then, for all $p \in [1 , +\infty)$, all $t >0$ and all $f \in \mathcal{C}_c^{\infty}(\mathbb{R}^d)$, 
\begin{align}\label{ineq:pseudo_poincar_frac_gradient}
\left\| P_t^\alpha(f) - f \right\|_{L^p(\mathbb{R}^d,dx)} \leq \frac{t^{1 - \frac{1}{\alpha}}}{\alpha-1} \left\| \frac{\nabla(p_\alpha)}{p_\alpha} \right\|_{L^p(\mu_\alpha)}  \left\| D^{\alpha-1}(f) \right\|_{L^p(\mathbb{R}^d,dx)}.
\end{align}
Moreover, for all $ p \in [1, \alpha)$, all $t >0$ and all $f \in \mathcal{C}_c^{\infty}(\mathbb{R}^d)$, 
\begin{align}\label{ineq:pseudo_poincar_gradient}
\left\| P_t^\alpha(f) - f \right\|_{L^p(\mathbb{R}^d,dx)} \leq t^{\frac{1}{\alpha}} \left(\bbe |Y_\alpha|^p \right)^{\frac{1}{p}} \| \sigma_\alpha \left(\nabla(f)\right) \|_{L^p(\mathbb{R}^d,dx)},
\end{align}
where $Y_\alpha$ is an $\alpha$-stable random variable with characteristic function given by $\bbe \exp\left(i \xi Y_\alpha \right) = \exp \left( - |\xi|^\alpha\right)$, for all $\xi \in \bbr$.
\end{prop}
\noindent
The proof of Proposition \ref{prop:pseudo_poincare} is presented in Section \ref{sec:not_prelim}. As previously mentioned, it relies on the Bismut-type formulas proved in \cite{AH20_4} (recalled in \eqref{eq:bismut_type_gradient} and \eqref{eq:bismut_type_frac_gradient} below) together with duality arguments. The pseudo-Poincar\'e inequalities \eqref{ineq:pseudo_poincar_frac_gradient} and \eqref{ineq:pseudo_poincar_gradient} for $p = 1$ lead us to introduce two definitions of perimeter for which isoperimetric inequalities can be studied: let $\mathcal{P}_{\operatorname{cl}}$ and $\mathcal{P}_{\operatorname{frac}}$ be defined, for all $E$ Borel measurable subset of $\mathbb{R}^d$ with finite $d$-dimensional Lebesgue measure, by
\begin{align}\label{eq:classical}
\mathcal{P}_{\operatorname{cl}}(E) =  \underset{t \rightarrow 0^+}{\lim} \| \sigma_{\alpha}\left(\nabla P^\alpha_t (\bbone_E) \right) \|_{L^1(\mathbb{R}^d,dx)},
\end{align}
\begin{align}\label{eq:perim_frac}
\mathcal{P}_{\operatorname{frac}}(E) =  \underset{t \rightarrow 0^+}{\lim} \| D^{\alpha-1} P^\alpha_t (\bbone_E)\|_{L^1(\mathbb{R}^d,dx)}.
\end{align}
Then, \eqref{ineq:pseudo_poincar_frac_gradient} and \eqref{ineq:pseudo_poincar_gradient} combined with the ultracontractivity of $(P^\alpha_t)_{t \geq 0}$ ensure the following isoperimetric-type inequalities. (Below and in the sequel, $\mathcal{L}_d(E)$ denotes the $d$-dimensional Lebesgue measure of the Borel set $E$ and $\|\cdot\|_{\infty}$ the supremum norm on $\mathbb{R}^d$.) 

\begin{thm}\label{thm:isoperimetric_type}
Let $\alpha \in (1,2)$ and let $\mu_\alpha$ be a non-degenerate symmetric $\alpha$-stable probability measure on $\mathbb{R}^d$, $d \geq 1$, with positive Lebesgue density $p_\alpha$. Then, for all $E\in \mathcal{B}(\mathbb{R}^d)$ with $\mathcal{L}_d(E)<+\infty$,
\begin{align}\label{ineq:perim_frac_isoperi}
 \mathcal{L}_d(E)^{\frac{d-(\alpha-1)}{d}}  \leq  C^1_{\alpha,d} \mathcal{P}_{\operatorname{frac}}(E), \quad  \mathcal{L}_d(E)^{\frac{d-1}{d}} \leq C^2_{\alpha,d} \mathcal{P}_{\operatorname{cl}}(E),
\end{align}
where $C^1_{\alpha,d}$ and $C^2_{\alpha,d}$ are given by 
\begin{align}\label{eq:c_1_alpha}
C^1_{\alpha,d} = \left(\frac{1}{2(\alpha-1)} (2d)^{\frac{\alpha-1}{d+ \alpha-1}} + \left(\frac{1}{2d}\right)^{\frac{d}{d+\alpha-1}}\right)^{\frac{d+\alpha-1}{d}} 2^{\frac{\alpha-1}{\alpha}} \|p_\alpha\|_{\infty}^{\frac{\alpha-1}{d}} \|\nabla(p_\alpha)\|_{L^1(\mathbb{R}^d,dx)},
\end{align}
\begin{align}\label{eq:c_2_alpha}
C^2_{\alpha,d} = \left( \frac{1}{2} (2d)^{\frac{1}{d+1}}+\frac{1}{(2d)^{\frac{d}{d+1}}}\right)^{\frac{d+1}{d}}2^{\frac{1}{\alpha}} \|p_\alpha\|_{\infty}^{\frac{1}{d}} \bbe |Y_\alpha|,
\end{align}
and where $Y_\alpha$ is an $\alpha$-stable random variable with characteristic function given by $\bbe \exp\left(i \xi Y_\alpha\right) = \exp \left(- |\xi|^\alpha\right)$, for all $\xi \in \bbr$.
\end{thm}
\noindent
The proof of Theorem \ref{thm:isoperimetric_type} is presented in Section \ref{sec:frac_class_iso_ineq} and relies on semigroup arguments and it fulfills the first main objective of our work. Once the Sobolev and isoperimetric inequalities associated with the semigroup $(P_t^\alpha)_{t \geq 0}$ have been identified, the next step is to investigate the optimal versions of these inequalities as well as the corresponding optimization problems. This is \textit{the second main objective} of these notes. For this purpose,   
it is necessary to set an appropriate functional framework in which the solutions to the optimization problems belong. This leads us to introduce and study homogeneous fractional Sobolev spaces based on the fractional gradient operator $D^{\alpha-1}$. However, at this point, $D^{\alpha-1}$ is only well-defined on functions belonging to $\mathcal{S}(\mathbb{R}^d)$, the Schwartz space of infinitely differentiable functions which, with their derivatives of any order, are rapidly decreasing. Moreover, since $D^{\alpha-1}$ is a fractional operator, one cannot proceed by duality to define $D^{\alpha-1}$ on the classical Lebesgue spaces. To overcome this difficulty, we employ a truncation argument at the beginning of Section \ref{sec:anisotropic_fractional_Sobolev_spaces} in order to define $D^{\alpha-1}(f)$, for all $p \in [1,+\infty)$ and all $f \in L^p(\mathbb{R}^d,dx)$, as a tempered distribution. Then, for $p \in (1,d/(\alpha-1))$, the homogeneous fractional Sobolev space is defined by
\begin{align}\label{def:Lp-version_homo_frac_Sob_space}
\dot{W}^{\alpha-1,p}\big(\mathbb{R}^d,dx\big) := \{f \in L^{p_\alpha^*}(\mathbb{R}^d,dx) : \, |f|_{\alpha-1,p}<+\infty\}, 
\end{align}
where $p_\alpha^* = pd/(d-p(\alpha-1))$ and $|\cdot|_{\alpha-1,p}$ is given, for all $f \in L^{p_\alpha^*}\left(\mathbb{R}^d,dx\right)$, by 
\begin{align}\label{eq:seminorm_Lp}
|f|_{\alpha-1,p} = \|D^{\alpha-1}(f)\|_{L^p(\mathbb{R}^d,dx)}. 
\end{align}
One of the main results of these notes is the following Meyers-Serrin identification of $\dot{W}^{\alpha-1,p}\left(\mathbb{R}^d,dx\right)$ as the completion of $\mathcal{C}^{\infty}_c(\mathbb{R}^d)$ with respect to the norm $|\cdot|_{\alpha-1,p}$ denoted by $\mathcal{D}^{\alpha-1,p}(\mathbb{R}^d,dx)$.

\begin{thm}\label{thm:Meyers-Serrin_theorem}
Let $d \geq 2$ be an integer, let $\alpha \in (1,2)$ and let $p \in (1,d/(\alpha-1))$. Let $\mu_\alpha$ be a non-degenerate symmetric $\alpha$-stable probability measure on $\mathbb{R}^d$ with L\'evy measure $\nu_\alpha$. Then, 
\begin{enumerate}
\item $(\dot{W}^{\alpha-1,p}(\mathbb{R}^d,dx),|\cdot|_{\alpha-1,p})$ is a normed vector space;
\item $\mathcal{D}^{\alpha-1,p}(\mathbb{R}^d,dx)=\dot{W}^{\alpha-1,p}(\mathbb{R}^d,dx)$. 
\end{enumerate}
\end{thm}
\noindent
Theorem \ref{thm:Meyers-Serrin_theorem} is proved in Section \ref{sec:Sobolev_embeddings}. Its proof relies on an approximation result by smooth and compactly supported functions which is Theorem \ref{thm:approx_scs_homo_Lp} below. Since the operator $D^{\alpha-1}$ commutes with translations, the approximation by smooth functions is performed by convolution with a sequence of standard mollifiers as it is classical in the theory of Sobolev spaces. Due to the non-locality of $D^{\alpha-1}$ and the critical integrability of functions in $\dot{W}^{\alpha-1,p}(\mathbb{R}^d,dx)$, the smooth truncation procedure requires a fine analysis of the fractional Leibniz rule satisfied by $D^{\alpha-1}$ (see \eqref{eq:frac_Leibniz} below) combined with a convexity trick based on the Banach-Saks property for $L^p$-spaces with $p \in (1,d/(\alpha-1))$. The method of proof shares similarities with the corresponding one regarding the homogeneous fractional Sobolev spaces based on the Gagliardo-Slobodeckij seminorms. This case has been recently investigated in \cite[Theorem $3.1$]{BG-CV_cv21}. 

Armed with this functional framework, it is then possible to properly set the optimization problem associated with the fractional Sobolev inequality \eqref{ineq:FSI_NDS_full}. Let $p \in (1, d/(\alpha-1))$ and $p_\alpha^* = pd/(d-p(\alpha-1))$. Let $S_{p,\alpha,d}$ be the non-negative constant defined by 
\begin{align}\label{eq:full_minimization_problem_frac_sobolev-p}
S_{p,\alpha,d} := \underset{f \in \dot{W}^{\alpha-1,p}(\mathbb{R}^d,dx)}{\inf} \{\|D^{\alpha-1}(f)\|^p_{L^p(\mathbb{R}^d,dx)} : \quad \|f\|_{L^{p^*_\alpha}\left(\mathbb{R}^d,dx\right)} = 1\}. 
\end{align}
Then, when studying \eqref{eq:full_minimization_problem_frac_sobolev-p}, one wants to establish the existence of non-trivial optimizers and compute the constant $S_{p,\alpha,d}$ with respect to the parameters. In Section \ref{sec:opt_FSI}, we perform the first part of this program for $p = 2$ and for $p \in (1,d/(\alpha-1))$. For the Hilbertian case, we adapt the method of the missing
mass, initially developed in \cite{Lieb_AM83}, to our fractional setting.   

\begin{thm}\label{thm:existence_non-trivial_optimizer}
Let $d \geq 2$ be an integer, let $\alpha \in (1,2)$ and let $2_\alpha^* = 2d/(d-2(\alpha-1))$.~Let $\mu_\alpha$ be a symmetric and non-degenerate $\alpha$-stable probability measure on $\mathbb{R}^d$ with positive Lebesgue density $p_\alpha$ satisfying \eqref{ineq:condint_logarithmic_derivative} with $p=2$.~Let $\dot{W}^{\alpha-1,2}\left(\mathbb{R}^d,dx\right)$ be the homogeneous fractional Sobolev space with $p=2$. Then, there exists $f \in \dot{W}^{\alpha-1,2}(\mathbb{R}^d,dx) \setminus \{0\}$ such that
\begin{align*}
S_{2,\alpha,d} = \dfrac{\|D^{\alpha-1}(f)\|^2_{L^2(\mathbb{R}^d,dx)}}{\|f\|^2_{L^{2_\alpha^*}\left(\mathbb{R}^d,dx\right)}}. 
\end{align*}   
\end{thm}
\noindent
The proof of Theorem \ref{thm:existence_non-trivial_optimizer} is deferred to Section \ref{sec:opt_FSI}. It is based on a refined fractional Sobolev inequality obtained in Lemma \ref{lem:extension_RFSI} in order to circumvent the lack of non-trivial compactness of minimizing sequences for the problem \eqref{eq:full_minimization_problem_frac_sobolev-p}. This method is by now classical and has been implemented in different settings (see, e.g., \cite{G_esaimcocv,FL_AM12,PP_14,Frank_24} and the references therein). Moreover, to use the Br\'ezis-Lieb lemma (see \cite[Theorem $1.9.$]{Lieb_Loss_book01}), we prove a technical result (Lemma \ref{lem:almost_everywhere_convergence}) which allows to get strong local convergence in $L^2(\mathbb{R}^d,dx)$ for minimizing sequences in $\dot{W}^{\alpha-1,2}(\mathbb{R}^d,dx)$. This lemma relies on the pseudo-Poincar\'e inequality \eqref{ineq:pseudo_poincar_frac_gradient} of Proposition \ref{prop:pseudo_poincare} and requires the condition \eqref{ineq:condint_logarithmic_derivative} to hold for $p = 2$. 

For the general case $p \in (1,d/(\alpha-1))$, the lack of an Hilbertian structure when $p \ne 2$ makes us change our strategy to prove the existence of non-trivial optimizers for \eqref{eq:full_minimization_problem_frac_sobolev-p}. More precisely, we combine a refined fractional Sobolev inequality (see Theorem \ref{thm:refined_frac_Sobolev_Lp} below) together with a fractional concentration-compactness principle (see Theorem \ref{thm:ccp_homo_frac_p} below) to reach existence. Moreover, to derive the refined fractional Sobolev inequality, we employ the method developed in \cite{PP_14} which relies on a potential theory pointwise formula and weighted inequalities for the Riesz potential operator (see \eqref{eq:def_Riesz_Potential_Operator} for a definition) due to \cite{SaWhe_AJM92}. In particular, for the pointwise formula in our setting, we need to assume that the L\'evy measure $\nu_\alpha$ is a $\gamma$-measure with $\gamma$ large enough. The definition of $\gamma$-measure has been introduced and studied in \cite{BS_SM07,Szt_MN10,BSK_DM20} in order to model the effects of anisotropy on pointwise estimates for the potential kernel of $(P_t^{\alpha})_{t \geq 0}$ (see \eqref{eq:potential_kernel_stable_measure} below) and its derivatives. Recall that a Borel measure $\nu$ is a $\gamma$-measure on $\mathbb{S}^{d-1}$, with $\gamma \in [1,d]$, if there exists a positive constant $c$ such that, for all $x \in \mathbb{R}^d$ with $\|x\| = 1$ and all $r \in (0,1/2)$, 
\begin{align}\label{eq:gamma-measure}
\nu\left(B(x,r)\right) \leq cr^\gamma,
\end{align}
where $B(x,r)$ is the Euclidean ball centered at $x$ with radius $r$. 

\begin{thm}\label{thm:existence_non-trivial_optimizer_p}
Let $d \geq 2$ be an integer, let $\alpha \in (1,2)$, let $p \in (1,d/(\alpha-1))$ and let $p_\alpha^* = pd/(d-p(\alpha-1))$.~~Let $\mu_\alpha$ be a symmetric and non-degenerate $\alpha$-stable probability measure on $\mathbb{R}^d$ with positive Lebesgue density $p_\alpha$ satisfying \eqref{ineq:condint_logarithmic_derivative}. Let the L\'evy measure $\nu_\alpha$ be a $\gamma$-measure with $\gamma \in [1,d]$ and $\gamma - d + 2\alpha > 1$. Then, there exists $f \in \dot{W}^{\alpha-1,p}(\mathbb{R}^d,dx) \setminus \{0\}$ such that
\begin{align*}
S_{p,\alpha,d} = \dfrac{\|D^{\alpha-1}(f)\|^p_{L^p(\mathbb{R}^d,dx)}}{\|f\|^p_{L^{p_\alpha^*}\left(\mathbb{R}^d,dx\right)}}.
\end{align*}
\end{thm}
\noindent
The proof of Theorem \ref{thm:existence_non-trivial_optimizer_p} is done in Section \ref{sec:opt_FSI}. In the terminology of \cite{Lions_RMI185,Lions_RMI285}, the refined fractional Sobolev inequality \eqref{eq:refined_frac_Sobolev_Lp} together with Lemma \ref{lem:there_is_something_somewhere_non-trivial_p} and Lemma \ref{lem:almost_everywhere_convergence_p} allow to get rid of vanishing along renormalized minimizing sequences for $\eqref{eq:full_minimization_problem_frac_sobolev-p}$. Finally, the reverse H\"older inequalities \eqref{ineq:rH_atomic} and \eqref{ineq:reverse_Holder_inf} of the concentration-compactness principle combined with a convexity trick allow to reach the conclusion of Theorem \ref{thm:existence_non-trivial_optimizer_p}. This method of proof shares similarities with the ones developed in \cite{HR_aihpc07,BSS_NoDEA18,CKW_NA23}. 

\subsection{Literature review on the Riesz fractional gradient}
Of particular interest is the case where the underlying non-degenerate and symmetric $\alpha$-stable probability measure is invariant under rotations. Let $\sigma_L$ denote the spherical part of the $d$-dimensional Lebesgue measure and let $\nu_{\alpha}^{\operatorname{rot}}$ be the L\'evy measure on $\mathbb{R}^d$ with polar decomposition 
\begin{align}\label{eq:Levy_Rot}
\nu_{\alpha}^{\operatorname{rot}}(du) = c_{\alpha,d} \bbone_{(0,+\infty)}(r) \bbone_{\mathbb{S}^{d-1}}(y)\dfrac{dr}{r^{\alpha+1}}\sigma_L(dy),
\end{align}
where
\begin{align}\label{eq:renorm}
c_{\alpha,d} = \dfrac{-\alpha (\alpha-1) \Gamma\left(\frac{\alpha+d}{2}\right)}{4 \cos \left(\frac{\alpha\pi}{2}\right)\Gamma\left(\frac{\alpha+1}{2}\right) \pi^{\frac{d-1}{2}}\Gamma(2-\alpha)}.
\end{align}
Let $\mu_{\alpha}^{\operatorname{rot}}$ be the rotationally invariant $\alpha$-stable probability measure on $\mathbb{R}^d$ with L\'evy measure given by \eqref{eq:Levy_Rot} and Fourier transform given, for all $\xi \in \mathbb{R}^d$, by
\begin{align}\label{eq:charac_rot}
\widehat{\mu_{\alpha}^{\operatorname{rot}}}\left(\xi\right) = \exp\left(-\frac{\|\xi\|^\alpha}{2}\right).
\end{align}
The measure $\mu_{\alpha}^{\operatorname{rot}}$ is absolutely continuous with respect to the $d$-dimensional Lebesgue measure and its Lebesgue density, denoted by $p_{\alpha}^{\operatorname{rot}}$, is infinitely differentiable and such that, for all $x\in\mathbb{R}^d$, 
\begin{align*}
\frac{C_2}{\left(1+\|x\|\right)^{\alpha+d}} \leq p_{\alpha}^{\operatorname{rot}}(x) \leq \dfrac{C_1}{\left(1+\|x\|\right)^{\alpha+d}},
\end{align*}
for some constants $C_1,C_2>0$ only depending on $\alpha$ and $d$. In this case, the operator $D^{\alpha-1}$ is named \textit{the Riesz fractional gradient} and denoted by $D^{\alpha-1, \operatorname{rot}}$. In the last decade, it has been the subject of several studies starting with \cite{Shieh_Spector_ACV15,SSS_17,Shieh_Spector_ACV18}. Indeed, in these references, $D^{\alpha-1, \operatorname{rot}}$ is used as the canonical object to develop a theory of fractional Sobolev spaces and of fractional partial differential equations as well as a fractional calculus of variations. In particular, in \cite{Shieh_Spector_ACV15}, the non-homogeneous fractional Sobolev spaces based on $D^{\alpha-1,\operatorname{rot}}$ with integrability exponent $p \in (1,+\infty)$ are introduced as the completion of $\mathcal{C}^{\infty}_c(\mathbb{R}^d)$ with respect to the norm 
\begin{align}\label{eq:non-homo_frac-norm_rotinv}
\|f\|_{\alpha-1,p,\operatorname{rot}} = \left(\|f\|^p_{L^p(\mathbb{R}^d,dx)} + \sum_{k = 1}^d \|D_k^{\alpha-1,\operatorname{rot}}(f)\|^p_{L^p(\mathbb{R}^d,dx)}\right)^{\frac{1}{p}}, \quad f \in \mathcal{C}^{\infty}_c(\mathbb{R}^d). 
\end{align}
In \cite[Theorem $1.7$]{Shieh_Spector_ACV15}, these fractional Sobolev spaces are identified with the classical Bessel potential spaces of order $\alpha-1$ denoted by $L^{\alpha-1,p}\left(\mathbb{R}^d\right)$ (see Remark \ref{rem:fractional_Sobolev-type_space_order1} below). Moreover, in \cite{comi_stefani,BCGS_CRM22}, the following distributional-type definition of the non-homogeneous fractional Sobolev spaces is introduced: for $p \in [1, +\infty)$, 
\begin{align}
W^{\alpha-1,p}_{\operatorname{rot}}(\mathbb{R}^d,dx) := \left\{f \in L^p(\mathbb{R}^d,dx) : \, \|D^{\alpha-1,\operatorname{rot}}(f)\|_{L^p(\mathbb{R}^d,dx)}<+\infty\right\},
\end{align}
where $D^{\alpha-1, \operatorname{rot}}(f)$ is defined in the distributional sense, for $f \in L^p(\mathbb{R}^d,dx)$. In particular, it is proved in \cite[Theorem $3.23$]{comi_stefani} for $p = 1$ and in \cite[Appendix A, Theorem $25$]{BCGS_CRM22} for $p \in (1, +\infty)$, that $W^{\alpha-1,p}_{\operatorname{rot}}(\mathbb{R}^d,dx)$ coincides with the non-homogeneous fractional Sobolev space obtained by the completion procedure and, consequently, with the classical Bessel potential space $L^{\alpha-1,p}(\mathbb{R}^d)$.   

In our setting and thanks to the truncation procedure already mentioned, the non-homogeneous fractional Sobolev spaces based on $D^{\alpha-1}$ are introduced in Section \ref{sec:anisotropic_fractional_Sobolev_spaces} (see Definition \ref{def:non-homogeneous_frac_Sobolev_spaces} and Remark \ref{rem:fractional_Sobolev-type_space_order1} below). In Proposition \ref{prop:density_compactsupport}, we prove that these two spaces coincide in this general situation. The proof is based on the Young inequality and the interpolation inequality of Lemma \ref{lem:interpolation_inequality} as in the rotationally invariant case. Moreover, a general theory of Bessel potential-type spaces has been developed in \cite[Chapter $3$, Section $3.3$]{NJ02_2} for convolution semigroups associated with continuous negative definite functions. In particular, the results in there apply to the stable heat semigroup $(P^{\alpha}_t)_{t \geq 0}$ associated with $\psi_\alpha$ given, for all $\xi \in \mathbb{R}^d$, by
\begin{align}\label{eq:multiplier_fractional_laplacian}
 \psi_\alpha(\xi) = - \int_{\mathbb{S}^{d-1}} |\langle \xi ; y \rangle|^\alpha  \lambda_1(dy).
\end{align} 
For $p \in (1,+\infty)$ and $s>0$, let $H_p^{\psi_\alpha,s}(\mathbb{R}^d)$ be the Bessel potential-type space defined by 
\begin{align}\label{eq:def_fractional_BPT_space_orders}
H_p^{\psi_\alpha,s}(\mathbb{R}^d) := \left\{ f\in L^p(\mathbb{R}^d,dx) : \quad \|f\|_{\psi_\alpha,s,p} : = \| \left(E - \mathcal{A}_{\alpha,p}\right)^{\frac{s}{2}}(f) \|_{L^p(\mathbb{R}^d,dx)}<+\infty\right\}, 
\end{align} 
where $\mathcal{A}_{\alpha,p}$ is the $L^p(\mathbb{R}^d,dx)$-generator of $(P^{\alpha}_t)_{t \geq 0}$. The set $H_p^{\psi_\alpha,s}(\mathbb{R}^d)$ is the $L^p(\mathbb{R}^d,dx)$-domain of $\left(E-\mathcal{A}_{\alpha,p}\right)^{s/2}$. In Proposition \ref{prop:converse_inclusion}, we complete the picture by proving the following Calder\'on-type result: for all $p \in (1, +\infty)$, the non-homogeneous fractional Sobolev space $W^{\alpha-1,p}(\mathbb{R}^d,dx)$ considered in Definition \ref{def:non-homogeneous_frac_Sobolev_spaces} coincides with $H_p^{\psi_\alpha,r(\alpha)}(\mathbb{R}^d)$ where $r(\alpha) = 2(\alpha-1)/\alpha$. The proof is partly based on the Riesz transform-type results of \cite{AH20_4}. 

Returning to the literature on the Riesz fractional gradient, a fractional Sobolev inequality is given in \cite[Theorem $1.8$]{Shieh_Spector_ACV15} for $p \in (1,d/(\alpha-1))$ with $D^{\alpha-1, \operatorname{rot}}$. The proof is based on \cite[Theorem $1.12$]{Shieh_Spector_ACV15} which can be understood as a fractional version of the fundamental theorem of calculus. In Section \ref{sec:not_prelim}, Lemma \ref{lem:l_alpha_ball} and Proposition \ref{prop:FFTC_L1_NDSstrong} provide extensions of this fractional fundamental theorem of calculus when the underlying L\'evy measure is given as in \eqref{eq:LevyIndAxes} below, or when $\nu_\alpha$ is a $\gamma$-measure in the sense of \eqref{eq:gamma-measure} with $\gamma \in [1,d]$ such that $\gamma - d + 2\alpha > 1$. Moreover, based on these two results, we obtain new proofs of the fractional Sobolev inequality \eqref{ineq:FSI_NDS_full} in Proposition \ref{prop:well-defined_continuity} and Proposition \ref{prop:SI_NDSstrong} for these two situations. This is reminiscent of the proof provided by S.~L.~Sobolev in \cite{Sobolev_MS38} for the classical Sobolev inequality. 

As for the geometric regime ($p=1$), the fractional Sobolev inequality with $D^{\alpha-1, \operatorname{rot}}$ has been proved in \cite[Theorem A']{SSS_17}. This inequality is used in turn to prove the fractional isoperimetric inequality with $D^{\alpha-1, \operatorname{rot}}$ in \cite[Theorem $4.4$]{comi_stefani} for Borel subsets of $\mathbb{R}^d$ with finite Lebesgue measure and finite fractional perimeter. In \cite[Definition $4.1$]{comi_stefani}, the notion of fractional perimeter is defined through a duality formula based on a fractional divergence operator designed from the definition of $D^{\alpha-1,\operatorname{rot}}$. Moreover, by analogy with the classical bounded variation space, a theory of fractional bounded variation space is developed in \cite[Section $3$]{comi_stefani}. In the Appendix, Proposition \ref{prop:variational_representation_fractional_perimeter} ensures that the semigroup definition of the fractional perimeter \eqref{eq:perim_frac} and the variational definition introduced in \cite{comi_stefani} coincide for the rotationally invariant case. Actually, using the fractional divergence operator introduced in Definition \ref{defi:fractional_divergence} below and Proposition \ref{prop:variational_representation_fractional_perimeter}, the two definitions of fractional perimeter do coincide for all non-degenerate and symmetric $\alpha$-stable L\'evy measures. The proof relies on a structure result for the associated fractional bounded variation spaces and the lower semicontinuity of the fractional variation norms. 

Finally, in \cite{comi_stefani_rmc22,BCGS_CRM22}, the asymptotics of the operator $D^{\alpha-1,\operatorname{rot}}$ is analyzed when $\alpha$ tends to $2$ from below or to $1$ from above in the spirit of the celebrated Bourgain-Br\'ezis-Mironescu and Maz’ya-Shaposhnikova results (see \cite{BBM_2001} and \cite{MS_02} respectively). More precisely, \cite[Theorem $4.11$]{comi_stefani_rmc22} ensures that as $\alpha$ tends to $2$, $D^{\alpha-1,\operatorname{rot}}(f)$ converges strongly in $L^p(\mathbb{R}^d,dx)$ to $\nabla(f)$, for all $p \in [1,+\infty)$ and all $f \in W^{1,p}(\mathbb{R}^d,dx)$, the classical Sobolev space of order $1$ and integrability exponent $p$. For the general non-degenerate and symmetric $\alpha$-stable situation, we obtain an analogous result in Theorem \ref{thm:prop:BBM_strong_convergence_NDS} after a suitable normalization. In particular, the limiting first order differential operator inherits the (an)isotropy of the underlying L\'evy measure $\nu_\alpha$. It is given, for all $f \in \mathcal{S}(\mathbb{R}^d)$ and all $x \in \mathbb{R}^d$, by 
\begin{align}\label{eq:anisotropic_local_gradient}
D_\sigma(f)(x) = \int_{\mathbb{S}^{d-1}} y \langle y ; \nabla(f)(x) \rangle \sigma(dy),
\end{align}
where $\sigma$ is the spherical component of $\nu_\alpha$. For the asymptotics when $\alpha$ tends to $1$ from above, \cite[Theorem $17$, (ii)]{BCGS_CRM22} ensures that $D^{\alpha-1,\operatorname{rot}}(f)$ converges strongly in $L^p(\mathbb{R}^d,dx)$ to $-\mathcal{R}(f)/2$, for all $p \in (1,+\infty)$ and all $f \in \cup_{\alpha \in (1,2)}W^{\alpha-1,p}_{\operatorname{rot}}(\mathbb{R}^d,dx)$. Here and in the sequel, $\mathcal{R}(f)$ denotes the classical vectorial Riesz transform of order $1$ with Fourier symbol $-i\xi /\| \xi\|$, for all $\xi \in \mathbb{R}^d \setminus \{0\}$. In particular, the proof of \cite[Theorem $17$, (ii)]{BCGS_CRM22} relies on sharp fractional interpolation inequalities which are stated in \cite[Theorem $13$, (i)]{BCGS_CRM22}. In Theorem \ref{thm:MS_strong_convergence_NDSstrong} below, we extend this strong convergence result to all non-degenerate and symmetric $\alpha$-stable L\'evy measures which are $\gamma$-measures in the sense of \eqref{eq:gamma-measure} with $\gamma > 2 + d - 2\alpha$. Moreover, the limiting operator which appears as $\alpha \rightarrow 1^+$ is a Fourier multiplier operator defined, for all $f \in \mathcal{C}^{\infty}_c(\mathbb{R}^d)$ and all $x \in \mathbb{R}^d$, by 
\begin{align}\label{eq:avRT}
\mathcal{R}_\sigma(f)(x) = \frac{1}{(2\pi)^d} \int_{\mathbb{R}^d} \mathcal{F}(f)(\xi) e^{i \langle x ; \xi \rangle } m_\sigma(\xi) d\xi,
\end{align}
where, for all $\xi \in \mathbb{R}^d \setminus \{0\}$, 
\begin{align}\label{eq:mulitplier_anisotropic_vRT}
m_\sigma(\xi) = i \frac{\pi}{2} \int_{\mathbb{S}^{d-1}} y \operatorname{sign}\left( \langle y ; \xi \rangle\right) \sigma(dy),
\end{align}
and $\sigma$ is the spherical component of $\nu_\alpha$. $\mathcal{R}_\sigma$ is an (an)isotropic extension of the usual vectorial Riesz transform (up to a multiplicative non-zero constant) and, using the method of rotations, it admits a bounded extension on the Lebesgue space $L^p(\mathbb{R}^d,dx)$, for all $p \in (1, +\infty)$ (see Lemma \ref{lem:Lp_continuity_avRT} below). The proof of Theorem \ref{thm:MS_strong_convergence_NDSstrong} is similar to the one of \cite[Theorem $17$, (ii)]{BCGS_CRM22} and also relies on sharp fractional interpolation inequalities which are derived in Proposition \ref{prop:fractional_interpolation_NDSstrong} in this general situation (see also respectively Proposition \ref{prop:fractional_interpolation_inequality2} and Proposition \ref{prop:fractional_interpolation_independent} for the rotationally invariant and the independent cases). The proof of these propositions is \textit{new} and based on several ingredients: the ``optimal" splitting point technique, the method of rotations to infer $L^p(\mathbb{R}^d,dx)$-boundedness from $L^p(\mathbb{R},dx)$-boundedness and the fractional versions of the fundamental theorem of calculus instanced in \eqref{eq:spatial_rep_comistefani}, \eqref{eq:formula_Stein_axes} and \eqref{eq:FFTC_NDS}.

\subsection{Side results}
Along the course of these notes, we obtain side results of independent interest which we now discuss. In Section \ref{sec:anisotropic_fractional_Sobolev_spaces}, we provide a fractional extension of order $\alpha$ of a well-known result due to A.~P.~Calder\'on (see \cite[Theorem $7$]{Cal_61} and also \cite[Chapter $V$, Section $3$, Theorem $3$]{Stein_1}). Indeed, Theorem \ref{thm:identification_domain_generator} ensures that the fractional Sobolev space of order $\alpha$ defined, for all $p \in (1,+\infty)$, by 
\begin{align*}
W^{\alpha,p}(\mathbb{R}^d,dx) = \{f \in L^p(\mathbb{R}^d,dx) :\forall (j,k) \in \{1, \dots, d\}^2,\, \partial_j(f), D^{\alpha-1}_k(f), \partial_j (D^{\alpha-1}_k(f)) \in L^p(\mathbb{R}^d,dx) \}
\end{align*} 
is equal to the $L^p(\mathbb{R}^d,dx)$-domain of the $L^p(\mathbb{R}^d,dx)$-generator of $(P^\alpha_t)_{t \geq 0}$. Combined with \cite[Theorem $3.3.11$]{NJ02_2}, this ensures the following identification: for all $p \in (1,+\infty)$, 
\begin{align}\label{eq:Calderon_order_alpha}
W^{\alpha,p}(\mathbb{R}^d,dx) = H_p^{\psi_\alpha,2}(\mathbb{R}^d).
\end{align} 

As a consequence of the sharp fractional interpolation inequality \eqref{ineq:fractional_interpolation_NDSstrong} of Proposition \ref{prop:fractional_interpolation_NDSstrong}, we obtain the following continuous embeddings: for all $p \in (1,+\infty)$ and all $1<\beta<\alpha<2$, 
\begin{align*}
W^{\alpha-1, p}\left(\mathbb{R}^d , dx\right) \hookrightarrow W^{\beta-1,p}\left(\mathbb{R}^d,dx\right) \Leftrightarrow H_p^{\psi_\alpha, r(\alpha)}(\mathbb{R}^d) \hookrightarrow H_p^{\psi_\beta, r(\beta)}(\mathbb{R}^d), 
\end{align*}
when the underlying L\'evy measures are $\gamma$-measures with $\gamma > d+2-2\alpha$. Moreover, from these continuous embeddings, Corollary \ref{cor:non-smooth_mutliplier_NDSstrong} provides $L^p(\mathbb{R}^d,dx)$-boundedness of Fourier multiplier operators whose symbols are defined, for all $\xi \in \mathbb{R}^d$, by
\begin{align*}
t_{\alpha,\beta}(\xi) = \dfrac{\left(1 - \psi_\beta(\xi)\right)^{\frac{\beta-1}{\beta}}}{\left(1 - \psi_\alpha(\xi) \right)^{\frac{\alpha-1}{\alpha}}},
\end{align*}   
where $\psi_\alpha$ is given by \eqref{eq:multiplier_fractional_laplacian}. Note that, in the general non-degenerate and symmetric $\alpha$-stable case, the function $\psi_\alpha$ lacks smoothness. Finally, this type of Fourier multipliers linked to L\'evy processes have been studied in \cite{BB_jfa07,BBB_bcp11} but the results there do not apply to the ones given by $t_{\alpha,\beta}$.   

As already observed, at the endpoint $\alpha=2$, the fractional gradient operator $(2-\alpha)D^{\alpha-1}$ becomes an (an)isotropic local gradient operator $D_\sigma$ given by \eqref{eq:anisotropic_local_gradient}. In Proposition \ref{prop:AL_FTC}, we provide an (an)isotropic fundamental theorem of calculus which is a generalization of \cite[Chapter $V$, Section $2$, identity $(18)$]{Stein_1}. Combined with Lemma \ref{lem:from_iso_to_aniso} and the classical Hardy-Littlewood-Sobolev inequality, this version of the fundamental theorem of calculus provides the (an)isotropic Sobolev inequalities for certain scalar product norms in Theorem \ref{thm:sobolev-inequality-anisotropic} (see also \cite{cordero_nazaret_villani} for the sharp versions of these inequalities by completely different methods). Still in relation to Sobolev-type embeddings with general norms, we obtain in Theorem \ref{thm:PL_H_implies_Sobolev} a sharp Sobolev-type inequality in the Lorentz scale. In particular, this embedding implies the optimal anisotropic isoperimetric inequality (see Remark \ref{long_rem}, (i), below). The proof relies on the anisotropic P\'olya-Szeg\"o principle together with a Hardy-type inequality and Lemma \ref{lemma_technical_Lorentz_Hardy} of the Appendix.

Finally, in Section \ref{sec:frac_class_iso_ineq}, we point out a natural connection between convex geometry and $p_\alpha$, the Lebesgue density of the non-degenerate symmetric $\alpha$-stable probability measures on $\mathbb{R}^d$. Lemma \ref{lem:special_rep_palpha} presents a new integral representation formula for $p_\alpha$ in \eqref{eq:Generalized_Bessel_rep} which is a direct extension of the well-known formula in the rotationally invariant case based on the Fourier transform of radial functions (see, e.g., \cite[Appendix $B.4$ and $B.5$]{G08}). Equation \eqref{eq:Generalized_Bessel_rep} makes use    
of a function which is defined (up to some constant) as the Fourier transform of the cone probability measure naturally associated with $p_\alpha$ (see \eqref{def:cone_measure_K_alpha} below). In Proposition \ref{prop:limit_discrete_case} and Proposition \ref{prop:limit_general_case}, we further explore this connection having the sharp anisotropic isoperimetric inequalities in mind. More precisely, we compute certain limits in \eqref{eq:limit_discrete_case} and \eqref{eq:limit_general_case} which imply the small time first order asymptotic expansions of the heat contents associated with the independent case (see \eqref{eq:StableIndAxes} below) and the general non-degenerate symmetric $\alpha$-stable case for some specific convex subsets of $\mathbb{R}^d$.

\subsection{Organization}
\noindent
Let us now describe the content of these notes. In the next section, we introduce the remaining notations and definitions used throughout the manuscript and prove preliminary results on the pseudo-Poincar\'e inequalities associated with $(P_t^{\alpha})_{t \geq 0}$, composition formulas for the fractional gradient and divergence operators and the (an)isotropic (non)local versions of the fundamental theorem of calculus. In Section \ref{sec:anisotropic_fractional_Sobolev_spaces}, we develop the fractional Sobolev spaces theory based on the extension of the fractional gradient operator $D^{\alpha-1}$ and study the asymptotic theories as $\alpha$ tends to $2$ from below or to $1$ from above. Section \ref{sec:Sobolev_embeddings} deals with the weak and strong Sobolev-type embeddings naturally linked to $(P_t^{\alpha})_{t\geq 0}$. A particular emphazis is put on the refined fractional Sobolev inequalities and the homogeneous fractional Sobolev spaces. These tools are then used in Section \ref{sec:opt_FSI} to study the minimization problem \eqref{eq:full_minimization_problem_frac_sobolev-p} for $p=2$ and for general $p \in (1,d/(\alpha-1))$. In Section \ref{sec:frac_class_iso_ineq}, we study the geometric inequalities associated with $(P^{\alpha}_t)_{t \geq 0}$ and, in particular, the isoperimetric inequalities \eqref{ineq:perim_frac_isoperi}. Finally, the manuscript ends with an appendix section gathering technical results used throughout.

\section{Notations and preliminaries}\label{sec:not_prelim}
\noindent
The standard Gaussian probability measure on $\mathbb{R}^d$ denoted by $\gamma$ is defined through its Fourier transform, for all $\xi \in \mathbb{R}^d$, by
\begin{align}\label{def:gauss}
\widehat{\gamma}(\xi) := \int_{\mathbb{R}^d} e^{i \langle y ; \xi\rangle} \gamma(dy) =  \exp\left(- \frac{\|\xi\|^2}{2}\right).
\end{align}
Next, let $\nu_{\alpha,1}$ be the L\'evy measure on $\bbr$ given by
\begin{align}\label{sym:1}
\nu_{\alpha,1}(du) = c_{\alpha}\frac{du}{|u|^{\alpha+1}},
\end{align}
with
\begin{align}
c_{\alpha} = \left( \dfrac{-\alpha (\alpha -1)}{2 \Gamma(2-\alpha)\cos\left(\frac{\alpha \pi}{2}\right)} \right),
\end{align}
and let $\mu_{\alpha,1}$ be the $\alpha$-stable probability measure on $\bbr$ with L\'evy measure $\nu_{\alpha,1}$ and with Fourier transform defined, for all $\xi \in \bbr$, by
\begin{align}
\widehat{\mu_{\alpha,1}}(\xi) = \exp \left(\int_{\bbr} \left(e^{i \langle u;\xi\rangle}-1-i \langle u;\xi \rangle\right) \nu_{\alpha,1}(du)\right) = \exp\left(-|\xi|^{\alpha}\right).
\end{align}
Finally, let $\mu_{\alpha,d} = \mu_{\alpha,1} \otimes \dots \otimes \mu_{\alpha,1}$ be the product probability measure on $\mathbb{R}^d$ whose Fourier transform is given, for all $\xi \in \mathbb{R}^d$, by
\begin{align}\label{eq:StableIndAxes}
\widehat{\mu_{\alpha,d}}(\xi) = \prod_{k=1}^d \widehat{\mu_{\alpha,1}}(\xi_k) = \exp\left(\int_{\mathbb{R}^d} \left(e^{i \langle \xi;u \rangle}-1-i\langle \xi;u \rangle\right)\nu_{\alpha,d}(du)\right),
\end{align}
with
\begin{align}\label{eq:LevyIndAxes}
\nu_{\alpha,d}(du) = \sum_{k=1}^d \delta_0(du_1) \otimes \dots \otimes \delta_0(du_{k-1}) \otimes \nu_{\alpha,1}(du_k) \otimes \delta_0(du_{k+1}) \otimes \dots \otimes \delta_0(du_d),
\end{align}
where $\delta_0$ is the Dirac measure at $0$. In the sequel, $\mathcal{F}$ is the Fourier transform operator given, for all $f \in \mathcal{S}(\mathbb{R}^d)$ and all $\xi \in \mathbb{R}^d$, by
\begin{align*}
\mathcal{F}(f)(\xi) = \int_{\mathbb{R}^d} f(x) e^{- i \langle x; \xi \rangle} dx.
\end{align*}
On $\mathcal{S}(\mathbb{R}^d)$, the Fourier transform is an isomorphism and the following well-known inversion formula holds
\begin{align*}
f(x) = \frac{1}{(2\pi)^d} \int_{\mathbb{R}^d} \mathcal{F}(f)(\xi)e^{i \langle \xi ; x\rangle} d\xi,\quad x\in \mathbb{R}^d.
\end{align*}
For $p \in [1, +\infty)$, $L^p(\mathbb{R}^d,dx)$ is the classical Lebesgue space where the reference measure is the Lebesgue measure. It is endowed with the norm $\|\cdot\|_{L^p(\mathbb{R}^d,dx)}$ defined, for all suitable $f$, by
\begin{align*}
\|f\|_{L^p(\mathbb{R}^d,dx)} := \left(\int_{\mathbb{R}^d} |f(x)|^p dx\right)^{\frac{1}{p}}.
\end{align*}
In the sequel, $(D^{\alpha-1})^*$ and $\mathbf{D}^{\alpha-1}$ denote the fractional and non-local operators respectively defined, for all $f \in \mathcal{S}(\mathbb{R}^d)$ and all $x\in \mathbb{R}^d$, by
\begin{align}\label{eq:perim_fracGradDual}
(D^{\alpha-1})^*(f)(x) := \int_{\mathbb{R}^d} (f(x-u)-f(x)) u \nu_\alpha(du),
\end{align}
\begin{align}\label{eq:perim_fracGradWB}
\mathbf{D}^{\alpha-1}(f)(x) :=\frac{1}{2}\left(D^{\alpha-1}\left(f\right)(x)-(D^{\alpha-1})^*\left(f\right)(x)\right).
\end{align}
Note that, for all $f \in \mathcal{S}(\mathbb{R}^d)$ and all $p \in [1, +\infty]$, $D^{\alpha-1}(f) \in L^p(\mathbb{R}^d,dx)$, as easily seen as a straightforward consequence of the Minkowski integral inequality. Let us introduce as well a notion of fractional divergence operator acting, for example, on $\mathbb{R}^d$-valued infinitely differentiable functions with compact support.
\begin{defi}\label{defi:fractional_divergence}
Let $d \geq 1$ be an integer, let $\alpha \in (1,2)$, let $\nu_\alpha$ be a non-degenerate symmetric L\'evy measure on $\mathbb{R}^d$ verifying \eqref{eq:scale} and let $D^{\alpha-1}$ and $(D^{\alpha-1})^*$ be the operators defined by \eqref{eq:perim_fracGrad} and \eqref{eq:perim_fracGradDual} respectively. Then,~the fractional divergence operator $\operatorname{div}_\alpha$ is defined, for all $w = (w_1, \dots, w_d) \in \mathcal{C}_c^{\infty}(\mathbb{R}^d , \mathbb{R}^d)$ and all $x \in \mathbb{R}^d$, by
\begin{align}\label{eq:perim_fractional_divergence}
\operatorname{div}_\alpha(w)(x) = - \sum_{k = 1}^d (D^{\alpha-1}_k)^*(w_k)(x) = \sum_{k = 1}^d (D^{\alpha-1}_k)(w_k)(x).
\end{align}
Moreover, for all $w = (w_1, \dots, w_d) \in \mathcal{C}_c^{\infty}(\mathbb{R}^d , \mathbb{R}^d)$ and all $f \in \mathcal{S}(\mathbb{R}^d)$, 
\begin{align}\label{eq:ipp}
\int_{\mathbb{R}^d} \langle w(x) ; D^{\alpha-1}(f)(x) \rangle dx = - \int_{\mathbb{R}^d} \operatorname{div}_\alpha(w)(x) f(x)dx.
\end{align}
\end{defi}
\noindent
Let us now recall the definition of the gamma transform of order $r>0$.~For $\mathcal{C}_0$-semigroups of contractions on a Banach space, $(P_t)_{t\geq 0}$, with generator $\mathcal{A}$, the gamma transform of order $r>0$ is defined, for all suitable $f$, by
\begin{align}\label{eq:Gam_Tr}
\left(E-\mathcal{A} \right)^{-\frac{r}{2}} f = \frac{1}{\Gamma(\frac{r}{2})} \int_0^{+\infty} \dfrac{e^{-t}}{t^{1-\frac{r}{2}}} P_t(f) dt,
\end{align}
where $E$ is the identity operator and the integral on the right-hand side is understood in the Bochner sense. Moreover, for all $\lambda >0$, all $r>0$ and all suitable $f$, 
\begin{align*}
\left(\lambda E-\mathcal{A} \right)^{-\frac{r}{2}} f = \frac{1}{\Gamma(\frac{r}{2})} \int_0^{+\infty} \dfrac{e^{-\lambda t}}{t^{1-\frac{r}{2}}} P_t(f) dt.
\end{align*}
In particular, as $\lambda$ tends to $0^+$ and when the quantities make sense,
\begin{align*}
\left(-\mathcal{A} \right)^{-\frac{r}{2}} f = \frac{1}{\Gamma(\frac{r}{2})} \int_0^{+\infty} t^{\frac{r}{2}-1} P_t(f) dt.
\end{align*}
In the Fourier domain, $\mathcal{A}_\alpha$ is a multiplier operator which admits the following representation: for all $f \in \mathcal{S}(\mathbb{R}^d)$ and all $x \in \mathbb{R}^d$, 
\begin{align}\label{eq:Stheatgen_FourRep}
\mathcal{A}_\alpha(f)(x) = \frac{1}{(2 \pi)^d} \int_{\mathbb{R}^d} \mathcal{F}(f)(\xi) e^{i \langle \xi ; x \rangle} \psi_\alpha(\xi) d\xi, 
\end{align} 
where $\psi_\alpha$ is given by \eqref{eq:multiplier_fractional_laplacian}. Note that since $\sigma$ is symmetric, 
\begin{align*}
D^{\alpha-1} = -(D^{\alpha-1})^*, \quad (P_t^{\alpha})^* = P_t^\alpha , \quad (\mathcal{A}_\alpha)^* = \mathcal{A}_\alpha, 
\end{align*}
where $*$ here denotes the adjoint operation. Moreover, for all $f,g$ smooth enough on $\mathbb{R}^d$, 
\begin{align}\label{eq:ipp_symmetry}
\langle f ; (- \mathcal{A}_\alpha)(g) \rangle_{L^2(\mathbb{R}^d,dx)} = \frac{1}{2} \int_{\mathbb{R}^d} \int_{\mathbb{R}^d} (f(x+u) - f(x))(g(x+u) - g(x)) \nu_\alpha(du) dx,
\end{align}
and so it is natural to introduce the following gradient-length: for all $x \in \mathbb{R}^d$ and all $f \in \mathcal{S}(\mathbb{R}^d)$, 
\begin{align}\label{ineq:perim_frac_gradient_length}
\nabla_{\nu_\alpha}(f)(x) = \left(\int_{\mathbb{R}^d} |f(x+u)-f(x)|^2 \nu_\alpha(du)\right)^{\frac{1}{2}}.
\end{align}
Note that this gradient-length is linked to the ``carr\'e du champs" operator associated with $\mathcal{A}_\alpha$ and which is given, for all $f,g$ smooth enough on $\mathbb{R}^d$ and all $x \in \mathbb{R}^d$, by 
\begin{align*}
\Gamma_\alpha(f,g)(x)& = \frac{1}{2} \left(\mathcal{A}_{\alpha}(fg)(x) - f(x) \mathcal{A}_\alpha(g)(x) - g(x) \mathcal{A}_{\alpha}(f)(x)\right) \\
&= \frac{1}{2} \int_{\mathbb{R}^d} \left(f(x+u)-f(x)\right) \left(g(x+u) - g(x)\right)\nu_\alpha(du),
\end{align*}
so that $\Gamma_\alpha(f,f) = \nabla_{\nu_\alpha}(f)^2/2$. Next, let us recall the two Bismut-type formulas associated with $(P^\alpha_t)_{t\geq 0}$ and already put forward in \cite[Proposition $2.1$ and Lemma $3.7$]{AH20_4}: for all $f \in \mathcal{S}(\mathbb{R}^d)$, all $x \in \mathbb{R}^d$ and all $t >0$, 
\begin{align}\label{eq:bismut_type_gradient}
\nabla P^{\alpha}_t (f)(x) = \frac{1}{t^{\alpha}} \int_{\mathbb{R}^d} \dfrac{- \nabla(p_\alpha)(y)}{p_\alpha(y)} f \left(x + t^{\frac{1}{\alpha}} y \right) \mu_\alpha(dy),
\end{align}
and
\begin{align}\label{eq:bismut_type_frac_gradient}
D^{\alpha - 1} P^\alpha_t(f)(x) = \frac{1}{t^{1 - \frac{1}{\alpha}}} \int_{\mathbb{R}^d} y f\left(x + t^{\frac{1}{\alpha}} y\right) \mu_\alpha(dy).
\end{align}
As an application of the Minkowski integral inequality, for all $f$ smooth enough on $\mathbb{R}^d$, all $x \in \mathbb{R}^d$ and all $t>0$,
\begin{align}\label{ineq:sub_commutation}
\nabla_{\nu_\alpha}\left(P_t^\alpha(f)\right)(x) \leq P^\alpha_t\left(\nabla_{\nu_\alpha}(f)\right)(x).
\end{align}
Finally, the following lemma will be useful.

\begin{lem}\label{lem:local_reverse_poincar}
Let $\alpha \in (1,2)$ and let $(P^\alpha_t)_{t \geq 0}$ be the stable heat semigroup given by \eqref{eq:StheatSM}.  Then, for all $f \in \mathcal{S}(\mathbb{R}^d)$, all $t>0$ and all $x \in \mathbb{R}^d$, 
\begin{align}\label{ineq:local_reverse_poincar}
P^{\alpha}_t(f^2)(x) - (P^{\alpha}_t(f)(x))^2 \geq t (\nabla_{\nu_\alpha}\left(P_t^\alpha(f)\right)(x))^2. 
\end{align}
\end{lem}

\begin{proof}
The proof is very classical in the framework of Markov diffusion semigroups in the sense of \cite{BGL14}.  Let $t>0$ and let $f \in \mathcal{S}(\mathbb{R}^d)$. Then,  for all $x \in \mathbb{R}^d$, and using \eqref{ineq:sub_commutation},
\begin{align*}
P^{\alpha}_t(f^2)(x) - (P^{\alpha}_t(f)(x))^2 & = \int_0^t  \dfrac{d}{ds} \left(P_s^\alpha \left(P^\alpha_{t-s}(f)\right)^2\right)(x)ds \\
& = 2 \int_0^t P^\alpha_s \left( \Gamma_\alpha \left(P^\alpha_{t-s}(f) , P^\alpha_{t-s}(f)\right)  \right)(x) ds  \\
& = \int_0^t P^\alpha_s \left( |\nabla_{\nu_\alpha}\left(P^\alpha_{t-s}(f)\right)|^2 \right)(x) ds \\
& \geq \int_0^t \left(P^\alpha_s \left( \nabla_{\nu_\alpha}\left(P^\alpha_{t-s}(f)\right) \right)(x)\right)^2 ds \\
& \geq t (\nabla_{\nu_\alpha}\left(P_t^\alpha(f)\right)(x))^2. 
\end{align*}
\end{proof}
\noindent
Let us now prove the pseudo-Poincar\'e inequalities of Proposition \ref{prop:pseudo_poincare} for $(P^\alpha_t)_{t \geq 0}$ based on the decomposition \eqref{eq:Stheatgen_nice_decomposition}.\\

\noindent
\textit{Proof of Proposition \ref{prop:pseudo_poincare}.}
The proof is rather elementary and is based on a duality argument together with the Bismut-type formulas \eqref{eq:bismut_type_gradient} and \eqref{eq:bismut_type_frac_gradient}.~Let us start with the case $p \in (1, +\infty)$.~Let $f \in \mathcal{S}(\mathbb{R}^d)$ and let $g \in \mathcal{C}^{\infty}_c(\mathbb{R}^d)$.~Then, by symmetry and standard semigroup arguments, 
\begin{align*}
\langle P^\alpha_t(f) - f ; g \rangle_{L^2(\mathbb{R}^d,dx)} =  \langle f ; P_t^\alpha(g) - g \rangle_{L^2(\mathbb{R}^d,dx)} = \int_0^t \langle f ; \mathcal{A}_\alpha P^\alpha_s(g) \rangle_{L^2(\mathbb{R}^d,dx)} ds. 
\end{align*}
Thus, by duality and the Bismut-type formula \eqref{eq:bismut_type_gradient}, 
\begin{align*}
\langle P^\alpha_t(f) - f ; g \rangle_{L^2(\mathbb{R}^d,dx)} & = \frac{1}{\alpha}  \int_0^t \langle (D^{\alpha - 1})^*f ; \nabla P^\alpha_s(g) \rangle_{L^2(\mathbb{R}^d,dx)} ds \\
& = - \int_0^t \frac{ds}{\alpha s^{\frac{1}{\alpha}}} \int_{\mathbb{R}^d} \int_{\mathbb{R}^d} \langle (D^{\alpha - 1})^*f(x) ;  \dfrac{\nabla(p_\alpha)(y)}{p_\alpha(y)} \rangle g \left(x + s^{\frac{1}{\alpha}} y \right) \mu_\alpha(dy)dx.
\end{align*}
Now, the H\"older's inequality with $p \in (1, +\infty)$ and $q = p/(p-1)$ gives (with $L^p(\bbr^{2d}, \mu_\alpha \otimes dx)$)
\begin{align*}
\left| \int_{\mathbb{R}^d} \int_{\mathbb{R}^d} \langle (D^{\alpha - 1})^*f(x) ; \dfrac{\nabla(p_\alpha)(y)}{p_\alpha(y)} \rangle g \left(x + s^{\frac{1}{\alpha}} y \right) \mu_\alpha(dy)dx \right| & \leq \|(D^{\alpha - 1})^*(f)\|_{L^p(\mathbb{R}^d,dx)} \left\| \dfrac{\nabla(p_\alpha)}{p_\alpha} \right\|_{L^p(\mu_\alpha)} \\
&\quad\quad \times \|g\|_{L^q(\mathbb{R}^d,dx)}.
\end{align*}
Inequality \eqref{ineq:pseudo_poincar_frac_gradient}, for $p \in (1, +\infty)$, follows by straightforward arguments.~Let us explain now how one can obtain the important case $p=1$ for the inequality \eqref{ineq:pseudo_poincar_frac_gradient}. Let $f \in \mathcal{S}(\mathbb{R}^d)$ and let $g \in L^{\infty}(\mathbb{R}^d,dx)$ with 
$\|g\|_{L^{\infty}(\mathbb{R}^d,dx)} = 1$. Then, for all $t>0$, 
\begin{align*}
\langle P^\alpha_t(f) - f ; g \rangle_{L^2(\mathbb{R}^d,dx)}&  = \int_0^t \langle \mathcal{A}_\alpha P^\alpha_s(f) ; g \rangle_{L^2(\mathbb{R}^d,dx)} ds  \\
& = \int_0^t \langle P^\alpha_s  \mathcal{A}_\alpha(f) ; g \rangle_{L^2(\mathbb{R}^d,dx)} ds  \\
& = \int_0^t \langle  \mathcal{A}_\alpha(f) ;  P^\alpha_s (g) \rangle_{L^2(\mathbb{R}^d,dx)} ds  \\
& = \frac{1}{\alpha} \int_0^t \langle  \nabla \cdot D^{\alpha-1}(f) ;  P^\alpha_s (g) \rangle_{L^2(\mathbb{R}^d,dx)} ds,
\end{align*}
where $\nabla \cdot D^{\alpha-1} = \sum_{k = 1}^d \partial_k D^{\alpha-1}_k$. Now, the function $P_s^\alpha(g)$ is the convolution of $g$ with the kernel $s^{-\frac{d}{\alpha}}p_\alpha(./s^{\frac{1}{\alpha}})$ so that it belongs to $\mathcal{C}^1(\mathbb{R}^d)$, the set of continuously differentiable real-valued functions on $\mathbb{R}^d$, for all $s >0$, since 
\begin{align*}
\int_{\mathbb{R}^d} \|\nabla p_\alpha (x)\|dx <+\infty,
\end{align*}
(see, e.g., \cite[Theorem $2.1$]{KS_13}).~Moreover, for all $s >0$ and all $x \in \mathbb{R}^d$, 
\begin{align*}
\nabla P_s^\alpha(g)(x) = \frac{1}{s^{\frac{1}{\alpha}}} \int_{\mathbb{R}^d} g(z) \frac{1}{s^{\frac{d}{\alpha}}} \nabla p_\alpha \left( \frac{x - z}{s^{\frac{1}{\alpha}}} \right) dz = - \frac{1}{s^{\frac{1}{\alpha}}} \int_{\mathbb{R}^d} g(x + s^{\frac{1}{\alpha}} y ) \frac{\nabla p_\alpha \left( y \right) }{p_\alpha(y)} \mu_\alpha(dy),
\end{align*}
where the second equality follows by straightforward manipulations. Now, by integration by parts, for all $s \in (0, t]$
\begin{align*}
\langle  \nabla \cdot D^{\alpha-1}(f) ;  P^\alpha_s (g) \rangle_{L^2(\mathbb{R}^d,dx)} & = - \langle D^{\alpha-1}(f) ; \nabla P^\alpha_s (g) \rangle_{L^2(\mathbb{R}^d,dx)} \\
& = \frac{1}{s^{\frac{1}{\alpha}}}  \int_{\mathbb{R}^d} \int_{\mathbb{R}^d} \langle D^{\alpha - 1}(f)(x) ;  \dfrac{\nabla(p_\alpha)(y)}{p_\alpha(y)} \rangle g \left(x + s^{\frac{1}{\alpha}} y \right) \mu_\alpha(dy)dx.
\end{align*}
The end of the case $p=1$ follows similarly.~As for the inequality \eqref{ineq:pseudo_poincar_gradient}, for all $p \in (1, \alpha)$, for all $f \in \mathcal{S}(\mathbb{R}^d)$, all $g \in \mathcal{C}_c^{\infty}(\mathbb{R}^d)$ and all $t>0$, 
\begin{align*}
\langle P_t^\alpha(f) - f ; g \rangle_{L^2(\mathbb{R}^d,dx)} & = - \frac{1}{\alpha} \int_0^t \langle \nabla(f) ; D^{\alpha - 1} P^\alpha_s (g)\rangle_{L^2(\mathbb{R}^d,dx)} ds \\
& = - \frac{1}{\alpha} \int_0^t  \frac{ds}{s^{1-\frac{1}{\alpha}}} \int_{\mathbb{R}^d} \int_{\mathbb{R}^d} \langle \nabla(f)(x) ; y \rangle g(x+s^{\frac{1}{\alpha}} y) \mu_\alpha(dy)dx. 
\end{align*}
Once again, by H\"older's inequality with $p \in (1, \alpha)$ and with $q = p/(p-1)$, 
\begin{align*}
\left| \int_{\mathbb{R}^d} \int_{\mathbb{R}^d} \langle \nabla(f)(x) ; y \rangle g \left(x + s^{\frac{1}{\alpha}} y \right) \mu_\alpha(dy)dx \right| \leq \| \langle \nabla(f)(x) ; y \rangle \|_{L^p(\bbr^{2d},\mu_\alpha \otimes dx)} \|g\|_{L^q(\mathbb{R}^d,dx)}.
\end{align*}
Now, observe that, under $\mu_\alpha$ and for fixed $x \in \mathbb{R}^d$, $\langle \nabla(f)(x) ; y \rangle$ is an $\alpha$-stable random variable with characteristic function given, for all $\xi \in \bbr$, by 
\begin{align*}
\int_{\mathbb{R}^d} e^{i \xi \langle \nabla(f)(x) ; y \rangle} \mu_\alpha(dy) = \exp \left( - |\xi|^\alpha \int_{\mathbb{S}^{d-1}} |\langle y ; \nabla(f)(x) \rangle|^\alpha \lambda_1(dy) \right).
\end{align*}
Thus, for all $x \in \mathbb{R}^d$ fixed, 
\begin{align*}
\left(\int_{\mathbb{R}^d} \left| \langle \nabla(f)(x) ; y \rangle \right|^p \mu_\alpha(dy)\right) \leq \left(\bbe |Y_\alpha|^p\right) \sigma_\alpha(\nabla(f)(x))^p. 
\end{align*}
Note that since $\sigma$ is non-degenerate in the sense of \eqref{eq:non_deg}, for all $x \in \mathbb{R}^d$, $\sigma_\alpha(\nabla(f)(x))$ is bounded from above and below by $\| \nabla(f)(x) \|$ (up to some constant). Then, integrating with respect to the $x$-variable gives
\begin{align*}
\left( \int_{\mathbb{R}^d} \int_{\mathbb{R}^d} |\langle \nabla(f)(x) ; y \rangle|^p  \mu_\alpha(dy) dx \right)^{\frac{1}{p}} \leq \left(\bbe |Y_\alpha|^p \right)^{\frac{1}{p}} \left(\int_{\mathbb{R}^d} \sigma_\alpha \left(\nabla(f)(x)\right)^p dx \right)^{\frac{1}{p}}. 
\end{align*}
Thus, for all $t >0$, 
\begin{align*}
\left| \langle P^\alpha_t(f) - f ; g \rangle_{L^2(\mathbb{R}^d,dx)} \right| \leq t^{\frac{1}{\alpha}} \left(\bbe |Y_\alpha|^p \right)^{\frac{1}{p}} \| \sigma_\alpha \left(\nabla(f)\right) \|_{L^p(\mathbb{R}^d,dx)} \|g\|_{L^q(\mathbb{R}^d,dx)}. 
\end{align*}
To conclude, let us consider the case $p = 1$ for the inequality \eqref{ineq:pseudo_poincar_gradient}.~Let $f \in \mathcal{S}(\mathbb{R}^d)$ and let $g \in L^\infty(\mathbb{R}^d,dx)$ be such that $\| g \|_{L^\infty(\mathbb{R}^d,dx)} = 1$.~As previously, by straightforward semigroup arguments, 
\begin{align*}
\langle P_t^\alpha(f) -f ; g \rangle_{L^2(\mathbb{R}^d,dx)} = \int_0^t \langle \mathcal{A}_\alpha(f)  ; P_s^\alpha (g) \rangle_{L^2(\mathbb{R}^d,dx)} ds = \frac{1}{\alpha}  \int_0^t \langle  D^{\alpha-1} \cdot \nabla(f)  ; P_s^\alpha (g) \rangle_{L^2(\mathbb{R}^d,dx)} ds. 
\end{align*}
Now, by duality, for all $s \in (0,t]$, 
\begin{align*}
\langle  D^{\alpha-1} \cdot \nabla(f)  ; P_s^\alpha (g) \rangle_{L^2(\mathbb{R}^d,dx)} = -  \langle  \nabla(f)  ; D^{\alpha-1} P_s^\alpha (g) \rangle_{L^2(\mathbb{R}^d,dx)}.
\end{align*}
Moreover, recall that, for all $x \in \mathbb{R}^d$, $D^{\alpha-1}(p_\alpha)(x) = -x p_\alpha(x)$. Thus, for all $x \in \mathbb{R}^d$, by a change of variables in the radial coordinate, 
\begin{align*}
D^{\alpha-1} P_s^\alpha (g)(x) & = \int_{\mathbb{R}^d} g(z) \left(\int_{\mathbb{R}^d} u \nu_\alpha(du) \left( p_\alpha\left( \frac{x+u-z}{s^{\frac{1}{\alpha}}} \right) - p_\alpha\left( \frac{x-z}{s^{\frac{1}{\alpha}}} \right) \right) \right) \frac{dz}{s^{\frac{d}{\alpha}}} \\
& = \frac{1}{s^{1 - \frac{1}{\alpha}}} \int_{\mathbb{R}^d} g(z) \left(\int_{\mathbb{R}^d} u \nu_\alpha(du) \left( p_\alpha\left( \frac{x-z}{s^{\frac{1}{\alpha}}} + u \right) - p_\alpha\left( \frac{x-z}{s^{\frac{1}{\alpha}}} \right) \right) \right) \frac{dz}{s^{\frac{d}{\alpha}}} \\
& = -  \frac{1}{s^{1 - \frac{1}{\alpha}}} \int_{\mathbb{R}^d} g(z)  \frac{x-z}{s^{\frac{1}{\alpha}}}  p_\alpha\left( \frac{x-z}{s^{\frac{1}{\alpha}}} \right) \frac{dz}{s^{\frac{d}{\alpha}}} \\
& = \frac{1}{s^{1 - \frac{1}{\alpha}}} \int_{\mathbb{R}^d} g(x + s^{\frac{1}{\alpha}}z)  z  p_\alpha\left( z \right) dz.
\end{align*}
The end of the proof of the proposition follows easily.$\qed$\\

\noindent 
Next, based on Lemma \ref{lem:local_reverse_poincar} and the inequality \eqref{ineq:sub_commutation}, let us investigate the pseudo-Poincar\'e inequality with the non-linear operator $\nabla_{\nu_\alpha}$. 

\begin{prop}\label{prop:pseudo_poincar_gradient_length}
Let $\alpha \in (1,2)$, let $\nu_\alpha$ be a non-degenerate symmetric $\alpha$-stable L\'evy measure on $\mathbb{R}^d$, $d \geq 1$, and let $(P^\alpha_t)_{t \geq 0}$ be the associated stable heat semigroup.~Then,  for all $p \in [1,2]$, all $f \in \mathcal{S}(\mathbb{R}^d)$ and all $t>0$, 
\begin{align}\label{ineq:pseudo_poincar_gradient_length}
\|f- P^\alpha_t(f)\|_{L^p(\mathbb{R}^d,dx)} \leq \sqrt{t} \| \nabla_{\nu_\alpha}(f) \|_{L^p(\mathbb{R}^d,dx)}.
\end{align}
\end{prop}

\begin{proof}
Let $f \in \mathcal{S}(\mathbb{R}^d)$ and let $t>0$. At first let $p \in (1,2]$. Then, take $g \in \mathcal{C}_c^{\infty}(\mathbb{R}^d)$.  Thus, by semigroup arguments and making use of \eqref{eq:ipp_symmetry}, 
\begin{align*}
\langle P^\alpha_t(f) - f ; g \rangle_{L^2(\mathbb{R}^d,dx)} & = \int_0^t \langle f ; \mathcal{A}_\alpha P^\alpha_s(g) \rangle_{L^2(\mathbb{R}^d,dx)} ds \\
& = - \int_0^t \left(\int_{\mathbb{R}^d} \Gamma_\alpha(f , P^\alpha_s(g))(x) dx \right) ds. 
\end{align*}
Thus, by the Cauchy-Schwarz inequality together with H\"older's inequality where $q = p/(p-1)$, 
\begin{align*}
|\langle P^\alpha_t(f) - f ; g \rangle_{L^2(\mathbb{R}^d,dx)}| \leq \frac{1}{2} \|\nabla_{\nu_\alpha}(f)\|_{L^p(\mathbb{R}^d,dx)} \int_0^t  \|\nabla_{\nu_\alpha}(P_s^\alpha(g))\|_{L^q(\mathbb{R}^d,dx)} ds.
\end{align*}
Now, Lemma \ref{lem:local_reverse_poincar} as well as convexity arguments ($q \geq 2$),  for all $s \in (0,t]$,
\begin{align*}
\|\nabla_{\nu_\alpha}(P_s^\alpha(g))\|^q_{L^q(\mathbb{R}^d,dx)} & \leq \frac{1}{s^{\frac{q}{2}}} \int_{\mathbb{R}^d} \left(P_s^\alpha(g^2)(x)\right)^{\frac{q}{2}} dx \\
& \leq \frac{1}{s^{\frac{q}{2}}} \int_{\mathbb{R}^d} P_s^\alpha(|g|^q)(x) dx = \frac{1}{s^{\frac{q}{2}}} \int_{\mathbb{R}^d} |g(x)|^q dx.
\end{align*}
The end of the proof for $p \in (1,2]$ follows easily.~Let us discuss now the case $p =1$.~Take $g \in L^\infty(\mathbb{R}^d,dx)$ with $\|g\|_{L^{\infty}(\mathbb{R}^d,dx)} = 1$.  Then,  
\begin{align*}
\langle P^\alpha_t(f) - f ; g \rangle_{L^2(\mathbb{R}^d,dx)} & = \int_0^t \langle \mathcal{A}_\alpha P_s^\alpha (f) ; g \rangle_{L^2(\mathbb{R}^d,dx)} ds \\
& =\int_0^t \langle P_s^\alpha \mathcal{A}_\alpha (f) ; g \rangle_{L^2(\mathbb{R}^d,dx)} ds \\
& =\int_0^t \langle \mathcal{A}_\alpha (f) ; P_s^\alpha (g) \rangle_{L^2(\mathbb{R}^d,dx)} ds \\
& = - \int_{0}^t \left( \int_{\mathbb{R}^d}  \Gamma_\alpha (f ,  P^\alpha_s(g))(x) dx \right) ds.
\end{align*}
The end of the proof for this case follows similarly. 
\end{proof}

\begin{rem}\label{rem:Lp_norm_gradient_length}
Let us discuss the finiteness of the quantity $\|\nabla_{\nu_\alpha}(f)\|_{L^p(\mathbb{R}^d,dx)}$, for $p \in [1, 2]$ and $f \in \mathcal{S}(\mathbb{R}^d)$.~For simplicity, let us consider a bump function $f$ with compact support contained in the Euclidean unit ball and the L\'evy measure $\nu_\alpha(du)= du/\|u\|^{\alpha+d}$.~Then, for all $x \in \mathbb{R}^d$ with $\|x\| \geq 3$, 
\begin{align*}
(\nabla_{\nu_\alpha}(f)(x))^2 = \int_{\mathbb{R}^d} |f(x+u)|^2 \frac{du}{\|u\|^{\alpha+d}} = \int_{B(0,1)} |f(v)|^2 \frac{dv}{\|x-v\|^{\alpha+d}} \geq \frac{c_f}{\|x\|^{\alpha+d}},
\end{align*}
for some $c_f>0$.~Then, for all $d \geq 2$ and all $p \in [1, 2d/(\alpha+d)]$, $\nabla_{\nu_\alpha}(f) \notin L^p(\mathbb{R}^d,dx)$.~Thus, \eqref{ineq:pseudo_poincar_gradient_length} is only relevant for $p \in (2d/(\alpha+d),2]$.~While for $d=1$, $\nabla_{\nu_\alpha}(f) \in L^p(\mathbb{R}^d,dx)$, for all $p \in [1,2]$.
\end{rem}

\begin{rem}\label{rem:nonlocal_Li_Yau}
As emphasized in \cite[Proof of Lemma $3$]{MLedoux_03}, in order to reach a pseudo-Poincar\'e inequality for the cases $p \in (2,+\infty)$ (in the absence of a Bismut-type formula),  a Li-Yau-type inequality is the main ingredient in the diffusive setting.~Indeed, the abstract integration by parts formula and the diffusion property (which ensures for the chain rule to hold) together with standard arguments allow to treat the cases $p \in(2, +\infty)$.~In the non-local setting,  not much seems to be known regarding Li-Yau-type inequalities.~For the isotropic fractional Laplacian, a certain form of this inequality has been recently obtained in \cite[Lemma $3.1$ and Theorem $3.2$]{WZ_22} despite the fact that this operator has infinite dimension (see \cite{SWZ_20}).
\end{rem} 

\begin{rem}\label{rem:pseudo_poincar_stable_heat_semigroup}
(i) In the Gaussian case, it is well-known that (see, e.g., \cite{MLedoux_03}), for all $p \in [1, +\infty)$, all $t>0$ and all $f \in \mathcal{S}(\mathbb{R}^d)$, 
\begin{align}\label{eq:pseudo_poincar_gaussien}
\| P^H_t(f) - f \|_{L^p(\mathbb{R}^d,dx)} \leq \sqrt{2t} \left(\bbe |X_1|^p \right)^{\frac{1}{p}} \| \nabla(f) \|_{L^p(\mathbb{R}^d,dx)},
\end{align}
where $X_1$ is a standard normal random variable.~Thus, inequalities \eqref{ineq:pseudo_poincar_frac_gradient} and \eqref{ineq:pseudo_poincar_gradient} are straightforward generalizations of inequality \eqref{eq:pseudo_poincar_gaussien} to the non-degenerate symmetric $\alpha$-stable case with $\alpha \in (1,2)$.~In particular, they take into account the fine structure of the generator $\mathcal{A}_\alpha$ as well as the polar decomposition of the L\'evy measure $\nu_\alpha$.~When the spectral measure $\lambda_1$ is invariant by rotation, for all $x \in \mathbb{R}^d$, 
\begin{align*}
\sigma_\alpha(\nabla(f)(x)) = c_{\alpha,d} \|\nabla(f)(x)\|,
\end{align*}
for some well-chosen $c_{\alpha,d}>0$ depending only on $\alpha$ and $d$.~Finally, let us assume that the spectral measure $\lambda_1$ is discrete and given by 
\begin{align}\label{eq:spectral_measure_discrete}
\lambda^{\operatorname{dis}}_1(dy) = \sum_{k = 1}^d w_k \delta_{\pm e_k}(dy),
\end{align}
where $e_k = (0, \dots, 0, 1, 0 , \dots, 0)^t$, $\delta_{e_k}$ is the Dirac measure at $e_k$ and $w_k$ are positive weights, for all $k \in \{1, \dots, d\}$. Then, for all $f \in \mathcal{S}(\mathbb{R}^d)$ and all $x \in \mathbb{R}^d$, 
\begin{align}\label{eq:norm_discrete}
\sigma^{\operatorname{dis}}_\alpha\left(\nabla(f)(x)\right) = 2^{\frac{1}{\alpha}} \left(\sum_{k=1}^d w_k \left| \partial_k(f)(x) \right|^{\alpha}\right)^{\frac{1}{\alpha}}. 
\end{align}
Thus, the quantity $\sigma_\alpha(\nabla(f)(x))$ encodes the anisotropy inherent to $\nu_\alpha$ into the upper bound of \eqref{ineq:pseudo_poincar_gradient}.\\
(ii) Let us further discuss the integrability condition \eqref{ineq:condint_logarithmic_derivative} in the rotationally invariant case. Then, for all $x \in \mathbb{R}^d$,
\begin{align*}
\dfrac{C_1}{\left(1+ \|x\|\right)^{\alpha+d}} \leq p_{\alpha}^{\operatorname{rot}}(x) \leq \dfrac{C_2}{\left(1+ \|x\|\right)^{\alpha+d}}, \quad \|\nabla(p_\alpha^{\operatorname{rot}})(x)\| \leq \dfrac{C_3}{\left(1+ \|x\|\right)^{\alpha + d + 1}},
\end{align*}
for some $C_1, C_2$ and $C_3$ positive constants depending on $\alpha$ and $d$ (see, e.g., \cite{W07}). Thus, in this case, $\|\nabla(p_\alpha)\|/p_\alpha \in L^p(\mu_\alpha)$, for all $p \in [1, +\infty)$. 
\end{rem}
\noindent
Now, recall that $(P^\alpha_t)_{t \geq 0}$ satisfies the following ultracontractive estimates: for all $t>0$ and all $f \in L^1(\mathbb{R}^d,dx)$, 
\begin{align}\label{ineq:ultracontractivity_stable_heat_L1}
\left \| P_t^\alpha(f)\right \|_{L^\infty(\mathbb{R}^d,dx)} \leq \frac{\| p _\alpha \|_{L^\infty(\mathbb{R}^d,dx)}}{t^{\frac{d}{\alpha}}} \| f \|_{L^1(\mathbb{R}^d,dx)}
\end{align}
and, for all $p \in [1,+\infty)$ and all $f \in L^p(\mathbb{R}^d,dx)$, 
\begin{align}\label{ineq:ultracontractivity_stable_heat_Lp}
\left \| P_t^\alpha(f)\right \|_{L^\infty(\mathbb{R}^d,dx)} \leq \frac{\| p _\alpha \|^{\frac{1}{p}}_{L^\infty(\mathbb{R}^d,dx)}}{t^{\frac{d}{\alpha p}}} \| f \|_{L^p(\mathbb{R}^d,dx)}.
\end{align}
In the sequel, we are interested in Sobolev-type inequalities associated with the non-degenerate symmetric $\alpha$-stable probability measures on $\mathbb{R}^d$ generalizing the classical Sobolev inequality and its refined versions.~The incoming proofs are based on the pseudo-Poincar\'e inequalities obtained in Proposition \ref{prop:pseudo_poincare} and Proposition \ref{prop:pseudo_poincar_gradient_length} and on the ultracontractive estimates \eqref{ineq:ultracontractivity_stable_heat_L1} and \eqref{ineq:ultracontractivity_stable_heat_Lp}.~Based on these ultracontractive estimates, let us define a Besov space specifically designed for $(P_t^\alpha)_{t \geq 0}$, namely, for all $s<0$ and all $\alpha \in (1,2)$, let $\mathcal{B}_{\infty, \infty}^{s,\alpha}$ be the set of function $f$ in $\cup_{p \in [1,+\infty]} L^p(\mathbb{R}^d, dx)$ such that
\begin{align*}
\| f \|_{\mathcal{B}_{\infty, \infty}^{s,\alpha}} : = \underset{t >0 }{\sup}\, t^{- \frac{s}{\alpha}} \|P^\alpha_t(f)\|_{L^{\infty}(\mathbb{R}^d,dx)}<+\infty.  
\end{align*}

Next, let us discuss the geometric content of the non-degenerate symmetric $\alpha$-stable probability measures on $\mathbb{R}^d$.  Under the non-degeneracy condition \eqref{eq:non_deg}, the mapping $\sigma_\alpha$ defined by \eqref{eq:rep_spectral_measure} is a norm on $\mathbb{R}^d$.  Let $K_\alpha$ be defined by
\begin{align}\label{def:unit_ball_stable}
K_\alpha = \{ z \in \mathbb{R}^d,\,  \sigma_\alpha(z) \leq 1 \}. 
\end{align}
$K_\alpha$ is a compact convex subset of $\mathbb{R}^d$ which contains the origin and is symmetric with respect to it (if $z \in K_\alpha$, then $-z \in K_\alpha$). Let $\mathring{K}_\alpha$ be its polar set given by
\begin{align}\label{def:polar_set_stable}
\mathring{K}_\alpha = \{ x \in \mathbb{R}^d , \langle x; z \rangle \leq 1, \,  \forall z \in K_\alpha \}
\end{align}
and with the corresponding dual norm given, for all $x \in \mathbb{R}^d$, by
\begin{align}\label{def:norm_dual_syn_ndeg}
\sigma^*_\alpha(x) : = \underset{z \in \mathbb{R}^d, \, \sigma_\alpha(z)\leq 1}{\sup} |\langle x; z \rangle|. 
\end{align} 
Let us describe briefly the sets $K_\alpha$ and $\mathring{K}_\alpha$ for specific choices of the spectral measure $\lambda_1$.  When $\alpha \in (1,2)$ and $\nu_\alpha = \nu_\alpha^{\operatorname{rot}}$, for all $z \in \mathbb{R}^d$, 
\begin{align}\label{eq:rot_inv_sigma}
\sigma_\alpha^{\operatorname{rot}}(z) = \frac{\|z\|}{2^{\frac{1}{\alpha}}},
\end{align}
so that $K^{\operatorname{rot}}_\alpha  = B(0 , 2^{\frac{1}{\alpha}})$ and $\mathring{K}^{\operatorname{rot}}_\alpha = B(0, 1/2^{1/\alpha})$.  Next, let us consider the discrete case in full generality. Then, for all $z \in \mathbb{R}^d$, 
\begin{align}\label{eq:discrete_case_full_generality}
\sigma^{\operatorname{dis}}_\alpha(z) = \left( \sum_{k=1}^d 2 w_k |z_k|^{\alpha}  \right)^{\frac{1}{\alpha}},  \quad \sigma_\alpha^{\operatorname{dis},*}(z) = \left(\sum_{k = 1}^d \dfrac{|z_k|^{\frac{\alpha}{\alpha-1}}}{(2w_k)^{\frac{1}{\alpha-1}}}\right)^{\frac{\alpha-1}{\alpha}},
\end{align}
where $(w_1, \dots, w_d) \in (\mathbb{R}_+^*)^d$. From the definition, 
\begin{align*}
K_\alpha^{\operatorname{dis}} \subset R_\alpha := [-1/(2w_1)^{1/\alpha}, 1/(2w_1)^{1/\alpha}] \times \cdots \times [-1/(2w_d)^{1/\alpha}, 1/(2w_d)^{1/\alpha}].
\end{align*}
In particular, when $w_k = 1/2$, for $k \in \{1, \dots, d\}$,  $K_\alpha^{\operatorname{dis}}$ is exactly the unit ball (denoted also by $\mathbb{B}_\alpha^d$) of the Minkowski space $\ell_\alpha^d$ and its polar set $\mathring{K}_\alpha^{\operatorname{dis}}$ is the unit ball of $\ell^d_{\alpha/(\alpha-1)}$. The $d$-dimensional Lebesgue measures of these convex bodies have been computed exactly (see, e.g., \cite{Pisier_book89} and the references therein) and are given, for all $\alpha \in (1,2)$, by 
\begin{align*}
\mathcal{L}_d(\mathbb{B}_\alpha^d) = \dfrac{(2 \Gamma(1+\frac{1}{\alpha}))^d}{\Gamma(1+ \frac{d}{\alpha})}, \quad \mathcal{L}_d(\mathbb{B}_{\alpha/(\alpha-1)}^d) = \dfrac{(2 \Gamma(1+\frac{\alpha-1}{\alpha}))^d}{\Gamma(1+ \frac{(\alpha-1)d}{\alpha})}.
\end{align*}
Let us briefly describe geometric functionals associated with the set $K_\alpha$ (and its boundary $\partial K_\alpha$).  First, by a change of variables in spherical coordinates,
\begin{align}\label{eq:formula_volume}
\mathcal{L}_d(K_\alpha) = \frac{1}{d} \int_{\mathbb{S}^{d-1}} \frac{\sigma_L(d\theta)}{(\sigma_\alpha(\theta))^d}. 
\end{align}
Let $\mu_{K_\alpha}$ be the cone measure on $\partial K_\alpha$ defined, for all $B$ Borel subset of $\partial K_\alpha$, by
\begin{align}\label{def:cone_measure_K_alpha}
\mu_{K_\alpha}(B) := \dfrac{\mathcal{L}_d([0,1]B)}{\mathcal{L}_d(K_\alpha)},
\end{align}
where $[0,1]B = \{ t\omega , \, t \in [0,1],\, \omega \in B \}$. Thanks to spherical coordinates, $x = r\theta$ belongs to $[0,1]B$ if and only if $(r,\theta/\sigma_\alpha(\theta)) \in [0, 1/\sigma_\alpha(\theta)]\times B$. Thus, for all $B \in \mathcal{B}(\partial K_\alpha)$, 
\begin{align}\label{eq:representation_formula_integral_cone_measure_1}
\mu_{K_\alpha}(B) = \frac{1}{d \mathcal{L}_d(K_\alpha)} \int_{\mathbb{S}^{d-1}} \bbone_{B} \left( \dfrac{\theta}{\sigma_\alpha(\theta)}   \right) \frac{\sigma_L(d\theta)}{(\sigma_\alpha(\theta))^d}. 
\end{align}
By a standard procedure, it is possible to extend the previous formula to all suitable function $f$ defined on $\partial K_\alpha$ by
\begin{align}\label{eq:representation_formula_integral_cone_measure_2}
\mu_{K_\alpha}(f) : = \int_{\partial K_\alpha} f(x) \mu_{K_\alpha}(dx) = \frac{1}{d \mathcal{L}_d(K_\alpha)} \int_{\mathbb{S}^{d-1}} f \left( \dfrac{\theta}{\sigma_\alpha(\theta)}   \right) \frac{\sigma_L(d\theta)}{(\sigma_\alpha(\theta))^d}.
\end{align}
Finally, recall the following polar-type coordinate integration formula (see, e.g., \cite[Proposition $1$]{Naor_Romik_2003}): for all $f \in L^1(\mathbb{R}^d,dx)$, 
\begin{align}\label{eq:polar_coordinate_cone_measure}
\int_{\mathbb{R}^d} f(x) dx = d \mathcal{L}_d(K_\alpha) \int_0^{+\infty} s^{d-1} \left(\int_{\partial K_\alpha} f(s \omega) \mu_{K_\alpha}(d\omega) \right) ds.
\end{align}
Based on these sets, it is natural to consider a notion of anisotropic perimeter related to the non-degenerate symmetric $\alpha$-stable probability measures on $\mathbb{R}^d$, with $\alpha \in (1,2)$.  First, let us define it by a variational representation formula: for all $E\in \mathcal{B}(\mathbb{R}^d)$ with $\mathcal{L}_d(E)<+\infty$, let   
\begin{align}\label{def:anisotropic_perimeter_stable}
\mathcal{P}_{\alpha}^{\operatorname{var}}(E) = \sup \left\{ \int_{E} \operatorname{div}(\phi(x))dx , \, \phi \in \mathcal{C}_c^1(\mathbb{R}^d, \mathbb{R}^d), \, \| \sigma_{\alpha}^*(\phi)\|_{\infty} \leq 1\right\}. 
\end{align}
The pseudo-Poincar\'e inequality \eqref{ineq:pseudo_poincar_gradient} encodes this (an)isotropic geometric feature which can be observed as well from the corresponding sharp strong Sobolev and sharp isoperimetric inequalities. Indeed, \cite{Milman_Schechtman_86,Lions_97,cordero_nazaret_villani,Van_Schaftingen_Jean_06,Cabre_CAM17},  for all suitable $f$ and $E$, 
\begin{align}\label{ineq:sharp_anisotropic_sobolev}
\int_{\mathbb{R}^d} \sigma_\alpha\left(\nabla(f)(x)\right) dx \geq d \mathcal{L}_d\left(\mathring{K}_\alpha\right)^{\frac{1}{d}} \left( \int_{\mathbb{R}^d} \left| f(x) \right|^{\frac{d}{d-1}} dx \right)^{\frac{d-1}{d}}, \quad d \geq 2,
\end{align}
and
\begin{align}\label{ineq:sharp_anisotropic_iso}
\mathcal{P}_{\alpha}^{\operatorname{var}}(E) \geq d \mathcal{L}_d\left(\mathring{K}_\alpha\right)^{\frac{1}{d}} \mathcal{L}_d(E)^{\frac{d-1}{d}}.
\end{align}
The equality case in \eqref{ineq:sharp_anisotropic_iso} is obtained when $E$ is (up to scaling and translation) equal to $\mathring{K}_\alpha$.  Note that one can retrieve the exact formula for the perimeter $\mathcal{P}_\alpha^{\operatorname{var}}(\mathring{K}_\alpha)$ using,
\begin{align*}
\mathcal{P}_\alpha^{\operatorname{var}}(\mathring{K}_\alpha) = \int_{\partial^* \mathring{K}_\alpha} \sigma_\alpha\left( \nu_{\mathring{K}_\alpha}(y)\right) \mathcal{H}^{d-1}(dy),
\end{align*}
where $\nu_{\mathring{K}_\alpha}(y)$ is the (measure theoretic) outer normal to $\mathring{K}_\alpha$ at $y \in \partial^* \mathring{K}_\alpha$, together with the explicit Radon-Nikod\'ym derivative of the surface measure on $\partial \mathring{K}_\alpha$ with respect to the cone measure $\mu_{\mathring{K}_\alpha}$ (see \cite[Lemma $1$]{Naor_Romik_2003}).

Finally, to end this preliminaries section, let us state and prove composition formulas as well as potential theoretic results which are used to establish Sobolev-type embeddings with the fractional gradient operator $D^{\alpha-1}$.~These composition formulas for fractional operators have been obtained in \cite{Hor_PAMS59} and, more recently, stated in \cite[Theorem $5.3$]{Silhavy_CMT20}, for the rotationally invariant case only. Let us recall the Riesz transform-type operators introduced in \cite{AH20_4}: for all $\alpha \in (1,2)$, let $\mathcal{R}_\alpha$ and $\mathcal{R}^\alpha$ be the linear operators defined, for all $f \in \mathcal{C}_c^{\infty}(\mathbb{R}^d)$, all $\varphi \in \mathcal{C}_c^{\infty}(\mathbb{R}^d, \mathbb{R}^d)$ and all $x \in \mathbb{R}^d$, by 
\begin{align}\label{eq:Riesz_transform_frac}
\mathcal{R}_\alpha(f)(x) = \frac{1}{(2\pi)^{d}} \int_{\mathbb{R}^d} \mathcal{F}(f)(\xi) m_\alpha(\xi) e^{i \langle x; \xi \rangle} d\xi,
\end{align}
and
\begin{align*}
\mathcal{R}^{\alpha}(\varphi)(x) = \frac{1}{(2\pi)^{d}} \int_{\mathbb{R}^d} \langle \mathcal{F}(\varphi)(\xi) ; m_\alpha(\xi)\rangle e^{i \langle x; \xi \rangle} d\xi,
\end{align*}
where, for all $\xi \in \mathbb{R}^d \setminus \{0\}$, 
\begin{align}\label{eq:symbol_Riesz_transform}
m_\alpha(\xi) = i \alpha \nabla(\sigma_\alpha)(\xi).
\end{align}

\begin{prop}\label{prop:composition_formulas}
Let $1<\alpha$, $\beta<2$. Then, for all $f \in \mathcal{C}^\infty_c(\mathbb{R}^d)$, $d \geq 1$, and all $x \in \mathbb{R}^d$, 
\begin{align}\label{eq:composition_formulas_NDS}
(\operatorname{div}_\beta \circ D^{\alpha-1})(f)(x) = \mathcal{R}^\alpha\left(\mathcal{R}_\beta \circ (- \mathcal{A}_\alpha)^{\frac{\alpha-1}{\alpha}} \circ (- \mathcal{A}_\beta)^{\frac{\beta-1}{\beta}}(f)\right)(x).
\end{align}
In particular, for all $f \in \mathcal{C}^\infty_c(\mathbb{R}^d)$ and all $x \in \mathbb{R}^d$, 
\begin{align}\label{eq:composition_formula_rot}
(\operatorname{div}^{\operatorname{rot}}_\beta \circ D^{\alpha-1, \operatorname{rot}})(f)(x) = - \frac{\alpha}{2} \frac{\beta}{2} \left(- \Delta\right)^{\frac{\alpha+\beta-2}{2}}(f)(x).
\end{align}
\end{prop}

\begin{proof}
Let us start with the proof of \eqref{eq:composition_formula_rot}.~Let $f \in \mathcal{C}_c^{\infty}(\mathbb{R}^d)$. By the Fourier inversion formula (see, e.g., \cite{AH20_4}), for all $x \in \mathbb{R}^d$, 
\begin{align*}
D^{\alpha-1, \operatorname{rot}}(f)(x) = \frac{1}{(2\pi)^d} \int_{\mathbb{R}^d} \mathcal{F}(f)(\xi) \tau_\alpha^{\operatorname{rot}}(\xi) e^{i \langle x;\xi \rangle} d\xi, 
\end{align*}
where, for all $\xi \in \mathbb{R}^d \setminus \{0\}$, 
\begin{align*}
 \tau_\alpha^{\operatorname{rot}}(\xi) = \frac{\alpha}{2} \dfrac{i \xi }{\|\xi\|} \|\xi\|^{\alpha-1}. 
\end{align*}
Similarly, for all $\varphi \in \mathcal{C}_c^{\infty}(\mathbb{R}^d , \mathbb{R}^d)$ and all $x \in \mathbb{R}^d$, 
\begin{align*}
\operatorname{div}_\beta^{\operatorname{rot}}(\varphi)(x) & = \int_{\mathbb{R}^d} \langle \varphi(x+u) - \varphi(x) ; u \rangle \nu_\beta^{\operatorname{rot}}(du) = \frac{1}{(2\pi)^d} \int_{\mathbb{R}^d} \langle \mathcal{F}(\varphi)(\xi) ; \tau_\beta^{\operatorname{rot}}(\xi)\rangle e^{i \langle x ; \xi \rangle} d\xi. 
\end{align*}
Then, for all $x \in \mathbb{R}^d$, 
\begin{align*}
(\operatorname{div}^{\operatorname{rot}}_\beta \circ D^{\alpha-1, \operatorname{rot}})(f)(x) & =  \frac{1}{(2\pi)^d} \int_{\mathbb{R}^d} \langle \mathcal{F}\left(D^{\alpha-1, \operatorname{rot}})(f)\right)(\xi) ; \tau_\beta^{\operatorname{rot}}(\xi)\rangle e^{i \langle x ; \xi \rangle} d\xi \\
& = \frac{1}{(2\pi)^d} \int_{\mathbb{R}^d} \mathcal{F}(f)(\xi) \langle \tau_\alpha^{\operatorname{rot}}(\xi) ; \tau_\beta^{\operatorname{rot}}(\xi)\rangle e^{i \langle x ; \xi \rangle} d\xi.
\end{align*}
Now, for all $\xi \in \mathbb{R}^d \setminus \{0\}$, 
\begin{align*}
\langle \tau_\alpha^{\operatorname{rot}}(\xi) ;  \tau_\beta^{\operatorname{rot}}(\xi) \rangle = - \frac{\alpha}{2} \frac{\beta}{2} \|\xi\|^{\alpha+ \beta -2}. 
\end{align*}
Thus, for all $f \in \mathcal{C}_c^{\infty}(\mathbb{R}^d)$ and all $x \in \mathbb{R}^d$, 
\begin{align*}
(\operatorname{div}^{\operatorname{rot}}_\beta \circ D^{\alpha-1, \operatorname{rot}})(f)(x) = - \frac{\alpha}{2} \frac{\beta}{2} \left(- \Delta\right)^{\frac{\alpha+\beta-2}{2}}(f)(x),
\end{align*}
proving \eqref{eq:composition_formula_rot}. For the general non-degenerate symmetric case, 
\begin{align*}
& D^{\alpha-1}(f)(x) = \frac{1}{(2\pi)^d} \int_{\mathbb{R}^d} \mathcal{F}(f)(\xi) \tau_\alpha(\xi) e^{i \langle x;\xi \rangle} d\xi, \\
& \operatorname{div}_\beta(\varphi)(x) =  \frac{1}{(2\pi)^d} \int_{\mathbb{R}^d} \langle \mathcal{F}(\varphi)(\xi) ; \tau_\beta(\xi)\rangle e^{i \langle x ; \xi \rangle} d\xi, 
\end{align*}
where, for all $\xi \in \mathbb{R}^d \setminus \{0\}$, 
\begin{align}\label{eq:tau_alpha_NDS}
\tau_\alpha(\xi) = \frac{1}{i} \dfrac{\nabla(\widehat{\mu_\alpha})(\xi)}{\widehat{\mu_\alpha}(\xi)} = i \alpha \nabla(\sigma_\alpha)(\xi) (\sigma_\alpha(\xi))^{\alpha-1},
\end{align}
with $\sigma_\alpha$ given by \eqref{eq:rep_spectral_measure}. Thus, for all $x \in \mathbb{R}^d$, 
\begin{align*}
\operatorname{div}_\beta \left(D^{\alpha-1}(f)\right)(x) & = \frac{1}{(2\pi)^d} \int_{\mathbb{R}^d} \mathcal{F}(f)(\xi) e^{i \langle x ; \xi \rangle} \langle \tau_\beta(\xi) ; \tau_\alpha(\xi) \rangle d\xi \\
& = \frac{1}{(2\pi)^d} \int_{\mathbb{R}^d} \mathcal{F}(f)(\xi) e^{i \langle \xi ; x \rangle} \sigma_{\alpha,\beta}(\xi) d\xi,
\end{align*}
where, for all $ \xi \in \mathbb{R}^d \setminus \{0\}$, 
\begin{align*}
\sigma_{\alpha,\beta}(\xi) & = - \alpha \beta \langle \nabla(\sigma_\alpha)(\xi) ; \nabla(\sigma_\beta)(\xi) \rangle (\sigma_{\alpha}(\xi))^{\alpha-1}(\sigma_{\beta}(\xi))^{\beta-1} \\
& = \langle m_\alpha(\xi) ; m_{\beta}(\xi) \rangle (\sigma_{\alpha}(\xi))^{\alpha-1}(\sigma_{\beta}(\xi))^{\beta-1}, 
\end{align*}
from which the end of the proof easily follows. 
\end{proof}
\noindent
The next lemma provides an explicit spatial-integral  representation of an identity contained in \cite[Theorem $1.12$]{Shieh_Spector_ACV15}. For this purpose, recall the Riesz potential operator of order $s \in (0,1)$ denoted by $I^s$ and defined, for all $f \in \mathcal{S}(\mathbb{R}^d)$ and all $x \in \mathbb{R}^d$, by 
\begin{align}\label{eq:def_Riesz_Potential_Operator}
I^s(f)(x) = \frac{1}{(2\pi)^d} \int_{\mathbb{R}^d} \mathcal{F}(f)(\xi) e^{i \langle \xi ;x \rangle} \frac{d\xi}{\|\xi\|^s} = \int_{\mathbb{R}^d} f(y) k^{\operatorname{rot}}_s(x-y)dy,
\end{align}
where $k^{\operatorname{rot}}_s$ is the Riesz kernel of order $s$ such that $k^{\operatorname{rot}}_s(x ) = \gamma_{s,d} \| x \|^{s-d}$, $x \in \mathbb{R}^d$ with $x \ne 0$, with 
\begin{align}\label{eq:normalizing_constant_Rieszkernel}
\gamma_{s,d} = \frac{1}{2^s \pi^{\frac{d}{2}}} \dfrac{\Gamma\left(\frac{d-s}{2}\right)}{\Gamma\left(\frac{s}{2}\right)}.
\end{align}

\begin{lem}\label{lem:Fourier_Transform_HD}
Let $d \geq 2$ be an integer, let $s \in (0,d)$ and let $k_s$ be defined, for all $x \in \mathbb{R}^d \setminus \{0\}$, by 
\begin{align}\label{eq:convolution_Riesz_Riesz}
k_s(x) = \dfrac{x}{\|x\|^{d+1-s}}. 
\end{align}   
Then, for all $\xi \in \mathbb{R}^d \setminus \{0\}$, 
\begin{align}\label{eq:TF_CRiesz_Riesz}
\mathcal{F}(k_s)(\xi) = 2^{s} \dfrac{\Gamma\left(\frac{s+1}{2}\right)}{\Gamma(\frac{d-s+1}{2})}\pi^{\frac{d}{2}} \frac{-i\xi}{\|\xi\|^{s+1}}.
\end{align}
Moreover, for all $f \in \mathcal{C}_c^{\infty}(\mathbb{R}^d)$ and all $x \in \mathbb{R}^d$, 
\begin{align}\label{eq:spatial_rep_comistefani}
f(x) & = \frac{2}{\alpha} \left[ I^{\alpha-1} \circ \left(\mathcal{R} \cdot D^{\alpha-1, \operatorname{rot}}\right) \right](f)(x) \nonumber \\
& = \frac{2}{\alpha} \left(2^{\alpha-1} \pi^{\frac{d}{2}} \dfrac{\Gamma\left(\frac{\alpha}{2}\right)}{\Gamma(\frac{d-\alpha+2}{2})}\right)^{-1} \int_{\mathbb{R}^d} \dfrac{1}{\|x-y\|^{d+2-\alpha}} \langle x-y; D^{\alpha-1, \operatorname{rot}}(f)(y) \rangle dy. 
\end{align}
\end{lem}

\begin{proof}
The proof of this technical lemma is very similar to the computation of the Fourier transform of the potential Riesz kernel. Let $R \geq 1$. Then, using spherical coordinates, for all $ \xi \in \mathbb{R}^d \setminus \{0\}$, 
\begin{align}\label{eq:Fourier_Transform_Truncation}
\mathcal{F}\left(k_{s,R}\right)(\xi) & = \int_{B(0,R)} k_s(x) e^{- i \langle x ; \xi\rangle} dx \nonumber \\
& = \int_{(0,R) \times \mathbb{S}^{d-1}} \frac{ \theta}{r^{1 - s}} e^{- r i \langle \theta; \xi \rangle} dr d\theta. 
\end{align} 
Next, let us fix $r \in (0,R)$ and $\xi \in \mathbb{R}^d \setminus \{0\}$. Then, 
\begin{align*}
\int_{\mathbb{S}^{d-1}} \theta e^{ - i r \langle \theta; \xi \rangle} d\theta & = -i \int_{\mathbb{S}^{d-1}} \theta \sin \left(r \langle \theta ; \xi \rangle\right) d\theta \\
& = \frac{i}{r} \nabla_\xi \left(\int_{\mathbb{S}^{d-1}} \cos\left(\langle r \theta ; \xi \rangle\right) d\theta\right).
\end{align*}
Now, by \cite[formula $B.4$]{G08}, for all $r>0$ and all $\xi \in \mathbb{R}^d \setminus \{0\}$, 
\begin{align}\label{eq:TF_sm}
\int_{\mathbb{S}^{d-1}} \cos \left( r \langle \xi ; \theta \rangle \right) d\theta = \frac{(2\pi)^{\frac{d}{2}}}{r^{\frac{d-2}{2}} \|\xi\|^{\frac{d-2}{2}}} J_{\frac{d}{2}-1}(r\|\xi\|). 
\end{align}
Also, for all $r>0$ and all $\xi \in \mathbb{R}^d \setminus \{0\}$, 
\begin{align*}
\nabla_\xi \left( \|r\xi\|^{1- \frac{d}{2}} J_{\frac{d}{2}-1}(r\|\xi\|) \right) = r \nabla_\xi \left(\|\xi\|\right) \frac{d}{dt} \left(t^{1-\frac{d}{2}}J_{\frac{d}{2}-1}(t)\right) \bigg|_{t = r\|\xi\|}.
\end{align*}
Moreover, by \cite[formula $B.2$,(1)]{G08}, for all $t>0$ and all $d>1$, 
\begin{align*}
\frac{d}{dt} \left(t^{1-\frac{d}{2}}J_{\frac{d}{2}-1}(t)\right)(t) = - t^{1-\frac{d}{2}} J_{\frac{d}{2}}(t),
\end{align*}
which gives, for all $r>0$ and all $\xi \in \mathbb{R}^d \setminus \{0\}$, 
\begin{align*}
\nabla_\xi \left( \|r\xi\|^{1- \frac{d}{2}} J_{\frac{d}{2}-1}(r\|\xi\|) \right) = - r \frac{\xi}{r^{\frac{d}{2}-1} \|\xi\|^{\frac{d}{2}}} J_{\frac{d}{2}}\left(r \|\xi\|\right).  
\end{align*}
Thus, for all $r\in (0,R)$ and $\xi \in \mathbb{R}^d \setminus \{0\}$, 
\begin{align*}
\int_{\mathbb{S}^{d-1}} \theta e^{ - i r \langle \theta; \xi \rangle} d\theta = (2\pi)^{\frac{d}{2}} \frac{-i\xi}{r^{\frac{d}{2}-1} \|\xi\|^{\frac{d}{2}}} J_{\frac{d}{2}}\left(r \|\xi\|\right). 
\end{align*}
Plugging this identity into \eqref{eq:Fourier_Transform_Truncation} gives, for all $R>0$ and all $\xi \in \mathbb{R}^d \setminus \{0\}$, 
\begin{align*}
\mathcal{F}\left(k_{s,R}\right)(\xi) = (2\pi)^{\frac{d}{2}} \frac{-i\xi}{\|\xi\|^{\frac{d}{2}}} \int_0^R J_{\frac{d}{2}}\left(r \|\xi\|\right) \frac{dr}{r^{\frac{d}{2}-s}}.
\end{align*}
Finally, by \cite[formula $10.22.43$]{OLBC_nist10} with $s < (d+1)/2$, 
\begin{align*}
\int_0^{+\infty} J_{\frac{d}{2}}\left(t\right) \frac{dt}{t^{\frac{d}{2}-s}} = 2^{s - \frac{d}{2}} \dfrac{\Gamma\left(\frac{s+1}{2}\right)}{\Gamma(\frac{d-s+1}{2})}.  
\end{align*}
So, for all $\xi \in \mathbb{R}^d \setminus \{0\}$, 
\begin{align*}
\mathcal{F}(k_s)(\xi) = 2^{s} \dfrac{\Gamma\left(\frac{s+1}{2}\right)}{\Gamma(\frac{d-s+1}{2})}\pi^{\frac{d}{2}} \frac{(-i\xi)}{\|\xi\|^{s+1}}.  
\end{align*}
Standard Fourier analysis arguments combined with \eqref{eq:TF_CRiesz_Riesz} give the integral representation \eqref{eq:spatial_rep_comistefani}. This concludes the proof of the lemma. 
\end{proof} 

\begin{rem}\label{rem:FFTC_dimension1}
For $d=1$, by Fourier analytic arguments, it is straightforward to prove that, for all $f \in \mathcal{C}_c^{\infty}(\mathbb{R})$ and all $x\in \mathbb{R}$, 
\begin{align}\label{eq:FFTC_dimension1}
f(x) = \dfrac{1}{\alpha \Gamma(\alpha-1) \sin \left(\frac{\pi(\alpha-1)}{2}\right)} \int_{\mathbb{R}} \dfrac{x-y}{|x-y|^{3-\alpha}} D^{\alpha-1,\operatorname{rot}}(f)(y)dy,
\end{align}
where $D^{\alpha-1,\operatorname{rot}}$ is the fractional gradient operator associated with $\nu_\alpha^{\operatorname{rot}}$. 
\end{rem}
\noindent
Next, let $m_{\alpha,d}$ be the singular vector-valued symbol defined by
\begin{align}\label{eq:anisotropic_composition_kernel}
m_{\alpha,d}(\xi) = \dfrac{-i \xi}{\|\xi\|_\alpha^\alpha}, \quad \xi \in \mathbb{R}^d \setminus \{0\}.
\end{align}
Note that the function $m_{\alpha,d}$ defines a homogeneous tempered distribution since $m_{\alpha,d} \in L^1_{\operatorname{loc}}(\mathbb{R}^d,dx)$ denoted by $M_{\alpha,d}$. The next technical lemma investigates the properties of the Fourier transform (in the sense of tempered distributions) of the kernel $m_{\alpha,d}$. 

\begin{lem}\label{lem:l_alpha_ball}
Let $d \geq 2$ be an integer and let $\alpha \in (1,2)$.~Let $k_{\alpha,d}$ be the inverse Fourier transform of the symbol $m_{\alpha,d}$. Then, 
\begin{enumerate}
\item for all $x \in \mathbb{R}^d$ with $x_k \ne 0$, $k \in \{1, \dots, d\}$, and all $j \in \{1, \dots, d\}$, 
\begin{align}\label{eq:int_representation_stable_density}
k^j_{\alpha,d}(x) = - \frac{\alpha}{(2\pi)^d} \int_0^{+\infty} t^{d-\alpha} \bigg[\prod_{k \ne j} \gamma_\alpha(t x_k) \bigg]\gamma'_{\alpha}(t x_j)dt,
\end{align}
where $\gamma_\alpha$ is the Fourier transform of the function $t \mapsto \exp\left(- |t|^\alpha\right)$, $t \in \bbr$; 
\item $k_{\alpha,d}$ is homogeneous of degree $\alpha-d-1$ and $k_{\alpha,d} \in L^1_{\operatorname{loc}}(\mathbb{R}^d,dx)$; 
\item for all $f \in \mathcal{C}_c^{\infty}(\mathbb{R}^d)$ and all $x \in \mathbb{R}^d$, 
\begin{align}\label{eq:formula_Stein_axes}
f(x) = \frac{1}{\alpha} \int_{\mathbb{R}^d} \langle k_{\alpha,d}(x-y) ; D^{\alpha-1,d}(f)(y) \rangle dy,
\end{align}
where $D^{\alpha-1,d}$ is the fractional gradient operator associated with the $\alpha$-stable probability measure $\mu_{\alpha,d}$ such that $\widehat{\mu_{\alpha,d}}(\xi) = \exp (- \|\xi\|_\alpha^\alpha)$, $\xi \in \mathbb{R}^d$. 
\end{enumerate}
\end{lem}

\begin{proof}
First, recall that the Fourier transform of $f_\alpha$ defined, for all $x \in \mathbb{R}^d \setminus \{0\}$, by 
\begin{align}\label{eq:Koldobsky_particular_case}
f_\alpha(x) = \frac{1}{\|x\|_\alpha^\alpha}, 
\end{align}
is given, for all $\xi \in \mathbb{R}^d$ with $\xi_k \ne 0$, $k \in \{1, \dots, d\}$, by 
\begin{align}\label{eq:FT_Koldobsky_particular_case}
\mathcal{F}\left(f_\alpha\right)(\xi) = \alpha \int_0^{+\infty} t^{d-\alpha-1} \prod_{k=1}^d \gamma_\alpha(t \xi_k) dt,
\end{align}
where $\gamma_\alpha$ is the Fourier transform of the function $t \mapsto \exp\left(- |t|^\alpha\right)$, $t \in \bbr$ (see \cite[Lemma $8$]{Koldo_AM98}).~Using standard calculus rules for the Fourier transform, for all $\xi \in \mathbb{R}^d$, $\xi_k \ne 0$, $k \in \{1, \dots, d\}$, 
\begin{align}\label{eq:FT_malphad}
\mathcal{F}\left(m_{\alpha,d}\right)(\xi) = \nabla_{\xi} \left(\mathcal{F} \left(f_\alpha\right)(\xi)\right). 
\end{align}
Then, thanks to \eqref{eq:FT_Koldobsky_particular_case}, for all $\xi \in \mathbb{R}^d$, $\xi_k \ne 0$, $k \in \{1, \dots, d\}$, and all $j \in \{1, \dots, d\}$, 
\begin{align*}
\partial_{\xi_j} \left(\mathcal{F}(f_\alpha)(\xi)\right) = \alpha \int_0^{+\infty} t^{d-\alpha} \bigg[\prod_{k \ne j} \gamma_\alpha(t \xi_k) \bigg]\gamma'_{\alpha}(t\xi_j)dt. 
\end{align*}
Note that the integral above is absolutely convergent since, for all $t \in \bbr$, 
\begin{align*}
\gamma_{\alpha}(t) \leq \frac{c_{\alpha,1}}{(1+|t|)^{\alpha+1}}, \quad |\gamma_\alpha'(t)| \leq \frac{c_{\alpha,2}}{(1+|t|)^{\alpha+2}}, \quad c_{\alpha,1},c_{\alpha,2}>0.
\end{align*}
Next, for all $j \in \{1, \dots, d\}$, let $k_{\alpha,d}^j$ be defined, for all $x \in \mathbb{R}^d$ with $x_k \ne 0$, $k \in \{1, \dots, d\}$, by 
\begin{align*}
k_{\alpha,d}^j(x) = - \frac{\alpha}{(2\pi)^d} \int_0^{+\infty} t^{d-\alpha} \bigg[\prod_{k \ne j} \gamma_\alpha(t x_k) \bigg]\gamma'_{\alpha}(t x_j)dt.
\end{align*}
From this integral representation, it is clear that $k_{\alpha,d} = ( k^1_{\alpha,d} , \dots, k^d_{\alpha,d})$ is homogeneous of degree $\alpha-d-1$. Namely, for all $\lambda>0$ and all $j \in \{1,\dots, d\}$, 
\begin{align}\label{eq:homo_kalphad}
k_{\alpha,d}^j(\lambda x) = \left(\frac{1}{\lambda}\right)^{d-\alpha+1} k_{\alpha,d}^j(x). 
\end{align}
Now, let $Q_1$ be the open unit cube of $\mathbb{R}^d$, i.e., $Q_1 = (0,1)^d$. Then, for all $t \in (0,+\infty)$ and all $x \in \mathbb{R}^d$ with $x_k \ne 0$, $k \in \{1, \dots, d\}$, 
\begin{align*}
\int_{Q_1} \prod_{k \ne j} \gamma_\alpha(t x_k) |\gamma_\alpha'(t x_j)| dx = \left(\int_0^1 \gamma_\alpha(t x) dx\right)^{d-1} \int_0^1 |\gamma_\alpha'(t x)| dx. 
\end{align*}
But, 
\begin{align*}
\int_0^1 |\gamma_\alpha'(t x)| dx \leq c_\alpha \int_0^1 \frac{dx}{\left(1+ t x\right)^{\alpha+2}} =  \frac{c_\alpha}{t(\alpha+1)} \left(1 - \frac{1}{(1+t)^{\alpha+1}}\right). 
\end{align*}
Moreover, 
\begin{align*}
\int_0^1 \gamma_\alpha(tx)dx \leq \frac{c_\alpha}{\alpha t} \left(1 - \frac{1}{(1+t)^\alpha}\right). 
\end{align*}
Thus, for $t<1$, the integral of $\prod_{k \ne j} \gamma_\alpha(t x_k) |\gamma_\alpha'(t x_j)|$ over the unit cube $Q_1$ is uniformly bounded in $t$ and, for $t>1$, 
\begin{align*}
\int_{Q_1} \prod_{k \ne j} \gamma_\alpha(t x_k) |\gamma_\alpha'(t x_j)| dx \leq \frac{c_{\alpha,d}}{t^d},
\end{align*} 
for some $c_{\alpha,d}>0$ depending only on $\alpha$ and $d$. Since $\alpha \in (1,2)$, for all $j \in \{1, \dots, d\}$, $k^j_{\alpha,d} \in L^1_{\operatorname{loc}}(\mathbb{R}^d,dx)$. Finally, the identity \eqref{eq:formula_Stein_axes} follows by standard Fourier arguments.
\end{proof}
\noindent
The next lemma provides find sharp pointwise bounds on $k_{\alpha,d}$ and on $\partial_\ell (k_{\alpha,d})$, with $\ell \in \{1,\dots,d\}$, based on the integral representation \eqref{eq:int_representation_stable_density}. 

\begin{lem}\label{lem:sharp_pointwise_bound}
Let $d \geq 2$ be an integer, let $\alpha \in (1,2)$ and let $k_{\alpha,d}$ be the inverse Fourier transform of the symbol $m_{\alpha,d}$. Then, 
\begin{enumerate}
\item for all $x \in \mathbb{R}^d$, $x_k \ne 0$, $k \in \{1, \dots, d\}$ and all $j \in \{1, \dots, d\}$, 
\begin{align}\label{ineq:pointwise_bound_tensorized_kernel}
\left| k^j_{\alpha,d}(x) \right| \leq c_{\alpha,d} \prod_{k = 1}^d |x_k|^{- \frac{d+1-\alpha}{d}},
\end{align} 
for some $c_{\alpha,d}>0$ depending on $\alpha$ and $d$ only;
\item for all $x \in \mathbb{R}^d$, $x_k \ne 0$, $k \in \{1, \dots, d\}$ and all $j,\ell \in \{1, \dots, d\}$, 
\begin{align}\label{ineq:pointwise_bound_tensorized_Dkernel}
\left| \partial_{\ell} k^j_{\alpha,d}(x) \right| \leq C_{\alpha,d} \prod_{k = 1}^d |x_k|^{- \frac{d+2-\alpha}{d}},
\end{align}
for some $C_{\alpha,d}>0$ depending on $\alpha$ and $d$ only; 
\item assuming that $\alpha < (d+1)/2$, for all $x \in \mathbb{R}^d$, $x_k \ne 0$, $k \in \{1, \dots, d\}$ and all $j,\ell \in \{1, \dots, d\}$, 
\begin{align}\label{ineq:pointwise_bound_mixed_Dkernel}
\left| \partial_{\ell} (k_{\alpha,d}^j )(x)  \right| \leq A_{\alpha,d} \frac{1}{\|x\|^{\frac{\alpha(d+1)-2}{d-1}}_1} \prod_{k = 1}^d |x_k|^{-\frac{d+1-2\alpha}{d-1}}, 
\end{align}
for some $A_{\alpha,d}>0$ depending on $\alpha$ and $d$ only; 
\item for $d = 2$, 
\begin{align}\label{ineq:L1_increments_mixedhyperrectangles}
\underset{\omega \in \mathbb{S}}{\sup} \int_{(0,2) \times (2,+\infty)} \|k_{\alpha,2}(z+ \omega) - k_{\alpha,2}(z) \| dz < + \infty;
\end{align}
\item for all $d \geq 3$, 
\begin{align}\label{eq:L1_increments_axes_fulld}
\sup_{\omega \in \mathbb{S}^{d-1}} \int_{\mathbb{R}^d} \left\| k_{\alpha,d}(\omega+z) - k_{\alpha,d}(z)  \right\| dz <+\infty.
\end{align}
\end{enumerate}
\end{lem}

\begin{proof}
The proofs of the inequalities \eqref{ineq:pointwise_bound_tensorized_kernel}, \eqref{ineq:pointwise_bound_tensorized_Dkernel} and \eqref{ineq:pointwise_bound_mixed_Dkernel} rely on a truncation-optimization argument.~In the sequel, $c_{\alpha,d}$ denotes a positive constant which might change from a line to another. Let $R>0$ be a parameter to be fixed later on. From Lemma \ref{lem:l_alpha_ball}, for all $x \in \mathbb{R}^d$, $x_k \ne 0$, $k \in \{1,\dots, d\}$, 
\begin{align*}
\left| k_{\alpha,d}^j(x) \right| \leq c_{\alpha,d} \int_0^{+\infty} t^{d-\alpha} \left| \partial_j(p_{\alpha,d})(tx)\right| dt, 
\end{align*}
where $p_{\alpha,d}$ is the Lebesgue density of the probability measure $\mu_{\alpha,d}$, and from the product structure of $\mu_{\alpha,d}$ and standard estimates on $p_{\alpha,1}$, for all $j \in \{1, \dots, d\}$, 
\begin{align*}
\underset{x \in \mathbb{R}^d}{\sup} \dfrac{\left| \partial_j(p_{\alpha,d})(x)\right|}{p_{\alpha,d}(x)} <+\infty. 
\end{align*} 
Then, for all $x \in \mathbb{R}^d$, $x_k \ne 0$, $k \in \{1,\dots, d\}$,
\begin{align*}
\left| k_{\alpha,d}^j(x) \right| &\leq c_{\alpha,d} \int_0^{+\infty} t^{d-\alpha} p_{\alpha,d}(tx) dt \\
& \leq c_{\alpha,d} \int_{0}^{+\infty} t^{d-\alpha} \dfrac{dt}{\prod_{k=1}^d \left(1+t |x_k|\right)^{\alpha+1}} \\
& \leq c_{\alpha,d} \int_0^{+\infty} t^{d-\alpha} \dfrac{dt}{\left(1+ \|x\|_1 t + \dots + \prod_k |x_k| t^d\right)^{\alpha+1}} \\
& \leq c_{\alpha,d} \|x\|^{\alpha-1-d}_1 \int_0^{+\infty} t^{d-\alpha} \dfrac{dt}{\left(1+ t + \dots + \frac{\prod_k |x_k|}{\|x\|^d_1} t^d\right)^{\alpha+1}},
\end{align*}
where $\|x\|_1 = \sum_k |x_k|$. Now, let us cut the integral on the right-hand side of the previous inequality using the truncation parameter $R$. Then, for all $x \in \mathbb{R}^d$, $x_k \ne 0$, $k \in \{1,\dots, d\}$, 
\begin{align*}
\left| k_{\alpha,d}^j(x) \right| & \leq c_{\alpha,d} \left( \|x\|^{\alpha-1-d}_1 \int_0^R t^{d-\alpha} dt + \|x\|^{\alpha-1-d}_1 \int_{R}^{+\infty} \dfrac{dt}{t^{\alpha(d+1)}} \dfrac{\|x\|^{d(\alpha+1)}_1}{\prod_k |x_k|^{\alpha+1}} \right) \\
& \leq c_{\alpha,d} \left(  \|x\|^{\alpha-1-d}_1 R^{d+1-\alpha} +  \dfrac{\|x\|^{d\alpha + \alpha-1}_1}{\prod_k |x_k|^{\alpha+1}} \frac{1}{R^{d\alpha+ \alpha-1}} \right). 
\end{align*}
Next, choose $R = \|x\|_1 / \prod_k |x_k|^{\frac{1}{d}}$.~The end of the proof of \eqref{ineq:pointwise_bound_tensorized_kernel} follows by standard computations.~Next, let us explain the main steps in the proofs of \eqref{ineq:pointwise_bound_tensorized_Dkernel} and \eqref{ineq:pointwise_bound_mixed_Dkernel}. Thanks to the integral representation \eqref{eq:int_representation_stable_density}, for all $x \in \mathbb{R}^d$, $x_k \ne 0$, $k \in \{1, \dots, d\}$, 
\begin{align*}
\left| \partial_\ell k^j_{\alpha,d}(x) \right| \leq c_{\alpha,d} \int_0^{+\infty} t^{d-\alpha+1} \left| \partial^2_{\ell,j}(p_{\alpha,d})(tx) \right| dt. 
\end{align*} 
Next, for all $t>0$ and all $x \in \mathbb{R}^d$, $x_k \ne 0$, $k \in \{1, \dots, d\}$, 
\begin{align*}
\left| \partial_{\ell, j}^2(p_{\alpha,d})(tx) \right| & \leq c_{\alpha,d} \prod_{k = 1}^d \dfrac{1}{\left(1+ t |x_k|\right)^{\alpha+1}} \\
& \leq \dfrac{c_{\alpha,d}}{\left(1+ t \|x\|_1 + \dots + \prod_k |x_k| t^d \right)^{\alpha+1}}. 
\end{align*}
Thus, for all $x \in \mathbb{R}^d$, $x_k \ne 0$, $k \in \{1, \dots, d\}$, 
\begin{align*}
\left| \partial_{\ell} k^j_{\alpha,d}(x) \right| & \leq c_{\alpha,d} \int_0^{+\infty} t^{d-\alpha+1} \dfrac{dt}{\left(1+ t \|x\|_1 + \dots + \prod_k |x_k| t^d \right)^{\alpha+1}} \\
& \leq c_{\alpha,d} \|x\|_1^{\alpha-d-2} \int_0^{+\infty} t^{d-\alpha+1} \dfrac{dt}{\left(1+ t + \dots + \frac{\prod_k |x_k|}{\|x\|^d_1} t^d \right)^{\alpha+1}}.
\end{align*}
As previously, let $R>0$ be a truncation parameter to be fixed later on.~Then, by standard computations, 
\begin{align*}
\left| \partial_{\ell} k^j_{\alpha,d}(x) \right| & \leq c_{\alpha,d} \left(\|x\|_1^{\alpha-d-2} R^{d-\alpha+2} + \frac{1}{R^{d \alpha + \alpha-2}} \dfrac{\|x\|^{\alpha(d+1)-2}_1}{\prod_k |x_k|^{\alpha+1}} \right).
\end{align*}
Now, choose again $R = \|x\|_1 / \prod_k |x_k|^{1/d}$. Then, for all $x \in \mathbb{R}^d$, $x_k \ne 0$, $k \in \{1, \dots, d\}$, 
\begin{align*}
\left| \partial_{\ell} k^j_{\alpha,d}(x) \right| \leq c_{\alpha,d} \dfrac{1}{\prod_k |x_k|^{\frac{d-\alpha+2}{d}}}. 
\end{align*}
To prove \eqref{ineq:pointwise_bound_mixed_Dkernel}, let us proceed as previously. Let $R>0$ to be fixed later on. Then, for all $x \in \mathbb{R}^d$, $x_k \ne 0$, $k \in \{1, \dots, d\}$, and all $j , \ell \in \{1, \dots, d\}$, 
\begin{align*}
|\partial_\ell k_{\alpha,d}^j(x) | & \leq c_{\alpha,d} \int_0^{+\infty} t^{d-\alpha+1} \dfrac{dt}{\left(1+ \|x\|_1 t + \dots + \prod_{k} |x_k| t^d\right)^{\alpha+1}} \\
& \leq c_{\alpha,d} \left( \|x\|_1^{-(\alpha+1)} \int_0^R \dfrac{dt}{t^{2\alpha-d}} + \prod_k |x_k|^{-\alpha-1}\int_R^{+\infty} \dfrac{dt}{t^{d\alpha+\alpha-1}} \right)  \\
& \leq c_{\alpha,d} \left(  \|x\|_1^{-(\alpha+1)}R^{d-2\alpha+1} + \prod_k |x_k|^{-(\alpha+1)} \frac{1}{R^{d\alpha+\alpha-2}}\right).
\end{align*}
Choosing $R = (\|x\|_1/\prod_k |x_k|)^{1/(d-1)}$, 
\begin{align*}
|\partial_\ell k_{\alpha,d}^j(x) | & \leq c_{\alpha,d} \dfrac{1}{\|x\|_1^{\frac{\alpha d+\alpha-2}{d-1}}} \dfrac{1}{\prod_k |x_k|^{\frac{d+1-2\alpha}{d-1}}}. 
\end{align*}
Next, let us prove the finiteness of the supremum \eqref{ineq:L1_increments_mixedhyperrectangles}. For all $\omega=(\omega_1 , \omega_2) \in \mathbb{S} := \mathbb{S}^{1}$, 
\begin{align*}
\int_{(0,2) \times (2,+\infty)} \|k_{\alpha,2}(x+\omega) - k_{\alpha,2}(x)\|dx &\leq \int_{(0,2) \times (2,+\infty)}  |k^1_{\alpha,2}(x+\omega) - k^1_{\alpha,2}(x)|dx \\
& \quad\quad +  \int_{(0,2) \times (2,+\infty)}  |k^2_{\alpha,2}(x+\omega) - k^2_{\alpha,2}(x)|dx.
\end{align*}
Without loss of generality, let us treat the first term of the sum on the right-hand side of the previous inequality.~The product structure of the probability measure $\mu_{\alpha,2}$ gives: 
\begin{align*}
\int_{(0,2) \times (2,+\infty)}  |k^1_{\alpha,2}(x+\omega) - k^1_{\alpha,2}(x)|dx \leq & \alpha \int_{(0,2) \times (2,+\infty)} \int_0^{+\infty} t^{2-\alpha} \bigg| \partial_1\left(p_{\alpha,1}\right)(tx_1+ t\omega_1) p_{\alpha,1}(t(x_2+\omega_2)) \\
& \quad\quad - \partial_1\left(p_{\alpha,1}\right)(tx_1) p_{\alpha,1}(tx_2) \bigg| dt dx_1dx_2 \\
& \leq A_1(\omega) + A_2(\omega), 
\end{align*}
where
\begin{align*}
& A_1(\omega) =  \alpha \int_{(0,2) \times (2,+\infty)} \int_0^{+\infty} t^{2-\alpha} \bigg| \partial_1\left(p_{\alpha,1}\right)(tx_1+ t\omega_1) p_{\alpha,1}(t(x_2+\omega_2)) \bigg|  dt dx_1dx_2,\\
& A_2(\omega) =  \alpha \int_{(0,2) \times (2,+\infty)} \int_0^{+\infty} t^{2-\alpha} \bigg| \partial_1\left(p_{\alpha,1}\right)(tx_1)p_{\alpha,1}(tx_2) \bigg|  dt dx_1dx_2.
\end{align*}
Let us start with $A_1$. So, for all $t \in (0,+\infty)$ and all $\omega=(\omega_1 , \omega_2) \in \mathbb{S}$, 
\begin{align*}
\int_{(2,+\infty)} \left| p_{\alpha,1}\left(tx_2 + t \omega_2\right) \right| dx_2 & \leq C_\alpha \int_{(2, +\infty)} \dfrac{dx_2}{\left(1+t |x_2+\omega_2|\right)^{\alpha+1}}  \\
& \leq C_\alpha \int_{(2, +\infty)} \dfrac{dx_2}{\left(1+t (x_2+\omega_2)\right)^{\alpha+1}}  \\
& \leq \frac{C_\alpha}{\alpha t} \dfrac{1}{\left(1+t \left(2+ \omega_2\right)\right)^\alpha}  \\
& \leq \frac{C_\alpha}{\alpha t} \dfrac{1}{\left(1+t\right)^\alpha}.
\end{align*}
Moreover, for all $t \in (0,+\infty)$ and all $\omega=(\omega_1 , \omega_2) \in \mathbb{S}$, 
\begin{align*}
\int_{(0,2)} dx_1 \left| \partial_1(p_{\alpha,1})\left(tx_1 + t\omega_1\right)\right|  & = \int_{(0,2)+\omega_1} dx_1 \left| \partial_1(p_{\alpha,1})\left(tx_1 \right)\right| \leq \int_{]-1,3[} dx_1 \left| \partial_1(p_{\alpha,1})\left(tx_1 \right)\right| \\
& \leq C_\alpha \int_{]-1,3[} \dfrac{dx_1}{\left(1+t |x_1|\right)^{\alpha+2}}  \\
& \leq  \dfrac{C_\alpha}{(\alpha+1)t} \left( 1- \dfrac{1}{ (1+t)^{\alpha+1}} + 1 - \dfrac{1}{ (1+3t)^{\alpha+1}} \right). 
\end{align*}
Then, for all $\omega= (\omega_1, \omega_2) \in \mathbb{S}$,  
\begin{align*}
A_1(\omega) &\leq C_\alpha \int_0^{+\infty} t^{-\alpha}\dfrac{1}{\left(1+t\right)^\alpha} \left( 1- \dfrac{1}{ (1+t)^{\alpha+1}} + 1 - \dfrac{1}{ (1+3t)^{\alpha+1}} \right) dt <+\infty. 
\end{align*}
Next, let us deal with $A_2$. So, for the integral with respect to the $x_1$-variable, for all $t>0$, 
\begin{align*}
\int_{(0,2)} \left| \partial_1(p_{\alpha,1})\left(tx_1\right) \right| dx_1 & \leq C_\alpha \int_{(0,2)} \dfrac{dx_1}{\left(1+t x_1\right)^{\alpha+2}} \\
& \leq \frac{C_\alpha}{t} \left(1- \dfrac{1}{\left(1+2t\right)^{\alpha+1}} \right).
\end{align*}
Finally, for the integral with respect to the $x_2$-variable, for all $t>0$,
\begin{align*}
\int_{(2, +\infty)} \left| p_{\alpha,1}\left(tx_2\right) \right| dx_2 & \leq C_\alpha \int_{(2,+\infty)} \dfrac{dx_2}{\left(1+t |x_2|\right)^{\alpha+1}} \\
& \leq \frac{C_\alpha}{t} \frac{1}{\left(1+2t\right)^{\alpha}}, 
\end{align*}
for some $C_\alpha>0$ depending only on $\alpha$.~Thus, for all $\omega = \left(\omega_1, \omega_2\right) \in \mathbb{S}$, 
\begin{align*} 
A_2(\omega) \leq C_\alpha \int_0^{+\infty} t^{- \alpha}\left(1- \dfrac{1}{\left(1+2t\right)^{\alpha+1}} \right)\frac{1}{\left(1+2t\right)^{\alpha}}dt <+\infty.
\end{align*}
To finish, let us prove the finiteness of the supremum \eqref{eq:L1_increments_axes_fulld}. Recall that $k^j_{\alpha,d} \in L^1_{\operatorname{loc}}(\mathbb{R}^d,dx)$, for all $j \in \{1, \dots, d\}$. In particular, with $Q_2 = (0,2)^d$, for all $\nu \in \mathbb{S}^{d-1}$, 
\begin{align*}
\int_{Q_2} \left\| k_{\alpha,d}(\omega +\nu) - k_{\alpha,d}(\omega) \right\| d\omega & \leq \int_{Q_2} \left\| k_{\alpha,d}(\omega) \right\| d\omega + \int_{Q_2} \left\| k_{\alpha,d}(\omega +\nu) \right\| d\omega  \\
& \leq \int_{Q_2} \left\| k_{\alpha,d}(\omega) \right\| d\omega + \int_{Q_2 + \nu} \left\| k_{\alpha,d}(\omega) \right\| d\omega.
\end{align*} 
Now, it is clear that there exists a compact subset $K \subset \mathbb{R}^d$ such that $Q_2 + \nu \subset K$, for all $\nu \in \mathbb{S}^{d-1}$. Thus,  
\begin{align*}
\int_{Q_2} \left\| k_{\alpha,d}(\omega +\nu) - k_{\alpha,d}(\omega) \right\| d\omega \leq \int_{Q_2} \left\| k_{\alpha,d}(\omega) \right\| d\omega + \int_{K} \left\| k_{\alpha,d}(\omega) \right\| d\omega < +\infty.
\end{align*}
Next, let us deal with, 
\begin{align*}
I : = \int_{\bbr_+^d \setminus Q_2} \left\| k_{\alpha,d}(\omega +\nu) - k_{\alpha,d}(\omega) \right\| d\omega. 
\end{align*}
The integral $I$ can be broken down into two terms:
\begin{align*}
I = I_1 + I_2, \quad I_1 = \int_{[2, +\infty[^d} \left\| k_{\alpha,d}(\omega +\nu) - k_{\alpha,d}(\omega) \right\| d\omega , \quad  I_2 = I - I_1. 
\end{align*}
Let us start with $I_1$. Using \eqref{ineq:pointwise_bound_tensorized_Dkernel}, for all $j \in \{1,\dots, d\}$, 
\begin{align*}
\int_{[2, +\infty[^d} \left| k^j_{\alpha,d}(\omega+\nu) - k^j_{\alpha,d}(\omega) \right| d\omega & \leq \int_0^1  \int_{[2,+\infty[^d}  \left| \langle \nabla (k^j_{\alpha,d})(\omega+s\nu) ; \nu \rangle  \right| d\omega ds  \\
& \leq \sum_{\ell}  \int_0^1  \int_{[2,+\infty[^d} \left| \partial_\ell (k_{\alpha,d}^j)(\omega+s\nu) \right| d\omega ds  \\
& \leq c_{\alpha,d} \sum_{\ell} \int_0^1 \int_{[2, +\infty[^d} \prod_{k = 1}^d |\omega_k +s \nu_k|^{- \frac{d-\alpha+2}{d}} d\omega ds  \\
& \leq c_{\alpha,d} \int_0^1 \left(\int_{[2, +\infty[} \frac{d\omega}{(\omega-s)^{\frac{d-\alpha+2}{d}}} \right)^d ds \\
& \leq c_{\alpha,d} \int_0^1 \dfrac{ds}{(2-s)^{2-\alpha}} < +\infty. 
\end{align*}
It remains to deal with $I_2$ which contains integrals over mixed hyperrectangles, where by mixed hyperrectangles, we mean hyperrectangles formed by intervals of the types $(0,2)$ and $(2, +\infty)$. Next, let us consider the $d$-dimensional rectangle $(0,2)^{d-1} \times (2, +\infty)$, with $d\geq 3$.~Then, using \eqref{ineq:pointwise_bound_mixed_Dkernel}, for all $j \in \{1, \dots, d\}$, 
\begin{align*}
\int_{(0,2)^{d-1} \times (2, +\infty)} \left| k^j_{\alpha,d}(\omega + \nu) - k^j_{\alpha,d}(\omega)\right| d\omega & \leq \int_0^1 \int_{(0,2)^{d-1} \times (2, +\infty)} \left| \langle \nabla(k^j_{\alpha,d})(\omega+s\nu) ; \nu \rangle \right| d\omega ds \\
& \leq \sum_{\ell} \int_0^1 \int_{(0,2)^{d-1} \times (2,+\infty)} \left| \partial_\ell (k_{\alpha,d}^j)(\omega+s\nu) \right| d\omega ds \\
& \leq c_{\alpha,d} \int_0^1 \bigg(\int_{(0,2)^{d-1} \times (2,+\infty)} \frac{1}{\| \omega + s\nu \|^{\frac{\alpha(d+1)-2}{d-1}}_1} \\
& \quad\quad\quad \times \prod_{k = 1}^d |\omega_k + s \nu_k|^{-\frac{d+1-2\alpha}{d-1}} d\omega\bigg) ds \\
& \leq c_{\alpha,d} \int_0^1 \left(\prod_{k = 1}^{d-1} \int_{(0,2)} \frac{d\omega}{|\omega + s \nu_k|^{\frac{d+1-2\alpha}{d-1}}} \right) \\
&\quad\quad\quad \times \left(\int_{[2, +\infty)} \dfrac{d\omega}{|\omega+s\nu_d|^{\alpha+1}}\right)ds \\
& \leq c_{\alpha,d} \int_0^1 \left(\int_{(-1,0)\cup(0,3)} \frac{d\omega}{|\omega|^{\frac{d+1-2\alpha}{d-1}}}\right)^{d-1} \\
&\quad\quad\quad \times \left(\int_{[1,+\infty)} \frac{d\omega}{|\omega|^{\alpha+1}}\right)ds <+\infty. 
\end{align*}
The other cases can be treated in a similar way. This concludes the proof of the lemma. 
\end{proof}
\noindent
Let us now consider the general non-degenerate symmetric $\alpha$-stable situation, with $\alpha \in (1,2)$. For all $T>0$, let $k_{\alpha,T}$ be the $\mathbb{R}^d$-valued function defined, for all $x \in \mathbb{R}^d$, by
\begin{align}\label{eq:truncated_convolution_NDSkernel}
k_{\alpha,T}(x) = -\alpha \int_0^T t^{d-\alpha} \nabla(p_\alpha)\left(tx\right)dt , \quad d \geq 2. 
\end{align}

\begin{lem}\label{lem:truncated_convolution_NDSkernel}
Let $d \geq 2$ be an integer, let $\alpha \in (1,2)$ and let $T>0$. Then, 
\begin{itemize}
\item for all $\xi \in \bbr^d$ such that $\xi \ne 0$, 
\begin{align}\label{eq:FT_truncated_convolution_NDSkernel}
\mathcal{F}\left(k_{\alpha,T}\right)(\xi) = \dfrac{(-i\xi)}{\sigma_\alpha(\xi)^\alpha} \exp \left( - \dfrac{\sigma_\alpha(\xi)^\alpha}{T^\alpha}\right), 
\end{align}
where $\sigma_\alpha$ is given by \eqref{eq:rep_spectral_measure};
\item $\mathcal{F}(k_{\alpha,T}) \in L^p(\bbr^d,dx)$, with $p \in [1,d/(\alpha-1))$;  
\item for all $f \in \mathcal{C}_c^{\infty}(\mathbb{R}^d)$ and all $x \in \mathbb{R}^d$, 
\begin{align}\label{eq:frac_FTC}
\underset{T \rightarrow +\infty}{\lim} \frac{1}{\alpha} \int_{\bbr^d} \langle k_{\alpha,T}(x-y) ; D^{\alpha-1}(f)(y) \rangle dy = f(x). 
\end{align}
\end{itemize}
\end{lem}

\begin{proof}
By Fubini's theorem and changing variables, for all $T>0$ and all $\xi \in \mathbb{R}^d$ such that $\xi \ne 0$, 
\begin{align*}
\mathcal{F} \left(k_{\alpha,T}\right)(\xi) & = -\alpha \int_0^T t^{d-\alpha} \mathcal{F}\left(\nabla(p_\alpha)\right)(\xi/t) \frac{dt}{t^d}\\
& = (-i\xi) \alpha \int_0^T t^{-\alpha-1} \exp \left(- \frac{\sigma_\alpha(\xi)^\alpha}{t^{\alpha}}\right) dt \\
& = (-i\xi)  \int_{\frac{1}{T^{\alpha}}}^{+\infty} \exp\left(-u \sigma_\alpha(\xi)^\alpha \right) du \\
& = \frac{(-i\xi)}{\sigma_\alpha(\xi)^\alpha}  \int_{\frac{\sigma_\alpha(\xi)^\alpha}{T^{\alpha}}}^{+\infty} \exp\left(-v \right) dv = \dfrac{(-i\xi)}{\sigma_\alpha(\xi)^\alpha} \exp \left( - \dfrac{\sigma_\alpha(\xi)^\alpha}{T^\alpha}\right). 
\end{align*}
Now, for all $p \in [1, + \infty)$ and all $\xi \in \mathbb{R}^d$ such that $\xi \ne 0$, 
\begin{align*}
\left\| \mathcal{F} \left(k_{\alpha,T}\right)(\xi) \right\|^p & \leq \dfrac{\|\xi\|^p}{ \|\xi\|^{\alpha p} \sigma_\alpha(e_\xi)^{\alpha p}} \exp \left( - p \frac{\|\xi\|^\alpha \sigma_\alpha(e_\xi)^\alpha}{T^\alpha}\right).
\end{align*}
Since $\mu_\alpha$ is non-degenerate, $\underset{e \in \mathbb{S}^{d-1}}{\inf} \sigma_\alpha(e)>0$. Thus, for all $p \in [1, + \infty)$ and all $\xi \in \mathbb{R}^d \setminus \{0\}$, 
\begin{align*}
\left\| \mathcal{F} \left(k_{\alpha,T}\right)(\xi) \right\|^p & \leq  C^1_{\alpha,d,p}\|\xi\|^{(1-\alpha)p} \exp\left( - \frac{C^2_{\alpha,d,p}  \|\xi\|^\alpha}{T^\alpha}\right), 
\end{align*}
for some $C^1_{\alpha,d,p},C^2_{\alpha,d,p}>0$ depending on $\alpha$, $p$ and $d$. Moreover, for all $T\in (0,+\infty)$, 
\begin{align*}
\int_{\mathbb{R}^d} \|\xi\|^{(1-\alpha)p} \exp\left( - \frac{C^2_{\alpha,d,p}  \|\xi\|^\alpha}{T^\alpha}\right) d\xi <+\infty \Leftrightarrow p < \frac{d}{\alpha-1}. 
\end{align*}
Finally, by the Fourier inversion formula, for all $f \in \mathcal{C}_c^{\infty}(\bbr^d)$, all $T>0$ and all $x \in \mathbb{R}^d$, 
\begin{align*}
\frac{1}{\alpha} \int_{\mathbb{R}^d} \langle k_{\alpha,T}(x-y)  ; D^{\alpha-1}(f)(y) \rangle dy & = \frac{1}{(2\pi)^d} \frac{1}{\alpha} \int_{\mathbb{R}^d} \mathcal{F}(f)(\xi) e^{i \langle \xi ; x \rangle } \langle (-i\xi) ; \tau_\alpha(\xi) \rangle \\
&\quad\quad\times \exp \left( - \dfrac{\sigma_\alpha(\xi)^\alpha}{T^\alpha}\right) \frac{d\xi}{\sigma_\alpha(\xi)^\alpha},
\end{align*}
where $\tau_\alpha$ is given, for all $\xi \in \mathbb{R}^d$, by 
\begin{align*}
\tau_\alpha(\xi) = \int_{\mathbb{R}^d} \left(e^{i \langle u ; \xi \rangle}-1\right)u\nu_\alpha(du). 
\end{align*}
Now, making use of polar coordinates and an integration by parts in the radial coordinate, for all $\xi \in \mathbb{R}^d$, 
\begin{align*}
\frac{1}{\alpha} \langle - i \xi ; \tau_\alpha(\xi) \rangle & = \frac{-i}{\alpha} \int_{\bbr^d} \left(e^{i \langle u ; \xi \rangle}-1\right) \langle \xi ;u \rangle \nu_\alpha(du)  \\
& =  \frac{-i}{\alpha} \int_{(0,+\infty) \times \mathbb{S}^{d-1}} \left(e^{i \langle ry ; \xi \rangle}-1\right) \langle y; \xi \rangle \frac{dr}{r^\alpha} \sigma(dy)  \\
& = \frac{-1}{\alpha} \int_{\mathbb{S}^{d-1}} \int_{(0,+\infty)} \dfrac{d}{dr} \left(e^{i \langle ry ; \xi \rangle}-1 - i  \langle r y ; \xi\rangle \right) \frac{dr}{r^\alpha} \sigma(dy)  \\
& = \frac{1}{\alpha} \int_{\mathbb{S}^{d-1}} \int_{(0,+\infty)} \left(e^{i \langle ry ; \xi \rangle}-1 - i  \langle r y ; \xi\rangle \right) \frac{d}{dr} \left(r^{-\alpha}\right) dr \sigma(dy)  \\
& = - \int_{\mathbb{R}^d} \left(e^{i \langle u ; \xi \rangle}-1 - i  \langle u ; \xi\rangle \right) \nu_\alpha(du)  \\
& = \sigma_\alpha(\xi)^{\alpha}. 
\end{align*}
Thus, for all $f \in \mathcal{C}_c^{\infty}(\bbr^d)$, all $T>0$ and all $x \in \mathbb{R}^d$, 
\begin{align*}
\frac{1}{\alpha} \int_{\mathbb{R}^d} \langle k_{\alpha,T}(x-y)  ; D^{\alpha-1}(f)(y) \rangle dy & = \frac{1}{(2\pi)^d} \int_{\mathbb{R}^d} \mathcal{F}(f)(\xi) e^{i \langle \xi ; x \rangle } \exp \left( - \dfrac{\sigma_\alpha(\xi)^\alpha}{T^\alpha}\right)d\xi \\
& = P^{\alpha}_{1/T^\alpha}(f)(x) \longrightarrow f(x),  
\end{align*}
as $T$ tends to $+\infty$. This concludes the proof of the lemma. 
\end{proof}
\noindent
Next, let us introduce standard objects from the potential theory of symmetric $d$-dimensional $\alpha$-stable L\'evy processes. The potential operator $\mathbb{V}_\alpha$ is defined, for all suitable $f$, by 
\begin{align}\label{eq:potential_operator_stable_measure}
\mathbb{V}_\alpha(f) : = \int_0^{+\infty} P^\alpha_t(f) dt.
\end{align}
The potential kernel $V_\alpha$ is defined, for all $x \in \bbr^d \setminus \{0\}$, by 
\begin{align}\label{eq:potential_kernel_stable_measure}
V_\alpha(x) := \int_{0}^{+\infty} p_\alpha\left(\frac{x}{t^{\frac{1}{\alpha}}}\right) \frac{dt}{t^{\frac{d}{\alpha}}}. 
\end{align} 
The regularity of this kernel function has been studied in \cite{BS_SM07,Szt_MN10}. By a standard change of variables, $V_\alpha$ admits the following representation: for all $x \in \mathbb{R}^d$ such that $x \ne 0$, 
\begin{align}\label{eq:potential_kernel_alternative_rep}
V_\alpha(x) = \alpha \int_0^{+\infty} t^{d - (\alpha+1)} p_\alpha\left(t x\right)dt. 
\end{align}
When the L\'evy measure $\nu_\alpha$ is a $\gamma$-measure in the sense of \cite[Definition 2.1]{Szt_MN10} with $\gamma \in [1,d]$ such that $\gamma - d + 2\alpha > 1$, then \cite[Theorem $1.1$]{Szt_MN10} asserts that $V_\alpha$ is differentiable on $\mathbb{R}^d \setminus \{0\}$ and, for all $x \in \mathbb{R}^d \setminus \{0\}$ and all $\beta \in \mathbb{N}_0^d$ such that $ 0\leq |\beta| < \gamma - d + 2\alpha$,  
\begin{align}\label{eq:derivatives_estimates_potential_kernel}
\left| D^{\beta}\left(V_\alpha\right)(x) \right| \leq c(\alpha,d,\gamma,\beta) \|x\|^{\alpha-d-|\beta|}. 
\end{align}
In the next technical lemma, let us specify the link between the potential kernel $V_\alpha$ and the generalized vector-valued kernel $k_\alpha$ defined, for all $\xi \in \mathbb{R}^d \setminus \{0\}$, through
\begin{align}\label{eq:generalized_kernel}
\mathcal{F}\left(k_\alpha\right)(\xi) = \dfrac{(-i\xi)}{\sigma_\alpha\left(\xi\right)^{\alpha}}. 
\end{align}  
\begin{lem}\label{lem:link_between_kernels}
Let $\alpha \in (1,2)$ and let $\nu_\alpha$ be a non-degenerate symmetric $\alpha$-stable L\'evy measure on $\mathbb{R}^d$, $d \geq 2$, such that $\nu_\alpha$ is a $\gamma$-measure with $\gamma \in [1,d]$ and $\gamma - d + 2\alpha > 1$. Then, for all $x \in \mathbb{R}^d \setminus \{0\}$, 
\begin{align}\label{eq:link_between_kernels}
k_\alpha(x) = -\nabla\left(V_\alpha\right)(x). 
\end{align}
In particular, for all $x \in \mathbb{R}^d \setminus \{0\}$, 
\begin{align}\label{ineq:pointwise_estimate_FFTCkernel}
\|k_\alpha(x)\| \leq C_{\alpha,d,\gamma} \|x\|^{\alpha-d-1},
\end{align}
for some constant $C_{\alpha,d,\gamma}>0$ depending on $\alpha$, $d$ and $\gamma$.  
\end{lem}

\begin{proof}
Since $\gamma - d + 2\alpha>1$, the potential kernel $V_\alpha$ is differentiable on $\mathbb{R}^d \setminus \{0\}$.~Moreover, as in the proof of \cite[Theorem $1.1$]{Szt_MN10}, for all $x \in \mathbb{R}^d \setminus \{0\}$, 
\begin{align}\label{eq:gradient_potential_kernel}
\nabla \left(V_\alpha\right)(x) = \alpha \int_0^{+\infty} t^{d-\alpha} \nabla(p_\alpha)\left(tx\right)dt. 
\end{align}
Let us compute in the tempered distributions sense the Fourier transform of the gradient of the potential kernel $V_\alpha$ based on \eqref{eq:gradient_potential_kernel}. Let $\psi \in \mathcal{S}(\mathbb{R}^d)$. Then, for all $j \in \{1, \dots, d\}$, 
\begin{align*}
\langle \mathcal{F}\left( \partial_j\left(V_\alpha\right)\right) ; \psi \rangle & = \int_{\mathbb{R}^d} \partial_j \left(V_\alpha\right)(x) \mathcal{F}(\psi)(x) dx \\
& = \alpha \int_{\mathbb{R}^d} \left(\int_0^{+\infty} t^{d-\alpha} \partial_j(p_\alpha)\left(tx\right)dt\right) \mathcal{F}(\psi)(x)dx. 
\end{align*}
Now, 
\begin{align*}
\int_{\mathbb{R}^d \times (0,+\infty)} t^{d-\alpha} \left|\partial_j(p_\alpha)\left(tx\right) \right| \left| \mathcal{F}(\psi)(x) \right| dxdt \leq C_{\alpha,d,\gamma} \int_{\mathbb{R}^d} \left| \mathcal{F}(\psi)(x) \right| \frac{dx}{\|x\|^{d+1-\alpha}}. 
\end{align*}
The integral on the right-hand side of the previous inequality is finite since $\psi \in \mathcal{S}(\mathbb{R}^d)$ and $\alpha \in (1,2)$. Thus, by Fubini's theorem, 
\begin{align*}
\langle \mathcal{F}\left( \partial_j\left(V_\alpha\right)\right) ; \psi \rangle = \alpha \int_0^{+\infty} t^{d-\alpha} \left( \int_{\mathbb{R}^d}\partial_j\left(p_\alpha\right)(tx) \mathcal{F}(\psi)(x)dx \right) dt.
\end{align*}
Now, by Fubini's theorem again, for all $t>0$, 
\begin{align*}
 \int_{\mathbb{R}^d}\partial_j\left(p_\alpha\right)(tx) \mathcal{F}(\psi)(x)dx & = \int_{\mathbb{R}^d} \mathcal{F} \left(\partial_j\left(p_\alpha\right)(t \cdot)\right)(\xi) \psi(\xi) d\xi \\
 & = \frac{1}{t^d}\int_{\mathbb{R}^d} \mathcal{F} \left(\partial_j\left(p_\alpha\right)\right)\left(\frac{\xi}{t}\right) \psi(\xi) d\xi \\
 & = \frac{1}{t^{d+1}} \int_{\mathbb{R}^d} \exp\left( - \frac{\sigma_\alpha(\xi)^\alpha}{t^\alpha} \right) \psi(\xi) (i \xi_j) d\xi. 
\end{align*}
Then,
\begin{align*}
\langle \mathcal{F}\left( \partial_j\left(V_\alpha\right)\right) ; \psi \rangle & = \alpha \int_0^{+\infty} \frac{1}{t^{\alpha+1}} \left(\int_{\mathbb{R}^d} \exp\left( - \frac{\sigma_\alpha(\xi)^\alpha}{t^\alpha} \right) \psi(\xi) (i \xi_j) d\xi\right) dt \\
& = \int_{\mathbb{R}^d} \dfrac{i \xi_j}{\sigma_\alpha\left(\xi\right)^\alpha} \psi\left(\xi\right) d\xi ,
\end{align*}
where we have applied Fubini's theorem again in the last line since $\nu_\alpha$ is non-degenerate and $\alpha \in (1,2)$. The inequality \eqref{ineq:pointwise_estimate_FFTCkernel} follows from \eqref{eq:derivatives_estimates_potential_kernel} and \eqref{eq:link_between_kernels}, concluding the proof of the lemma. 
\end{proof}
\noindent
The next proposition investigates the fractional fundamental theorem of calculus and the $L^1$-increments for the vector-valued function $k_\alpha$. 

\begin{prop}\label{prop:FFTC_L1_NDSstrong}
Let $\alpha \in (1,2)$ and let $\nu_\alpha$ be a non-degenerate symmetric $\alpha$-stable L\'evy measure on $\mathbb{R}^d$, $d \geq 2$, such that $\nu_\alpha$ is a $\gamma$-measure with $\gamma \in [1,d]$ and $\gamma - d + 2\alpha > 1$. Then, for all $f \in \mathcal{C}^{\infty}_c(\mathbb{R}^d)$ and all $x \in \mathbb{R}^d$, 
\begin{align}\label{eq:FFTC_NDS}
f(x) = \frac{1}{\alpha} \int_{\mathbb{R}^d} \langle k_{\alpha}(x-y) ; D^{\alpha-1}(f)(y) \rangle dy. 
\end{align}
Moreover, if $\gamma - d + 2\alpha > 2$, 
\begin{align}\label{ineq:uniform_bound_L1-increments}
\underset{\omega \in \mathbb{S}^{d-1}}{\sup} \int_{\mathbb{R}^d} \left\| k_\alpha(\omega +z ) - k_\alpha(z) \right\| dz < + \infty.  
\end{align} 
\end{prop}

\begin{proof}
By the Fourier inversion formula and the proof of Lemma \ref{lem:truncated_convolution_NDSkernel}, for all $f \in \mathcal{C}^{\infty}_c(\mathbb{R}^d)$ and all $x \in \mathbb{R}^d$, 
\begin{align*}
f(x) & = \frac{1}{(2\pi)^d} \int_{\mathbb{R}^d} \mathcal{F}(f)(\xi) e^{i \langle x ;\xi  \rangle} d\xi \\
& = \frac{1}{\left(2\pi\right)^d} \int_{\mathbb{R}^d} \mathcal{F}(f)(\xi) \dfrac{\sigma_\alpha(\xi)^\alpha}{\sigma_\alpha(\xi)^\alpha} e^{i \langle x ; \xi \rangle} d\xi \\
& =  \frac{1}{\left(2\pi\right)^d}\int_{\mathbb{R}^d} \frac{1}{\alpha} \bigg\langle \frac{(-i \xi)}{\sigma_\alpha(\xi)^\alpha} ; \tau_\alpha(\xi) \bigg\rangle \mathcal{F}(f)(\xi) e^{i \langle x ; \xi \rangle} d\xi \\
& = \frac{1}{\alpha} \frac{1}{(2\pi)^d} \int_{\mathbb{R}^d} \big\langle \mathcal{F}\left(k_\alpha\right)(\xi) ; \mathcal{F}\left(D^{\alpha-1}(f)\right)(\xi)\big\rangle e^{i \langle x ; \xi \rangle} d\xi \\
& = \frac{1}{\alpha} \int_{\mathbb{R}^d} \langle k_\alpha(x-y) ; D^{\alpha-1}(f)(y) \rangle dy,   
\end{align*}
where in the last line we have used standard Fourier analysis together with Lemma \ref{lem:link_between_kernels}.~Now, \cite[Lemma $2.4$]{BSK_DM20} asserts that, for all $j,k \in \{1, \dots, d\}^2$ and all $x \in \mathbb{R}^d$, 
\begin{align*}
\left| \partial^2_{j,k}\left(p_\alpha\right)(x) \right| \leq \dfrac{C_{\alpha,d,\gamma}}{\left(1+\|x\|\right)^{\alpha+\gamma}},
\end{align*}
for some $C_{\alpha,d,\gamma}>0$ depending only on $\alpha$, $d$ and $\gamma$. This estimate implies the following one on the first order partial derivatives of $k_\alpha$: for all $j,k \in \{1,\dots, d\}^2$ and all $x \in \mathbb{R}^d \setminus \{0\}$, 
\begin{align}\label{ineq:pointwise_inequality_diffkalpha}
\left| \partial_k\left(k^j_\alpha\right)(x) \right| \leq \dfrac{C_{\alpha,d,\gamma}}{\|x\|^{d+2-\alpha}}, 
\end{align}
for some $C_{\alpha,d,\gamma}>0$ depending on $\alpha$, $d$ and $\gamma$. Next, let us perform an analysis similar to the one in the proof of \cite[Proposition $3.14$]{comi_stefani} with the isotropic estimates at our disposal. Let $\omega \in \mathbb{S}^{d-1}$ and let $B(0,2)$ be the Euclidean ball centered at the origin with radius $2$. Then, 
\begin{align*}
 \int_{B(0,2)} \left\| k_\alpha(\omega +z ) - k_\alpha(z) \right\| dz \leq \int_{B(0,2)} \| k_{\alpha}(\omega +z) \| dz + \int_{B(0,2)} \| k_{\alpha}(z) \| dz. 
\end{align*}
Thanks to \eqref{ineq:pointwise_estimate_FFTCkernel}, 
\begin{align*}
\int_{B(0,2)} \| k_{\alpha}(z) \| dz \leq C_{\alpha,d,\gamma} \int_{B(0,2)} \|z\|^{\alpha-d-1} dz <+\infty, 
\end{align*}
since $\alpha \in (1,2)$.~Moreover,
\begin{align*}
 \int_{B(0,2)} \left\| k_\alpha(\omega +z ) \right\| dz \leq C_{\alpha,d,\gamma} \int_{B(0,2)} \|\omega+z\|^{\alpha-d-1} dz \leq  C_{\alpha,d,\gamma} \int_{B(0,3)} \|y\|^{\alpha-d-1} dy<+\infty. 
\end{align*} 
Finally, thanks to \eqref{ineq:pointwise_inequality_diffkalpha}, 
\begin{align*}
 \int_{\mathbb{R}^d \setminus B(0,2)} \left\| k_\alpha(\omega +z ) - k_\alpha(z) \right\| dz & \leq \sum_{\ell=1}^d \int_{\mathbb{R}^d \setminus B(0,2)} \int_0^1 \left| \langle \nabla(k^\ell_\alpha)(z+t \omega) ;  \omega \rangle \right| dt dz \\
 & \leq \sum_{\ell=1}^d \int_0^1 \left(\int_{\mathbb{R}^d \setminus B(0,2)} \left\| \nabla(k^\ell_\alpha)(z+t \omega) \right\| dz\right) dt \\
 & \leq C_{\alpha,d,\gamma} \sum_{\ell=1}^d \int_0^1 \left(\int_{\mathbb{R}^d \setminus B(0,2)} \|z + t\omega\|^{\alpha-d-2} dz\right) dt \\
 & \leq C_{\alpha,d,\gamma} \sum_{\ell=1}^d \int_0^1 \left(\int_{\mathbb{R}^d \setminus B(0,1)} \|y\|^{\alpha-d-2} dy\right) dt <+\infty,
\end{align*}
since $\alpha<2$. All the previous estimates hold uniformly over $\mathbb{S}^{d-1}$. This concludes the proof of the proposition. 
\end{proof}
\noindent
Next, let us study anisotropic version of the fundamental theorem of calculus which puts into play the anisotropic local gradient operator $D_\sigma$ defined in \eqref{eq:anisotropic_local_gradient}.~For this purpose, let us introduce an asymptotic non-degeneracy condition for the spherical component of the L\'evy measure:
\begin{align}\label{eq:AsymNDC}
\inf_{e \in \mathbb{S}^{d-1}} \int_{\mathbb{S}^{d-1}} |\langle e ; y \rangle|^2 \sigma(dy)>0.
\end{align}
Moreover, let us introduce the following notations: for all $x \in \mathbb{R}^d$, $\tau_x$ is the translation operator defined, for all $f$ and all $y \in \mathbb{R}^d$, by $\tau_x(f)(y) = f(x+y)$. Moreover, the operator $\breve{}$ is defined, for all $f$ and all $y \in \mathbb{R}^d$, by $\breve{f}(y) = f(-y)$.   

\begin{prop}\label{prop:AL_FTC}
Let $d \geq 2$ be an integer and let $\sigma$ be a symmetric finite positive measure on $\mathbb{S}^{d-1}$ verifying \eqref{eq:AsymNDC}. Then, for all $f \in \mathcal{C}^{\infty}_c(\mathbb{R}^d)$ and all $x \in \mathbb{R}^d$, 
\begin{align*}
f(x) = \langle \tau_x \circ \breve{k_\sigma} ; D_\sigma(f) \rangle,
\end{align*}
where $k_\sigma$ is the tempered distribution characterized, for all $\xi \in \mathbb{R}^d \setminus \{0\}$, by 
\begin{align}\label{eq:anisotropic_kernel}
\mathcal{F}\left(k_\sigma\right)(\xi) : = \dfrac{-i \xi}{\int_{\mathbb{S}^{d-1}} |\langle\xi ; y \rangle|^2 \sigma(dy)}. 
\end{align}
\end{prop} 

\begin{proof}
First, since $\sigma$ satisfies \eqref{eq:AsymNDC}, for all $\alpha \in (1,2)$, 
\begin{align*}
\inf_{e\in \mathbb{S}^{d-1}} \int_{\mathbb{S}^{d-1}} |\langle e ; y \rangle|^{\alpha} \sigma(dy) \geq \inf_{e\in \mathbb{S}^{d-1}} \int_{\mathbb{S}^{d-1}} |\langle e ; y \rangle|^{2} \sigma(dy)>0.
\end{align*}
Let $f \in \mathcal{C}^{\infty}_c(\mathbb{R}^d)$. By the Fourier inversion formula, for all $x \in \mathbb{R}^d$, 
\begin{align*}
f(x) & = \frac{1}{(2\pi)^d} \int_{\mathbb{R}^d} \mathcal{F}(f)(\xi) e^{i \langle x ; \xi \rangle} d\xi  \\
& = \frac{1}{\left(2\pi\right)^d} \int_{\mathbb{R}^d} \mathcal{F}(f)(\xi) \dfrac{\sigma_\alpha(\xi)^\alpha}{\sigma_\alpha(\xi)^\alpha} e^{i \langle x ; \xi \rangle} d\xi  \\
& = \frac{1}{\left(2\pi\right)^d}\int_{\mathbb{R}^d} \frac{1}{\alpha} \bigg\langle \frac{(-i \xi)}{(2-\alpha)\sigma_\alpha(\xi)^\alpha} ; (2-\alpha)\tau_\alpha(\xi) \bigg\rangle \mathcal{F}(f)(\xi) e^{i \langle x ; \xi \rangle} d\xi.
\end{align*}
Now, for all $\xi \in \mathbb{R}^d \setminus \{0\}$, 
\begin{align*}
(2-\alpha) \tau_\alpha(\xi) & = i \frac{\Gamma(3-\alpha) \cos\left(\frac{\alpha \pi}{2}\right)}{1- \alpha} \int_{\mathbb{S}^{d-1}} y |\langle y ; \xi \rangle|^{\alpha-1} \operatorname{sign}\left( \langle y ; \xi \rangle\right) \sigma(dy) \\
& \longrightarrow i \int_{\mathbb{S}^{d-1}} y \langle y ; \xi \rangle \sigma(dy), 
\end{align*}
as $\alpha$ tends to $2^-$. Moreover, for all $\xi \in \mathbb{R}^d \setminus \{0\}$, 
\begin{align*}
(2- \alpha) \sigma_\alpha(\xi)^\alpha & = (2-\alpha) \int_{\mathbb{S}^{d-1}} |\langle \xi ; y  \rangle|^\alpha \lambda_1(dy) \\
& =  \dfrac{\Gamma(3-\alpha) |\cos \left(\frac{\alpha\pi}{2}\right)|}{\alpha(\alpha-1)}\int_{\mathbb{S}^{d-1}} |\langle \xi ; y  \rangle|^\alpha \sigma(dy) \\
& \longrightarrow \frac{1}{2} \int_{\mathbb{S}^{d-1}} |\langle \xi ; y  \rangle|^2 \sigma(dy),  
\end{align*}
as $\alpha$ tends to $2^-$. Thus, for all $\xi \in \mathbb{R}^d \setminus \{0\}$, 
\begin{align*}
\underset{\alpha \rightarrow 2^-}{\lim} \frac{1}{\alpha} \bigg\langle \frac{(-i \xi)}{(2-\alpha)\sigma_\alpha(\xi)^\alpha} ; (2-\alpha)\tau_\alpha(\xi) \bigg\rangle = \bigg\langle \dfrac{(-i \xi)}{ \int_{\mathbb{S}^{d-1}} |\langle \xi ; y  \rangle|^2 \sigma(dy)} ; i \int_{\mathbb{S}^{d-1}} y \langle y ; \xi \rangle \sigma(dy) \bigg\rangle. 
\end{align*}
Moreover, for all $\alpha \in (1,2)$ and all $\xi \in \mathbb{R}^d \setminus \{0\}$, 
\begin{align*}
 \bigg| \frac{1}{\alpha} \bigg\langle \frac{(-i \xi)}{(2-\alpha)\sigma_\alpha(\xi)^\alpha} ; (2-\alpha)\tau_\alpha(\xi) \bigg\rangle \bigg| \leq 1.
\end{align*}
By the Lebesgue dominated convergence theorem, 
\begin{align*}
f(x) & = \frac{1}{(2\pi)^d} \int_{\mathbb{R}^d} \mathcal{F}(f)(\xi) \bigg\langle \dfrac{(-i \xi)}{ \int_{\mathbb{S}^{d-1}} |\langle \xi ; y  \rangle|^2 \sigma(dy)} ; i \int_{\mathbb{S}^{d-1}} y \langle y ; \xi \rangle \sigma(dy) \bigg\rangle e^{i \langle x ; \xi \rangle} d\xi. 
\end{align*}
The end of the proof easily follows by standard Fourier analysis. 
\end{proof}
\noindent
Next, let us study the tempered distribution $k_\sigma$ defined by \eqref{eq:anisotropic_kernel}. For this purpose, let $p_{2,\sigma}$ be defined, for all $\xi \in \mathbb{R}^d$, by
\begin{align}\label{eq:anisotropic_Gaussian}
\mathcal{F}\left(p_{2,\sigma}\right)\left(\xi\right) = \exp \left( - \int_{\mathbb{S}^{d-1}} |\langle \xi ; y \rangle|^2 \sigma(dy)\right). 
\end{align} 
But, for all $\xi \in \mathbb{R}^d$, 
\begin{align*}
\int_{\mathbb{S}^{d-1}} |\langle \xi ; y \rangle|^2 \sigma(dy) = \sum_{\ell,k=1}^d \xi_\ell \xi_k \sigma_{\ell,k}, \quad \sigma_{\ell,k} = \int_{\mathbb{S}^{d-1}} y_\ell y_k \sigma(dy). 
\end{align*}
Let $\mathbf{\Sigma}$ be the $d \times d$ matrix whose coefficient at the position $(\ell,k)$ is given by $\sigma_{\ell,k}$. Then, for all $\xi \in \mathbb{R}^d$, 
\begin{align}\label{eq:from_measure_to_matrix}
\int_{\mathbb{S}^{d-1}} |\langle \xi ; y \rangle|^2 \sigma(dy) = \big\langle \xi ; \mathbf{\Sigma}\xi \big\rangle. 
\end{align}
$\mathbf{\Sigma}$ is a real symmetric positive definite matrix since $\sigma$ satisfies \eqref{eq:AsymNDC}. Thus, it is clear that $p_{2, \sigma} \in \mathcal{S}(\mathbb{R}^d)$, $p_{2,\sigma}(x) \geq 0$, for $x \in \mathbb{R}^d$, and $\int_{\mathbb{R}^d} p_{2,\sigma}(x)dx = 1$. 

\begin{lem}\label{lem:from_iso_to_aniso}
Let $d \geq 2$ be an integer, let $\sigma$ be a symmetric finite positive measure on $\mathbb{S}^{d-1}$ satisfying \eqref{eq:AsymNDC} and let $\mathbf{\Sigma}$ be the symmetric positive definite matrix defined in \eqref{eq:from_measure_to_matrix}. Let $k_\sigma$ be defined by \eqref{eq:anisotropic_kernel}. Then, for all $x \in \mathbb{R}^d \setminus \{0\}$,
\begin{align}\label{eq:spatial_rep_anisokernel}
k_{\sigma}(x) = \dfrac{1}{\sqrt{\operatorname{det}\left(\mathbf{\Sigma}\right)}} \left(\dfrac{2}{\Gamma(\frac{d}{2})}\pi^{\frac{d}{2}}\right)^{-1}  \dfrac{\mathbf{\Sigma}^{-1}x}{\| \mathbf{\Sigma}^{-\frac{1}{2}}x\|^{d}}.
\end{align}
\end{lem}

\begin{proof}
The proof of \eqref{eq:spatial_rep_anisokernel} follows from well-chosen changes of variables.~Indeed, for all $\varphi \in \mathcal{S}(\mathbb{R}^d, \mathbb{R}^d)$, 
\begin{align*}
\langle k_\sigma ; \mathcal{F}(\varphi) \rangle & = \int_{\mathbb{R}^d} \left\langle \dfrac{-i\xi}{\| \mathbf{\Sigma}^{\frac{1}{2}} \xi \|^2} ; \varphi(\xi) \right\rangle d\xi  \\
& =  \int_{\mathbb{R}^d}  \left\langle \dfrac{-i \mathbf{\Sigma}^{-\frac{1}{2}} \xi}{\| \xi \|^2} ;  \varphi\left( \mathbf{\Sigma}^{-\frac{1}{2}} \xi\right) \right\rangle \frac{d\xi}{\sqrt{\operatorname{det}\left(\mathbf{\Sigma}\right)}}  \\
& =   \int_{\mathbb{R}^d}  \left\langle \dfrac{-i\xi}{\| \xi \|^2} ;  \mathbf{\Sigma}^{-\frac{1}{2}} \varphi\left( \mathbf{\Sigma}^{-\frac{1}{2}} \xi\right) \right\rangle \frac{d\xi}{\sqrt{\operatorname{det}\left(\mathbf{\Sigma}\right)}}  \\
& = \frac{1}{\sqrt{\operatorname{det}\left(\mathbf{\Sigma}\right)}}  \langle  \mathbf{\Sigma}^{-\frac{1}{2}} \tilde{k_1} ;   \mathcal{F}\left(\varphi \circ \mathbf{\Sigma}^{-\frac{1}{2}}\right)\rangle  \\
& = \frac{1}{\sqrt{\operatorname{det}\left(\mathbf{\Sigma}\right)}}  \langle  \mathbf{\Sigma}^{-\frac{1}{2}} \tilde{k_1}\circ  \mathbf{\Sigma}^{-\frac{1}{2}} ;   \mathcal{F}\left(\varphi\right)\rangle ,
\end{align*} 
where $\tilde{k_1}$ is given, for all $x\in \mathbb{R}^d\setminus \{0\}$, 
\begin{align*}
\tilde{k_1}(x) = \left(\dfrac{2}{\Gamma(\frac{d}{2})}\pi^{\frac{d}{2}}\right)^{-1}  \dfrac{x}{\|x\|^{d}}. 
\end{align*}
This concludes the proof of the lemma. 
\end{proof}

\section{(An)isotropic fractional Sobolev spaces}\label{sec:anisotropic_fractional_Sobolev_spaces}
\noindent
Let us recall the definition of the fractional gradient operator $D^{\alpha-1}$: for all $f \in \mathcal{S}(\mathbb{R}^d)$ and all $x \in \mathbb{R}^d$, 
\begin{align*}
D^{\alpha-1} (f)(x) = \int_{\mathbb{R}^d} (f(x+u) - f(x)) u \nu_\alpha(du) = \frac{1}{(2\pi)^d} \int_{\mathbb{R}^d} \mathcal{F}(f)(\xi) e^{i \langle x; \xi \rangle } \tau_{\alpha}(\xi) d\xi, 
\end{align*}
with, for all $\xi \in \mathbb{R}^d$,
\begin{align*}
\tau_\alpha(\xi) = \int_{\mathbb{R}^d} \left(e^{i \langle u; \xi \rangle}-1\right) u \nu_\alpha(du). 
\end{align*}
In general, the multiplier $\xi \mapsto \tau_\alpha(\xi)$ is not smooth so that the operator $D^{\alpha-1}$ does not map $\mathcal{S}(\mathbb{R}^d)$ into $\mathcal{S}(\mathbb{R}^d , \mathbb{R}^d)$ and does not extend to $\mathcal{S}'(\mathbb{R}^d)$ a priori. To circumvent this fact, inspired by \cite[Chapter $3$, Section $3.3$, pages $270$ - $275$]{NJ02_2}, let us introduce a truncation procedure: observe that, for all $R>0$, all $f \in \mathcal{S}(\mathbb{R}^d)$ and all $x \in \mathbb{R}^d$, 
\begin{align*}
D^{\alpha-1}(f)(x) = \mathcal{F}^{-1} \left( \tau_{\alpha, R} \mathcal{F}(f)\right)(x) + \int_{B(0,R)^c} \left(f(x+u) - f(x)\right) u \nu_\alpha(du), 
\end{align*}
where $B(0,R) = \{x \in \mathbb{R}^d : \, \|x\| < R\}$ and where $\tau_{\alpha,R}$ is defined, for all $\xi \in \mathbb{R}^d$, by
\begin{align*}
\tau_{\alpha,R}\left( \xi \right) = \int_{B(0,R)}  \left(e^{i \langle \xi  ; u \rangle}-1\right) u\nu_\alpha(du). 
\end{align*}
Clearly, $\xi \mapsto \tau_{\alpha, R}(\xi)$ is infinitely differentiable on $\mathbb{R}^d$. Moreover, for all $\beta \in \mathbb{N}_0^d$ with $|\beta| := \beta_1 + \cdots + \beta_d \geq 1$, 
\begin{align*}
D^{\beta} \left(\tau_{\alpha, R}\right)(\xi) = \int_{B(0,R)} i^{|\beta|} \prod_{k = 1}^d u_k^{\beta_k} e^{i \langle u ; \xi \rangle } u \nu_\alpha(du).
\end{align*}
Then, for all $\beta \in \mathbb{N}_0^d$ with $|\beta| \geq 1$ and all $\xi \in \mathbb{R}^d$, 
\begin{align*}
\left\| D^{\beta} \left(\tau_{\alpha, R}\right)(\xi)  \right\| \leq \int_{B(0,R)} \prod_{k = 1}^d |u_k|^{\beta_k} \|u\| \nu_\alpha(du). 
\end{align*}
Thus, $\tau_{\alpha,R} \in \mathcal{C}^{\infty}_p(\mathbb{R}^d , \mathbb{R}^d)$ which implies that the operator $D^{\alpha-1}_R$, with $R>0$, defined, for all $f \in \mathcal{S}(\mathbb{R}^d)$, by
\begin{align}\label{eq:frac_gradient_truncated}
D^{\alpha-1}_{R}(f) = \mathcal{F}^{-1} \left(\tau_{\alpha,R}\mathcal{F}(f)\right),
\end{align}
maps $\mathcal{S}(\mathbb{R}^d)$ into $\mathcal{S}(\mathbb{R}^d, \mathbb{R}^d)$. Then, $D^{\alpha-1}_R$ extends by duality to $S'(\mathbb{R}^d)$. In particular, it is well-defined for $f \in L^p(\mathbb{R}^d, dx)$, with $p \in [1,+\infty)$, as an $\mathbb{R}^d$-valued tempered distribution.~Moreover, by Minkowski's integral inequality, for all $f \in \mathcal{S}(\mathbb{R}^d)$, 
\begin{align*}
\left\| \int_{B(0,R)^c} \left(f(\cdot+u) - f(\cdot)\right) u \nu_\alpha(du) \right\|_{L^p(\mathbb{R}^d, dx)} \leq 2 \|f\|_{L^p(\mathbb{R}^d,dx)} \int_{B(0,R)^c} \|u\| \nu_\alpha(du), 
\end{align*}
which is finite since $\alpha \in (1,2)$. This upper bound allows, by a classical density argument, to extend the previous linear integral operator to any $f \in L^p(\mathbb{R}^d,dx)$. Then, for all $R>0$ and all $f \in L^p(\mathbb{R}^d, dx)$, $D^{\alpha-1}(f)$ is a well-defined element of $S'(\mathbb{R}^d, \mathbb{R}^d)$. Next, let us prove that this definition of $D^{\alpha-1}(f)$ does not depend on the truncation parameter $R>0$: let $R_1 > R_2 > 0$ and let us denote for now $D^{\alpha-1, R_1} (f)$ and $D^{\alpha-1,R_2}(f)$, for $f \in L^p(\mathbb{R}^d,dx)$, the fractional gradients associated with the truncation parameters $R_1$ and $R_2$ respectively. Now, for all $k \in \{1, \dots, d\}$ and all $\varphi \in \mathcal{S}(\mathbb{R}^d)$, 
\begin{align*}
\langle \mathcal{F}^{-1} \left( \tau_{\alpha,R_1,k} \mathcal{F}(f) \right) ; \varphi \rangle = \langle f ; (D^{\alpha-1}_{R_1,k})^*(\varphi) \rangle, \quad (D^{\alpha-1}_{R_1,k})^*(\varphi)(x) = \int_{B(0,R_1)} \Delta_{-u}(\varphi)(x) u_k \nu_\alpha(du). 
\end{align*} 
Then, 
\begin{align*}
\langle  \mathcal{F}^{-1} \left( \tau_{\alpha,R_1,k} \mathcal{F}(f) \right) ; \varphi \rangle -\langle  \mathcal{F}^{-1} \left( \tau_{\alpha,R_2,k} \mathcal{F}(f) \right) ; \varphi \rangle & =  \langle f ; (D^{\alpha-1}_{R_1,k})^*(\varphi) - (D^{\alpha-1}_{R_2,k})^*(\varphi) \rangle  \\
& =  \left\langle f ; \int_{B(0,R_1) \setminus B(0,R_2)} \Delta_{-u}(\varphi) u_k \nu_\alpha(du) \right\rangle  \\
& =  \left\langle \int_{B(0,R_1) \setminus B(0,R_2)} \Delta_{u}(f) u_k \nu_\alpha(du) ; \varphi \right\rangle.
\end{align*}
Thus, in $S'(\mathbb{R}^d , \mathbb{R}^d)$, 
\begin{align*}
\mathcal{F}^{-1} \left( \tau_{\alpha,R_1} \mathcal{F}(f) \right)  - \mathcal{F}^{-1} \left( \tau_{\alpha,R_2} \mathcal{F}(f) \right)  =  \int_{B(0,R_1) \setminus B(0,R_2)} \Delta_{u}(f) u\nu_\alpha(du). 
\end{align*}
From this equality, it easily follows that $D^{\alpha-1 , R_1}(f) = D^{\alpha-1,R_2}(f)$, for all $R_1>R_2>0$ and all $f \in L^p(\mathbb{R}^d,dx)$. Since, for all $p \in [1, +\infty)$ and all $f \in L^p(\mathbb{R}^d,dx)$, $D^{\alpha-1}(f)$ is a well-defined $\mathbb{R}^d$-valued tempered distribution, it is natural to identify its action on Schwartz functions. This is the purpose of the next lemma. 

\begin{lem}\label{lem:weak_fractional_derivative}
Let $d \geq 1$ be an integer, let $\alpha \in (1,2)$ and let $p \in [1, +\infty)$.~Let $f \in L^p(\mathbb{R}^d,dx) $. Then, for all $\varphi \in \mathcal{S}(\mathbb{R}^d)$ and all $k \in \{1, \dots, d\}$, 
\begin{align}\label{eq:weak_fractional_derivative}
\langle D^{\alpha-1}_k(f) ; \varphi \rangle = \langle f ; (D^{\alpha-1}_k)^*(\varphi) \rangle,
\end{align}
where, for all $x \in \mathbb{R}^d$, 
\begin{align*}
(D^{\alpha-1}_k)^*(\varphi)(x) = \int_{\mathbb{R}^d} \left(\varphi(x-u) -\varphi(x)\right) u_k \nu_{\alpha}(du). 
\end{align*}
\end{lem}

\begin{proof}
Let $R>0$. Then, for all $\varphi \in \mathcal{S}(\mathbb{R}^d)$ and all $k \in \{1, \dots, d\}$, 
\begin{align}\label{eq:decomposition_dualitybracket}
\langle D^{\alpha-1}_k(f) ; \varphi \rangle = \langle D^{\alpha-1}_{R,k}(f) ; \varphi \rangle + \left\langle \int_{B(0,R)^c} \left(f(\cdot + u) - f(\cdot) \right) u_k \nu_\alpha(du) ; \varphi \right\rangle.
\end{align}
Let $(f_n)_{n \geq 1}$ be a sequence of Schwartz functions which converges toward $f$ in $L^p(\mathbb{R}^d,dx)$. Then, in $L^p(\mathbb{R}^d,dx)$, 
\begin{align*}
 \int_{B(0,R)^c} \left(f_n(\cdot + u) - f_n(\cdot) \right) u_k \nu_\alpha(du) \longrightarrow  \int_{B(0,R)^c} \left(f(\cdot + u) - f(\cdot) \right) u_k \nu_\alpha(du), \quad n \rightarrow +\infty.
\end{align*}
Moreover, by Fubini's theorem, for all $n \geq 1$, 
\begin{align*}
\left\langle \int_{B(0,R)^c} \left(f_n(\cdot + u) - f_n(\cdot) \right) u_k \nu_\alpha(du) ; \varphi \right\rangle = \left\langle f_n ; \int_{B(0,R)^c} \left( \varphi(\cdot -u) - \varphi(\cdot)\right) u_k \nu_\alpha(du) \right\rangle.
\end{align*}
Thus, 
\begin{align*}
\left\langle \int_{B(0,R)^c} \left(f(\cdot + u) - f(\cdot) \right) u_k \nu_\alpha(du) ; \varphi \right\rangle = \left\langle f ; \int_{B(0,R)^c} \left( \varphi(\cdot -u) - \varphi(\cdot)\right) u_k \nu_\alpha(du) \right\rangle.
\end{align*}
For the first term of the decomposition \eqref{eq:decomposition_dualitybracket}, 
\begin{align*}
\langle D^{\alpha-1}_{R,k}(f) ; \varphi \rangle & = \langle \mathcal{F}^{-1} \left(\tau_{\alpha,R,k} \mathcal{F}(f)\right) ; \varphi \rangle   \\
& =\frac{1}{(2\pi)^d} \langle \tau_{\alpha,R,k} \mathcal{F}(f) ; \mathcal{F}(\varphi) \rangle  \\
& = \frac{1}{(2\pi)^d} \langle \mathcal{F}(f) ; \overline{\tau_{\alpha,R,k}}\mathcal{F}(\varphi) \rangle \\
& = \langle f ; \mathcal{F}^{-1} \left( \overline{\tau_{\alpha,R,k}}\mathcal{F}(\varphi) \right) \rangle, 
\end{align*}
where $\tau_{\alpha,R,k}(\xi) = \langle \tau_{\alpha,R}(\xi) ; e_k \rangle$, for all $\xi \in \mathbb{R}^d$. Now, by the Fourier inversion formula, for all $x \in \mathbb{R}^d$,  
\begin{align*}
\mathcal{F}^{-1} \left( \overline{\tau_{\alpha,R,k}}\mathcal{F}(\varphi)\right) (x) & = \frac{1}{(2\pi)^d} \int_{\mathbb{R}^d} \mathcal{F}(\varphi)(\xi)  \overline{\tau_{\alpha,R,k}}(\xi) e^{i \langle x ; \xi \rangle} d\xi  \\
& = \frac{1}{(2\pi)^d} \int_{\mathbb{R}^d} \mathcal{F}(\varphi)(\xi) \left(\int_{B(0,R)} (e^{- i\langle u;\xi \rangle}- 1) u_k \nu_\alpha(du)\right) e^{i \langle x ; \xi \rangle} d\xi  \\
& = \int_{B(0,R)} \left( \varphi(x-u) - \varphi(x) \right) u_k \nu_\alpha(du). 
\end{align*}
The conclusion of the proof easily follows. 
\end{proof}

\begin{rem}\label{rem:comparison_weak_alpha-gradient}
A weak fractional gradient operator is introduced in \cite[Definition $3.19$]{comi_stefani} and defined by duality when the underlying L\'evy measure is invariant by rotation.  Based on Lemma \ref{lem:weak_fractional_derivative}, it is clear that this definition of weak fractional gradient operator coincides with $D^{\alpha-1,\operatorname{rot}}(f)$, for $f \in L^p(\mathbb{R}^d,dx)$, with $p \in [1,+\infty)$.   
\end{rem}

\begin{defi}\label{def:non-homogeneous_frac_Sobolev_spaces}
Let $d \geq 1$ be an integer, let $\alpha \in (1,2)$ and let $p \in [1,+\infty)$. The non-homogeneous fractional Sobolev space of order $\alpha-1$ is defined by 
\begin{align}\label{eq:def_fractional_Sobolev-type_space_order1}
W^{\alpha-1,p}(\mathbb{R}^d,dx) : = \{f \in L^p(\mathbb{R}^d,dx):\, \forall k \in \{1, \dots, d\},\, D^{\alpha-1}_k(f) \in L^p(\mathbb{R}^d,dx) \}. 
\end{align}
This set is endowed with the positively one-homogeneous functional $\|\cdot\|_{\alpha-1,p}$ defined, for all $f \in W^{\alpha-1,p}(\mathbb{R}^d,dx)$, by 
\begin{align}\label{eq:def_fractional_Sobolev-type_norm_order1}
\|f\|_{\alpha-1,p} : = \left(\|f\|^p_{L^p(\mathbb{R}^d,dx)} + \sum_{k = 1}^d \left\|D_k^{\alpha-1}(f)\right\|^p_{L^p(\mathbb{R}^d,dx)}\right)^{\frac{1}{p}},
\end{align}
and $\left(W^{\alpha-1,p}(\mathbb{R}^d,dx) , \|\cdot\|_{\alpha-1,p} \right)$ becomes a normed vector space. 
\end{defi}
\noindent
With the help of a classical argument combined with Lemma \ref{lem:weak_fractional_derivative}, $\left(W^{\alpha-1,p}(\mathbb{R}^d,dx) , \|\cdot\|_{\alpha-1,p}\right)$ is a Banach space.

\begin{rem}\label{rem:fractional_Sobolev-type_space_order1}
Let $W_0^{\alpha-1,p}(\mathbb{R}^d,dx)$ be the set of functions defined by
\begin{align}\label{eq:def_fractional_Sobolev-type_space_order1_completion}
W_0^{\alpha-1,p}(\mathbb{R}^d,dx) : = \overline{\mathcal{C}_c^{\infty}\left(\mathbb{R}^d\right)}^{\| \cdot \|_{\alpha-1,p}}. 
\end{align}
By definition, $W_0^{\alpha-1,p}\left(\mathbb{R}^d,dx\right) \subset W^{\alpha-1,p}\left(\mathbb{R}^d,dx\right)$.~In the rotationally invariant case, \cite[Theorem $3.23$]{comi_stefani} for $p = 1$ and \cite[Appendix A, Theorem $25$]{BCGS_CRM22} for $p \in (1, +\infty)$ prove that both sets are actually equal. For all $p \in (1,+\infty)$, recall that $\mathcal{A}_{\alpha,p}$ denotes the $L^p(\mathbb{R}^d,dx)$-generator of $(P^{\alpha}_t)_{t\geq 0}$ and, for all $s>0$, $H_p^{\psi_\alpha,s}(\mathbb{R}^d)$ denotes the Bessel potential-type space defined by \eqref{eq:def_fractional_BPT_space_orders}. These spaces have been introduced and deeply analyzed in \cite[Chapter $3$, Sections $3.2$ and $3.3$]{NJ02_2}. In particular, \cite[Proposition $3.3.14$]{NJ02_2} ensures that, for all $s>0$,
\begin{align*}
\overline{\mathcal{S}\left(\mathbb{R}^d\right)}^{\|\cdot\|_{\psi_\alpha,s,p}} = H_p^{\psi_\alpha,s}(\mathbb{R}^d).
\end{align*}
In the sequel, we prove that the following equalities (with equivalence of corresponding norms) hold true: for all $p \in (1, +\infty)$, 
\begin{align}\label{eq:Calderon_result_order1}
W^{\alpha-1,p}(\mathbb{R}^d,dx) = W_0^{\alpha-1,p}(\mathbb{R}^d,dx) = H_p^{\psi_\alpha,r(\alpha)}(\mathbb{R}^d), \quad r(\alpha) = \dfrac{2(\alpha-1)}{\alpha}. 
\end{align}
Again, in the rotationally invariant case, \cite[Theorem $1.7.$]{Shieh_Spector_ACV15} provides
\begin{align*}
W_{0, \operatorname{rot}}^{\alpha-1,p}(\mathbb{R}^d,dx) = L^{\alpha-1,p}(\mathbb{R}^d), 
\end{align*}
where $L^{\alpha-1,p}(\mathbb{R}^d)$ is the classical fractional Bessel potential space of order $\alpha-1$. 
\end{rem}

\begin{prop}\label{prop:consq_Riesz_transform}
Let $\alpha \in (1,2)$, let $p \in (1,+\infty)$, let $r(\alpha) = 2(\alpha-1)/\alpha$ and let $d \geq 1$ be an integer. Then, $H^{\psi_\alpha, r(\alpha)}_p(\mathbb{R}^d) \subset W^{\alpha-1,p}(\mathbb{R}^d,dx)$ and, for all $f \in H_p^{\psi_\alpha,r(\alpha)}(\mathbb{R}^d)$ and all $k \in \{1, \dots, d\}$, 
\begin{align}\label{eq:continuous_embedding1}
\| D^{\alpha-1}_k(f) \|_{L^p(\mathbb{R}^d,dx)} \leq C_{\alpha,p,d} \|f\|_{\psi_\alpha, r(\alpha), p},  
\end{align}
for some $C_{\alpha,p,d}>0$ depending on $\alpha$, $p$ and $d$. 
\end{prop}

\begin{proof}
Let $f \in \mathcal{S}(\mathbb{R}^d)$. Then, from \cite[Theorem $3.2$]{AH20_4}, for all $k \in \{1, \dots, d\}$, 
\begin{align*}
\left\| D_{k}^{\alpha-1}(f)  \right\|_{L^p(\mathbb{R}^d,dx)} & = \left\| D_{k}^{\alpha-1}\circ \left(E - \mathcal{A}_{\alpha,p}\right)^{-\frac{\alpha-1}{\alpha}} \circ \left(E - \mathcal{A}_{\alpha,p}\right)^{\frac{\alpha-1}{\alpha}}(f) \right\|_{L^p(\mathbb{R}^d,dx)} \\
& \leq C_{\alpha,p,d} \|  \left(E - \mathcal{A}_{\alpha,p}\right)^{\frac{\alpha-1}{\alpha}}(f) \|_{L^p(\mathbb{R}^d,dx)}\\
& \leq C_{\alpha,p,d} \|f\|_{\psi_\alpha,r(\alpha),p}. 
\end{align*} 
Now, using \cite[Proposition $3.3.14$]{NJ02_2} combined with a standard argument, for all $f \in H^{\psi_\alpha, r(\alpha)}_p(\mathbb{R}^d)$ there exists a sequence $(f_n)_{n \geq 1}$ of Schwartz functions such that $(D^{\alpha-1}_k(f_n))_{n \geq 1}$ converges to some $g_k$ in $L^p(\mathbb{R}^d,dx)$, for all $k \in \{1,\dots, d\}$. Moreover, thanks to Lemma \ref{lem:weak_fractional_derivative}, for all $k \in \{1, \dots, d\}$ and all $\varphi \in \mathcal{S}(\mathbb{R}^d)$, 
\begin{align*}
\langle g_k ; \varphi \rangle =\underset{n \rightarrow +\infty}{\lim} \langle D^{\alpha-1}_k(f_n) ; \varphi \rangle = \underset{n \rightarrow +\infty}{\lim} \langle f_n ; (D^{\alpha-1}_k)^*(\varphi) \rangle = \langle  D^{\alpha-1}_k(f) ; \varphi  \rangle. 
\end{align*}
Then, $f \in W^{\alpha-1,p}(\mathbb{R}^d,dx)$ and inequality \eqref{eq:continuous_embedding1} is obtained by passing to the limit. 
\end{proof}
\noindent
In the next lemma, let us provide a sufficient condition which ensures that the $\mathbb{R}^d$-valued tempered distribution $D^{\alpha-1}(f)$, with $f \in L^p(\mathbb{R}^d,dx)$, belongs to $L^p(\mathbb{R}^d,dx)$. For this purpose, recall that $W^{1,p}(\mathbb{R}^d,dx)$ denotes the classical first order Sobolev space with integrability exponent $p$. 

\begin{lem}\label{lem:interpolation_inequality}
Let $\alpha \in (1,2)$, let $p \in [1,+\infty)$ and let $d \geq 1$ be an integer.~Then, for all $f \in W^{1,p}(\mathbb{R}^d,dx)$, $D^{\alpha-1}_k(f) \in L^p(\mathbb{R}^d,dx)$, for all $k \in \{1, \dots, d\}$. Moreover, for all $k \in \{1, \dots, d\}$, 
\begin{align}\label{ineq:interpolation_inequality1}
\left\| D^{\alpha-1}_k(f) \right\|_{L^p(\mathbb{R}^d,dx)} \leq c_{\alpha,d} \|f\|^{2-\alpha}_{L^p(\mathbb{R}^d,dx)} \|\nabla(f)\|^{\alpha-1}_{L^p(\mathbb{R}^d,dx)},
\end{align}
for some $c_{\alpha,d}>0$ depending on $\alpha$, $d$ and the spherical part of the L\'evy measure $\nu_\alpha$. 
\end{lem}

\begin{proof}
Let $f \in \mathcal{S}(\mathbb{R}^d)$ and let $R>0$ to be chosen later on.~Then, for all $x \in \mathbb{R}^d$ and all $k \in \{1, \dots, d\}$, 
\begin{align*}
\left| D^{\alpha-1}_k(f)(x) \right| \leq \int_{B(0,R)} |f(x+u)-f(x)| |u_k| \nu_\alpha(du) + \int_{B(0,R)^c} |f(x+u)-f(x)| |u_k| \nu_\alpha(du).
\end{align*}
Then, by the triangle inequality, 
\begin{align*}
\left\| D^{\alpha-1}_k(f) \right\|_{L^p(\mathbb{R}^d,dx)} \leq & \left\| \int_{B(0,R)} |f(\cdot+u)-f(\cdot)| |u_k| \nu_\alpha(du)\right\|_{L^p(\mathbb{R}^d,dx)} \\
& \quad \quad + \left\| \int_{B(0,R)^c} |f(\cdot+u)-f(\cdot)| |u_k| \nu_\alpha(du) \right\|_{L^p(\mathbb{R}^d,dx)}.
\end{align*}
For the second term on the right-hand side of the previous inequality, by Minkowski's integral inequality, 
\begin{align*}
\left\| \int_{B(0,R)^c} |f(\cdot+u)-f(\cdot)| |u_k| \nu_\alpha(du) \right\|_{L^p(\mathbb{R}^d,dx)} & \leq 2 \|f\|_{L^p(\mathbb{R}^d,dx)} \int_{B(0,R)^c} \|u\| \nu_\alpha(du) \\
& \leq 2 \sigma\left(\mathbb{S}^{d-1}\right) \int_{r \geq R} \frac{dr}{r^{\alpha}} \|f\|_{L^p(\mathbb{R}^d,dx)} \\
& \leq 2 \sigma\left(\mathbb{S}^{d-1}\right) \frac{ \|f\|_{L^p(\mathbb{R}^d,dx)} }{(\alpha-1)R^{\alpha-1}} .
\end{align*} 
For the other term, using Taylor's formula since $f \in \mathcal{S}(\mathbb{R}^d)$, for all $x,u \in \mathbb{R}^d$, 
\begin{align*}
f(x+u) - f(x) = \int_0^1 \langle \nabla(f)(x+tu) ; u \rangle dt. 
\end{align*}
Then, by Minkowski's integral inequality again and translation invariance, 
\begin{align*}
\left\| \int_{B(0,R)} |f(\cdot+u)-f(\cdot)| |u_k| \nu_\alpha(du)\right\|_{L^p(\mathbb{R}^d,dx)} & \leq \| \nabla(f) \|_{L^p(\mathbb{R}^d,dx)} \int_{B(0,R)} \|u\|^2 \nu_\alpha(du) \\
& \leq \sigma\left(\mathbb{S}^{d-1}\right) \int_0^R r^{1- \alpha} dr \| \nabla(f) \|_{L^p(\mathbb{R}^d,dx)} \\
&\leq  \sigma\left(\mathbb{S}^{d-1}\right) \frac{R^{2-\alpha}}{2-\alpha} \| \nabla(f) \|_{L^p(\mathbb{R}^d,dx)}. 
\end{align*}
Now, setting $R =  \|f\|_{L^p(\mathbb{R}^d,dx)} / \| \nabla(f) \|_{L^p(\mathbb{R}^d,dx)}$, 
\begin{align*}
\left\| D^{\alpha-1}_k(f) \right\|_{L^p(\mathbb{R}^d,dx)}  \leq c_{\alpha,d} \|f\|^{2-\alpha}_{L^p(\mathbb{R}^d,dx)} \|\nabla(f)\|^{\alpha-1}_{L^p(\mathbb{R}^d,dx)},
\end{align*}
for some $c_{\alpha,d}>0$ depending on $\alpha$, $\sigma$ and $d$. Next, let $f \in W^{1,p}(\mathbb{R}^d,dx)$ and let $(f_n)_{n \geq 1}$ be a sequence of Schwartz functions converging to $f$ in $W^{1,p}(\mathbb{R}^d,dx)$. Then, for all $n,m \geq 1$ and all $k \in \{1, \dots, d\}$, 
\begin{align*}
\left\| D^{\alpha-1}_k(f_m-f_n) \right\|_{L^p(\mathbb{R}^d,dx)}  \leq c_{\alpha,d} \|f_m-f_n\|^{2-\alpha}_{L^p(\mathbb{R}^d,dx)} \|\nabla(f_m-f_n)\|^{\alpha-1}_{L^p(\mathbb{R}^d,dx)}.
\end{align*}
Since $(f_n)_{n \geq 1}$ is convergent in $W^{1,p}(\mathbb{R}^d,dx)$, it is fundamental in $W^{1,p}(\mathbb{R}^d,dx)$.~The previous inequality implies that $\left(D^{\alpha-1}_k(f_n)\right)_{n \geq 1}$ is fundamental in $L^p(\mathbb{R}^d,dx)$ so that there exists $g_k \in L^p(\mathbb{R}^d,dx)$ such that $D^{\alpha-1}_k(f_n) \rightarrow g_k$ in $L^p(\mathbb{R}^d,dx)$, as $n \rightarrow +\infty$. To conclude the proof of the lemma, it remains to prove that $D^{\alpha-1}_k(f) = g_k$: for all $\varphi \in \mathcal{S}(\mathbb{R}^d)$ and all $k \in \{1, \dots, d\}$,
\begin{align*}
\langle g_k ; \varphi \rangle = \underset{n \rightarrow +\infty}{\lim} \langle D^{\alpha-1}_k(f_n) ; \varphi \rangle =  \underset{n \rightarrow +\infty}{\lim} \langle f_n ; (D^{\alpha-1}_k)^*(\varphi) \rangle = \langle f ; (D^{\alpha-1}_k)^*(\varphi) \rangle = \langle D^{\alpha-1}_k(f) ; \varphi \rangle.
\end{align*}
This concludes the proof of the lemma.  
\end{proof}

\begin{prop}\label{prop:density_compactsupport}
Let $d \geq 1$ be an integer, let $\alpha \in (1,2)$ and let $p \in (1, +\infty)$. Then, 
\begin{align}\label{eq:density_result1}
W^{\alpha-1,p}(\mathbb{R}^d,dx) = \overline{\mathcal{C}^{\infty}_c(\mathbb{R}^d)}^{\|\cdot\|_{\alpha-1,p}}. 
\end{align}
\end{prop}

\begin{proof}
Let $f \in W^{\alpha-1,p}(\mathbb{R}^d,dx)$. Let $\rho$ be a symmetric non-negative function in $\mathcal{C}_c^{\infty} \left( \mathbb{R}^d\right)$ such that $\int_{\mathbb{R}^d} \rho(x)dx=1$. Let $(\rho_{\varepsilon})_{\varepsilon>0}$ be defined, for all $\varepsilon>0$ and all $x \in \mathbb{R}^d$, by 
\begin{align*}
\rho_{\varepsilon}(x) = \frac{1}{\varepsilon^d} \rho \left(\frac{x}{\varepsilon}\right). 
\end{align*} 
Now, let $(f_{\varepsilon})_{\varepsilon>0}$ be defined, for all $\varepsilon>0$, by $f_{\varepsilon} = f \ast \rho_{\varepsilon}$ where $\ast$ stands for the convolution product. It is clear that $f_\varepsilon \in \mathcal{C}^{\infty}(\mathbb{R}^d)$, for all $\varepsilon >0$, and that $f_\varepsilon  \longrightarrow f$ in $L^p(\mathbb{R}^d,dx)$ as $\varepsilon$ tends to $0^+$. Moreover, for all $\varepsilon>0$, all $x \in \mathbb{R}^d$ and all $k \in \{1, \dots, d\}$, 
\begin{align*}
D^{\alpha-1}_k(f_\varepsilon)(x) = \int_{\mathbb{R}^d} \left(f_{\varepsilon}(x+u) - f_{\varepsilon}(x) \right) u_k \nu_\alpha(du).
\end{align*}
Next, let $\varphi \in \mathcal{S}(\mathbb{R}^d)$. Then, for all $\varepsilon>0$, 
\begin{align*}
\langle D^{\alpha-1}_k(f_\varepsilon) ; \varphi \rangle & = \langle f_\varepsilon ; (D^{\alpha-1}_k)^*(\varphi) \rangle = \int_{\mathbb{R}^d} \int_{\mathbb{R}^d} f(y) \rho_{\varepsilon} (x-y) (D^{\alpha-1}_k)^*(\varphi)(x) dx dy \\
& =  \int_{\mathbb{R}^d} \int_{\mathbb{R}^d} f(y) \rho_{\varepsilon} (y-x) (D^{\alpha-1}_k)^*(\varphi)(x) dx dy \\
& =\int_{\mathbb{R}^d} f(y) (D^{\alpha-1}_k)^*(\varphi) \ast \rho_{\varepsilon}(y)dy.
\end{align*}
Moreover, for all $y \in \mathbb{R}^d$
\begin{align*}
(D^{\alpha-1}_k)^*(\varphi) \ast \rho_{\varepsilon}(y) & = \int_{\mathbb{R}^d} (D^{\alpha-1}_k)^*(\varphi)(z) \rho_{\varepsilon}(y-z) dz \\
& = \int_{\mathbb{R}^d} \left( \int_{\mathbb{R}^d} (\varphi(z-u) - \varphi(z)) u_k \nu_\alpha(du) \right)  \rho_{\varepsilon}(y-z) dz \\
& = \int_{\mathbb{R}^d} (\varphi \ast \rho_{\varepsilon}(y-u) - \varphi \ast \rho_{\varepsilon}(y)) u_k\nu_\alpha(du) = (D^{\alpha-1}_k)^*(\varphi \ast \rho_{\varepsilon})(y).  
\end{align*}
Thus, for all $\varepsilon>0$, 
\begin{align*}
\langle D^{\alpha-1}_k(f_\varepsilon) ; \varphi \rangle & = \langle f ; (D^{\alpha-1}_k)^*(\varphi \ast \rho_{\varepsilon}) \rangle = \langle D^{\alpha-1}_k(f) ; \varphi \ast \rho_\varepsilon \rangle \\
& =  \langle D^{\alpha-1}_k(f) \ast \rho_{\varepsilon} ; \varphi \rangle. 
\end{align*}
By standard mollifier theory, $D^{\alpha-1}_k(f) \ast \rho_{\varepsilon} \rightarrow  D^{\alpha- 1}_k(f)$ in $L^p(\mathbb{R}^d,dx)$, as $\varepsilon$ tends to $0^+$. Now, since $\rho_\varepsilon \in \mathcal{C}^{\infty}_c(\mathbb{R}^d)$, it follows that $\partial_j(\rho_\varepsilon) \in L^1(\mathbb{R}^d,dx)$, for all $j \in \{1, \dots, d\}$ and all $\varepsilon>0$. Moreover,
\begin{align*}
\nabla \left( f \ast \rho_\varepsilon \right) = f \ast \nabla \left(\rho_\varepsilon\right).
\end{align*}
Finally, by Young's inequality, for all $\varepsilon>0$ and all $j \in \{1, \cdots, d\}$,
\begin{align*}
\| \partial_j(f_\varepsilon)\|_{L^p(\mathbb{R}^d,dx)} \leq \| \partial_j \left(\rho_\varepsilon\right) \|_{L^1(\mathbb{R}^d,dx)} \|f\|_{L^p(\mathbb{R}^d,dx)}.
\end{align*}
Thus, $\mathcal{C}^{\infty}(\mathbb{R}^d) \cap W^{1,p}(\mathbb{R}^d,dx)$ is dense in $W^{\alpha-1,p}\left(\mathbb{R}^d,dx\right)$. Without loss of generality, let us assume that $f \in W^{\alpha-1, p }(\mathbb{R}^d,dx) \cap W^{1,p} \left(\mathbb{R}^d,dx\right)$. Now, recall that $\mathcal{C}_c^{\infty}(\mathbb{R}^d)$ is dense in $W^{1,p}(\mathbb{R}^d,dx)$ with respect to the classical Sobolev norm $\|\cdot\|_{1,p}$.~Then, let $(f_n)_{n \geq 1}$ be a sequence of functions in $\mathcal{C}^{\infty}_c(\mathbb{R}^d)$ which converges to $f$ in $W^{1,p}(\mathbb{R}^d,dx)$. From the interpolation inequality \eqref{ineq:interpolation_inequality1}, for all $n \geq 1$, 
\begin{align*}
\left\| D^{\alpha-1}_k(f-f_n) \right\|_{L^p(\mathbb{R}^d,dx)} \leq c_{\alpha,d} \|f-f_n\|^{2-\alpha}_{L^p(\mathbb{R}^d,dx)} \|\nabla(f-f_n)\|^{\alpha-1}_{L^p(\mathbb{R}^d,dx)} \longrightarrow 0,
\end{align*}
as $n$ tends to $+\infty$. Then, $(f_n)_{n \geq 1}$ converges to $f$ in the norm $\|\cdot\|_{\alpha-1,p}$.
\end{proof}
\noindent

\begin{prop}\label{prop:converse_inclusion}
Let $d \geq 1$ be an integer, let $\alpha \in (1,2)$, let $p \in (1, +\infty)$ and let $r(\alpha) = 2(\alpha-1)/\alpha$. Then, 
\begin{align}
W^{\alpha-1,p}(\mathbb{R}^d,dx) = H^{\psi_\alpha, r(\alpha)}_p(\mathbb{R}^d). 
\end{align}
\end{prop}

\begin{proof}
Proposition \ref{prop:consq_Riesz_transform} ensures that $H^{\psi_\alpha, r(\alpha)}_p\left(\mathbb{R}^d\right) \subset W^{\alpha-1,p}\left(\mathbb{R}^d,dx\right)$ and the embedding is continuous. Now, thanks to \cite[Proposition 3.9]{AH20_4}, for all $f \in \mathcal{S}(\mathbb{R}^d)$, 
\begin{align*}
\| \left( E - \mathcal{A}_\alpha \right)^{\frac{\alpha-1}{\alpha}}(f)\|_{L^p(\mathbb{R}^d,dx)} \leq C_{\alpha,p,d} \left( \|f\|_{L^p(\mathbb{R}^d,dx)} + \| D^{\alpha-1}(f)\|_{L^p(\mathbb{R}^d,dx)}\right),
\end{align*}
for some $C_{\alpha,p,d}>0$ depending on $\alpha$, $p$ and $d$. Next, let $f \in W^{\alpha-1,p}(\mathbb{R}^d,dx)$. Then, Proposition \ref{prop:density_compactsupport} asserts that there exists a sequence $(f_n)_{n \geq 1}$ of functions belonging to $\mathcal{C}^{\infty}_c(\mathbb{R}^d)$ which tends to $f$ in $W^{\alpha-1,p}(\mathbb{R}^d,dx)$. Then, for all $n,m \geq 1$, 
\begin{align*}
\| \left( E - \mathcal{A}_\alpha \right)^{\frac{\alpha-1}{\alpha}}(f_n - f_m)\|_{L^p(\mathbb{R}^d,dx)} \leq C_{\alpha,p,d} \left( \|f_n- f_m\|_{L^p(\mathbb{R}^d,dx)} + \| D^{\alpha-1}(f_n -f_m)\|_{L^p(\mathbb{R}^d,dx)}\right).
\end{align*} 
Thus, $( \left( E - \mathcal{A}_\alpha \right)^{\frac{\alpha-1}{\alpha}}(f_n))_{n \geq 1}$ is fundamental in $L^{p}\left(\mathbb{R}^d,dx\right)$ which is complete by the Riesz-Fischer theorem. Then, there exists $g \in L^p(\mathbb{R}^d,dx)$ such that 
\begin{align*}
\left( E - \mathcal{A}_\alpha \right)^{\frac{\alpha-1}{\alpha}}(f_n) \longrightarrow g, \quad n \rightarrow +\infty,
\end{align*}
in $L^p(\mathbb{R}^d,dx)$. The closedness of $\left( E - \mathcal{A}_{\alpha,p} \right)^{\frac{\alpha-1}{\alpha}}$ concludes the proof. 
\end{proof}
\noindent
Based on Definition \ref{def:non-homogeneous_frac_Sobolev_spaces}, it is natural to introduce the following fractional Sobolev spaces of order $\alpha$: for all $p \in [1, +\infty)$,  let $W^{\alpha,p}(\mathbb{R}^d,dx)$ be the set of functions $f \in L^p(\mathbb{R}^d, dx)$ such that $\partial_j(f)$, $D^{\alpha-1}_k(f)$ and $\partial_j \circ D^{\alpha-1}_{k}(f)$ belong to $L^p(\mathbb{R}^d,dx)$, for all $j,k \in \{1, \dots, d\}$. On this set of functions, let us introduce a positively one-homogeneous functional $\| \cdot \|_{\alpha,p}$ defined by 
\begin{align}\label{eq:norm_fractional_sobolev_space}
\|f\|^p_{\alpha,p} & : = \|f\|^p_{L^p(\mathbb{R}^d,dx)} + \sum_{j = 1}^d \|\partial_j(f)\|^p_{L^p(\mathbb{R}^d,dx)} + \sum_{k = 1}^d \| D^{\alpha-1}_k(f)\|^p_{L^p(\mathbb{R}^d,dx)} \nonumber \\
& \quad\quad + \sum_{j ,k =1}^d \| \partial_j D_k^{\alpha-1} (f)\|^p_{L^p(\mathbb{R}^d,dx)}. 
\end{align}
$\|\cdot\|_{\alpha,p}$ is a norm on $W^{\alpha,p}(\mathbb{R}^d,dx)$ and $\left(W^{\alpha,p}(\mathbb{R}^d,dx) , \| \cdot \|_{\alpha,p}\right)$ is a Banach space.~In the sequel, let us investigate the following equalities (with equivalence of corresponding norms): for all $p \in (1, +\infty)$, 
\begin{align}\label{eq:Calderon_result_order2}
W^{\alpha,p}(\mathbb{R}^d,dx) = \overline{\mathcal{C}^{\infty}_c(\mathbb{R}^d)}^{\|\cdot\|_{\alpha,p}} = H_p^{\psi_\alpha,2}(\mathbb{R}^d) = D(\mathcal{A}_{\alpha,p}) = \{f \in L^p(\mathbb{R}^d,dx) : \, \mathcal{A}_{\alpha}(f) \in L^p(\mathbb{R}^d,dx)\}.
\end{align}
\noindent
As shown next, the known equality $ H_p^{\psi_\alpha,2}(\mathbb{R}^d) = D\left(\mathcal{A}_{\alpha,p}\right)$ (see, e.g., \cite[Chapter $3$, Section $3.3$, Theorem $3.3.11$]{NJ02_2}) will imply that the inclusion $W^{\alpha,p}(\mathbb{R}^d,dx) \subset H_p^{\psi_\alpha,2}(\mathbb{R}^d)$ becomes clear once it is proved that $\mathcal{S}(\mathbb{R}^d)$ is dense in $W^{\alpha,p}(\mathbb{R}^d,dx)$ and that, for all $f \in W^{\alpha,p}(\mathbb{R}^d,dx)$ and all $t>0$, 
\begin{align*}
\mathcal{A}_\alpha (P^\alpha_t(f)) = P_t^{\alpha}\left(\mathcal{A}_\alpha(f)\right). 
\end{align*} 
Indeed, let us assume for now that the two previous properties are true.~Then, let $(f_n)_{n \geq 1}$ be a sequence of Schwartz functions such that $f_n \longrightarrow f$ in $W^{\alpha,p}(\mathbb{R}^d,dx)$, as $n \rightarrow +\infty$. Then, for all $t>0$ and all $n \geq 1$, 
\begin{align*}
\left\| \dfrac{P^\alpha_t(f)-f}{t} - \mathcal{A}_\alpha(f) \right\|_{L^p(\mathbb{R}^d,dx)} &\leq \left\| \dfrac{P^\alpha_t(f)-f}{t} -\dfrac{P^\alpha_t(f_n)-f_n}{t} \right\|_{L^p(\mathbb{R}^d,dx)} +\left\| \mathcal{A}_\alpha(f_n) - \mathcal{A}_\alpha(f) \right\|_{L^p(\mathbb{R}^d,dx)} \\
& \quad\quad + \left\| \dfrac{P^\alpha_t(f_n)-f_n}{t} - \mathcal{A}_\alpha(f_n) \right\|_{L^p(\mathbb{R}^d,dx)}.
\end{align*}
Since $(f_n)_{n \geq 1}$ tends to $f$ in $W^{\alpha,p}(\mathbb{R}^d,dx)$, it is clear that $\mathcal{A}_\alpha(f_n) \rightarrow \mathcal{A}_\alpha(f)$ in $L^p(\mathbb{R}^d,dx)$, as $n \rightarrow +\infty$. Moreover, regarding the first term, for all $t>0$ and all $n \geq 1$, 
\begin{align*}
\left\| \dfrac{P^\alpha_t(f)-f}{t} -\dfrac{P^\alpha_t(f_n)-f_n}{t} \right\|_{L^p(\mathbb{R}^d,dx)} & = \left\| \frac{1}{t} \int_0^t \left[\mathcal{A}_\alpha(P_s^\alpha(f_n)) - \mathcal{A}_\alpha(P_s^\alpha(f)) \right] ds\right\|_{L^p(\mathbb{R}^d,dx)}  \\
& = \left\| \frac{1}{t} \int_0^t \left[P_s^\alpha(\mathcal{A}_\alpha(f_n)) - P_s^\alpha(\mathcal{A}_\alpha(f)) \right] ds\right\|_{L^p(\mathbb{R}^d,dx)}  \\
& \leq \left\| \mathcal{A}_\alpha(f_n) - \mathcal{A}_\alpha(f) \right\|_{L^p(\mathbb{R}^d,dx)},
\end{align*}
since $(P_t^\alpha)_{t \geq 0}$ is a semigroup of contractions on $L^p(\mathbb{R}^d,dx)$, for every $p \in [1, +\infty]$.~Furthermore, the previous upper bound is uniform in $t$. Finally, since $\mathcal{S}(\mathbb{R}^d) \subset D \left( \mathcal{A}_{\alpha,p}\right)$ and $\mathcal{A}_{\alpha,p} = \mathcal{A}_\alpha$ on $\mathcal{S}(\mathbb{R}^d)$, for all $n \geq 1$, 
\begin{align*}
\left\| \dfrac{P^\alpha_t(f_n)-f_n}{t} - \mathcal{A}_\alpha(f_n) \right\|_{L^p(\mathbb{R}^d,dx)} \underset{t \rightarrow 0^+}{\longrightarrow} 0. 
\end{align*}
Let us prove the validity of the two properties mentioned above. 

\begin{lem}\label{lem:interpolation_inequality2_density_result2}
Let $d \geq 1$ be an integer, let $\alpha \in (1,2)$ and let $p \in (1,+\infty)$. Then, $W^{2,p}\left(\mathbb{R}^d,dx\right) \subset W^{\alpha,p}\left(\mathbb{R}^d,dx\right)$ and, for all $f \in W^{2,p}\left(\mathbb{R}^d,dx\right)$, all $j \in \{1, \dots, d\}$ and all $k \in \{1, \dots, d\}$,
\begin{align}\label{ineq:interpolation_inequality2}
\| \partial_j D_k^{\alpha-1}(f)\|_{L^p(\mathbb{R}^d,dx)} \leq c_{\alpha,d} \| \partial_j(f) \|^{2-\alpha}_{L^p(\mathbb{R}^d,dx)} \| \nabla(\partial_j(f))\|^{\alpha-1}_{L^p(\mathbb{R}^d,dx)},
\end{align}
for some $c_{\alpha,d}>0$ depending on $\alpha$, $d$ and the spherical part of the L\'evy measure $\nu_\alpha$. Moreover,
\begin{align}\label{eq:density_result2}
W^{\alpha,p}(\mathbb{R}^d,dx) = \overline{\mathcal{C}^{\infty}_c(\mathbb{R}^d)}^{\|\cdot\|_{\alpha,p}}. 
\end{align}
\end{lem}

\begin{proof}
Let $f \in \mathcal{S}(\mathbb{R}^d)$. Then, for all $j,k \in \{1, \cdots, d\}$ and all $x \in \mathbb{R}^d$, 
\begin{align}\label{eq:commutation_derivatives}
\partial_j D^{\alpha-1}_k(f)(x) = D_k^{\alpha-1}\left(\partial_j(f)\right)(x).
\end{align}
Moreover, since $\partial_j(f) \in \mathcal{S}(\mathbb{R}^d)$, one can apply inequality \eqref{ineq:interpolation_inequality1} to $\partial_j(f)$.~So, the inequality \eqref{ineq:interpolation_inequality2} is proved on $\mathcal{S}(\mathbb{R}^d)$. Now, let $f \in W^{2,p}(\mathbb{R}^d,dx)$ and let $(f_n)_{n \geq 1}$ be a sequence of Schwartz functions such that $f_n \rightarrow f$ in $W^{2,p}(\mathbb{R}^d,dx)$, as $n$ tends to $+\infty$. Then, for all $n,m \geq 1$, 
\begin{align*}
\| \partial_j D_k^{\alpha-1}(f_n - f_m)\|_{L^p(\mathbb{R}^d,dx)} \leq c_{\alpha,d} \| \partial_j(f_n - f_m) \|^{2-\alpha}_{L^p(\mathbb{R}^d,dx)} \| \nabla(\partial_j(f_n - f_m))\|^{\alpha-1}_{L^p(\mathbb{R}^d,dx)}. 
\end{align*}
This implies that $(\partial_j D_k^{\alpha-1}(f_n))_{n \geq 1}$ converges to some $g_{j,k}$ in $L^p(\mathbb{R}^d,dx)$, for all $j,k \in \{1, \dots, d\}$. Then, for all $\varphi \in \mathcal{S}(\mathbb{R}^d)$ and all $j,k \in \{1, \dots, d\}$, 
\begin{align*}
\langle g_{j,k} ; \varphi \rangle = \underset{n \rightarrow +\infty}{\lim} \langle  \partial_j D_k^{\alpha-1}(f_n) ; \varphi \rangle = - \langle f ; (D^{\alpha-1}_k)^* \partial_j (\varphi)  \rangle = \langle \partial_j D_k^{\alpha-1}(f) ; \varphi \rangle, 
\end{align*}
which proves that $f \in W^{\alpha,p}\left( \mathbb{R}^d,dx \right)$. Next, let us prove equality \eqref{eq:density_result2}.~As in the proof of Proposition \ref{prop:density_compactsupport}, let $f \in W^{\alpha,p}(\mathbb{R}^d,dx)$, let $(\rho_{\varepsilon})_{\varepsilon>0}$ be a sequence of standard mollifiers and let $(f_{\varepsilon})_{\varepsilon>0}$ be defined, for all $\varepsilon>0$, by $f_{\varepsilon} = f \ast \rho_{\varepsilon}$. It is clear that $f_\varepsilon \in \mathcal{C}^{\infty}(\mathbb{R}^d)$, for all $\varepsilon >0$, and that $f_\varepsilon  \longrightarrow f$ in $L^p(\mathbb{R}^d,dx)$ as $\varepsilon$ tends to $0^+$. Moreover, for all $j \in \{1, \dots, d\}$, 
\begin{align*}
\partial_j (f \ast \rho_{\varepsilon}) = \partial_j(f) \ast \rho_{\varepsilon} \longrightarrow \partial_j(f),
\end{align*} 
in $L^{p}(\mathbb{R}^d,dx)$, as $\varepsilon$ tends to $0^+$. Also, as in the proof of Proposition \ref{prop:density_compactsupport}, $D^{\alpha-1}_k(f_\varepsilon) \rightarrow  D^{\alpha- 1}_k(f)$ in $L^p(\mathbb{R}^d,dx)$, as $\varepsilon$ tends to $0^+$. Finally, for all $j,k \in \{1, \dots, d\}$ and all $\varepsilon>0$, 
\begin{align*}
\partial_j D^{\alpha-1}_k(f_\varepsilon) = \partial_j D^{\alpha-1}_k(f) \ast \rho_{\varepsilon} \longrightarrow \partial_j D^{\alpha-1}_k(f), 
\end{align*} 
in $L^p(\mathbb{R}^d,dx)$. Then, based on the previous analysis, $f_{\varepsilon}$ tends to $f$ in $W^{\alpha,p}(\mathbb{R}^d,dx)$ as $\varepsilon \rightarrow 0^+$. Next, let us prove that $f_\varepsilon \in W^{2,p}\left(\mathbb{R}^d,dx\right)$, for all $\varepsilon>0$. Clearly, for all $\varepsilon>0$, $f_\varepsilon \in W^{1,p}(\mathbb{R}^d,dx)$. Moreover, for all $j,k \in \{1, \dots, d\}$ and all $\varepsilon>0$, 
\begin{align*}
\partial^2_{j,k} (f_{\varepsilon}) = \partial_j(f) \ast \partial_k(\rho_\varepsilon) \Rightarrow \| \partial^2_{j,k} (f_{\varepsilon}) \|_{L^p(\mathbb{R}^d,dx)} \leq \| \partial_k(\rho_\varepsilon) \|_{L^1(\mathbb{R}^d,dx)} \| \partial_j(f) \|_{L^p(\mathbb{R}^d,dx)} <+\infty.
\end{align*}
Thus, without loss of generality, let us assume that $f \in W^{\alpha,p}(\mathbb{R}^d,dx) \cap W^{2,p}\left(\mathbb{R}^d,dx\right) \cap \mathcal{C}^{\infty}(\mathbb{R}^d)$. Then, there exists a sequence of functions $(f_n)_{n \geq 1}$ belonging to $\mathcal{C}^{\infty}_c(\mathbb{R}^d)$ such that $f_n \rightarrow f$ in $W^{2,p}(\mathbb{R}^d,dx)$, as $n$ tends to $+\infty$. By the interpolation inequalities \eqref{ineq:interpolation_inequality1} and \eqref{ineq:interpolation_inequality2}, for all $n \geq 1$ and all $j,k \in \{1, \cdots, d\}$, 
\begin{align*}
&\| \partial_j D_k^{\alpha-1}(f-f_n)\|_{L^p(\mathbb{R}^d,dx)} \leq c_{\alpha,d} \| \partial_j(f-f_n) \|^{2-\alpha}_{L^p(\mathbb{R}^d,dx)} \| \nabla(\partial_j(f-f_n))\|^{\alpha-1}_{L^p(\mathbb{R}^d,dx)} , \\
&\| D_k^{\alpha-1}(f-f_n)\|_{L^p(\mathbb{R}^d,dx)} \leq c_{\alpha,d} \| f-f_n \|^{2-\alpha}_{L^p(\mathbb{R}^d,dx)} \| \nabla(f-f_n)\|^{\alpha-1}_{L^p(\mathbb{R}^d,dx)},
\end{align*}
which imply that $(f_n)_{n \geq 1}$ converges to $f$ in $W^{\alpha,p}\left(\mathbb{R}^d,dx\right)$. This concludes the proof of the lemma. 
\end{proof}
\noindent
For the second property, we have:

\begin{lem}\label{lem:commutation_result}
Let $\alpha \in (1,2)$ and let $p \in (1,+\infty)$.~Then, for all $f \in W^{\alpha,p}(\mathbb{R}^d,dx)$, $d \geq 1$, and all $t \geq 0$, $P_t^\alpha(f) \in W^{\alpha,p}(\mathbb{R}^d,dx)$ and
\begin{align*}
\mathcal{A}_\alpha (P^\alpha_t(f)) = P^{\alpha}_t(\mathcal{A}_\alpha(f)),
\end{align*}
where the equality holds in $L^p(\mathbb{R}^d,dx)$.
\end{lem}

\begin{proof}
Let $f \in W^{\alpha,p}(\mathbb{R}^d,dx)$ and let $(f_n)_{n \geq 1}$ be a sequence of Schwartz functions which converges to $f$ in $W^{\alpha,p}(\mathbb{R}^d,dx)$ (which exists thanks to Lemma \ref{lem:interpolation_inequality2_density_result2}). Recall that, for all $t>0$ and all $x \in \mathbb{R}^d$, 
\begin{align*}
P^\alpha_t(f)(x) = \int_{\mathbb{R}^d} f(y) p_\alpha \left(\frac{x-y}{t^{\frac{1}{\alpha}}}\right) \frac{dy}{t^{\frac{d}{\alpha}}},
\end{align*}
where $p_\alpha$ is the Lebesgue density of the non-degenerate and symmetric $\alpha$-stable probability measure $\mu_\alpha$. Then, for all $\varphi \in \mathcal{S}(\mathbb{R}^d)$, all $t \geq 0$ and all $k \in \{1, \dots, d\}$,  
\begin{align*}
\langle D^{\alpha-1}_k(P_t^{\alpha}(f)) ; \varphi \rangle & = \langle P^\alpha_t(f) ; \left(D^{\alpha-1}_k\right)^*(\varphi) \rangle \\
& = \langle f ; P^{\alpha}_t \circ \left(D^{\alpha-1}_k\right)^*(\varphi) \rangle  \\
& = \langle f ; (D^{\alpha-1}_k)^* \circ P^\alpha_t\left(\varphi\right) \rangle  \\
&= \underset{n\rightarrow+\infty}{\lim} \langle f_n ;   (D^{\alpha-1}_k)^* \circ P^\alpha_t\left(\varphi\right) \rangle  \\
&= \underset{n\rightarrow+\infty}{\lim} \langle D^{\alpha-1}_k(f_n) ; P^\alpha_t\left(\varphi\right) \rangle  \\
& = \langle D_k^{\alpha-1}(f) ; P^\alpha_t(\varphi) \rangle \\
& = \langle P_t^\alpha \left(D_k^{\alpha-1}(f)\right) ; \varphi \rangle.
\end{align*}
The last equality ensures that $D^{\alpha-1}_k(P_t^{\alpha}(f)) \in L^p(\mathbb{R}^d,dx)$ and that $D^{\alpha-1}_k(P_t^{\alpha}(f)) = P_t^\alpha \left(D_k^{\alpha-1}(f) \right)$, for all $k \in \{1, \dots, d\}$. Similarly, since $\partial_j \circ P_t^\alpha = P^\alpha_t \circ \partial_j$, $\partial_j P_t^\alpha(f)$ belong to $L^p(\mathbb{R}^d,dx)$, for all $j \in \{1, \dots, d\}$. To conclude, it remains to prove that $\partial_j D^{\alpha-1}_k \left(P_t^\alpha(f)\right) = P_t^\alpha \left( \partial_j D^{\alpha-1}_k(f)\right)$, for all $j,k \in \{1, \dots, d\}$. But, for all $j,k \in \{1, \cdots, d\}$, all $\varphi \in \mathcal{S}(\mathbb{R}^d)$ and all $t \geq 0$,  
\begin{align*}
\langle \partial_j D^{\alpha-1}_k(P_t^\alpha(f)) ; \varphi \rangle & = - \langle f ; P^\alpha_t \circ \partial_j \circ (D^{\alpha-1}_k)^{*} \left(\varphi\right)   \rangle \\
& = - \langle f ; (D^{\alpha-1}_k)^* \circ \partial_j P^\alpha_t\left(\varphi\right) \rangle \\
& = - \underset{n \rightarrow +\infty}{\lim} \langle f_n ; (D^{\alpha-1}_k)^* \circ \partial_j P^\alpha_t\left(\varphi\right) \rangle \\
& = \underset{n \rightarrow +\infty}{\lim} \langle \partial_j D^{\alpha-1}_k(f_n) ; P^{\alpha}_t(\varphi) \rangle \\
& = \langle \partial_j D^{\alpha-1}_k(f) ; P^{\alpha}_t(\varphi) \rangle = \langle P^{\alpha}_t (\partial_j D^{\alpha-1}_k(f)) ; \varphi \rangle. 
\end{align*}
This concludes the proof of the lemma.
\end{proof}
\noindent
These two lemmas and the arguments preceding them prove that $W^{\alpha,p}(\mathbb{R}^d,dx) \subset H_p^{\psi_\alpha,2}(\mathbb{R}^d)$. Next, it remains to prove the converse inclusion which is based on the $L^p(\mathbb{R}^d,dx)$-continuity of the $\alpha$-order fractional Riesz transform analyzed in \cite{AH20_4}.  

\begin{thm}\label{thm:identification_domain_generator}
Let $d \geq 1$ be an integer, let $\alpha \in (1,2)$ and let $p \in (1,+\infty)$.~Then, 
\begin{align}\label{eq:identification_domain_generator}
W^{\alpha,p}(\mathbb{R}^d,dx) = D \left(\mathcal{A}_{\alpha,p}\right). 
\end{align}
\end{thm}

\begin{proof}
Based on previous observations, it is clear that $W^{\alpha,p}\left(\mathbb{R}^d,dx\right) \subset D\left(\mathcal{A}_{\alpha,p}\right)$. Next, let $f \in D\left(\mathcal{A}_{\alpha,p}\right) = H_p^{\psi_\alpha,2}(\mathbb{R}^d)$. Then, since $H_p^{\psi_\alpha,2}(\mathbb{R}^d) \subset H^{\psi_\alpha, r(\alpha)}_p(\mathbb{R}^d)$ (with continuous embedding) and since $H^{\psi_\alpha, r(\alpha)}_p(\mathbb{R}^d) = W^{\alpha-1,p}(\mathbb{R}^d,dx)$ by Proposition \ref{prop:converse_inclusion}, $D_{k}^{\alpha-1}(f)$ belongs to $L^p(\mathbb{R}^d,dx)$, for all $k \in \{1, \cdots, d\}$. Moreover, \cite[Inequality (3.40)]{AH20_4}, for all $j \in \{1, \dots, d\}$,
\begin{align*}
\| \partial_j(f)\|_{L^p(\mathbb{R}^d,dx)} & = \| \partial_j \circ \left(E- \mathcal{A}_{\alpha,p}\right)^{- \frac{1}{\alpha}} \circ \left(E- \mathcal{A}_{\alpha,p}\right)^{\frac{1}{\alpha}} (f)\|_{L^p(\mathbb{R}^d,dx)} \\
& \leq C_{\alpha,p,d} \left( \|f\|_{L^p(\mathbb{R}^d,dx)} + \|(-\mathcal{A}_{\alpha,p})^{\frac{1}{\alpha}}(f)\|_{L^p(\mathbb{R}^d,dx)} \right)<+\infty,
\end{align*}
since $f \in D\left(\mathcal{A}_{\alpha,p}\right)$. Now, let $(f_n)_{n \geq 1}$ be a sequence of Schwartz functions such that $f_n \longrightarrow f$ in $D\left(\mathcal{A}_{\alpha,p}\right)$ (such a sequence exists \cite[Theorem $3.3.11.$]{NJ02_2}). Then, by \cite[Remark $3.10$]{AH20_4}, for all $n,m \geq 1$ and all $j,k \in \{1, \cdots, d\}$, 
\begin{align}\label{ineq:calderon_zygmund_ineq}
\| \partial_j D_k^{\alpha-1}\left(f_n-f_m\right)\|_{L^p(\mathbb{R}^d,dx)} \leq C_{\alpha,p,d} \left(  \|f_n - f_m\|_{L^p(\mathbb{R}^d,dx)} + \| \mathcal{A}_{\alpha,p}(f_n - f_m)\|_{L^p(\mathbb{R}^d,dx)} \right), 
\end{align}
for some constant $C_{\alpha,p,d}>0$ depending on $\alpha$, $p$ and $d$. For all $j,k \in \{1, \cdots, d\}$, $(\partial_j D_k^{\alpha-1}(f_n))_{n \geq 1}$ converges to some $g_{j,k}$ in $L^p(\mathbb{R}^d,dx)$. Finally, for all $\varphi \in \mathcal{S}(\mathbb{R}^d)$ and all $j,k \in \{1, \cdots, d\}$, 
\begin{align*}
\langle g_{j,k} ; \varphi \rangle & = \underset{n \rightarrow +\infty}{\lim} \langle \partial_j D^{\alpha-1}_k(f_n)  ; \varphi \rangle = - \underset{n \rightarrow +\infty}{\lim} \langle f_n ; (D^{\alpha-1}_k)^* \circ \partial_j \left(\varphi\right) \rangle \\
& = - \langle f ; (D^{\alpha-1}_k)^* \circ \partial_j \left(\varphi\right) \rangle =  \langle \partial_j D^{\alpha-1}_k(f) ; \varphi \rangle.
\end{align*}
This concludes the proof of the theorem. 
\end{proof}
\noindent 
In the sequel, let us analyze the behavior of the (an)isotropic fractional gradient operator $D^{\alpha-1}$ when $\alpha$ tends to $2$ after an appropriate normalization procedure. For this purpose, let us start with a pointwise result and a $L^2$-result with test function belonging to the Schwartz space $\mathcal{S}(\mathbb{R}^d)$. Lemmas \ref{lem:pointwise_L2_BBM_result} and \ref{lem:uniform_inequality_H1} below are similar to the results obtained in \cite[Section $5$]{BS_JMAA18}.   

\begin{lem}\label{lem:pointwise_L2_BBM_result}
Let $\alpha \in (1,2)$, let $\sigma$ be the spherical part of the L\'evy measure $\nu_\alpha$ and let $D_\sigma$ be the anisotropic local gradient operator given by \eqref{eq:anisotropic_local_gradient}.~Then, for all $f \in \mathcal{S}(\mathbb{R}^d)$ and all $x \in \mathbb{R}^d$, $d\geq 1$,
\begin{align}\label{eq:pointwise_limit_BBManisotropic}
\left(2-\alpha\right) D^{\alpha-1}(f)(x) \underset{\alpha \rightarrow 2^-}{\longrightarrow} D_\sigma(f)(x).
\end{align}
Moreover, for all $f \in \mathcal{S}(\mathbb{R}^d)$, 
\begin{align}\label{eq:L2_limit_BBManisotropic}
\left\| \left(2-\alpha\right) D^{\alpha-1}(f) \right\|_{L^2(\mathbb{R}^d,dx)} \underset{\alpha \rightarrow 2^-}{\longrightarrow} \|D_\sigma(f) \|_{L^2(\mathbb{R}^d,dx)}. 
\end{align}
\end{lem}

\begin{proof}
By the Fourier inversion formula, for all $f \in \mathcal{S}(\mathbb{R}^d)$ and all $x \in \mathbb{R}^d$, 
\begin{align*}
D^{\alpha-1}(f)(x) = \frac{1}{(2\pi)^d} \int_{\mathbb{R}^d} \mathcal{F}(f)(\xi) e^{i \langle x; \xi \rangle} \tau_\alpha(\xi) d\xi, \quad \tau_\alpha(\xi) = \int_{\mathbb{R}^d} u \left(e^{i \langle u ; \xi \rangle }-1\right) \nu_\alpha(du). 
\end{align*}
Next, let us compute the vectorial symbol $\tau_\alpha$.~Passing to spherical coordinates, using Fubini's theorem and the symmetry of $\sigma$, for all $\xi \in \mathbb{R}^d \setminus \{0\}$, 
\begin{align*}
\tau_\alpha(\xi) & = \int_{(0, + \infty) \times \mathbb{S}^{d-1}} y \left(e^{ir \langle y ; \xi \rangle} - 1\right) \frac{dr}{r^{\alpha}} \sigma(dy)  \\
& = \int_{\mathbb{S}^{d-1}} y \left( \int_{(0,+\infty)} \left(e^{i r \langle y ; \xi \rangle}-1\right) \frac{dr}{r^{\alpha}}\right)  \sigma(dy)  \\
& = i\int_{\mathbb{S}^{d-1}} y \left(\int_{(0,+\infty)} \sin \left(r \langle y ; \xi \rangle\right)  \frac{dr}{r^\alpha}\right) \sigma(dy)  \\
& = i\int_{\mathbb{S}^{d-1}} y |\langle y ; \xi \rangle|^{\alpha-1} \left(\int_{(0,+\infty)} \sin \left(r \operatorname{sign}\left( \langle y ; \xi \rangle\right) \right)\frac{dr}{r^\alpha}\right) \sigma(dy)  \\
& = i\int_{\mathbb{S}^{d-1}} y |\langle y ; \xi \rangle|^{\alpha-1} \operatorname{sign}\left( \langle y ; \xi \rangle\right) \sigma(dy) \int_0^{+\infty} \sin(r) \frac{dr}{r^\alpha}  \\
& = i \frac{\Gamma(2-\alpha) \cos\left(\frac{\alpha \pi}{2}\right)}{1- \alpha} \int_{\mathbb{S}^{d-1}} y |\langle y ; \xi \rangle|^{\alpha-1} \operatorname{sign}\left( \langle y ; \xi \rangle\right) \sigma(dy),
\end{align*}
where 
\begin{align*}
\int_0^{+\infty} \sin(r) \frac{dr}{r^\alpha} = \dfrac{\Gamma(2- \alpha) \cos \left(\frac{\alpha \pi}{2}\right)}{1- \alpha}. 
\end{align*}
Thus, for all $\xi \in \mathbb{R}^d \setminus \{0\}$, 
\begin{align*}
\left(2- \alpha \right) \tau_\alpha(\xi) & = i \frac{\Gamma(3-\alpha) \cos\left(\frac{\alpha \pi}{2}\right)}{1- \alpha} \int_{\mathbb{S}^{d-1}} y |\langle y ; \xi \rangle|^{\alpha-1} \operatorname{sign}\left( \langle y ; \xi \rangle\right) \sigma(dy) \\
& \longrightarrow \int_{\mathbb{S}^{d-1}} y  \langle y ; i \xi \rangle \sigma(dy), \quad \alpha \longrightarrow 2.
\end{align*}
Moreover, for all $\xi \in \mathbb{R}^d \setminus \{0\}$ and all $\alpha \in (1,2)$, 
\begin{align}\label{ineq:domination_assumption_Asymptotic1}
\left(2- \alpha \right) \|\tau_\alpha(\xi)\| & \leq \frac{\Gamma(3-\alpha) |\cos\left(\frac{\alpha \pi}{2}\right)|}{\alpha-1} \left(1+ \|\xi\|\right) \sigma\left(\mathbb{S}^{d-1}\right) \nonumber \\
& \leq \underset{\alpha \in (1,2)}{\sup}  \left(\frac{\Gamma(3-\alpha) |\cos\left(\frac{\alpha \pi}{2}\right)|}{\alpha-1} \right) \left(1+ \|\xi\|\right) \sigma\left(\mathbb{S}^{d-1}\right). 
\end{align}
The Fourier inversion formula concludes the proof of Equation \eqref{eq:pointwise_limit_BBManisotropic}.~Finally, the Parseval-Plancherel identity combined with the pointwise convergence of the Fourier multiplier and with the bound \eqref{ineq:domination_assumption_Asymptotic1} gives Equation \eqref{eq:L2_limit_BBManisotropic}.  
\end{proof}
\noindent
Let $H^1(\mathbb{R}^d,dx)$ be the set (of equivalence classes) of functions defined via 
\begin{align}\label{eq:definition_H1Sobolev}
H^1(\mathbb{R}^d,dx) = \{f \in L^2(\mathbb{R}^d,dx) : \, \nabla(f) \in L^2(\mathbb{R}^d,\mathbb{R}^d,dx)\}.
\end{align}
This vector space endowed with the norm, 
\begin{align}\label{eq:H1_norm}
\|f\|_{H^1(\mathbb{R}^d,dx)} :=\left( \|f\|^2_{L^2(\mathbb{R}^d,dx)} + \|\nabla(f)\|^2_{L^2(\mathbb{R}^d,dx)}\right)^{\frac{1}{2}},
\end{align}
is a Hilbert space and is also denoted by $W^{1,2}(\mathbb{R}^d,dx)$.  

\begin{lem}\label{lem:uniform_inequality_H1}
Let $\alpha_0 \in (1,2)$. Let $\sigma$ be the spherical component of the L\'evy measure $\nu_\alpha$. Then, for all $\alpha \in (\alpha_0, 2)$ and all $f \in H^1(\mathbb{R}^d,dx)$, $d \geq 1$,
\begin{align}\label{ineq:uniform_inequality_H1}
\| (2-\alpha) D^{\alpha-1}(f)\|_{L^2(\mathbb{R}^d,dx)} \leq C(\alpha_0, \sigma,d) \|f\|_{H^1(\mathbb{R}^d,dx)}, 
\end{align}
for some $C(\alpha_0, \sigma,d)>0$ which depends on $\alpha_0$, $\sigma$ and $d$. Moreover, for all $f \in H^1(\mathbb{R}^d,dx)$, 
\begin{align}\label{eq:target_space_H1}
\| (2-\alpha) D^{\alpha-1}(f)\|_{L^2(\mathbb{R}^d,dx)} \underset{\alpha \rightarrow 2^-}{\longrightarrow} \|D_\sigma(f)\|_{L^2(\mathbb{R}^d,dx)}. 
\end{align}
\end{lem}

\begin{proof}
Let $\alpha_0 \in (1,2)$, let $f \in H^1(\mathbb{R}^d,dx)$ and let $(f_n)_{n \geq 1}$ be a sequence of functions in the Schwartz space which converges, as $n$ tends to $+\infty$, to $f$ in $H^1(\mathbb{R}^d,dx)$.~Lemma \ref{lem:interpolation_inequality} gives $D^{\alpha-1}(f) \in L^2(\mathbb{R}^d,\mathbb{R}^d,dx)$ and
\begin{align*}
\| D^{\alpha-1}(f) \|_{L^2(\mathbb{R}^d,dx)} & \leq \left( \frac{2}{\alpha-1} + \frac{1}{2-\alpha} \right) \sigma\left(\mathbb{S}^{d-1}\right) \|f\|^{2-\alpha}_{L^2(\mathbb{R}^d,dx)} \|\nabla(f)\|^{\alpha-1}_{L^2(\mathbb{R}^d,dx)} \\
& \leq \left( \frac{2}{\alpha-1} + \frac{1}{2-\alpha} \right) \sigma\left(\mathbb{S}^{d-1}\right) \|f\|_{H^1(\mathbb{R}^d,dx)}. 
\end{align*} 
Thus, for all $\alpha \in (\alpha_0,2)$, 
\begin{align*}
\| (2-\alpha) D^{\alpha-1}(f) \|_{L^2(\mathbb{R}^d,dx)} & \leq \underset{\alpha \in (\alpha_0,2)}{\sup} \left( \frac{2(2-\alpha)}{\alpha-1} + 1 \right)  \sigma\left(\mathbb{S}^{d-1}\right) \|f\|_{H^1(\mathbb{R}^d,dx)},
\end{align*}
so that the constant $C(\alpha_0, \sigma,d)$ in inequality \eqref{ineq:uniform_inequality_H1} is given by 
\begin{align*}
C(\alpha_0, \sigma,d) : = \underset{\alpha \in (\alpha_0,2)}{\sup} \left( \frac{2(2-\alpha)}{\alpha-1} + 1 \right)  \sigma\left(\mathbb{S}^{d-1}\right). 
\end{align*}
Now, for all $n \geq 1$ and all $\alpha \in (\alpha_0 , 2)$, 
\begin{align*}
\left| \| (2-\alpha) D^{\alpha-1}(f) \|_{L^2(\mathbb{R}^d,dx)} - \|D_\sigma(f)\|_{L^2(\mathbb{R}^d,dx)} \right| & \leq \bigg| \| (2-\alpha) D^{\alpha-1}(f_n) \|_{L^2(\mathbb{R}^d,dx)} \\
&\quad\quad -\| (2-\alpha) D^{\alpha-1}(f) \|_{L^2(\mathbb{R}^d,dx)} \bigg| \\
&\quad\quad +\bigg| \|D_\sigma(f_n)\|_{L^2(\mathbb{R}^d,dx)} \\
&\quad\quad -\| (2-\alpha) D^{\alpha-1}(f_n) \|_{L^2(\mathbb{R}^d,dx)} \bigg| \\
&\quad\quad +\bigg|  \|D_\sigma(f)\|_{L^2(\mathbb{R}^d,dx)} - \|D_\sigma(f_n)\|_{L^2(\mathbb{R}^d,dx)}\bigg|. 
\end{align*}
But, \eqref{ineq:uniform_inequality_H1} ensures that, for all $n \geq 1$ and all $\alpha \in (\alpha_0,2)$, 
\begin{align*}
 \bigg| \| (2-\alpha) D^{\alpha-1}(f_n) \|_{L^2(\mathbb{R}^d,dx)}-\| (2-\alpha) D^{\alpha-1}(f) \|_{L^2(\mathbb{R}^d,dx)} \bigg| & \leq (2-\alpha) \| D^{\alpha-1}(f_n-f) \|_{L^2(\mathbb{R}^d,dx)} \\
 & \leq C(\alpha_0,d,\sigma) \|f_n -f\|_{H^1(\mathbb{R}^d,dx)}.
\end{align*}
Moreover, for all $n \geq 1$, 
\begin{align*}
\bigg| \| D_\sigma(f_n) \|_{L^2(\mathbb{R}^d,dx)} - \|D_\sigma(f)\|_{L^2(\mathbb{R}^d,dx)} \bigg| \leq \| D_\sigma(f_n) - D_\sigma(f) \|_{L^2(\mathbb{R}^d,dx)} \leq \sigma\left(\mathbb{S}^{d-1}\right)\|f_n -f\|_{H^1(\mathbb{R}^d,dx)}.
\end{align*}
Thus, thanks to Lemma \ref{lem:pointwise_L2_BBM_result}, for all $n \geq 1$
\begin{align*}
\underset{\alpha \rightarrow 2^-}{\limsup} \bigg| \| (2-\alpha) D^{\alpha-1}(f) \|_{L^2(\mathbb{R}^d,dx)} - \|D_\sigma(f)\|_{L^2(\mathbb{R}^d,dx)} \bigg| & \leq C(\alpha_0,d,\sigma) \|f_n - f\|_{H^1(\mathbb{R}^d,dx)} \\
& \quad\quad + \sigma\left(\mathbb{S}^{d-1}\right) \|f_n -f\|_{H^1(\mathbb{R}^d,dx)}. 
\end{align*}
Passing to the limit in the previous inequality as $n$ tends to $+\infty$ concludes the proof of the lemma. 
\end{proof}
\noindent 
To analyze the asymptotic behavior of the operator $D^{\alpha-1}$ when $\alpha$ tends to $1^+$, we study the boundedness of $\mathcal{R}_\sigma$ given by \eqref{eq:avRT} on the classical Lebesgue spaces $L^p(\mathbb{R}^d,dx)$, with $p \in (1,+\infty)$. This linear operator can be understood as an (an)isotropic version of the vectorial Riesz transform based on the following spherical representation of its multiplier: for all $\xi \in \mathbb{R}^d \setminus \{0\}$, 
\begin{align}\label{eq:spherical_representation_vRT_multiplier}
- \frac{i\xi}{\|\xi\|} = - i \frac{\pi}{2} \dfrac{\Gamma\left(\frac{d+1}{2}\right)}{\pi^{\frac{d+1}{2}}} \int_{\mathbb{S}^{d-1}} y \operatorname{sign}\left( \langle y ; \xi \rangle\right) \sigma_L(dy).
\end{align}
Using the method of rotations as in \cite{DuoRubio_CRAS85}, one extends continuously the linear operator $\mathcal{R}_\sigma$ to $L^p(\mathbb{R}^d,dx)$, for all $p \in (1,+\infty)$.

\begin{lem}\label{lem:Lp_continuity_avRT}
Let $p \in (1, +\infty)$. Let $\sigma$ be a symmetric finite positive and non-degenerate measure on $\mathbb{S}^{d-1}$, $d \geq 1$. Let $\mathcal{R}_\sigma$ be the linear Fourier multiplier operator defined by \eqref{eq:avRT}. Then, for all $f \in \mathcal{S}(\mathbb{R}^d)$, 
\begin{align}\label{ineq:Lp_continuity_avRT}
\|\mathcal{R}_\sigma\left(f\right)\|_{L^p(\mathbb{R}^d,dx)} \leq C_{p,\sigma,d} \|f\|_{L^p(\mathbb{R}^d,dx)}, 
\end{align} 
where $C_{p,\sigma,d}$ is given by 
\begin{align}\label{eq:upper_bound_avRT}
C_{p,\sigma,d} := \frac{\pi}{2} \sigma(\mathbb{S}^{d-1}) C_p,  
\end{align}
with 
\begin{align}\label{eq:Pichorides_constant}
C_p : = \left\{
    \begin{array}{ll}
        \tan(\frac{\pi}{2p}), & p \in (1,2],\\
        \cot(\frac{\pi}{2p}), & p\geq 2.\\
    \end{array}
\right.
\end{align}
\end{lem}

\begin{proof}
By Fubini's theorem and the Fourier inversion formula, for all $f \in \mathcal{S}(\mathbb{R}^d)$ and all $x \in \mathbb{R}^d$, 
\begin{align*}
\mathcal{R}_\sigma(f)(x) & =  i \frac{\pi}{2} \frac{1}{(2\pi)^d} \int_{\mathbb{S}^{d-1}} y \left(\int_{\mathbb{R}^d} \mathcal{F}(f)(\xi) e^{i \langle \xi ; x\rangle} \operatorname{sign}\left( \langle y ; \xi \rangle\right) d\xi \right) \sigma(dy) \\
& = - \frac{\pi}{2} \int_{\mathbb{S}^{d-1}} y H_y(f)(x) \sigma(dy), 
\end{align*}
where $H_y(f)$ is the directional Hilbert transform of $f$ in the direction $y \in \mathbb{S}^{d-1}$ with Fourier multiplier given by $- i \operatorname{sign}\left(\langle y ; \xi \rangle\right)$, for all $\xi \in \mathbb{R}^d$. Next, let $\tilde{\sigma}$ be the probability measure on $\mathbb{S}^{d-1}$ defined by $\tilde{\sigma}(dy) = \sigma(dy)/ \sigma(\mathbb{S}^{d-1})$. By H\"older's inequality, for all $f \in \mathcal{S}(\mathbb{R}^d)$ and all $x \in \mathbb{R}^d$, 
\begin{align*}
\left\| \mathcal{R}_\sigma(f)(x) \right\| &\leq \frac{\pi}{2} \sigma\left(\mathbb{S}^{d-1}\right) \left(\int_{\mathbb{S}^{d-1}} |\langle y ; \lambda \rangle|^q \tilde{\sigma}(dy)\right)^{\frac{1}{q}} \left(\int_{\mathbb{S}^{d-1}} |H_y(f)(x)|^p \tilde{\sigma}(dy)\right)^{\frac{1}{p}} \\ 
& \leq \frac{\pi}{2} \sigma\left(\mathbb{S}^{d-1}\right) \left(\int_{\mathbb{S}^{d-1}} |H_y(f)(x)|^p \tilde{\sigma}(dy)\right)^{\frac{1}{p}},
\end{align*}
for some $\lambda \in \mathbb{S}^{d-1}$.~Taking the $p$-th power and integrating with respect to the $d$-dimensional Lebesgue measure, 
\begin{align*}
\int_{\mathbb{R}^d} \left\| \mathcal{R}_\sigma(f)(x) \right\|^p dx \leq \left(\frac{\pi}{2} \sigma\left(\mathbb{S}^{d-1}\right)\right)^p \int_{\mathbb{S}^{d-1}} \left(\int_{\mathbb{R}^d}|H_y(f)(x)|^p dx \right) \tilde{\sigma}(dy).
\end{align*}
Now, thanks to \cite[Theorem $4.1$]{Pichorides_SM72}, for all $y \in \mathbb{S}^{d-1}$
\begin{align*}
\|H_y(f)\|_{L^p(\mathbb{R}^d,dx)} \leq C_p \|f\|_{L^p(\mathbb{R}^d,dx)}, 
\end{align*}
where $C_p$ is given by \eqref{eq:Pichorides_constant}. Then, for all $f \in \mathcal{S}(\mathbb{R}^d)$, 
\begin{align*}
\|\mathcal{R}_{\sigma}(f)\|_{L^p(\mathbb{R}^d,dx)} \leq \frac{\pi}{2} \sigma (\mathbb{S}^{d-1}) C_p \|f\|_{L^p(\mathbb{R}^d,dx)}.
\end{align*}
This concludes the proof of the lemma. 
\end{proof}
\noindent 
Based on the Fourier multiplier of $D^{\alpha-1}$ and the previous discussion, we are now ready to analyze the pointwise behavior and the $L^2$-behavior of $D^{\alpha-1}(f)$, for $f \in \mathcal{S}(\mathbb{R}^d)$, when $\alpha$ tends to $1^+$. 

\begin{lem}\label{lem:pointwise_L2_MS_result}
Let $\alpha \in (1,2)$ and let $\sigma$ be the spherical part of the L\'evy measure $\nu_\alpha$. Then, for all $f \in \mathcal{S}(\mathbb{R}^d)$ and all $x \in \mathbb{R}^d$, $d \geq 1$, 
\begin{align}\label{eq:pointwise_limit_MSanisotropic}
D^{\alpha-1}(f)(x) \underset{\alpha \rightarrow 1^+}{\longrightarrow} \mathcal{R}_\sigma(f)(x).
\end{align}
Moreover, for all $f \in \mathcal{S}(\mathbb{R}^d)$, 
\begin{align}\label{eq:L2_limit_MSanisotropic}
\left\| D^{\alpha-1}(f) \right\|_{L^2(\mathbb{R}^d,dx)} \underset{\alpha \rightarrow 1^+}{\longrightarrow}  \left\|\mathcal{R}_\sigma(f)\right\|_{L^2(\mathbb{R}^d,dx)}.
\end{align}
\end{lem}

\begin{proof}
As in the proof of Lemma \ref{lem:pointwise_L2_BBM_result}, for $\alpha \in (1,2)$, all $f \in \mathcal{S}(\mathbb{R}^d)$ and all $x \in \mathbb{R}^d$, 
\begin{align*}
D^{\alpha-1}(f)(x) = \frac{1}{(2\pi)^d} \int_{\mathbb{R}^d} \mathcal{F}(f)(\xi) e^{i \langle x ;  \xi \rangle} \tau_\alpha(\xi) d\xi, 
\end{align*}
where, for all $\xi \in \mathbb{R}^d$, 
\begin{align*}
\tau_\alpha(\xi) = i \frac{\Gamma(2-\alpha) \cos\left(\frac{\alpha \pi}{2}\right)}{1- \alpha} \int_{\mathbb{S}^{d-1}} y |\langle y ; \xi \rangle|^{\alpha-1} \operatorname{sign}\left( \langle y ; \xi \rangle\right) \sigma(dy). 
\end{align*}
Since $\cos (\alpha \pi/2) / (1-\alpha) \rightarrow \pi/2$ as $\alpha$ tends to $1^+$, for all $\xi \in \mathbb{R}^d \setminus \{0\}$, 
\begin{align*}
\tau_\alpha(\xi) \longrightarrow  i \frac{\pi}{2} \int_{\mathbb{S}^{d-1}} y \operatorname{sign}\left( \langle y ; \xi \rangle\right) \sigma(dy).
\end{align*}
Let $\alpha_1 \in (1,2)$. Moreover, for all $\xi \in \mathbb{R}^d \setminus \{0\}$ and all $\alpha \in (1,\alpha_1)$,
\begin{align}\label{ineq:domination_assumption_asymptotic2}
\left\| \tau_\alpha(\xi) \right\| & \leq \frac{\Gamma(2-\alpha) |\cos\left(\frac{\alpha \pi}{2}\right)|}{\alpha-1} \|\xi\|^{\alpha-1} \sigma\left(\mathbb{S}^{d-1}\right) \nonumber \\
& \leq \sup_{\alpha \in (1,\alpha_1)} \left(\frac{\Gamma(2-\alpha) |\cos\left(\frac{\alpha \pi}{2}\right)|}{\alpha-1}\right) \left(1+ \|\xi\|\right) \sigma \left(\mathbb{S}^{d-1}\right). 
\end{align} 
A standard application of the Lebesgue dominated convergence theorem together with the Fourier inversion formula concludes the proof of Equation \eqref{eq:pointwise_limit_MSanisotropic}. Finally, the Parseval-Plancherel identity combined with the pointwise convergence of the Fourier multiplier and with the bound \eqref{ineq:domination_assumption_asymptotic2} give Equation \eqref{eq:L2_limit_MSanisotropic}.  
\end{proof}
\noindent
Next, we want to extend the previous weak convergence results to strong convergence in $L^p(\mathbb{R}^d, \mathbb{R}^d,dx)$, for all $p \in (1, +\infty)$. It is straightforward to check that, for all $\xi \in \mathbb{R}^d \setminus \{0\}$, 
\begin{align*}
\tau_\alpha(\xi) : = \int_{\mathbb{R}^d} \left(e^{i \langle u ; \xi \rangle}-1\right) u \nu_\alpha(du) = i \alpha \nabla(\sigma_\alpha)(\xi) \sigma_\alpha(\xi)^{\alpha-1},
\end{align*}
which implies, for all $f \in \mathcal{C}_c^{\infty}(\mathbb{R}^d)$ and all $x \in \mathbb{R}^d$, 
\begin{align*}
D^{\alpha-1}(f)(x) = \mathcal{R}_\alpha \circ \left(-\mathcal{A}_\alpha \right)^{\frac{\alpha-1}{\alpha}}(f)(x). 
\end{align*}
Then, for all $p \in (1,+\infty)$ and all $f \in \mathcal{C}_c^{\infty}(\mathbb{R}^d)$, 
\begin{align*}
\left\| \mathcal{R}_\alpha \circ \left(-\mathcal{A}_\alpha \right)^{\frac{\alpha-1}{\alpha}}(f) - \mathcal{R}_\sigma(f) \right\|_{L^p(\mathbb{R}^d,dx)} &\leq \left\| \mathcal{R}_\alpha \right\|_{L^p(\mathbb{R}^d,dx)\rightarrow L^p(\mathbb{R}^d,dx)} \left\|\left(-\mathcal{A}_\alpha \right)^{\frac{\alpha-1}{\alpha}}(f) - f \right\|_{L^p(\mathbb{R}^d,dx)} \\
&\quad\quad + \|\mathcal{R}_\alpha(f)-\mathcal{R}_\sigma(f)\|_{L^p(\mathbb{R}^d,dx)}.
\end{align*}
In the rotationally invariant case, the Riesz transform-type operator $\mathcal{R}_\alpha^{\operatorname{rot}}$ admits the following Fourier multiplier: for all $\xi \in \mathbb{R}^d$ such that $\xi \ne 0$, 
\begin{align*}
m_\alpha^{\operatorname{rot}}(\xi) = \frac{\alpha}{2^{\frac{1}{\alpha}}} \frac{i \xi}{\|\xi\|}, 
\end{align*} 
which agrees with the following normalization: for all $\xi \in \mathbb{R}^d$, 
\begin{align*}
\int_{\mathbb{S}^{d-1}} \left| \langle \xi ;y \rangle \right|^\alpha \lambda_1(dy) = \frac{\|\xi\|^\alpha}{2}. 
\end{align*}
Thus, for all $p \in (1, +\infty)$,  
\begin{align*}
\|\mathcal{R}^{\operatorname{rot}}_\alpha\|_{L^p\left(\mathbb{R}^d,dx\right)\rightarrow L^p\left(\mathbb{R}^d,dx\right)} = \frac{\alpha}{2^{\frac{1}{\alpha}}} \|\mathcal{R}\|_{L^p\left(\mathbb{R}^d,dx\right)\rightarrow L^p\left(\mathbb{R}^d,dx\right)} \longrightarrow \frac{1}{2} \|\mathcal{R}\|_{L^p\left(\mathbb{R}^d,dx\right)\rightarrow L^p\left(\mathbb{R}^d,dx\right)},
\end{align*} 
as $\alpha$ tends to $1^+$. The previous choice of normalization corresponds to the L\'evy measure $\nu_\alpha^{\operatorname{rot}}$ given by \eqref{eq:Levy_Rot} and with spherical part equal to $\sigma(dy) = c_{\alpha,d} \sigma_L(dy)$. Thus, in this situation, the Fourier multiplier of the anisotropic Riesz transform $\mathcal{R}_{\sigma_\ell}$ is given, for all $\xi \in \mathbb{R}^d\setminus \{0\}$, by 
\begin{align*}
m_{\sigma_\ell}(\xi) = \frac{i}{2} \frac{\pi}{2} \dfrac{\Gamma\left(\frac{1+d}{2}\right)}{\pi^{\frac{d+1}{2}}} \int_{\mathbb{S}^{d-1}} y \operatorname{sign}\left( \langle y ; \xi \rangle\right) \sigma_L(dy) = \frac{1}{2} \frac{i\xi}{\|\xi\|}. 
\end{align*}
Hence, for all $p \in (1,+\infty)$ and all $f \in \mathcal{C}_c^{\infty}(\mathbb{R}^d)$, 
\begin{align*}
 \|\mathcal{R}^{\operatorname{rot}}_\alpha(f)-\mathcal{R}_{\sigma_\ell}(f)\|_{L^p(\mathbb{R}^d,dx)} \longrightarrow 0,
\end{align*}
as $\alpha$ tends to $1^+$. Moreover, using an adaptation of \cite[Lemma $21$]{BCGS_CRM22}, for all $f \in \mathcal{C}^{\infty}_c(\mathbb{R}^d)$ and all $p \in (1, +\infty)$, 
\begin{align*}
\left\|\left(-\mathcal{A}^{\operatorname{rot}}_\alpha \right)^{\frac{\alpha-1}{\alpha}}(f) - f \right\|_{L^p(\mathbb{R}^d,dx)} \longrightarrow 0,
\end{align*}
as $\alpha$ tends to $1^+$ (since the Fourier symbol of $(- \mathcal{A}^{\operatorname{rot}}_\alpha)^{(\alpha-1)/\alpha}$ is given by $\|\xi\|^{\alpha-1}/2^{(\alpha-1)/\alpha}$, for all $\xi \in \mathbb{R}^d$). This retrieves the conclusion of \cite[Lemma 21, Equation $(88)$]{BCGS_CRM22}. The next proposition provides an upper bound, for all $p \in (1, +\infty)$, on 
\begin{align*}
\|\mathcal{R}_\alpha\|_{L^p(\mathbb{R}^d,dx) \rightarrow L^p(\bbr,dx)} := \underset{f \in L^p(\mathbb{R}^d,dx) \setminus \{0\}}{\sup} \dfrac{\|\mathcal{R}_\alpha(f)\|_{L^p(\mathbb{R}^d,dx)}}{\|f\|_{L^p(\mathbb{R}^d,dx)}}, 
\end{align*}
complementing the results of \cite{AH20_4}. 

\begin{prop}\label{prop:Riesz_Lp_forall}
Let $\alpha \in (1,2)$, let $d \geq 1$ be an integer, let $\sigma$ be the spherical component of the L\'evy measure $\nu_\alpha$, let $\mathcal{R}_{\alpha}$ be the fractional Riesz transform-type operator given by \eqref{eq:Riesz_transform_frac} and let $\sigma_\alpha$ be given by \eqref{eq:rep_spectral_measure}. Then, for all $p \in (1, +\infty)$, 
\begin{align}\label{ineq:Riesz_Lp_forall}
\left\| \mathcal{R}_\alpha\right\|_{L^p(\mathbb{R}^d,dx) \rightarrow L^p(\mathbb{R}^d,dx)} \leq \dfrac{\alpha C_p}{2} \frac{1}{\pi^{\frac{d-1}{2}}} \Gamma\left(\dfrac{d+1}{2}\right)\int_{\mathbb{S}^{d-1}} \sigma_\alpha(y) \sigma_L(dy),
\end{align}
where $C_p>0$ is given by \eqref{eq:Pichorides_constant}. In particular, for all $p \in (1,+\infty)$,
\begin{align}\label{ineq:limsup_Riesz_Lp_forall}
\underset{\alpha \rightarrow 1^+}{\limsup} \left\| \mathcal{R}_\alpha\right\|_{L^p(\mathbb{R}^d,dx) \rightarrow L^p(\mathbb{R}^d,dx)} \leq  \dfrac{\pi C_p}{2} \sigma(\mathbb{S}^{d-1}).
\end{align}
\end{prop}

\begin{proof}
From \cite[Lemma $5.3$]{AH20_4}, for all $f \in \mathcal{C}^{\infty}_c(\mathbb{R}^d)$ and all $x \in \mathbb{R}^d$, 
\begin{align}\label{eq:superposition_dirSIO}
\mathcal{R}_\alpha(f)(x) = \frac{1}{2\Gamma\left(1 - \frac{1}{\alpha}\right)} \int_{\mathbb{R}^d} y T^{y}_\alpha(f)(x) \mu_\alpha(dy), 
\end{align}
where, for all $y \in \mathbb{R}^d$, 
\begin{align}\label{eq:dirSIO}
T^y_\alpha(f)(x) = \int_{(0, +\infty)} \left(f(x+t^{\frac{1}{\alpha}}y) - f(x - t^{\frac{1}{\alpha}}y)\right) \frac{dt}{t}. 
\end{align}
Then, by Minkowski's integral inequality, by transference and by \cite[Theorem $4.1$]{Pichorides_SM72}, for all $p \in (1, +\infty)$ and all $f \in \mathcal{C}_c^{\infty}(\mathbb{R}^d)$,  
\begin{align*}
\left\| \mathcal{R}_\alpha(f) \right\|_{L^p(\mathbb{R}^d,dx)} \leq \dfrac{\alpha \pi C_p}{2 \Gamma\left(1- \frac{1}{\alpha}\right)} \left(\int_{\mathbb{R}^d} \|y\| \mu_\alpha(dy)\right) \|f\|_{L^p(\mathbb{R}^d,dx)},
\end{align*}
where $C_p>0$ is given by \eqref{eq:Pichorides_constant}. Moreover, \cite[Theorem $6.1$]{Molchanov_JMA09}, 
\begin{align}\label{eq:radial_mean_NDS}
\int_{\mathbb{R}^d} \|y\| \mu_\alpha(dy) = \frac{1}{\pi^{\frac{d+1}{2}}} \Gamma\left(\dfrac{d+1}{2}\right) \Gamma\left(\frac{\alpha-1}{\alpha}\right) \int_{\mathbb{S}^{d-1}} \sigma_\alpha(y) \sigma_L(dy). 
\end{align} 
The upper bounds \eqref{ineq:Riesz_Lp_forall} and \eqref{ineq:limsup_Riesz_Lp_forall} easily follow.
\end{proof}
\noindent
The next lemma studies the pointwise convergence of $(-\mathcal{A}_\alpha)^{(\alpha-1)/\alpha}(f)$, for $f \in \mathcal{S}(\mathbb{R}^d)$, as $\alpha$ tends to $1^+$.

\begin{lem}\label{lem:pointwise_fracpow_generator}
Let $\alpha \in (1,2)$. Then, for all $f \in \mathcal{S}(\mathbb{R}^d)$, $d\geq 1$, and all $x \in \mathbb{R}^d$, 
\begin{align}\label{eq:pointwise_fracpow_generator}
\left(- \mathcal{A}_\alpha\right)^{\frac{\alpha-1}{\alpha}}(f)(x) \longrightarrow f(x), 
\end{align}
as $\alpha$ tends to $1^+$. 
\end{lem}

\begin{proof}
The proof of this lemma is based on the Fourier inversion formula and on the Lebesgue dominated convergence theorem.  Indeed, using \eqref{eq:Stheatgen_FourRep}, for all $f \in \mathcal{S}(\mathbb{R}^d)$ and all $x \in \mathbb{R}^d$, 
\begin{align*}
\left(- \mathcal{A}_\alpha\right)^{\frac{\alpha-1}{\alpha}}(f)(x) = \frac{1}{(2\pi)^d} \int_{\mathbb{R}^d} \mathcal{F}(f)(\xi) e^{i \langle x ; \xi \rangle} (\sigma_\alpha(\xi))^{\alpha-1} d\xi. 
\end{align*}
Now, for all $\xi \in \mathbb{R}^d$ such that $\xi \ne 0$, 
\begin{align*}
\sigma_\alpha(\xi)^{\alpha-1} \longrightarrow 1, 
\end{align*}
as $\alpha$ tends to $1^+$.  Let $\alpha_1 \in (1,2)$.  For all $\xi \in \mathbb{R}^d \setminus \{0\}$ and all $\alpha \in (1, \alpha_1)$, 
\begin{align*}
\sigma_\alpha(\xi)^{\alpha-1} & \leq \left(1+ \|\xi\|\right) \left( - \dfrac{\Gamma(2-\alpha)\cos\left(\alpha \frac{\pi}{2}\right)}{\alpha(\alpha-1)} \right)^{\frac{\alpha-1}{\alpha}} \sigma\left(\mathbb{S}^{d-1}\right)^{\frac{\alpha-1}{\alpha}} \\
& \leq \left(1+ \|\xi\|\right) \underset{\alpha \in (1, \alpha_1)}{\sup} \left( \left(- \dfrac{\Gamma(2-\alpha)\cos\left(\alpha \frac{\pi}{2}\right)}{\alpha(\alpha-1)} \right)^{\frac{\alpha-1}{\alpha}} \sigma\left(\mathbb{S}^{d-1}\right)^{\frac{\alpha-1}{\alpha}}\right).
\end{align*}
The Lebesgue dominated convergence theorem concludes the proof of the lemma. 
\end{proof}
\noindent
In the next proposition, let us investigate $L^p(\mathbb{R}^d,dx)$-bound on $\left(-\mathcal{A}_{\alpha}\right)^{\frac{\alpha-1}{\alpha}}(f)$ based on the fractional Riesz transform-type results of \cite[Proposition $3.9$]{AH20_4}. For this purpose, let us recall the second fractional Riesz transform-type operator introduced in \cite{AH20_4}: let $\mathbf{R}^\alpha$ be the linear operator defined, for all $f \in \mathcal{S}(\mathbb{R}^d)$ and all $x \in \mathbb{R}^d$, by 
\begin{align}\label{eq:second_Riesz_transform}
\mathbf{R}^\alpha(f)(x) = \nabla \circ \left(- \mathcal{A}_\alpha\right)^{-\frac{1}{\alpha}}(f)(x) = \frac{1}{2\Gamma\left(\frac{1}{\alpha}\right)} \int_{\mathbb{R}^d} F_\alpha(y) T^y_\alpha(f)(x) \mu_\alpha(dy), 
\end{align}   
where $T^y_\alpha(f)$ is given by \eqref{eq:dirSIO} and where $F_\alpha$ is defined, for all $y \in \mathbb{R}^d$, by 
\begin{align*}
F_\alpha(y) = \dfrac{- \nabla(p_\alpha)(y)}{p_\alpha(y)}. 
\end{align*}
Combining Minkowski's integral inequality, a transference argument and \cite[Theorem $4.1$]{Pichorides_SM72}, for all $p \in (1, +\infty)$ and all $f \in \mathcal{S}(\mathbb{R}^d)$, 
\begin{align}\label{inerq:Lp_second_Riesz_transform}
\| \mathbf{R}^\alpha(f) \|_{L^p(\mathbb{R}^d,dx)} \leq \dfrac{\alpha C_p}{2 \Gamma\left(\frac{1}{\alpha}\right)} \left(\int_{\mathbb{R}^d} \|F_\alpha(y)\| \mu_\alpha(dy)\right) \|f\|_{L^p(\mathbb{R}^d,dx)},
\end{align} 
where $C_p$ is given by \eqref{eq:Pichorides_constant}. 

\begin{lem}\label{lem:Lp_bound_FractionalPower}
Let $\alpha \in (1,2)$. Then, for all $p \in (1,+\infty)$ and all $f \in \mathcal{C}_c^{\infty}(\mathbb{R}^d)$, $d \geq 1$,
\begin{align}\label{ineq:Lp_bound_FractionalPower}
\left\| \left(- \mathcal{A}_\alpha\right)^{\frac{\alpha-1}{\alpha}}(f) \right\|_{L^p(\mathbb{R}^d,dx)} \leq \dfrac{C_q}{2\Gamma\left(\frac{1}{\alpha}\right)} \left( \int_{\mathbb{R}^d} \|F_\alpha(y)\| \mu_\alpha(dy) \right) \|D^{\alpha-1}(f)\|_{L^p(\mathbb{R}^d,dx)},
\end{align}
where $C_q$ is given by \eqref{eq:Pichorides_constant} with $q = p/(p-1)$.  
\end{lem}

\begin{proof}
The proof is similar to the one of \cite[Proposition 3.9]{AH20_4}.~Let $f,g \in \mathcal{C}_c^{\infty}(\mathbb{R}^d)$ and let $p,q \in (1, +\infty)$ be such that $1/p+1/q=1$. Then, by Fourier analytic arguments and \eqref{eq:Stheatgen_nice_decomposition}, 
\begin{align*}
 \int_{\mathbb{R}^d} (-\mathcal{A}_\alpha)^{\frac{\alpha-1}{\alpha}}(f)(x) \left(-\mathcal{A}_\alpha\right)^{\frac{1}{\alpha}}(g)(x) dx = \int_{\mathbb{R}^d} (-\mathcal{A}_\alpha)(f)(x) g(x) dx = \frac{1}{\alpha} \int_{\mathbb{R}^d} \langle D^{\alpha-1}(f)(x) ; \nabla(g)(x) \rangle dx.
\end{align*} 
Then, by H\"older's inequality, 
\begin{align*}
\left|  \int_{\mathbb{R}^d} (-\mathcal{A}_\alpha)^{\frac{\alpha-1}{\alpha}}(f)(x) \left(-\mathcal{A}_\alpha\right)^{\frac{1}{\alpha}}(g)(x) dx\right| & \leq \frac{1}{\alpha} \|D^{\alpha-1}(f)\|_{L^p(\mathbb{R}^d,dx)} \|\nabla(g)\|_{L^q(\mathbb{R}^d,dx)} \\
& \leq \dfrac{C_q}{2\Gamma\left(\frac{1}{\alpha}\right)} \left(\int_{\mathbb{R}^d} \|F_\alpha(y)\| \mu_\alpha(dy) \right)  \|D^{\alpha-1}(f)\|_{L^p(\mathbb{R}^d,dx)} \\
&\quad\quad \times \|\left(-\mathcal{A}_\alpha\right)^{\frac{1}{\alpha}}(g)\|_{L^q(\mathbb{R}^d,dx)}.
\end{align*}
Now, from \cite[Lemma $1$]{Russ_MS00}, the set $\mathcal{E}_\alpha$, defined by 
\begin{align*}
\mathcal{E}_\alpha = \left\{ (-\mathcal{A}_\alpha)^{\frac{1}{\alpha}}(g):\, g \in L^q(\mathbb{R}^d,dx) \cap D\left((-\mathcal{A}_\alpha)^{\frac{1}{\alpha}}\right)\right\},
\end{align*}
is dense in $L^q(\mathbb{R}^d,dx)$. Then, by duality, 
\begin{align*}
\left\| \left(- \mathcal{A}_\alpha\right)^{\frac{\alpha-1}{\alpha}}(f)\right\|_{L^p(\mathbb{R}^d,dx)} & = \underset{g^* \in \mathcal{E}_\alpha,\, \|g^*\|_{L^q(\mathbb{R}^d,dx)} \leq 1}{\sup} \left| \langle \left(- \mathcal{A}_\alpha\right)^{\frac{\alpha-1}{\alpha}}(f) ; g^* \rangle_{L^2(\mathbb{R}^d,dx)} \right| \\
& \leq \dfrac{C_q}{2\Gamma\left(\frac{1}{\alpha}\right)} \left(\int_{\mathbb{R}^d} \|F_\alpha(y)\| \mu_\alpha(dy) \right)  \|D^{\alpha-1}(f)\|_{L^p(\mathbb{R}^d,dx)}.
\end{align*}
\end{proof}

\begin{lem}\label{lem:technical-lemma_pointwise_decay}
Let $\alpha_1 \in (1,2)$, let $\alpha \in (1,\alpha_1)$ and let $\nu_\alpha$ be a non-degenerate symmetric $\alpha$-stable L\'evy measure on $\mathbb{R}^d$, $d\geq 1$, with spherical component $\sigma$. Let $D^{\alpha-1}$ be given by \eqref{eq:perim_fracGrad}. Then, for all $f \in \mathcal{S}(\mathbb{R}^d)$, all $x \in \mathbb{R}^d$ and all $m > 1$ integer,  
\begin{align}\label{eq:pointwise_decay}
\left\| D^{\alpha-1}(f)(x) \right\| & \leq C_{d,f,m,\alpha_1}\bigg[ \int_{\{|\langle x;y \rangle| \leq 1|\}} \dfrac{\sigma(dy)}{\left(1+ \|x_w\|_1 \right)^{m-1}} \nonumber \\
& \quad\quad + \int_{\{| \langle y;x \rangle |>1\}} \bigg( \frac{1}{(1+|x_y|)} \dfrac{1}{\left(1+ \|x_w\|_1\right)^{m-1}} + \frac{\|x\|_1}{\left(1+ \|x\|_1\right)^m}\bigg) \sigma(dy)\bigg],
\end{align}
where $\|x\|_1 = |\langle x;y \rangle| + \|x_w\|_1$ with $x_w \in y^{\perp}$, for $y \in \mathbb{S}^{d-1}$, and where $C_{d,f,m,\alpha_1}>0$ depends on $d$, $f$, $m$ and $\alpha_1$. Moreover, the right-hand side of inequality \eqref{eq:pointwise_decay} belongs to $L^p(\mathbb{R}^d,dx)$, for all $p \in (1, +\infty)$, as soon as $m$ is large enough.  
\end{lem}

\begin{proof}
Making use of the polar decomposition of the L\'evy measure $\nu_\alpha$ and the symmetry of $\sigma$, for all $f \in \mathcal{S}(\mathbb{R}^d)$ and all $x \in \mathbb{R}^d$, 
\begin{align}\label{eq:superposition_non-singular_int_operator}
D^{\alpha-1}(f)(x) & = \int_{\mathbb{R}^d} \left(f(x+u) - f(x)\right) u \nu_\alpha(du) \nonumber\\
& = \frac{1}{2}\int_{\mathbb{S}^{d-1}} y \left(\int_{(0,+\infty)}\left(f(x+ry) - f(x-ry)\right) \frac{dr}{r^{\alpha}}\right) \sigma(dy) \nonumber \\
& = \frac{1}{2} \int_{\mathbb{S}^{d-1}} y H^y_\alpha(f)(x) \sigma(dy), 
\end{align}
where, for all $y \in \mathbb{S}^{d-1}$, 
\begin{align*}
H^y_\alpha(f)(x) = \int_{(0,+\infty)} \left(f(x+ry)-f(x-ry)\right)\frac{dr}{r^\alpha}. 
\end{align*}
Let $y \in \mathbb{S}^{d-1}$ and let $x \in \mathbb{R}^d \setminus \{0\}$ be such that $x = x_y y + x_w $ where $x_w \in y^{\perp}$ is the orthogonal projection of $x$ on $y^{\perp}$. First, let us assume that $|x_y|>1$. Then, 
\begin{align}\label{ineq:cut-3}
\bigg| H^y_\alpha(f)(x) \bigg| & = \bigg| \int_{(0,+\infty)} \left(f(x+ry)-f(x-ry)\right)\frac{dr}{r^\alpha} \bigg| \nonumber \\
& \leq \int_{(0, |x_y|/2]} \left| f(x+ry)-f(x-ry) \right| \frac{dr}{r^\alpha} + \int_{(|x_y|/2,|x_y|)} \left| f(x+ry)-f(x-ry) \right| \frac{dr}{r^\alpha} \nonumber \\
& \quad\quad + \int_{[|x_y|, +\infty)} \left| f(x+ry)-f(x-ry) \right| \frac{dr}{r^\alpha}. 
\end{align}
Let us start with the first integral in the right-hand side of the previous inequality. As, 
\begin{align*}
f(x+ry)-f(x-ry) = \int_{-1}^{1} \frac{d}{dt} \left(f(x+try)\right) dt = \int_{-1}^{1} \langle \nabla(f)(x+try) ; ry \rangle dt, 
\end{align*}
then
\begin{align*}
\int_{(0, |x_y|/2)} \left| f(x+ry)-f(x-ry) \right| \frac{dr}{r^\alpha} &\leq \int_0^{\frac{|x_y|}{2}} \int_{-1}^{1} \| \nabla(f)(x+try) \| dt \frac{dr}{r^{\alpha-1}} \\
& \leq \int_0^{\frac{|x_y|}{2}} \int_{-1}^1 \dfrac{C_{d,f,m}}{\left(1+ \|x_w\|_1 + |x_y +tr|\right)^m}dt \frac{dr}{r^{\alpha-1}},
\end{align*}
for some integer $m$ strictly larger than $1$ and for some $C_{d,f,m}>0$ depending on $d$, $f$ and $m$. Moreover, 
\begin{align*}
|x_y +tr| \geq ||x_y|-|t|r| = |x_y| - |t|r \geq |x_y| - r \geq \frac{|x_y|}{2}. 
\end{align*}
Thus, 
\begin{align}\label{ineq:integral1}
\int_{(0, |x_y|/2)} \left| f(x+ry)-f(x-ry) \right| \frac{dr}{r^\alpha} &\leq 2^{m+1} \frac{C_{d,f,m}}{\left(1+ \|x\|_1\right)^m} \dfrac{\|x\|^{2-\alpha}_1}{2-\alpha}. 
\end{align}
Next, for the middle integral in the right-hand side of the inequality \eqref{ineq:cut-3}, 
\begin{align}\label{ineq:B1+B2}
 \int_{(|x_y|/2,|x_y|)} \left| f(x+ry)-f(x-ry) \right| \frac{dr}{r^\alpha} \leq B_1 + B_2,
\end{align}
where
\begin{align*}
B_1 :=  \int_{(|x_y|/2,|x_y|)} \left| f(x+ry)\right| \frac{dr}{r^\alpha}, \quad B_2 :=  \int_{(|x_y|/2,|x_y|)} \left|f(x-ry)\right| \frac{dr}{r^\alpha}. 
\end{align*}
Without loss of generality, let us treat the term $B_1$ only. Then, 
\begin{align*}
 \int_{(|x_y|/2,|x_y|)} \left| f(x+ry)\right| \frac{dr}{r^\alpha} & \leq \frac{2^\alpha}{|x_y|^\alpha}  \int_{(|x_y|/2,|x_y|)} \left| f(x+ry)\right| dr \\
 & \leq  \frac{2^\alpha}{|x_y|^\alpha} \int_{\frac{|x_y|}{2}}^{|x_y|} \dfrac{C_{d,f,m}}{\left(1 + \|x_w\|_1 + |x_y+r|\right)^m} dr. 
\end{align*}
Now, $|x_y+r| \geq \left| |x_y| -r \right| = |x_y| - r$, and so, 
\begin{align*}
 \int_{(|x_y|/2,|x_y|)} \left| f(x+ry)\right| \frac{dr}{r^\alpha} & \leq  \frac{2^\alpha}{|x_y|^\alpha} C_{d,f,m} \int_{\frac{|x_y|}{2}}^{|x_y|} \dfrac{dr}{\left(1 + \|x_w\|_1 +  |x_y| - r \right)^m} \\
 & \leq \frac{2^\alpha}{|x_y|^\alpha} \frac{C_{d,f,m}}{m-1}  \dfrac{1}{\left(1+ \|x_w\|_1 \right)^{m-1}}.
\end{align*}
Since $|x_y|>1$, 
\begin{align}\label{ineq:B1-bound}
 \int_{(|x_y|/2,|x_y|)} \left| f(x+ry)\right| \frac{dr}{r^\alpha}  \leq C'_{d,f,m}\dfrac{4^\alpha}{(1+ |x_y|)^\alpha}\dfrac{1}{\left(1+ \|x_w\|_1 \right)^{m-1}},
\end{align}
for some $C'_{d,f,m}>0$ depending on $d$, $f$ and $m$. Finally, let us deal with the last integral in \eqref{ineq:cut-3}. 
\begin{align}\label{ineq:C1+C2}
 \int_{[|x_y|, +\infty)} \left| f(x+ry)-f(x-ry) \right| \frac{dr}{r^\alpha} \leq C_1 + C_2,
\end{align}
where, 
\begin{align*}
C_1:=  \int_{[|x_y|, +\infty)} \left| f(x+ry)\right| \frac{dr}{r^\alpha}, \quad C_2 :=   \int_{[|x_y|, +\infty)} \left| f(x-ry)\right| \frac{dr}{r^\alpha}. 
\end{align*}
Without loss of generality, let us upper-bound $C_1$ since $C_2$ can be bounded similarly. Thus, 
\begin{align*}
C_1 & \leq \frac{1}{|x_y|^\alpha}  \int_{[|x_y|, +\infty)} \left| f(x+ry)\right| dr  \\
& \leq C_{d,f,m} \frac{1}{|x_y|^\alpha} \int_{[|x_y|, +\infty)} \dfrac{dr}{\left(1+ \|x_w\|_1 + |x_y+r|\right)^m}.
\end{align*}
Now, $|x_y+r| \geq | |x_y|- r| = r-|x_y|$, and so, 
\begin{align*}
C_1 &\leq C_{d,f,m} \frac{1}{|x_y|^\alpha} \int_{|x_y|}^{+\infty} \dfrac{dr}{\left(1+ \|x_w\|_1 + r-|x_y|\right)^m} \\
&\leq  \frac{C'_{d,f,m}}{|x_y|^\alpha} \dfrac{1}{\left(1+ \|x_w\|_1\right)^{m-1}},
\end{align*}
for some $C'_{d,f,m}>0$ depending on $d$, $f$ and $m$. Since $|x_y|>1$, 
\begin{align}\label{ineq:C1-bound}
C_1 \leq \frac{2^{\alpha}C'_{d,f,m}}{(1+|x_y|)^\alpha} \dfrac{1}{\left(1+ \|x_w\|_1\right)^{m-1}}. 
\end{align}
Combining the inequalities \eqref{ineq:cut-3}--\eqref{ineq:C1-bound}, for all $y \in \mathbb{S}^{d-1}$ and all $x \in \mathbb{R}^d$ such that $|x_y|>1$, 
\begin{align}\label{ineq:fractional_directional_hilbert}
\left| H_\alpha^y(f)(x)\right| &\leq C_{d,f,m} \bigg( \frac{2^{\alpha+1}}{(1+|x_y|)^\alpha} \dfrac{1}{\left(1+ \|x_w\|_1\right)^{m-1}}+\dfrac{4^{\alpha+1}}{(1+ |x_y|)^\alpha}\dfrac{1}{\left(1+ \|x_w\|_1 \right)^{m-1}} \nonumber \\ 
& \quad\quad\quad\quad + \frac{2^{m+1}}{\left(1+ \|x\|_1\right)^m} \dfrac{\|x\|^{2-\alpha}_1}{2-\alpha} \bigg).
\end{align}
Note that since $f \in \mathcal{S}(\mathbb{R}^d)$, $m$ can be chosen as large as we want. Next, let us study the case when $|x_y|\leq 1$. As previously, 
\begin{align}\label{ineq:cut-2}
\left|  H_\alpha^y(f)(x) \right| & \leq \int_0^{1} |f(x+ry) - f(x-ry)| \frac{dr}{r^{\alpha}} + \int_{1}^{+\infty} |f(x+ry) - f(x-ry)| \frac{dr}{r^{\alpha}}. 
\end{align}
Now, 
\begin{align}\label{ineq:integral2}
\int_0^{1} |f(x+ry) - f(x-ry)| \frac{dr}{r^{\alpha}} & \leq C_{d,f,m} \int_0^1 \int_{-1}^{1} \dfrac{dt}{\left(1+ \|x_w\|_1 + |x_y+rt| \right)^m} \frac{dr}{r^{\alpha-1}} \nonumber \\
& \leq \frac{2 C_{d,f,m}}{(2-\alpha)\left(1+ \|x_w\|_1\right)^m}.
\end{align}
Moreover, 
\begin{align}\label{ineq:D1+D2}
\int_{1}^{+\infty} |f(x+ry) - f(x-ry)| \frac{dr}{r^{\alpha}} \leq D_1 + D_2, 
\end{align}
where 
\begin{align*}
D_1 := \int_{1}^{+\infty} |f(x+ry)| \frac{dr}{r^{\alpha}}, \quad D_2 := \int_{1}^{+\infty} |f(x-ry)| \frac{dr}{r^{\alpha}}. 
\end{align*}
Then, 
\begin{align}\label{ineq:D1-bound}
D_1 & \leq C_{d,f,m} \int_1^{+\infty} \dfrac{dr}{\left(1+ \|x_w\|_1 + |x_y+r|\right)^m} \nonumber \\
& \leq \dfrac{C'_{d,f,m}}{\left(1+ \|x_w\|_1 \right)^{m-1}}. 
\end{align}
Combining \eqref{ineq:cut-2}--\eqref{ineq:D1-bound}, for all $y \in \mathbb{S}^{d-1}$ and all $x \in \mathbb{R}^d$ such that $|x_y| \leq 1$,
\begin{align*}
\left| H_\alpha^y(f)(x)\right| \leq C_{d,f,m} \dfrac{3-\alpha}{2-\alpha} \dfrac{1}{\left(1+ \|x_w\|_1 \right)^{m-1}}.
\end{align*}
Finally, for all $x \in \mathbb{R}^d$, 
\begin{align*}
\|D^{\alpha-1}(f)(x)\| & \leq \int_{\mathbb{S}^{d-1}} \left| H_\alpha^y(f)(x)\right| \sigma(dy) \\
& \leq \int_{\{|\langle x;y \rangle| \leq 1\}} \left| H_\alpha^y(f)(x)\right| \sigma(dy) + \int_{\{|\langle x;y \rangle|>1\}} \left| H_\alpha^y(f)(x)\right| \sigma(dy) \\
& \leq C_{d,f,m} \bigg[ \frac{3-\alpha}{2-\alpha}  \int_{\{|\langle x;y \rangle| \leq 1|\}} \dfrac{\sigma(dy)}{\left(1+ \|x_w\|_1 \right)^{m-1}} \\
& \quad\quad + \int_{\{| \langle y;x \rangle |>1\}} \bigg( \frac{1}{(1+|x_y|)^\alpha} \dfrac{1}{\left(1+ \|x_w\|_1\right)^{m-1}} + \frac{1}{\left(1+ \|x\|_1\right)^m} \dfrac{\|x\|^{2-\alpha}_1}{2-\alpha} \bigg) \sigma(dy) \bigg]. 
\end{align*}
Next, let $\alpha_1 \in (1,2)$. Then, for all $\alpha \in (1, \alpha_1)$ and all $x \in \mathbb{R}^d$, 
\begin{align*}
\|D^{\alpha-1}(f)(x)\| & \leq C_{d,f,m,\alpha_1}\bigg[ \int_{\{|\langle x;y \rangle| \leq 1|\}} \dfrac{\sigma(dy)}{\left(1+ \|x_w\|_1 \right)^{m-1}} \\
& \quad\quad + \int_{\{| \langle y;x \rangle |>1\}} \bigg( \frac{1}{(1+|x_y|)} \dfrac{1}{\left(1+ \|x_w\|_1\right)^{m-1}} + \frac{\|x\|_1}{\left(1+ \|x\|_1\right)^m}\bigg) \sigma(dy)\bigg].
\end{align*}
This concludes the proof of the lemma. 
\end{proof}
\noindent
Combining Lemma \ref{lem:pointwise_L2_MS_result} with Lemma \ref{lem:technical-lemma_pointwise_decay} and the Lebesgue dominated convergence theorem, for all $p \in (1,+\infty)$ and all $f \in \mathcal{S}(\mathbb{R}^d)$, 
\begin{align}\label{eq:strong_convergence_MS_Schwartz}
\underset{\alpha \rightarrow 1, \, \alpha >1}{\lim} \left\| D^{\alpha-1}(f) - R_\sigma(f) \right\|_{L^p(\mathbb{R}^d,dx)} = 0. 
\end{align}
In order to extend the previous result to the target functional space $\cup_{\alpha \in (1,2)} W^{\alpha-1,p}(\mathbb{R}^d,dx)$, an intermediary question to resolve is the following: let $1<\beta < \alpha <2$. Do the following continuous embeddings hold true (recall also Proposition \ref{prop:converse_inclusion}):
\begin{align}\label{eq:equivalence_continuous_embedding}
W^{\alpha-1, p}\left(\mathbb{R}^d , dx\right) \hookrightarrow W^{\beta-1,p}\left(\mathbb{R}^d,dx\right) \Leftrightarrow H_p^{\psi_\alpha, r(\alpha)}(\mathbb{R}^d) \hookrightarrow H_p^{\psi_\beta, r(\beta)}(\mathbb{R}^d)\quad ?
\end{align}
Let $T_{\alpha,\beta}$ be the linear operator defined, for all $f \in \mathcal{S}(\mathbb{R}^d)$ and all $x \in \mathbb{R}^d$, by 
\begin{align}\label{eq:Fourier_multiplier_ratio}
T_{\alpha,\beta}(f)(x) = \frac{1}{(2\pi)^d} \int_{\mathbb{R}^d} \mathcal{F}(f)(\xi) e^{i \langle \xi ; x \rangle} t_{\alpha,\beta}(\xi) d\xi, 
\end{align} 
where, for all $\xi \in \mathbb{R}^d$, 
\begin{align}\label{eq:fourier_multiplier_ratio}
t_{\alpha,\beta}(\xi) = \dfrac{\left(1 - \psi_\beta(\xi)\right)^{\frac{\beta-1}{\beta}}}{\left(1 - \psi_\alpha(\xi) \right)^{\frac{\alpha-1}{\alpha}}},
\end{align}
with $\psi_\alpha$ given by \eqref{eq:multiplier_fractional_laplacian}. Note that by \cite[Theorem $3.3.28$]{NJ02_2}, the second embedding in \eqref{eq:equivalence_continuous_embedding} is equivalent to the $L^p(\mathbb{R}^d,dx)$-boundedness of $T_{\alpha,\beta}$. 

\begin{lem}\label{lem:fractional_calculus-type_result}
Let $1<\beta<\alpha<2$ and let $I^{\alpha-\beta}$ be the classical Riesz potential operator defined in \eqref{eq:def_Riesz_Potential_Operator} with $s = \alpha-\beta$. Then, for all $f \in \mathcal{C}^{\infty}_c(\mathbb{R}^d)$, all $x \in \mathbb{R}^d$, $d \geq 1$, and all $k \in \{1, \dots, d\}$,
\begin{align}\label{lem:fractional_calculus-type_result}
D_k^{\beta-1, \operatorname{rot}}(f)(x) = \frac{\beta}{\alpha} I^{\alpha-\beta}\left(D_k^{\alpha-1, \operatorname{rot}}(f)\right)(x).
\end{align}
\end{lem} 

\begin{proof}
The proof follows from Fourier manipulation at the level of the multiplier of $D_k^{\beta-1, \operatorname{rot}}$. Indeed, for all $f \in \mathcal{C}^{\infty}_c(\mathbb{R}^d)$, all $\xi \in \mathbb{R}^d \setminus \{0\}$ and all $k \in \{1, \dots, d\}$,
\begin{align*}
\mathcal{F} \left(D_k^{\beta-1, \operatorname{rot}}(f)\right)(\xi) & = \frac{i \beta}{2} \dfrac{\xi_k}{\|\xi\|} \|\xi\|^{\beta-1} \mathcal{F}(f)(\xi) \\
& = \frac{\beta}{\alpha} \frac{i \alpha}{2} \dfrac{\xi_k}{\|\xi\|}\|\xi\|^{\alpha-1} \|\xi\|^{\beta-\alpha} \mathcal{F}(f)(\xi).
\end{align*}
This concludes the proof.
\end{proof}

\begin{prop}\label{prop:fractional_interpolation_inequality}
Let $1<\beta <\alpha<2$, let $D^{\alpha-1,\operatorname{rot}}$ be the fractional gradient operator with associated L\'evy measure $\nu_\alpha^{\operatorname{rot}}$ and let $p \in [1,+\infty)$. Then, for all $f \in \mathcal{C}^{\infty}_c(\mathbb{R}^d)$, $d \geq 2$, 
\begin{align}\label{ineq:fractional_interpolation_inequality}
\| D^{\beta-1,\operatorname{rot}}(f) \|_{L^p(\mathbb{R}^d,dx)} \leq c_{\alpha,\beta,d} \|f\|^{\frac{\alpha-\beta}{\alpha-1}}_{L^p(\mathbb{R}^d,dx)} \|D^{\alpha-1,\operatorname{rot}}(f)\|^{\frac{\beta-1}{\alpha-1}}_{L^p(\mathbb{R}^d,dx)}, 
\end{align}
where $c_{\alpha,\beta,d}>0$ depends on $\alpha$, $\beta$ and $d$. 
\end{prop}

\begin{proof}
The proof of this proposition is similar to the proof of the interpolation inequality of Lemma \ref{lem:interpolation_inequality} using Lemma \ref{lem:Fourier_Transform_HD} instead of the classical Taylor's formula. Let $R>0$. First, for all $x \in \mathbb{R}^d$, 
\begin{align*}
\left\| D^{\beta-1, \operatorname{rot}}(f)(x) \right\| \leq \int_{\|u\|\leq R} |f(x+u)-f(x)| \|u\| \nu^{\operatorname{rot}}_{\beta}(du) + \int_{\|u\| \geq R} | f(x+u) -f(x)| \|u \| \nu^{\operatorname{rot}}_{\beta}(du). 
\end{align*}
By the triangle inequality and by Minkowski's integral inequality,
\begin{align*}
\|D^{\beta-1, \operatorname{rot}}(f)\|_{L^p(\bbr^d,dx)} & \leq \int_{\|u\| \leq R} \|f(.+u)-f(.)\|_{L^p(\bbr^d,dx)} \|u\| \nu^{\operatorname{rot}}_{\beta}(du) \\
& \quad\quad + \int_{\| u \| \geq R} \|f(.+u)-f(.)\|_{L^p(\bbr^d,dx)} \|u \| \nu^{\operatorname{rot}}_{\beta}(du). 
\end{align*} 
For the right-most integral, as in the proof of lemma \ref{lem:interpolation_inequality}, 
\begin{align*}
\int_{\| u \| \geq R} \|f(.+u)-f(.)\|_{L^p(\bbr^d,dx)} \|u \| \nu^{\operatorname{rot}}_{\beta}(du) & \leq 2 \|f\|_{L^p(\bbr^d,dx)} \int_{\| u \| \geq R} \|u \| \nu^{\operatorname{rot}}_{\beta}(du) \\
& \leq  2 \|f\|_{L^p(\bbr^d,dx)} \dfrac{\sigma_L\left(\mathbb{S}^{d-1}\right) c_{\beta,d}}{(\beta-1)R^{\beta-1}},
\end{align*}
where $c_{\beta,d}>0$ is given by \eqref{eq:renorm}. Now, thanks to \eqref{eq:spatial_rep_comistefani}, for all $x \in \mathbb{R}^d$ and all $u \in B(0,R)$, 
\begin{align*}
f(x+u) - f(x) & = \frac{2}{\alpha} \left(2^{\alpha-1} \pi^{\frac{d}{2}} \dfrac{\Gamma\left(\frac{\alpha}{2}\right)}{\Gamma(\frac{d-\alpha+2}{2})}\right)^{-1} \int_{\mathbb{R}^d} \bigg\langle \dfrac{x+u-y}{\|x+u-y\|^{d+2-\alpha}} \\
& \quad\quad -\dfrac{x-y}{\|x-y\|^{d+2-\alpha}} ; D^{\alpha-1, \operatorname{rot}}(f)(y) \bigg\rangle dy \\
& = \frac{2}{\alpha} \left(2^{\alpha-1} \pi^{\frac{d}{2}} \dfrac{\Gamma\left(\frac{\alpha}{2}\right)}{\Gamma(\frac{d-\alpha+2}{2})}\right)^{-1} \int_{\mathbb{R}^d} \bigg\langle \dfrac{u-z}{\|u-z\|^{d+2-\alpha}} \\
& \quad\quad +\dfrac{z}{\|z\|^{d+2-\alpha}} ; D^{\alpha-1, \operatorname{rot}}(f)(x+z) \bigg\rangle dz.
\end{align*}
Then, for all $x \in \mathbb{R}^d$ and all $u \in B(0,R)$, 
\begin{align*}
\left| \Delta_u(f)(x) \right| \leq \frac{2}{\alpha} \left(2^{\alpha-1} \pi^{\frac{d}{2}} \dfrac{\Gamma\left(\frac{\alpha}{2}\right)}{\Gamma(\frac{d-\alpha+2}{2})}\right)^{-1}\int_{\mathbb{R}^d} \bigg\| \dfrac{u-z}{\|u-z\|^{d+2-\alpha}}+\dfrac{z}{\|z\|^{d+2-\alpha}} \bigg\| \bigg\| D^{\alpha-1, \operatorname{rot}}(f)(x+z)\bigg\|dz.
\end{align*}
Again, by Minkowski's integral inequality, for all $u \in B(0,R)$, 
\begin{align*}
\| \Delta_u(f)\|_{L^p(\bbr^d,dx)} \leq \frac{2}{\alpha} \left(2^{\alpha-1} \pi^{\frac{d}{2}} \dfrac{\Gamma\left(\frac{\alpha}{2}\right)}{\Gamma(\frac{d-\alpha+2}{2})}\right)^{-1} \left(\int_{\mathbb{R}^d} \bigg\| \dfrac{u-z}{\|u-z\|^{d+2-\alpha}}+\dfrac{z}{\|z\|^{d+2-\alpha}} \bigg\|dz \right) \bigg\| D^{\alpha-1, \operatorname{rot}}(f)\bigg\|_{L^p(\bbr^d,dx)}. 
\end{align*}
Using the rotational invariance of the $d$-dimensional Lebesgue measure as in \cite[proof of Proposition $3.14$]{comi_stefani}, 
\begin{align*}
\left(\int_{\mathbb{R}^d} \bigg\| \dfrac{u-z}{\|u-z\|^{d+2-\alpha}}+\dfrac{z}{\|z\|^{d+2-\alpha}} \bigg\|dz \right) & \leq \|u\|^{\alpha-1}  \left(\int_{\mathbb{R}^d} \bigg\| \dfrac{e_1-z}{\|e_1-z\|^{d+2-\alpha}}+\dfrac{z}{\|z\|^{d+2-\alpha}} \bigg\|dz \right) \\
& \leq C(\alpha,d) \|u\|^{\alpha-1},
\end{align*}
where $e_1 = (1,0, \dots, 0)^t$. Thus, 
\begin{align*}
\int_{\|u\|\leq R} \| \Delta_u(f)\|_{L^p(\bbr^d,dx)} \|u\| \nu_{\beta}^{\operatorname{rot}}(du)
& \leq  C(\alpha,d) \frac{2}{\alpha} \left(2^{\alpha-1} \pi^{\frac{d}{2}} \dfrac{\Gamma\left(\frac{\alpha}{2}\right)}{\Gamma(\frac{d-\alpha+2}{2})}\right)^{-1} \\
& \quad\quad \times \left(\int_{\|u\|\leq R}  \|u\|^{\alpha} \nu_{\beta}^{\operatorname{rot}}(du) \right) \bigg\| D^{\alpha-1, \operatorname{rot}}(f)\bigg\|_{L^p(\bbr^d,dx)} \\
& \leq C(\alpha,d) \frac{2}{\alpha} \left(2^{\alpha-1} \pi^{\frac{d}{2}} \dfrac{\Gamma\left(\frac{\alpha}{2}\right)}{\Gamma(\frac{d-\alpha+2}{2})}\right)^{-1} \\
& \quad\quad \times c_{\beta,d} \frac{R^{\alpha-\beta}}{\alpha-\beta} \sigma_L\left(\mathbb{S}^{d-1}\right) \bigg\| D^{\alpha-1, \operatorname{rot}}(f)\bigg\|_{L^p(\bbr^d,dx)}.
\end{align*}
Taking $R = \left(\|f\|_{L^p(\bbr^d,dx)}/ \|D^{\alpha-1}(f)\|_{L^p(\bbr^d,dx)}\right)^{\frac{1}{\alpha-1}}$, 
\begin{align*}
\left\| D^{\beta-1, \operatorname{rot}}(f) \right\|_{L^p(\bbr^d,dx)} & \leq 2 \dfrac{\sigma_L\left(\mathbb{S}^{d-1}\right) c_{\beta,d}}{(\beta-1)} \|f\|^{\frac{\alpha-\beta}{\alpha-1}}_{L^p(\bbr^d,dx)} \| D^{\alpha-1, \operatorname{rot}}(f) \|^{\frac{\beta-1}{\alpha-1}}_{L^p(\bbr^d,dx)} \\
&\quad\quad + C(\alpha,d) \frac{2}{\alpha} \left(2^{\alpha-1} \pi^{\frac{d}{2}} \dfrac{\Gamma\left(\frac{\alpha}{2}\right)}{\Gamma(\frac{d-\alpha+2}{2})}\right)^{-1} \frac{c_{\beta,d}}{\alpha-\beta} \sigma_L\left(\mathbb{S}^{d-1}\right) \\
&\quad\quad \times \|f\|^{\frac{\alpha-\beta}{\alpha-1}}_{L^p(\bbr^d,dx)} \| D^{\alpha-1, \operatorname{rot}}(f) \|^{\frac{\beta-1}{\alpha-1}}_{L^p(\bbr^d,dx)}.
\end{align*}
This concludes the proof of the proposition. 
\end{proof}

\begin{prop}\label{prop:fractional_interpolation_inequality2}
Let $1<\beta <\alpha<2$, let $D^{\alpha-1,\operatorname{rot}}$ be the fractional gradient operator with associated L\'evy measure $\nu_\alpha^{\operatorname{rot}}$ and let $p \in (1,+\infty)$. Then, for all $f \in \mathcal{C}^{\infty}_c(\mathbb{R}^d)$, $d \geq 2$,
\begin{align}\label{ineq:fractional_interpolation_inequality2}
\| D^{\beta-1,\operatorname{rot}}(f) \|_{L^p(\mathbb{R}^d,dx)} \leq C_{\alpha,\beta,d,p} \|f\|^{\frac{\alpha-\beta}{\alpha-1}}_{L^p(\mathbb{R}^d,dx)} \|D^{\alpha-1,\operatorname{rot}}(f)\|^{\frac{\beta-1}{\alpha-1}}_{L^p(\mathbb{R}^d,dx)}, 
\end{align}
where
\begin{align}\label{eq:cste_fractional_interpolation_inequality2}
C_{\alpha,\beta,d,p} &:= \frac{c_{\beta,d}}{2} \sigma_L \left(\mathbb{S}^{d-1}\right) C \max\left(p, (p-1)^{-1}\right) \left( 2 + \beta\left(4 +  \frac{2}{2-\beta} \right)\right) \nonumber \\
& \quad\quad + C(\alpha,d) \frac{2}{\alpha} \left(2^{\alpha-1} \pi^{\frac{d}{2}} \dfrac{\Gamma\left(\frac{\alpha}{2}\right)}{\Gamma(\frac{d-\alpha+2}{2})}\right)^{-1} \frac{c_{\beta,d}}{\alpha-\beta} \sigma_L\left(\mathbb{S}^{d-1}\right),
\end{align}
and where $C>0$ is a numerical constant. 
\end{prop}

\begin{proof}
As in the proof of Lemma \ref{lem:technical-lemma_pointwise_decay}, for all $f \in \mathcal{C}^{\infty}_c(\bbr^d)$, all $R>0$ and all $x \in \bbr^d$, 
\begin{align}\label{eq:decomposition_method_rotation}
&D^{\beta-1, \operatorname{rot}}(f)(x) = \int_{B(0,R)} (f(x+u)-f(x)) u \nu_\beta^{\operatorname{rot}}(du) + \int_{B(0,R)^c} (f(x+u)-f(x)) u \nu_\beta^{\operatorname{rot}}(du) \nonumber \\
& = \int_{B(0,R)} (f(x+u)-f(x)) u \nu_\beta^{\operatorname{rot}}(du) + \frac{1}{2} c_{\beta,d} \int_{\mathbb{S}^{d-1}} y \left(\int_{R}^{+\infty} \left(f(x+ry) - f(x-ry)\right) \frac{dr}{r^\beta}\right) \sigma_L(dy), 
\end{align}
where $c_{\beta,d}>0$ is given by \eqref{eq:renorm}.~Next, for all $y \in \mathbb{S}^{d-1}$, let $H^y_{\beta,R}$ be the truncated fractional directional Hilbert transform defined, for all $f \in\mathcal{C}_c^{\infty}(\bbr^d)$ and all $x \in \bbr^d$, by 
\begin{align}\label{eq:trunc_frac_direc_HT}
H^y_{\beta,R}(f)(x) = \int_{R}^{+\infty} \left(f(x+ry) - f(x-ry)\right) \frac{dr}{r^\beta}. 
\end{align}
Let $e_1=(1, 0, \cdots, 0)^t$ and let $O$ be a $d \times d$ rotation matrix. Then, for all $f \in \mathcal{C}^\infty_c(\bbr^d)$ and all $x \in \mathbb{R}^d$, 
\begin{align}\label{eq:reduction_e1}
H^{O(e_1)}_{\beta,R}(f)(x) =  \int_{R}^{+\infty} \left(f\left(x+r O(e_1)\right) - f\left(x-r O(e_1)\right)\right) \frac{dr}{r^\beta} = H^{e_1}_{\beta,R}(f\circ O)\left(O^{-1}x\right). 
\end{align} 
Using the previous equality and the rotational invariance of the $d$-dimensional Lebesgue measure, let us study the mapping property of the linear operator $H^{e_1}_{\beta,R}$ on the Lebesgue spaces $L^p(\mathbb{R}^d,dx)$, for all $p \in (1, +\infty)$. For all $x = (x_1, \dots, x_d) \in \bbr^d$, let $\underline{x_1} = (x_2, \dots, x_d)$. Then, for all $f \in \mathcal{C}^{\infty}_c(\bbr^d)$ and all $x \in \mathbb{R}^d$, 
\begin{align*}
H^{e_1}_{\beta,R}(f)\left(x\right) & = \int_R^{+\infty} \left(f(x+re_1) - f(x-re_1)\right) \frac{dr}{r^\beta} \\
& = \int_R^{+\infty} \left(f\left(x_1+r, x_2, \dots, x_d\right) - f\left(x_1 - r, x_2,\dots, x_d\right)\right) \frac{dr}{r^\beta} \\
& = \int_R^{+\infty} \left(f_{\underline{x_1}}\left(x_1+r\right) - f_{\underline{x_1}}\left(x_1 - r\right)\right) \frac{dr}{r^\beta} , \quad f_{\underline{x_1}}(x_1) = f(x_1, x_2, \dots, x_d). 
\end{align*}
From these last observations, let us study the one-dimensional multiplier operator defined, for all $f \in \mathcal{S}(\bbr)$ and all $x \in \bbr$, by
\begin{align}\label{eq:one_dimensional_trunc_frac_Hilbert}
\mathcal{H}^\beta_R(f)(x) = \int_R^{+\infty} (f(x+r)-f(x-r)) \frac{dr}{r^\beta}. 
\end{align} 
By the Fourier inversion formula, for all $f \in \mathcal{S}(\bbr)$ and all $x \in \bbr$, 
\begin{align*}
\mathcal{H}^\beta_R(f)(x) & = \frac{1}{2\pi} \int_{\bbr} \mathcal{F}(f)(\xi) e^{ix\xi} \left(\int_R^{+\infty} (e^{i \xi r} - e^{-i\xi r}) \frac{dr}{r^\beta}\right) d\xi \\
& =  \frac{2i}{2\pi} \int_{\bbr} \mathcal{F}(f)(\xi) e^{ix\xi} \left(\int_R^{+\infty} \sin\left(r\xi\right) \frac{dr}{r^\beta}\right) d\xi \\
& =  \frac{1}{2\pi} \int_{\bbr} \mathcal{F}(f)(\xi) e^{ix\xi} h_R^\beta(\xi) d\xi, \quad  h_R^\beta(\xi) = 2i\int_R^{+\infty} \sin\left(r\xi\right) \frac{dr}{r^\beta}.
\end{align*}
Now, for all $R>0$ and all $\xi \in \bbr \setminus \{0\}$, 
\begin{align*}
h_R^\beta(\xi) & = 2i\int_R^{+\infty} \sin\left(r\xi\right) \frac{dr}{r^\beta} = 2i \operatorname{sign}\left(\xi\right) \int_R^{+\infty} \sin\left(r|\xi|\right) \frac{dr}{r^\beta} \\
& = 2i \operatorname{sign}\left(\xi\right) |\xi|^{\beta-1} \int_{R|\xi|}^{+\infty} \sin\left(r\right) \frac{dr}{r^\beta}. 
\end{align*}
In order to precisely estimate the $L^p(\bbr,dx)$-norm, $p \in (1,+\infty)$, of the operator $\mathcal{H}^\beta_R$ and to apply \cite[Theorem $5.2.7$]{G08}, let us estimate the supremum norm of the multiplier $h_R^\beta$ and of its first derivative. Then, for all $\xi \in \mathbb{R} \setminus \{0\}$, 
\begin{align*}
\left| h_R^\beta\left(\xi\right) \right| = 2 |\xi|^{\beta-1}  \left| \int_{R|\xi|}^{+\infty} \sin\left(r\right) \frac{dr}{r^\beta} \right|. 
\end{align*}
By an integration by parts, for all $\xi>0$, 
\begin{align*}
\int_{R|\xi|}^{+\infty} \sin\left(r\right) \frac{dr}{r^\beta} & = \frac{\cos\left(R|\xi|\right)}{(R|\xi|)^\beta} - \beta \int_{R|\xi|}^{+\infty} \cos(r)\frac{dr}{r^{1+\beta}}. 
\end{align*}
Thus, for all $\xi R \geq 1$, 
\begin{align*}
\left| h_R^\beta(\xi) \right|  = 2 |\xi|^{\beta-1} \left| \frac{\cos\left(R|\xi|\right)}{(R|\xi|)^\beta} - \beta \int_{R|\xi|}^{+\infty} \cos(r)\frac{dr}{r^{1+\beta}} \right| & \leq 2 |\xi|^{\beta-1} \left( \dfrac{1}{(R|\xi|)^\beta} + \beta \int_{R|\xi|}^{+\infty} \frac{dr}{r^{1+\beta}}\right) \\
& \leq \frac{4}{R^{\beta-1}} |R\xi|^{\beta-1} \dfrac{1}{(R|\xi|)^\beta} \leq \frac{4}{R^{\beta-1}}. 
\end{align*}
Finally, for all $\xi R \leq 1$, 
\begin{align*}
\left| h_R^{\beta}\left(\xi\right) \right| & = 2 |\xi|^{\beta-1} \left| \int_{R|\xi|}^1 \sin(r) \frac{dr}{r^\beta} + \int_{1}^{+\infty} \sin(r) \frac{dr}{r^\beta} \right| \leq 2 |\xi|^{\beta-1} \left(\int_{R |\xi|}^1 \frac{dr}{r^{\beta-1}} + 2 \right) \\
& \leq 2 |\xi|^{\beta-1} \left(\frac{1}{2-\beta}+ 2 \right) \leq \frac{2}{R^{\beta-1}}\left(\frac{1}{2-\beta}+ 2 \right). 
\end{align*}
Then, for all $R>0$, 
\begin{align}\label{ineq:estimate_supremum_norm}
\| h_R^{\beta} \|_{\infty , \bbr} := \sup_{\xi \in \mathbb{R}} \left|h_R^{\beta}(\xi)\right| \leq \frac{1}{R^{\beta-1}} \left(4+ \frac{2}{2-\beta}\right). 
\end{align}
It remains to estimate the first order derivative of the multiplier $h_R^\beta$. Without loss of generality, let us perform this analysis on $(0,+\infty)$.~Then, by standard computations, for all $\xi>0$, 
\begin{align*}
\frac{d}{d\xi} \left(h_R^{\beta}(\xi)\right) & = 2i \left(\beta -1\right) \xi^{\beta-2}\int_{R\xi}^{+\infty} \sin\left(r\right) \frac{dr}{r^\beta} - 2i \xi^{\beta-1} \dfrac{R\sin(R\xi)}{(R\xi)^\beta}.
\end{align*}
Thus, for all $\xi>0$, 
\begin{align*}
\left| \frac{d}{d\xi} \left(h_R^{\beta}(\xi)\right) \right| & \leq \frac{2(\beta-1)}{|\xi|} |\xi|^{\beta-1}  \left| \int_{R\xi}^{+\infty} \sin\left(r\right) \frac{dr}{r^\beta} \right| + \frac{2}{|\xi| R^{\beta-1} } \\
& \leq \frac{(\beta-1)}{|\xi|} \left(4+ \frac{2}{2-\beta}\right) \frac{1}{R^{\beta-1}} + \frac{2}{|\xi| R^{\beta-1} }.
\end{align*}
Then, \cite[Theorem $5.2.7$, (a)]{G08} in dimension $1$ ensures that the linear operator $\mathcal{H}^\beta_R$ is bounded from $L^p(\bbr,dx)$ to $L^p(\bbr,dx)$, with $p \in (1, +\infty)$, and, for all $f \in \mathcal{C}^{\infty}_c(\bbr)$, 
\begin{align*}
\left\| \mathcal{H}^\beta_R(f) \right\|_{L^p(\bbr,dx)} \leq C \max\left(p, (p-1)^{-1}\right) \left( 2 + \beta\left(4 +  \frac{2}{2-\beta} \right)\right) \frac{\|f\|_{L^p(\bbr,dx)}}{R^{\beta-1}}, 
\end{align*}
for some $C>0$. Then, by Fubini's theorem, for all $f \in \mathcal{C}^{\infty}_c(\bbr^d)$, 
\begin{align*}
\left\| H^{e_1}_{\beta,R}(f) \right\|_{L^p(\bbr^d,dx)} \leq C \max\left(p, (p-1)^{-1}\right) \left( 2 + \beta\left(4 +  \frac{2}{2-\beta} \right)\right) \frac{\|f\|_{L^p(\bbr^d,dx)}}{R^{\beta-1}},
\end{align*}  
which implies the boundedness of the operator $H^{y}_{\beta,R}$, for all $y \in \mathbb{S}^{d-1}$, from $L^p(\bbr^d,dx)$ to $L^p(\bbr^d,dx)$, with $p \in (1, +\infty)$, using the representation \eqref{eq:reduction_e1}. Thus, by Minkowski's integral inequality, 
for all $f \in \mathcal{C}^{\infty}_c(\bbr^d)$ and all $R>0$, 
\begin{align*}
\left\| \int_{B(0,R)^c} (f(\cdot+u)-f(\cdot)) u \nu_\beta^{\operatorname{rot}}(du) \right\|_{L^p(\bbr^d,dx)} & \leq \frac{c_{\beta,d}}{2} \sigma_L \left(\mathbb{S}^{d-1}\right) C \max\left(p, (p-1)^{-1}\right)  \\
& \quad\quad \times \left( 2 + \beta\left(4 +  \frac{2}{2-\beta} \right)\right) \frac{\|f\|_{L^p(\bbr^d,dx)}}{R^{\beta-1}}. 
\end{align*}
The end of the proof is similar to the one of Proposition \ref{prop:fractional_interpolation_inequality}. 
\end{proof}
\noindent
When the underlying symmetric non-degenerate $\alpha$-stable L\'evy measure is $\nu_{\alpha,d}$, the previous result becomes: 

\begin{prop}\label{prop:fractional_interpolation_independent}
Let $d \geq 2$ be an integer, let $1 < \beta < \alpha< 2$ and let $D^{\alpha-1,d}$ be the fractional gradient operator associated with the $\alpha$-stable probability measure $\mu_{\alpha,d}$ characterized by \eqref{eq:StableIndAxes}.~Let $p \in (1, +\infty)$. Then, for all $f \in \mathcal{C}_c^{\infty}(\mathbb{R}^d)$,
\begin{align}\label{ineq:fractional_interpolation_independent}
\|D^{\beta-1,d}(f)\|_{L^p(\mathbb{R}^d,dx)} \leq C^{ind}_{\alpha,\beta,p,d} \|f\|^{\frac{\alpha-\beta}{\alpha-1}}_{L^p(\mathbb{R}^d,dx)} \| D^{\alpha-1,d}(f) \|^{\frac{\beta-1}{\alpha-1}}_{L^p(\mathbb{R}^d,dx)},
\end{align}
where,
\begin{align}\label{cste:fractional_interpolation_independent}
C^{ind}_{\alpha,\beta,p,d} & = \frac{\sigma_d \left(\mathbb{S}^{d-1}\right)}{2} C \max\left(p, (p-1)^{-1}\right)\left( 2 + \beta\left(4 +  \frac{2}{2-\beta} \right)\right)+ \frac{c_{\alpha,d}}{\alpha}\dfrac{\sigma_d\left(\mathbb{S}^{d-1}\right)}{\alpha-\beta},
\end{align}
where $C>0$ is a numerical constant, $\sigma_d$ is the spherical component of the L\'evy measure of $\mu_{\alpha,d}$ and
\begin{align*}
c_{\alpha,d} : =  \sup_{\omega \in \mathbb{S}^{d-1}} \int_{\mathbb{R}^d} \left\| k_{\alpha,d}(\omega+z) - k_{\alpha,d}(z)  \right\| dz,
\end{align*}
with $k_{\alpha,d}$ defined in Lemma \ref{lem:l_alpha_ball}. 
\end{prop}

\begin{proof}
The proof is similar to the proof of the previous proposition. So, for all $x \in \mathbb{R}^d$ and all $R>0$, 
\begin{align*}
D^{\beta-1,d}(f)(x) = \int_{B(0,R)} (f(x+u) - f(x)) u \nu_{\beta,d}(du) + \int_{B(0,R)^c} (f(x+u) - f(x)) u \nu_{\beta,d}(du). 
\end{align*}
From the polar decomposition of the L\'evy measure $\nu_{\beta,d}$, 
\begin{align*}
 \int_{B(0,R)^c} (f(x+u) - f(x)) u \nu_{\beta,d}(du) & = \frac{1}{2} \int_{\mathbb{S}^{d-1}} y \sigma_d(dy) \int_{R}^{+\infty} \left(f(x+ry)-f(x-ry)\right) \frac{dr}{r^\beta} \\
 & = \frac{1}{2} \int_{\mathbb{S}^{d-1}} y \sigma_d(dy) H^y_{\beta,R}(f)(x). 
\end{align*}
Thus, as in the proof of Proposition \ref{prop:fractional_interpolation_inequality2} and by Minkowski's integral inequality, for all $R>0$, 
\begin{align*}
\left\|  \int_{B(0,R)^c} (f(.+u) - f(.)) u \nu_{\beta,d}(du) \right\|_{L^p(\mathbb{R}^d,dx)} & \leq \frac{1}{2} \int_{\mathbb{S}^{d-1}} \sigma_d(dy) \|H^y_{\beta,R}(f)\|_{L^p(\mathbb{R}^d,dx)} \\
& \leq \frac{\sigma_d \left(\mathbb{S}^{d-1}\right)}{2} C \max\left(p, (p-1)^{-1}\right)  \\
& \quad\quad \times \left( 2 + \beta\left(4 +  \frac{2}{2-\beta} \right)\right) \frac{\|f\|_{L^p(\bbr^d,dx)}}{R^{\beta-1}}. 
\end{align*}
Next, thanks to \eqref{eq:formula_Stein_axes}, for all $x \in \mathbb{R}^d$ and all $R>0$, 
\begin{align*}
 \int_{B(0,R)} (f(x+u) - f(x)) u \nu_{\beta,d}(du) & = \frac{1}{\alpha} \int_{B(0,R)} \bigg(\int_{\mathbb{R}^d} \langle k_{\alpha,d}(x+u-y) \\
 & \quad\quad\quad - k_{\alpha,d}(x-y) ; D^{\alpha-1,d}(f)(y) \rangle dy \bigg) u \nu_{\beta,d}(du) \\
 & = \frac{1}{\alpha} \int_{B(0,R)} \bigg(\int_{\mathbb{R}^d} \langle k_{\alpha,d}(u-z) \\
 & \quad\quad\quad - k_{\alpha,d}(-z) ; D^{\alpha-1,d}(f)(x+z) \rangle dy \bigg) u \nu_{\beta,d}(du).
\end{align*}
Then, by Minkowski's integral inequality, 
\begin{align*}
\bigg\|  \int_{B(0,R)} (f(.+u) - f(.)) u \nu_{\beta,d}(du) \bigg\|_{L^p(\mathbb{R}^d,dx)} & \leq \frac{1}{\alpha} \int_{B(0,R)}  \left(\int_{\mathbb{R}^d} \left\| k_{\alpha,d}(u-z) - k_{\alpha,d}(-z)  \right\| dz \right) \|u\| \nu_{\beta,d}(du) \\
& \quad\quad \times \|D^{\alpha-1,d}(f)\|_{L^p(\mathbb{R}^d,dx)}. 
\end{align*}
Changing variables and using Lemma \ref{lem:l_alpha_ball}, for all $u \in B(0,R) \setminus \{0\}$, 
\begin{align*}
\int_{\mathbb{R}^d} \left\| k_{\alpha,d}(u-z) - k_{\alpha,d}(-z)  \right\| dz = \|u\|^{\alpha-1}\int_{\mathbb{R}^d} \left\| k_{\alpha,d}(e_u-z) - k_{\alpha,d}(-z)  \right\| dz, \quad e_u = \frac{u}{\|u\|}. 
\end{align*}
Moreover, from Lemma \ref{lem:sharp_pointwise_bound}, for all $d \geq 2$, 
\begin{align*}
c_{\alpha,d} = \sup_{\omega \in \mathbb{S}^{d-1}} \int_{\mathbb{R}^d} \left\| k_{\alpha,d}(\omega+z) - k_{\alpha,d}(z)  \right\| dz <+\infty. 
\end{align*}
Thus, 
\begin{align*}
\bigg\|  \int_{B(0,R)} (f(.+u) - f(.)) u \nu_{\beta,d}(du) \bigg\|_{L^p(\mathbb{R}^d,dx)} & \leq \frac{c_{\alpha,d}}{\alpha}\sigma_d\left(\mathbb{S}^{d-1}\right) \dfrac{R^{\alpha-\beta}}{\alpha-\beta} \|D^{\alpha-1,d}(f)\|_{L^p(\mathbb{R}^d,dx)}. 
\end{align*}
The end of the proof follows by taking 
\begin{align*}
R = \left(\dfrac{\|f\|_{L^p(\mathbb{R}^d,dx)}}{\|D^{\alpha-1,d}(f)\|_{L^p(\mathbb{R}^d,dx)}}\right)^{\frac{1}{\alpha-1}}. 
\end{align*}
\end{proof}
\noindent

\begin{prop}\label{prop:fractional_interpolation_NDSstrong}
Let $d \geq 2$ be an integer, let $1 < \beta < \alpha< 2$ and let $\gamma \in [1,d]$ be such that $\gamma > d+2-2\alpha$. Let $D^{\alpha-1}$ be the fractional gradient operator associated with a symmetric non-degenerate $\alpha$-stable probability measure $\mu_{\alpha}$ with L\'evy measure $\nu_\alpha$ which is assumed to be a $\gamma$-measure in the sense of \eqref{eq:gamma-measure}.~Let $p \in (1, +\infty)$. Then, for all $f \in \mathcal{C}_c^{\infty}(\mathbb{R}^d)$, 
\begin{align}\label{ineq:fractional_interpolation_NDSstrong}
\|D^{\beta-1}(f)\|_{L^p(\mathbb{R}^d,dx)} \leq C^{\gamma}_{\alpha,\beta,p,d} \|f\|^{\frac{\alpha-\beta}{\alpha-1}}_{L^p(\mathbb{R}^d,dx)} \| D^{\alpha-1}(f) \|^{\frac{\beta-1}{\alpha-1}}_{L^p(\mathbb{R}^d,dx)},
\end{align}
where
\begin{align}\label{cste:fractional_interpolation_NDSstrong}
C^{\gamma}_{\alpha,\beta,p,d} & = \frac{\sigma \left(\mathbb{S}^{d-1}\right)}{2} C \max\left(p, (p-1)^{-1}\right)\left( 2 + \beta\left(4 +  \frac{2}{2-\beta} \right)\right)+ \frac{c_{\alpha,d,\gamma}}{\alpha}\dfrac{\sigma\left(\mathbb{S}^{d-1}\right)}{\alpha-\beta},
\end{align}
where $C>0$ is a numerical constant, $\sigma$ is the spherical component of the L\'evy measure $\nu_\alpha$ and
\begin{align*}
c_{\alpha,d,\gamma} : =  \sup_{\omega \in \mathbb{S}^{d-1}} \int_{\mathbb{R}^d} \left\| k_{\alpha}(\omega+z) - k_{\alpha}(z)  \right\| dz,
\end{align*}
with $k_{\alpha}$ defined in Fourier by \eqref{eq:generalized_kernel}. 
\end{prop}

\begin{proof}
The proof is similar to the ones of Proposition \ref{prop:fractional_interpolation_inequality2} and Proposition \ref{prop:fractional_interpolation_independent}. So, for all $x \in \mathbb{R}^d$ and all $R>0$, 
\begin{align*}
D^{\beta-1}(f)(x) = \int_{B(0,R)} (f(x+u) - f(x)) u \nu_{\beta}(du) + \int_{B(0,R)^c} (f(x+u) - f(x)) u \nu_{\beta}(du). 
\end{align*}
From the polar decomposition of the L\'evy measure $\nu_{\beta}$, 
\begin{align*}
 \int_{B(0,R)^c} (f(x+u) - f(x)) u \nu_{\beta}(du) & = \frac{1}{2} \int_{\mathbb{S}^{d-1}} y \sigma(dy) \int_{R}^{+\infty} \left(f(x+ry)-f(x-ry)\right) \frac{dr}{r^\beta} \\
 & = \frac{1}{2} \int_{\mathbb{S}^{d-1}} y \sigma(dy) H^y_{\beta,R}(f)(x). 
\end{align*}
Thus, as in the proof of Proposition \ref{prop:fractional_interpolation_inequality2} and by Minkowski's integral inequality, for all $R>0$, 
\begin{align*}
\left\|  \int_{B(0,R)^c} (f(.+u) - f(.)) u \nu_{\beta,d}(du) \right\|_{L^p(\mathbb{R}^d,dx)} & \leq \frac{1}{2} \int_{\mathbb{S}^{d-1}} \sigma(dy) \|H^y_{\beta,R}(f)\|_{L^p(\mathbb{R}^d,dx)} \\
& \leq \frac{\sigma \left(\mathbb{S}^{d-1}\right)}{2} C \max\left(p, (p-1)^{-1}\right)  \\
& \quad\quad \times \left( 2 + \beta\left(4 +  \frac{2}{2-\beta} \right)\right) \frac{\|f\|_{L^p(\bbr^d,dx)}}{R^{\beta-1}}. 
\end{align*}
Next, from Proposition \ref{prop:FFTC_L1_NDSstrong}, for all $x \in \mathbb{R}^d$ and all $R>0$, 
\begin{align*}
 \int_{B(0,R)} (f(x+u) - f(x)) u \nu_{\beta}(du) & = \frac{1}{\alpha} \int_{B(0,R)} \bigg(\int_{\mathbb{R}^d} \langle k_{\alpha}(x+u-y) \\
 & \quad\quad - k_{\alpha}(x-y) ; D^{\alpha-1}(f)(y) \rangle dy \bigg) u \nu_{\beta}(du) \\
 & = \frac{1}{\alpha} \int_{B(0,R)} \bigg(\int_{\mathbb{R}^d} \langle k_{\alpha}(u-z) \\
 & \quad\quad - k_{\alpha}(-z) ; D^{\alpha-1}(f)(x+z) \rangle dy \bigg) u \nu_{\beta}(du).
\end{align*}
Then, 
\begin{align*}
\bigg\|  \int_{B(0,R)} (f(.+u) - f(.)) u \nu_{\beta}(du) \bigg\|_{L^p(\mathbb{R}^d,dx)} & \leq \frac{1}{\alpha} \int_{B(0,R)}  \left(\int_{\mathbb{R}^d} \left\| k_{\alpha}(u-z) - k_{\alpha}(-z)  \right\| dz \right) \|u\| \nu_{\beta}(du) \\
& \quad\quad \times \|D^{\alpha-1}(f)\|_{L^p(\mathbb{R}^d,dx)}. 
\end{align*}
Changing variables and by homogeneity of the potential kernel $V_\alpha$, for all $u \in B(0,R) \setminus \{0\}$, 
\begin{align*}
\int_{\mathbb{R}^d} \left\| k_{\alpha}(u-z) - k_{\alpha}(-z)  \right\| dz = \|u\|^{\alpha-1}\int_{\mathbb{R}^d} \left\| k_{\alpha}(e_u-z) - k_{\alpha}(-z)  \right\| dz, \quad e_u = \frac{u}{\|u\|}. 
\end{align*}
Moreover, from Proposition \ref{prop:FFTC_L1_NDSstrong}, for all $d \geq 2$, 
\begin{align*}
c_{\alpha,d,\gamma} = \sup_{\omega \in \mathbb{S}^{d-1}} \int_{\mathbb{R}^d} \left\| k_{\alpha}(\omega+z) - k_{\alpha}(z)  \right\| dz <+\infty. 
\end{align*}
Thus, 
\begin{align*}
\bigg\|  \int_{B(0,R)} (f(.+u) - f(.)) u \nu_{\beta}(du) \bigg\|_{L^p(\mathbb{R}^d,dx)} & \leq \frac{c_{\alpha,d,\gamma}}{\alpha}\sigma\left(\mathbb{S}^{d-1}\right) \dfrac{R^{\alpha-\beta}}{\alpha-\beta} \|D^{\alpha-1}(f)\|_{L^p(\mathbb{R}^d,dx)}. 
\end{align*}
The end of the proof follows by taking 
\begin{align*}
R = \left(\dfrac{\|f\|_{L^p(\mathbb{R}^d,dx)}}{\|D^{\alpha-1}(f)\|_{L^p(\mathbb{R}^d,dx)}}\right)^{\frac{1}{\alpha-1}}. 
\end{align*}
\end{proof}

\begin{cor}\label{cor:non-smooth_mutliplier_NDSstrong}
Let $d \geq 2$ be an integer, let $1 < \beta < \alpha< 2$ and let $\gamma \in [1,d]$ be such that $\gamma > d+2-2\alpha$.~Let $\mu_\alpha$ be a symmetric non-degenerate $\alpha$-stable probability measure on $\mathbb{R}^d$ with L\'evy measure $\nu_\alpha$ which is assumed to be a $\gamma$-measure and L\'evy-Khintchine exponent $\psi_\alpha$. Let $T_{\alpha,\beta}$ be the linear multiplier operator defined by \eqref{eq:Fourier_multiplier_ratio}. Then, for all $p \in (1,+\infty)$, $T_{\alpha,\beta}$ is $L^p(\mathbb{R}^d,dx)$-bounded.   
\end{cor}

\begin{proof}
Combining Inequality \eqref{ineq:fractional_interpolation_NDSstrong} and Proposition \ref{prop:converse_inclusion} leads to 
\begin{align*}
W^{\alpha-1, p}\left(\mathbb{R}^d , dx\right) \hookrightarrow W^{\beta-1,p}\left(\mathbb{R}^d,dx\right) \Leftrightarrow H_p^{\psi_\alpha, r(\alpha)}(\mathbb{R}^d) \hookrightarrow H_p^{\psi_\beta, r(\beta)}(\mathbb{R}^d).
\end{align*} 
\cite[Theorem $3.3.28$]{NJ02_2} concludes the proof of the corollary. 
\end{proof}
\noindent 
In the next theorem, let us prove a generalization of \cite[Theorem $17$]{BCGS_CRM22} for the general symmetric non-degenerate $\alpha$-stable case, with $\alpha \in (1,2)$, when the associated L\'evy measure $\nu_\alpha$ is a $\gamma$-measure on $\mathbb{S}^{d-1}$ with $\gamma > 2 + d - 2\alpha$.  

\begin{thm}\label{thm:MS_strong_convergence_NDSstrong}
Let $d \geq 2$ be an integer, let $\alpha \in (1,2)$ and let $\gamma \in [1,d]$ be such that $\gamma > 2 + d - 2\alpha$. Let $\mu_\alpha$ be a non-degenerate symmetric $\alpha$-stable probability measure on $\mathbb{R}^d$ with L\'evy measure $\nu_\alpha$ which is assumed to be a $\gamma$-measure with spherical component $\sigma$. Let $D^{\alpha-1}$ be the fractional gradient operator associated with $\nu_\alpha$ and let $R_\sigma$ be the linear operator defined by \eqref{eq:avRT}. Then, for all $p \in (1, +\infty)$ and all $f \in \cup_{\alpha \in (1,2)} W^{\alpha-1,p}(\mathbb{R}^d,dx)$, 
\begin{align}\label{eq:MS_strong_convergence_NDSstrong}
\underset{\alpha \rightarrow 1^+}{\lim} \|D^{\alpha-1}(f)-R_\sigma(f)\|_{L^p(\mathbb{R}^d,dx)} =0. 
\end{align}
\end{thm}

\begin{proof}
Let $\alpha \in (1,2)$, let $p \in (1, +\infty)$ and let $f \in W^{\alpha-1,p}(\mathbb{R}^d,dx)$. By Proposition \ref{prop:density_compactsupport}, there exists a sequence of Schwartz functions $(f_k)_{k \geq 1}$ which converges to $f$ in $W^{\alpha-1,p}(\mathbb{R}^d,dx)$. Let $\beta \in (1, \alpha)$. Then, combining Lemma \ref{lem:Lp_continuity_avRT} and Proposition \ref{prop:fractional_interpolation_NDSstrong}, for all integer $k \geq 1$, 
\begin{align*}
 \|D^{\beta-1}(f)-R_\sigma(f)\|_{L^p(\mathbb{R}^d,dx)} & \leq \| D^{\beta-1}(f_k)-R_\sigma(f_k) \|_{L^p(\mathbb{R}^d,dx)} + \| D^{\beta-1}(f_k) - D^{\beta-1}(f) \|_{L^p(\mathbb{R}^d,dx)} \\
 &\quad\quad + \|R_\sigma(f_k) - R_\sigma(f)\|_{L^p(\mathbb{R}^d,dx)} \\
 & \leq \| D^{\beta-1}(f_k)-R_\sigma(f_k) \|_{L^p(\mathbb{R}^d,dx)} + C_{p,\sigma,d} \|f-f_k\|_{L^p(\mathbb{R}^d,dx)} \\
 & \quad\quad+ C^{\gamma}_{\alpha,\beta,p,d} \|f-f_k\|^{\frac{\alpha-\beta}{\alpha-1}}_{L^p(\mathbb{R}^d,dx)} \| D^{\alpha-1}(f_k-f) \|^{\frac{\beta-1}{\alpha-1}}_{L^p(\mathbb{R}^d,dx)},
\end{align*} 
where $C_{p,\sigma,d}>0$ is given by \eqref{eq:upper_bound_avRT}. Now, 
\begin{align}\label{eq:lim_constante_interpolation_NDSstrong}
C^{\gamma}_{\alpha,p,d} : = \underset{\beta \rightarrow 1^+}{\lim} C^{\gamma}_{\alpha,\beta,p,d} = 4\sigma \left(\mathbb{S}^{d-1}\right) C \max\left(p, (p-1)^{-1}\right)+ \frac{c_{\alpha,d,\gamma}}{\alpha}\dfrac{\sigma\left(\mathbb{S}^{d-1}\right)}{\alpha-1}.
\end{align} 
Thus, by the limit \eqref{eq:strong_convergence_MS_Schwartz}, 
\begin{align*}
\underset{\beta \rightarrow 1^+}{\limsup} \|D^{\beta-1}(f)-R_\sigma(f)\|_{L^p(\mathbb{R}^d,dx)} & \leq \underset{\beta \rightarrow 1^+}{\limsup} \| D^{\beta-1}(f_k)-R_\sigma(f_k) \|_{L^p(\mathbb{R}^d,dx)} + C_{p,\sigma,d} \|f-f_k\|_{L^p(\mathbb{R}^d,dx)} \\
&\quad\quad + C^{\gamma}_{\alpha,p,d} \|f-f_k\|_{L^p(\mathbb{R}^d,dx)} \\
& \leq C_{p,\sigma,d} \|f-f_k\|_{L^p(\mathbb{R}^d,dx)} + C^{\gamma}_{\alpha,p,d} \|f-f_k\|_{L^p(\mathbb{R}^d,dx)}. 
\end{align*}
Taking the limit as $k$ tends to $+\infty$ concludes the proof of the theorem. 
\end{proof}
\noindent
In the sequel, let us study strong convergence results in the spirit of the Bourgain-Br\'ezis-Mironescu asymptotics (see, e.g., \cite{BBM_2001,Lu_14_2,comi_stefani_rmc22}), for all $p \in [1, +\infty)$. 

\begin{prop}\label{prop:BBM_strong_convergence_NDS_densesubspace}
Let $d \geq 2$ be an integer and let $\alpha \in (1,2)$. Let $\mu_\alpha$ be a non-degenerate symmetric $\alpha$-stable probability measure on $\mathbb{R}^d$ with L\'evy measure $\nu_\alpha$ and spherical component $\sigma$. Let $D^{\alpha-1}$ be the associated fractional gradient operator and let $D_\sigma$ be the anisotropic local gradient operator given by \eqref{eq:anisotropic_local_gradient}. Then, for all $p \in [1, +\infty)$ and all $f \in \mathcal{C}^{\infty}_c(\mathbb{R}^d)$, 
\begin{align}\label{eq:BBM_strong_convergence_NDS_densesubspace}
\underset{\alpha \rightarrow 2^-}{\lim} \| (2-\alpha) D^{\alpha-1}(f) - D_\sigma(f)\|_{L^p(\mathbb{R}^d,dx)} = 0. 
\end{align}
\end{prop}

\begin{proof}
From \eqref{eq:superposition_non-singular_int_operator} and \eqref{eq:anisotropic_local_gradient}, for all $f \in \mathcal{C}^{\infty}_c(\mathbb{R}^d)$ and all $x \in \mathbb{R}^d$, 
\begin{align*}
(2-\alpha)D^{\alpha-1}(f)(x) = \frac{2-\alpha}{2}\int_{\mathbb{S}^{d-1}} y \left(\int_{(0,+\infty)} (f(x+ry)- f(x-ry))\frac{dr}{r^\alpha}\right) \sigma(dy),
\end{align*}
and 
\begin{align*}
D_\sigma(f)(x) = \int_{\mathbb{S}^{d-1}}y \langle y ; \nabla(f)(x)\rangle \sigma(dy). 
\end{align*}
Thus, for all $f \in \mathcal{C}_c^{\infty}(\mathbb{R}^d)$ and all $x \in \mathbb{R}^d$, 
\begin{align*}
(2-\alpha)D^{\alpha-1}(f)(x) - D_\sigma(f)(x)  & = \int_{\mathbb{S}^{d-1}} y \bigg[  \frac{2-\alpha}{2} \int_{(0,+\infty)} (f(x+ry)- f(x-ry))\frac{dr}{r^\alpha} \\
& \quad\quad\quad -\langle y ; \nabla(f)(x)\rangle\bigg] \sigma(dy).
\end{align*}
Let $y \in \mathbb{S}^{d-1}$ be fixed. Then,
\begin{align*}
\frac{2-\alpha}{2} \int_{(0,+\infty)} & (f(x+ry)- f(x-ry))\frac{dr}{r^\alpha} -\langle y ; \nabla(f)(x)\rangle = (2-\alpha) \int_0^1 \frac{dr}{r^\alpha} \bigg[ \dfrac{f(x+ry)-f(x-ry)}{2} \\
& \quad\quad - \langle ry ; \nabla(f)(x) \rangle \bigg] + \frac{2-\alpha}{2} \int_{1}^{+\infty} \frac{dr}{r^\alpha}\left( f(x+ry)-f(x-ry) \right). 
\end{align*}
For the second term of the sum appearing on the right-hand side of the previous equality and for all $p \in [1, +\infty)$, 
\begin{align*}
\left\| \frac{2-\alpha}{2} \int_{1}^{+\infty} \left( f(\cdot+ry)-f(\cdot-ry) \right) \frac{dr}{r^\alpha} \right\|_{L^p(\mathbb{R}^d,dx)} & \leq (2-\alpha) \|f\|_{L^p(\mathbb{R}^d,dx)} \int_{1}^{+\infty} \frac{dr}{r^\alpha} \\
& \leq \dfrac{2-\alpha}{\alpha-1} \|f\|_{L^p(\mathbb{R}^d,dx)} \underset{\alpha \rightarrow 2^-}{\longrightarrow} 0.
\end{align*}
Now, by an integration by parts, 
\begin{align*}
(2-\alpha) \int_0^1 \frac{dr}{r^\alpha}  \bigg[ \dfrac{f(x+ry)-f(x-ry)}{2}  - \langle ry ; \nabla(f)(x) \rangle \bigg]  & = \frac{2-\alpha}{1-\alpha} \int_{0}^1 \frac{d}{dr} \left(r^{1-\alpha}\right)  \\
& \quad\quad \times \bigg[ \dfrac{f(x+ry)-f(x-ry)}{2} \\
& \quad\quad - \langle ry ; \nabla(f)(x) \rangle \bigg] dr  \\
& =  \frac{2-\alpha}{2(\alpha-1)} \int_0^1 r^{1-\alpha} \bigg\langle y ; \nabla(f)(x+ry)  \\
& \quad\quad + \nabla(f)(x-ry) - 2\nabla(f)(x) \bigg\rangle dr \\
& \quad\quad + \frac{2-\alpha}{2(1-\alpha)} \bigg[ f(x+y)-f(x-y) \\
&\quad\quad - 2 \langle y ; \nabla(f)(x)\rangle \bigg].
\end{align*} 
Note that the second term of the sum of the previous equality vanishes when $\alpha$ tends to $2$. By another integration by parts, 
\begin{align*}
 & \frac{2-\alpha}{2(\alpha-1)} \int_0^1 r^{1-\alpha} \langle y ; \nabla(f)(x+ry)+ \nabla(f)(x-ry) - 2 \nabla(f)(x) \rangle dr \\
 & = \frac{1}{2(\alpha-1)} \int_0^1 \frac{d}{dr} \left(r^{2-\alpha}\right) \langle y ; \nabla(f)(x+ry)+ \nabla(f)(x-ry) - 2 \nabla(f)(x) \rangle dr \\
 & = - \frac{1}{2(\alpha-1)} \int_0^1 r^{2-\alpha} \bigg\langle y ; \frac{d}{dr} \bigg( \nabla(f)(x+ry)+ \nabla(f)(x-ry) \bigg) \bigg\rangle dr \\
 &\quad\quad\quad + \frac{1}{2(\alpha-1)} \langle y ; \nabla(f)(x+y) + \nabla(f)(x-y) - 2\nabla(f)(x) \rangle. 
\end{align*}
Thus,
\begin{align}\label{eq:smart_decomposition_Ball}
(2-\alpha) & \int_0^1  \bigg[ \dfrac{f(x+ry)-f(x-ry)}{2} - \langle ry ; \nabla(f)(x) \rangle \bigg] \frac{dr}{r^\alpha} = \frac{2-\alpha}{2(1-\alpha)} \bigg[ f(x+y)-f(x-y) - 2 \langle y ; \nabla(f)(x)\rangle \bigg] \nonumber  \\
& + \frac{1}{2(\alpha-1)} \langle y ; \nabla(f)(x+y) + \nabla(f)(x-y) - 2\nabla(f)(x) \rangle \nonumber \\
&  - \frac{1}{2(\alpha-1)} \int_0^1 r^{2-\alpha} \bigg\langle y ; \frac{d}{dr} \bigg( \nabla(f)(x+ry)+ \nabla(f)(x-ry) \bigg) \bigg\rangle dr.
\end{align}
From the previous decomposition, 
\begin{align*}
(2-\alpha) \int_0^1  \bigg[ \dfrac{f(x+ry)-f(x-ry)}{2} - \langle ry ; \nabla(f)(x) \rangle \bigg] \frac{dr}{r^\alpha} \longrightarrow 0,
\end{align*}
as $\alpha$ tends to $2$. Next, let $\alpha_0 \in (1,2)$. Then, for all $\alpha \in (\alpha_0,2)$, 
\begin{align*}
\bigg| (2-\alpha) \int_0^1  \bigg[ \dfrac{f(x+ry)-f(x-ry)}{2} & - \langle ry ; \nabla(f)(x) \rangle \bigg] \frac{dr}{r^\alpha} \bigg| \leq \underset{\alpha \in (\alpha_0,2)}{\sup}\left(\frac{2-\alpha}{2(\alpha-1)}\right) \\
& \quad \times \bigg| f(x+y)-f(x-y) - 2 \langle y ; \nabla(f)(x)\rangle \bigg| \\
& \quad + \frac{1}{2(\alpha_0-1)} \\
& \quad \times \| \nabla(f)(x+y) + \nabla(f)(x-y) - 2\nabla(f)(x) \| \\
& \quad + \frac{1}{2(\alpha_0-1)} \\
& \quad \times \int_0^1 \left\| \dfrac{d}{dr} (\nabla(f)(x+ry)+ \nabla(f)(x-ry))\right\| dr. 
\end{align*}
The right-hand sides of the previous inequality clearly belong to $L^p(\mathbb{R}^d,dx)$, for all $p \in [1, +\infty)$. An application of the Lebesgue dominated convergence theorem concludes the proof of the proposition. 
\end{proof}
\noindent
In the next theorem, the limit \eqref{eq:BBM_strong_convergence_NDS_densesubspace} is extended to the target functional space $W^{1,p}\left(\mathbb{R}^d,dx\right)$, $p \in [1, +\infty)$, as in \cite[Theorem $4.11$]{comi_stefani_rmc22}.

\begin{thm}\label{thm:prop:BBM_strong_convergence_NDS}
Let $\alpha \in (1,2)$. Let $\mu_\alpha$ be a non-degenerate symmetric $\alpha$-stable probability measure on $\mathbb{R}^d$, $d \geq 2$, with spherical component $\sigma$ and with L\'evy measure $\nu_\alpha$.~Let $D^{\alpha-1}$ be the associated fractional gradient operator and let $D_\sigma$ be the anisotropic local gradient operator. Then, for all $p \in [1, +\infty)$ and all $f \in W^{1,p}\left(\mathbb{R}^d,dx\right)$, 
\begin{align}\label{eq:BBM_strong_convergence_NDS}
\underset{\alpha \rightarrow 2^-}{\lim} \| (2-\alpha) D^{\alpha-1}(f) - D_\sigma(f)\|_{L^p(\mathbb{R}^d,dx)} = 0. 
\end{align}
\end{thm}

\begin{proof}
Let $p \in [1, +\infty)$ and let $f \in W^{1,p}\left(\mathbb{R}^d,dx\right)$. Let $(f_k)_{k \geq 1}$ be a sequence of functions in $\mathcal{C}_c^{\infty} \left(\mathbb{R}^d\right)$ such that $(f_k)_{k \geq 1}$ converges to $f$ in $W^{1,p}\left(\mathbb{R}^d,dx\right)$, i.e., $(f_k)_{k\geq 1}$ converges to $f$ in $L^p(\mathbb{R}^d,dx)$ and $(\nabla(f_k))_{k \geq 1}$ converges to $\nabla(f)$ in $L^p(\mathbb{R}^d,\mathbb{R}^d,dx)$. Moreover,
\begin{align*}
\left\| D_{\sigma}(f_k) - D_{\sigma}(f) \right\|_{L^p(\mathbb{R}^d,dx)} & = \left\| \int_{\mathbb{S}^{d-1}} y \langle \nabla(f_k - f)  ; y \rangle \sigma(dy) \right\|_{L^p(\mathbb{R}^d,dx)} \\
&  \leq \sigma\left(\mathbb{S}^{d-1}\right) \|\nabla(f_k)-\nabla(f)\|_{L^p(\mathbb{R}^d,dx)}. 
\end{align*} 
Thus, $(D_\sigma(f_k))_{k \geq 1}$ converges to $D_\sigma(f)$ in $L^p(\mathbb{R}^d,\mathbb{R}^d,dx)$. Furthermore, by the interpolation inequality \eqref{ineq:interpolation_inequality1}, for all $\alpha \in (1,2)$ and all $k \geq 1$, 
\begin{align*}
\left\| (2-\alpha) D^{\alpha-1}(f) - (2-\alpha) D^{\alpha}(f_k) \right\|_{L^p(\mathbb{R}^d,dx)} &\leq \left( \dfrac{2(2-\alpha) \sigma\left(\mathbb{S}^{d-1}\right)}{\alpha-1} + \sigma\left(\mathbb{S}^{d-1}\right) \right) \|f-f_k\|^{2-\alpha}_{L^p(\mathbb{R}^d,dx)} \\
&\quad\quad\quad  \times \| \nabla(f)-\nabla(f_k)\|^{\alpha-1}_{L^p(\mathbb{R}^d,dx)}.  
\end{align*} 
Then, for all $\alpha \in (1,2)$ and all $k \geq 1$, 
\begin{align*}
 \| (2-\alpha) D^{\alpha-1}(f) - D_\sigma(f)\|_{L^p(\mathbb{R}^d,dx)} & \leq \left\| (2-\alpha) D^{\alpha-1}(f) - (2-\alpha) D^{\alpha}(f_k) \right\|_{L^p(\mathbb{R}^d,dx)} \\
 & \quad\quad + \| (2-\alpha) D^{\alpha-1}(f_k) - D_\sigma(f_k)\|_{L^p(\mathbb{R}^d,dx)} \\
 & \quad\quad + \left\| D_{\sigma}(f_k) - D_{\sigma}(f) \right\|_{L^p(\mathbb{R}^d,dx)} \\
 & \leq \left( \dfrac{2(2-\alpha) \sigma\left(\mathbb{S}^{d-1}\right)}{\alpha-1} + \sigma\left(\mathbb{S}^{d-1}\right) \right) \|f-f_k\|^{2-\alpha}_{L^p(\mathbb{R}^d,dx)} \\
&\quad\quad  \times \| \nabla(f)-\nabla(f_k)\|^{\alpha-1}_{L^p(\mathbb{R}^d,dx)} \\
& \quad\quad +  \| (2-\alpha) D^{\alpha-1}(f_k) - D_\sigma(f_k)\|_{L^p(\mathbb{R}^d,dx)} \\
& \quad\quad + \left\| D_{\sigma}(f_k) - D_{\sigma}(f) \right\|_{L^p(\mathbb{R}^d,dx)}. 
\end{align*}
Then, by Proposition \ref{prop:BBM_strong_convergence_NDS_densesubspace}, for all $k \geq 1$, 
\begin{align*}
\underset{\alpha \rightarrow 2^-}{\limsup} \| (2-\alpha) D^{\alpha-1}(f) - D_\sigma(f)\|_{L^p(\mathbb{R}^d,dx)} & \leq \sigma\left(\mathbb{S}^{d-1}\right) \| \nabla(f)-\nabla(f_k)\|_{L^p(\mathbb{R}^d,dx)}\\
& \quad\quad+\left\| D_{\sigma}(f_k) - D_{\sigma}(f) \right\|_{L^p(\mathbb{R}^d,dx)}. 
\end{align*}
Letting $k$ tend to $+\infty$ in the previous inequality concludes the proof of the theorem. 
\end{proof}
\noindent
Let us now study the fractional bounded variation spaces associated with the symmetric non-degenerate $\alpha$-stable probability measures with $\alpha \in (1,2)$.~For this purpose, let $| \cdot |_\alpha$ be the mapping defined, for all $f$ good enough, by
\begin{align}\label{eq:bounded_fractional_snorm_general}
|f|_\alpha : = \sup \left\{ \int_{\mathbb{R}^d} f(x) \operatorname{div}_\alpha(\varphi)(x)dx: \, \varphi \in \mathcal{C}_c^\infty(\mathbb{R}^d,\mathbb{R}^d), \|\varphi\|_{\infty} \leq 1 \right\}, 
\end{align}
where $\operatorname{div}_\alpha$ is given by Definition \ref{defi:fractional_divergence} and let $BV_\alpha(\mathbb{R}^d)$ be the functional space given by
\begin{align}\label{def:bounded_fractional_space_general}
BV_\alpha(\mathbb{R}^d) = \left\{ f \in L^1(\mathbb{R}^d, dx): \, \left| f\right|_\alpha<+\infty\right\}. 
\end{align}
Thanks to Step $3$ in the proof of Proposition \ref{prop:variational_representation_fractional_perimeter}, a structure theorem for $BV_\alpha(\mathbb{R}^d)$ holds true: $f \in BV_\alpha (\mathbb{R}^d)$ if and only if there exists a finite $\mathbb{R}^d$-valued Radon measure $d_{\alpha-1}(f)$ on $\mathcal{B}(\mathbb{R}^d)$ such that, for all $\varphi \in \mathcal{C}_c^{\infty}(\mathbb{R}^d, \mathbb{R}^d)$, 
\begin{align}\label{eq:frac_ipp_NDS}
\int_{\mathbb{R}^d} f(x) \operatorname{div}_\alpha(\varphi)(x)dx = - \int_{\mathbb{R}^d} \langle \varphi(x) ; d_{\alpha-1}(f)(dx) \rangle. 
\end{align}
Moreover, the semi-norm $|\cdot|_{\alpha}$ is lower semicontinuous:
\begin{align*}
\left| f\right|_{\alpha} \leq \underset{n \rightarrow +\infty}{\liminf} \left|f_n \right|_{\alpha} <+\infty,
\end{align*}
where $(f_n)_{n\geq 1}$ is a sequence of functions belonging to $BV_\alpha\left(\mathbb{R}^d\right)$ and which converges to $f$ in $L^1\left(\mathbb{R}^d,dx\right)$. Finally, Proposition \ref{prop:approximation_general_case_BV} of the Appendix shows that $\mathcal{C}_c^{\infty}\left(\mathbb{R}^d\right)$ is dense in $BV_\alpha(\mathbb{R}^d)$ in the following sense: for all $f \in BV_\alpha\left(\mathbb{R}^d\right)$ there exists a sequence of functions $(f_k)_{k \geq 1}$ in $\mathcal{C}_c^{\infty}\left(\mathbb{R}^d\right)$ such that $f_k \longrightarrow f$ in $L^1\left(\mathbb{R}^d,dx\right)$ and $|f_k|_{\alpha} \longrightarrow |f|_{\alpha}$, as $k$ tends to $+\infty$. Now, let us prove a compactness result in the fractional bounded variation space associated with the probability measure $\mu_{\alpha,d}$ given by \eqref{eq:StableIndAxes}, with $\alpha \in (1,2)$ and $d \geq 2$. Let $BV_{\alpha,d}(\mathbb{R}^d)$ be defined by 
\begin{align}\label{eq:fractional_bounded_variation_space_axes}
BV_{\alpha,d}(\mathbb{R}^d) = \left\{f \in L^1(\mathbb{R}^d,dx): \, \left| f\right|_{\alpha,d}<+\infty\right\},
\end{align}  
with, 
\begin{align}\label{eq:fractional_seminorm_axes}
 \left| f\right|_{\alpha,d} : =  \sup \left\{ \int_{\mathbb{R}^d} f(x) \operatorname{div}_{\alpha,d}(\varphi)(x)dx: \, \varphi \in \mathcal{C}_c^\infty(\mathbb{R}^d,\mathbb{R}^d), \|\varphi\|_{\infty} \leq 1 \right\},
\end{align}
and with, for all $\varphi \in \mathcal{C}_c^{\infty}(\mathbb{R}^d, \mathbb{R}^d)$ and all $x \in \mathbb{R}^d$, 
\begin{align*}
\operatorname{div}_{\alpha,d}(\varphi)(x) = \sum_{k = 1}^d D^{\alpha-1,d}_k(\varphi_k)(x). 
\end{align*}
The next technical lemma provides a quantitative error bound in $L^1$-norm for the translations of functions in $\mathcal{C}^{\infty}_c(\mathbb{R}^d)$ and is similar to \cite[Proposition $3.14$]{comi_stefani}. 

\begin{prop}\label{prop:L1_estimate_translation_axes}
Let $d \geq 2$ be an integer, let $\alpha \in (1,2)$ and let $D^{\alpha-1,d}$ be the fractional gradient operator associated with the $\alpha$-stable probability measure $\mu_{\alpha,d}$ given by \eqref{eq:StableIndAxes}. Then, for all $f \in \mathcal{C}_c^{\infty}(\mathbb{R}^d)$ and all $y \in \mathbb{R}^d$, 
\begin{align}\label{ineq:translationL1}
\int_{\mathbb{R}^d} \left| f(x+y) - f(x) \right| dx \leq c_{\alpha,d} \|y\|^{\alpha-1} \left\| D^{\alpha-1,d}(f) \right\|_{L^1(\mathbb{R}^d,dx)}, 
\end{align}
for some $c_{\alpha,d}>0$ depending only on $\alpha$ and $d$. 
\end{prop}

\begin{proof}
Let $f \in \mathcal{C}^{\infty}_c(\mathbb{R}^d)$ and let $y \in \mathbb{R}^d \setminus \{0\}$.~In the sequel, $c_{\alpha,d}$ denotes a positive constant which might change from a line to another. Then, by Lemma \ref{lem:l_alpha_ball}, for all $x \in \mathbb{R}^d$, 
\begin{align}\label{eq:fractional_fundamental_theorem}
f(x+y)- f(x) = \frac{1}{\alpha} \int_{\mathbb{R}^d} \langle k_{\alpha,d}(x+y-z)- k_{\alpha,d}(x-z) ; D^{\alpha-1,d}(f)(z) \rangle dz.  
\end{align}
Thus, by invariance by translations and homogeneity, 
\begin{align*}
\int_{\mathbb{R}^d} |f(x+y) - f(x)| dx & \leq \frac{1}{\alpha} \int_{\mathbb{R}^d} \int_{\mathbb{R}^d} \left\| k_{\alpha,d}(y-\omega)- k_{\alpha,d}(-\omega) \right\| \left\| D^{\alpha-1,d}(f)(\omega+x) \right\| d\omega dx  \\
& \leq \frac{1}{\alpha} \left\| D^{\alpha-1,d}(f) \right\|_{L^1(\mathbb{R}^d,dx)} \int_{\mathbb{R}^d} \left\| k_{\alpha,d}(y+\omega)- k_{\alpha,d}(\omega) \right\|  d\omega \\
& \leq \dfrac{\|y\|^{\alpha-1}}{\alpha} \left\| D^{\alpha-1,d}(f) \right\|_{L^1(\mathbb{R}^d,dx)} \int_{\mathbb{R}^d} \left\| k_{\alpha,d}(\omega+ e_y)- k_{\alpha,d}(\omega) \right\|  d\omega,
\end{align*}
with $e_y = y/ \|y\|$. \eqref{eq:L1_increments_axes_fulld} concludes the proof of the proposition. 
\end{proof}
\noindent
A straightforward consequence of the previous proposition is the following corollary which contains a quantitative $L^1$-error estimate in the space $BV_{\alpha,d}(\mathbb{R}^d)$. 

\begin{cor}\label{cor:L1error_BVaxes}
Let $d \geq 2$, let $\alpha \in (1,2)$ and let $D^{\alpha-1,d}$ be the fractional gradient operator associated with the $\alpha$-stable probability measure $\mu_{\alpha,d}$ given by \eqref{eq:StableIndAxes}.~Let $(\rho_\varepsilon)_{\varepsilon>0}$ be a sequence of standard mollifiers. Then, for all $f \in BV_{\alpha,d}(\mathbb{R}^d)$ and all $\varepsilon>0$,
\begin{align}\label{ineq:errorL1}
\|f \ast \rho_\varepsilon - f\|_{L^1(\mathbb{R}^d,dx)} \leq c_{\alpha,d} \varepsilon^{\alpha-1} \left| f \right|_{\alpha,d},
\end{align}
for some positive constant $c_{\alpha,d}$ which only depends on $\alpha$ and $d$. 
\end{cor}

\begin{proof}
The proof is similar to the proof of \cite[Corollary $3.15$]{comi_stefani} by combining the estimate \eqref{ineq:translationL1} together with the definition of the convolution operation and with the approximation result Proposition \ref{prop:approximation_general_case_BV}. 
\end{proof}
\noindent
We are now ready to state the compactness result in the space $BV_{\alpha,d}(\mathbb{R}^d)$.

\begin{thm}\label{thm:compactness_BVfrac_axes}
Let $d \geq 2$ and let $\alpha \in (1,2)$. Let $(f_k)_{k \geq 1}$ be a sequence of functions in $BV_{\alpha,d}(\mathbb{R}^d)$ such that 
\begin{align}\label{eq:uniformly_bounded}
\underset{k \geq 1}{\sup}\left( \|f_k\|_{L^1(\mathbb{R}^d,dx)}+ | f_k |_{\alpha,d} \right) <+\infty. 
\end{align}
Then, there exists a subsequence $(f_{n_k})_{k \geq 1}$ and a function $f \in L^1(\mathbb{R}^d,dx)$ such that
\begin{align}\label{eq:convergenceL1}
f_{n_k} \longrightarrow f \, \operatorname{in} \, L^1_{\operatorname{loc}}(\mathbb{R}^d,dx), 
\end{align}
as $k$ tends to $+\infty$.
\end{thm}

\begin{proof}
The proof is similar to the one of \cite[Theorem $3.16$]{comi_stefani} and is left to the interested reader. 
\end{proof}
\noindent

\begin{prop}\label{prop:L1_estimate_translation_anisotropic}
Let $d \geq 2$ be an integer, let $\alpha \in (1,2)$ and let $\gamma \in [1,d]$ be such that $\gamma > 2+d-2\alpha$. Let $D^{\alpha-1}$ be the fractional gradient operator associated with the symmetric non-degenerate $\alpha$-stable probability measure $\mu_{\alpha}$ with L\'evy measure $\nu_\alpha$ which is assumed to be a $\gamma$-measure and has spherical component $\sigma$.~Then, for all $f \in \mathcal{C}_c^{\infty}(\mathbb{R}^d)$ and all $y \in \mathbb{R}^d$, 
\begin{align}\label{ineq:L1_estimate_translation_anisotropic}
\int_{\mathbb{R}^d} \left| f(x+y) - f(x) \right| dx \leq c_{\alpha,d,\sigma} \|y\|^{\alpha-1} \left\| D^{\alpha-1}(f) \right\|_{L^1(\mathbb{R}^d,dx)}, 
\end{align}
for some $c_{\alpha,d,\sigma}>0$ depending only on $\alpha$, $d$ and $\sigma$.  
\end{prop}

\begin{proof}
With the help of Proposition \ref{prop:FFTC_L1_NDSstrong}, for all $f \in \mathcal{C}_c^{\infty}\left(\mathbb{R}^d\right)$ and all $x,y \in \mathbb{R}^d$, 
\begin{align*}
f(x+y) - f(x) = \frac{1}{\alpha} \int_{\mathbb{R}^d} \langle k_{\alpha}\left(x+y-z\right) - k_\alpha(x-z); D^{\alpha-1}(f)(z) \rangle dz.
\end{align*}
By a change of variables and the Cauchy-Schwarz inequality, 
\begin{align*}
\left| f(x+y) - f(x) \right| \leq \frac{1}{\alpha} \int_{\mathbb{R}^d} \left\| k_{\alpha}(y-z) - k_{\alpha}(-z) \right\| \|D^{\alpha-1}(f)(x+z)\| dz. 
\end{align*}
Then, by the Fubini theorem and the translations invariance, for all $y \in \bbr^d$, 
\begin{align*}
\int_{\mathbb{R}^d} \left| f(x+y)-f(x) \right| dx \leq \frac{1}{\alpha} \int_{\mathbb{R}^d} \|k_{\alpha}(y-z) - k_\alpha(-z)\| dz \|D^{\alpha-1}(f)\|_{L^1(\mathbb{R}^d,dx)}. 
\end{align*} 
By homogeneity of the potential kernel and Proposition \ref{prop:FFTC_L1_NDSstrong}, 
\begin{align*}
\int_{\mathbb{R}^d} \left| f(x+y)-f(x) \right| dx \leq c_{\alpha,d, \sigma} \|y\|^{\alpha-1} \|D^{\alpha-1}(f)\|_{L^1(\mathbb{R}^d,dx)}, 
\end{align*}
where $c_{\alpha,d, \sigma}$ is given by 
\begin{align*}
c_{\alpha,d, \sigma} : = \dfrac{1}{\alpha} \underset{\omega \in \mathbb{S}^{d-1}}{\sup} \left( \int_{\mathbb{R}^d} \|k_{\alpha}(\omega+z) - k_\alpha(z)\| dz \right) <+\infty. 
\end{align*}
This concludes the proof of the proposition. 
\end{proof}

\begin{cor}\label{cor:L1error_BVanisotropic}
Let $d \geq 2$ be an integer, let $\alpha \in (1,2)$ and let $\gamma \in [1,d]$ be such that $\gamma > 2+d-2\alpha$. Let $D^{\alpha-1}$ be the fractional gradient operator associated with the symmetric non-degenerate $\alpha$-stable probability measure $\mu_{\alpha}$ with L\'evy measure $\nu_\alpha$ which is assumed to be a $\gamma$-measure and has spherical component $\sigma$.~Let $(\rho_\varepsilon)_{\varepsilon>0}$ be a sequence of standard mollifiers.~Then, for all $f \in BV_{\alpha}(\mathbb{R}^d)$ and all $\varepsilon>0$,
\begin{align}\label{ineq:L1error_BVanisotropic}
\|f \ast \rho_\varepsilon - f\|_{L^1(\mathbb{R}^d,dx)} \leq c_{\alpha,d,\sigma} \varepsilon^{\alpha-1} \left| f \right|_{\alpha},
\end{align}
for some positive constant $c_{\alpha,d,\sigma}$ which depends on $\alpha$, $d$ and $\sigma$. 
\end{cor}

\begin{proof}
The proof is similar to the proof of \cite[Corollary $3.15$]{comi_stefani} by combining the estimate \eqref{ineq:L1_estimate_translation_anisotropic} together with the definition of the convolution operation and with the approximation result Proposition \ref{prop:approximation_general_case_BV} of the Appendix. 
\end{proof}

\begin{thm}\label{thm:compactness_BVanisotropic}
Let $d \geq 2$ be an integer, let $\alpha \in (1,2)$ and let $\gamma \in [1,d]$ be such that $\gamma > 2+d-2\alpha$. Let $(f_k)_{k \geq 1}$ be a sequence of functions in $BV_{\alpha}(\mathbb{R}^d)$ such that 
\begin{align}\label{eq:uniformly_bounded_BVanisotropic}
\underset{k \geq 1}{\sup}\left( \|f_k\|_{L^1(\mathbb{R}^d,dx)}+ | f_k |_{\alpha} \right) <+\infty. 
\end{align}
Then, there exists a subsequence $(f_{n_k})_{k \geq 1}$ and a function $f \in L^1(\mathbb{R}^d,dx)$ such that
\begin{align}\label{eq:convergenceL1_BVanisotropic}
f_{n_k} \longrightarrow f \, \operatorname{in} \, L^1_{\operatorname{loc}}(\mathbb{R}^d,dx), 
\end{align}
as $k$ tends to $+\infty$.
\end{thm}

\begin{proof}
The proof is similar to the one of \cite[Theorem $3.16$]{comi_stefani} and is left to the interested reader.
\end{proof}

\section{Sobolev embeddings}\label{sec:Sobolev_embeddings}
\noindent
In this section, we are mainly interested in the fractional Sobolev-type inequalities which bring into play the operator $D^{\alpha - 1}$ in the general non-degenerate symmetric $\alpha$-stable case.~Formally, for all $\alpha \in (1,2)$ and all $f \in \mathcal{S}(\mathbb{R}^d)$, 
\begin{align}\label{ineq:strong_Sobolev_Frac}
\| f \|_{L^{p_\alpha^*}(\mathbb{R}^d,dx)} \leq C \| D^{\alpha-1}(f) \|_{L^p(\mathbb{R}^d,dx)},
\end{align}
for some positive constant $C$ depending on $p,d$ and $\alpha$ and with
\begin{align}\label{eq:relation_frac_sobolev}
p_\alpha^* = \dfrac{pd}{d - p(\alpha-1)}, \quad p(\alpha-1)<d.
\end{align}
First, let us consider a Nash-type inequality involving the fractional gradient operator $D^{\alpha-1}$.~The proof is an adaptation of the Fourier argument contained in the classical paper \cite{Nash}.

\begin{prop}\label{prop:Nash_type_inequality}
Let $\alpha \in (1,2)$, let $\nu_\alpha$ be a non-degenerate symmetric $\alpha$-stable L\'evy measure on $\mathbb{R}^d \setminus \{0\}$, $d \geq 1$, verifying \eqref{eq:scale} with spectral measure $\lambda_1$ and let $D^{\alpha-1}$ be defined by \eqref{eq:perim_fracGrad}.~Then, for all $f \in \mathcal{S}(\mathbb{R}^d)$, 
\begin{align}\label{eq:Nash_type_Inequality}
\|f\|^{1+\frac{2(\alpha-1)}{d}}_{L^2(\mathbb{R}^d,dx)} \leq C_{\alpha,d} \| D^{\alpha-1}(f)\|_{L^2(\mathbb{R}^d,dx)} \|f\|^{\frac{2(\alpha-1)}{d}}_{L^1(\mathbb{R}^d,dx)},
\end{align}
for some $C_{\alpha,d}>0$ depending on $\alpha$, $d$ and the spectral measure $\lambda_1$. 
\end{prop}

\begin{proof}
First, let us compute the Fourier symbol of the linear operator $D^{\alpha-1}$.~By the Fourier inversion formula, for all $\xi \in \mathbb{R}^d$, 
\begin{align}\label{eq:fourier_frac_gradient}
\mathcal{F}(D^{\alpha-1}(f))(\xi) =\left( \int_{\mathbb{R}^d} u \left(e^{i \langle u ; \xi \rangle}-1\right) \nu_\alpha(du) \right) \mathcal{F}(f)(\xi). 
\end{align}
Now, for all $j \in \{1, \dots, d\}$ and all $\xi \in \mathbb{R}^d$, 
\begin{align*}
\int_{\mathbb{R}^d} u_j \left(e^{i \langle u ; \xi \rangle}-1\right) \nu_\alpha(du) =\frac{1}{i} \partial_j \int_{\mathbb{R}^d} \left( e^{i \langle u ; \xi \rangle}-1- i \langle u ; \xi \rangle\right) \nu_\alpha(du) = i \partial_j \left( \int_{\mathbb{S}^{d-1}} |\langle y ; \xi \rangle|^{\alpha} \lambda_1(dy)\right), 
\end{align*}
where $\lambda_1$ is the spectral measure associated with $\nu_\alpha$.~Then, by straightforward computations (since $\alpha \in (1,2)$), 
\begin{align*}
\int_{\mathbb{R}^d} u_j \left(e^{i \langle u ; \xi \rangle}-1\right) \nu_\alpha(du) = i \alpha \int_{\mathbb{S}^{d-1}} y_j |\langle y ; \xi \rangle|^{\alpha-1} \operatorname{sign}(\langle \xi ; y \rangle)\lambda_1(dy).
\end{align*}
In the sequel, let us denote by $\eta_\alpha$ the $\mathbb{R}^d$-valued function defined, for all $\xi \in \mathbb{R}^d$, by
\begin{align*}
\eta_\alpha(\xi) =  \int_{\mathbb{S}^{d-1}} y |\langle y ; \xi \rangle|^{\alpha-1} \operatorname{sign}(\langle \xi ; y \rangle)\lambda_1(dy).
\end{align*}
Take $R >0$ and let $\xi \in \mathbb{R}^d \setminus B(0,R)$.~Then,
\begin{align*}
\| \eta_\alpha(\xi) \| = \| \xi \|^{\alpha-1} \left\|  \int_{\mathbb{S}^{d-1}} y |\langle y ; e_\xi \rangle|^{\alpha-1} \operatorname{sign}(\langle e_\xi ; y \rangle)\lambda_1(dy) \right\|,
\end{align*}
where $e_\xi = \xi /\|\xi\| \in \mathbb{S}^{d-1}$. Now, by duality, 
\begin{align*}
\left\|  \int_{\mathbb{S}^{d-1}} y |\langle y ; e_\xi \rangle|^{\alpha-1} \operatorname{sign}(\langle e_\xi ; y \rangle)\lambda_1(dy) \right\| & = \underset{v \in \mathbb{S}^{d-1}}{\sup} \left|\langle v ;  \int_{\mathbb{S}^{d-1}} y |\langle y ; e_\xi \rangle|^{\alpha-1} \operatorname{sign}(\langle e_\xi ; y \rangle)\lambda_1(dy)  \rangle \right| \\
&\geq \left|\langle e_\xi ;  \int_{\mathbb{S}^{d-1}} y |\langle y ; e_\xi \rangle|^{\alpha-1} \operatorname{sign}(\langle e_\xi ; y \rangle)\lambda_1(dy)  \rangle \right| \\
&\geq  \int_{\mathbb{S}^{d-1}} |\langle y ; e_\xi \rangle|^{\alpha} \lambda_1(dy) \\
&\geq \underset{e \in \mathbb{S}^{d-1}}{\inf} \int_{\mathbb{S}^{d-1}} |\langle y ; e \rangle|^{\alpha} \lambda_1(dy).
\end{align*} 
Thus, 
\begin{align*}
\underset{e \in \mathbb{S}^{d-1}}{\inf}\left\|  \int_{\mathbb{S}^{d-1}} y |\langle y ; e \rangle|^{\alpha-1} \operatorname{sign}(\langle e ; y \rangle)\lambda_1(dy) \right\| &\geq \underset{e \in \mathbb{S}^{d-1}}{\inf} \int_{\mathbb{S}^{d-1}} |\langle y ; e \rangle|^{\alpha} \lambda_1(dy) >0,
\end{align*}
which implies that
\begin{align*}
\| \eta_\alpha(\xi) \| \geq  \| \xi \|^{\alpha-1} c(\alpha,d) \geq R^{\alpha-1} c(\alpha,d),
\end{align*}
for some $c(\alpha,d)>0$.~Next, by the Parseval-Plancherel formula, for all $f \in \mathcal{S}(\mathbb{R}^d)$, 
\begin{align*}
\|f\|^2_{L^2(\mathbb{R}^d,dx)} & = \frac{1}{(2\pi)^d} \int_{\mathbb{R}^d} |\mathcal{F}(f)(\xi)|^2 d\xi \\
& = \frac{1}{(2\pi)^d} \int_{B(0,R)} |\mathcal{F}(f)(\xi)|^2 d\xi +  \frac{1}{(2\pi)^d} \int_{\mathbb{R}^d \setminus B(0,R)} |\mathcal{F}(f)(\xi)|^2 d\xi.
\end{align*}  
Hence, 
\begin{align*}
 \int_{\mathbb{R}^d \setminus B(0,R)} |\mathcal{F}(f)(\xi)|^2 d\xi & \leq \frac{1}{R^{2(\alpha-1)} c(\alpha,d)^2} \int_{\mathbb{R}^d} \| \eta_\alpha(\xi) \|^2 |\mathcal{F}(f)(\xi)|^2 d\xi \\ 
 & \leq \frac{C_{1,\alpha,d}}{R^{2(\alpha-1)}} \int_{\mathbb{R}^d} |\mathcal{F}(D^{\alpha-1}(f))(\xi)|^2 d\xi,
\end{align*}
for some positive constant $C_{1,\alpha,d}$ depending on $\alpha$, $d$ and $\lambda_1$.~Moreover, 
\begin{align*}
\int_{B(0,R)} |\mathcal{F}(f)(\xi)|^2 d\xi \leq C_{2,d} R^d \|f\|^2_{L^1(\mathbb{R}^d,dx)}.
\end{align*}
The end of the proof follows by a standard optimization in $R$.
\end{proof}

\begin{rem}\label{rem:classical_Nash_equivalent_Sobolev}
In the classical case with the standard gradient, it is now well-known that the Nash inequality is equivalent to the (a priori stronger) Sobolev inequality: for all $d >2$ and all $f \in W^{1,2}(\mathbb{R}^d,dx)$, 
\begin{align}\label{ineq:Sobolev_classical}
\|f\|_{L^q(\mathbb{R}^d,dx)} \leq C_{d} \| \nabla(f) \|_{L^2(\mathbb{R}^d,dx)},
\end{align}
 for some $C_{d}>0$ and for $q = 2d/(d-2)$.~The fact that the inequality \eqref{ineq:Sobolev_classical} implies the Nash inequality can be seen by an application of H\"older's inequality.~The converse follows by an application of the truncation technique developed and generalized in \cite{BCLS}.
\end{rem}
\noindent 
Next, let us discuss the fractional Sobolev inequality \eqref{ineq:strong_Sobolev_Frac} in the rotationally invariant case.~Recall that the Fourier symbol of $D^{\alpha-1, \operatorname{rot}}$ defined by \eqref{eq:perim_fracGrad} with $\nu_\alpha = \nu_\alpha^{\operatorname{rot}}$ is given, for all $\xi \in \mathbb{R}^d \setminus \{0\}$, by 
\begin{align}\label{eq:Fourier_symbol_frac_grad_rotinv}
\tau_{\alpha}^{\operatorname{rot}}(\xi) =  \frac{\alpha}{2} \frac{i \xi}{\|\xi\|} \|\xi\|^{\alpha-1}.
\end{align}
Then, for all $f \in \mathcal{S}(\mathbb{R}^d)$ and all $x \in \mathbb{R}^d$, 
\begin{align}\label{eq:factorization_rot_invariant}
D^{\alpha-1, \operatorname{rot}}(f)(x) = \frac{\alpha}{2} \mathcal{R} \circ (-(- \Delta)^{\frac{\alpha-1}{2}})(f)(x) = \frac{\alpha}{2} (-(- \Delta)^{\frac{\alpha-1}{2}})(\mathcal{R}(f))(x),
\end{align}
where $(-(- \Delta)^{\frac{\alpha-1}{2}})$ denotes the fractional Laplacian of order $\alpha-1$ with Fourier symbol given by $-\| \xi \|^{\alpha-1}$, for all $\xi \in \mathbb{R}^d$. It is well-known that there is a fractional version of the Sobolev embedding theorem using the Gagliardo semi-norm and which involves the fractional Laplacian when $p=2$. By \cite[Theorem $6.5$]{NPV_12}, for all $s \in (0,1)$, all $p \in [1,d/s)$ and all $f \in \mathcal{C}_c^{\infty}(\mathbb{R}^d)$, 
\begin{align}\label{ineq:classical_fractional_Sobolev_embedding_theorem}
\|f\|_{L^q(\mathbb{R}^d,dx)} \leq C_{s,d,p} \left(\int_{\mathbb{R}^d} \int_{\mathbb{R}^d} \dfrac{|f(x) - f(y)|^p}{\|x- y\|^{d+ps}}dxdy\right)^{\frac{1}{p}},
\end{align}
for some $C_{s,d,p}>0$ depending only on $s$, $d$ and $p$ and for $q = pd/(d-ps)$.~Now, by Fourier arguments, for all $f \in \mathcal{C}_c^{\infty}(\mathbb{R}^d)$,
\begin{align*}
\int_{\mathbb{R}^d} \int_{\mathbb{R}^d} \dfrac{|f(x) - f(y)|^2}{\|x- y\|^{d+2s}}dxdy = \tilde{C}_{s,d} \int_{\mathbb{R}^d} \|\xi\|^{2s} \left| \mathcal{F}(f)(\xi)\right|^2 \frac{d \xi}{(2\pi)^d} =  \tilde{C}_{s,d} \| (-(- \Delta)^{\frac{s}{2}})(f)  \|^2_{L^2(\mathbb{R}^d,dx)},
\end{align*}
for some normalizing constant $\tilde{C}_{s,d} >0$ depending only on $s$ and $d$.~By \eqref{ineq:classical_fractional_Sobolev_embedding_theorem}, for all $d \geq 2$, all $\alpha \in (1,2)$ and all $f \in \mathcal{C}^{\infty}_c(\mathbb{R}^d)$, 
\begin{align*}
\| f \|_{L^{2_\alpha^*}(\mathbb{R}^d,dx)} \leq C_{\alpha-1,d,2}  (\tilde{C}_{\alpha-1,d})^{\frac{1}{2}} \| (-(- \Delta)^{\frac{\alpha-1}{2}})(f)  \|_{L^2(\mathbb{R}^d,dx)},
\end{align*}
with $2_\alpha^* = 2d/(d-2(\alpha-1))$.~Thus, since the classical Riesz transform is an isomorphism between $L^p(\mathbb{R}^d,dx)$-spaces, with $p \in (1, +\infty)$, for all $\ell \in \{1, \dots, d\}$, 
\begin{align}\label{ineq:perim_fractional_Sobolev_embedding_square_rot}
\|f\|_{L^{2_\alpha^*}(\mathbb{R}^d,dx)} \leq C_{\alpha,d} \|D^{\alpha-1, \operatorname{rot}}_{\ell}(f)\|_{L^2(\mathbb{R}^d,dx)},
\end{align}
with $C_{\alpha,d} >0$.~Next, let us discuss the general case $p \in [1,d/(\alpha-1))$, $d \geq 2$ and $\alpha \in (1,2)$ for $\nu_\alpha = \nu_\alpha^{\operatorname{rot}}$. By the Hardy-Littlewood-Sobolev theorem (see \cite[Chapter $V$, Section $1.2$, Theorem $1$]{Stein_1}), for all $d \geq 1$, all $s \in (0,1)$, all $p \in (1, d/s)$ and all $f \in L^p(\mathbb{R}^d,dx)$, 
\begin{align}\label{ineq:HLS_isotropic}
\| I^s(f)  \|_{L^q(\mathbb{R}^d,dx)} \leq C_{s,p,d} \|f\|_{L^p(\mathbb{R}^d,dx)},
\end{align}
for some $C_{s,p,d}>0$ depending on $s$, $p$ and $d$ and for $q = pd/(d-ps)$. Thanks to the inequality \eqref{ineq:HLS_isotropic} together with \eqref{eq:factorization_rot_invariant} and $I^{\alpha-1} \circ (-\Delta)^{(\alpha-1)/2} = E$, for all $f \in \mathcal{C}^{\infty}_c(\mathbb{R}^d)$ and all $\ell \in \{1, \dots, d\}$,
\begin{align*}
\|f\|_{L^{p_\alpha^*}(\mathbb{R}^d,dx)} \leq \tilde{C}_{\alpha,p,d} \| D^{\alpha-1, \operatorname{rot}}_{\ell}(f) \|_{L^p(\mathbb{R}^d,dx)},
\end{align*}
for some constant $\tilde{C}_{\alpha,p,d}>0$ depending only on $\alpha$, $p$ and $d$. 

Actually, the above discussion and observations are contained in \cite{Shieh_Spector_ACV15,Shieh_Spector_ACV18}. Indeed, \cite[Theorem $1.8$]{Shieh_Spector_ACV15} implies, for all $p \in (1,d/(\alpha-1))$ and all $f \in \mathcal{C}_c^{\infty}(\mathbb{R}^d)$, 
\begin{align}\label{ineq:strong_Sobolev_inequality_rotinv}
\|f\|_{L^{p_\alpha^*}(\mathbb{R}^d,dx)} \leq c_{\alpha,d,p} \|D^{\alpha-1, \operatorname{rot}}(f)\|_{L^p(\mathbb{R}^d,dx)}, 
\end{align} 
for some $c_{\alpha,d,p}>0$ depending on $\alpha$, $d$ and $p$.~The proof of this inequality is based on a composition formula (see \cite[Theorem $1.12$]{Shieh_Spector_ACV15}) named the fractional fundamental theorem of calculus. The fractional fundamental theorem of calculus of \cite[Theorem $1.12$]{Shieh_Spector_ACV15} is a generalization of an identity due to E.~M.~Stein \cite[Chapter $V$, Section $2$, identity $(17)$]{Stein_1}.~Actually, this identity is compared to the following one due to Sobolev (\cite{Sobolev_MS38}): for all $f \in \mathcal{C}_c^{\infty}(\mathbb{R}^d)$ and all $x \in \mathbb{R}^d$, 
\begin{align}\label{eq:identity18_Sobolev_method}
f(x) & = \left(2 \pi^{\frac{d}{2}} \dfrac{1}{\Gamma(\frac{d}{2})}\right)^{-1} \int_{\mathbb{R}^d} \dfrac{1}{\|x-y\|^{d}} \langle x-y; \nabla(f)(y) \rangle dy,
\end{align}
 which is \cite[Chapter $V$, Section $2$, identity $(18)$]{Stein_1}. 
 
 At this point, it is clear that one can extend Sobolev's method to deduce the strong Sobolev inequality \eqref{ineq:strong_Sobolev_Frac} with $D^{\alpha-1, \operatorname{rot}}$ using the integral representation \eqref{eq:spatial_rep_comistefani} together with the mapping properties of the standard Riesz potential operator $I^{\alpha-1}$ (recall \eqref{ineq:HLS_isotropic}). Indeed, for all $f \in \mathcal{C}_c^{\infty}(\mathbb{R}^d)$ and all $x \in \mathbb{R}^d$, 
\begin{align}\label{ineq:pointwise_domination_rotinv}
|f(x)| & \leq  \frac{2}{\alpha} \left(2^{\alpha-1} \pi^{\frac{d}{2}} \dfrac{\Gamma\left(\frac{\alpha}{2}\right)}{\Gamma(\frac{d-\alpha+2}{2})}\right)^{-1} \int_{\mathbb{R}^d} \dfrac{1}{\|x-y\|^{d+1-\alpha}} \| D^{\alpha-1,\operatorname{rot}}(f)(y) \| dy \nonumber \\
& \leq c_{\alpha,d} I^{\alpha-1}\left( \| D^{\alpha-1,\operatorname{rot}}(f) \| \right)(x),
\end{align}
for some $c_{\alpha,d}>0$ depending on $\alpha$ and $d$ (which can be made explicit). By Remark \ref{rem:FFTC_dimension1} and the continuity properties of the Riesz potential operator with $d=1$ (see \cite[Chapter $2$, Section $12.1$]{SKM93}), the Sobolev method can be applied to deduce that, for all $f \in \mathcal{C}^{\infty}_c(\mathbb{R})$ and all $p \in (1,1/(\alpha-1))$, 
\begin{align}\label{ineq:fractional_Sobolev_ineq_dimension1}
\|f\|_{L^{p_\alpha^*}(\mathbb{R},dx)} \leq C_{\alpha,p} \|D^{\alpha-1,\operatorname{rot}}(f)\|_{L^{p}(\mathbb{R},dx)}, 
\end{align}
for some $C_{\alpha,p}>0$ depending on $\alpha$ and $p$ and with $p_\alpha^* = p /(1-p(\alpha-1))$. 

Let us perform a similar analysis when the reference $\alpha$-stable probability measure is $\mu_{\alpha,d}$ given by \eqref{eq:StableIndAxes}. Then, for all $\xi \in \mathbb{R}^d$, $\widehat{\mu_{\alpha,d}}(\xi) = \exp\left(- \|\xi\|_\alpha^\alpha\right)$. By Fourier arguments, for all $f \in \mathcal{C}_c^{\infty}(\mathbb{R}^d)$ and all $x \in \mathbb{R}^d$, 
\begin{align}\label{eq:fractional_discrete_case}
f(x) & = \left(- \mathcal{A}_\alpha\right)^{-1} \circ \left(- \mathcal{A}_\alpha\right)(f)(x) \nonumber \\
& = - \frac{1}{\alpha} \left(-\mathcal{A}_\alpha\right)^{-1} \circ (\nabla \cdot D^{\alpha-1,d})(f)(x) \nonumber \\
& = \frac{1}{\alpha} I_\sigma^{\alpha-1} \circ R_\sigma \cdot D^{\alpha-1,d}(f)(x), 
\end{align}
where the operators $I_\sigma^{\alpha-1}$ and $R_\sigma$ are given, for all $f \in \mathcal{C}_c^{\infty}(\mathbb{R}^d)$ and all $x \in \mathbb{R}^d$, by
\begin{align}\label{eq:anisotropic_vectorial_Riesz_transform2}
R_\sigma(f)(x) & = \frac{1}{(2\pi)^d} \int_{\mathbb{R}^d} \mathcal{F}(f)(\xi) e^{i \langle x ; \xi \rangle} \frac{-i\xi}{\|\xi\|_\alpha}d\xi, 
\end{align} 
and
\begin{align}\label{eq:anisotropic_Riesz_potential}
I_\sigma^{\alpha-1}(f)(x) & = \frac{1}{(2\pi)^d} \int_{\mathbb{R}^d} \mathcal{F}(f)(\xi) e^{i \langle \xi ; x \rangle} \frac{d\xi}{\|\xi\|_\alpha^{\alpha-1}}. 
\end{align} 
\noindent
From Lemma \ref{lem:l_alpha_ball} and Lemma \ref{lem:sharp_pointwise_bound}, it is natural to introduce a Riesz polypotential operator thanks to the pointwise bound \eqref{ineq:pointwise_bound_tensorized_kernel}. This operator already appears in the work \cite{Ok_JLMS69}. Let $I^{\alpha-1,d}$ be the linear operator defined, for all $f \in \mathcal{S}(\mathbb{R}^d)$ and all $x \in \mathbb{R}^d$, by 
\begin{align}\label{eq:Rpotential_product}
I^{\alpha-1,d}(f)(x) = (\gamma_{(\alpha-1)/d,1})^d \int_{\mathbb{R}^d} \prod_{k=1}^d |x_k - y_k|^{- \frac{d+1-\alpha}{d}} f(y) dy,
\end{align} 
where $\gamma_{(\alpha-1)/d,1}$ is given by \eqref{eq:normalizing_constant_Rieszkernel}. In particular, for all $f = f_1 \otimes \dots \otimes f_d$, with $f_i \in \mathcal{C}_c^{\infty}(\mathbb{R}^d)$, $i \in \{1, \dots, d\}$,  
\begin{align}\label{eq:product_structure_Riesz_operator}
I^{\alpha-1,d}(f_1 \otimes \dots \otimes f_d) = I^{\frac{\alpha-1}{d}}(f_1) \otimes \dots \otimes  I^{\frac{\alpha-1}{d}}(f_d),
\end{align}
The next proposition investigates the mapping properties of the linear operator $I^{\alpha-1,d}$ on the classical Lebesgue spaces $L^p(\mathbb{R}^d,dx)$, with $p > 1$.

\begin{prop}\label{prop:well-defined_continuity}
Let $d \geq 2$ be an integer, let $\alpha \in (1,2)$ and let $I^{\alpha-1,d}$ be the linear operator defined by \eqref{eq:Rpotential_product}. Then, 
\begin{enumerate}
\item for all $p \in (1, d/(\alpha-1))$ and all $f \in L^p(\mathbb{R}^d,dx)$, 
\begin{align}\label{ineq:Riesz_LpLq_ineq}
\| I^{\alpha-1,d}(f)  \|_{L^{p_\alpha^*}(\mathbb{R}^d,dx)} \leq A_{p,\alpha,d} \|f\|_{L^p(\mathbb{R}^d,dx)},
\end{align}
for some $A_{p,\alpha,d}>0$ depending on $p$, $\alpha$ and $d$;
\item for all $f \in \mathcal{C}^{\infty}_c(\mathbb{R}^d)$ and all $p \in (1, d/(\alpha-1))$, 
\begin{align}\label{ineq:strong_sobolev_axes}
\|f\|_{L^{p_\alpha^*}(\mathbb{R}^d,dx)} \leq C_{p,\alpha,d} \|D^{\alpha-1,d}(f)\|_{L^p(\mathbb{R}^d,dx)},
\end{align}
for some $C_{p, \alpha,d}>0$ depending on $p$, $\alpha$ and $d$. 
\end{enumerate} 
\end{prop}

\begin{proof}
Inequality \eqref{ineq:Riesz_LpLq_ineq} is a direct consequence of the product structure of the convolution kernel of $I^{\alpha-1,d}$ together with the classical Hardy-Littlewood-Sobolev inequality \eqref{ineq:HLS_isotropic} in dimension $1$ (see, e.g., \cite[Section $24.12$, Theorem $24.9$]{SKM93}).~Next, let us prove inequality \eqref{ineq:strong_sobolev_axes}.~Then, by \eqref{eq:formula_Stein_axes} of Lemma \ref{lem:l_alpha_ball} and \eqref{ineq:pointwise_bound_tensorized_kernel} of Lemma \ref{lem:sharp_pointwise_bound}, for all $f \in \mathcal{C}_c^{\infty}(\mathbb{R}^d)$ and all $x \in \mathbb{R}^d$, 
\begin{align}\label{ineq:pointwise_domination_axes}
|f(x)| & = \left|  \frac{1}{\alpha} \int_{\mathbb{R}^d} \langle k_{\alpha,d}(x-y) ; D^{\alpha-1,d}(f)(y) \rangle dy \right| \nonumber \\
& \leq \frac{1}{\alpha} \int_{\mathbb{R}^d} \| k_{\alpha,d}(x-y) \| \left\| D^{\alpha-1,d}(f)(y) \right\| dy \nonumber \\
& \leq c_{\alpha,d} I^{\alpha-1,d}\left( \left\| D^{\alpha-1,d}(f) \right\|\right)(x),
\end{align}
for some $c_{\alpha,d}>0$ depending on $\alpha$ and $d$. Then, inequality \eqref{ineq:strong_sobolev_axes} is a straightforward consequence of  inequality \eqref{ineq:Riesz_LpLq_ineq} together with the pointwise bound \eqref{ineq:pointwise_domination_axes}. 
\end{proof}
\noindent
In the next proposition, let us investigate the fractional Sobolev inequality \eqref{ineq:strong_Sobolev_Frac} in the general non-degenerate symmetric $\alpha$-stable situation for which the isotropic pointwise estimates are strong enough. 

\begin{prop}\label{prop:SI_NDSstrong}
Let $\alpha \in (1,2)$ and let $\nu_\alpha$ be a non-degenerate symmetric $\alpha$-stable L\'evy measure on $\mathbb{R}^d$, $d \geq 2$, such that $\nu_\alpha$ is a $\gamma$-measure with $\gamma \in [1,d]$ and $\gamma - d + 2\alpha > 1$. Then, for all $p \in (1,d/(\alpha-1))$ and all $f \in \mathcal{C}^{\infty}_c(\mathbb{R}^d)$, 
\begin{align}\label{ineq:FSI_NDS}
\|f\|_{L^{p_\alpha^*}(\mathbb{R}^d,dx)} \leq C_{\alpha,p,d,\gamma} \|D^{\alpha-1}(f)\|_{L^p(\mathbb{R}^d,dx)}, 
\end{align}
where $C_{\alpha,p,d,\gamma}>0$ depends on $\alpha$, $p$, $d$ and $\gamma$.
\end{prop}

\begin{proof}
Thanks to \eqref{ineq:pointwise_estimate_FFTCkernel}, for all $f \in \mathcal{C}_c^{\infty}(\mathbb{R}^d)$ and all $x \in \mathbb{R}^d$,
\begin{align}\label{ineq:pointwise_inequality_NDS}
\left| f(x) \right| & \leq \frac{1}{\alpha} \int_{\mathbb{R}^d} \|k_{\alpha}(x-y)\| \|D^{\alpha-1}(f)(y)\| dy \nonumber \\
& \leq C_{\alpha,d,\gamma} \int_{\mathbb{R}^d} \|x-y\|^{\alpha-d-1} \|D^{\alpha-1}(f)(y)\| dy \nonumber \\
& \leq C_{\alpha,d,\gamma} I^{\alpha-1}\left(\|D^{\alpha-1}(f)\|\right)(x). 
\end{align} 
The classical Hardy-Littlewood-Sobolev theorem ensures that inequality \eqref{ineq:FSI_NDS} holds true. 
\end{proof}
\noindent
We have now arrived at the proof of Theorem \ref{thm:FSI_NDS_full} for the whole family of non-degenerate symmetric $\alpha$-stable probability measures, with $\alpha \in (1,2)$.\\

\noindent   
\textit{Proof of Theorem \ref{thm:FSI_NDS_full}.}
\textbf{Step 1 :} Let us prove that, for all $1 \leq r < s \leq +\infty$, all $f \in L^r(\mathbb{R}^d,dx)$ and all $t>0$, 
\begin{align}\label{ineq:ultracontractive_rs}
\|P^\alpha_t(f)\|_{L^s(\mathbb{R}^d,dx)} \leq t^{ - \frac{d}{\alpha} (\frac{1}{r} - \frac{1}{s})} \|f\|_{L^r(\mathbb{R}^d,dx)} \|p_\alpha\|^{1-\frac{r}{s}}_{L^{r'}(\mathbb{R}^d,dx)}. 
\end{align}
This follows from a standard interpolation argument combined with ultracontractivity. Indeed, for all $1 \leq r < s \leq +\infty$, all $f \in L^r(\mathbb{R}^d,dx)$ and all $t>0$, 
\begin{align*}
\|P^\alpha_t(f)\|_{L^s(\mathbb{R}^d,dx)}^s & = \int_{\mathbb{R}^d} |P_t^\alpha(f)(x)|^s dx = \int_{\mathbb{R}^d}|P_t^\alpha(f)(x)|^{s-r} |P_t^\alpha(f)(x)|^{r} dx \\
& \leq \|P_t^\alpha(f)\|^{s-r}_{L^{\infty}(\mathbb{R}^d,dx)}\int_{\mathbb{R}^d}|P_t^\alpha(f)(x)|^{r} dx \\
& \leq \|P_t^\alpha(f)\|^{s-r}_{L^{\infty}(\mathbb{R}^d,dx)} \|f\|^r_{L^r(\mathbb{R}^d,dx)}, 
\end{align*}
since $(P_t^{\alpha})_{t \geq 0}$ is a contraction semigroup on $L^r(\mathbb{R}^d,dx)$. Moreover, by the H\"older inequality with $r>1$, for all $t>0$ and all $x \in \mathbb{R}^d$, 
\begin{align*}
|P_t^{\alpha}(f)(x)| & \leq \|f\|_{L^r(\mathbb{R}^d,dx)} \left(\int_{\mathbb{R}^d}p_\alpha \left(\dfrac{x-z}{t^{\frac{1}{\alpha}}}\right)^{r'}\frac{dz}{t^{\frac{dr'}{\alpha}}}\right)^{\frac{1}{r'}} \\
& \leq \|f\|_{L^r(\mathbb{R}^d,dx)} t^{\frac{d}{\alpha}\left(\frac{1}{r'}-1\right)} \|p_\alpha\|_{L^{r'}(\mathbb{R}^d,dx)} \\
& \leq t^{- \frac{d}{\alpha r}} \|f\|_{L^r(\mathbb{R}^d,dx)}\|p_\alpha\|_{L^{r'}(\mathbb{R}^d,dx)},
\end{align*}
with $r' = r/(r-1)$. A similar pointwise bound holds true when $r=1$ since $\|p_\alpha\|_{\infty}<+\infty$. Then, for all $1 \leq r < s \leq +\infty$, all $f \in L^r(\mathbb{R}^d,dx)$ and all $t>0$, 
\begin{align*}
\|P^\alpha_t(f)\|_{L^s(\mathbb{R}^d,dx)}^s & \leq t^{- \frac{d}{\alpha r}(s-r)} \|f\|^{s-r}_{L^r(\mathbb{R}^d,dx)}\|p_\alpha\|^{s-r}_{L^{r'}(\mathbb{R}^d,dx)} \|f\|^r_{L^r(\mathbb{R}^d,dx)} \\
& \leq t^{- \frac{d}{\alpha}(\frac{s}{r}-1)} \|f\|^s_{L^r(\mathbb{R}^d,dx)}\|p_\alpha\|^{s-r}_{L^{r'}(\mathbb{R}^d,dx)}.
\end{align*}
This concludes the proof of \eqref{ineq:ultracontractive_rs}.\\ 
\textbf{Step 2 :} Let $\beta  = (\alpha-1)/\alpha$ and $0< \delta <R < +\infty$. Let $T_{\delta,R}$ be the linear operator defined, for all $f \in \mathcal{C}^{\infty}_c(\mathbb{R}^d)$, by
\begin{align}\label{eq:truncated_fractional_power}
T_{\delta,R}(f) = \frac{1}{\Gamma(\beta)} \int_{\delta}^R t^{\beta-1} P^{\alpha}_{t}(f)dt. 
\end{align}
Let us prove that this truncated operator extends continuously from $L^p(\mathbb{R}^d,dx)$ to $L^{p_\alpha^*,\infty}(\mathbb{R}^d,dx)$, for all $p \in (1, d/(\alpha-1))$, with a constant not depending on the truncation parameters. Here and in the sequel, $L^{p_\alpha^*,\infty}(\mathbb{R}^d,dx)$ denotes the weak Lebesgue space of exponent $p_\alpha^*$ (see, e.g., \cite[Volume $I$, Chapter $1$, Section $1.1$]{G08}). Let $a \in (\delta, R)$. By a standard union bound, 
\begin{align*}
\mathcal{L}_d \left(|T_{\delta,R}(f)(x)| > \rho\right) \leq \mathcal{L}_d \left(|T_{\delta,a}(f)(x)| > \frac{\rho}{2}\right) + \mathcal{L}_d \left(|T_{a,R}(f)(x)| > \frac{\rho}{2}\right).
\end{align*} 
Next, let $p<s_0 < p_\alpha^* < s_1 < + \infty$. Then, from \eqref{ineq:ultracontractive_rs}, \eqref{eq:truncated_fractional_power} and the Minkowski integral inequality, for all $f \in \mathcal{C}_c^{\infty}(\mathbb{R}^d)$, 
\begin{align*}
\|T_{\delta,a}(f)\|_{L^{s_0}(\mathbb{R}^d,dx)} & \leq \frac{1}{\Gamma(\beta)} \int_{\delta}^a t^{\beta-1} \|P_t^{\alpha}(f)\|_{L^{s_0}(\mathbb{R}^d,dx)} dt \\
& \leq \frac{1}{\Gamma(\beta)} \int_{\delta}^{a} t^{\beta-1} t^{- \frac{d}{\alpha}(\frac{1}{p} - \frac{1}{s_0})} dt \|f\|_{L^p(\mathbb{R}^d,dx)} \|p_\alpha\|^{1-\frac{p}{s_0}}_{L^{p'}(\mathbb{R}^d,dx)} \\
& \leq \frac{1}{\Gamma(\beta)}\int_{\delta}^{a} t^{\frac{d}{\alpha}\left(\frac{1}{s_0}-\frac{1}{p_{\alpha}^*}\right)} \frac{dt}{t} \|f\|_{L^p(\mathbb{R}^d,dx)} \|p_\alpha\|^{1-\frac{p}{s_0}}_{L^{p'}(\mathbb{R}^d,dx)} \\
& \leq C_0(\alpha,d,p,s_0) a^{\frac{d}{\alpha}\left(\frac{1}{s_0}-\frac{1}{p_\alpha^*}\right)} \|f\|_{L^p(\mathbb{R}^d,dx)}, 
\end{align*}
where we have used the facts that $\beta = d(1/p - 1/p_\alpha^*)/\alpha$ and $s_0 < p_\alpha^*$ and where $C_0(\alpha,d,p,s_0)>0$ depends on $\alpha$, $d$, $p$ and $s_0$. Similarly, for all $f \in \mathcal{C}_c^{\infty}(\mathbb{R}^d)$,
\begin{align*}
\|T_{a,R}(f)\|_{L^{s_1}(\mathbb{R}^d,dx)} & \leq \frac{1}{\Gamma(\beta)} \int_{a}^R t^{\beta-1} \|P_t^{\alpha}(f)\|_{L^{s_1}(\mathbb{R}^d,dx)} dt \\
& \leq \frac{1}{\Gamma(\beta)} \int_{a}^{R} t^{\beta-1} t^{- \frac{d}{\alpha}(\frac{1}{p} - \frac{1}{s_1})} dt \|f\|_{L^p(\mathbb{R}^d,dx)} \|p_\alpha\|^{1-\frac{p}{s_1}}_{L^{p'}(\mathbb{R}^d,dx)} \\
& \leq \frac{1}{\Gamma(\beta)}\int_{a}^{R} t^{\frac{d}{\alpha}\left(\frac{1}{s_1}-\frac{1}{p_{\alpha}^*}\right)} \frac{dt}{t} \|f\|_{L^p(\mathbb{R}^d,dx)} \|p_\alpha\|^{1-\frac{p}{s_1}}_{L^{p'}(\mathbb{R}^d,dx)} \\
& \leq C_1(\alpha,d,p,s_1) a^{-\frac{d}{\alpha}\left(\frac{1}{p_\alpha^*}-\frac{1}{s_1}\right)} \|f\|_{L^p(\mathbb{R}^d,dx)},
\end{align*}
using that $s_1 > p_\alpha^*$ and where $C_1(\alpha,d,p,s_1)>0$ depends on $\alpha$, $d$, $p$ and $s_1$. Now, by Chebyshev's inequality with $\|f\|_{L^p(\mathbb{R}^d,dx)} = 1$ and taking $a = \rho^{ - \frac{\alpha p_\alpha^*}{d}}$,  
\begin{align*}
\mathcal{L}_d \left(|T_{\delta,R}(f)(x)| > \rho\right) & \leq \dfrac{2^{s_1}}{\rho^{s_1}} \|T_{a,R}(f)\|^{s_1}_{L^{s_1}(\mathbb{R}^d,dx)} + \dfrac{2^{s_0}}{\rho^{s_0}} \|T_{\delta,a}(f)\|^{s_0}_{L^{s_0}(\mathbb{R}^d,dx)} \\
& \leq C(\alpha,d,p,s_0,s_1) \left( \frac{1}{\rho^{s_1}} \|f\|^{s_1}_{L^p(\mathbb{R}^d,dx)} a^{ - \frac{s_1 d}{\alpha} \left(\frac{1}{p_\alpha^*}-\frac{1}{s_1}\right)} + \frac{1}{\rho^{s_0}} \|f\|^{s_0}_{L^p(\mathbb{R}^d,dx)} a^{\frac{s_0 d}{\alpha} \left(\frac{1}{s_0} - \frac{1}{p_\alpha^*}\right)}  \right) \\
& \leq 2 C(\alpha,d,p,s_0,s_1) \rho^{ - p_\alpha^*},
\end{align*}
where $C(\alpha,d,p,s_0,s_1)>0$ depends on $\alpha$, $d$, $p$, $s_0$ and $s_1$. This concludes the proof of the second step.\\
\textbf{Step 3 :} Now, let $p_0$ and $p_1$ be such that $1< p_0 < p < p_1 < d/(\alpha-1)$. Observe that $p_{\alpha,0}^* = p_0 d / (d - p_0(\alpha-1))$ and $p_{\alpha,1}^* = p_1 d / (d - p_1(\alpha-1))$ are such that $p_{\alpha,0}^* < p_\alpha^* < p_{\alpha,1}^*$. Then, by the Marcinkiewicz interpolation theorem, for all $0<\delta<R< +\infty$, $T_{\delta,R}$ is strong-type $(p,p_\alpha^*)$ with a constant which does not depend on $\delta$ neither on $R$. Namely, for all $f \in L^p(\mathbb{R}^d,dx)$ and all $0<\delta <R<+\infty$,
\begin{align}\label{ineq:strong_type_uniform}
\|T_{\delta,R}(f)\|_{L^{p_\alpha^*}(\mathbb{R}^d,dx)} \leq C \|f\|_{L^p(\mathbb{R}^d,dx)}, 
\end{align}
where $C>0$ does not depend on $\delta$ neither on $R$. Next, let $f \in \mathcal{C}_c^{\infty}(\mathbb{R}^d)$ and $(-\mathcal{A}_\alpha)^{\frac{\alpha-1}{\alpha}}(f)$ be defined, for all $x \in \mathbb{R}^d$, by 
\begin{align*}
(-\mathcal{A}_\alpha)^{\frac{\alpha-1}{\alpha}}(f)(x) = \frac{1}{(2\pi)^d} \int_{\mathbb{R}^d} \mathcal{F}(f)(\xi) (- \psi_\alpha(\xi))^{\frac{\alpha-1}{\alpha}} e^{i \langle x ; \xi \rangle} d\xi. 
\end{align*}
By Lemma \ref{lem:Lp_bound_FractionalPower}, $(-\mathcal{A}_\alpha)^{\frac{\alpha-1}{\alpha}}(f) \in L^p(\mathbb{R}^d,dx)$. Thus, for all $f \in \mathcal{C}^{\infty}_c(\mathbb{R}^d)$ and all $0<\delta<R<+\infty$, 
\begin{align*}
\|T_{\delta,R}((-\mathcal{A}_\alpha)^{\frac{\alpha-1}{\alpha}}(f))\|_{L^{p_\alpha^*}(\mathbb{R}^d,dx)} \leq C \|(-\mathcal{A}_\alpha)^{\frac{\alpha-1}{\alpha}}(f)\|_{L^p(\mathbb{R}^d,dx)}. 
\end{align*}
Moreover, by standard Fourier analysis and the non-degeneracy condition \eqref{eq:non_deg}, for all $f \in \mathcal{C}_c^{\infty}(\mathbb{R}^d)$ and all $x \in \mathbb{R}^d$, 
\begin{align*}
\underset{\delta \rightarrow 0^+, R\rightarrow +\infty}{\lim}T_{\delta,R}((-\mathcal{A}_\alpha)^{\frac{\alpha-1}{\alpha}}(f))(x) = f(x).  
\end{align*}
Then, by Fatou's lemma, for all $f \in \mathcal{C}^{\infty}_c(\mathbb{R}^d)$, 
\begin{align*}
\|f\|_{L^{p_\alpha^*}(\mathbb{R}^d,dx)} \leq C \|(-\mathcal{A}_\alpha)^{\frac{\alpha-1}{\alpha}}(f)\|_{L^p(\mathbb{R}^d,dx)}. 
\end{align*}
Inequality \eqref{ineq:Lp_bound_FractionalPower} concludes the proof of the theorem.$\qed$

\begin{thm}\label{thm:sobolev-inequality-anisotropic}
Let $d \geq 2$ be an integer, let $\sigma$ be a symmetric finite positive measure on $\mathbb{S}^{d-1}$ verifying the non-degeneracy condition \eqref{eq:AsymNDC} and let $\mathbf{\Sigma}$ be the symmetric positive definite matrix defined in \eqref{eq:from_measure_to_matrix}. Let $D_\sigma$ be the anisotropic local gradient operator as in \eqref{eq:anisotropic_local_gradient}.~Let $p \in (1,d)$ and let $q = pd/(d-p)$. Then, for all $f \in \mathcal{C}^{\infty}_c(\mathbb{R}^d)$, 
\begin{align}\label{ineq:sobolev-inequality-anisotropic}
\|f\|_{L^q(\mathbb{R}^d,dx)} \leq C_{\sigma,d,p} \| \mathbf{\Sigma}^{-\frac{1}{2}} D_\sigma(f)\|_{L^p(\mathbb{R}^d,dx)}, 
\end{align}
for some $C_{\sigma,d,p}>0$ depending on $d$, $\sigma$ and $p$. 
\end{thm}

\begin{proof}
Let $f \in \mathcal{C}^{\infty}_c(\mathbb{R}^d)$. Then, by Proposition \ref{prop:AL_FTC} and Lemma \ref{lem:from_iso_to_aniso}, for all $x \in \mathbb{R}^d$, 
\begin{align*}
f(x) & =  \dfrac{1}{\sqrt{\operatorname{det}\left(\mathbf{\Sigma}\right)}} \left(\dfrac{2}{\Gamma(\frac{d}{2})}\pi^{\frac{d}{2}}\right)^{-1} \int_{\mathbb{R}^d} \left\langle \dfrac{\mathbf{\Sigma}^{-1}(x-y)}{\| \mathbf{\Sigma}^{-\frac{1}{2}}(x-y)\|^{d}} ; D_\sigma(f)(y)\right\rangle dy \\
& =  \dfrac{1}{\sqrt{\operatorname{det}\left(\mathbf{\Sigma}\right)}} \left(\dfrac{2}{\Gamma(\frac{d}{2})}\pi^{\frac{d}{2}}\right)^{-1} \int_{\mathbb{R}^d} \left\langle \dfrac{\mathbf{\Sigma}^{-\frac{1}{2}}(x-y)}{\| \mathbf{\Sigma}^{-\frac{1}{2}}(x-y)\|^{d}} ; \mathbf{\Sigma}^{-\frac{1}{2}} D_\sigma(f)(y)\right\rangle dy.
\end{align*}
Thus, by the Cauchy-Schwarz inequality, for all $x \in \mathbb{R}^d$, 
\begin{align*}
\left| f(x)\right| \leq  \dfrac{1}{\sqrt{\operatorname{det}\left(\mathbf{\Sigma}\right)}} \left(\dfrac{2}{\Gamma(\frac{d}{2})}\pi^{\frac{d}{2}}\right)^{-1} \int_{\mathbb{R}^d} \dfrac{1}{\| \mathbf{\Sigma}^{-\frac{1}{2}}(x-y)\|^{d-1}} \left\| \mathbf{\Sigma}^{-\frac{1}{2}} D_\sigma(f)(y)\right\| dy. 
\end{align*}
Moreover, for all $z \in \mathbb{R}^d$, 
\begin{align*}
\| \mathbf{\Sigma}^{-\frac{1}{2}} z \| \geq \frac{1}{\sqrt{\lambda_{\operatorname{max}}\left(\mathbf{\Sigma}\right)}} \|z\|,
\end{align*}
where $\lambda_{\operatorname{max}}\left(\mathbf{\Sigma}\right)$ is the largest eigenvalue of the matrix $\mathbf{\Sigma}$. Thus, for all $x \in \mathbb{R}^d$, 
\begin{align*}
\left| f(x)\right| & \leq  \dfrac{\lambda_{\operatorname{max}}\left(\mathbf{\Sigma}\right)^{\frac{d-1}{2}}}{\sqrt{\operatorname{det}\left(\mathbf{\Sigma}\right)}} \left(\dfrac{2}{\Gamma(\frac{d}{2})}\pi^{\frac{d}{2}}\right)^{-1} \int_{\mathbb{R}^d} \dfrac{1}{\| (x-y)\|^{d-1}} \left\| \mathbf{\Sigma}^{-\frac{1}{2}} D_\sigma(f)(y)\right\| dy \\
& \leq C_{\sigma,d,p} I^{1} \left(  \left\| \mathbf{\Sigma}^{-\frac{1}{2}} D_\sigma(f)\right\| \right)(x). 
\end{align*}
The Hardy-Littlewood-Sobolev inequality concludes the proof of the theorem.
\end{proof}
\noindent
Now, let $H_\sigma$ be the mapping defined, for all $x \in \mathbb{R}^d$, by 
\begin{align}\label{eq:Gauge_Function_Aniso}
H_\sigma(x) = \left\| \mathbf{\Sigma}^{-\frac{1}{2}} \int_{\mathbb{S}^{d-1}} y \langle x ; y \rangle \sigma(dy) \right\|,
\end{align}
where $\mathbf{\Sigma}$ is given in \eqref{eq:from_measure_to_matrix} and $\sigma$ satisfies \eqref{eq:AsymNDC}.~The next technical lemma provides standard properties of $H_\sigma$.

\begin{lem}\label{lem:prop_Gauge_Function_Aniso}
Let $\sigma$ be a symmetric positive finite measure on $\mathbb{S}^{d-1}$, $d \geq 2$, verifying \eqref{eq:AsymNDC}. Then, 
\begin{enumerate}
\item $H_\sigma$ is a norm on $\mathbb{R}^d$;
\item the dual norm $\overset{\circ}{H}_\sigma$ is given, for all $x \in \mathbb{R}^d$, by 
\begin{align}\label{eq:Dual_Gauge_Function_Aniso}
\overset{\circ}{H}_\sigma(x) = \left\| \mathbf{\Sigma}^{-\frac{1}{2}}x \right\|.
\end{align}
\end{enumerate}
\end{lem}

\begin{proof}
The triangle inequality and the non-negative $1$-homogeneity of $H_\sigma$ are straightforward from \eqref{eq:Gauge_Function_Aniso}. Let us prove that 
\begin{align}\label{eq:positiveness}
H_\sigma(x) = 0 \Rightarrow x = 0. 
\end{align}
Since $\|\cdot\|$ is the Euclidean norm on $\mathbb{R}^d$ and since the matrix $\mathbf{\Sigma}$ is non-degenerate,
\begin{align*}
H_\sigma(x) = 0 \Rightarrow \int_{\mathbb{S}^{d-1}} y \langle y; x \rangle \sigma(dy) = 0.
\end{align*}
But, the non-degeneracy condition ensures that 
\begin{align*}
\int_{\mathbb{S}^{d-1}} \langle y ; x \rangle^2 \sigma(dy) \geq \|x\|^{2} \underset{e \in \mathbb{S}^{d-1}}{\inf} \int_{\mathbb{S}^{d-1}} \left| \langle y ; e \rangle \right|^2 \sigma(dy). 
\end{align*}
Thus, \eqref{eq:positiveness} is proved. To conclude the proof of the lemma, let us obtain a simple formula for the norm $H_\sigma$. Note that the mapping $L_\sigma$ from $\mathbb{R}^d$ to $\mathbb{R}^d$ defined, for all $x \in \mathbb{R}^d$, by 
\begin{align*}
L_\sigma(x) := \int_{\mathbb{S}^{d-1}} y \langle x; y \rangle \sigma(dy),
\end{align*}
is linear. Moreover, by standard linear algebra arguments, the associated $d \times d$ matrix is exactly the symmetric positive definite matrix $\mathbf{\Sigma}$. Thus, for all $x \in \mathbb{R}^d$, 
\begin{align*}
H_\sigma(x) = \left\| \mathbf{\Sigma}^{- \frac{1}{2}} \mathbf{\Sigma} x \right\| = \left\| \mathbf{\Sigma}^{\frac{1}{2}}x \right\|. 
\end{align*}
By duality, 
\begin{align*}
\overset{\circ}{H}_{\sigma}(x) : = \underset{y \in \mathbb{R}^d \setminus \{0\},\, H_\sigma(y)=1}{\sup} \left| \langle y;x \rangle \right| = \left\| \mathbf{\Sigma}^{-\frac{1}{2}}x \right\|.
\end{align*}
This conclude the proof the lemma. 
\end{proof}
\noindent 
Sharp anisotropic Sobolev inequalities as well as their optimizers have been obtained in \cite{cordero_nazaret_villani}. In particular, \cite[Theorem $2$, (ii)]{cordero_nazaret_villani}, the sharp version of  \eqref{ineq:sobolev-inequality-anisotropic} is given, for all $f \in \mathcal{C}_c^{\infty}\left(\mathbb{R}^d\right)$ and all $p \in (1,d)$, by 
\begin{align}\label{ineq:sharp-sobolev-inequality-anisotropic}
\|\mathbf{\Sigma}^{\frac{1}{2}} \nabla(f)\|_{L^p(\mathbb{R}^d,dx)} \geq \| \mathbf{\Sigma}^{\frac{1}{2}} \nabla(h_\sigma)\|_{L^p(\mathbb{R}^d,dx)} \|f\|_{L^q(\mathbb{R}^d,dx)}, 
\end{align} 
where $q = pd /(d-p)$ and where, for all $x \in \mathbb{R}^d$, 
\begin{align}\label{eq:optimizer-sobolev-inequality-anisotropic}
h_\sigma(x) = \dfrac{1}{\left(\lambda_p + \| \mathbf{\Sigma}^{-\frac{1}{2}} x  \|^{p'}\right)^{\frac{d-p}{p}}},
\end{align}
with $p' = p/(p-1)$ and with $\lambda_p$ chosen in such a way that $\|h_\sigma\|_{L^q(\mathbb{R}^d,dx)} = 1$.

Next, let us investigate refined fractional Sobolev-type inequalities in order to study sharp versions of \eqref{ineq:strong_Sobolev_Frac}. 

\begin{thm}\label{thm:Refined_Fractional_Sobolev-type_ineq}
Let $\alpha \in (1,2)$ and let $\mu_\alpha$ be a symmetric non-degenerate $\alpha$-stable probability measure on $\mathbb{R}^d$, $d \geq 2$, with associated L\'evy measure $\nu_\alpha$. Let $2_\alpha^* = 2d / (d-2(\alpha-1))$, let $D^{\alpha-1}$ be the fractional gradient operator associated with $\nu_\alpha$, and let $\mathcal{A}_\alpha$ be the generator of the stable heat semigroup $(P_t^\alpha)_{t \geq 0}$. Let $\chi$ be an infinitely differentiable function defined on $\mathbb{R}_+$ with compact support such that $\chi$ is identically equal to $1$ in a neighborhood of the origin. Then, for all $f \in \mathcal{C}_c^{\infty}(\mathbb{R}^d)$, 
\begin{align}\label{ineq:Refined_Fractional_Sobolev-type_ineq}
\|f\|_{L^{2_\alpha^*}(\mathbb{R}^d,dx)} \leq C(\alpha,d,\chi) \left(\underset{t>0}{\sup}\, t^{\frac{d-2(\alpha-1)}{2\alpha}} \|\chi \left(-t \mathcal{A}_\alpha\right)(f)\|_{\infty} \right)^{1-2/2_\alpha^*} \|D^{\alpha-1}(f)\|^{2/2_\alpha^*}_{L^2(\mathbb{R}^d,dx)},
\end{align}
where $C(\alpha,d,\chi)>0$ depends on $\alpha$, $d$ and $\chi$. Moreover, for all $f \in \mathcal{C}_c^{\infty}\left(\mathbb{R}^d\right)$, 
\begin{align}\label{ineq:Fractional_Sobolev-type_ineq_p2}
\|f\|_{L^{2_\alpha^*}(\mathbb{R}^d,dx)} \leq A(\alpha,d,\chi)\|D^{\alpha-1}(f)\|_{L^2(\mathbb{R}^d,dx)}, 
\end{align} 
where $A(\alpha,d,\chi)>0$ depends on $\alpha$, $d$ and $\chi$.  
\end{thm}

\begin{proof}
The proof of this theorem is very similar to the proof of \cite[Proposition $1.2$]{Frank_24}. Let $d$, $\alpha$ and $q$ be as in the statement of the theorem. Let $\beta$ be given by 
\begin{align*}
\beta = \frac{d-2(\alpha-1)}{2\alpha}.
\end{align*}
By the layer cake representation formula, for all $f \in \mathcal{C}^{\infty}_c(\mathbb{R}^d)$, 
\begin{align*}
\|f\|^{2_\alpha^*}_{L^{2_\alpha^*}(\mathbb{R}^d,dx)} = 2_\alpha^* \int_0^{+\infty} t^{2_\alpha^*-1} \mathcal{L}_d \left( |f(x)| > t \right) dt, 
\end{align*}
where $\mathcal{L}_d$ is the $d$-dimensional Lebesgue measure. Now, by a classical union bound, for all $t>0$ and all $s>0$, 
\begin{align*}
\mathcal{L}_d \left( |f(x)| > t \right) \leq \mathcal{L}_d \left( \left|\chi \left(-s \mathcal{A}_\alpha\right)(f)(x) \right| > t/2 \right) + \mathcal{L}_d \left( \left|(E-\chi \left(-s \mathcal{A}_\alpha\right))(f)(x)\right| > t/2 \right). 
\end{align*}
Next, let $s(t)>0$ be defined by 
\begin{align}\label{eq:critical_parameter}
\frac{t}{2} = s(t)^{- \beta} \|f\|_{\infty,\infty,\alpha}, 
\end{align}
with
\begin{align*}
\|f\|_{\infty,\infty,\alpha} := \underset{t>0}{\sup}\, t^{\beta} \|\chi \left(-t \mathcal{A}_\alpha\right)(f)\|_{\infty}. 
\end{align*}
Then, 
\begin{align*}
\mathcal{L}_d \left( \left|\chi \left(-s(t) \mathcal{A}_\alpha\right)(f)(x) \right| > t/2 \right) = 0. 
\end{align*}
Moreover, by Chebyshev's inequality, 
\begin{align*}
\mathcal{L}_d \left( \left|(E-\chi \left(-s \mathcal{A}_\alpha\right))(f)(x)\right| > t/2 \right) \leq \frac{4}{t^2} \|(E-\chi \left(-s(t) \mathcal{A}_\alpha\right))(f)\|^2_{L^2(\mathbb{R}^d,dx)}.
\end{align*}
By the Parseval-Plancherel theorem and Fubini's theorem, 
\begin{align*}
\|f\|^{2_\alpha^*}_{L^{2_\alpha^*}(\mathbb{R}^d,dx)} & \leq 4 2_\alpha^* \int_{0}^{+\infty} t^{2_\alpha^*-3} \|(E-\chi \left(-s(t) \mathcal{A}_\alpha\right))(f)\|^2_{L^2(\mathbb{R}^d,dx)} dt \\
& \leq C(\alpha,d) \int_{0}^{+\infty} t^{2_\alpha^*-3} \left( \int_{\mathbb{R}^d} \left| \mathcal{F}(f)(\xi) \right|^2 \left| 1 - \chi \left(s(t) \sigma_\alpha(\xi)^{\alpha}\right) \right|^2 d\xi \right) dt \\
& \leq C(\alpha,d) \int_{\mathbb{R}^d} \left| \mathcal{F}(f)(\xi) \right|^2 \left( \int_{0}^{+\infty} t^{2_\alpha^*-3} \left| 1 - \chi \left(s(t) \sigma_\alpha(\xi)^{\alpha}\right) \right|^2 dt \right) d\xi. 
\end{align*}
Using \eqref{eq:critical_parameter} and changing variables in the integral with respect to $t$, 
\begin{align*}
\|f\|^{2_\alpha^*}_{L^{2_\alpha^*}(\mathbb{R}^d,dx)} & \leq C(\alpha,d) \left(\int_{\mathbb{R}^d} \left| \mathcal{F}(f)(\xi) \right|^2 \sigma_\alpha\left(\xi\right)^{\alpha \beta (2_\alpha^*-2)} d\xi \right) \|f\|^{2_\alpha^*-2}_{\infty,\infty,\alpha} \int_{0}^{+\infty} s^{2_\alpha^*-3} |1 - \chi \left(\frac{1}{s^{1/\beta}}\right)|^2 ds. 
\end{align*}
Note that the integral with respect to the $s$ variable is finite since $\chi$ is compactly supported on $\mathbb{R}_+$, and is identically equal to $1$ in a right neighborhood of the origin and since $2_\alpha^* = 2d / (d-2(\alpha-1))$. Moreover, using the explicit expressions of $\beta$ and $2_\alpha^*$ and recalling \eqref{eq:Stheatgen_FourRep},  
\begin{align*}
\int_{\mathbb{R}^d}  \left| \mathcal{F}(f)(\xi) \right|^2 \sigma_\alpha\left(\xi\right)^{\alpha \beta (2_\alpha^*-2)} d\xi  & = \int_{\mathbb{R}^d} \left| \mathcal{F}(f)(\xi) \right|^2  \sigma_\alpha\left(\xi\right)^{2(\alpha-1)} d\xi \\
& = \int_{\mathbb{R}^d} \left| \mathcal{F}((- \mathcal{A}_\alpha)^{\frac{\alpha-1}{\alpha}}(f))(\xi) \right|^2 d\xi. 
\end{align*}  
Lemma \ref{lem:Lp_bound_FractionalPower} concludes the proof of the inequality \eqref{ineq:Refined_Fractional_Sobolev-type_ineq}. Thus, by the Fourier inversion formula and the Cauchy-Schwarz inequality, for all $f \in \mathcal{C}_c^{\infty}\left(\mathbb{R}^d\right)$, all $x \in \mathbb{R}^d$ and all $t>0$, 
\begin{align*}
\left| \chi \left(-t \mathcal{A}_\alpha\right)(f)(x) \right| & = \left|  \frac{1}{\left(2\pi\right)^d} \int_{\mathbb{R}^d} \mathcal{F}(f)(\xi)  \chi \left(t \sigma_\alpha(\xi)^\alpha\right) e^{i \langle x ; \xi \rangle} d\xi \right| \\
& \leq \frac{1}{\left(2\pi\right)^d} \int_{\mathbb{R}^d} \left| \mathcal{F}(f)(\xi) \right| \left| \chi \left(t \sigma_\alpha(\xi)^\alpha\right) \right| d\xi \\
& \leq \frac{1}{\left(2\pi\right)^d} \int_{\mathbb{R}^d} \left| \mathcal{F}(f)(\xi) \right| \sigma_\alpha(\xi)^{\alpha-1} \dfrac{\left| \chi \left(t \sigma_\alpha(\xi)^\alpha\right) \right|}{\sigma_\alpha(\xi)^{\alpha-1}} d\xi \\
& \leq c_d \left( \int_{\mathbb{R}^d} \left| \mathcal{F}(f)(\xi) \right|^2 \sigma_\alpha(\xi)^{2(\alpha-1)} d\xi \right)^{\frac{1}{2}} \left(\int_{\mathbb{R}^d}  \dfrac{\left| \chi \left(t \sigma_\alpha(\xi)^\alpha\right) \right|^2}{\sigma_\alpha(\xi)^{2(\alpha-1)}} d\xi\right)^{\frac{1}{2}} \\
& \leq c_d \left( \int_{\mathbb{R}^d} \left| \mathcal{F}\left(\left(-\mathcal{A}_\alpha\right)^{\frac{\alpha-1}{\alpha}}(f)\right)(\xi) \right|^2 d\xi \right)^{\frac{1}{2}} \left(\int_{\mathbb{R}^d}  \dfrac{\left| \chi \left(t \sigma_\alpha(\xi)^\alpha\right) \right|^2}{\sigma_\alpha(\xi)^{2(\alpha-1)}} d\xi\right)^{\frac{1}{2}} .
\end{align*}
Now, changing variables, for all $t>0$, 
\begin{align*}
\int_{\mathbb{R}^d}  \dfrac{\left| \chi \left(t \sigma_\alpha(\xi)^\alpha\right) \right|^2}{\sigma_\alpha(\xi)^{2(\alpha-1)}} d\xi = t^{- \frac{d-2(\alpha-1)}{\alpha}} \int_{\mathbb{R}^d}  \dfrac{\left| \chi \left(\sigma_\alpha(\xi)^\alpha\right) \right|^2}{\sigma_\alpha(\xi)^{2(\alpha-1)}} d\xi.
\end{align*}
Since $\mu_\alpha$ is non-degenerate in the sense of \eqref{eq:non_deg}, $\sigma_\alpha(\omega)$ is comparable to $\|\omega\|$ and the integral on the right-hand side of the previous equality is finite since $\alpha \in (1,2)$, $d\geq 2$ and $\chi$ is compactly supported on $\mathbb{R}_+$. Then, for all $t>0$, 
\begin{align*}
 t^{\frac{d-2(\alpha-1)}{2\alpha}}\left\| \chi \left(-t \mathcal{A}_\alpha\right)(f)\right \|_{\infty} \leq c_{d,\alpha,\chi}  \| (-\mathcal{A}_\alpha)^{\frac{\alpha-1}{\alpha}}(f) \|_{L^2(\mathbb{R}^d,dx)}, 
\end{align*}
which implies, for all $f \in \mathcal{C}_c^{\infty} \left(\mathbb{R}^d\right)$, that
\begin{align*}
\left(\underset{t>0}{\sup}\, t^{\frac{d-2(\alpha-1)}{2\alpha}}\left\| \chi \left(-t \mathcal{A}_\alpha\right)(f)\right \|_{\infty}\right)^{1-\frac{2}{2_\alpha^*}} \leq c_{d,\alpha,\chi} \| (-\mathcal{A}_\alpha)^{\frac{\alpha-1}{\alpha}}(f) \|^{1-\frac{2}{2_\alpha^*}}_{L^2(\mathbb{R}^d,dx)}.
\end{align*}
Lemma \ref{lem:Lp_bound_FractionalPower} concludes the proof.  
\end{proof}
\noindent
\begin{rem}\label{rem:refined_sobolev_ineq}
The non-degeneracy condition \eqref{eq:non_deg} implies that, for all $\xi \in \mathbb{R}^d$, 
\begin{align}\label{ineq:equivalence_norms}
c_{\alpha,d,\sigma} \|\xi\| \leq \sigma_\alpha\left(\xi\right) \leq C_{\alpha,d,\sigma} \|\xi\|, 
\end{align}
for some constants $c_{\alpha,d,\sigma}, C_{\alpha,d,\sigma}>0$ depending on $\alpha$, $d$ and $\sigma$ (the spherical component of the L\'evy measure $\nu_\alpha$). Then, by a proof completely similar to the one given to obtain the inequality \eqref{ineq:Refined_Fractional_Sobolev-type_ineq}, the following refined fractional Sobolev inequalities hold true: for all $f \in \mathcal{C}_c^{\infty}(\mathbb{R}^d)$, 
\begin{align}\label{ineq:refined_frac_sobolev_ineq-2}
\|f\|_{L^{2_\alpha^*}(\mathbb{R}^d,dx)} \leq C_1(\alpha,d,\chi) \left(\underset{t>0}{\sup}\, t^{\frac{d-2(\alpha-1)}{2\alpha}} \|\chi \left(t(-\Delta)^{\frac{\alpha}{2}}\right)(f)\|_{\infty} \right)^{1-\frac{2}{2_\alpha^*}} \|D^{\alpha-1}(f)\|^{\frac{2}{2_\alpha^*}}_{L^2(\mathbb{R}^d,dx)},
\end{align}  
where $C_1(\alpha,d,\chi)>0$ depends on $\alpha$, $d$ and $\chi$ and where $(-\Delta)^{\frac{\alpha}{2}}$ stands for the fractional power of $(-\Delta)$ of order $\frac{\alpha}{2}$ with symbol $\|\xi\|^{\alpha}$, for all $\xi \in \mathbb{R}^d$. Moreover, for all $f \in \mathcal{C}^{\infty}_c(\mathbb{R}^d)$, 
\begin{align}\label{ineq:refined_frac_sobolev_ineq-3}
\|f\|_{L^{2_\alpha^*}(\mathbb{R}^d,dx)} \leq C_2(\alpha,d,\chi) \left(\underset{t>0}{\sup}\, t^{\frac{d-2(\alpha-1)}{4}} \|\chi \left(-t\Delta\right)(f)\|_{\infty} \right)^{1-\frac{2}{2_\alpha^*}} \|D^{\alpha-1}(f)\|^{\frac{2}{2_\alpha^*}}_{L^2(\mathbb{R}^d,dx)},
\end{align} 
where $C_2(\alpha,d,\chi)>0$ depends on $\alpha$, $d$ and $\chi$.
\end{rem}
\noindent
In order to pursue our investigations towards an optimal form for the fractional Sobolev-type inequality \eqref{ineq:FSI_NDS_full}, let us introduce and recall some notations and definitions. Let $p \in (1,d/(\alpha-1))$ and let $p_\alpha^*$ be defined by 
\begin{align}\label{eq:critical-Sobolev-exponent}
p_\alpha^* := \dfrac{pd}{d-p(\alpha-1)}.
\end{align}
Let $\dot{W}^{\alpha-1,p}\left(\mathbb{R}^d,dx\right)$ be the set of functions defined by \eqref{def:Lp-version_homo_frac_Sob_space} and let $|\cdot|_{\alpha-1,p}$ be defined by \eqref{eq:seminorm_Lp}. Thanks to the Riesz transform-type results of \cite{AH20_4} and Lemma \ref{lem:Lp_bound_FractionalPower}, for all $f \in \mathcal{C}_c^{\infty}\left(\mathbb{R}^d\right)$, 
\begin{align}\label{ineq:two-sided_Lpestimates}
\|D^{\alpha-1}(f)\|_{L^p(\mathbb{R}^d,dx)} \asymp \|(-\mathcal{A}_\alpha)^{\frac{\alpha-1}{\alpha}}(f)\|_{L^p(\mathbb{R}^d,dx)}. 
\end{align}
Thanks to the non-degeneracy condition, $|\cdot|_{\alpha-1,p}$ is a norm on $\mathcal{C}_c^{\infty}\left(\mathbb{R}^d\right)$. Let $\mathcal{D}^{\alpha-1,p}(\mathbb{R}^d,dx)$ be the completion of the normed vector space $\left(\mathcal{C}_c^{\infty}\left(\mathbb{R}^d\right) , |\cdot|_{\alpha-1,p}\right)$. By construction, $\mathcal{D}^{\alpha-1,p}(\mathbb{R}^d,dx)$ is a Banach space with dense subset $\mathcal{C}_c^{\infty}\left(\mathbb{R}^d\right)$.~In the sequel, let us prove that, for all $p \in (1, d/(\alpha-1))$, 
\begin{align}\label{eq:Meyer-Serrin-result-Lp}
\dot{W}^{\alpha-1,p}(\mathbb{R}^d,dx) = \mathcal{D}^{\alpha-1,p}(\mathbb{R}^d,dx).
\end{align}
As in the proof of Proposition \ref{prop:density_compactsupport}, it is clear that $\mathcal{C}^{\infty}(\mathbb{R}^d) \cap \dot{W}^{\alpha-1,p}(\mathbb{R}^d,dx) \cap W^{1,\infty}(\mathbb{R}^d,dx)$ is ``dense" in $\dot{W}^{\alpha-1,p}(\mathbb{R}^d,dx)$ with respect to $|\cdot|_{\alpha-1,p}$ by convolution with a sequence of standard mollifiers and Young's inequality (see, e.g., \cite[Theorem $4.33$]{B_FA11}). Moreover, for such a regularization, for all $x \in \mathbb{R}^d$, 
\begin{align*}
D^{\alpha-1}(f)(x) = \int_{\mathbb{R}^d}(f(x+u)-f(x))u\nu_\alpha(du).
\end{align*}
Let us start by an approximation result based on the Banach-Saks property (see \cite[Proposition $10.8.$]{Ponce_book16}). 

\begin{thm}\label{thm:approx_scs_homo_Lp}
Let $\alpha \in (1,2)$, let $d \geq 2$ be an integer and let $p \in (1,d/(\alpha-1))$. Let $f \in \dot{W}^{\alpha-1,p}(\mathbb{R}^d,dx)$. Then, there exists a sequence of functions $(f_n)_{n \geq 1}$ in $\mathcal{C}^{\infty}_c(\mathbb{R}^d)$ such that 
\begin{align*}
f_n \longrightarrow f, \quad n \longrightarrow +\infty,
\end{align*}
in $L^{p^*_\alpha}(\mathbb{R}^d,dx)$, and, 
\begin{align*}
D^{\alpha-1}(f_n) \longrightarrow D^{\alpha-1}(f), \quad n \longrightarrow +\infty,
\end{align*}
in $L^{p}(\mathbb{R}^d,\mathbb{R}^d,dx)$. 
\end{thm}

\begin{proof}
Without loss of generality and from the proof of Proposition \ref{prop:density_compactsupport}, let us assume that $f \in \mathcal{C}^{\infty}(\mathbb{R}^d) \cap \dot{W}^{\alpha-1,p}(\mathbb{R}^d,dx) \cap W^{1,\infty}(\mathbb{R}^d,dx)$. Let $\eta$ be a compactly supported infinitely differentiable function on $\mathbb{R}^d$ such that $\eta(-x) = \eta(x)$, for all $x \in \mathbb{R}^d$, $\operatorname{Supp}(\eta) \subset B(0,2)$, $0\leq \eta \leq 1$ and $\eta(x) = 1$, for all $x \in B(0,1)$. For all $j \geq 1$ an integer, let $\eta_j$ be defined, for all $x \in \mathbb{R}^d$, by 
\begin{align}\label{eq:smooth_cut-off}
\eta_j(x)=\eta\left(\frac{x}{j}\right). 
\end{align} 
Using \cite[Equation (3.8)]{AH23}, for all $x \in \mathbb{R}^d$ and all $j \in \{1,\dots,d\}$, 
\begin{align}\label{eq:frac_Leibniz}
D^{\alpha-1}\left(f \eta_j\right)(x) = \eta_j(x) D^{\alpha-1}(f)(x) + f(x) D^{\alpha-1}(\eta_j)(x) + R^\alpha(\eta_j,f)(x),
\end{align} 
where, 
\begin{align}\label{eq:frac_remainder}
R^\alpha(\eta_j,f)(x) = \int_{\mathbb{R}^d} \left(\eta_j(x+u)-\eta_j(x)\right) \left(f(x+u) - f(x)\right)u\nu_\alpha(du). 
\end{align}
Then, by Minkowski's inequality, for all $k \in \{1,\dots,d\}$ and all $j \geq 1$, 
\begin{align}\label{eq:consequence_triangle_ineq}
\|D^{\alpha-1}_k(\eta_jf)-D^{\alpha-1}_k(f)\|_{L^p(\mathbb{R}^d,dx)} &\leq \|(1-\eta_j)D^{\alpha-1}_k(f)\|_{L^p(\mathbb{R}^d,dx)} + \|f D^{\alpha-1}_k(\eta_j) \nonumber \\ 
&\quad\quad  + R^\alpha_k(\eta_j,f)\|_{L^p(\mathbb{R}^d,dx)}. 
\end{align}
It is clear that the first term appearing on the right-hand side of the previous inequality tends to $0$ as $j$ tends to $+\infty$; this follows from a standard application of the Lebesgue dominated convergence theorem. Next, let us prove that, for all $k \in \{1, \dots, d\}$,  
\begin{align}\label{eq:uniform_boundedness1_general_p}
\sup_{j \geq 1} \|f D^{\alpha-1}_k(\eta_j)\|_{L^p(\mathbb{R}^d,dx)} < + \infty.
\end{align}
First, by scale invariance, for all $j \geq 1$ and all $x \in \mathbb{R}^d$, 
\begin{align*}
D^{\alpha-1}(\eta_j)(x) = \frac{1}{j^{\alpha-1}} D^{\alpha-1}\left(\eta\right)\left(\frac{x}{j}\right). 
\end{align*}
Moreover, since $\eta \in \mathcal{C}^{\infty}_c(\mathbb{R}^d)$, $D^{\alpha-1}(\eta)\in L^r(\mathbb{R}^d,dx)$, for all $r \in [1,+\infty]$.~Now, by H\"older's inequality with $r = d/(d-p(\alpha-1))$ and $q = r/(r-1)$, 
\begin{align*}
\|f D^{\alpha-1}(\eta_j)\|^p_{L^p(\mathbb{R}^d,dx)} & = \int_{\mathbb{R}^d} |f(x)|^p \|D^{\alpha-1}(\eta_j)(x)\|^p dx \\
& = \frac{1}{j^{p(\alpha-1)}} \int_{\mathbb{R}^d} |f(x)|^p \left\|D^{\alpha-1}(\eta)\left(\frac{x}{j}\right)\right\|^p dx \\
& \leq \frac{1}{j^{p(\alpha-1)}} \left(\int_{\mathbb{R}^d}|f(x)|^{p_\alpha^*} dx\right)^{\frac{1}{r}} \left(\int_{\mathbb{R}^d}\left\|D^{\alpha-1}(\eta)\left(\frac{x}{j}\right)\right\|^{pq} dx\right)^{\frac{1}{q}} \\
& \leq \left(\int_{\mathbb{R}^d}|f(x)|^{p_\alpha^*} dx\right)^{\frac{1}{r}} \left(\int_{\mathbb{R}^d}\left\|D^{\alpha-1}(\eta)\left(x\right)\right\|^{pq} dx\right)^{\frac{1}{q}},
\end{align*}
from which \eqref{eq:uniform_boundedness1_general_p} follows. Next, let us prove that, for all $k \in \{1, \dots, d\}$, 
\begin{align}\label{eq:uniform_boundedness_remainder2_general_p}
\sup_{j \geq 1}\left\| R_k^\alpha(\eta_j,f)\right\|_{L^p(\mathbb{R}^d,dx)}<+\infty. 
\end{align}
By Minkowski's integral inequality, for all integer $j \geq 1$, 
\begin{align*}
\left\| \int_{\mathbb{R}^d}\Delta_u(\eta_j)(.)\Delta_u(f)(.) u\nu_\alpha(du)\right\|_{L^p(\mathbb{R}^d,dx)} & \leq \int_{\mathbb{R}^d} \|u\|\nu_\alpha(du) \|\Delta_u(\eta_j)(.)\Delta_u(f)(.)\|_{L^{p}(\mathbb{R}^d,dx)} \\
& \leq \int_{\mathbb{R}^d} \|u\|\nu_\alpha(du) \bigg(\|\Delta_u(\eta_j)(.)f(.+u)\|_{L^{p}(\mathbb{R}^d,dx)} \\
& \quad\quad + \|\Delta_u(\eta_j)(.)f(.)\|_{L^{p}(\mathbb{R}^d,dx)}\bigg) \\
& \leq \frac{1}{j^{\alpha-1}}\int_{\mathbb{R}^d} \|u\|\nu_\alpha(du) \bigg(\bigg\|\Delta_u(\eta)\left(\frac{.}{j}\right)f(.+u)\bigg\|_{L^{p}(\mathbb{R}^d,dx)} \\
& \quad\quad + \bigg\|\Delta_u(\eta)\left(\frac{.}{j}\right)f(.)\bigg\|_{L^{p}(\mathbb{R}^d,dx)}\bigg).
\end{align*}
Again, by H\"older's inequality with $r = d/(d-p(\alpha-1))$ and $q = r/(r-1)=d/(p(\alpha-1))$,
\begin{align*}
\left\|\Delta_u(\eta)\left(\frac{.}{j}\right)f(.)\right\|^p_{L^{p}(\mathbb{R}^d,dx)} & = \int_{\mathbb{R}^d} f(x)^p \left| \eta\left(\frac{x}{j}+u\right) - \eta\left(\frac{x}{j}\right) \right|^pdx\\
& \leq \left(\int_{\mathbb{R}^d} |f(x)|^{pr} dx \right)^{\frac{1}{r}} \left(\int_{\mathbb{R}^d}\left|\eta\left(\frac{x}{j}+u\right) - \eta\left(\frac{x}{j}\right) \right|^{pq}dx\right)^{\frac{1}{q}}\\
& \leq j^{\frac{d}{q}} \left(\int_{\mathbb{R}^d} |f(x)|^{pr} dx \right)^{\frac{1}{r}} \left(\int_{\mathbb{R}^d}\left|\eta\left(x+u\right) - \eta\left(x\right) \right|^{pq}dx\right)^{\frac{1}{q}}.
\end{align*}
Hence, 
\begin{align*}
\left\| \int_{\mathbb{R}^d}\Delta_u(\eta_j)(.)\Delta_u(f)(.) u\nu_\alpha(du)\right\|_{L^p(\mathbb{R}^d,dx)} & \leq
2\frac{j^{\frac{d}{pq}}}{j^{\alpha-1}} \int_{\mathbb{R}^d} \|u\|\nu_{\alpha}(du) \left(\int_{\mathbb{R}^d} |f(x)|^{pr} dx \right)^{\frac{1}{rp}} \\
& \quad\quad \times \left(\int_{\mathbb{R}^d}\left|\eta\left(x+u\right) - \eta\left(x\right) \right|^{pq}dx\right)^{\frac{1}{pq}}.
\end{align*}
Moreover, for all $x \in \mathbb{R}^d$, 
\begin{align*}
\left\| R^\alpha\left(\eta_j, f\right)(x) \right\| & \leq \int_{B(0,1)} |\Delta_u(\eta_j)(x)||\Delta_u(f)(x)| \|u\| \nu_\alpha(du) \\
& \quad\quad\quad + \int_{B(0,1)^c} |\Delta_u(\eta_j)(x)||\Delta_u(f)(x)| \|u\| \nu_\alpha(du) \\
& \leq 2 \|f\|_{\infty} \int_{B(0,1)} |\Delta_u(\eta_j)(x)| \|u\| \nu_\alpha(du) \\
& \quad\quad\quad + 4 \|f\|_{\infty} \|\eta\|_{\infty} \int_{B(0,1)^c} \|u\| \nu_\alpha(du) \\
& \leq 2 \|f\|_{\infty} \| \nabla(\eta_j)\|_{\infty} \int_{B(0,1)} \|u\|^2 \nu_\alpha(du)
\\
& \quad\quad\quad + 4 \|f\|_{\infty} \|\eta\|_{\infty} \int_{B(0,1)^c} \|u\| \nu_\alpha(du) \\
& \leq \frac{2}{j} \|f\|_{\infty} \| \nabla(\eta)\|_{\infty} \int_{B(0,1)} \|u\|^2 \nu_\alpha(du)
\\
& \quad\quad\quad + 4 \|f\|_{\infty} \|\eta\|_{\infty} \int_{B(0,1)^c} \|u\| \nu_\alpha(du),
\end{align*}
where $\Delta_u(f)(x) = f(x+u) - f(x)$, for all $x,u \in \mathbb{R}^d$. Next, for all $x \in \mathbb{R}^d$ and all $u \in \mathbb{R}^d$, 
\begin{align}\label{eq:pointwise_limit_first_order}
\underset{j \rightarrow +\infty}{\lim} (\eta_j(x+u)-\eta_j(x)) = 0. 
\end{align}
Moreover, for all $x\in \mathbb{R}^d$ and all $u \in B(0,1)$, 
\begin{align*}
\left|\Delta_u(\eta_j)(x)\right| \left|\Delta_u(f)(x)\right| \leq 2 \|\nabla(\eta)\|_{\infty} \|f\|_{\infty} \|u\|, 
\end{align*}
and, for all $x\in \mathbb{R}^d$ and all $u \in B(0,1)^c$,
\begin{align*}
\left|\Delta_u(\eta_j)(x)\right| \left|\Delta_u(f)(x)\right| \leq 4 \|\eta\|_{\infty} \|f\|_{\infty}.
\end{align*}
Thus, by the Lebesgue dominated convergence theorem, for all $x \in \mathbb{R}^d$, 
\begin{align}\label{eq:pointwise_limit_frac_remainder-2}
\underset{j \rightarrow +\infty}{\lim} R^\alpha\left(\eta_j,f\right)(x) = 0.
\end{align}
Combining \eqref{eq:uniform_boundedness1_general_p} and \eqref{eq:uniform_boundedness_remainder2_general_p} leads to: 
\begin{align}\label{eq:uniform_boundedness_both_general_p}
\sup_{j \geq 1} \left\|f D^{\alpha-1}(\eta_j) + R^\alpha(\eta_j,f)\right\|_{L^p(\mathbb{R}^d,dx)} < +\infty.
\end{align}
Thanks to \cite[Theorem $3.18$]{B_FA11}, there exists a subsequence $(f D^{\alpha-1}(\eta_{j_k})+R^\alpha\left(\eta_{j_k},f\right))_{k \geq 1}$ which converges weakly to some element of $L^{p}\left(\mathbb{R}^d,\mathbb{R}^d,dx\right)$. Now, for all $\varphi \in \mathcal{S}(\mathbb{R}^d)$ and all $k \in \{1,\dots,d\}$, 
\begin{align*}
\underset{\ell \rightarrow +\infty}{\lim} \langle fD_k^{\alpha-1}\left(\eta_{j_\ell}\right) ; \varphi\rangle = \underset{\ell \rightarrow +\infty}{\lim} \frac{1}{j_\ell^{\alpha-1}}\langle D_k^{\alpha-1}\left(\eta\right)(\frac{.}{j_\ell}) ; f\varphi\rangle = 0,
\end{align*}
where the last equality follows by the Lebesgue dominated convergence theorem. Moreover, combined with \eqref{eq:pointwise_limit_frac_remainder-2}, $(f D^{\alpha-1}(\eta_{j_k})+R^\alpha\left(\eta_{j_k},f\right))_{k \geq 1}$ converges weakly to $0$ in $L^{p}\left(\mathbb{R}^d,\mathbb{R}^d,dx\right)$. Finally, by the Banach-Saks property (see \cite[Proposition $10.8$]{Ponce_book16}), up to a subsequence, 
\begin{align}\label{eq:Banach-Saks}
f D^{\alpha-1}\left(\frac{1}{N}\sum_{k=1}^N\eta_{j_k}\right)+R^\alpha\left(\frac{1}{N}\sum_{k=1}^N\eta_{j_k},f\right) \longrightarrow 0, \quad N \longrightarrow +\infty,
\end{align}
in $L^p(\mathbb{R}^d,\mathbb{R}^d,dx)$. Now, let $(f_N)_{N \geq 1}$ be defined, for all $N \geq 1$ and all $x \in \mathbb{R}^d$, by 
\begin{align}\label{eq:convex_smooth_cut-off}
f_N(x) = \frac{1}{N} \sum_{k=1}^{N}\eta_{j_k}(x)f(x). 
\end{align}
Then, $0 \leq \sum_{k=1}^N \eta_{j_k}(x)/N \leq 1$ and, for all $N \geq 1$, $\sum_{k = 1}^{N} \eta_{j_k}/N \in \mathcal{C}_c^{\infty}\left(\mathbb{R}^d\right)$. Finally, for all $x \in \mathbb{R}^d$, 
\begin{align*}
\frac{1}{N} \sum_{k=1}^N \eta_{j_k}(x) \longrightarrow 1, \quad N \longrightarrow +\infty.
\end{align*} 
This concludes the proof of the theorem.
\end{proof}
\noindent
Next, let us prove Theorem \ref{thm:Meyers-Serrin_theorem} which provides standard properties of the vector space $\dot{W}^{\alpha-1,p}(\mathbb{R}^d,dx)$ endowed with the mapping $|\cdot|_{\alpha-1,p}$ and identifies it with $\mathcal{D}^{\alpha-1,p}(\mathbb{R}^d,dx)$.\\

\noindent
\textit{Proof of Theorem \ref{thm:Meyers-Serrin_theorem}.}
Let us start with the proof of the first statement. Since the operator $D^{\alpha-1}$ is linear and since $\|\cdot\|_{L^p(\mathbb{R}^d,dx)}$ is a norm, we only need to prove: if $f \in L^{p^\star_\alpha}(\mathbb{R}^d,dx)$ is such that $\|D^{\alpha-1}(f)\|_{L^p(\mathbb{R}^d,dx)}=0$, then $f=0$ almost everywhere.  By Theorem \ref{thm:approx_scs_homo_Lp}, let $(f_n)_{n\geq 1}$ be a sequence of smooth and compactly supported functions defined on $\mathbb{R}^d$ such that 
\begin{align*}
D^{\alpha-1}(f_n) \longrightarrow D^{\alpha-1}(f),\quad n \rightarrow+\infty,
\end{align*}
in $L^p(\mathbb{R}^d, \mathbb{R}^d,dx)$, and such that 
\begin{align*}
f_n \longrightarrow f, \quad n \rightarrow+\infty, 
\end{align*}
in $L^{p^*_\alpha}(\mathbb{R}^d,dx)$. Now, by the fractional Sobolev inequality \eqref{ineq:FSI_NDS_full}, for all $n \geq 1$, 
\begin{align*}
\|f_n\|_{L^{p^*_\alpha}(\mathbb{R}^d,dx)} \leq C(\alpha,d,p) \|D^{\alpha-1}(f_n)\|_{L^{p}(\mathbb{R}^d,dx)},
\end{align*}
which implies, 
\begin{align*}
\|f\|_{L^{p^*_\alpha}(\mathbb{R}^d,dx)} \leq C(\alpha,d,p) \|D^{\alpha-1}(f)\|_{L^{p}(\mathbb{R}^d,dx)}.
\end{align*}
Thus, if $f \in \dot{W}^{\alpha-1,p}(\mathbb{R}^d,dx)$ is such that $|f|_{\alpha-1,p}=0$, then $f=0$ almost everywhere. Let $F$ be an element of $\mathcal{D}^{\alpha-1,p}(\mathbb{R}^d,dx)$ and let $(f_m)_{m \geq 1}$ be a representative which is a fundamental sequence of functions in $\mathcal{C}_c^{\infty}(\mathbb{R}^d)$ with respect to the norm $|\cdot|_{\alpha-1,p}$. Thanks to \eqref{ineq:FSI_NDS_full}, $(f_m)_{m\geq 1}$ is then a fundamental sequence in $L^{p_\alpha^*}(\mathbb{R}^d,dx)$ which is complete by the Fischer-Riesz theorem. Let us denote by $f$ its limit. Moreover, since $(f_m)_{m\geq 1}$ is a fundamental sequence with respect to the norm $|\cdot|_{\alpha-1,p}$, there exists $g_\alpha \in L^p(\mathbb{R}^d, \mathbb{R}^d,dx)$ such that $(D^{\alpha-1}(f_m))_{m\geq 1}$ converges to $g_\alpha$ in $L^p(\mathbb{R}^d,\mathbb{R}^d,dx)$ as $m$ tends to $+\infty$. Thus, for all $\varphi \in \mathcal{S}(\mathbb{R}^d)$ and all $k \in \{1, \dots, d\}$, 
\begin{align*}
\langle g_{\alpha,k} ; \varphi \rangle & = \underset{m \rightarrow +\infty}{\lim} \langle D^{\alpha-1}_k(f_m) ; \varphi\rangle \\
& = \underset{m \rightarrow +\infty}{\lim} \langle f_m ; (D^{\alpha-1}_k)^*(\varphi)\rangle \\
& = \langle f ; (D^{\alpha-1}_k)^*(\varphi)\rangle \\
& = \langle D^{\alpha-1}_k(f) ; \varphi\rangle.
\end{align*}
Then, $f$ belongs to $\dot{W}^{\alpha-1,p}(\mathbb{R}^d,dx)$. Let $\mathcal{J}$ be the linear map defined, for all $F \in \mathcal{D}^{\alpha-1,p}(\mathbb{R}^d,dx)$, by 
\begin{align*}
\mathcal{J}(F) = f,
\end{align*}
where $f \in \dot{W}^{\alpha-1,p}(\mathbb{R}^d,dx)$ is the function obtained previously. Note that this map is well-defined (it does not depend on the representative of $F$) and is an isometry between $\mathcal{D}^{\alpha-1,p}(\mathbb{R}^d,dx)$ and $\dot{W}^{\alpha-1,p}(\mathbb{R}^d,dx)$ since, 
\begin{align*}
|F|_{\alpha-1,p} = \underset{m \rightarrow + \infty}{\lim} |f_m|_{\alpha-1,p} = |f|_{\alpha-1,p}.  
\end{align*}
Finally, let $f \in \dot{W}^{\alpha-1,p}(\mathbb{R}^d,dx)$. Then, from Theorem \ref{thm:approx_scs_homo_Lp}, there exists a sequence of functions $(f_n)_{n \geq 1}$ in $\mathcal{C}^{\infty}_c(\mathbb{R}^d)$ such that 
\begin{align*}
\left| f_n - f \right|_{\alpha-1,p} \longrightarrow 0,
\end{align*}
as $n$ tends to $+\infty$. This converging sequence of functions is fundamental. So, 
\begin{align*}
\mathcal{J}\left((f_n)_{n \geq 1}\right) = f.
\end{align*}
This concludes the proof of the theorem.$\qed$\\

\noindent
The next technical lemma extends the refined fractional Sobolev inequality \eqref{ineq:refined_frac_sobolev_ineq-3} to the functional space $\dot{W}^{\alpha-1,2}\left(\mathbb{R}^d,dx\right)$.

\begin{lem}\label{lem:extension_RFSI}
Let $\alpha \in (1,2)$, let $d \geq 2$ be an integer and let $2^*_\alpha=2d/(d-2(\alpha-1))$.~Let $\chi$ be an infinitely differentiable function defined on $\mathbb{R}_+$ with compact support such that $\chi$ is identically equal to $1$ in a neighborhood of the origin. Then, for all $f \in \dot{W}^{\alpha-1,2}\left(\mathbb{R}^d,dx\right)$,
\begin{align}\label{ineq:extension_RFSI}
\|f\|_{L^{2^*_\alpha}(\mathbb{R}^d,dx)} \leq C_2(\alpha,d,\chi) \left(\underset{t>0}{\sup}\, t^{\frac{d-2(\alpha-1)}{4}} \|\chi \left(-t\Delta\right)(f)\|_{\infty} \right)^{1-\frac{2}{2^*_\alpha}} \|D^{\alpha-1}(f)\|^{\frac{2}{2^*_\alpha}}_{L^2(\mathbb{R}^d,dx)},
\end{align} 
where $C_2(\alpha,d,\chi)>0$ depends on $\alpha$, $d$ and $\chi$. 
\end{lem}

\begin{proof}
Let $f \in \dot{W}^{\alpha-1,2}\left(\mathbb{R}^d,dx\right)$. By Theorem \ref{thm:approx_scs_homo_Lp}, there exists a sequence of functions $(f_n)_{n \geq 1}$ in $\mathcal{C}^{\infty}_c(\mathbb{R}^d)$ such that 
\begin{align*}
f_n \longrightarrow f, \quad n \longrightarrow +\infty,  
\end{align*}
in $L^{2^*_\alpha}\left(\mathbb{R}^d,dx\right)$, and 
\begin{align*}
D^{\alpha-1}(f_n) \longrightarrow D^{\alpha-1}(f), \quad n \longrightarrow +\infty,  
\end{align*}
in $L^{2}\left(\mathbb{R}^d,\mathbb{R}^d,dx\right)$. Now, by the Fourier inversion formula, for all $x \in \mathbb{R}^d$, all $t>0$ and all $n \geq 1$,  
\begin{align*}
\chi\left(- t \Delta\right)(f_n)(x) & = \frac{1}{(2\pi)^d} \int_{\mathbb{R}^d} \mathcal{F}(f_n)(\xi) \chi\left(t \|\xi\|^2\right) e^{i \langle x;\xi\rangle} d\xi \\
& = \int_{\mathbb{R}^d} f_n(y) g_t(x-y) dy,
\end{align*}
where,
\begin{align*}
g_t(x) = \frac{1}{(2\pi)^d} \int_{\mathbb{R}^d} \chi\left(t \|\xi\|^2\right) e^{i \langle x;\xi\rangle}d\xi = \frac{1}{t^{\frac{d}{2}}}g_1\left(\frac{x}{\sqrt{t}}\right). 
\end{align*}
Similarly, for all $x \in \mathbb{R}^d$ and all $t>0$, 
\begin{align*}
\chi\left(-t\Delta\right)(f)(x) = \int_{\mathbb{R}^d} f(y) g_t(x-y)dy.
\end{align*}
Note that the previous integral is well-defined since $f \in L^{2^*_\alpha}\left(\mathbb{R}^d,dx\right)$ and $g_1\in L^p(\mathbb{R}^d,dx)$, for all $p \in [1,+\infty]$. Thus, by H\"older's inequality, for all $x \in \mathbb{R}^d$, all $t>0$ and all $n \geq 1$,  
\begin{align*}
\left|\chi\left(-t\Delta\right)(f_n)(x)-\chi\left(-t\Delta\right)(f)(x)\right|& \leq \int_{\mathbb{R}^d} |f_n(y)-f(y)| |g_t(x-y)| dy\\
& \leq \|f_n-f\|_{L^{2_\alpha^*}(\mathbb{R}^d,dx)} \|g_t\|_{L^{(2_\alpha^*)'}(\mathbb{R}^d,dx)} \\
& \leq t^{-\frac{d-2(\alpha-1)}{4}}\|f_n-f\|_{L^{2_\alpha^*}(\mathbb{R}^d,dx)} \|g_1\|_{L^{(2_\alpha^*)'}(\mathbb{R}^d,dx)},
\end{align*}
where $(2_\alpha^*)' =2_\alpha^*/(2_\alpha^*-1)$. This implies that, 
\begin{align*}
\underset{t>0}{\sup}\, t^{\frac{d-2(\alpha-1)}{4}} \|\chi\left(-t\Delta\right)(f_n)-\chi\left(-t\Delta\right)(f)\|_{\infty} \leq \|f_n-f\|_{L^{2_\alpha^*}(\mathbb{R}^d,dx)} \|g_1\|_{L^{(2_\alpha^*)'}(\mathbb{R}^d,dx)}.
\end{align*}
Thus,
\begin{align*}
\underset{n \rightarrow +\infty}{\lim} \underset{t>0}{\sup}\, t^{\frac{d-2(\alpha-1)}{4}} \|\chi\left(-t\Delta\right)(f_n)-\chi\left(-t\Delta\right)(f)\|_{\infty} = 0.
\end{align*}
Now, for all $n \geq 1$, 
\begin{align*}
\|f_n\|_{L^{2^*_\alpha}(\mathbb{R}^d,dx)} \leq C_2(\alpha,d,\chi) \left(\underset{t>0}{\sup}\, t^{\frac{d-2(\alpha-1)}{4}} \|\chi \left(-t\Delta\right)(f_n)\|_{\infty} \right)^{1-\frac{2}{2^*_\alpha}} \|D^{\alpha-1}(f_n)\|^{\frac{2}{2^*_\alpha}}_{L^2(\mathbb{R}^d,dx)}.
\end{align*}
Passing to the limit in the previous inequality concludes the proof of the lemma. 
\end{proof}
\noindent
In the sequel, let us introduce the Morrey spaces; they provide a suitable functional setting for a refined version of the fractional Sobolev inequality \eqref{ineq:FSI_NDS_full} in the general $L^p$ case, with $p \in (1,d/(\alpha-1))$. This approach originates in \cite{PP_14} where refined versions of the classical Sobolev inequality have been obtained based on weighted inequalities for the Riesz potential operator of \cite{SaWhe_AJM92}. For $r \in [1, + \infty)$ and $\beta \in [0,d]$, let $L^{r,\beta}(\mathbb{R}^d,dx)$ be the set of Borel measurable real-valued functions defined on $\mathbb{R}^d$ such that 
\begin{align}\label{eq:Morrey_norm}
\|f\|_{r,\beta} : = \left(\sup_{R>0,x \in \mathbb{R}^d} R^{\beta-d}\int_{B(x,R)} |f(y)|^rdy\right)^{\frac{1}{r}}<+\infty,
\end{align}
where $B(x,R)$ is the Euclidean open ball centered at $x$ and of radius $R$.~It becomes a normed vector space when quotiented by the almost everywhere equality with respect to the $d$-dimensional Lebesgue measure and endowed with the norm \eqref{eq:Morrey_norm}. It is still denoted by $L^{r,\beta}(\mathbb{R}^d,dx)$ in the sequel. Next, let $W$ and $V$ be two non-negative Borel measurable functions defined on $\mathbb{R}^d$. Let $s \in (0,d)$ and $1<r \leq t< +\infty$ be such that, for all cubes $Q \subset \mathbb{R}^d$, 
\begin{align}\label{eq:weights_condition}
|Q|^{\frac{s}{d}+\frac{1}{t}-\frac{1}{r}} \left(\frac{1}{|Q|}\int_Q W(x)^{\tau} dx\right)^{\frac{1}{t \tau}}\left(\frac{1}{|Q|}\int_Q V(x)^{(1-r')\tau} dx\right)^{\frac{1}{r'\tau}} \leq C_{\tau},
\end{align}
where $r' = r/(r-1)$, for some $\tau>1$ and some $C_{\tau}>0$. Then, thanks to \cite[Theorem $1$, $(A)$]{SaWhe_AJM92}, for all $f \in L^r(\mathbb{R}^d,V(x)dx)$ non-negative,
\begin{align}\label{eq:weighted_inequality_Riesz_op}
\left(\int_{\mathbb{R}^d} [I^s(f)(x)]^{t} W(x)dx \right)^{\frac{1}{t}} \leq c C_{\tau} \left(\int_{\mathbb{R}^d} f(x)^r V(x)dx\right)^{\frac{1}{r}},
\end{align} 
for some $c>0$ depending on $d$, $r$ and $t$. 

\begin{thm}\label{thm:refined_frac_Sobolev_Lp}
Let $\alpha \in (1,2)$, let $d \geq 2$ be an integer, let $p \in (1,d/(\alpha-1))$ and let $p_\alpha^* = pd/(d-p(\alpha-1))$.~Let $\nu_\alpha$ be a non-degenerate symmetric $\alpha$-stable L\'evy measure on $\mathbb{R}^d$ such that $\nu_\alpha$ is a $\gamma$-measure with $\gamma \in [1,d]$ and $\gamma-d+2\alpha>1$. Let $t\in[p,p^*_\alpha)$ be such that $t>p_{\alpha}^*-1$. Then, for all $f \in \mathcal{C}^{\infty}_c\left(\mathbb{R}^d\right)$, 
\begin{align}\label{eq:refined_frac_Sobolev_Lp}
\|f\|_{L^{p_\alpha^*}\left(\mathbb{R}^d,dx\right)} \leq \tilde{C}_{\alpha,d,p,\gamma,t} \|f\|^{1-\frac{t}{p_\alpha^*}}_{1, \frac{d-p(\alpha-1)}{p}}  \|D^{\alpha-1}(f)\|_{L^p(\mathbb{R}^d,dx)}^{\frac{t}{p_{\alpha}^*}},
\end{align}
for some $\tilde{C}_{\alpha,d,p,\gamma,t}>0$ depending on $\alpha$, $d$, $p$, $\gamma$ and $t$. Moreover, for all $f \in \dot{W}^{\alpha-1,p}(\mathbb{R}^d,dx)$, 
\begin{align}\label{eq:refined_frac_Sobolev_Lp_full}
\|f\|_{L^{p_\alpha^*}\left(\mathbb{R}^d,dx\right)} \leq \tilde{C}_{\alpha,d,p,\gamma,t} \|f\|^{1-\frac{t}{p_\alpha^*}}_{1, \frac{d-p(\alpha-1)}{p}}  \|D^{\alpha-1}(f)\|_{L^p(\mathbb{R}^d,dx)}^{\frac{t}{p_{\alpha}^*}}.
\end{align}
\end{thm}

\begin{proof}
Let $f \in \mathcal{C}^{\infty}_c(\mathbb{R}^d)$. Then, for all $x \in \mathbb{R}^d$, 
\begin{align*}
|f(x)|^{p_\alpha^*} = |f(x)|^{t} W(x), \quad W(x) = |f(x)|^{p_\alpha^*-t},
\end{align*}
for some $t \in [p, p_\alpha^*)$. Next, using \eqref{ineq:pointwise_estimate_FFTCkernel}, for all $x \in \mathbb{R}^d$,
\begin{align}
\left| f(x) \right| \leq C_{\alpha,d,\gamma} I^{\alpha-1}\left(\|D^{\alpha-1}(f)\|\right)(x),
\end{align} 
for some $C_{\alpha,d,\gamma}>0$ depending on $\alpha$, $d$ and $\gamma$. Then, for all $x \in \mathbb{R}^d$, 
\begin{align*}
|f(x)|^{p_\alpha^*} \leq C_{\alpha,d,\gamma,t} [I^{\alpha-1}\left(\|D^{\alpha-1}(f)\|\right)(x)]^t W(x). 
\end{align*}
In order to apply \eqref{eq:weighted_inequality_Riesz_op}, let us check that the condition \eqref{eq:weights_condition} is satisfied with $r = p$, $s = \alpha-1$, $W(x) = |f(x)|^{p_\alpha^*-t}$ and $V(x) = 1$, for all $x \in \mathbb{R}^d$. To do so, set $\tau = 1/(p_\alpha^*-t)$ and choose $t$ such that $\tau>1$. Then, 
\begin{align*}
|Q|^{\frac{\alpha-1}{d}+\frac{1}{t}-\frac{1}{p}} \left(\frac{1}{|Q|}\int_Q W(x)^{\frac{1}{p_\alpha^*-t}} dx\right)^{\frac{p_\alpha^*-t}{t}} & = |Q|^{\frac{\alpha-1}{d}+\frac{1}{t}-\frac{1}{p}} \left(\frac{1}{|Q|}\int_Q |f(x)| dx\right)^{\frac{p_\alpha^*-t}{t}} \\
& \leq C R^{\alpha-1 + \frac{d}{t}-\frac{d}{p}} \left( \frac{1}{R^d}\int_{B(x,R)} |f(y)| dy\right)^{\frac{p_\alpha^*-t}{t}} \\
& \leq C \left( R^{\frac{d-p(\alpha-1)}{p}
-d}\int_{B(x,R)} |f(y)|dy\right)^{\frac{p_\alpha^*-t}{t}}\\
& \leq C \|f\|^{\frac{p_\alpha^*-t}{t}}_{1,\frac{d-p(\alpha-1)}{p}} = : C_\tau.
\end{align*}
Thus, from \eqref{eq:weighted_inequality_Riesz_op},
\begin{align*}
\int_{\mathbb{R}^d} |f(x)|^{p_\alpha^*}dx & \leq C_{\alpha,d,\gamma,t} \int_{\mathbb{R}^d} [I^{\alpha-1}\left(\|D^{\alpha-1}(f)\|\right)(x)]^t W(x)dx \\
& \leq C_{\alpha,d,\gamma,t} (c C_{\tau})^t \left(\int_{\mathbb{R}^d} \|D^{\alpha-1}(f)(x)\|^p dx\right)^{\frac{t}{p}} \\
& \leq C_{\alpha,d,\gamma,t} (cC)^t \|f\|^{p_\alpha^*-t}_{1, \frac{d-p(\alpha-1)}{p}} \left(\int_{\mathbb{R}^d} \|D^{\alpha-1}(f)(x)\|^p dx\right)^{\frac{t}{p}}.
\end{align*}
Next, let $f \in \dot{W}^{\alpha-1,p}(\mathbb{R}^d,dx)$ and, by Theorem \ref{thm:approx_scs_homo_Lp}, let $(f_n)_{n \geq 1}$ be a sequence of infinitely differentiable and compactly supported functions such that
\begin{align*}
D^{\alpha-1}(f_n) \longrightarrow D^{\alpha-1}(f),\quad n \rightarrow+\infty,
\end{align*}
in $L^p(\mathbb{R}^d, \mathbb{R}^d,dx)$, and such that 
\begin{align*}
f_n \longrightarrow f, \quad n \rightarrow+\infty, 
\end{align*}
in $L^{p^*_\alpha}(\mathbb{R}^d,dx)$. Moreover, by H\"older's inequality, for all $R>0$ and all $x \in \mathbb{R}^d$, 
\begin{align*}
R^{\frac{d-p(\alpha-1)}{p}-d} \int_{B(x,R)} |f(y)|dy \leq C \left(\int_{\mathbb{R}^d} |f(x)|^{p_\alpha^*} dx\right)^{\frac{1}{p_{\alpha}^*}},
\end{align*}
for some $C>0$. Thus, for all $n \geq 1$, 
\begin{align*}
\|f_n-f\|_{1,\frac{d-p(\alpha-1)}{p}} \leq C \|f_n-f\|_{L^{p_\alpha^*}\left(\mathbb{R}^d,dx\right)}. 
\end{align*}
Passing to the limit in \eqref{eq:refined_frac_Sobolev_Lp} concludes the proof of the theorem.
\end{proof}
\noindent 
The next lemma proves a weak version of the Sobolev-Besov-type inequality. 

\begin{lem}\label{lem:weak_Sobolev_Besov_frac_gradient}
Let $d \geq 1$ be an integer, let $\alpha \in (1,2)$, let $\nu_\alpha$ be a non-degenerate symmetric L\'evy measure on $\mathbb{R}^d$ verifying \eqref{eq:scale}, let $\mu_\alpha$ be the corresponding $\alpha$-stable probability measure on $\mathbb{R}^d$ characterized by \eqref{stable:characteristic} and let $p_\alpha$ be its positive Lebesgue density on $\mathbb{R}^d$. Finally, assume that, for all $p \in [1 , +\infty)$, 
\begin{align*}
\left\| \frac{\nabla(p_\alpha)}{p_\alpha} \right\|_{L^p(\mu_\alpha)} < +\infty. 
\end{align*}
Then, for all $1 \leq p < q <+\infty$ and all $f \in \mathcal{S}(\mathbb{R}^d)$, 
\begin{align}\label{ineq:weak_sobolev_besov}
\| f \|_{L^{q,\infty}(\mathbb{R}^d,dx)} \leq C_{\alpha,p,q,d} \|D^{\alpha-1}(f)\|^{\theta}_{L^p(\mathbb{R}^d,dx)} \| f \|^{1-\theta}_{\mathcal{B}_{\infty, \infty}^{s_\alpha,\alpha}} ,
\end{align}
where $\theta =p /q$, $s_\alpha = (\alpha-1) p/(p-q)$ and 
\begin{align}\label{eq:constante_Sobolev_Besov_frac}
C_{\alpha,p,q,d} = \frac{2}{(\alpha-1)^{\frac{p}{q}}} \left\| \frac{\nabla(p_\alpha)}{p_\alpha} \right\|^{\frac{p}{q}}_{L^p(\mu_\alpha)}.
\end{align}
\end{lem}

\begin{proof}
Recall that, for all $q \in (1, +\infty)$, 
\begin{align*}
\| f \|_{L^{q,\infty}(\mathbb{R}^d,dx)} : = \underset{t > 0 }{\sup}\, \left(t \mathcal{L}_d \{x \in \mathbb{R}^d:\, |f(x)| > t \}^{\frac{1}{q}}\right).
\end{align*}
Let $f \in \mathcal{S}(\mathbb{R}^d)$ and let us assume that $ \| f \|_{\mathcal{B}_{\infty, \infty}^{s_\alpha,\alpha}} \leq 1$ by homogeneity. Thus, for all $t >0$ and all $x \in \mathbb{R}^d$, 
\begin{align*}
|P_t^\alpha(f)(x)| \leq t^{\frac{s_\alpha}{\alpha}}. 
\end{align*}
Next, for all $u >0$, set $t(u) = u^{\frac{\alpha}{s_\alpha}}$. Then, by Chebyshev's inequality, for all $u>0$, 
\begin{align*}
u^q \mathcal{L}_d \{x \in \mathbb{R}^d:\, |f(x)| > 2u \} & \leq u^q \mathcal{L}_d \{x \in \mathbb{R}^d:\, |f(x) - P_{t(u)}^\alpha(f)(x) | > u \} \\
& \leq u^{q-p} \int_{\mathbb{R}^d} \left| f(x)-  P_{t(u)}^\alpha(f)(x)  \right|^p dx. 
\end{align*}
From the pseudo-Poincar\'e inequality \eqref{ineq:pseudo_poincar_frac_gradient}, for all $u>0$, 
\begin{align*}
u^q \mathcal{L}_d \{x \in \mathbb{R}^d:\, |f(x)| > 2u \} & \leq u^{q-p} \frac{t(u)^{p - \frac{p}{\alpha}}}{(\alpha-1)^p} \left\| \frac{\nabla(p_\alpha)}{p_\alpha} \right\|^p_{L^p(\mu_\alpha)}  \left\| (D^{\alpha-1})(f) \right\|^p_{L^p(\mathbb{R}^d,dx)} \\
& \leq u^{q - p +  p\frac{\alpha}{s_\alpha}(1-\frac{1}{\alpha})} \frac{1}{(\alpha-1)^p} \left\| \frac{\nabla(p_\alpha)}{p_\alpha} \right\|^p_{L^p(\mu_\alpha)}  \left\| (D^{\alpha-1})(f) \right\|^p_{L^p(\mathbb{R}^d,dx)}.
\end{align*}
Now, taking $s_\alpha = (\alpha-1) p/(p-q)$ gives
\begin{align*}
u \mathcal{L}_d \{x \in \mathbb{R}^d:\, |f(x)| > 2u \}^{\frac{1}{q}} \leq \frac{1}{(\alpha-1)^{\frac{p}{q}}} \left\| \frac{\nabla(p_\alpha)}{p_\alpha} \right\|^{\frac{p}{q}}_{L^p(\mu_\alpha)}  \left\| (D^{\alpha-1})(f) \right\|^{\frac{p}{q}}_{L^p(\mathbb{R}^d,dx)}.
\end{align*}
This concludes the proof of the lemma.
\end{proof}

\begin{rem}\label{rem:weak_Sobolev_inequality}
A ``weak" version of the Sobolev-type inequality \eqref{ineq:strong_Sobolev_Frac} can be obtained from Lemma \ref{lem:weak_Sobolev_Besov_frac_gradient} and the ultracontractive estimate \eqref{ineq:ultracontractivity_stable_heat_Lp}.~Indeed, by the definition of the Besov norm $\| . \|_{\mathcal{B}^{s_\alpha, \alpha}_{\infty, \infty}}$ and the ultracontractive estimate \eqref{ineq:ultracontractivity_stable_heat_Lp}, for all $t>0$, 
\begin{align*}
t^{-s_\alpha/\alpha} \| P_t^\alpha(f)\|_{L^{\infty}(\mathbb{R}^d,dx)} & \leq t^{- \frac{\alpha-1}{\alpha} \frac{p}{p-q} - \frac{d}{\alpha q}} \|p_\alpha\|^{\frac{1}{q}}_{L^\infty(\mathbb{R}^d,dx)} \|f\|_{L^q(\mathbb{R}^d,dx)}\\
& \leq \|p_\alpha\|^{\frac{1}{q}}_{L^\infty(\mathbb{R}^d,dx)} \|f\|_{L^{q}(\mathbb{R}^d,dx)},
\end{align*}
with $q = dp /(d- p(\alpha-1))$.~Thus, for all $p \in [1,+\infty)$ such that $p < d/(\alpha-1)$ and all $f \in \mathcal{S}(\mathbb{R}^d)$, 
\begin{align}\label{ineq:weak_Sobolev_frac_gradient}
\|f\|_{L^{q, \infty}(\mathbb{R}^d,dx)} \leq  C_{\alpha,p,q,d} \|p_\alpha\|^{\frac{1-\theta}{q}}_{L^{\infty}(\mathbb{R}^d,dx)} \|D^{\alpha-1}(f)\|^{\theta}_{L^p(\mathbb{R}^d,dx)} \| f \|^{1-\theta}_{L^q(\mathbb{R}^d,dx)}, 
\end{align}
where $C_{\alpha, p,q,d}$ is given by \eqref{eq:constante_Sobolev_Besov_frac}.
\end{rem}
\noindent 
Next, let us study the consequences of the pseudo-Poincar\'e inequality \eqref{ineq:pseudo_poincar_gradient}.~Note that since it brings into play the local gradient operator, the strong Sobolev-type inequality and its refined versions are easier to deal with.~Let us start with an analogue of Lemma \ref{lem:weak_Sobolev_Besov_frac_gradient}. 

\begin{lem}\label{lem:weak_Sobolev_Besov_local_gradient}
Let $\alpha \in (1,2)$ and let $\nu_\alpha$ be a non-degenerate symmetric L\'evy measure on $\mathbb{R}^d$, $d \geq 1$, verifying \eqref{eq:scale} with spectral measure $\lambda_1$.  Let $\sigma_\alpha$ be defined, for all $\xi \in \mathbb{R}^d$, by 
\begin{align*}
\sigma_\alpha(\xi) = \left(\int_{\mathbb{S}^{d-1}} |\langle y ; \xi\rangle|^\alpha \lambda_1(dy)\right)^{\frac{1}{\alpha}},
\end{align*}
and let $Y_\alpha$ be a symmetric $\alpha$-stable random variable with characteristic function given by $\exp \left( - |\xi|^\alpha \right)$, for all $\xi \in \bbr$.  Then, for all $(p,q)$ such that $p \in [1, \alpha) \cap [1,q) $ and all $f \in \mathcal{S}(\mathbb{R}^d)$, 
\begin{align}
\|f\|_{L^{q, \infty}(\mathbb{R}^d,dx)} \leq 2 \left( \bbe |Y_\alpha|^p \right)^{\frac{1}{q}} \| \sigma_\alpha\left(\nabla(f)\right) \|^{\theta}_{L^p(\mathbb{R}^d,dx)} \| f  \|^{1-\theta}_{\mathcal{B}^{\theta/(\theta-1),\alpha}_{\infty,\infty}},
\end{align}
with $\theta = p/q$. 
\end{lem}

\begin{proof}
The proof is completely similar to the one of Lemma \ref{lem:weak_Sobolev_Besov_frac_gradient}.~Recall that from Chebyshev's inequality, for all $u>0$, 
\begin{align*}
u^q \mathcal{L}_d \{x \in \mathbb{R}^d:\, |f(x)| > 2u \} & \leq u^q \mathcal{L}_d \{x \in \mathbb{R}^d:\, |f(x) - P_{t(u)}^\alpha(f)(x) | > u \} \\
& \leq u^{q-p} \int_{\mathbb{R}^d} \left| f(x)-  P_{t(u)}^\alpha(f)(x)  \right|^p dx,
\end{align*}
where $t(u) = u^{\frac{\alpha}{s}}$, and $f \in \mathcal{S}(\mathbb{R}^d)$ such that $\|f\|_{\mathcal{B}_{\infty,\infty}^{s,\alpha}} \leq 1$ with $s<0$ to be chosen in the sequel.  Then, from the pseudo-Poincar\'e inequality \eqref{ineq:pseudo_poincar_gradient},  for all $u>0$, 
\begin{align*}
u^q \mathcal{L}_d \{x \in \mathbb{R}^d:\, |f(x)| > 2u \} & \leq u^{q-p} t(u)^{\frac{p}{\alpha}} \bbe |Y_\alpha|^p \| \sigma_\alpha\left(\nabla(f)\right) \|^p_{L^p(\mathbb{R}^d,dx)} \\
& \leq u^{q-p + \frac{p}{s}}  \bbe |Y_\alpha|^p \| \sigma_\alpha\left(\nabla(f)\right) \|^p_{L^p(\mathbb{R}^d,dx)}.
\end{align*}
Taking $s = p/(p-q) = \theta/(\theta - 1)$, for all $u >0$, 
\begin{align*}
u^q \mathcal{L}_d \{x \in \mathbb{R}^d:\, |f(x)| > 2u \} & \leq \bbe |Y_\alpha|^p \| \sigma_\alpha\left(\nabla(f)\right) \|^p_{L^p(\mathbb{R}^d,dx)}.
\end{align*}
The conclusion easily follows.
\end{proof}
\noindent
Under the non-degeneracy condition \eqref{eq:non_deg}, the function $\sigma_\alpha$ is a norm on $\mathbb{R}^d$. Then, it is natural to retrieve the couple $(\theta,s)$ associated with the classical heat semigroup.~As a straightforward consequence of the truncation technique contained in \cite[Theorem $1$]{MLedoux_03} one obtains the following strong version of Lemma \ref{lem:weak_Sobolev_Besov_local_gradient}.~The proof is given here for the sake of completeness.

\begin{prop}\label{prop:refinement_sobolev_besov_gradient}
Let $\alpha \in (1,2)$ and let $\nu_\alpha$ be a non-degenerate symmetric L\'evy measure on $\mathbb{R}^d$, $d \geq 1$, verifying \eqref{eq:scale} with spectral measure $\lambda_1$.~Let $\sigma_\alpha$ be defined, for all $\xi \in \mathbb{R}^d$, by 
\begin{align*}
\sigma_\alpha(\xi) = \left(\int_{\mathbb{S}^{d-1}} |\langle y ; \xi\rangle|^\alpha \lambda_1(dy)\right)^{\frac{1}{\alpha}},
\end{align*}
and let $Y_\alpha$ be a symmetric $\alpha$-stable random variable with characteristic function given by $\exp \left( - |\xi|^\alpha \right)$, for all $\xi \in \bbr$.  Then, for all $(p,q)$ such that $p \in [1, \alpha) \cap [1,q) $ and all $f \in W^{1,p}(\mathbb{R}^d,dx)$, 
\begin{align}\label{ineq:refinement_sobolev_besov_gradient}
\|f\|_{L^{q}(\mathbb{R}^d,dx)} \leq C_{\alpha,p,q} \| \sigma_\alpha\left(\nabla(f)\right) \|^{\theta}_{L^p(\mathbb{R}^d,dx)} \| f  \|^{1-\theta}_{\mathcal{B}^{\theta/(\theta-1),\alpha}_{\infty,\infty}},
\end{align}
with $\theta = p/q$ and with $C_{\alpha,p,q}>0$ depending only on $\alpha,p$ and $q$.
\end{prop}

\begin{proof}
The proof is an adaptation of an argument contained in \cite[Proof of Theorem $1$, Steps $2$ and $3$]{MLedoux_03}.  First, let us assume that the function $f$ belongs to $L^q(\mathbb{R}^d,dx)$ as well.~As previously, assume that $\| f \|_{\mathcal{B}^{s_\alpha,\alpha}_{\infty,\infty}}\leq 1$, for some $s_\alpha<0$ to be chosen later on and let $t(u) = u^{\alpha/s_\alpha}$,  for all $u>0$.  Now, let us introduce a truncation parameter $R \geq 10$.  By the layer cake representation formula, 
\begin{align*}
\frac{1}{5^q}\|f\|^q_{L^q(\mathbb{R}^d,dx)} = \int_{0}^{+\infty} \mathcal{L}_d\{ x\in \mathbb{R}^d: \, |f(x)| \geq 5 u \} qu^{q-1} du.
\end{align*}
For all $u >0$, let $f_{R,u}$ be the truncated function defined, for a.e. $x \in \mathbb{R}^d$, by
\begin{align*}
f_{R,u}(x)  = (f(x)- u )_+ \wedge ((R-1)u) +(f(x)+ u )_{-} \vee (-(R-1)u),
\end{align*}
where $(c)_+ = \max (c, 0)$ and $(c)_{-} = \min(c, 0)$.~It is clear that $\{ x \in \mathbb{R}^d, \,  |f(x)| \geq 5 u\} \subset \{ x \in \mathbb{R}^d, \,  |f_{R,u}(x)| \geq 4 u\}$. Then,
\begin{align}\label{ineq:union_bound}
\int_0^{+\infty}\mathcal{L}_d\{ x\in \mathbb{R}^d: \, |f(x)| \geq 5 u \} qu^{q-1} du & \leq \int_0^{+\infty}\mathcal{L}_d\{ x\in \mathbb{R}^d: \, |f_{R,u}(x)| \geq 4 u \} qu^{q-1} du  \nonumber \\
& \leq \int_0^{+\infty} \mathcal{L}_d\{ x\in \mathbb{R}^d: \, | f_{R,u}(x) - P^\alpha_{t(u)}(f_{R,u})(x) | \geq u \}q u^{q-1} du \nonumber \\
&\quad + \int_0^{+\infty}\mathcal{L}_d\{ x\in \mathbb{R}^d: \,  P^\alpha_{t(u)}(|f - f_{R,u}|)(x) \geq 2 u \} q u^{q-1} du,
\end{align}
where we have used the triangle inequality together with the fact that $\| f \|_{\mathcal{B}^{s_\alpha,\alpha}_{\infty,\infty}}\leq 1$ to get, for a.e. $x \in \mathbb{R}^d$, 
\begin{align*}
| f_{R,u}(x) | \leq |  f_{R,u}(x) - P^\alpha_{t(u)}(f_{R,u})(x)  | + P^\alpha_{t(u)}(| f - f_{R,u}| )(x) + u.
\end{align*}
Now, using the pseudo-Poincar\'e inequality \eqref{ineq:pseudo_poincar_gradient} to deal with the first term on the right-hand side of \eqref{ineq:union_bound}, 
\begin{align*}
\mathcal{L}_d\{ x\in \mathbb{R}^d: \, | f_{R,u}(x) - P^\alpha_{t(u)}(f_{R,u})(x) | \geq u \} & \leq u^{-p} \int_{\mathbb{R}^d} | f_{R,u}(x) - P^\alpha_{t(u)}(f_{R,u})(x) |^p dx\\
& \leq u^{-p} (t(u))^{\frac{p}{\alpha}} \bbe |Y_\alpha|^p \| \sigma_\alpha \left( \nabla(f_{R,u}) \right)\|^p_{L^p(\mathbb{R}^d,dx)}\\
& \leq u^{\frac{p}{s_\alpha} - p } \bbe |Y_\alpha|^p \| \sigma_\alpha \left( \nabla(f_{R,u}) \right)\|^p_{L^p(\mathbb{R}^d,dx)}\\
& \leq u^{-q} \bbe |Y_\alpha|^p \| \sigma_\alpha \left( \nabla(f_{R,u}) \right)\|^p_{L^p(\mathbb{R}^d,dx)},
\end{align*}
where we have chosen $s_\alpha = p/(p-q)$ in the last line. Now, by the very definition of $f_{R,u}$, the gradient of $f_{R,u}$ is dominated by the one of $f$ on $\{x \in \mathbb{R}^d:\, u \leq |f(x)| \leq Ru\}$ and it is null otherwise.~More precisely, the gradient of $f_{R,u}$ is null on $\{x  \in \mathbb{R}^d:\,  |f(x)| > Ru \}$ and on $\{x  \in \mathbb{R}^d:\,  |f(x)| < u\}$ and it is equal to the gradient of $f$ on $\{x \in \mathbb{R}^d:\, u \leq |f(x)| \leq Ru\}$.~Thus, 
\begin{align*}
\mathcal{L}_d\{ x\in \mathbb{R}^d: \, | f_{R,u}(x) - P^\alpha_{t(u)}(f_{R,u})(x) | \geq u \} \leq u^{-q} \bbe |Y_\alpha|^p \int_{u \leq |f(x)| \leq Ru} |\sigma_\alpha \left( \nabla(f)(x) \right)|^p dx. 
\end{align*}
Integrating with respect to $q u^{q-1} du$ gives
\begin{align*}
\int_0^{+\infty} \mathcal{L}_d\{ x\in \mathbb{R}^d: \, | f_{R,u}(x) - P^\alpha_{t(u)}(f_{R,u})(x) | \geq u \}& d(u^q) \leq q \bbe |Y_\alpha|^p \\
&\quad \times \int_0^{+\infty} \left( \int_{u \leq |f(x)| \leq Ru} |\sigma_\alpha \left( \nabla(f)(x) \right)|^p dx \right) \frac{du}{u}. 
\end{align*}
Fubini's theorem ensures that, 
\begin{align*}
\int_0^{+\infty} \mathcal{L}_d\{ x\in \mathbb{R}^d: \, | f_{R,u}(x) - P^\alpha_{t(u)}(f_{R,u})(x) | \geq u \} d(u^q) \leq q \bbe |Y_\alpha|^p \log(R) \int_{\mathbb{R}^d}  |\sigma_\alpha \left( \nabla(f)(x) \right)|^p dx. 
\end{align*}
Now, the semigroup $(P^\alpha_t)_{t \geq 0}$ is positivity preserving and mass conservative (thanks to \eqref{eq:StheatSM}) and the $d$-dimensional Lebesgue measure is an invariant measure for $(P^\alpha_t)_{t \geq 0}$ (by translation invariance). Then, 
\begin{align*}
 \int_0^{+\infty}\mathcal{L}_d\{ x\in \mathbb{R}^d: \,  P^\alpha_{t(u)}(|f - f_{R,u}|)(x) \geq 2 u \} d(u^q) \leq \frac{q}{q-1} \frac{1}{R^{q-1}} \|f\|^q_{L^q(\mathbb{R}^d,dx)}. 
\end{align*}
Taking $R$ large enough ensures that 
\begin{align*}
\|f\|_{L^q(\mathbb{R}^d,dx)} \leq C_{\alpha, p,q}  \| \sigma_\alpha \left( \nabla(f)\right) \|_{L^p(\mathbb{R}^d,dx)},
\end{align*}
for some $C_{\alpha,p,q}>0$ only depending on $\alpha,p$ and $q$. The end of the proof then follows as in \cite[Proof of Theorem $1$, Step $3$]{MLedoux_03}. 
\end{proof}

\begin{rem}\label{rem:dim_free_sob_gagliardo_nirenberg}
(i) At this point, it is interesting to note that the positive constant $C_{\alpha,p,q}$ in \eqref{ineq:refinement_sobolev_besov_gradient} does not depend on the dimension.~This follows from the proof of Proposition \ref{prop:refinement_sobolev_besov_gradient} and since the constant in the inequality \eqref{ineq:pseudo_poincar_gradient} is dimension-free. \\
(ii) As a straightforward consequence of this refined inequality together with the definition of the Besov-type space $\mathcal{B}^{s, \alpha}_{\infty,\infty}$, for $s<0$, and the ultracontractive estimate \eqref{ineq:ultracontractivity_stable_heat_Lp}, for all $f \in W^{1,p}(\mathbb{R}^d,dx)$, 
\begin{align}\label{ineq:strong_sobolev_inequality_gradient}
\|f\|_{L^q(\mathbb{R}^d,dx)} \leq \left(C_{\alpha,p,q} \|p_\alpha\|^{\frac{1-\theta}{q}}_{L^\infty(\mathbb{R}^d,dx)}   \right)^{\frac{1}{\theta}}  \| \sigma_\alpha(\nabla(f)) \|_{L^p(\mathbb{R}^d,dx)},
\end{align}
with $\theta = p /q$, $q = dp/(d-p)$, $p \in [1, \alpha\wedge q)$ such that $d>p$ ($d \geq 2$) and $C_{\alpha,p,q}>0$ the same constant as in the inequality \eqref{ineq:refinement_sobolev_besov_gradient}.\\
(iii) Similarly, for all $f \in \mathcal{S}(\mathbb{R}^d)$
\begin{align}\label{ineq:strong_GN_inequality_gradient}
\|f\|_{L^q(\mathbb{R}^d,dx)} \leq C_{\alpha,p,q} \|p_\alpha\|^{\frac{1-\theta}{r}}_{L^{\infty}(\mathbb{R}^d,dx)} \|\sigma_\alpha(\nabla(f))\|^\theta_{L^p(\mathbb{R}^d,dx)} \|f\|^{1-\theta}_{L^r(\mathbb{R}^d,dx)},  
\end{align}
with $\theta = p/q$, $1/q = 1/p - r/(dq)$, $p \in [1, \alpha\wedge q)$ and $C_{\alpha,p,q}>0$ the same constant as in \eqref{ineq:refinement_sobolev_besov_gradient}. 
\end{rem}

\section{Fractional and classical isoperimetric inequalities}\label{sec:frac_class_iso_ineq}
\noindent
In this section, we want to explore the natural notions of perimeter associated with the semigroup $(P^{\alpha}_t)_{t \geq 0}$ and investigate consequences at the level 
of isoperimetric-type inequalities, co-area formulas and strong Sobolev embeddings.~In the classical Euclidean case, there is an equivalence between the strong Sobolev inequality and the isoperimetric inequality of purely geometric nature.~The key to this equivalence is the co-area formula.~We expect that parts of the classical theory will result as a particular case of the one developed here and based on the non-degenerate symmetric $\alpha$-stable 
probability measures with $\alpha \in (1,2)$.~In \cite{De_Giorgi}, the author introduced the following notion of perimeter based on the classical heat semigroup: for all Borel measurable set $E$ with finite $d$-dimensional Lebesgue measure,
\begin{align}\label{eq:perimeter_classical}
\mathcal{P}(E) =  \underset{t \rightarrow 0^+}{\lim} \| \nabla P^H_t (\bbone_E)\|_{L^1(\mathbb{R}^d,dx)} = \underset{t>0}{\sup}  \| \nabla P^H_t (\bbone_E)\|_{L^1(\mathbb{R}^d,dx)}, 
\end{align}
where $(P^H_t)_{t \geq 0}$ is the heat semigroup with generator $\Delta$ given by \eqref{eq:classical_heat_sg}. The fact that the previous limit exists (finite or infinite) is based on the commutation relation and the contraction property of $(P^H_t)_{t \geq 0}$.~Moreover, this definition of perimeter coincides with the one based on the variational formula which is given, for all Borel measurable set $E$ with finite $d$-dimensional Lebesgue measure, by
\begin{align}\label{eq:variational_representation_formula}
\mathcal{P}(E) = \underset{\phi \in \Phi }{\sup} \int_{\mathbb{R}^d} \bbone_E(x) \operatorname{div}(\phi)(x)dx , 
\end{align}
where $\Phi = \{ \phi \in \mathcal{C}_c^1(\mathbb{R}^d , \mathbb{R}^d), \, \|\phi\|_{\infty} \leq 1 \}$.~Thus, it is natural to consider the notions of perimeter \eqref{eq:classical} and \eqref{eq:perim_frac} in our context.~Recall also that the following commutation relations hold true:~for all $f$ smooth enough on $\mathbb{R}^d$, all $t \geq0$ and all $x \in \mathbb{R}^d$, 
\begin{align*}
\nabla P^\alpha_t(f)(x) = P^{\alpha}_t(\nabla(f))(x) , \quad D^{\alpha-1} P^\alpha_t(f)(x) = P^\alpha_t(D^{\alpha-1}(f))(x). 
\end{align*}
\noindent
Before moving on, let us discuss the reasoning behind the notation $\mathcal{P}_{\operatorname{frac}}$ for \eqref{eq:perim_frac}.~From standard scaling arguments, for all $\lambda>0$, all $E\in \mathcal{B}(\mathbb{R}^d)$ with $\mathcal{L}_d(E)<+\infty$, and all $x \in \mathbb{R}^d$, 
\begin{align*}
P^\alpha_t\left(\bbone_{\lambda E}\right)(\lambda x) = P^\alpha_{\frac{t}{\lambda^\alpha}}(\bbone_E)(x).
\end{align*}
Then, using the scaling property of the L\'evy measure $\nu_\alpha$, 
\begin{align*}
\left\| D^{\alpha-1} P^\alpha_{t}(\bbone_{\lambda E}) \right\|_{L^1(\mathbb{R}^d,dx)} = \lambda^{d- (\alpha-1)} \left\| D^{\alpha-1} P^\alpha_{\frac{t}{\lambda^\alpha}}(\bbone_E) \right\|_{L^1(\mathbb{R}^d,dx)},
\end{align*}
which implies
\begin{align*}
\mathcal{P}_{\operatorname{frac}}(\lambda E) = \lambda^{d-(\alpha-1)} \mathcal{P}_{\operatorname{frac}}(E). 
\end{align*}
Moreover, this set function is invariant by translation:~namely, $\mathcal{P}_{\operatorname{frac}}(x+E) = \mathcal{P}_{\operatorname{frac}}(E)$, for all $x \in \mathbb{R}^d$ and all $E \in \mathcal{B}(\mathbb{R}^d)$.~These two facts justify the terminology of fractional perimeter for $\mathcal{P}_{\operatorname{frac}}$.~Now, let us discuss the existence of $\| \sigma_\alpha \left( \nabla P^\alpha_t (\bbone_E) \right) \|_{L^1(\mathbb{R}^d,dx)}$ and $\| D^{\alpha-1} P^\alpha_t (\bbone_E)\|_{L^1(\mathbb{R}^d,dx)}$, for all $t>0$.~This is the purpose of the first technical lemma.

\begin{lem}\label{lem:well_defined_perimeters}
Let $\alpha \in (1,2)$, let $\mu_\alpha$ be a non-degenerate symmetric $\alpha$-stable probability measure on $\mathbb{R}^d$, $d \geq 1$, and let $(P^\alpha_t)_{t\geq 0}$ be the associated stable heat semigroup.~Then, for all $E\in \mathcal{B}(\mathbb{R}^d)$ with $\mathcal{L}_d(E)<+\infty$ and all $t>0$, 
\begin{align}\label{ineq:bound_non-asymp_pcl}
\| \sigma_\alpha \left( \nabla P^\alpha_t (\bbone_E) \right)\|_{L^1(\mathbb{R}^d,dx)} \leq \frac{1}{t^{\frac{1}{\alpha}}} \mathcal{L}_d(E) \left( \int_{\mathbb{R}^d} \sigma_\alpha(\nabla(p_\alpha)(z))  dz \right),
\end{align}
\begin{align}\label{ineq:bound_non-asymp_pfrac}
\| D^{\alpha-1} P^\alpha_t (\bbone_E)\|_{L^1(\mathbb{R}^d,dx)} \leq  \frac{1}{t^{1-\frac{1}{\alpha}}} \mathcal{L}_d(E) \left( \int_{\mathbb{R}^d} \| z\| \mu_\alpha(dz) \right), 
\end{align}
where $p_\alpha$ is the positive Lebesgue density of $\mu_\alpha$. 
\end{lem}
\begin{proof}
Recalling the proof of Proposition \eqref{prop:pseudo_poincare}, for all $g \in \cap_{p \geq 1} L^p(\mathbb{R}^d,dx)$, all $t>0$ and all $x \in \mathbb{R}^d$, 
\begin{align*}
\nabla P_t^\alpha(g)(x) =  \frac{1}{t^{\frac{1}{\alpha}}} \int_{\mathbb{R}^d} g(z) \nabla p_\alpha \left( \frac{x - z}{t^{\frac{1}{\alpha}}} \right) \frac{dz}{t^{\frac{d}{\alpha}}}. 
\end{align*}
Thus, for all $t>0$ and all $x \in \mathbb{R}^d$, 
\begin{align*}
\sigma_\alpha \left( \nabla P_t^\alpha(g)(x) \right) \leq \frac{1}{t^{\frac{1}{\alpha}}} \int_{\mathbb{R}^d} |g(z)| \sigma_\alpha\left(\nabla p_\alpha \left( \frac{x - z}{t^{\frac{1}{\alpha}}} \right) \right)\frac{dz}{t^{\frac{d}{\alpha}}}.
\end{align*}
Hence, for all $t>0$,
\begin{align*}
\int_{\mathbb{R}^d} \sigma_\alpha \left( \nabla P_t^\alpha(g)(x) \right) dx & \leq \frac{1}{t^{\frac{1}{\alpha}}} \int_{\mathbb{R}^d} \int_{\mathbb{R}^d} | g(z) |  \frac{1}{t^{\frac{d}{\alpha}}} \sigma_\alpha\left(\nabla p_\alpha \left( \frac{x - z}{t^{\frac{1}{\alpha}}} \right) \right) dzdx \\
& \leq \frac{1}{t^{\frac{1}{\alpha}}}  \left(\int_{\mathbb{R}^d} | g(z) | dz \right) \left(\int_{\mathbb{R}^d} \sigma_\alpha\left(\nabla p_\alpha \left(x\right) \right) dx\right),
\end{align*}
and so \eqref{ineq:bound_non-asymp_pcl} easily follows.~For \eqref{ineq:bound_non-asymp_pfrac}, first note that, for all $x \in \mathbb{R}^d$, all $t>0$ and all $g \in \cap_{p \geq 1} L^p(\mathbb{R}^d,dx)$, 
\begin{align*}
D^{\alpha - 1} P_t^\alpha(g)(x) = \frac{1}{t^{1 - \frac{1}{\alpha}}} \int_{\mathbb{R}^d} g(z) D^{\alpha - 1}(p_\alpha)\left(\frac{x-z}{t^{\frac{1}{\alpha}}}\right) \frac{dz}{t^{\frac{d}{\alpha}}}. 
\end{align*}
Thus, by a similar reasoning, for all $t>0$,
\begin{align*}
\int_{\mathbb{R}^d} \left\| D^{\alpha - 1} P_t^\alpha(g)(x) \right\| dx & \leq \frac{1}{t^{1- \frac{1}{\alpha}}} \left(\int_{\mathbb{R}^d} | g(z) | dz \right) \left(\int_{\mathbb{R}^d} \left\| D^{\alpha-1} (p_\alpha) \left( x \right) \right\| dx\right). 
\end{align*}
This concludes the proof of the lemma. 
\end{proof}
\noindent
Note that the two limits in \eqref{eq:classical} and \eqref{eq:perim_frac} are well-defined thanks to the commutation relations, the contraction property of $(P^\alpha_t)_{t \geq 0}$ on $L^1(\mathbb{R}^d,dx)$ and the fact that $\sigma_\alpha$ is a norm on $\mathbb{R}^d$.~Based on these definitions of perimeter together with the pseudo-Poincar\'e inequalities for $(P_t^{\alpha})_{t \geq 0}$, it is clear that, for all $t>0$ and all $E$ Borel measurable subset of $\mathbb{R}^d$ with finite $d$-dimensional Lebesgue measure, 
\begin{align*}
\|P_t^\alpha\left( \bbone_E\right) - \bbone_E\|_{L^1(\mathbb{R}^d,dx)} \leq \frac{t^{1 - \frac{1}{\alpha}}}{\alpha-1} \left\| \frac{\nabla(p_\alpha)}{p_\alpha} \right\|_{L^1(\mu_\alpha)} \mathcal{P}_{\operatorname{frac}}(E),
\end{align*}
and, 
\begin{align*}
\|P_t^\alpha\left( \bbone_E\right) - \bbone_E\|_{L^1(\mathbb{R}^d,dx)} \leq t^{\frac{1}{\alpha}} \left(\bbe |Y_\alpha| \right) \mathcal{P}_{\operatorname{cl}}(E).
\end{align*}
In order to pursue our investigation of these new perimeters,~let us lower bound the quantity $\|P_t^\alpha\left( \bbone_E\right) - \bbone_E\|_{L^1(\mathbb{R}^d,dx)}$.~From the convolution structure and the symmetry of the semigroup $(P_t^{\alpha})_{t \geq 0}$, it is clear that:~for all $t \geq 0$ and all $E$ Borel measurable subset of $\mathbb{R}^d$ with finite $d$-dimensional Lebesgue measure, 
\begin{align}\label{eq:easy_one}
\|P_t^\alpha\left( \bbone_E\right) - \bbone_E\|_{L^1(\mathbb{R}^d,dx)} = 2 \left( \mathcal{L}_d(E) - \int_E P^{\alpha}_t(\bbone_E)(x)dx \right).
\end{align}
Finally, the next step before getting the (non-sharp) isoperimetric inequalities is to upper bound
\begin{align}\label{def:A_t}
A_t : =  \int_E P^{\alpha}_t(\bbone_E)(x)dx = \langle P^{\alpha}_t(\bbone_E) ; \bbone_E \rangle_{L^2(\mathbb{R}^d,dx)}, \quad t>0.
\end{align}
This is the purpose of the next technical lemma.

\begin{lem}\label{lem:ultra}
Let $\alpha \in (1,2)$, let $\mu_\alpha$ be a non-degenerate symmetric $\alpha$-stable probability measure on $\mathbb{R}^d$, $d \geq 1$, with positive Lebesgue density $p_\alpha$ and let $(P^\alpha_t)_{t\geq 0}$ be the associated stable heat semigroup.~Then, for all $E\in \mathcal{B}(\mathbb{R}^d)$ with $\mathcal{L}_d(E)<+\infty$ and all $t>0$, 
\begin{align*}
 \int_E P^{\alpha}_t(\bbone_E)(x)dx \leq 2^{\frac{d}{\alpha}}  \dfrac{\|p_\alpha\|_{\infty}}{t^{\frac{d}{\alpha}}} \mathcal{L}_d(E)^2.
\end{align*}
\end{lem}

\begin{proof}
To start,  recall that, for all $f \in L^1(\mathbb{R}^d,dx)$ and all $t>0$, 
\begin{align*}
\left\| P_t^{\alpha}(f) \right\|_{L^\infty(\mathbb{R}^d,dx)} \leq \frac{\|p_\alpha\|_{\infty}}{t^{\frac{d}{\alpha}}} \|f\|_{L^1(\mathbb{R}^d,dx)} , \quad \left\| P_t^{\alpha}(f) \right\|_{L^1(\mathbb{R}^d,dx)} \leq \|f\|_{L^1(\mathbb{R}^d,dx)}.
\end{align*}
Thus, for all $t>0$, 
\begin{align}\label{ineq:L2_L1}
\| P_t^{\alpha}(f)  \|_{L^2(\mathbb{R}^d,dx)} \leq \dfrac{\sqrt{\|p_\alpha\|_{\infty}}}{t^{\frac{d}{2\alpha}}}\|f\|_{L^1(\mathbb{R}^d,dx)}.
\end{align}
Now, by the semigroup property and symmetry, for all $t \geq 0$, 
\begin{align*}
A_t  = \langle P^\alpha_t(\bbone_E) ;  \bbone_E \rangle_{L^2(\mathbb{R}^d,dx)} = \langle P^\alpha_{\frac{t}{2}}(\bbone_E); P^\alpha_{\frac{t}{2}}(\bbone_E) \rangle_{L^2(\mathbb{R}^d,dx)} = \left \| P^\alpha_{\frac{t}{2}}(\bbone_E) \right\|^2_{L^2(\mathbb{R}^d,dx)}.
\end{align*}
Then, \eqref{ineq:L2_L1} implies that, for all $t>0$, 
\begin{align*}
A_t \leq 2^{\frac{d}{\alpha}} \dfrac{\|p_\alpha\|_{\infty}}{t^{\frac{d}{\alpha}}} \mathcal{L}_d(E)^2.
\end{align*}
This concludes the proof of the lemma.
\end{proof}
\noindent
Hence, combining Lemma \ref{lem:ultra} together with \eqref{eq:easy_one} leads to the following lower bound, for all $t>0$ and all $E\in\mathcal{B}(\mathbb{R}^d)$ with $\mathcal{L}_d(E)<+\infty$, 
\begin{align*}
\|P_t^\alpha\left( \bbone_E\right) - \bbone_E\|_{L^1(\mathbb{R}^d,dx)} \geq 2 \left( \mathcal{L}_d(E) - 2^{\frac{d}{\alpha}}  \dfrac{\|p_\alpha\|_{\infty}}{t^{\frac{d}{\alpha}}} \mathcal{L}_d(E)^2 \right).
\end{align*}
We are now ready to prove the isoperimetric inequalities associated with the perimeters $\mathcal{P}_{\operatorname{frac}}(E)$ and $ \mathcal{P}_{\operatorname{cl}}(E)$ of Theorem \ref{thm:isoperimetric_type}.\\

\noindent
\textit{Proof of Theorem \ref{thm:isoperimetric_type}.}
Let $E$ be a Borel measurable subset of $\mathbb{R}^d$ with finite Lebesgue measure such that $0<\mathcal{P}_{\operatorname{frac}}(E)<+\infty$ and $0<\mathcal{P}_{\operatorname{cl}}(E)<+\infty$.~Without loss of generality, one can always assume that these perimeters are different of $0$:~indeed, if $\mathcal{P}_{\operatorname{frac}}(E) = 0$, then, for all $t>0$, 
\begin{align*}
\mathcal{L}_d(E) - 2^{\frac{d}{\alpha}}  \dfrac{\|p_\alpha\|_{\infty}}{t^{\frac{d}{\alpha}}} \mathcal{L}_d(E)^2 \leq 0,
\end{align*} 
which implies, by letting $t \rightarrow +\infty$, $\mathcal{L}_d(E) = 0$ (with a similar argument for $\mathcal{P}_{\operatorname{cl}}$).~Then, for all $t>0$,
\begin{align}\label{ineq:perim_frac_t}
 \mathcal{L}_d(E)  \leq \frac{t^{1 - \frac{1}{\alpha}}}{2(\alpha-1)} \left\| \frac{\nabla(p_\alpha)}{p_\alpha} \right\|_{L^1(\mu_\alpha)} \mathcal{P}_{\operatorname{frac}}(E)+2^{\frac{d}{\alpha}}  \dfrac{\|p_\alpha\|_{\infty}}{t^{\frac{d}{\alpha}}} \mathcal{L}_d(E)^2
\end{align}
and 
\begin{align}\label{ineq:cl_t}
\mathcal{L}_d(E) \leq \frac{t^{\frac{1}{\alpha}}}{2} \left(\bbe |Y_\alpha| \right) \mathcal{P}_{\operatorname{cl}}(E) + 2^{\frac{d}{\alpha}}  \dfrac{\|p_\alpha\|_{\infty}}{t^{\frac{d}{\alpha}}} \mathcal{L}_d(E)^2.
\end{align}
Let us start with the optimization of the inequality \eqref{ineq:perim_frac_t}.~For this purpose, let $t_0>0$ be defined by
\begin{align*}
t_0 := \left( \dfrac{2d \|p_\alpha\|_{\infty} 2^{\frac{d}{\alpha}}}{\| \nabla(p_\alpha)\|_{L^1(\mathbb{R}^d,dx)}} \right)^{\frac{\alpha}{d + \alpha-1}} \left( \dfrac{\mathcal{L}_d(E)^2}{\mathcal{P}_{\operatorname{frac}}(E)} \right)^{\frac{\alpha}{d+\alpha-1}}. 
\end{align*}
Then, for the first term of inequality \eqref{ineq:perim_frac_t},
\begin{align*}
(I) & := \frac{t_0^{1 - \frac{1}{\alpha}}}{2(\alpha-1)} \left\|\nabla(p_\alpha)\right\|_{L^1(\mathbb{R}^d,dx)} \mathcal{P}_{\operatorname{frac}}(E) \\
& = \frac{1}{2(\alpha-1)} (2d 2^{\frac{d}{\alpha}} \|p_\alpha\|_{\infty})^{\frac{\alpha-1}{d+ \alpha-1}} \|\nabla(p_\alpha)\|_{L^1(\mathbb{R}^d,dx)}^{\frac{d}{d+\alpha-1}}\mathcal{L}_d(E)^{\frac{2(\alpha-1)}{d+\alpha-1}} \mathcal{P}_{\operatorname{frac}}(E)^{\frac{d}{d+ \alpha-1}}. 
\end{align*} 
Moreover, for the second term of the inequality \eqref{ineq:perim_frac_t},
\begin{align*}
(II) & := 2^{\frac{d}{\alpha}}  \dfrac{\|p_\alpha\|_{\infty}}{t_0^{\frac{d}{\alpha}}} \mathcal{L}_d(E)^2 \\
& = \left(\frac{1}{2d}\right)^{\frac{d}{d+\alpha-1}}  (2^{\frac{d}{\alpha}} \|p_\alpha\|_{\infty})^{\frac{\alpha-1}{d+ \alpha-1}} \|\nabla(p_\alpha)\|_{L^1(\mathbb{R}^d,dx)}^{\frac{d}{d+\alpha-1}} \mathcal{L}_d(E)^{\frac{2(\alpha-1)}{d+\alpha-1}} \mathcal{P}_{\operatorname{frac}}(E)^{\frac{d}{d+ \alpha-1}}.
\end{align*}
Combining these expressions for $(I)$ and $(II)$, 
\begin{align*}
 \mathcal{L}_d(E) & \leq \left(\frac{1}{2(\alpha-1)} (2d)^{\frac{\alpha-1}{d+ \alpha-1}} + \left(\frac{1}{2d}\right)^{\frac{d}{d+\alpha-1}}\right) (2^{\frac{d}{\alpha}} \|p_\alpha\|_{\infty})^{\frac{\alpha-1}{d+ \alpha-1}} \|\nabla(p_\alpha)\|_{L^1(\mathbb{R}^d,dx)}^{\frac{d}{d+\alpha-1}} \\
 &\quad\quad \times \mathcal{L}_d(E)^{\frac{2(\alpha-1)}{d+\alpha-1}} \mathcal{P}_{\operatorname{frac}}(E)^{\frac{d}{d+ \alpha-1}},
\end{align*}
which implies the first inequality in \eqref{ineq:perim_frac_isoperi} with the constant $C^1_{\alpha,d}$ given by \eqref{eq:c_1_alpha}.~Next, for \eqref{ineq:cl_t}, set 
\begin{align*}
t_1 := \left( \dfrac{2d \|p_\alpha\|_{\infty} 2^{\frac{d}{\alpha}}}{\bbe |Y_\alpha|} \right)^{\frac{\alpha}{d+1}} \left( \dfrac{\mathcal{L}_d(E)^2}{\mathcal{P}_{\operatorname{cl}}(E)}\right)^{\frac{\alpha}{d+1}}.
\end{align*}
Then, for the first term of the inequality \eqref{ineq:cl_t}, 
\begin{align*}
A : = \dfrac{t_1^{\frac{1}{\alpha}}}{2} \bbe |Y_\alpha| \mathcal{P}_{\operatorname{cl}}(E) = \frac{1}{2} \left( 2d \|p_\alpha\|_{\infty} 2^{\frac{d}{\alpha}}\right)^{\frac{1}{d+1}} (\bbe |Y_\alpha|)^{\frac{d}{d+1}} \mathcal{P}_{\operatorname{cl}}(E)^{\frac{d}{d+1}} \mathcal{L}_d(E)^{\frac{2}{d+1}}. 
\end{align*}
For the second term of the inequality \eqref{ineq:cl_t}, 
\begin{align*}
B : = 2^{\frac{d}{\alpha}}  \dfrac{\|p_\alpha\|_{\infty}}{t_1^{\frac{d}{\alpha}}} \mathcal{L}_d(E)^2 = \frac{1}{(2d)^{\frac{d}{d+1}}}\left( 2^{\frac{d}{\alpha}} \|p_\alpha\|_{\infty}\right)^{\frac{1}{d+1}} (\bbe |Y_\alpha|)^{\frac{d}{d+1}} \mathcal{P}_{\operatorname{cl}}(E)^{\frac{d}{d+1}} \mathcal{L}_d(E)^{\frac{2}{d+1}}.
\end{align*}
Plugging these identities in \eqref{ineq:cl_t} leads to,
\begin{align*}
\mathcal{L}_d(E) \leq \left( \frac{1}{2} (2d)^{\frac{1}{d+1}}+\frac{1}{(2d)^{\frac{d}{d+1}}}\right)\left( 2^{\frac{d}{\alpha}} \|p_\alpha\|_{\infty}\right)^{\frac{1}{d+1}} (\bbe |Y_\alpha|)^{\frac{d}{d+1}} \mathcal{P}_{\operatorname{cl}}(E)^{\frac{d}{d+1}} \mathcal{L}_d(E)^{\frac{2}{d+1}},
\end{align*} 
which implies the second inequality in \eqref{ineq:perim_frac_isoperi} with the constant $C^2_{\alpha,d}$ given by \eqref{eq:c_2_alpha}.~This concludes the proof of the theorem.$\qed$

\begin{rem}\label{rem:Comi_et_al_Frank_et_al}
(i) In \cite[Theorem 4.4]{comi_stefani}, a fractional isoperimetric-type inequality is obtained in the rotationally invariant case.~Namely, for all $s \in (0,1)$, all $d \geq 2$ and all $E\in \mathcal{B}(\mathbb{R}^d)$ with $\mathcal{L}_d(E)<+\infty$ such that $|D^s (\bbone_E)|(\mathbb{R}^d)< +\infty$, 
\begin{align}\label{ineq:perim_fractional_iso_comi}
\mathcal{L}_d(E)^{\frac{d-s}{d}} \leq c_{s,d} |D^s (\bbone_E)|(\mathbb{R}^d),
\end{align}
for some constant $c_{s,d}>0$ and where $D^s (\bbone_E)$ is defined in \cite[Definition 4.1]{comi_stefani} by duality with the fractional gradient operator $D^{s,\operatorname{rot}}$.~In order to obtain inequality \eqref{ineq:perim_fractional_iso_comi},~the authors use a previous result regarding strong Sobolev-type embedding for the operator $D^{s,\operatorname{rot}}$ with $p = 1$ established in \cite[Theorem A']{SSS_17} (see also \cite{Spec_19}).~The restriction on the dimension $d$ is a consequence of this embedding theorem.~Moreover, the constant $c_{s,d}$ is not explicitely given.~Theorem \ref{thm:isoperimetric_type} provides explicit constants and holds for all non-degenerate symmetric $\alpha$-stable L\'evy measures $\nu_\alpha$ on $\mathbb{R}^d$ verifying \eqref{eq:scale} with $\alpha \in (1,2)$.~Finally, the equality cases, even when $\nu_\alpha = \nu_\alpha^{\operatorname{rot}}$, are open (see \cite[Section $1.4.$ Future Developments]{comi_stefani}).\\
(ii) In \cite[Theorem $4.1$]{FS_08},~the authors obtain the following version of a fractional isoperimetric-type inequality: for all $d \geq 1$, all $s \in (0,1)$ and all $E\in \mathcal{B}(\mathbb{R}^d)$ with $\mathcal{L}_d(E)<+\infty$, 
\begin{align}\label{ineq:perim_fractional_iso_frank}
\mathcal{L}_d(E)^{\frac{d-s}{d}} \leq C_{s,d} \mathcal{P}_s(E),
\end{align}
where $C_{s,d}>0$ is an explicit sharp constant depending on $d$ and $s$ and $\mathcal{P}_s(E)$ is the fractional perimeter of $E$ defined by
\begin{align}\label{def:frac_perimeter_frank}
\mathcal{P}_s(E)  := \int_{E}\int_{E^c} \dfrac{dxdy}{\|x-y\|^{d+s}}.
\end{align}
Moreover, there is equality if and only if $E$ is an Euclidean ball.~For an explicit computation of the constant $C_{s,d}$ see \cite{G_RLMA20}.\\
(iii) Let $d \geq 1$, let $\alpha \in (1,2)$ and let $\nu_\alpha(du) = du /\| u \|^{\alpha+d}$.~Then, for all $E \in \mathcal{B}(\mathbb{R}^d)$ with $\mathcal{L}_d(E) < + \infty$ such that $\mathcal{P}_{\alpha-1}(E)<+\infty$,  
\begin{align*}
\left\| D^{\alpha - 1}\left(\bbone_E\right) \right\|_{L^1(\mathbb{R}^d,dx)} & \leq \int_{\mathbb{R}^d} \left(\int_{\mathbb{R}^d} \frac{du}{\|u\|^{d+ \alpha - 1}} \left|\bbone_E(x+u) - \bbone_E(x)\right| \right) dx\\
&\leq \int_{\mathbb{R}^d} \int_{\mathbb{R}^d} \dfrac{\left|\bbone_E(u) - \bbone_E(x)\right|}{\|x - u\|^{d+ \alpha - 1}} dxdu = 2 \mathcal{P}_{\alpha-1}(E). 
\end{align*}
Thus, for these choices of L\'evy measure and of $E$, 
\begin{align}\label{ineq:simple_compineq}
\mathcal{P}_{\operatorname{frac}}\left(E\right) \leq 2 \mathcal{P}_{\alpha-1}(E) < +\infty. 
\end{align}
The previous reasoning continues to hold if the spherical component of $\nu_\alpha$ is only dominated by a uniform measure on $\mathbb{S}^{d-1}$.\\
(iv) \cite[Proposition $4.8$, Inequality $(4.10)$]{comi_stefani} compares these two notions of fractional perimeter (up to some normalization constant).~Moreover, it is asked in \cite{comi_stefani} for which sets the equality case holds true. This question is fully answered in \cite[Corollary $2.11$, (i)]{comi_stefani_rmc22} where it is proved that for $E \subset \mathbb{R}^d$ such that $\mathcal{P}_{\alpha-1}(E)<+\infty$ and such that $\mathcal{L}_d(E)\mathcal{L}_d(E^c)>0$, inequality \eqref{ineq:simple_compineq} is strict. 
\end{rem}
\noindent
From Remark \ref{rem:Comi_et_al_Frank_et_al}, (iii), we have the following non-asymptotic representation formula for the fractional perimeter $\mathcal{P}_{\operatorname{frac}}$ for a certain collection of Borel subsets of $\mathbb{R}^d$. 

\begin{lem}\label{lem:non_asymptotic_representation}
Let $\alpha \in (1,2)$, let $\nu_\alpha$ be a non-degenerate symmetric $\alpha$-stable L\'evy measure on $\mathbb{R}^d$, $d \geq 1$, such that its spherical part 
is dominated by a uniform measure on $\mathbb{S}^{d-1}$.~Let $E \in \mathcal{B}(\mathbb{R}^d)$ with finite Lebesgue measure be such that $\mathcal{P}_{\alpha-1}(E)<+\infty$. Then, 
\begin{align*}
\mathcal{P}_{\operatorname{frac}}(E) = \left\| D^{\alpha-1} (\bbone_E) \right\|_{L^1(\mathbb{R}^d,dx)}.
\end{align*}  
\end{lem}

\begin{proof}
From Remark \ref{rem:Comi_et_al_Frank_et_al}, (iii), $\left\| D^{\alpha-1} (\bbone_E) \right\|_{L^1(\mathbb{R}^d,dx)} < +\infty$.~Moreover, by commutation and the contraction property of $(P^\alpha_t)_{t \geq 0}$ on $L^1(\mathbb{R}^d,dx)$,
\begin{align*}
\underset{t \rightarrow 0^+}{\lim} \left\| D^{\alpha-1} P^\alpha_t(\bbone_E) \right\|_{L^1(\mathbb{R}^d,dx)} \leq \left\| D^{\alpha-1} (\bbone_E) \right\|_{L^1(\mathbb{R}^d,dx)}.
\end{align*}
Finally, by standard regularization arguments (see, e.g., \cite{Stein_1}), the semigroup $(P^\alpha_t)_{t \geq 0}$ is strongly continuous on $L^1(\mathbb{R}^d,dx)$. This concludes the proof of the lemma.
\end{proof}
\noindent
From the proof of the isoperimetric inequalities,~a crucial quantity to analyze is $A_t$, $t \geq 0$, defined by \eqref{def:A_t}.~More generally, it is interesting to study the following variant:~for all $t \geq 0$ and all $A,B$ Borel measurable subsets of $\mathbb{R}^d$, 
\begin{align}\label{eq:definition_K_functional}
K^\alpha_t(A,B) = \langle P^\alpha_t(\bbone_A) ; \bbone_B \rangle_{L^2(\mathbb{R}^d,dx)}. 
\end{align}  
Note that, for all $t \geq 0$ and all $A\in \mathcal{B}(\mathbb{R}^d)$ with $\mathcal{L}_d(A)<+\infty$, 
\begin{align*}
\mathcal{L}_d(A) - \int_A P^\alpha_t(\bbone_A) dx = \int_{\mathbb{R}^d} P^\alpha_t(\bbone_A) dx - \int_AP^\alpha_t(\bbone_A) dx = K^\alpha_t(A,A^c). 
\end{align*}
For the classical heat semigroup, the functional $(K^2_t)_{t \geq 0}$ already appears in the works \cite{MLedoux_94,Preun_04} in connection with the isoperimetric inequality.~In \cite{Valverde_17}, for $\alpha \in (0,2)$ and when $\nu_\alpha = \nu_\alpha^{\operatorname{rot}}$, the corresponding functional $(K_t^{\alpha,\operatorname{rot}})_{t \geq 0}$ has been studied but not directly in connection with isoperimetric problems.~In \cite{Valverde_17}, the author is interested in the small time asymptotic of $K_t^{\alpha,\operatorname{rot}}(A,A^c)$, when the boundary of $A$ possesses some form of regularity.~In particular, \cite[Theorem 1.2, inequality (1.17)]{Valverde_17} is completely analogous to the one obtained in \cite[Proposition 1.1]{MLedoux_94} (see, also \cite[Proposition 8]{Preun_04}) in the Gaussian case.~Thanks to the pseudo-Poincar\'e inequalities, for all $t \geq 0$ and all $A\in \mathcal{B}(\mathbb{R}^d)$ with $\mathcal{L}_d(A)<+\infty$,
\begin{align}\label{ineq:heat_content_first_ineq}
K^\alpha_t(A,A^c)  \leq \frac{t^{\frac{1}{\alpha}}}{2} \left(\bbe |Y_\alpha| \right) \mathcal{P}_{\operatorname{cl}}(A)
\end{align}
and
\begin{align}\label{ineq:heat_content_second_ineq}
K^\alpha_t(A,A^c) \leq  \dfrac{t^{1 - \frac{1}{\alpha}}}{2(\alpha-1)} \left\| \dfrac{\nabla(p_\alpha)}{p_\alpha}\right\|_{L^1(\mu_\alpha)} \mathcal{P}_{\operatorname{frac}}(A),
\end{align}
where $Y_\alpha$ is a random variable such that $ \bbe \exp \left( i \xi Y_\alpha \right) = \exp\left( - |\xi|^\alpha\right)$, for $\xi \in \bbr$, and $p_\alpha$ is the positive Lebesgue density of the probability measure $\mu_\alpha$.~In the next proposition, for $\nu_\alpha = \nu_\alpha^{\operatorname{rot}}$, let us retrieve the sharp classical isoperimetric inequality based on \eqref{ineq:heat_content_first_ineq} and the Riesz rearrangement inequality. 

\begin{thm}\label{thm:sharp_classical_isoperimetric_ineq}
Let $\alpha \in (1,2)$ and let $\mu_\alpha^{\operatorname{rot}}$ be the rotationally invariant $\alpha$-stable probability measure on $\mathbb{R}^d$, $d \geq 1$, with positive Lebesgue density $p_\alpha^{\operatorname{rot}}$ and L\'evy measure $\nu_\alpha^{\operatorname{rot}}$. Then, for all $E\in \mathcal{B}(\mathbb{R}^d)$ with $\mathcal{L}_d(E)<+\infty$ such that $\mathcal{P}(E)<+\infty$ and all $t > 0$, 
\begin{align*}
K^{\alpha, \operatorname{rot}}_t(B,B^c) \leq K^{\alpha, \operatorname{rot}}_t(E,E^c)\leq \dfrac{t^{\frac{1}{\alpha}}}{2^{1+\frac{1}{\alpha}}} \bbe |Y_\alpha| \mathcal{P}(E),
\end{align*}
where $B$ is an Euclidean ball centered at the origin such that $\mathcal{L}_d(B) = \mathcal{L}_d(E)$. Moreover,
\begin{align}\label{eq:limit_ball}
\underset{t \rightarrow 0^+}{\lim} \frac{K^{\alpha, \operatorname{rot}}_t(B,B^c)}{t^{\frac{1}{\alpha}}} = \frac{1}{2^{\frac{1}{\alpha}}} \frac{\Gamma \left(1- \frac{1}{\alpha}\right) }{\pi} \mathcal{P}(B).
\end{align}
Finally, 
\begin{align*}
\mathcal{P}(E) \geq \mathcal{P}(B). 
\end{align*}
\end{thm}

\begin{proof}
First, let us prove that when $\nu_\alpha = \nu_\alpha^{\operatorname{rot}}$, the perimeter $\mathcal{P}_{\operatorname{cl}}$ coincides (up to some constant) with the one introduced in \cite{De_Giorgi} and based on the classical heat semigroup. First, for all $z \in \mathbb{R}^d$, 
\begin{align*}
\sigma_\alpha^{\operatorname{rot}}(z)^\alpha = \int_{\mathbb{S}^{d-1}} |\langle  y ;  z \rangle|^\alpha \lambda(dy) = \frac{\|z\|^\alpha}{2}. 
\end{align*}
Then, for all $t>0$ and all $E\in \mathcal{B}(\mathbb{R}^d)$ with $\mathcal{L}_d(E)<+\infty$, 
\begin{align*}
\left\| \sigma_{\alpha}\left(\nabla P^{\alpha, \operatorname{rot}}_t (\bbone_E) \right) \right\|_{L^1(\mathbb{R}^d,dx)} = \frac{1}{2^{\frac{1}{\alpha}}} \| \nabla P^{\alpha, \operatorname{rot}}_t (\bbone_E) \|_{L^1(\mathbb{R}^d,dx)},
\end{align*}
so that, 
\begin{align}
\mathcal{P}_{\operatorname{cl}} (E) = \frac{1}{2^{\frac{1}{\alpha}}} \underset{t \rightarrow 0^+}{\lim} \| \nabla P^{\alpha, \operatorname{rot}}_t (\bbone_E) \|_{L^1(\mathbb{R}^d,dx)},
\end{align}
(which might or might not be finite) and where $(P^{\alpha, \operatorname{rot}}_t)_{t \geq 0}$ is given, for all $f \in \mathcal{S}(\mathbb{R}^d)$, all $t \geq 0$ and all $x \in \mathbb{R}^d$, by
\begin{align*}
P^{\alpha, \operatorname{rot}}_t(f)(x) = \frac{1}{(2\pi)^d} \int_{\mathbb{R}^d} \mathcal{F}(f)(\xi) e^{i \langle x ; \xi \rangle} \exp\left( - t \frac{\|\xi\|^\alpha}{2}\right) d\xi. 
\end{align*}
Now, let $E\in \mathcal{B}(\mathbb{R}^d)$ with $\mathcal{L}_d(E)<+\infty$ be such that 
\begin{align*}
\mathcal{P}(E) = \underset{\phi \in \Phi }{\sup} \int_{\mathbb{R}^d} \bbone_E(x) \operatorname{div}(\phi)(x)dx <+\infty,
\end{align*}
where $\Phi = \{ \phi \in \mathcal{C}_c^{1}(\mathbb{R}^d, \mathbb{R}^d) , \, \|\phi\|_{\infty} \leq 1\}$ and $\operatorname{div}$ is the classical divergence operator. With a slight abuse of notations, we note $\mathcal{P}(E) = \|\nabla \bbone_E\|_{L^1(\mathbb{R}^d,dx)}$. Now, the semigroup $(P^{\alpha, \operatorname{rot}}_t)_{t \geq 0}$ is strongly continuous on $L^1(\mathbb{R}^d,dx)$, i.e.,
\begin{align*}
P^{\alpha, \operatorname{rot}}_t(\bbone_E) \rightarrow \bbone_E, \quad L^1(\mathbb{R}^d,dx), 
\end{align*}
as $t$ tends to $0^+$. This implies that, 
\begin{align*}
\|\nabla \bbone_E\|_{L^1(\mathbb{R}^d,dx)} \leq \underset{t \rightarrow 0^+}{\liminf} \| \nabla P^{\alpha, \operatorname{rot}}_t (\bbone_E) \|_{L^1(\mathbb{R}^d,dx)} = \underset{t \rightarrow 0^+}{\lim} \| \nabla P^{\alpha, \operatorname{rot}}_t (\bbone_E) \|_{L^1(\mathbb{R}^d,dx)}.
\end{align*}
Hence, for all $t>0$ and all $\phi \in \Phi$,  
\begin{align*}
\int_{\mathbb{R}^d} P^{\alpha, \operatorname{rot}}_t(\bbone_E) \operatorname{div}(\phi)(x) dx = \int_E \operatorname{div} \left(P^{\alpha,\operatorname{rot}}_t(\phi)\right)(x)dx. 
\end{align*}
But, the semigroup $(P^{\alpha,\operatorname{rot}}_t)_{t \geq 0}$ is a contraction on $L^{\infty}(\mathbb{R}^d,dx)$ so that, 
\begin{align*}
\| P^{\alpha,\operatorname{rot}}_t (\phi)\|_{L^\infty(\mathbb{R}^d,dx)} \leq \|\phi\|_{L^\infty(\mathbb{R}^d,dx)} \leq 1.
\end{align*}
Moreover, $P^{\alpha,\operatorname{rot}}_t (\phi) \in \mathcal{C}^1(\mathbb{R}^d, \mathbb{R}^d)$ and $\| \operatorname{div} P^{\alpha,\operatorname{rot}}_t (\phi) \|_{L^1(\mathbb{R}^d,dx)} <+\infty$. Thus, 
\begin{align*}
\int_{\mathbb{R}^d} P^{\alpha, \operatorname{rot}}_t(\bbone_E) \operatorname{div}(\phi)(x) dx \leq \sup \left\{\int_E \operatorname{div}(\phi)dx , \, \phi \in \Phi \right\} = \mathcal{P}(E),
\end{align*}
which provides the converse inequality, and so, $\mathcal{P}_{\operatorname{cl}}(E) = \mathcal{P}(E) / 2^{1/\alpha}$. Thus, for all $t \geq 0$ and all $E\in \mathcal{B}(\mathbb{R}^d)$ with $\mathcal{L}_d(E) <+\infty$ such that $\mathcal{P}(E) <+\infty$,
\begin{align*}
K^{\alpha, \operatorname{rot}}_t(E,E^c)  \leq \frac{t^{\frac{1}{\alpha}}}{2} \left(\bbe |Y_\alpha| \right) \mathcal{P}_{\operatorname{cl}}(E) =  \frac{t^{\frac{1}{\alpha}}}{2^{1+\frac{1}{\alpha}}} \left(\bbe |Y_\alpha| \right)\mathcal{P}(E).
\end{align*}
Moreover,
\begin{align*}
\langle P^{\alpha, \operatorname{rot}}_t\left(\bbone_E\right) ; \bbone_E \rangle_{L^2(\mathbb{R}^d,dx)} = \int_{\mathbb{R}^d} \int_{\mathbb{R}^d} \bbone_E(x) \bbone_E(y) p^{\operatorname{rot}}_{\alpha} \left(\frac{x -y}{t^{\frac{1}{\alpha}}}\right) \frac{dxdy}{t^{\frac{d}{\alpha}}}. 
\end{align*}
Thanks to \cite[Theorem $3.7$]{Lieb_Loss_book01}, 
\begin{align*}
\langle P^{\alpha, \operatorname{rot}}_t\left(\bbone_E\right) ; \bbone_E \rangle_{L^2(\mathbb{R}^d,dx)}  \leq \int_{\mathbb{R}^d} \int_{\mathbb{R}^d} (\bbone_E)^*(x) (\bbone_E)^*(y) (p^{\operatorname{rot}}_{\alpha})^* \left(\frac{x -y}{t^{\frac{1}{\alpha}}}\right) \frac{dxdy}{t^{\frac{d}{\alpha}}},
\end{align*}
where $f^*$ is the symmetric non-increasing rearrangement of the function $f$. Now, $(\bbone_E)^* = \bbone_B$ where $B$ is an Euclidean ball centered at the origin such that $ \mathcal{L}_d(E) = \mathcal{L}_d(B)$. Moreover, $(p^{\operatorname{rot}}_{\alpha})^* = p^{\operatorname{rot}}_{\alpha}$ since $p^{\operatorname{rot}}_{\alpha}$ is the inverse Fourier transform of a radially symmetric function and, by subordination, for all $x \in \mathbb{R}^d$,
\begin{align*}
p^{\operatorname{rot}}_{\alpha}(x) = \int_0^{+\infty} p_2\left(\frac{x}{\sqrt{s}}\right) \frac{\eta_1^\alpha(ds)}{s^{\frac{d}{2}}},
\end{align*} 
where $\eta_1^\alpha$ is the density of an $\alpha/2$-stable subordinator (at time $1$) with Laplace transform given, for all $\lambda>0$, by 
\begin{align*}
\int_0^{+\infty} \exp\left( - s \lambda\right) \eta_1^\alpha(ds) = \exp\left( - \lambda^{\frac{\alpha}{2}}\right),
\end{align*}
and $p_2$ is the centered multivariate Gaussian density with appropriate variance. Then, for all $t \geq 0$, 
\begin{align*}
\langle P^{\alpha, \operatorname{rot}}_t\left(\bbone_E\right) ; \bbone_E \rangle_{L^2(\mathbb{R}^d,dx)} \leq  \langle P^{\alpha, \operatorname{rot}}_t\left(\bbone_B\right) ; \bbone_B \rangle_{L^2(\mathbb{R}^d,dx)}, 
\end{align*}
which implies, for all $t \geq 0$,
\begin{align*}
K^{\alpha, \operatorname{rot}}_t(B,B^c) \leq K^{\alpha, \operatorname{rot}}_t(E,E^c)  \leq \frac{t^{\frac{1}{\alpha}}}{2^{1+\frac{1}{\alpha}}} \left(\bbe |Y_\alpha| \right)\mathcal{P}(E).
\end{align*}
Next, using the argument in \cite[Proof of Theorem $1.2$, (i)]{Valverde_17}, 
\begin{align*}
\underset{t \rightarrow 0^+}{\lim} \frac{K^{\alpha, \operatorname{rot}}_t(B,B^c)}{t^{\frac{1}{\alpha}}} = \frac{1}{2^{\frac{1}{\alpha}}} \frac{\Gamma \left(1- \frac{1}{\alpha}\right) }{\pi} \mathcal{P}(B).
\end{align*}
Finally,  the previous limit together with \eqref{ineq:heat_content_first_ineq} implies the sharp isoperimetric inequality:
\begin{align*}
\mathcal{P}(E) \geq \frac{2}{\bbe |Y_\alpha|} \frac{\Gamma \left(1- \frac{1}{\alpha}\right) }{\pi} \mathcal{P}(B) =  \mathcal{P}(B),
\end{align*}
since $\bbe |Y_\alpha| = 2 \Gamma(1-1/\alpha)/\pi$, thanks to subordination. This concludes the proof of the theorem. 
\end{proof}
\noindent
\begin{rem}\label{rem:rotationally_invariant_Pcl}
(i) The stable heat kernel of the semigroup $(P^{\alpha,\operatorname{rot}}_t)_{t \geq 0}$ displays spatial properties (in term of localization) different from the Gaussian heat kernel ones. Both are radially symmetric and decreasing with respect to the radial coordinate but $p_\alpha^{\operatorname{rot}}$ satisfies,
\begin{align*}
\underset{\|x\| \rightarrow +\infty}{\lim} \|x\|^{\alpha+d} p_\alpha^{\operatorname{rot}}(x) = A_{\alpha,d},
\end{align*}
for some constant $A_{\alpha,d}>0$ depending on $\alpha$ and $d$ which can be computed explicitely (see, e.g., \cite[Theorem $2.1$]{BG_60}). The fact that one can retrieve the sharp classical isoperimetric inequality from the semigroup $(P^{\alpha , \operatorname{rot}}_t)_{t \geq 0}$ can be seen as a complementary result to the ones put forward in the series of works \cite{BBM_2001,JD_02,MRT_19} where general radially symmetric kernels are used to characterize sets of finite perimeter and functions of bounded variation.\\
(ii) The limiting behavior \eqref{eq:limit_ball} is actually true for all open set of $\mathbb{R}^d$ with finite Lebesgue measure and with finite perimeter. Indeed, for all such set $\Omega$ and \cite[Theorem $2$]{CG_17}, 
\begin{align}\label{eq:limit_open_rotationally}
\underset{t \rightarrow 0^+}{\lim} \frac{K^{\alpha, \operatorname{rot}}_t(\Omega,\Omega^c)}{t^{\frac{1}{\alpha}}} = \frac{1}{2^{\frac{1}{\alpha}}} \frac{\Gamma \left(1- \frac{1}{\alpha}\right) }{\pi} \mathcal{P}(\Omega).
\end{align}
\end{rem}
\noindent
Next, we are interested in the case where the spectral measure is discrete.~Namely,
\begin{align*}
\lambda^{\operatorname{dis}}_1 (dy) = \sum_{k =1}^d w_k \delta_{\pm e_k}(dy),
\end{align*} 
where $\left( w_1 , \dots, w_d\right) \in (\bbr^*_+)^d$. Then, for all $z \in \mathbb{R}^d$, 
\begin{align*}
\sigma^{\operatorname{dis}}_\alpha(z) = \left(2\sum_{k=1}^d w_k |\langle z; e_k \rangle|^\alpha\right)^{\frac{1}{\alpha}}. 
\end{align*}
By analogy with the rotationally invariant case, the first question to investigate is to find: for which Borel measurable subsets of $\mathbb{R}^d$ with finite Lebesgue measure, the following limit 
\begin{align}\label{eq:limit_K_discrete}
L^{\operatorname{dis}}(\Omega) = \underset{t \rightarrow 0^+}{\lim} \frac{K^{\alpha, \operatorname{dis}}_t(\Omega, \Omega^c)}{t^{\frac{1}{\alpha}}},
\end{align}
exists and can be computed explicitely. The next technical lemma is a straightforward consequence of the fact that the map $\sigma_\alpha$ is a norm under the non-degeneracy condition \eqref{eq:non_deg}.

\begin{lem}\label{lem:limsup_non-degeneracy}
Let $\alpha \in (1,2)$ and let $(K_t^\alpha)_{t \geq 0}$ be the functional defined by \eqref{eq:definition_K_functional}. Then, for all $E\in \mathcal{B}(\mathbb{R}^d)$, $d \geq 1$, with $\mathcal{L}_d(E)<+\infty$ such that $\mathcal{P}(E)<+\infty$, 
\begin{align*}
\underset{t \rightarrow 0^+}{\limsup} \dfrac{K^{\alpha}_t(E,E^c)}{t^{\frac{1}{\alpha}}} <+\infty. 
\end{align*}
\end{lem}

\begin{proof}
First, since $\sigma_\alpha$ is a norm on $\mathbb{R}^d$, for all $z \in \mathbb{R}^d$, 
\begin{align*}
C_{1,\alpha,d} \| z \| \leq \sigma_\alpha(z) \leq C_{2,\alpha,d} \| z \|.
\end{align*}
for some $C_{1,\alpha,d}, C_{2,\alpha,d}>0$. Then, for all $t >0$, 
\begin{align}
C_{1,\alpha,d}  \| \nabla P^\alpha_t \bbone_E \|_{L^1(\mathbb{R}^d,dx)} \leq \| \sigma_\alpha\left(\nabla P^\alpha_t \bbone_E\right) \|_{L^1(\mathbb{R}^d,dx)} \leq C_{2,\alpha,d} \| \nabla P^\alpha_t \bbone_E \|_{L^1(\mathbb{R}^d,dx)}. 
\end{align}
Now, as previously, the semigroup $(P^\alpha_t)_{t \geq 0}$ is strongly continuous on $L^1(\mathbb{R}^d,dx)$. Then, 
\begin{align*}
\mathcal{P}(E) \leq \underset{t \rightarrow 0^+}{\liminf} \| \nabla P^\alpha_t \bbone_E \|_{L^1(\mathbb{R}^d,dx)} = \underset{t \rightarrow 0^+}{\lim}  \| \nabla P^\alpha_t \bbone_E \|_{L^1(\mathbb{R}^d,dx)}.
\end{align*}
Moreover, using the fact that the semigroup $(P^\alpha_t)_{t \geq 0}$ is a contraction on $L^{\infty}(\mathbb{R}^d,dx)$, that $P^\alpha_t(\varphi) \in \mathcal{C}^1(\mathbb{R}^d,\mathbb{R}^d)$, for $\varphi \in \mathcal{C}_c^1(\mathbb{R}^d,\mathbb{R}^d)$, and that
\begin{align*}
\| \operatorname{div}\left(P^\alpha_t(\varphi)\right) \|_{L^1(\mathbb{R}^d, dx)} \leq \| \operatorname{div}(\varphi) \|_{L^1(\mathbb{R}^d,dx)} <+\infty,
\end{align*}
one gets 
\begin{align*}
\underset{t \rightarrow 0^+}{\lim}  \| \nabla P^\alpha_t \bbone_E \|_{L^1(\mathbb{R}^d,dx)} = \mathcal{P}(E) <+\infty.
\end{align*}
Then, for all $t>0$, 
\begin{align*}
\dfrac{K^{\alpha}_t(E,E^c)}{t^{\frac{1}{\alpha}}} \leq \frac{\bbe |Y_\alpha|}{2}C_{2, \alpha,d} \mathcal{P}(E).
\end{align*}
The conclusion follows easily.
\end{proof}
\noindent
Next, let us recall some basic facts about the stable heat kernel associated with the L\'evy measure whose spectral measure is given by \eqref{eq:spectral_measure_discrete}. In this situation, the Fourier transform of the associated $\alpha$-stable probability measure is given, for all $\xi \in \mathbb{R}^d$, by
\begin{align}\label{eq:characteristic_discrete}
\varphi^w_{\alpha,d}(\xi) = \exp\left( - 2 \sum_{k = 1}^{d} w_k |\langle \xi ; e_k\rangle|^{\alpha} \right), 
\end{align}
so that the $\alpha$-stable Lebesgue density tensorizes: for all $x \in \mathbb{R}^d$, 
\begin{align}\label{eq:stable_density_discrete}
p^w_{\alpha,d}\left(x \right) : = \prod_{k = 1}^d \frac{1}{(2w_k)^{\frac{1}{\alpha}} }p_{\alpha,1}\left( \frac{x_k}{(2 w_k)^{\frac{1}{\alpha}}} \right),
\end{align}
where $p_{\alpha,1}$ is defined, for all $x \in \bbr$, by
\begin{align*}
p_{\alpha,1}(x) = \int_{\bbr} e^{i x \xi}\exp\left( - |\xi|^\alpha \right) \frac{d\xi}{2\pi}. 
\end{align*}
Thus, the corresponding stable heat kernel is given, for all $x,y \in \mathbb{R}^d$ and all $t>0$, by
\begin{align}\label{eq:stable_heat_discrete}
p^w_{\alpha,d}(x,y,t) := \frac{1}{t^{\frac{d}{\alpha}}}p^w_{\alpha,d}\left( \frac{x-y}{t^{\frac{1}{\alpha}}} \right) = \prod_{k = 1}^d \frac{1}{(2w_k t)^{\frac{1}{\alpha}} }p_{\alpha,1}\left( \frac{x_k-y_k}{(2 w_k t)^{\frac{1}{\alpha}}} \right).
\end{align} 
With the help of the product structure of $p^w_{\alpha,d}$, let us perform our analysis on Borel measurable subsets of $\mathbb{R}^d$ which can be written as product of one-dimensional Borel measurable subsets of $\bbr$ with finite length. Then, for such a set $E = E_1 \times \cdots \times E_d$ and all $t > 0$, 
\begin{align}\label{eq:tensorization_formula}
K_t^{\alpha,\operatorname{dis}} \left(E , E^c\right) = \mathcal{L}_d(E) - \int_E P^{\alpha,\operatorname{dist}}_t(\bbone_E)(x)dx = \mathcal{L}_d(E) - \int_{E}\int_{E} p^w_{\alpha,d}\left( \frac{x-y}{t^{\frac{1}{\alpha}}} \right) \dfrac{dxdy}{t^{\frac{d}{\alpha}}}.
\end{align}
In particular, 
\begin{align*}
\int_{E}\int_{E} p^w_{\alpha,d}\left( \frac{x-y}{t^{\frac{1}{\alpha}}} \right) \dfrac{dxdy}{t^{\frac{d}{\alpha}}} & = \int_{E_1 \times \cdots \times E_d} \int_{E_1 \times \cdots \times E_d}  \prod_{k = 1}^d \frac{1}{(2w_k t)^{\frac{1}{\alpha}} }p_{\alpha,1}\left( \frac{x_k-y_k}{(2 w_k t)^{\frac{1}{\alpha}}} \right) dxdy\\
& = \prod_{k = 1}^d \int_{E_k} \int_{E_k} \frac{1}{(2w_k t)^{\frac{1}{\alpha}} }p_{\alpha,1}\left( \frac{x_k-y_k}{(2 w_k t)^{\frac{1}{\alpha}}} \right) dx_kdy_k. 
\end{align*}
Then, by the Riesz rearrangement inequality for $d=1$ (see \cite[Lemma $3.6$]{Lieb_Loss_book01}),
\begin{align*}
\int_{E}\int_{E} p^w_{\alpha,d}\left( \frac{x-y}{t^{\frac{1}{\alpha}}} \right) \dfrac{dxdy}{t^{\frac{d}{\alpha}}} \leq \int_{R}\int_{R} p^w_{\alpha,d}\left( \frac{x-y}{t^{\frac{1}{\alpha}}} \right) \dfrac{dxdy}{t^{\frac{d}{\alpha}}},
\end{align*}
where $R = R _1 \times \cdots \times R_d $ with $R_k$ an interval of finite length centered at the origin such that $\mathcal{L}_1(R_k) = \mathcal{L}_1(E_k)$, for all $k \in \{1, \dots, d\}$. Thus, for all such $E$ and all $t \geq 0$, 
\begin{align}\label{ineq:discrete_case}
K_t^{\alpha,\operatorname{dis}} \left(E , E^c\right) \geq K_t^{\alpha,\operatorname{dis}} \left(R , R^c\right). 
\end{align}
The next lemma computes the limit $L^{\operatorname{dis}}(R)$  in \eqref{eq:limit_K_discrete} using the one-dimensional result. 

\begin{lem}\label{lem:limit_discrete_K}
Let $\alpha\in (1,2)$ and let $\nu_\alpha^{\operatorname{dis}}$ be the L\'evy measure on $\mathbb{R}^d$, $d \geq 1$, verifying \eqref{eq:scale} with spectral measure given by \eqref{eq:spectral_measure_discrete}. Then, for all $R = R _1 \times \cdots \times R_d $, where, for $k \geq 1$, $R_k = [-a_k , a_k]$ with $0<a_k<+\infty$,
\begin{align}\label{eq:limit_discrete_K}
\underset{t \rightarrow 0^+}{\lim} \dfrac{K_t^{\alpha,\operatorname{dis}} \left(R , R^c\right)}{t^{\frac{1}{\alpha}}} = \frac{2}{\pi} \Gamma \left(1 - \frac{1}{\alpha}\right) \sum_{i = 1}^d (2w_i)^{\frac{1}{\alpha}} \prod_{k =1, k \ne i}^d 2a_k. 
\end{align}
\end{lem}

\begin{proof}
Let us use the product structure of $R$ and the density $p^w_{\alpha,d}$ to obtain the limit \eqref{eq:limit_discrete_K}.  First,  for all $t>0$, 
\begin{align*}
K_t^{\alpha, \operatorname{dis}} (R,R^c) & = \mathcal{L}_d(R) - \int_{R} P^{\alpha,\operatorname{dis}}_t(\bbone_R)(x)dx \\
& = \int_{\mathbb{R}^d} P^{\alpha,\operatorname{dis}}_t(\bbone_R)(x)dx - \int_{R} P^{\alpha,\operatorname{dis}}_t(\bbone_R)(x)dx \\
& = \int_{R^c} P^{\alpha,\operatorname{dis}}_t(\bbone_R)(x)dx  \\
& = \int_{R^c} \int_{R} p^w_{\alpha,d} \left(\frac{x-y}{t^{\frac{1}{\alpha}}}\right) \frac{dxdy}{t^{\frac{d}{\alpha}}}.
\end{align*}
Now,  $R = R_1 \times \cdots \times R_d$,  so that, 
\begin{align*}
R^c = \{ x \in \mathbb{R}^d, \, \exists i \in \{1, \dots, d\}, \ x_i \in R_i^c \} = \cup_{i  =1}^d  \tilde{R}_i, 
\end{align*}
$\tilde{R}_i = \bbr \times \cdots \times \bbr \times R_i^c \times \bbr \times \cdots \times \bbr $.  It is clear that $R^c$ can be written as the disjoint union of the $2^d-1$ sets of the form 
\begin{align*}
\tilde{R}_{\operatorname{elem}} = S_1 \times \cdots \times S_d,
\end{align*}
where $S_i = R_i$ or $R_i^c$, for all $i  \in \{1, \dots, d\}$, with at least one $i$ for which $S_i$ is equal to $R_i^c$ (The set $R$ is excluded from the previous list). Then, it is natural to encode the structure of such sets using vectors with values in $\{0,1\}^d$. Then, let $S_u$ be defined, for all $u \in \{0,1\}^d$, by 
\begin{align*}
S_u = S_{1,u(1)} \times \cdots \times S_{d,u(d)},
\end{align*}  
with $S_{i,0} = R_i^c$ and $S_{i,1} = R_i$, for all $i \in \{1, \dots, d\}$.  Then, 
\begin{align*}
R^c = \cup_{u \in \{0,1\}^d,\, u \ne (1, \dots, 1)} S_u, 
\end{align*}
and $S_u \cap S_v = \emptyset$,  for all $u , v \in \{0,1\}^d \setminus \{ (1, \cdots,1) \}$ with $v \ne u$. Then, for all $t>0$, 
\begin{align*}
K_t^{\alpha, \operatorname{dis}}(R,R^c) = \sum_{u \in \{0;1\}^d, \, u \ne (1, \cdots,1)} \int_{S_u} \int_R p^w_{\alpha,d} \left(\frac{x-y}{t^{\frac{1}{\alpha}}}\right) \frac{dxdy}{t^{\frac{d}{\alpha}}}.
\end{align*}
Now, using the product structure of $p^w_{\alpha,d}$,  for all $u \in \{0,1\}^d \setminus \{(1, \dots, d)\}$, 
\begin{align*}
 \int_{S_u} \int_R p^w_{\alpha,d} \left(\frac{x-y}{t^{\frac{1}{\alpha}}}\right) \frac{dxdy}{t^{\frac{d}{\alpha}}} = \prod_{k  =1}^d \int_{S_{k,u(k)}} \int_{R_k} p_{\alpha,1} \left(\frac{x_k -y_k}{(2 w_kt)^{\frac{1}{\alpha}}}\right) \frac{dx_kdy_k}{(2 w_kt)^{\frac{1}{\alpha}}}. 
\end{align*}
Next, let us identify the terms for which the contribution in the limit \eqref{eq:limit_discrete_K} is non-zero. This analysis is based on a recent one-dimensional result \cite[Theorem $1.1$]{Valverde_17}. First, assume that there are two distinct indices $i,j \in \{1, \dots, d\}$ for which $u(i) = u(j)=0$.  Then, 
\begin{align*}
\int_{S_u} \int_R p^w_{\alpha,d} \left(\frac{x-y}{t^{\frac{1}{\alpha}}}\right) \frac{dxdy}{t^{\frac{d}{\alpha}}} & = \prod_{k  =1,\,  k \notin \{i,j\}}^d \int_{S_{k,u(k)}} \int_{R_k} p_{\alpha,1} \left(\frac{x_k -y_k}{(2 w_kt)^{\frac{1}{\alpha}}}\right) \frac{dx_kdy_k}{(2 w_kt)^{\frac{1}{\alpha}}} \\
&\quad\quad \times \int_{R_i^c} \int_{R_i} p_{\alpha,1} \left(\frac{x_i -y_i}{(2 w_it)^{\frac{1}{\alpha}}}\right) \frac{dx_idy_i}{(2 w_it)^{\frac{1}{\alpha}}} \\
&\quad\quad \times \int_{R_j^c} \int_{R_j} p_{\alpha,1} \left(\frac{x_j -y_j}{(2 w_jt)^{\frac{1}{\alpha}}}\right) \frac{dx_jdy_j}{(2 w_jt)^{\frac{1}{\alpha}}}.
\end{align*}
Set,  for all $t>0$, 
\begin{align*}
A_i(t) =\int_{R_i^c} \int_{R_i} p_{\alpha,1} \left(\frac{x_i -y_i}{(2 w_it)^{\frac{1}{\alpha}}}\right) \frac{dx_idy_i}{(2 w_it)^{\frac{1}{\alpha}}}, \quad A_j(t) = \int_{R_j^c} \int_{R_j} p_{\alpha,1} \left(\frac{x_j -y_j}{(2w_jt)^{\frac{1}{\alpha}}}\right) \frac{dx_jdy_j}{(2w_jt)^{\frac{1}{\alpha}}}.
\end{align*}
Since both terms $A_i$ and $A_j$ tend to $0$ as $t$ tends to $0^+$ at the rate $t^{\frac{1}{\alpha}}$, it is clear that,  for such a $u \in \{0,1\}^d \setminus \{(1, \cdots, 1)\}$, 
\begin{align*}
\underset{t  \rightarrow 0^+}{\lim} \frac{1}{t^{\frac{1}{\alpha}}} \int_{S_u} \int_R p^w_{\alpha,d} \left(\frac{x-y}{t^{\frac{1}{\alpha}}}\right) \frac{dxdy}{t^{\frac{d}{\alpha}}} = 0.
\end{align*}
Then, it remains to analyze the terms with only one $0$.  Fix $i \in \{1,\cdots, d\}$ and consider the vector $u$ such that $u(i)=0$ and $u(k) = 1$, for all $k \in \{1, \cdots, d\}$ with $k \ne i$.  Then, 
\begin{align*}
\int_{S_u} \int_R p^w_{\alpha,d} \left(\frac{x-y}{t^{\frac{1}{\alpha}}}\right) \frac{dxdy}{t^{\frac{d}{\alpha}}} & = \prod_{k  =1,\,  k \ne i}^d \int_{R_k} \int_{R_k} p_{\alpha,1} \left(\frac{x_k -y_k}{(2 w_kt)^{\frac{1}{\alpha}}}\right) \frac{dx_kdy_k}{(2 w_kt)^{\frac{1}{\alpha}}} \\
&\quad\quad \times \int_{R_i^c} \int_{R_i} p_{\alpha,1} \left(\frac{x_i -y_i}{(2 w_it)^{\frac{1}{\alpha}}}\right) \frac{dx_idy_i}{(2 w_it)^{\frac{1}{\alpha}}}.
\end{align*}
Now, thanks to \cite[Theorem $1.1$]{Valverde_17}, 
\begin{align*}
\underset{t \rightarrow 0^+}{\lim} \frac{1}{t^{\frac{1}{\alpha}}}\int_{R_i^c} \int_{R_i} p_{\alpha,1} \left(\frac{x_i -y_i}{(2 w_it)^{\frac{1}{\alpha}}}\right) \frac{dx_idy_i}{(2 w_it)^{\frac{1}{\alpha}}} = (2w_i)^{\frac{1}{\alpha}} \frac{2}{\pi} \Gamma \left(1 - \frac{1}{\alpha}\right). 
\end{align*}
Moreover, 
\begin{align*}
\underset{t \rightarrow 0^+}{\lim} \prod_{k  =1,\,  k \ne i}^d \int_{R_k} \int_{R_k} p_{\alpha,1} \left(\frac{x_k -y_k}{(2 w_kt)^{\frac{1}{\alpha}}}\right) \frac{dx_kdy_k}{(2 w_kt)^{\frac{1}{\alpha}}} = 2^{d-1} \prod_{k =1, k \ne i}^d a_k.
\end{align*}
Thus,  for such a  vector $u$, 
\begin{align*}
\underset{t \rightarrow 0^+}{\lim} \frac{1}{t^{\frac{1}{\alpha}}} \int_{S_u} \int_R p^w_{\alpha,d} \left(\frac{x-y}{t^{\frac{1}{\alpha}}}\right) \frac{dxdy}{t^{\frac{d}{\alpha}}} = (2w_i)^{\frac{1}{\alpha}} \frac{2}{\pi} \Gamma \left(1 - \frac{1}{\alpha}\right) 2^{d-1} \prod_{k =1, k \ne i}^d a_k.
\end{align*}
Hence, 
\begin{align*}
\underset{t \rightarrow 0^+}{\lim} \dfrac{K_t^{\alpha,\operatorname{dis}} \left(R , R^c\right)}{t^{\frac{1}{\alpha}}} = \frac{2^d}{\pi} \Gamma \left(1 - \frac{1}{\alpha}\right) \sum_{i = 1}^d (2w_i)^{\frac{1}{\alpha}} \prod_{k =1, k \ne i}^d a_k. 
\end{align*}
This concludes the proof of the lemma.
\end{proof}
\noindent
The next lemma finds $\mathcal{P}^{\operatorname{dis}}_{\operatorname{cl}}(R)$ when $R$ is a hyperrectangle.

\begin{lem}\label{lem:discrete_perimeter}
Let $d \geq 1$ be an integer, let $\alpha \in (1,2)$ and let $\mathcal{P}^{\operatorname{dis}}_{\operatorname{cl}}$ be defined by \eqref{eq:classical} with $\sigma^{\operatorname{dis}}_\alpha$ given by \eqref{eq:norm_discrete}. Then, for all $R = R_1 \times \cdots \times R_d$ with $R_k = [-a_k , a_k]$, $a_k>0$, for $k \in \{1, \cdots, d\}$, 
\begin{align}\label{eq:perimeter_hyperrectangle_discrete}
\mathcal{P}^{\operatorname{dis}}_{\operatorname{cl}}(R) =2 \sum_{i = 1}^d (2w_i)^{\frac{1}{\alpha}} \prod_{k =1, k \ne i}^d 2a_k.
\end{align}  
\end{lem}

\begin{proof}
To prove \eqref{eq:perimeter_hyperrectangle_discrete}, let us use the variational representation formula for $\mathcal{P}^{\operatorname{dis}}_{\operatorname{cl}}$ provided by Proposition \ref{prop:nondeg_symmetric_equality} of the Appendix. Indeed, Proposition \ref{prop:nondeg_symmetric_equality} shows that $\mathcal{P}^{\operatorname{dis}}_{\alpha}(R) = \mathcal{P}^{\operatorname{dis}}_{\operatorname{cl}}(R)$ where
\begin{align}\label{eq:variational_representation_formula_discrete}
\mathcal{P}^{\operatorname{dis}}_{\alpha}(R) = \sup \left\{ \int_R \operatorname{div}(\varphi(x))dx, \, \varphi \in \mathcal{C}^1_c(\mathbb{R}^d , \mathbb{R}^d), \, \| \sigma^{\operatorname{dis},*}_\alpha(\varphi) \|_{\infty} \leq 1 \right\}. 
\end{align}
It thus remains to prove that $\mathcal{P}^{\operatorname{dis}}_{\alpha}(R)$ is given by the right-hand side of \eqref{eq:perimeter_hyperrectangle_discrete}. First, take $\varphi = (\varphi_1 , \cdots, \varphi_d)$ in $\mathcal{C}_c^{1} \left( \mathbb{R}^d, \mathbb{R}^d\right)$ such that $\| \sigma^{\operatorname{dis},*}_\alpha(\varphi)\|_{\infty} \leq 1$. This constraint implies that, for all $x \in \mathbb{R}^d$, 
\begin{align*}
\sum_{k = 1}^d \dfrac{|\varphi_k(x)|^{\frac{\alpha}{\alpha-1}}}{(2w_k)^{\frac{1}{\alpha-1}}} \leq 1 , 
\end{align*}
which implies, for all $x \in \mathbb{R}^d$ and all $ k \in \{1, \cdots, d\}$, 
\begin{align}\label{ineq:contraint_k}
|\varphi_k(x)| \leq \left(2w_k\right)^{\frac{1}{\alpha}}.
\end{align}
Now,
\begin{align*}
\int_{R} \operatorname{div}\left(\varphi(x)\right) dx = \sum_{k =1}^d \int_{R} \partial_k (\varphi_k)(x)dx & = \sum_{k = 1}^d \int_{-a_1}^{a_1} \dots \int_{-a_d}^{a_d} \bigg(\varphi_k(x_1, \dots, x_{k-1} , a_k , x_{k+1} , \dots, x_d) \\
&\quad\quad- \varphi_k(x_1, \dots, x_{k-1} , -a_k , x_{k+1} , \dots, x_d) \bigg) \\
&\quad\quad dx_1 \dots dx_{k -1} dx_{k+1} \dots dx_d. 
\end{align*}
So, using \eqref{ineq:contraint_k}, 
\begin{align*}
\left| \int_{R} \operatorname{div}\left(\varphi(x)\right) dx \right| \leq \sum_{k = 1}^d 2 \left(2 w_k\right)^{\frac{1}{\alpha}} \prod_{j =1\, j \ne k}^d \left(2 a_j\right). 
\end{align*}
Taking the supremum proves that 
\begin{align*}
\mathcal{P}^{\operatorname{dis}}_{\alpha}(R) \leq \sum_{k = 1}^d 2 \left(2 w_k\right)^{\frac{1}{\alpha}} \prod_{j =1\, j \ne k}^d \left(2 a_j\right). 
\end{align*}
For the reverse inequality, it is sufficient to take functions $(\varphi_1 , \dots, \varphi_d) \in \mathcal{C}_c^{1}(\mathbb{R}^d,\mathbb{R}^d)$ such that $\|\sigma^{\operatorname{dis},*}_\alpha(\varphi)\|_{\infty} \leq 1$ with, for all $ k \in \{1, \dots, d\}$ and all $x \in \mathbb{R}^d$, 
\begin{align*}
&\varphi_k \left(x_1 ,\dots,x_{k-1},a_k ,x_{k+1},\dots, x_d\right) = \varphi_k(a_k) = (2w_k)^{\frac{1}{\alpha}} , \\
&\varphi_k \left(x_1 ,\dots,x_{k-1}, -a_k ,x_{k+1},\dots, x_d\right) = \varphi_k(- a_k) = - (2w_k)^{\frac{1}{\alpha}}. 
\end{align*}
Then, 
\begin{align*}
\int_{R} \operatorname{div}\left(\varphi(x)\right) dx = \sum_{k =1}^d 2 \left(2 w_k\right)^{\frac{1}{\alpha}} \prod_{j = 1, \, j \ne k}^d (2a_j).
\end{align*}
The reverse inequality follows by taking, once again, the supremum. This concludes the proof of the lemma.
\end{proof}
\noindent
From the previous analysis, combining Lemma \ref{lem:limit_discrete_K} and Lemma \ref{lem:discrete_perimeter} together with the inequality \eqref{ineq:discrete_case}, one obtains the following sharp isoperimetric-type inequality in the symmetric non-degenerate $\alpha$-stable case, $\alpha \in (1,2)$, with a discrete spectral measure.  

\begin{prop}\label{prop:isoperimetry_discrete}
Let $\alpha \in (1,2)$ and let $\mathcal{P}^{\operatorname{dis}}_{\operatorname{cl}}$ be defined by \eqref{eq:classical} with $\sigma^{\operatorname{dis}}_\alpha$ given by \eqref{eq:norm_discrete}. Then, for all $E = E_1\times \dots \times E_d$ Borel subset of $\mathbb{R}^d$, $d \geq 1$, with finite Lebesgue measure and with $\mathcal{P}^{\operatorname{dis}}_{\operatorname{cl}}(E) <+ \infty$, 
\begin{align}
\mathcal{P}^{\operatorname{dis}}_{\operatorname{cl}}(E) \geq \mathcal{P}^{\operatorname{dis}}_{\operatorname{cl}}(R), 
\end{align}
where $R = R_1\times \dots \times R_d$ is such that $\mathcal{L}_d(E_k) = \mathcal{L}_d(R_k)$ and $R_k = [-a_k , a_k]$, with $a_k >0$, for all $k \in \{1, \dots, d\}$.
\end{prop}

\begin{proof}
The proof is now a direct application of the pseudo-Poincar\'e inequality \eqref{ineq:pseudo_poincar_gradient} with $p = 1$. Indeed, for all $t>0$, 
\begin{align*}
\mathcal{P}^{\operatorname{dis}}_{\operatorname{cl}}(E) \frac{ \bbe |Y_\alpha|}{2} \geq \dfrac{K_t^{\alpha,\operatorname{dis}} \left(E , E^c\right)}{t^{\frac{1}{\alpha}}} \geq \dfrac{K_t^{\alpha,\operatorname{dis}} \left(R , R^c\right)}{t^{\frac{1}{\alpha}}}. 
\end{align*}
Thus, by passing to the limit, 
\begin{align*}
\mathcal{P}^{\operatorname{dis}}_{\operatorname{cl}}(E) \geq \frac{2}{\bbe |Y_\alpha|} \underset{t \rightarrow 0^+ }{\lim}\dfrac{K_t^{\alpha,\operatorname{dis}} \left(R , R^c\right)}{t^{\frac{1}{\alpha}}} =  \frac{2}{\bbe |Y_\alpha|} \frac{1}{\pi} \Gamma \left(1 - \frac{1}{\alpha}\right) \mathcal{P}^{\operatorname{dis}}_{\operatorname{cl}}(R) = \mathcal{P}^{\operatorname{dis}}_{\operatorname{cl}}(R). 
\end{align*}
This concludes the proof of the proposition. 
\end{proof}
\noindent
At this point, in analogy with the previous analysis and in view of the sharp anisotropic isoperimetric inequality \eqref{ineq:sharp_anisotropic_iso}, it seems reasonable to wonder whether or not, for all $t>0$ and all $E\in \mathcal{B}(\mathbb{R}^d)$ such that $\mathcal{L}_d(E) = \mathcal{L}_d(\mathring{K}_\alpha)$, 
\begin{align}\label{ineq:first_conjecture}
K_t^\alpha(E,E^c) \geq K_t^\alpha(\mathring{K}_\alpha,\mathring{K}_\alpha^c)
\end{align}
and
\begin{align}\label{ineq:second_conjecture}
 \underset{t \rightarrow 0^+}{\lim} \dfrac{K_t^\alpha(\mathring{K}_\alpha,\mathring{K}_\alpha^c)}{t^{\frac{1}{\alpha}}} = \frac{\Gamma(1- \frac{1}{\alpha})}{\pi} \mathcal{P}_\alpha^{\operatorname{var}} (\mathring{K}_\alpha),
\end{align}
for all non-degenerate symmetric $\alpha$-stable probability measure $\mu_\alpha$ with $\alpha \in (1,2)$. From the proof of \eqref{ineq:first_conjecture}, in the rotationally invariant case,  an intermediate question could be the following:\\ 
\\
Is the Lebesgue density $p_\alpha$ of the non-degenerate symmetric $\alpha$-stable probability measure $\mu_\alpha$ convex symmetric with respect to the norm $\sigma_\alpha^*$ in the sense of \cite{Lions_97,Van_Schaftingen_Jean_06}?\\
\\
First, let us recall basic facts from convex symmetrization (see \cite{Lions_97,Van_Schaftingen_Jean_06}). Let $\sigma_\alpha^*$ be the dual norm for which the convex body $\mathring{K}_\alpha$ is the unit ball.  Then, for all $E\in\mathcal{B}(\mathbb{R}^d)$, let $E^\ast$ be the set defined by 
\begin{align}\label{eq:convex_symmetrization_set}
E^\ast = \bigg\{ x \in \mathbb{R}^d, \, \sigma_\alpha^*(x) \leq  \left(\frac{\mathcal{L}_d(E)}{\mathcal{L}_d(\mathring{K}_\alpha)}  \right)^{\frac{1}{d}}\bigg\}. 
\end{align}
Clearly, $\mathcal{L}_d(E^\ast) = \mathcal{L}_d(E)$ and $(\mathring{K}_\alpha)^\ast = \mathring{K}_\alpha$.  In particular, if $E \in \mathcal{B}(\mathbb{R}^d)$ is such that $\mathcal{L}_d(E) = \mathcal{L}_d(\mathring{K}_\alpha)$ then $E^\ast = \mathring{K}_\alpha$.  At the level of functions, for all Borel measurable $f$ from $\mathbb{R}^d$ to $\bbr$, let $f^\ast$ be defined, for all $x \in \mathbb{R}^d$, by 
\begin{align}\label{eq:convex_symmetrization_function}
f^\ast(x) = \sup \left\{ c \in \bbr_+ \cup \{ + \infty\}: \,  x \in \{|f| >c\}^{\ast} \right\},
\end{align}
where $\{ |f| > c\} = \{x \in \mathbb{R}^d, \ |f(x)|> c\}$. From the definition, it is clear that $(\bbone_E)^\ast = \bbone_{E^\ast}$. As already put forward, this type of symmetry appears naturally when one is interested in the sharp Sobolev and isoperimetric inequalities with non-Euclidean norms. Indeed,  the extremal functions for these inequalities enjoy this type of symmetry (see \cite[Theorems $2$ and $3$ as well as equation $(13)$ in there]{cordero_nazaret_villani}). Now, note that the Fourier transform of $p_\alpha$ is a fixed point of the convex symmetrization associated with the norm $\sigma_\alpha$ since, for all $\xi \in \mathbb{R}^d$, 
\begin{align*}
\mathcal{F}(p_\alpha)(\xi) = \exp \left( - \sigma_\alpha(\xi)^\alpha\right). 
\end{align*}
Moreover, $p_\alpha$ admits a particular integral representation.

\begin{lem}\label{lem:special_rep_palpha}
Let $\alpha \in (1,2)$, let $p_\alpha$ be the positive Lebesgue density of the non-degenerate symmetric $\alpha$-stable probability measure $\mu_\alpha$ on $\mathbb{R}^d$, $d \geq 1$, and let $\mu_{K_\alpha}$ be the cone measure on $\partial K_\alpha$ associated with the convex set $K_\alpha$. Then, for all $x \in \mathbb{R}^d$, 
\begin{align}\label{eq:Generalized_Bessel_rep}
p_\alpha(x) =  \int_{0}^{+\infty} \exp(-r^\alpha) J_{\alpha}(r,x)r^{d-1} \frac{dr}{(2 \pi)^d},
\end{align}
where, for all $r>0$ and all $x \in \mathbb{R}^d$, 
\begin{align}\label{eq:representation_formula_fourier_cone_measure}
J_\alpha(r,x) = d \mathcal{L}_d(K_\alpha) \int_{\partial K_\alpha} e^{i r \langle x ; y \rangle} \mu_{K_\alpha}(dy).
\end{align}
\end{lem}

\begin{proof}
Recall that the Fourier transform of $p_\alpha$ is given, for all $\xi \in \mathbb{R}^d$, by
\begin{align}\label{eq:fourier_transform_p_alpha}
\mathcal{F}(p_\alpha)(\xi) = \exp \left( - \sigma_\alpha(\xi)^\alpha\right) = u_\alpha\left(\sigma_\alpha(\xi)\right), 
\end{align} 
where $u_\alpha$ is the one-dimensional decreasing function defined, for all $r \geq 0$, by
\begin{align}
u_\alpha(r) = \exp\left( - r^\alpha\right). 
\end{align}
From \eqref{eq:fourier_transform_p_alpha}, the Fourier transform of $p_\alpha$ is convex symmetric with respect to the norm $\sigma_\alpha$. Using the Fourier inversion formula, polar coordinates and a change of variables, for all $x \in \mathbb{R}^d$, 
\begin{align*}
p_\alpha(x) & = \int_{\mathbb{R}^d} \exp\left( - \sigma_\alpha(\xi)^\alpha \right) \exp \left( i \langle x ; \xi \rangle\right) \frac{d \xi}{(2 \pi)^d} \\
& = \int_{(0, +\infty) \times \mathbb{S}^{d-1}} \exp\left( - r^\alpha \sigma_\alpha(\theta)^\alpha \right)\exp \left( i r \langle x ; \theta \rangle\right) \frac{r^{d-1} dr \sigma_L(d \theta)}{(2 \pi)^d} \\
& = \int_{0}^{+\infty} u_\alpha(r) \left( \int_{\mathbb{S}^{d-1}}  \exp \left( i r \left\langle x ; \frac{\theta}{\sigma_\alpha(\theta)} \right\rangle\right) \dfrac{\sigma_L(d \theta)}{(\sigma_\alpha(\theta))^d} \right) \frac{r^{d-1}dr}{(2 \pi)^d}.
\end{align*}
Next, let $J_\alpha$ be the function defined, for all $x \in \mathbb{R}^d$ and all $r>0$, by
\begin{align*}
J_\alpha(r,x) = \int_{\mathbb{S}^{d-1}}  \exp \left( i r \left\langle x ; \frac{\theta}{\sigma_\alpha(\theta)} \right\rangle\right) \dfrac{\sigma_L(d \theta)}{(\sigma_\alpha(\theta))^d}.
\end{align*}
Using \eqref{eq:representation_formula_integral_cone_measure_2}, for all $x\in \mathbb{R}^d$ and all $r>0$, 
\begin{align}\label{eq:representation_formula_fourier_cone_measure}
J_\alpha(r,x) = d \mathcal{L}_d(K_\alpha) \int_{\partial K_\alpha} e^{i r \langle x ; y \rangle} \mu_{K_\alpha}(dy).
\end{align}
This concludes the proof.
\end{proof}
\noindent
In certain instances of the spectal measure $\lambda_1$, the cone measure $\mu_{K_\alpha}$ admits a useful probabilistic representation (see, \cite{SZ_90,RR_91,PTT}).  When $\lambda_1 = \sum_{k = 1}^d \delta_{\pm e_k}/2$, the convex body $K_\alpha^{\operatorname{dis}}$ is $\mathbb{B}_\alpha^d$ whose geometry has been investigated in depth in the literature.  In particular, the cone measure can be realized as the law of the normalized random vector $Z_\alpha$ defined by 
\begin{align}
Z_\alpha : = \left( \frac{V_1}{\sigma_\alpha(V)} , \dots, \frac{V_d}{\sigma_\alpha(V)} \right),
\end{align} 
where $V = (V_1, \dots, V_d)$ and $\{ V_1, \dots, V_d \}$ is a collection of independent and identically distributed random variables with law given by the density function $h_\alpha$ defined, for all $v \in \bbr$,  by
\begin{align}\label{eq:density_V}
h_\alpha(v) = \dfrac{1}{2 (\alpha)^{\frac{1}{\alpha}} \Gamma(1+\frac{1}{\alpha})} \exp \left( - \frac{|v|^\alpha}{\alpha}\right). 
\end{align}
Thus, in this discrete situation, for all $r>0$ and all $x \in \mathbb{R}^d$, 
\begin{align*}
J_{\alpha}(r,x) = d \mathcal{L}_d(\mathbb{B}_\alpha^d) \bbe e^{i r \langle x ; Z_\alpha \rangle},
\end{align*}
so that, $p_{\alpha,d}$ is given, for all $x \in \mathbb{R}^d$, by 
\begin{align*}
p_{\alpha,d}(x) = \frac{d \mathcal{L}_d(\mathbb{B}_\alpha^d)}{(2\pi)^d} \int_0^{+\infty} e^{-r^\alpha} \bbe e^{i r \langle x ; Z_\alpha \rangle} r^{d-1} dr.
\end{align*}
Recall that we are interested in the computation of the following limit:
\begin{align*}
L(\mathbb{B}_{\alpha/(\alpha-1)}^d) = \underset{t \rightarrow 0^+}{\lim} \dfrac{K^{\alpha, \operatorname{dis}}_t(\mathbb{B}_{\alpha/(\alpha-1)}^d,(\mathbb{B}_{\alpha/(\alpha-1)}^d)^c)}{t^{\frac{1}{\alpha}}}.
\end{align*}
Next, let us give a proof of the classical Gaussian result which relies only on the stability property of the Gaussian probability measure without using invariance by rotations neither its product structure.  

\begin{prop}\label{prop:gaussian_short_proof}
Let $\tilde{\gamma}$ be the centered Gaussian probability measure on $\mathbb{R}^d$, $d \geq 1$, with covariance matrix $2I_d$.  Then, 
\begin{align*}
\underset{t \rightarrow 0^+}{\lim} \dfrac{K^{2, \operatorname{dis}}_t(\mathbb{B}_{2}^d,(\mathbb{B}_{2}^d)^c)}{t^{\frac{1}{2}}} = \frac{d \mathcal{L}_d( \mathbb{B}_2^d)}{\sqrt{\pi}}.
\end{align*}
\end{prop}

\begin{proof}
For all $t>0$, 
\begin{align}\label{eq:def_K_Functional_discrete}
K^{2, \operatorname{dis}}_t(\mathbb{B}_{2}^d,(\mathbb{B}_{2}^d)^c) = \int_{(\mathbb{B}_{2}^d)^c} \left( \int_{\mathbb{R}^d} \tilde{\gamma}\left(\frac{x-y}{t^{\frac{1}{2}}}\right) \bbone_{\mathbb{B}_{2}^d}(y) \frac{dy}{t^{\frac{d}{2}}} \right) dx.
\end{align}
Let $\mu_{\mathbb{B}_{2}^d}$ be the cone probability measure defined on $\partial \mathbb{B}_{2}^d$. Using \eqref{eq:polar_coordinate_cone_measure}, for all $t>0$, 
\begin{align*}
K^{2, \operatorname{dis}}_t(\mathbb{B}_{2}^d,(\mathbb{B}_{2}^d)^c) & = d \mathcal{L}_d(\mathbb{B}_{2}^d) \int_{\partial \mathbb{B}_{2}^d \times (1,+\infty)}  \left( \int_{\mathbb{R}^d} \tilde{\gamma}\left(\frac{s\omega-y}{t^{\frac{1}{2}}}\right) \bbone_{\mathbb{B}_{2}^d}(y) \frac{dy}{t^{\frac{d}{2}}} \right)\\
&\quad\quad \times s^{d-1} ds \mu_{\mathbb{B}_{2}^d}(d\omega) \\
& =  d \mathcal{L}_d(\mathbb{B}_{2}^d) \int_{\partial \mathbb{B}_{2}^d} \mu_{\mathbb{B}_{2}^d}(d\omega)\\
&\quad\quad \times \left(\int_{1}^{+\infty} \tilde{\gamma}\left(z:\,  \|s \omega + t^{\frac{1}{2}} z \| \leq 1 \right)  s^{d-1} ds\right).
\end{align*}
Next, fix $\omega \in \partial \mathbb{B}_{2}^d$ and $z \in \mathbb{R}^d$. Let us compute the following limit:
\begin{align*}
g_2(\omega, z) = \underset{t \rightarrow 0^+}{\lim} \frac{1}{t^{\frac{1}{2}}} \int_{1}^{+\infty} \bbone_{\|s \omega + t^{\frac{1}{2}} z \| \leq 1} s^{d-1} ds.
\end{align*}
First, observe that if $\langle \omega; z \rangle \geq 0$, then $g_2(\omega, z) = 0$.  Indeed, by the Pythagorean theorem, for all $s \in [1,+\infty)$ and all $t>0$, 
\begin{align*}
\|s \omega + t^{\frac{1}{2}} z \|^2 = s^2 + 2 s \langle \omega; z \rangle t^{\frac{1}{2}} + t \|z\|^2 \geq 1. 
\end{align*}
Thus, assume that $\langle \omega ; z \rangle<0$ and let $F_2$ be the function defined on $[1, +\infty)$ by 
\begin{align}\label{def:quadratic_function}
F_2(s) := s^2 + 2 s \langle \omega; z \rangle t^{\frac{1}{2}} + t \|z\|^2 - 1. 
\end{align} 
Let $s_\pm(t)$ be defined by
\begin{align*}
s_\pm(t) = - \langle \omega ; z \rangle \sqrt{t} \pm \sqrt{1+t \left(\langle \omega  ; z \rangle^2 - \|z\|^2 \right)}.
\end{align*}
Recall that $\omega$ and $z$ are fixed and that we are interested in $t$ tending to $0^+$.  Thus, $s_\pm(t)$ are well-defined, with $s_-(t)$ negative and $s_+(t)$ bigger than $1$.  Then, for $t$ small enough, 
\begin{align*}
\int_{1}^{+\infty} \bbone_{\|s \omega + t^{\frac{1}{2}} z \| \leq 1} s^{d-1} ds = \int_{1}^{s_+(t)} s^{d-1} ds = \frac{1}{d} \left(s_+(t)^d-1\right).  
\end{align*}
By standard computations, for all $\omega \in \partial \mathbb{B}^d_2$ and all $z \in \mathbb{R}^d$, 
\begin{align*}
g_2(\omega,z) = - \langle \omega ; z \rangle \bbone_{\langle \omega ; z \rangle \leq 0}. 
\end{align*}
Now, for all $\omega \in \partial \mathbb{B}_2^d$, under the probability measure $\tilde{\gamma}$, $\langle \omega ; z \rangle$ is a centered Gaussian random variable with variance $2$. Thus, 
\begin{align*} 
\int_{\mathbb{R}^d} g_2(\omega,z) \tilde{\gamma}(dz) = -  \sqrt{2} \int_{-\infty}^0 x \exp \left(- \frac{x^2}{2}\right) \frac{dx}{\sqrt{2\pi}} = \frac{1}{\sqrt{\pi}}.
\end{align*}
This concludes the proof of the proposition. 
\end{proof}

\begin{prop}\label{prop:limit_discrete_case}
Let $\alpha \in (1,2)$ and let $\mu_{\alpha,d}$ be the non-degenerate symmetric $\alpha$-stable product probability measure on $\mathbb{R}^d$, $d \geq 1$. Then, 
\begin{align}\label{eq:limit_discrete_case}
\underset{t \rightarrow 0^+}{\lim} \dfrac{K_t^{\alpha, \operatorname{dis}}(\mathbb{B}^d_{\alpha/(\alpha-1)}, (\mathbb{B}^d_{\alpha/(\alpha-1)})^c)}{t^{\frac{1}{\alpha}}} = \dfrac{\Gamma(1- \frac{1}{\alpha})}{\pi} d \mathcal{L}_d(\mathbb{B}^d_{\alpha/(\alpha-1)}) = \dfrac{\Gamma(1- \frac{1}{\alpha})}{\pi} \mathcal{P}_\alpha^{\operatorname{var}} \left(\mathbb{B}^d_{\alpha/(\alpha-1)}\right).
\end{align}
\end{prop}

\begin{proof}
For all $t>0$, 
\begin{align}\label{eq:def_K_Functional_discrete}
K^{\alpha, \operatorname{dis}}_t(\mathbb{B}_{\alpha/(\alpha-1)}^d,(\mathbb{B}_{\alpha/(\alpha-1)}^d)^c) = \int_{(\mathbb{B}^d_{\alpha/(\alpha-1)})^c} \left( \int_{\mathbb{R}^d} p_{\alpha,d}\left(\frac{x-y}{t^{\frac{1}{\alpha}}}\right) \bbone_{\mathbb{B}^d_{\alpha/(\alpha-1)}}(y) \frac{dy}{t^{\frac{d}{\alpha}}} \right) dx,
\end{align}
where $p_{\alpha,d}$ is the Lebesgue density of the probability measure $\mu_{\alpha,d}$. Let $\mu_{\mathbb{B}_{\alpha/(\alpha-1)}^d}$ be the cone probability measure defined on $\partial \mathbb{B}_{\alpha/(\alpha-1)}^d$.  Recall that the following polar-type coordinate integration formula holds true (see, e.g., \cite[Proposition $1$]{Naor_Romik_2003}): for all $f \in L^1(\mathbb{R}^d,dx)$, 
\begin{align*}
\int_{\mathbb{R}^d} f(x) dx = d \mathcal{L}_d(\mathbb{B}_{\alpha/(\alpha-1)}^d) \int_0^{+\infty} s^{d-1} \left(\int_{\partial \mathbb{B}_{\alpha/(\alpha-1)}^d} f(s \omega) \mu_{\mathbb{B}_{\alpha/(\alpha-1)}^d}(d\omega) \right) ds.
\end{align*}
Thus, for all $t>0$, 
\begin{align*}
K^{\alpha, \operatorname{dis}}_t(\mathbb{B}_{\alpha/(\alpha-1)}^d,(\mathbb{B}_{\alpha/(\alpha-1)}^d)^c) & = d \mathcal{L}_d(\mathbb{B}_{\alpha/(\alpha-1)}^d) \int_{\partial \mathbb{B}_{\alpha/(\alpha-1)}^d \times (1,+\infty)}  \left( \int_{\mathbb{R}^d} p_{\alpha,d}\left(\frac{s\omega-y}{t^{\frac{1}{\alpha}}}\right) \bbone_{\mathbb{B}^d_{\alpha/(\alpha-1)}}(y) \frac{dy}{t^{\frac{d}{\alpha}}} \right)\\
&\quad\quad \times s^{d-1} ds \mu_{\mathbb{B}_{\alpha/(\alpha-1)}^d}(d\omega) \\
& =  d \mathcal{L}_d(\mathbb{B}_{\alpha/(\alpha-1)}^d) \int_{\partial \mathbb{B}_{\alpha/(\alpha-1)}^d} \mu_{\mathbb{B}_{\alpha/(\alpha-1)}^d}(d\omega)\\
&\quad\quad \times \left(\int_{1}^{+\infty} \mu_{\alpha,d}\left(z:\,  \|s \omega + t^{\frac{1}{\alpha}} z \|_{\alpha/(\alpha-1)} \leq 1 \right)  s^{d-1} ds\right).
\end{align*}
Next, fix $\omega \in \partial \mathbb{B}^d_{\alpha/(\alpha-1)}$ and $z \in \mathbb{R}^d$, and let us compute the following limit:
\begin{align*}
g_\alpha(\omega, z) = \underset{t \rightarrow 0^+}{\lim} \frac{1}{t^{\frac{1}{\alpha}}} \int_{1}^{+\infty} \bbone_{\|s \omega + t^{\frac{1}{\alpha}} z \|_{\alpha/(\alpha-1)} \leq 1} s^{d-1} ds.
\end{align*}
Now, recall that $\alpha/(\alpha-1) \in (2,+\infty)$ since $\alpha \in (1,2)$.  Let $\nu_{\mathbb{B}^d_{\alpha/(\alpha-1)}}$ be a unit normal vector given, for all $\omega \in \partial\mathbb{B}^d_{\alpha/(\alpha-1)}$, by 
\begin{align*}
\nu_{\mathbb{B}^d_{\alpha/(\alpha-1)}}(\omega) = \dfrac{\nabla \left(\| \cdot \|_{\alpha/(\alpha-1)}\right)(\omega)}{\| \nabla \left(\| \cdot \|_{\alpha/(\alpha-1)}\right)(\omega)\|} = \left( \dfrac{\omega_1 |\omega_1|^{\frac{2-\alpha}{\alpha-1}}}{\left(\sum_{k = 1}^d |\omega_k|^{\frac{2}{\alpha-1}}\right)^{\frac{1}{2}}} ,\dots,  \dfrac{\omega_d |\omega_d|^{\frac{2-\alpha}{\alpha-1}}}{\left(\sum_{k = 1}^d |\omega_k|^{\frac{2}{\alpha-1}}\right)^{\frac{1}{2}}} \right),
\end{align*}
and let $\Omega = (\omega_1 |\omega_1|^{\frac{2-\alpha}{\alpha-1}}, \dots, \omega_d |\omega_d|^{\frac{2-\alpha}{\alpha-1}})$.  Note that $\|\Omega\|_{\alpha} = 1$. Now, for all $\omega \in \partial \mathbb{B}^d_{\alpha/(\alpha-1)}$, all $s \in [1, +\infty)$, all $t >0$ and all $z \in \mathbb{R}^d$, 
\begin{align}\label{eq:Taylor_order_one}
\| s \omega + t^{\frac{1}{\alpha}} z  \|_{\frac{\alpha}{\alpha-1}} = s +  t^{\frac{1}{\alpha}} \langle z ; \Omega \rangle + r\left(t^{\frac{1}{\alpha}}z\right), 
\end{align}
where $r$ is such that
\begin{align*}
\underset{\delta \rightarrow 0^+}{\lim} \underset{0<\|z\|\leq \delta}{\sup} \frac{|r(z)|}{\|z\|} = 0. 
\end{align*}
In order to establish \eqref{eq:Taylor_order_one},  note that, for all $x \in \mathbb{R}^d$ with $x \ne 0$, 
\begin{align*}
\nabla \left(  \| \cdot \|_{\alpha/(\alpha-1)}  \right)(x) =  \left(  \frac{x_1 |x_1|^{\frac{2- \alpha}{\alpha-1}}}{\left( \sum_{k = 1}^d |x_k|^{\frac{\alpha}{\alpha-1}}\right)^{\frac{1}{\alpha}}}  , \dots , \frac{x_d |x_d|^{\frac{2- \alpha}{\alpha-1}}}{\left( \sum_{k = 1}^d |x_k|^{\frac{\alpha}{\alpha-1}}\right)^{\frac{1}{\alpha}}}\right).
\end{align*}
As a first approximation, let $F_\alpha$ be defined, for all $s \in [1,+\infty)$,  by
\begin{align*}
F_\alpha(s) := s +  t^{\frac{1}{\alpha}} \langle z ; \Omega \rangle. 
\end{align*}
Clearly, if $\langle \Omega ; z \rangle \geq 0$,  then $F_\alpha(s) \geq 1$ for all $s \in [1,+\infty)$. Next, let us assume that $\langle \Omega ; z \rangle < 0$.  Then, let $s_\alpha(t)$ be given, for all $t>0$, by
\begin{align*}
s_\alpha(t) =1-t^{\frac{1}{\alpha}} \langle z ; \Omega\rangle >1.
\end{align*}
Using this first order approximation of $\| s \omega + t^{\frac{1}{\alpha}} z  \|_{\frac{\alpha}{\alpha-1}}$,  
\begin{align*}
\tilde{g}_{\alpha}(\omega,z) = \underset{t \rightarrow 0^+}{\lim} \frac{1}{t^{\frac{1}{\alpha}}} \int_1^{+\infty} \bbone_{\{F_\alpha(s) \leq 1\}} s^{d-1}ds = - \langle z ; \Omega \rangle \bbone_{\langle z ; \Omega \rangle \leq 0}.  
\end{align*}
Thus, 
\begin{align*}
\tilde{L}(\mathbb{B}^d_{\alpha/(\alpha-1)}) : = - d \mathcal{L}_d(\mathbb{B}^d_{\alpha / (\alpha-1)}) \int_{\partial \mathbb{B}^d_{\alpha / (\alpha-1)}} \mu_{\mathbb{B}^d_{\alpha / (\alpha-1)}}(d\omega) \int_{\mathbb{R}^d} \langle z ; \Omega \rangle \bbone_{\langle z ; \Omega \rangle \leq 0}\mu_{\alpha,d}(dz).
\end{align*}
Now, under $\mu_{\alpha,d}$, $\langle z ; \Omega \rangle$ is a symmetric $\alpha$-stable random variable with characteristic function given, for all $t \in \bbr$, by
\begin{align*}
\int_{\mathbb{R}^d} \exp\left( i  t \langle \Omega  ; z \rangle \right) \mu_{\alpha,d}(dz) = \exp \left( -  \|\Omega\|^\alpha_\alpha   |t|^\alpha  \right) = \exp \left( -   |t|^\alpha  \right).
\end{align*}
Then,  by subordination, 
\begin{align*}
\tilde{L}(\mathbb{B}^d_{\alpha/(\alpha-1)}) =  - d \mathcal{L}_d(\mathbb{B}^d_{\alpha/(\alpha-1)}) \bbe Y_\alpha \bbone_{Y_\alpha \leq 0}  = \frac{\Gamma(1  -\frac{1}{\alpha})}{\pi} d \mathcal{L}_d(\mathbb{B}^d_{\alpha/(\alpha-1)}).
\end{align*}
In order to finish the proof, let us show that one can indeed replace $\| s \omega + t^{\frac{1}{\alpha}} z  \|_{\frac{\alpha}{\alpha-1}}$ by $F_\alpha(s)$ in the computation of $g_\alpha(\omega,z)$.  Namely, let us prove that 
\begin{align*}
\underset{t \rightarrow 0^+}{\lim} \frac{1}{t^{\frac{1}{\alpha}}} \int_{1}^{+\infty} s^{d-1} \left( \bbone_{\| s \omega + t^{\frac{1}{\alpha} }z\|_{\alpha/(\alpha-1)} \leq 1} - \bbone_{F_\alpha(s) \leq 1} \right)ds = 0.
\end{align*}
Towards this goal, let us compute the Hessian matrix of the map $x \mapsto \|x\|_{\alpha/(\alpha-1)}$. For all $x \in \mathbb{R}^d$ with $x \ne 0$ and all $j,k \in \{1, \dots, d\}$ with $j \ne k$, 
\begin{align*}
\partial_j \partial_k \left( \| \cdot \|_{\alpha/(\alpha-1)}\right)(x) = - \frac{1}{\alpha-1} \dfrac{x_j |x_j|^{\frac{2-\alpha}{\alpha-1}}x_k |x_k|^{\frac{2-\alpha}{\alpha-1}}}{\|x\|^{\frac{\alpha+1}{\alpha-1}}_{\alpha/(\alpha-1)}},
\end{align*}
and, for all $k \in \{1, \dots, d\}$, 
\begin{align*}
\partial^2_k \left( \| \cdot \|_{\alpha/(\alpha-1)}\right)(x) = \frac{1}{\alpha-1} \left( \dfrac{|x_k|^{\frac{2-\alpha}{\alpha-1}}}{\|x\|^{\frac{1}{\alpha-1}}_{\alpha/(\alpha-1)}} - \dfrac{(x_k |x_k|^{\frac{2-\alpha}{\alpha-1}})^2}{\|x\|^{\frac{\alpha+1}{\alpha-1}}_{\alpha/(\alpha-1)}}  \right). 
\end{align*}
Now, from Taylor's formula of order $2$, 
\begin{align*}
\|s \omega + t^{\frac{1}{\alpha}} z  \|_{\alpha/(\alpha-1)} = s +t^{\frac{1}{\alpha}} \langle z ; \Omega \rangle + t^{\frac{2}{\alpha}} \int_0^1 \int_0^1 \langle  z  ; \operatorname{Hess}(\|\cdot\|_{\alpha/(\alpha-1)})\left(s \omega + t^{\frac{1}{\alpha}} z r u\right) z \rangle rdr du. 
\end{align*}
But, for all $s \in [1,+\infty)$, all $\omega \in \partial \mathbb{B}^d_{\alpha/(\alpha-1)}$, all $z \in \mathbb{R}^d$ and all $r,u \in [0,1]$, 
\begin{align*}
\left\langle  z  ; \operatorname{Hess}(\|\cdot\|_{\alpha/(\alpha-1)})\left(s \omega + t^{\frac{1}{\alpha}} z r u\right) z \right\rangle & = \frac{1}{(\alpha-1) \|s \omega + t^{\frac{1}{\alpha}} z r u\|^{\frac{1}{\alpha-1}}_{\alpha/(\alpha-1)}} \bigg( \sum_{j = 1}^d z_j^2 |s \omega_j + t^{\frac{1}{\alpha}} z_j r u|^{\frac{2-\alpha}{\alpha-1}}  \\
&\quad\quad-  \dfrac{\bigg(\sum_{j=1}^d z_j (s \omega_j + t^{\frac{1}{\alpha}} z_j r u) |s \omega_j + t^{\frac{1}{\alpha}} z_j r u|^{\frac{2-\alpha}{\alpha-1}}\bigg)^2}{\|s \omega + t^{\frac{1}{\alpha}} z r u\|^{\frac{\alpha}{\alpha-1}}_{\alpha/(\alpha-1)}}  \bigg) .
\end{align*}
Set $x = s \omega + t^{\frac{1}{\alpha}} z r u$. Then, by the Cauchy-Schwarz inequality, 
\begin{align*}
\left| \sum_{k=1}^d z_k x_k |x_k|^{\frac{2-\alpha}{\alpha-1}}   \right|^2 & \leq \left( \sum_{k = 1}^d z_k^2 |x_k|^{\frac{2-\alpha}{\alpha-1}} \right) \left( \sum_{k=1}^d |x_k|^{\frac{\alpha}{\alpha-1}} \right) \\
& \leq \left( \sum_{k = 1}^d z_k^2 |x_k|^{\frac{2-\alpha}{\alpha-1}} \right) \|x\|^{\frac{\alpha}{\alpha-1}}_{\alpha/(\alpha-1)}. 
\end{align*}
Thus,  
\begin{align*}
0 \leq \left\langle  z  ; \operatorname{Hess}(\|\cdot\|_{\alpha/(\alpha-1)})\left(s \omega + t^{\frac{1}{\alpha}} z r u\right) z \right\rangle \leq \frac{1}{(\alpha-1) \|x\|^{\frac{1}{\alpha-1}}_{\alpha/(\alpha-1)}} \sum_{k=1}^d z_k^2 |x_k|^{\frac{2- \alpha}{\alpha-1}}.
\end{align*}
Hence, for $0<t < t_0(\alpha,z)$ and $s \in [1,+\infty)$, 
\begin{align}\label{eq:error_bound_t_small_enough}
| \|s\omega + t^{\frac{1}{\alpha}}z\|_{\alpha/(\alpha-1)} - F_\alpha(s) | \leq c\frac{t^{\frac{2}{\alpha}}}{s} \|z\|^2 \leq c t^{\frac{2}{\alpha}} \|z\|^2,
\end{align}
where $c>0$ is a numerical constant.  Then,  for all $s \in [1, +\infty)$, 
\begin{align*}
\bbone_{ \{   F_\alpha(s) + c t^{\frac{2}{\alpha}} \|z\|^2 \leq 1 \} } \leq \bbone_{\{  \|s\omega + t^{\frac{1}{\alpha}}z\|_{\alpha/(\alpha-1)} \leq 1  \}} \leq \bbone_{  \{ F_\alpha(s) - c t^{\frac{2}{\alpha}} \|z\|^2 \leq 1  \} }.
\end{align*}
But, 
\begin{align*}
g_{\alpha,\pm}(\omega , z)= \underset{t \rightarrow 0^+}{\lim} \frac{1}{t^{\frac{1}{\alpha}}} \int_1^{+\infty} \bbone_{ \{ F_\alpha(s) \pm c t^{\frac{2}{\alpha}} \|z\|^2 \leq 1 \}} s^{d-1}ds  = - \langle z ; \Omega \rangle \bbone_{\langle z ; \Omega \rangle \leq 0}.
\end{align*}
This concludes the proof of the proposition. 
\end{proof}
\noindent
Now, let us investigate the limit in the general non-degenerate symmetric case.  The main difference with the previous discrete case is that one lacks (a priori) an analytical formula for the dual norm $\sigma_\alpha^*$. 

\begin{prop}\label{prop:limit_general_case}
Let $\alpha \in (1,2)$ and let $\mu_\alpha$ be a non-degenerate symmetric $\alpha$-stable probability measure on $\mathbb{R}^d$, $d \geq 1$. Then, 
\begin{align}\label{eq:limit_general_case}
\underset{t \rightarrow 0^+}{\lim} \dfrac{K_t^{\alpha}(\mathring{K}_\alpha, (\mathring{K}_\alpha)^c)}{t^{\frac{1}{\alpha}}} = \dfrac{\Gamma(1- \frac{1}{\alpha})}{\pi} d \mathcal{L}_d(\mathring{K}_\alpha) = \dfrac{\Gamma(1- \frac{1}{\alpha})}{\pi} \mathcal{P}_\alpha^{\operatorname{var}} \left(\mathring{K}_\alpha\right).
\end{align}
\end{prop}

\begin{proof}
The beginning of the proof is similar to the one of Proposition \ref{prop:limit_discrete_case}. Thanks to the polar coordinate integration formula \eqref{eq:polar_coordinate_cone_measure},
\begin{align*}
K_t^{\alpha}(\mathring{K}_\alpha, (\mathring{K}_\alpha)^c) = d \mathcal{L}_d(\mathring{K}_\alpha) \int_{\partial \mathring{K}_\alpha \times (1,+\infty)} \mu_{\mathring{K}_\alpha}(dy) s^{d-1} ds \int_{\mathbb{R}^d} p_\alpha(z)dz \bbone_{\{ \sigma_\alpha^*(s y + t^{\frac{1}{\alpha}}z) \leq 1 \}}.
\end{align*}
Let us prove that, for a.e. $\omega \in \partial \mathring{K}_\alpha$ and all $z \in \mathbb{R}^d$, 
\begin{align}\label{eq:pointwise_limit_general_case}
g_\alpha(\omega,z) = \underset{t\rightarrow 0^+}{\lim} \frac{1}{t^{\frac{1}{\alpha}}} \int_{1}^{+\infty} \bbone_{\{ \sigma_\alpha^*(s y + t^{\frac{1}{\alpha}}z) \leq 1 \}} s^{d-1}ds = - \langle z; \nabla(\sigma_\alpha^*)(\omega) \rangle \bbone_{\langle z; \nabla(\sigma_\alpha^*)(\omega) \rangle \leq 0}. 
\end{align}
Once this limit is obtained the remainder of the proof goes as follows: under $\mu_\alpha$, $\langle z ; \nabla(\sigma_\alpha^*)(\omega) \rangle$ is a symmetric $\alpha$-stable random variable with characteristic function given, for all $t \in \bbr$, by
\begin{align*}
\int_{\mathbb{R}^d} \exp\left( i  t \langle z ; \nabla(\sigma_\alpha^*)(\omega) \rangle \right)  \mu_\alpha(dz) = \exp \left( - |t|^\alpha \sigma_\alpha \left(\nabla(\sigma_\alpha^*)(\omega)\right)^\alpha\right).
\end{align*}
But, $\sigma_\alpha (\nabla(\sigma_\alpha^*)(\omega)) = 1$, for a.e. $\omega \in \partial \mathring{K}_\alpha$, since
\begin{align*}
\langle \omega ; \nabla(\sigma_\alpha^*)(\omega) \rangle = 1.
\end{align*}
Thus, after integrating with respect to $\mu_\alpha$, 
\begin{align*}
\underset{t \rightarrow 0^+}{\lim} \frac{K_t^{\alpha}(\mathring{K}_\alpha, (\mathring{K}_\alpha)^c)}{t^{\frac{1}{\alpha}}} = -d \mathcal{L}_d(\mathring{K}_\alpha) \bbe Y_\alpha \bbone_{Y_\alpha \leq 0} = \frac{\Gamma(1- \frac{1}{\alpha})}{\pi} d \mathcal{L}_d(\mathring{K}_\alpha). 
\end{align*}
Then, it remains to prove the second equality in \eqref{eq:pointwise_limit_general_case} as in the proof of Proposition \ref{prop:limit_discrete_case}. Now, let us consider the dual norm $\sigma_\alpha^*$. Since $\sigma_\alpha^*$ is a norm on $\mathbb{R}^d$, it is convex, positively $1$-homogeneous and Lipschitz on $\mathbb{R}^d$.  Namely, for all $r>0$ and all $x \in \mathbb{R}^d$,
\begin{align*}
\sigma_\alpha^*(rx) = r \sigma_\alpha^*(x). 
\end{align*}
Then, $\nabla \left( \sigma_\alpha^* \right)$ is positively homogeneous of degree $0$ at each point of differentiability of $\sigma_\alpha^*$, i.e., for a.e. $\omega\in \partial \mathring{K}_\alpha$ and all $s \in [1,+\infty)$, 
\begin{align*}
\nabla \left( \sigma_\alpha^* \right)(s \omega) = \nabla \left( \sigma_\alpha^* \right)( \omega).
\end{align*}
As previously, let us consider the first order Taylor approximation at $0$ of $t \mapsto \sigma_\alpha^*(s \omega + t^{\frac{1}{\alpha}} z )$ defined, for a.e. $\omega \in \partial \mathring{K}_\alpha$ and all $s \in [1, +\infty)$, by
\begin{align*}
F_\alpha\left(\omega, s\right) = s + t^{\frac{1}{\alpha}} \langle z ; \nabla \left(\sigma_\alpha^*\right)(\omega)\rangle. 
\end{align*}
From the above part of the proof,
\begin{align*}
\underset{t \rightarrow 0^+ }{\lim} \frac{1}{t^{\frac{1}{\alpha}}} \int_{1}^{+\infty}  \bbone_{F_\alpha(\omega, s) \leq 1}  s^{d-1}ds = -\langle z ; \nabla \left(\sigma_\alpha^*\right)(\omega) \rangle \bbone_{\langle z ; \nabla \left(\sigma_\alpha^*\right)(\omega) \rangle \leq 0}. 
\end{align*}
So it remains to compute the limit with $\sigma_\alpha^*(s\omega + t^{\frac{1}{\alpha}} z)$ instead of $F_\alpha(\omega,s)$. Let us study the second order differentiability property of $x \mapsto \sigma_\alpha^*(x)$. By the convexity of $\sigma_\alpha^*$ and at each point of differentiability $\omega\in \partial \mathring{K}_\alpha$, for all $s \in [1,+\infty)$, all $t>0$ and all $z \in \mathbb{R}^d$, 
\begin{align*}
\sigma_\alpha^*(s\omega + t^{\frac{1}{\alpha}} z ) \geq s + t^{\frac{1}{\alpha}} \langle \nabla \left(\sigma_\alpha^*\right)(\omega) ; z \rangle = F_\alpha\left( \omega,s \right). 
\end{align*}
which implies that,
\begin{align*}
\bbone_{\sigma_\alpha^*(s\omega + t^{\frac{1}{\alpha}} z ) \leq 1} \leq \bbone_{F_\alpha\left( \omega,s \right) \leq 1}.
\end{align*} 
Now, by Aleksandrov's theorem, $\sigma_\alpha^*$ is second order differentiable in the following sense: for a.e.  $x \in \mathbb{R}^d$, there exists a symmetric positive semi-definite matrix $H_\alpha(x)$ such that 
\begin{align*}
\sigma_\alpha^*(x + w) = \sigma_\alpha^*(x) + \langle \nabla(\sigma_\alpha^*)(x) ; w \rangle + \frac{1}{2} \langle H_\alpha(x)w ; w \rangle + o \left(\|w\|^2\right), \quad \|w\| \longrightarrow 0^+. 
\end{align*} 
With the help of the positive $1$-homogeneity and the convexity of $\sigma_\alpha^*$, let us prove that, at each point $x \in \mathbb{R}^d$ where $H_\alpha(x)$ exists,  $H_\alpha(sx)$,  with $s>0$, exists and is given by 
\begin{align*}
H_{\alpha}(sx) = \frac{1}{s} H_{\alpha}(x). 
\end{align*}
But,  for all such $x \in \mathbb{R}^d$ and all $s>0$, 
\begin{align*}
&\underset{\|w\| \rightarrow 0}{\lim} \dfrac{\sigma_\alpha^*(sx + w) -\sigma_\alpha^*(sx) - \langle \nabla(\sigma_\alpha^*)(x) ; w \rangle - \frac{1}{2s} \langle H_\alpha(x)w ; w \rangle}{\|w\|^2} \\
&= \underset{\|w\| \rightarrow 0}{\lim} \dfrac{s\sigma_\alpha^*(x + \frac{w}{s}) -s \sigma_\alpha^*(x) - \langle \nabla(\sigma_\alpha^*)(x) ; w \rangle - \frac{1}{2s} \langle H_\alpha(x)w ; w \rangle}{\|w\|^2} \\
& = \frac{1}{s} \underset{\|w\| \rightarrow 0}{\lim}  \dfrac{\sigma_\alpha^*(x + \frac{w}{s}) - \sigma_\alpha^*(x) - \langle \nabla(\sigma_\alpha^*)(x) ; \frac{w}{s} \rangle - \frac{1}{2s^2} \langle H_\alpha(x)w ; w \rangle}{\|\frac{w}{s}\|^2} = 0, 
\end{align*}
which proves the claim. Then, for a.e. $\omega \in \partial \mathring{K}_\alpha$,  all $s \in [1,+\infty)$, all $z \in \mathbb{R}^d$ and all $t>0$ small enough, 
\begin{align}\label{eq:second_order_Taylor_expansion}
\sigma_\alpha^*(s\omega + t^{\frac{1}{\alpha}}z) = s +t^{\frac{1}{\alpha}} \langle \nabla(\sigma_\alpha^*)(\omega) ; z \rangle + \frac{t^{\frac{2}{\alpha}}}{2s} \langle H_\alpha(\omega)z ; z \rangle + o \left(t^{\frac{2}{\alpha}}\right). 
\end{align}
The conclusion follows.
\end{proof}
\noindent
Next,  let us discuss the second interrogation \eqref{ineq:first_conjecture}. Equivalently,  do we have
\begin{align}\label{ineq:Riesz_rearrangement_type}
\int_{\mathbb{R}^d} \int_{\mathbb{R}^d} \bbone_E(x) \bbone_E(y) p_\alpha\left( \frac{x-y}{t^{\frac{1}{\alpha}}}\right) \frac{dxdy}{t^{\frac{d}{\alpha}}} \leq \int_{\mathbb{R}^d} \int_{\mathbb{R}^d} \bbone_{\mathring{K}_\alpha}(x) \bbone_{\mathring{K}_\alpha}(y) p_\alpha\left( \frac{x-y}{t^{\frac{1}{\alpha}}}\right) \frac{dxdy}{t^{\frac{d}{\alpha}}},
\end{align}
for all Borel measurable subsets $E$ of $\mathbb{R}^d$ such that $\mathcal{L}_d(E) = \mathcal{L}_d(\mathring{K}_\alpha)$? \\
At this point, let us remark that one cannot use a Riesz rearrangement-type inequality as in the proof of Theorem \ref{thm:sharp_classical_isoperimetric_ineq} in order to answer positively the previous question.~Indeed,  \cite[Proposition $4.2$]{Van_Schaftingen_Jean_06} ensures that if a Riesz rearrangement-type inequality holds true for an anisotropic symmetrization, then the associated gauge function is necessarily Euclidean (see also \cite[Proposition $4.4$]{Van_Schaftingen_Jean_06} for a similar conclusion with a weaker-type of inequality). However, the sharp anisotropic isoperimetric inequality is re-proved using this anisotropic symmetrization (see \cite[Theorem $7.2$]{Van_Schaftingen_Jean_06}). Let us restate inequality \eqref{ineq:Riesz_rearrangement_type} in a semigroup form. For this pupose, let $\psi \in \mathcal{C}_c^{\infty}(\mathbb{R}^d, \bbr_+)$ and let $\psi^*$ be its convex symmetrization with respect to the norm $\sigma_\alpha^*$. Then, one can ask if the following inequality holds: for all $t > 0$, 
\begin{align*}
\int_{\mathbb{R}^d} \int_{\mathbb{R}^d} \psi(x) \psi(y) p_\alpha\left( \frac{x-y}{t^{\frac{1}{\alpha}}}\right) \frac{dxdy}{t^{\frac{d}{\alpha}}} \leq \int_{\mathbb{R}^d} \int_{\mathbb{R}^d} \psi^*(x) \psi^*(y) p_\alpha\left( \frac{x-y}{t^{\frac{1}{\alpha}}}\right) \frac{dxdy}{t^{\frac{d}{\alpha}}},
\end{align*}
which is equivalent to, for all $t \geq 0$, 
\begin{align}\label{ineq:comparison_concentration_frac}
\int_{\mathbb{R}^d} \left| P_t^\alpha(\psi)(x) \right|^2 dx \leq \int_{\mathbb{R}^d} \left| P^\alpha_t(\psi^*)(x)\right|^2 dx. 
\end{align}
Note that the functions $u(t,x) = P^\alpha_t(\psi)(x)$ and $u^*(t,x) = P^\alpha_t(\psi^*)(x)$ are the respective solutions to the linear fractional parabolic equations:
\begin{align*}
(*) : \left\{
    \begin{array}{ll}
        \frac{\partial}{\partial t}(u)(t,x) = \mathcal{A}_\alpha(u)(t,x), &  \mbox{for } (t,x) \in (0, +\infty) \times \mathbb{R}^d, \\
         u(0,x) = \psi(x), &  \mbox{for } x \in \mathbb{R}^d,
    \end{array}
\right.
\end{align*}
and
\begin{align*}
(**) : \left\{
    \begin{array}{ll}
        \frac{\partial}{\partial t}(u^*)(t,x) = \mathcal{A}_\alpha(u)(t,x), &  \mbox{for } (t,x) \in (0, +\infty) \times \mathbb{R}^d, \\
         u^*(0,x) = \psi^*(x), &  \mbox{for } x \in \mathbb{R}^d,
    \end{array}
\right.
\end{align*}
where $\mathcal{A}_\alpha$ is defined by \eqref{eq:Stheatgen} on smooth enough function.~Note that comparison results based on symmetrization such as \eqref{ineq:comparison_concentration_frac} have been studied in relatively recent times for the isotropic fractional Laplacian in (non)-linear elliptic and parabolic partial differential equations (see, e.g., \cite{Biaso_Volzone_JDE12,VV_14,FV_21}) but not that much seems to be known in the general non-degenerate and symmetric situation despite regularity results (\cite{ROS_16,FRRO_17}).~Such topics will be investigated elsewhere.

To conclude our discussion regarding the sharp fractional and classical isoperimetric inequalities related to the perimeters $\mathcal{P}_{\operatorname{frac}}$ and $\mathcal{P}_{\operatorname{cl}}$, let us first present a proof of an anisotropic Sobolev-type inequality based on the P\'olya-Szeg\"o principle, the sharp anisotropic Hardy inequality and Lemma \ref{lemma_technical_Lorentz_Hardy} of the Appendix.~The idea of the proof is already contained in \cite{FS_08} but for the rotationally invariant case only.  

\begin{thm}\label{thm:PL_H_implies_Sobolev}
Let $d \geq 2$ be an integer and let $p \in [1,d)$.~Let $H$ be a norm on $\mathbb{R}^d$ and let $\mathring{H}$ be the associated dual norm. Then, for all $f \in \mathcal{C}_c^{\infty}(\mathbb{R}^d)$ 
\begin{align}\label{ineq:sobolev_anisotrop_Lorentz_scale}
\left(\int_{\mathbb{R}^d} H\left(\nabla(f)(x)\right)^p dx\right)^{\frac{1}{p}} \geq  \left(\frac{d-p}{p}\right) \mathcal{L}_d(B_{\mathring{H}})^{\frac{1}{d}} \|f\|_{p^*,p},
\end{align}
where $p^*=pd/(d-p)$, $\|.\|_{p^*,p}$ is the Lorentzian quasi-norm $(p^*,p)$ and $B_{\mathring{H}}$ is the unit ball with respect to the dual norm $\mathring{H}$. 
\end{thm}

\begin{proof}
Let $f \in \mathcal{C}_c^{\infty}(\mathbb{R}^d)$. Then, by \cite[Theorem $3.1$]{Lions_97}, 
\begin{align*}
\int_{\mathbb{R}^d} H\left(\nabla(f)(x)\right)^p dx \geq \int_{\mathbb{R}^d} H\left(\nabla(f^*)(x)\right)^p dx,
\end{align*} 
where $f^*$ is the convex symmetrization of $f$ with respect to the dual norm $\mathring{H}$. Then, by \cite[Proposition $7.5$]{Van_Schaftingen_Jean_06}, 
\begin{align*}
\int_{\mathbb{R}^d} H\left(\nabla(f^*)(x)\right)^p dx \geq \left(\frac{d-p}{p}\right)^p \int_{\mathbb{R}^d} \dfrac{|f^*(x)|^p}{\mathring{H}(x)^p}dx. 
\end{align*}
Thus, using Lemma \ref{lemma_technical_Lorentz_Hardy}, 
\begin{align*}
\int_{\mathbb{R}^d} H\left(\nabla(f^*)(x)\right)^p dx \geq \left(\frac{d-p}{p}\right)^p \mathcal{L}_d(B_{\mathring{H}})^{\frac{p}{d}}   \|f^*\|^p_{p^*,p} = \left(\frac{d-p}{p}\right)^p \mathcal{L}_d(B_{\mathring{H}})^{\frac{p}{d}}   \|f\|^p_{p^*,p}.
\end{align*}
This concludes the proof of the theorem. 
\end{proof}

\begin{rem}\label{long_rem}
(i) When $p=1$ and $f = \bbone_E$ where $E$ is a Borel subset of $\mathbb{R}^d$ with finite Lebesgue measure and finite anisotropic perimeter with respect to the norm $H$, the inequality \eqref{ineq:sobolev_anisotrop_Lorentz_scale} boils down to, 
\begin{align*}
\mathcal{P}_H(E) \geq d \mathcal{L}_d(B_{\mathring{H}})^{\frac{1}{d}} \mathcal{L}_d(E)^{\frac{d-1}{d}},
\end{align*}
where,
\begin{align}\label{eq:anisotropic_perimeter_H}
\mathcal{P}_H(E) = \sup \bigg\{ \int_{E} \operatorname{div}(\phi(x))dx,\, \phi \in \mathcal{C}^1_c(\mathbb{R}^d, \mathbb{R}^d),\, \| \mathring{H}(\phi) \|_{\infty} \leq 1 \bigg\}. 
\end{align}
(ii) Let us discuss, in the non-local cases, functional inequalities such as \eqref{ineq:sobolev_anisotrop_Lorentz_scale} as well as Hardy-type inequalities. When one considers the Gagliardo semi-norm $[\cdot]_{W^{s,p}(\mathbb{R}^d)}$ defined, for all $s \in (0,1)$ and all $1 \leq p < d/s$, by
\begin{align}\label{eq:gagliardo_seminorm}
[f]_{W^{s,p}(\mathbb{R}^d)} := \left(\int_{\mathbb{R}^d} \int_{\mathbb{R}^d} \dfrac{|f(x) - f(y)|^p}{\|x-y\|^{d+sp}} dxdy\right)^{\frac{1}{p}},
\end{align}
as a measure of fractional regularity, the analogous results are proved in \cite{FS_08}. Indeed, the authors first developed the above strategy to obtain sharp fractional Sobolev-type inequality on the Lorentz scale (\cite[Theorem 4.1]{FS_08}) based on a sharp version of Hardy inequality for fractional Sobolev spaces (\cite[Theorem 1.1]{FS_08}) and on a corresponding P\'olya-Szeg\"o principle (\cite[Theorem A.1]{FS_08}). When one is rather interested in the fractional gradient operator $D^{\alpha-1}$, results have been obtained recently for the rotationally invariant case. Indeed, starting with the cases $d \geq 2$, $s \in (0,1)$ and $p \in (1,d/s)$, \cite[Theorem 1.8]{Shieh_Spector_ACV15} provides the fractional Sobolev inequalities with the fractional gradient operator (on the Lebesgue scale) whereas \cite[Theorem 1.9]{Shieh_Spector_ACV15} furnishes the corresponding fractional Hardy-type inequalities. Regarding the geometric regime, the fractional Sobolev inequality is established in \cite[Theorem A']{SSS_17} while the fractional Hardy inequality is proved in \cite[Theorem 1.2]{Shieh_Spector_ACV18} for non-negative functions and in \cite[Theorem 1.4]{Spec_JFA20} for functions of arbitrary sign.~Note, in particular, that \cite[Theorem 1.2]{Shieh_Spector_ACV18} is sharp and that the question of optimality for this fractional Hardy-type inequality seems open. Finally, still for the $L^1$ regime, the fractional Sobolev inequality on the Lorentz scale is obtained in \cite[Inequality (1.5)]{Spec_19} (see also the \cite[proof of Theorem $1.4$ pages 21 - 22, Inequality (4.1)]{Spec_JFA20}). Let $f \in \mathcal{C}^{\infty}_c(\mathbb{R}^d, \bbr_+)$ and let $f^*$ be the convex symmetrization of $f$ with respect to the Euclidean norm $\|\cdot\|$. Applying \cite[Theorem 1.2]{Shieh_Spector_ACV18} to $f^*$ with $s = \alpha -1$, $\alpha \in (1,2)$, gives 
\begin{align}\label{ineq:optimal_Hardy_inequality_frac_gradient}
\int_{\mathbb{R}^d} \|\mathcal{D}^{\alpha-1}(f^*)(x)\| dx \geq \dfrac{(d-1) \Gamma\left(\frac{\alpha-1}{2}\right) \Gamma\left(\frac{d-1}{2}\right)}{\pi^{\frac{3-\alpha}{2}} 2^{2-\alpha} \Gamma\left(\frac{d+1-\alpha}{2}\right) } \int_{\mathbb{R}^d} \dfrac{f^*(x)}{\|x\|^{\alpha-1}}dx,
\end{align}
where $\mathcal{D}^{\alpha-1} = \nabla \circ I^{2-\alpha}$ (note that $\mathcal{D}^{\alpha-1}$ and $D^{\alpha-1, \operatorname{rot}}$ are proportional to each other as easily seen from standard Fourier analysis arguments). Then, using a variant of Lemma \ref{lemma_technical_Lorentz_Hardy} of the appendix section, 
\begin{align*}
\int_{\mathbb{R}^d} \|\mathcal{D}^{\alpha-1}(f^*)(x)\| dx \geq \dfrac{(d-1) \Gamma\left(\frac{\alpha-1}{2}\right) \Gamma\left(\frac{d-1}{2}\right)}{\pi^{\frac{3-\alpha}{2}} 2^{2-\alpha} \Gamma\left(\frac{d+1-\alpha}{2}\right) } \mathcal{L}_d(B(0,1))^{\frac{\alpha-1}{d}} \| f \|_{\frac{d}{d-\alpha+1},1},
\end{align*}
where $B(0,1)$ is the Euclidean unit ball of $\mathbb{R}^d$ and where $\| \cdot \|_{\frac{d}{d-\alpha+1},1}$ is the $(d/(d-\alpha+1),1)$-Lorentzian quasi-norm. If a $L^1$-P\'olya-Szeg\"o principle holds true for $\mathcal{D}^{\alpha-1}$, then, for all $f \in \mathcal{C}^{\infty}_c(\mathbb{R}^d)$, 
\begin{align}\label{ineq:sharp_frac_L1Sobolev_Lorentz}
\int_{\mathbb{R}^d} \|\mathcal{D}^{\alpha-1}(f)(x)\| dx \geq \dfrac{(d-1) \Gamma\left(\frac{\alpha-1}{2}\right) \Gamma\left(\frac{d-1}{2}\right)}{\pi^{\frac{3-\alpha}{2}} 2^{2-\alpha} \Gamma\left(\frac{d+1-\alpha}{2}\right) } \mathcal{L}_d(B(0,1))^{\frac{\alpha-1}{d}} \| f \|_{\frac{d}{d-\alpha+1},1}. 
\end{align}
Note that, for $p =2$, the $L^2$-P\'olya-Szeg\"o principle for $\mathcal{D}^{\alpha-1}$ follows from \cite[Theorem A.1]{FS_08}.\\
(iii) For the general non-degenerate and symmetric $\alpha$-stable case, $\alpha \in (1,2)$, nothing analogous seems to be known.
\end{rem}
\noindent
Next, let us investigate the co-area formulas naturally associated with the fractional perimeter $\mathcal{P}_{\operatorname{frac}}$.~First, for $\nu_\alpha(du) = du / \| u \|^{\alpha + d}$, based on the inequality \eqref{ineq:simple_compineq} and a known co-area formula for the fractional perimeter $\mathcal{P}_{\alpha-1}$ defined in \eqref{def:frac_perimeter_frank} (see, e.g., \cite[Lemma $4.3.$]{PS_20}), we have the following proposition which is a particular case of \eqref{ineq:classical_fractional_Sobolev_embedding_theorem} for $p=1$.

\begin{prop}\label{prop:strong_sobolev_embedding}
Let $d \geq 1$ be an integer, let $\alpha \in (1,2)$, let $\nu_\alpha(du) = du/\|u\|^{\alpha+d}$ and let $q = d/(d-(\alpha-1))$. Then, for all $f \in \mathcal{C}_c^{\infty}(\mathbb{R}^d)$, 
\begin{align}
\|f\|_{L^q(\mathbb{R}^d,dx)} \leq C^1_{\alpha,d} \int_{\mathbb{R}^d} \int_{\mathbb{R}^d} \dfrac{ | f(x) - f(y) |}{\|x -y\|^{d+\alpha-1}}dxdy,
\end{align}
where $C^1_{\alpha,d}$ is given by Theorem \ref{thm:isoperimetric_type}. 
\end{prop}

\begin{proof}
Let $f \in \mathcal{C}^{\infty}_c(\mathbb{R}^d)$, let $E_\sigma = \{x \in \mathbb{R}^d,\, :\, |f(x)|> \sigma\}$, for all $\sigma>0$, and let $q = d/(d- (\alpha - 1))$.~Then, by the layer cake representation formula, 
\begin{align*}
\| f \|^q_{L^q(\mathbb{R}^d,dx)} = \int_{\mathbb{R}^d} |f(x)|^q dx = \int_0^{+\infty} \mathcal{L}_d(E_\sigma) d(\sigma^q) =  \int_{0}^{+\infty} G(\sigma)^{q} d(\sigma^q),
\end{align*}
where $G(\sigma) = \mathcal{L}_d(E_\sigma)^{(d-(\alpha-1))/d}$.~Now, the function $G$ is non-increasing on $(0, +\infty)$, so that, 
\begin{align*}
\int_{0}^{+\infty} G(\sigma)^{q} d(\sigma^q) & = q \int_0^{+\infty} \left( \sigma G(\sigma)\right)^{q-1} G(\sigma) d\sigma  \\
& \leq q \int_{0}^{+\infty} \left(\int_0^\sigma G(r) dr \right)^{q-1} G(\sigma)d\sigma \\
& \leq \left(\int_0^{+\infty} G(\sigma) d\sigma \right)^q.
\end{align*}
Thus, by the isoperimetric inequality with $\mathcal{P}_{\operatorname{frac}}$ and the inequality \eqref{ineq:simple_compineq}, 
\begin{align*}
\| f \|_{L^q(\mathbb{R}^d,dx)}  \leq \int_0^{+\infty} \mathcal{L}_d(E_\sigma)^{(d-(\alpha-1))/d} d\sigma & \leq C^1_{\alpha,d}  \int_0^{+\infty}\mathcal{P}_{\operatorname{frac}}(E_\sigma) d\sigma \\
& \leq 2 C^1_{\alpha,d} \int_0^{+\infty}\mathcal{P}_{\alpha-1}(E_\sigma) d\sigma. 
\end{align*}
Now, for all $x, y \in \mathbb{R}^d$ with $x \ne y$, 
\begin{align*}
\int_0^{+\infty} \left| \bbone_{E_\sigma}(x) - \bbone_{E_\sigma}(y) \right| d\sigma = ||f(x)| - |f(y)|| \leq |f(x) - f(y)|. 
\end{align*}
Moreover, 
\begin{align*}
\mathcal{P}_{\alpha-1}(E_\sigma) = \frac{1}{2} \int_{\mathbb{R}^d} \int_{\mathbb{R}^d} \dfrac{ | \bbone_{E_\sigma}(x) - \bbone_{E_\sigma}(y) | }{\|x -y\|^{d+\alpha-1}}dxdy.
\end{align*}
Thus, combining the two previous identities leads to, 
\begin{align*}
\| f \|_{L^q(\mathbb{R}^d,dx)}  \leq C^1_{\alpha,d} \int_{\mathbb{R}^d} \int_{\mathbb{R}^d} \dfrac{ | f(x) - f(y) |}{\|x -y\|^{d+\alpha-1}}dxdy.
\end{align*}
The proof of the proposition is complete. 
\end{proof}
\noindent
Using Definition \ref{defi:fractional_divergence}, let us prove a coarea-type inequality.

\begin{lem}\label{lem:lower_bound}
Let $d \geq 1$ be an integer, let $\alpha \in (1,2)$, let $\nu_\alpha$ be a non-degenerate symmetric $\alpha$-stable L\'evy measure on $\mathbb{R}^d$ such that its spherical part 
is dominated by a uniform measure on $\mathbb{S}^{d-1}$.~Let $f \in \mathcal{C}_c^{\infty}(\mathbb{R}^d)$ and let $E_{\sigma} = \{x \in \mathbb{R}^d,\, |f(x)| > \sigma\}$, for all $\sigma>0$.~Then, 
\begin{align}\label{eq:easy_ineq}
\int_{0}^{+\infty} \mathcal{P}_{\operatorname{frac}}\left(E_{\sigma}\right) d\sigma \geq \int_{\mathbb{R}^d} \left\| D^{\alpha-1}(f)(x)\right\|dx. 
\end{align}
\end{lem}

\begin{proof}
Let us treat the case $f \in \mathcal{C}^\infty_c(\mathbb{R}^d)$ with $f \geq 0$.~The general case follows similarly. First, since
\begin{align*}
\int_{0}^{+\infty} \mathcal{P}_{\alpha-1}\left(E_\sigma\right) d\sigma \leq \frac{1}{2} \int_{\mathbb{R}^d} \int_{\mathbb{R}^d} \dfrac{|f(x)-f(y)|}{\|x - y\|^{d+\alpha-1}} dxdy <+\infty,   
\end{align*}
then, $\mathcal{P}_{\alpha-1}\left(E_\sigma\right)<+\infty$, for a.e.~$\sigma\in (0,+\infty)$.~Take such a $\sigma \in (0,+\infty)$.~Then, by Remark \ref{rem:Comi_et_al_Frank_et_al}, (iii),
\begin{align*}
\mathcal{P}_{\operatorname{frac}}(E_\sigma) \leq c \mathcal{P}_{\alpha-1}\left(E_\sigma\right) <+\infty, 
\end{align*}  
for some constant $c>0$ which might depend on $\alpha$ and on $d$.~Moreover, 
\begin{align*}
\mathcal{P}_{\operatorname{frac}}(E_\sigma) = \| D^{\alpha-1} \left( \bbone_{E_\sigma}\right) \|_{L^1(\mathbb{R}^d,dx)}. 
\end{align*}
Now, take $g \in \mathcal{C}_c^1(\mathbb{R}^d , \mathbb{R}^d)$ such that $\|g\|_{\infty} \leq 1$.~Then, 
\begin{align*}
\| D^{\alpha-1} \left( \bbone_{E_\sigma}\right) \|_{L^1(\mathbb{R}^d,dx)} \geq \int_{E_\sigma} \operatorname{div}_\alpha(g)(x) dx.
\end{align*}
Integrating with respect to $\sigma$ on $(0,+\infty)$ and using Fubini's theorem,
\begin{align*}
\int_0^{+\infty} \| D^{\alpha-1} \left( \bbone_{E_\sigma}\right) \|_{L^1(\mathbb{R}^d,dx)} d\sigma \geq \int_0^{+\infty} \left(\int_{E_\sigma} \operatorname{div}_\alpha(g)(x) dx\right) d\sigma = \int_{\mathbb{R}^d} f(x) \operatorname{div}_\alpha(g)(x) dx.
\end{align*} 
Thus, by duality and taking the supremum over all $g \in \mathcal{C}^1_c(\mathbb{R}^d, \mathbb{R}^d)$ with $\|g\|_\infty \leq 1$, 
\begin{align*}
\int_0^{+\infty} \| D^{\alpha-1} \left( \bbone_{E_\sigma}\right) \|_{L^1(\mathbb{R}^d,dx)} d\sigma \geq \int_{\mathbb{R}^d} \|D^{\alpha-1}(f)(x)\| dx. 
\end{align*}
This concludes the proof of the lemma.
\end{proof}
\noindent

\begin{rem}\label{remark:converse_trace_inequality}
(i) To finish, let us discuss a converse to \eqref{eq:easy_ineq}.~Namely, for all $f \in \mathcal{C}_c^{\infty}(\mathbb{R}^d)$, 
\begin{align}\label{ineq:converse_trace_inequality}
\int_{0}^{+\infty} \mathcal{P}_{\operatorname{frac}}\left(E_\sigma\right) d\sigma \leq C \int_{\mathbb{R}^d} \| D^{\alpha-1}(f)(x)\| dx,
\end{align} 
for some $C>0$.~This kind of inequality is called a $L^1$-trace inequality in the harmonic analysis community (whereas it is called a co-area inequality in the PDE/probability community, see, e.g., \cite{Bob_Houdr_97}).~In \cite[Theorem 1.3.]{Spec_19}, a negative result regarding this type of inequality for the Riesz fractional gradient operator has been presented and it is recalled.~More precisely, let $d \geq 2$, let $\alpha \in (1,2)$ and let $\mathcal{D}^{\alpha-1}$ be the fractional gradient operator defined, for all $f \in \mathcal{C}^\infty_c(\mathbb{R}^d)$, by
\begin{align}\label{ineq:Spector_fractional_Gradient}
\mathcal{D}^{\alpha-1} (f) = \nabla I^{1-(\alpha-1)} (f),
\end{align}
where $I^{2-\alpha}$ is the classical Riesz potential operator of order $2- \alpha \in (0,1)$.~Now, let $\mathcal{H}_\infty^{d-(\alpha-1)}$ be the Hausdorff content of dimension $d-\alpha+1$ defined, for all $A \subset \mathbb{R}^d$, by
\begin{align}\label{eq:hausdorff_content_dimension_d_alpha}
\mathcal{H}_\infty^{d-(\alpha-1)} \left( A \right) : = \inf \left\{ \sum_{k = 0}^{+\infty} \omega_{d-(\alpha-1)} r_k^{d-\alpha+1},\, : A \subset \cup_{k = 0}^{+\infty} B(x_k , r_k) \right\},
\end{align}
where $\omega_{d-\alpha+1} = \pi^{(d-\alpha+1)/2} / \Gamma((d-\alpha)/2+3/2)$.~Then, there is no universal constant, $C_{\alpha,d}>0$, such that
\begin{align*}
\int_0^{+\infty} \mathcal{H}_\infty^{d-(\alpha-1)} \left( E_\sigma \right) d\sigma \leq C_{\alpha,d} \int_{\mathbb{R}^d} \| \mathcal{D}^{\alpha-1} (f)(x)\| dx, 
\end{align*} 
for all $f \in L^q(\mathbb{R}^d,dx) \cap \mathcal{C}(\mathbb{R}^d)$ for some $q \in [1, d/(2-\alpha))$ such that $\mathcal{D}^{\alpha-1} (f) \in L^1(\mathbb{R}^d,dx)$~(with $E_\sigma = \{x \in \mathbb{R}^d,\, |f(x)| > \sigma \}$, for $\sigma>0$).~Moreover, this negative result when combined with \cite[Theorem 1.2.]{PS_20} implies that one cannot expect a trace inequality of the type: for all $f \in L^q(\mathbb{R}^d,dx) \cap \mathcal{C}(\mathbb{R}^d)$ such that $\mathcal{D}^{\alpha-1} (f) \in L^1(\mathbb{R}^d,dx)$~ 
\begin{align*}
\int_0^{+\infty} \mathcal{P}_{\alpha-1} \left( E_\sigma \right) d\sigma \leq C_{\alpha,d} \int_{\mathbb{R}^d} \| \mathcal{D}^{\alpha-1} (f)(x)\| dx, 
\end{align*}
for some $C_{\alpha,d}>0$, some $q \in [1, d/(2-\alpha))$ and where $\mathcal{P}_{\alpha-1}$ is the fractional perimeter defined by \eqref{ineq:perim_fractional_iso_frank}.~Indeed, according to \cite[Theorem 1.2.]{PS_20}, for all $E \in \mathcal{B}(\mathbb{R}^d)$ with $\mathcal{L}_d(E)<+\infty$ such that $\mathcal{P}_{\alpha-1}(E)<+\infty$, 
\begin{align}\label{ineq:lower_bound_inequality}
\mathcal{H}_{\infty}^{d-(\alpha-1)}(E) \leq C (\alpha-1)(2- \alpha) \mathcal{P}_{\alpha-1}(E), 
\end{align}
for some constant $C>0$.\\
(ii) In \cite[Corollary $5.6$]{comi_stefani}, it is proved that there exists function in the fractional bounded variation space such that the left-hand side
 of inequality \eqref{ineq:converse_trace_inequality} (when $\nu_\alpha = \nu_\alpha^{\operatorname{rot}}$) is equal to $+\infty$ which implies that this inequality does not hold for all functions in this space.~It is noted in \cite[Remark $5.7$]{comi_stefani}
 that the authors do not know if there can be equality in \eqref{eq:easy_ineq} for some functions in the fractional bounded variation space.
\end{rem}

\section{Optimal fractional Sobolev embeddings}\label{sec:opt_FSI}
\noindent
This section studies the optimal form of the fractional Sobolev inequality \eqref{ineq:FSI_NDS_full} as well as the existence of the non-trivial optimizers. In order to do so, let us properly define the minimization problem under consideration.~As before, let $\alpha \in (1,2)$, let $d \geq 2$ be an integer and let $2_\alpha^*$ be given by
\begin{align}\label{eq:critical_Sobolev_exponent}
2_\alpha^* : = \dfrac{2d}{d-2(\alpha-1)}.  
\end{align}
Again, let $\nu_\alpha$ be an $\alpha$-stable symmetric L\'evy measure which is non-degenerate in the sense of \eqref{eq:non_deg} and let $D^{\alpha-1}$ be the associated fractional gradient operator as defined in \eqref{eq:perim_fracGrad}. Let $\dot{W}^{\alpha-1,2}(\mathbb{R}^d,dx)$ be the set of (equivalence classes of) functions defined in \eqref{def:Lp-version_homo_frac_Sob_space} with $p=2$. The minimization problem we are interested in is the following:
\begin{align}\label{eq:minimization_problem_frac_sobolev-p2}
S_{2,\alpha,d} = \underset{f \in \dot{W}^{\alpha-1,2}(\mathbb{R}^d,dx)}{\inf} \{\|D^{\alpha-1}(f)\|^2_{L^2(\mathbb{R}^d,dx)} : \quad \|f\|_{L^{2^*_\alpha}\left(\mathbb{R}^d,dx\right)} = 1\}. 
\end{align}  
In the first part of this section, we study the existence of non-trivial minimizers following the approach presented in \cite[Section $1.2$]{Frank_24}. This approach shares similarities with the concentration-compactness principle developed in \cite{Lions_RMI185,Lions_RMI285} and who  proved the existence of non-trivial optimizers for the fractional Sobolev inequality expressed with the Gagliardo-Slobodeckij (semi)-norm (recall \eqref{ineq:classical_fractional_Sobolev_embedding_theorem}).

\begin{lem}\label{lem:there_is_something_somewhere_non-trivial}
Let $\alpha \in (1,2)$ and let $d \geq 2$ be an integer. Let $(f_n)_{n \geq 1}$ be a sequence of functions in $\dot{W}^{\alpha-1,2}\left(\mathbb{R}^d,dx\right)$ such that, 
\begin{align*}
\|f_n\|_{L^{2^*_\alpha}\left(\mathbb{R}^d,dx\right)} = 1, \quad n \geq 1,
\end{align*}  
and, 
\begin{align*}
\underset{n \rightarrow +\infty}{\lim} \|D^{\alpha-1}(f_n)\|_{L^2(\mathbb{R}^d,dx)} \in (0,+\infty).
\end{align*}
Then, up to symmetries, there exist a subsequence of $(f_n)_{n\geq 1}$ and a function $f \in \dot{W}^{\alpha-1,2}\left(\mathbb{R}^d,dx\right) \setminus\{0\}$ such that
\begin{align}\label{eq:weak_convergence_something_non-trivial}
f_{n_k} \longrightarrow f , \quad k \longrightarrow +\infty,
\end{align}
weakly in $\dot{W}^{\alpha-1,2}\left(\mathbb{R}^d,dx\right)$.
\end{lem}

\begin{proof}
Let $(f_n)_{n \geq 1}$ be a sequence of functions in $\dot{W}^{\alpha-1,2}\left(\mathbb{R}^d,dx\right)$ as in the statement. Thanks to Lemma \ref{lem:extension_RFSI}, for all $n \geq 1$, 
\begin{align*}
1 \leq C_2(\alpha,d,\chi) \left(\underset{t>0}{\sup}\, t^{\frac{d-2(\alpha-1)}{4}} \|\chi \left(-t\Delta\right)(f_n)\|_{\infty} \right)^{1-2/2^*_\alpha} \|D^{\alpha-1}(f_n)\|^{2/2^*_\alpha}_{L^2(\mathbb{R}^d,dx)},
\end{align*}
for some $C_2(\alpha,d,\chi)>0$ depending on $\alpha$, $d$ and $\chi$. Then, 
\begin{align*}
\underset{n \rightarrow +\infty}{\limsup}\, \underset{t>0}{\sup}\, t^{\frac{d-2(\alpha-1)}{4}} \|\chi \left(-t\Delta\right)(f_n)\|_{\infty}>0. 
\end{align*}
Next, for all $n \geq 1$ a fixed integer, let $(t_n, x_n)\in \mathbb{R}_+^* \times \mathbb{R}^d$ be such that 
\begin{align*}
t_n^{\frac{d-2(\alpha-1)}{4}} |\chi \left(-t_n\Delta\right)(f_n)(x_n)|\geq \frac{1}{3}\underset{t>0,x\in \mathbb{R}^d}{\sup}\, t^{\frac{d-2(\alpha-1)}{4}} |\chi \left(-t\Delta\right)(f_n)(x)|.
\end{align*}
Then, 
\begin{align*}
\underset{n \rightarrow +\infty}{\limsup}\, t_n^{\frac{d-2(\alpha-1)}{4}} |\chi \left(-t_n\Delta\right)(f_n)(x_n)|>0.
\end{align*}
Now, 
\begin{align*}
t_n^{\frac{d-2(\alpha-1)}{4}}\chi \left(-t_n\Delta\right)(f_n)(x_n) &=t_n^{\frac{d-2(\alpha-1)}{4}} \int_{\mathbb{R}^d} f_n(y) g_1\left(\frac{x_n-y}{\sqrt{t_n}}\right) \frac{dy}{t_n^{\frac{d}{2}}} \\
& = t_n^{\frac{d-2(\alpha-1)}{4}} \int_{\mathbb{R}^d} f_n(-y\sqrt{t_n}+x_n) g_1\left(y\right)dy \\
& = \int_{\mathbb{R}^d} \tilde{f}_n(y)g_1(y)dy,
\end{align*}
where, for all $y \in \mathbb{R}^d$,  
\begin{align}\label{eq:translated_dilated_sequence}
\tilde{f}_n(y) = t_n^{\frac{d-2(\alpha-1)}{4}} f_n\left(x_n-y\sqrt{t_n}\right).
\end{align}
So,
\begin{align*}
\|\tilde{f}_n\|_{L^{2_\alpha^*}\left(\mathbb{R}^d,dx\right)} = 1
\end{align*}
and 
\begin{align*}
\|D^{\alpha-1}\left(\tilde{f}_n\right)\|_{L^2(\mathbb{R}^d,dx)} = \|D^{\alpha-1}\left(f_n\right)\|_{L^2(\mathbb{R}^d,dx)} \Rightarrow \lim_{n \rightarrow +\infty}\|D^{\alpha-1}\left(\tilde{f}_n\right)\|_{L^2(\mathbb{R}^d,dx)} \in (0,+\infty). 
\end{align*}
By \cite[Theorem $3.18$]{B_FA11}, there exist a subsequence $(\tilde{f}_{n_k})_{k \geq 1}$ and $f \in L^{2_\alpha^*}(\mathbb{R}^d,dx)$ such that 
\begin{align*}
\tilde{f}_{n_k} \longrightarrow f , \quad k \longrightarrow +\infty, 
\end{align*}
weakly in $L^{2_\alpha^*}\left(\mathbb{R}^d,dx\right)$. Using \cite[Theorem $3.18$]{B_FA11} again, and by extracting a further subsequence, there exists $g_\alpha \in L^2(\mathbb{R}^d,\mathbb{R}^d,dx)$ such that
\begin{align*}
D^{\alpha-1}\left(\tilde{f}_{n_k}\right) \longrightarrow g_\alpha , \quad k \longrightarrow +\infty, 
\end{align*}
weakly in $L^{2}\left(\mathbb{R}^d,\mathbb{R}^d,dx\right)$. Next, for all $\varphi \in \mathcal{S}(\mathbb{R}^d)$ and all $\ell \in \{1, \dots, d\}$, 
\begin{align*}
\langle g_{\alpha,\ell} ; \varphi \rangle & = \underset{k \rightarrow +\infty}{\lim} \langle D^{\alpha-1}_\ell(\tilde{f}_{n_k}) ; \varphi \rangle \\
& = \underset{k \rightarrow +\infty}{\lim} \langle \tilde{f}_{n_k} ; (D^{\alpha-1}_\ell)^*(\varphi)\rangle \\
& = \langle f ; (D^{\alpha-1}_\ell)^*(\varphi)\rangle \\
& = \langle D^{\alpha-1}_\ell(f) ;\varphi\rangle.
\end{align*}
Thus, $D^{\alpha-1}(f) = g_\alpha$ and, so, $f \in \dot{W}^{\alpha-1,2}\left(\mathbb{R}^d,dx\right)$. By weak convergence, 
\begin{align*}
\left| \int_{\mathbb{R}^d} f(y)g_1(y) dy \right| = \lim_{k \rightarrow +\infty}\left| \int_{\mathbb{R}^d} \tilde{f}_{n_k}(y)g_1(y) dy \right| = \lim_{k \rightarrow +\infty} t_{n_k}^{\frac{d-2(\alpha-1)}{4}}\chi \left(-t_{n_k}\Delta\right)(f_{n_k})(x_{n_k}) >0.
\end{align*}
This concludes the proof of the lemma.
\end{proof}
\noindent
The next result is a fractional counterpart of \cite[Lemma $4.5$]{FL_AM12} and \cite[Theorems $8.6$ and $8.7$]{Lieb_Loss_book01} based on an extension of the pseudo-Poincar\'e inequality of Proposition \ref{prop:pseudo_poincare}.

\begin{lem}\label{lem:almost_everywhere_convergence}
Let $\alpha \in (1,2)$, let $d\geq 2$ be an integer and let $2_\alpha^*$ be given by \eqref{eq:critical_Sobolev_exponent}.~Let $\mu_\alpha$ be a symmetric and non-degenerate $\alpha$-stable probability measure on $\mathbb{R}^d$ with positive Lebesgue density $p_\alpha$ satisfying \eqref{ineq:condint_logarithmic_derivative} with $p=2$. Let $f \in \dot{W}^{\alpha-1,2}\left(\mathbb{R}^d,dx\right)$ and let $(f_n)_{n \geq 1}$ be a sequence of functions in $\dot{W}^{\alpha-1,2}\left(\mathbb{R}^d,dx\right)$ such that
\begin{align*}
\underset{n \geq 1}{\sup} \|D^{\alpha-1}(f_n)\|_{L^2(\mathbb{R}^d,dx)}<+\infty, 
\end{align*}
and
\begin{align*}
f_n \longrightarrow f, \quad n \longrightarrow +\infty, 
\end{align*}
weakly in $L^{2_\alpha^*}\left(\mathbb{R}^d,dx\right)$. Then, for all $K$ compact subset of $\mathbb{R}^d$, 
\begin{align}\label{eq:local_strong_convergence_L-2}
\chi_K f_n \longrightarrow \chi_K f, \quad n \longrightarrow +\infty,
\end{align}  
strongly in $L^2(\mathbb{R}^d,dx)$.~Finally, there exists a subsequence of $(f_n)_{n \geq 1}$ which converges almost everywhere to $f$. 
\end{lem}

\begin{proof}
From Proposition \ref{prop:pseudo_poincare}, for all $f \in \mathcal{C}^{\infty}_c(\mathbb{R}^d)$ and all $t>0$, 
\begin{align*}
\left\| P_t^\alpha(f) - f \right\|_{L^2(\mathbb{R}^d,dx)} \leq \frac{t^{1 - \frac{1}{\alpha}}}{\alpha-1} \left\| \frac{\nabla(p_\alpha)}{p_\alpha} \right\|_{L^2(\mu_\alpha)}  \left\| (D^{\alpha-1})(f) \right\|_{L^2(\mathbb{R}^d,dx)}. 
\end{align*}
By Theorem \ref{thm:approx_scs_homo_Lp} and Fatou's lemma, for all $f \in \dot{W}^{\alpha-1,2}\left(\mathbb{R}^d,dx\right)$ and all $t>0$, 
\begin{align*}
\left\| P_t^\alpha(f) - f \right\|_{L^2(\mathbb{R}^d,dx)} \leq \frac{t^{1 - \frac{1}{\alpha}}}{\alpha-1} \left\| \frac{\nabla(p_\alpha)}{p_\alpha} \right\|_{L^2(\mu_\alpha)}  \left\| (D^{\alpha-1})(f) \right\|_{L^2(\mathbb{R}^d,dx)}. 
\end{align*}
Now, let $(f_n)_{n\geq 1}$ and $f$ be as in the statement of the lemma and let $K$ be a compact subset of $\mathbb{R}^d$. Then, by the triangle inequality, for all $n \geq 1$ and all $t>0$, 
\begin{align*}
\| \chi_K f_n - \chi_K f\|_{L^2(\mathbb{R}^d,dx)} & \leq \|\chi_K \left(P_t^\alpha(f_n) - f_n\right)\|_{L^2(\mathbb{R}^d,dx)} + \|\chi_K \left(P_t^\alpha(f) - f\right)\|_{L^2(\mathbb{R}^d,dx)} \\
& \quad\quad + \|\chi_K \left(P_t^\alpha(f) - P_t^{\alpha}(f_n)\right)\|_{L^2(\mathbb{R}^d,dx)} \\
& \leq \frac{t^{1 - \frac{1}{\alpha}}}{\alpha-1} \left\| \frac{\nabla(p_\alpha)}{p_\alpha} \right\|_{L^2(\mu_\alpha)}  \left\| (D^{\alpha-1})(f) \right\|_{L^2(\mathbb{R}^d,dx)} \\
& \quad\quad + \frac{t^{1 - \frac{1}{\alpha}}}{\alpha-1} \left\| \frac{\nabla(p_\alpha)}{p_\alpha} \right\|_{L^2(\mu_\alpha)}  \left\| (D^{\alpha-1})(f_n) \right\|_{L^2(\mathbb{R}^d,dx)} \\
& \quad\quad + \|\chi_K \left(P_t^\alpha(f) - P_t^{\alpha}(f_n)\right)\|_{L^2(\mathbb{R}^d,dx)} \\
& \leq C(\alpha,d) \left( \underset{n \geq 1}{\sup} \|D^{\alpha-1}(f_n)\|_{L^2(\mathbb{R}^d,dx)} + \|D^{\alpha-1}(f)\|_{L^2(\mathbb{R}^d,dx)} \right) t^{1- \frac{1}{\alpha}}  \\
& \quad\quad + \|\chi_K \left(P_t^\alpha(f) - P_t^{\alpha}(f_n)\right)\|_{L^2(\mathbb{R}^d,dx)},
\end{align*}
for some constant $C(\alpha,d)>0$ depending on $\alpha$ and $d$. To conclude, it remains to prove that, for all $t>0$, 
\begin{align}\label{eq:strong_local_convergence_reg}
\underset{n \rightarrow +\infty}{\lim} \| \chi_K \left(P^\alpha_t(f)-P^{\alpha}_t(f_n)\right)\|_{L^2(\mathbb{R}^d,dx)} = 0. 
\end{align}
By H\"older's inequality and the fractional Sobolev inequality, for all $n \geq 1$, all $t>0$ and all $x \in \mathbb{R}^d$, 
\begin{align*}
\chi_K(x) |P_t^{\alpha}(f_n)(x)| &\leq \dfrac{\|f_n\|_{L^{2^*_\alpha}(\mathbb{R}^d,dx)}}{t^{\frac{d}{\alpha}}} \left\|p_\alpha\left(\frac{.}{t^{\frac{1}{\alpha}}}\right)\right\|_{L^{(2_\alpha^*)'}(\mathbb{R}^d,dx)}\chi_K(x) \\
&\leq C_2(\alpha,d) \dfrac{\sup_{n \geq 1}\|D^{\alpha-1}(f_n)\|_{L^{2}(\mathbb{R}^d,dx)}}{t^{\frac{d}{\alpha}}}\left\|p_\alpha\left(\frac{.}{t^{\frac{1}{\alpha}}}\right)\right\|_{L^{(2_\alpha^*)'}(\mathbb{R}^d,dx)}\chi_K(x),
\end{align*}
for some $C_2(\alpha,d)>0$.~The right-hand side of the previous inequality belongs to $L^2(\mathbb{R}^d,dx)$ since $K$ is a compact subset of $\mathbb{R}^d$. Moreover, for all $t>0$ and all $x \in \mathbb{R}^d$, 
\begin{align*}
P^{\alpha}_t(f_n)(x) = \int_{\mathbb{R}^d} f_n(y) p_\alpha \left(\dfrac{x-y}{t^{\frac{1}{\alpha}}}\right) \frac{dy}{t^{\frac{d}{\alpha}}} \longrightarrow \int_{\mathbb{R}^d} f(y) p_\alpha \left(\dfrac{x-y}{t^{\frac{1}{\alpha}}}\right) \frac{dy}{t^{\frac{d}{\alpha}}} = P_t^{\alpha}(f)(x),
\end{align*}
as $n$ tends to $+\infty$, since $(f_n)_{n \geq 1}$ converges weakly to $f$ in $L^{2_\alpha^*}\left(\mathbb{R}^d,dx\right)$.~An application of the Lebesgue dominated convergence theorem ensures \eqref{eq:strong_local_convergence_reg}. A standard diagonal argument concludes the proof of the lemma.  
\end{proof}
\noindent
With the help of the previous results, we are now in position to prove Theorem \ref{thm:existence_non-trivial_optimizer}. The proof is an extension to this fractional setting of the missing mass method developed in \cite{Lieb_AM83}. \\

\noindent
\textit{Proof of Theorem \ref{thm:existence_non-trivial_optimizer}.}
Let $(f_n)_{n \geq 1}$ be an optimizing sequence of the minimization problem \eqref{eq:minimization_problem_frac_sobolev-p2}, i.e., for all $n \geq 1$ integer, 
\begin{align*}
f_n \in \dot{W}^{\alpha-1,2}\left(\mathbb{R}^d,dx\right), \quad \|f_n\|_{L^{2^*_\alpha}\left(\mathbb{R}^d,dx\right)} = 1 
\end{align*}
and
\begin{align*}
\underset{n\longrightarrow+\infty}{\lim} \|D^{\alpha-1}(f_n)\|^2_{L^2(\mathbb{R}^d,dx)} = S_{2,\alpha,d}. 
\end{align*}
Using Lemma \ref{lem:there_is_something_somewhere_non-trivial}, up to symmetries, there exists a subsequence $(f_{n_k})_{k \geq 1}$ which converges weakly to $f \in \dot{W}^{\alpha-1,2}(\mathbb{R}^d,dx) \setminus \{0\}$ in $L^{2^*_\alpha}\left(\mathbb{R}^d,dx\right)$ and for which $(D^{\alpha-1}(f_{n_k}))_{k \geq 1}$ converges weakly to $D^{\alpha-1}(f)$ in $L^2(\mathbb{R}^d, \mathbb{R}^d,dx)$. Moreover, from Lemma \ref{lem:almost_everywhere_convergence} and, without loss of generality, $(f_{n_k})_{k \geq 1}$ converges to $f$ almost everywhere. Now, let $(r_{n_k})_{k \geq 1}$ be defined, for all $k \geq 1$ integer, by $r_{n_k} = f_{n_k}-f$. Then,
\begin{align*}
\|D^{\alpha-1}(r_{n_k})\|^2_{L^2(\mathbb{R}^d,dx)} \longrightarrow S_{2, \alpha, d} - \|D^{\alpha-1}(f)\|^2_{L^2(\mathbb{R}^d,dx)}, 
\end{align*}
as $k$ tends to $+\infty$. Since $\|f_{n_k}\|_{L^{2_\alpha^*}\left(\mathbb{R}^d,dx\right)} = 1$, for all $k \geq 1$,
\begin{align*}
\sup_{k \geq 1} \|f_{n_k}\|_{L^{2_\alpha^*}\left(\mathbb{R}^d,dx\right)} < + \infty. 
\end{align*}
By the Br\'ezis-Lieb lemma (see, e.g., \cite[Theorem 1.9.]{Lieb_Loss_book01}), 
\begin{align*}
\underset{k \rightarrow +\infty}{\lim} \int_{\mathbb{R}^d} \left| \left| f_{n_k}(x) \right|^{2^*_\alpha} -\left| f_{n_k}(x) - f(x)  \right|^{2^*_\alpha} -\left| f(x) \right|^{2^*_\alpha} \right| dx = 0. 
\end{align*} 
Thus, 
\begin{align*}
\underset{k \rightarrow +\infty}{\lim} \|r_{n_k}\|^{2_\alpha^*}_{L^{2^*_\alpha}(\mathbb{R}^d,dx)} = 1 - \|f\|^{2_\alpha^*}_{L^{2^*_\alpha}(\mathbb{R}^d,dx)}. 
\end{align*}
By the fractional Sobolev inequality, for all $k \geq 1$, 
\begin{align*}
S_{2, \alpha, d} \|r_{n_k}\|^2_{L^{2_\alpha^*}\left(\mathbb{R}^d,dx\right)} \leq \|D^{\alpha-1}(r_{n_k})\|^2_{L^2(\mathbb{R}^d,dx)}. 
\end{align*}
Passing to the limit in the previous inequality leads to, 
\begin{align*}
S_{2,\alpha,d} \left(1-\|f\|^{2_\alpha^*}_{L^{2_\alpha^*}(\mathbb{R}^d,dx)}\right)^{\frac{2}{2^*_{\alpha}}} \leq S_{2,\alpha,d}-\|D^{\alpha-1}(f)\|^2_{L^2(\mathbb{R}^d,dx)}.
\end{align*}
Moreover, since $2^*_\alpha>2$, for all $a,b \geq 0$, 
\begin{align*}
\left(a+b\right)^{\frac{2}{2_\alpha^*}} \leq a^{\frac{2}{2_\alpha^*}} + b^{\frac{2}{2_\alpha^*}}. 
\end{align*}
Thus, 
\begin{align*}
S_{2,\alpha,d}\left(1 - \|f\|^2_{L^{2_\alpha^*}\left(\mathbb{R}^d,dx\right)}\right) \leq S_{2,\alpha,d} - \|D^{\alpha-1}(f)\|^2_{L^2(\mathbb{R}^d,dx)}, 
\end{align*}
which implies that $S_{2,\alpha,d} \|f\|^2_{L^{2_\alpha^*}\left(\mathbb{R}^d,dx\right)} = \|D^{\alpha-1}(f)\|^2_{L^2(\mathbb{R}^d,dx)}$ since $f \ne 0$. This concludes the proof of the theorem.$\qed$\\

\noindent
Next, let us perform a similar analysis for $p \in (1, d/(\alpha-1))$ when the reference $\alpha$-stable L\'evy measure $\nu_\alpha$ given by \eqref{eq:polar} is non-degenerate symmetric and is a $\gamma$-measure with $\gamma \in [1,d]$ such that $\gamma - d + 2\alpha > 1$. Let $p_\alpha^*$ be given by 
\begin{align}\label{eq:critical_Sobolev_exponent_p}
p_\alpha^* = \dfrac{pd}{d-p(\alpha-1)}. 
\end{align}
Let $S_{p,\alpha,d,\gamma}$ be the non-negative constant defined by 
\begin{align}\label{eq:minimization_problem_frac_sobolev-p}
S_{p,\alpha,d,\gamma} = \underset{f \in \dot{W}^{\alpha-1,p}(\mathbb{R}^d,dx)}{\inf} \{\|D^{\alpha-1}(f)\|^p_{L^p(\mathbb{R}^d,dx)} : \quad \|f\|_{L^{p^*_\alpha}\left(\mathbb{R}^d,dx\right)} = 1\}. 
\end{align}

\begin{lem}\label{lem:there_is_something_somewhere_non-trivial_p}
Let $\alpha \in (1,2)$, let $d \geq 2$ be an integer and let $p \in (1,d/(\alpha-1))$.~Let $\mu_\alpha$ be a non-degenerate symmetric $\alpha$-stable probability measure on $\mathbb{R}^d$ with L\'evy measure $\nu_\alpha$ such that $\nu_\alpha$ is a $\gamma$-measure with $\gamma \in [1,d]$ and $\gamma - d + 2\alpha > 1$. Let $(f_n)_{n \geq 1}$ be a sequence of functions in $\dot{W}^{\alpha-1,p}\left(\mathbb{R}^d,dx\right)$ such that, 
\begin{align*}
\|f_n\|_{L^{p^*_\alpha}\left(\mathbb{R}^d,dx\right)} = 1, \quad n \geq 1
\end{align*}  
and
\begin{align*}
\underset{n \rightarrow +\infty}{\lim} \|D^{\alpha-1}(f_n)\|_{L^p(\mathbb{R}^d,dx)} \in (0,+\infty).
\end{align*}
Then, up to symmetries, there exist a subsequence of $(f_n)_{n\geq 1}$ and a function $f \in \dot{W}^{\alpha-1,p}\left(\mathbb{R}^d,dx\right)$ such that
\begin{align}\label{eq:weak_convergence_something_non-trivial_p}
f_{n_k} \longrightarrow f , \quad k \longrightarrow +\infty,
\end{align}
weakly in $\dot{W}^{\alpha-1,p}\left(\mathbb{R}^d,dx\right)$.
\end{lem}

\begin{proof}
Let $(f_n)_{n \geq 1}$ be a sequence of functions in $\dot{W}^{\alpha-1,p}\left(\mathbb{R}^d,dx\right)$ as in the statement of the lemma. Theorem \ref{thm:refined_frac_Sobolev_Lp} asserts that, for all $n \geq 1$, 
\begin{align*}
1 \leq \tilde{C}_{\alpha,d,p,\gamma,t} \|f_n\|^{1-\frac{t}{p_\alpha^*}}_{1, \frac{d-p(\alpha-1)}{p}}  \|D^{\alpha-1}(f_n)\|_{L^p(\mathbb{R}^d,dx)}^{\frac{t}{p_{\alpha}^*}},
\end{align*}
for some $t \in [p,p_\alpha^*)$ such that $t >p_\alpha^*-1$ and $\tilde{C}_{\alpha,d,p,\gamma,t}>0$. Thus, 
\begin{align*}
\underset{n \rightarrow +\infty}{\limsup}\, \|f_n\|_{1, \frac{d-p(\alpha-1)}{p}} >0. 
\end{align*}
Next, for any $n \geq 1$, let $(R_n, x_n)\in \mathbb{R}_+^* \times \mathbb{R}^d$ be such that 
\begin{align*}
R_n^{\frac{d-p(\alpha-1)}{p}-d} \int_{B(x_n,R_n)} |f_n(y)|dy\geq \frac{1}{3} \|f_n\|_{1, \frac{d-p(\alpha-1)}{p}}. 
\end{align*}
Then, 
\begin{align*}
\underset{n \rightarrow +\infty}{\limsup}\, R_n^{\frac{d-p(\alpha-1)}{p}-d} \int_{B(x_n,R_n)} |f_n(y)|dy >0.
\end{align*}
Next, for all $n \geq 1$, 
\begin{align*}
R_n^{\frac{d-p(\alpha-1)}{p}-d} \int_{B(x_n,R_n)} |f_n(y)|dy &= R_n^{\frac{d-p(\alpha-1)}{p}-d} \int_{B(0,R_n)} |f_n(x_n+z)| dz \\
&= R_n^{\frac{d-p(\alpha-1)}{p}} \int_{B(0,1)} |f_{n}(x_n+R_n y) | dy \\
& = \int_{B(0,1)} |\tilde{f}_n(y)| dy, 
\end{align*}
where, for all $y \in \mathbb{R}^d$,  
\begin{align}\label{eq:translated_dilated_sequence_p}
\tilde{f}_n(y) = R_n^{\frac{d-p(\alpha-1)}{p}}f_{n}(x_n+R_n y)  .
\end{align}
Now, for all $n \geq 1$, 
\begin{align*}
\|\tilde{f}_n\|_{L^{p_\alpha^*}\left(\mathbb{R}^d,dx\right)} = 1
\end{align*}
and
\begin{align*}
\|D^{\alpha-1}\left(\tilde{f}_n\right)\|_{L^p(\mathbb{R}^d,dx)} = \|D^{\alpha-1}\left(f_n\right)\|_{L^p(\mathbb{R}^d,dx)} \Rightarrow \lim_{n \rightarrow +\infty}\|D^{\alpha-1}\left(\tilde{f}_n\right)\|_{L^p(\mathbb{R}^d,dx)} \in (0,+\infty). 
\end{align*}
By \cite[Theorem $3.18$]{B_FA11}, there exist a subsequence $(\tilde{f}_{n_k})_{k \geq 1}$ and $f \in L^{p_\alpha^*}(\mathbb{R}^d,dx)$ such that 
\begin{align*}
\tilde{f}_{n_k} \longrightarrow f , \quad k \longrightarrow +\infty, 
\end{align*}
weakly in $L^{p_\alpha^*}\left(\mathbb{R}^d,dx\right)$. Thanks to \cite[Theorem $3.18$]{B_FA11} again and by extracting a further subsequence, there exists $g_\alpha \in L^p(\mathbb{R}^d,\mathbb{R}^d,dx)$ such that
\begin{align*}
D^{\alpha-1}\left(\tilde{f}_{n_k}\right) \longrightarrow g_\alpha , \quad k \longrightarrow +\infty, 
\end{align*}
weakly in $L^{p}\left(\mathbb{R}^d,\mathbb{R}^d,dx\right)$. Now, for all $\varphi \in \mathcal{S}(\mathbb{R}^d)$ and all $\ell \in \{1, \dots, d\}$, 
\begin{align*}
\langle g_{\alpha,\ell} ; \varphi \rangle & = \underset{k \rightarrow +\infty}{\lim} \langle D^{\alpha-1}_\ell(\tilde{f}_{n_k}) ; \varphi \rangle \\
& = \underset{k \rightarrow +\infty}{\lim} \langle \tilde{f}_{n_k} ; (D^{\alpha-1}_\ell)^*(\varphi)\rangle \\
& = \langle f ; (D^{\alpha-1}_\ell)^*(\varphi)\rangle \\
& = \langle D^{\alpha-1}_\ell(f) ;\varphi\rangle.
\end{align*}
Thus, $D^{\alpha-1}(f) = g_\alpha$ and, so, $f \in \dot{W}^{\alpha-1,p}\left(\mathbb{R}^d,dx\right)$.
\end{proof}
\noindent
To prove the non-triviality of the limit $f$, we need a strong convergence result similar to the one obtained in Lemma \ref{lem:almost_everywhere_convergence} in the $L^p$-setting, with $p \in (1, d/(\alpha-1))$.  

\begin{lem}\label{lem:almost_everywhere_convergence_p}
Let $\alpha \in (1,2)$, let $d\geq 2$ be an integer and let $p \in (1,d/(\alpha-1))$.~Let $\mu_\alpha$ be a symmetric and non-degenerate $\alpha$-stable probability measure on $\mathbb{R}^d$ with positive Lebesgue density $p_\alpha$ satisfying \eqref{ineq:condint_logarithmic_derivative} with $p$ and L\'evy measure $\nu_\alpha$. Let $f \in \dot{W}^{\alpha-1,p}\left(\mathbb{R}^d,dx\right)$ and let $(f_n)_{n \geq 1}$ be a sequence of functions in $\dot{W}^{\alpha-1,p}\left(\mathbb{R}^d,dx\right)$ such that
\begin{align*}
\underset{n \geq 1}{\sup} \|D^{\alpha-1}(f_n)\|_{L^p(\mathbb{R}^d,dx)}<+\infty 
\end{align*}
and
\begin{align*}
f_n \longrightarrow f, \quad n \longrightarrow +\infty, 
\end{align*}
weakly in $L^{p_\alpha^*}\left(\mathbb{R}^d,dx\right)$. Then, for all $K$ compact subsets of $\mathbb{R}^d$, 
\begin{align}\label{eq:local_strong_convergence_L-p}
\chi_K f_n \longrightarrow \chi_K f, \quad n \longrightarrow +\infty,
\end{align}  
strongly in $L^p(\mathbb{R}^d,dx)$.~Finally, there exists a subsequence of $(f_n)_{n \geq 1}$ which converges almost everywhere to $f$. 
\end{lem}

\begin{proof}
By Proposition \ref{prop:pseudo_poincare}, for all $f \in \mathcal{C}^{\infty}_c(\mathbb{R}^d)$ and all $t>0$, 
\begin{align*}
\left\| P_t^\alpha(f) - f \right\|_{L^p(\mathbb{R}^d,dx)} \leq \frac{t^{1 - \frac{1}{\alpha}}}{\alpha-1} \left\| \frac{\nabla(p_\alpha)}{p_\alpha} \right\|_{L^p(\mu_\alpha)}  \left\| (D^{\alpha-1})(f) \right\|_{L^p(\mathbb{R}^d,dx)}. 
\end{align*}
By Theorem \ref{thm:approx_scs_homo_Lp} and Fatou's lemma, for all $f \in \dot{W}^{\alpha-1,p}\left(\mathbb{R}^d,dx\right)$ and all $t>0$, 
\begin{align*}
\left\| P_t^\alpha(f) - f \right\|_{L^p(\mathbb{R}^d,dx)} \leq \frac{t^{1 - \frac{1}{\alpha}}}{\alpha-1} \left\| \frac{\nabla(p_\alpha)}{p_\alpha} \right\|_{L^p(\mu_\alpha)}  \left\| (D^{\alpha-1})(f) \right\|_{L^p(\mathbb{R}^d,dx)}. 
\end{align*}
Now, let $(f_n)_{n\geq 1}$ and $f$ be as in the statement of the lemma and let $K$ be a compact subset of $\mathbb{R}^d$. Then, by the triangle inequality, for all $n \geq 1$ an integer and all $t>0$, 
\begin{align*}
\| \chi_K f_n - \chi_K f\|_{L^p(\mathbb{R}^d,dx)} & \leq \|\chi_K \left(P_t^\alpha(f_n) - f_n\right)\|_{L^p(\mathbb{R}^d,dx)} + \|\chi_K \left(P_t^\alpha(f) - f\right)\|_{L^p(\mathbb{R}^d,dx)} \\
& \quad\quad + \|\chi_K \left(P_t^\alpha(f) - P_t^{\alpha}(f_n)\right)\|_{L^p(\mathbb{R}^d,dx)} \\
& \leq \frac{t^{1 - \frac{1}{\alpha}}}{\alpha-1} \left\| \frac{\nabla(p_\alpha)}{p_\alpha} \right\|_{L^p(\mu_\alpha)}  \left\| (D^{\alpha-1})(f) \right\|_{L^p(\mathbb{R}^d,dx)} \\
& \quad\quad + \frac{t^{1 - \frac{1}{\alpha}}}{\alpha-1} \left\| \frac{\nabla(p_\alpha)}{p_\alpha} \right\|_{L^p(\mu_\alpha)}  \left\| (D^{\alpha-1})(f_n) \right\|_{L^p(\mathbb{R}^d,dx)} \\
& \quad\quad + \|\chi_K \left(P_t^\alpha(f) - P_t^{\alpha}(f_n)\right)\|_{L^p(\mathbb{R}^d,dx)} \\
& \leq C(\alpha,d,p) \left( \underset{n \geq 1}{\sup} \|D^{\alpha-1}(f_n)\|_{L^p(\mathbb{R}^d,dx)} + \|D^{\alpha-1}(f)\|_{L^p(\mathbb{R}^d,dx)} \right) t^{1- \frac{1}{\alpha}}  \\
& \quad\quad + \|\chi_K \left(P_t^\alpha(f) - P_t^{\alpha}(f_n)\right)\|_{L^p(\mathbb{R}^d,dx)},
\end{align*}
for some constant $C(\alpha,d,p)>0$ depending on $\alpha$, $d$ and $p$. To conclude, it remains to prove that, for all $t>0$, 
\begin{align}\label{eq:strong_local_convergence_reg_p}
\underset{n \rightarrow +\infty}{\lim} \| \chi_K \left(P^\alpha_t(f)-P^{\alpha}_t(f_n)\right)\|_{L^p(\mathbb{R}^d,dx)} = 0. 
\end{align}
By H\"older's inequality and the fractional Sobolev inequality, for all $n \geq 1$, all $t>0$ and all $x \in \mathbb{R}^d$, 
\begin{align*}
\chi_K(x) |P_t^{\alpha}(f_n)(x)| &\leq \dfrac{\|f_n\|_{L^{p^*_\alpha}(\mathbb{R}^d,dx)}}{t^{\frac{d}{\alpha}}} \left\|p_\alpha\left(\frac{.}{t^{\frac{1}{\alpha}}}\right)\right\|_{L^{(p_\alpha^*)'}(\mathbb{R}^d,dx)}\chi_K(x) \\
&\leq C_2(\alpha,d,p,\gamma) \dfrac{\sup_{n \geq 1}\|D^{\alpha-1}(f_n)\|_{L^{p}(\mathbb{R}^d,dx)}}{t^{\frac{d}{\alpha}}}\left\|p_\alpha\left(\frac{.}{t^{\frac{1}{\alpha}}}\right)\right\|_{L^{(p_\alpha^*)'}(\mathbb{R}^d,dx)}\chi_K(x),
\end{align*}
for some $C_2(\alpha,d,p,\gamma)>0$.~The right-hand side of the previous inequality belongs to $L^p(\mathbb{R}^d,dx)$ since $K$ is a compact subset of $\mathbb{R}^d$. Moreover, for all $t>0$ and all $x \in \mathbb{R}^d$, 
\begin{align*}
P^{\alpha}_t(f_n)(x) = \int_{\mathbb{R}^d} f_n(y) p_\alpha \left(\dfrac{x-y}{t^{\frac{1}{\alpha}}}\right) \frac{dy}{t^{\frac{d}{\alpha}}} \longrightarrow \int_{\mathbb{R}^d} f(y) p_\alpha \left(\dfrac{x-y}{t^{\frac{1}{\alpha}}}\right) \frac{dy}{t^{\frac{d}{\alpha}}} = P_t^{\alpha}(f)(x),
\end{align*}
as $n$ tends to $+\infty$, since $(f_n)_{n \geq 1}$ converges weakly to $f$ in $L^{p_\alpha^*}\left(\mathbb{R}^d,dx\right)$.~An application of the Lebesgue dominated convergence theorem ensures \eqref{eq:strong_local_convergence_reg_p}. A standard diagonal argument concludes the proof of the lemma.  
\end{proof}
\noindent
Combining Lemma \ref{lem:there_is_something_somewhere_non-trivial_p} and Lemma \ref{lem:almost_everywhere_convergence_p}, we can conclude that the limit $f \in \dot{W}^{\alpha-1,p}\left(\mathbb{R}^d,dx\right)$ of Lemma \ref{lem:there_is_something_somewhere_non-trivial_p} is non-trivial; namely, $f \ne 0$. Indeed, by strong convergence in $L^1_{\operatorname{loc}}(\mathbb{R}^d,dx)$, 
\begin{align*}
\int_{B(0,1)} |f(y)|dy = \underset{n \rightarrow +\infty}{\lim} \int_{B(0,1)} |\tilde{f}_n(y)|dy >0, 
\end{align*}
where $\tilde{f}_n$ is given by \eqref{eq:translated_dilated_sequence_p} (up to a weakly convergent subsequence). With the help of the strict inequality above, let us prove that vanishing cannot occur in the terminology of \cite{Lions_RMI185,Lions_RMI285}.

\begin{lem}\label{lem:exclusion_vanishing}
Let $(\tilde{f}_n)_{n \geq 1}$ and $f$ be as in Lemma \ref{lem:there_is_something_somewhere_non-trivial_p} and Lemma \ref{lem:almost_everywhere_convergence_p} such that 
\begin{align}\label{eq:non-trivial_p}
\int_{B(0,1)} |f(y)|dy = \underset{n \rightarrow +\infty}{\lim} \int_{B(0,1)} |\tilde{f}_n(y)|dy >0. 
\end{align}
Then, 
\begin{align}\label{ineq:exclusion_vanishing}
\underset{n \longrightarrow +\infty}{\liminf} \underset{x \in \mathbb{R}^d}{\sup} \int_{B(x,1)} |\tilde{f}_n(y)|^{p_\alpha^*} dy>0.
\end{align}
\end{lem}
 
\begin{proof}
Using \eqref{eq:non-trivial_p}, 
\begin{align*}
\underset{n \longrightarrow +\infty}{\liminf} \underset{x \in \mathbb{R}^d}{\sup} \int_{B(x,1)} |\tilde{f}_n(y)| dy \geq \underset{n \longrightarrow +\infty}{\lim} \int_{B(0,1)} |\tilde{f}_n(y)| dy>0.
\end{align*}
Let us fix $x \in \mathbb{R}^d$. Then, by H\"older's inequality with $q = p_\alpha^*$ and $q'= q/(q-1)$, 
\begin{align*}
\int_{B(x,1)} |\tilde{f}_n(y)| dy \leq \mathcal{L}_d\left(B(0,1)\right)^{\frac{1}{q'}} \left(\int_{B(x,1)} |\tilde{f}_n(y)|^{p_\alpha^*} dy  \right)^{\frac{1}{p_\alpha^*}}.
\end{align*}
Thus, 
\begin{align*}
\underset{x \in \mathbb{R}^d}{\sup} \int_{B(x,1)} |\tilde{f}_n(y)|^{p_\alpha^*} dy \geq \dfrac{1}{\mathcal{L}_d\left(B(0,1)\right)^{\frac{p_\alpha^*}{q'}}} \left(\underset{x \in \mathbb{R}^d}{\sup} \int_{B(x,1)} |\tilde{f}_n(y)|dy\right)^{p_\alpha^*}.
\end{align*}
Passing to the limit in the previous inequality concludes the proof of the lemma. 
\end{proof}
\noindent
Let us now state and prove a concentration-compactness principle in this fractional setting. 

\begin{thm}\label{thm:ccp_homo_frac_p}
Let $\alpha \in (1,2)$, let $d \geq 2$ be an integer and let $p \in (1,d/(\alpha-1))$.~~Let $\mu_\alpha$ be a symmetric and non-degenerate $\alpha$-stable probability measure on $\mathbb{R}^d$ with positive Lebesgue density $p_\alpha$ satisfying \eqref{ineq:condint_logarithmic_derivative} with $p$ and L\'evy measure $\nu_\alpha$.~Let $f\in \dot{W}^{\alpha-1,p}\left(\mathbb{R}^d,dx\right)$ and let $(f_n)_{n\geq 1}$ be a sequence of functions in $\dot{W}^{\alpha-1,p}\left(\mathbb{R}^d,dx\right)$ such that, 
\begin{align*}
f_n \longrightarrow f, \quad n \longrightarrow +\infty,
\end{align*}
weakly in $L^{p_\alpha^*}\left(\mathbb{R}^d,dx\right)$ and
\begin{align*}
D^{\alpha-1}(f_n) \longrightarrow D^{\alpha-1}(f), \quad n \longrightarrow +\infty,
\end{align*}
weakly in $L^p\left(\mathbb{R}^d,\mathbb{R}^d,dx\right)$.~Let $\tau(dx)$ and $\rho(dx)$ be the weak-$\ast$ limits of the sequences of positive finite Borel measures $(|f_n(x)|^{p_\alpha^*} dx)_{n\geq 1}$ and $(\|D^{\alpha-1}(f_n)(x)\|^p dx)_{n \geq 1}$ respectively (up to a subsequence). Then, there exist a set $I$ at most enumerable, nonnegative reals $\rho_i, \tau_i$, $i \in I$, and distinct points $x_i$, $i \in I$, in $\mathbb{R}^d$ such that: 
\begin{enumerate}
\item $\|D^{\alpha-1}(f_n)(x)\|^p dx \longrightarrow \rho(dx) \geq \|D^{\alpha-1}(f)(x)\|^p dx + \sum_{i \in I} \rho_i \delta_{x_i}(dx)$ weakly-$\ast$ in the sense of measure;
\item $|f_n(x)|^{p_\alpha^*} dx \longrightarrow |f(x)|^{p_\alpha^*} dx + \sum_{i \in I} \tau_i \delta_{x_i}(dx)$ weakly-$\ast$ in the sense of measure;
\item for all $i \in I$, 
\begin{align}\label{ineq:rH_atomic}
S_{p,\alpha,d,\gamma} \tau_i^{\frac{p}{p_\alpha^*}} \leq \rho_i. 
\end{align}  
\end{enumerate}
Moreover, let $\tau_\infty$ and $\rho_\infty$ be defined by
\begin{align*}
&\tau_\infty : = \underset{R \rightarrow +\infty}{\lim}\, \underset{n \longrightarrow +\infty}{\limsup} \int_{\|x\| > R} |f_n(x)|^{p_\alpha^*} dx , \\
&\rho_\infty : =  \underset{R \rightarrow +\infty}{\lim}\, \underset{n \longrightarrow +\infty}{\limsup} \int_{\|x\| > R} \|D^{\alpha-1}(f_n)(x)\|^{p} dx. 
\end{align*}
Then, 
\begin{enumerate}
\item $\underset{n \rightarrow +\infty}{\limsup} \int_{\mathbb{R}^d} \|D^{\alpha-1}(f_n)(x)\|^pdx = \rho\left(\mathbb{R}^d\right) + \rho_\infty$ ; 
\item $\underset{n \rightarrow +\infty}{\limsup} \int_{\mathbb{R}^d} |f_n(x)|^{p_\alpha^*}dx = \int_{\mathbb{R}^d} |f(x)|^{p_\alpha^*}dx + \sum_{i \in I} \tau_i + \tau_\infty$.
\end{enumerate}
\end{thm}

\begin{proof}
\textbf{Step 1:} Let $(\rho_n)_{n \geq 1}$ and $(\tau_n)_{n \geq 1}$ be the sequences of finite Borel measures on $\mathbb{R}^d$ defined, for all integers $n \geq 1$, by 
\begin{align*}
\rho_n(dx) = \|D^{\alpha-1}(f_n)(x)\|^p dx, \quad \tau_n(dx) = |f_n(x)|^{p_\alpha^*}dx. 
\end{align*}
By assumptions, these sequences of finite Borel measures weakly-$\ast$ converge to $\rho$ and $\tau$ respectively. Next, let us prove the first three points in the statement of Theorem \ref{thm:ccp_homo_frac_p}. Following \cite{Lions_RMI185,BSS_NoDEA18}, let us treat the case $f=0$ first. Now, let us prove the following reverse H\"older inequality: for all $\varphi \in \mathcal{C}_c^{\infty}\left(\mathbb{R}^d\right)$, 
\begin{align}\label{ineq:reverse-Holder}
S_{p, \alpha,d,\gamma}^{\frac{1}{p}} \|\varphi\|_{L^{p_\alpha^*}(\mathbb{R}^d,\tau(dx))} \leq \|\varphi\|_{L^{p}(\mathbb{R}^d,\rho(dx))}.
\end{align}
By the fractional Sobolev inequality, for all $n \geq 1$ and all $\varphi \in \mathcal{C}_c^{\infty}\left(\mathbb{R}^d\right)$, 
\begin{align*}
S^{\frac{1}{p}}_{p,\alpha,d,\gamma} \|\varphi f_n\|_{L^{p_\alpha^*}\left(\mathbb{R}^d,dx\right)} \leq \|D^{\alpha-1}\left(\varphi f_n\right)\|_{L^p\left(\mathbb{R}^d,dx\right)}. 
\end{align*}
Now, 
\begin{align*}
\|\varphi f_n\|^{p_\alpha^*}_{L^{p_\alpha^*}\left(\mathbb{R}^d,dx\right)} = \int_{\mathbb{R}^d} |\varphi(x)|^{p_\alpha^*} |f_n(x)|^{p_\alpha^*} dx \longrightarrow \int_{\mathbb{R}^d} |\varphi(x)|^{p_\alpha^*} \tau(dx),
\end{align*}
as $n$ tends to $+\infty$. It remains to study $\underset{n \longrightarrow +\infty} \lim \|D^{\alpha-1}\left(\varphi f_n\right)\|_{L^p\left(\mathbb{R}^d,dx\right)}$.
But, for all $n \geq 1$ and a.e.-$x \in \mathbb{R}^d$, 
\begin{align*}
D^{\alpha-1}(\varphi f_n)\left(x\right) = \varphi(x) D^{\alpha-1}\left(f_n\right)(x) + f_n(x) D^{\alpha-1}\left(\varphi\right)(x) + R^{\alpha}\left(\varphi,f_n\right)(x),
\end{align*}
where
\begin{align*}
R^{\alpha}\left(\varphi,f_n\right)(x) = \int_{\mathbb{R}^d} \Delta_u(\varphi)(x) \Delta_u(f_n)(x) u \nu_\alpha(du). 
\end{align*} 
Note that, for all $n \geq 1$ and a.e.-$x \in \mathbb{R}^d$, 
\begin{align*}
f_n(x) D^{\alpha-1}\left(\varphi\right)(x) + R^{\alpha}\left(\varphi,f_n\right)(x) = \int_{\mathbb{R}^d}f_n(x+u) \Delta_u(\varphi)(x) u \nu_\alpha(du). 
\end{align*}
By Minkowski's inequality, for all $n \geq 1$, 
\begin{align*}
\left\| D^{\alpha-1}(\varphi f_n) \right\|_{L^p\left(\mathbb{R}^d,dx\right)} \leq \|\varphi D^{\alpha-1}\left(f_n\right)\|_{L^p\left(\mathbb{R}^d,dx\right)}+ \|f_n D^{\alpha-1}\left(\varphi\right) + R^{\alpha}\left(\varphi,f_n\right)\|_{L^p\left(\mathbb{R}^d,dx\right)}.
\end{align*}
Now, for all $n \geq 1$, 
\begin{align*}
\left\|f_n D^{\alpha-1}\left(\varphi\right) + R^{\alpha}\left(\varphi,f_n\right) \right\|_{L^p\left(\mathbb{R}^d,dx\right)} & = \left(\int_{\mathbb{R}^d} \left\| \int_{\mathbb{R}^d} f_n(x+u)\left(\varphi(x+u) - \varphi(x)\right) u\nu_\alpha(du)\right\|^p dx \right)^{\frac{1}{p}} \\
& \leq \int_{\mathbb{R}^d} \|u\|\nu_\alpha(du) \|f_n(\cdot + u)(\varphi(\cdot+u)-\varphi(\cdot))\|_{L^p(\mathbb{R}^d,dx)}. 
\end{align*}
But, 
\begin{align*}
\underset{n \longrightarrow+\infty}{\lim} \|\varphi D^{\alpha-1}\left(f_n\right)\|_{L^p\left(\mathbb{R}^d,dx\right)} = \|\varphi\|_{L^p(\mathbb{R}^d,\rho(dx))}. 
\end{align*}
To conclude the proof of \eqref{ineq:reverse-Holder}, it remains to prove that 
\begin{align*}
\underset{n \longrightarrow +\infty}{\lim} \int_{\mathbb{R}^d} \|u\|\nu_\alpha(du) \|f_n(\cdot + u)(\varphi(\cdot+u)-\varphi(\cdot))\|_{L^p(\mathbb{R}^d,dx)} = 0.
\end{align*}
Observe that 
\begin{align}\label{eq:cut_two}
 \int_{\mathbb{R}^d} \|u\|\nu_\alpha(du) &\|f_n(\cdot + u)(\varphi(\cdot+u)-\varphi(\cdot))\|_{L^p(\mathbb{R}^d,dx)} =  \int_{\|u\| > 1} \|u\|\nonumber \\
 &\quad \times \|f_n(\cdot + u)(\varphi(\cdot+u)-\varphi(\cdot))\|_{L^p(\mathbb{R}^d,dx)}\nu_\alpha(du)
+ \int_{\|u\|\leq 1} \|u\|\nonumber \\
&\quad \times  \|f_n(\cdot + u)(\varphi(\cdot+u)-\varphi(\cdot))\|_{L^p(\mathbb{R}^d,dx)}\nu_\alpha(du).
\end{align}
Now, let us deal with the first term of the sum on the right-hand side of the previous equality. For all $u \in \mathbb{R}^d$ such that $\|u\|>1$ and all $R>0$, 
\begin{align*}
\|f_n(\cdot + u)(\varphi(\cdot+u)-\varphi(\cdot))\|^p_{L^p(\mathbb{R}^d,dx)} & = \int_{\mathbb{R}^d} |f_n(x+u)|^p |\varphi(x+u) - \varphi(x)|^p dx \\
& = \int_{\mathbb{R}^d} |f_n(x)|^p |\varphi(x) - \varphi(x-u)|^p dx \\
& = \int_{\|x\| \leq R} |f_n(x)|^p |\varphi(x) - \varphi(x-u)|^p dx \\
& \quad \quad + \int_{\|x\|>R} |f_n(x)|^p |\varphi(x) - \varphi(x-u)|^p dx.
\end{align*}
By Lemma \ref{lem:almost_everywhere_convergence_p},
\begin{align*}
\int_{\|x\| \leq R} |f_n(x)|^p |\varphi(x) - \varphi(x-u)|^p dx \leq 2^p \|\varphi\|^p_{\infty} \int_{\|x\| \leq R} |f_n(x)|^p dx \longrightarrow 0, 
\end{align*}
as $n$ tends to $+\infty$. Moreover, by the H\"older inequality with $r = d/(d-p(\alpha-1))$ and $r' = d/(p(\alpha-1))$, 
\begin{align*}
\int_{\|x\|>R} |f_n(x)|^p |\varphi(x) - \varphi(x-u)|^p dx & \leq \left(\int_{\|x\|>R} |f_n(x)|^{p_\alpha^*}dx\right)^{\frac{1}{r}} \left(\int_{\|x\|>R} |\varphi(x) - \varphi(x-u)|^{pr'}dx\right)^{\frac{1}{r'}} \\
& \leq \sup_{n\geq 1} \|f_n\|^p_{L^{p_\alpha^*}(\mathbb{R}^d,dx)} \bigg( \|\varphi\|_{L^{\frac{d}{\alpha-1}}(B(0,R)^c,dx)} \\
& \quad\quad + \|\varphi(\cdot-u)\|_{L^{\frac{d}{\alpha-1}}(B(0,R)^c,dx)} \bigg)^p.
\end{align*}
Then, for all $u \in \mathbb{R}^d$ such that $\|u\|>1$, 
\begin{align*}
\underset{R \longrightarrow +\infty}{\lim} \sup_{n \geq 1} \int_{\|x\|>R} |f_n(x)|^p |\varphi(x) - \varphi(x-u)|^p dx = 0.
\end{align*}
Thus, for all $u \in \mathbb{R}^d$ such that $\|u\|>1$,
\begin{align*}
\underset{n \longrightarrow +\infty}{\lim}\|f_n(\cdot + u)(\varphi(\cdot+u)-\varphi(\cdot))\|_{L^p(\mathbb{R}^d,dx)} = 0.
\end{align*}
Finally, by the H\"older inequality again, 
\begin{align*}
\|f_n(\cdot + u)(\varphi(\cdot+u)-\varphi(\cdot))\|^p_{L^p(\mathbb{R}^d,dx)} & \leq \left(\int_{\mathbb{R}^d} |f_n(x)|^{p_\alpha^*}dx\right)^{\frac{1}{r}} \left(\int_{\mathbb{R}^d} |\varphi(x) - \varphi(x-u)|^{pr'}dx\right)^{\frac{1}{r'}} \\
& \leq 2^p\sup_{n\geq 1} \|f_n\|^p_{L^{p_\alpha^*}(\mathbb{R}^d,dx)} \|\varphi\|^p_{L^{\frac{d}{\alpha-1}}(\mathbb{R}^d,dx)}.
\end{align*}
The Lebesgue dominated convergence theorem ensures that 
\begin{align*}
\underset{n \longrightarrow +\infty}{\lim} \int_{\|u\| > 1} \|u\|\nu_\alpha(du) \|f_n(\cdot + u)(\varphi(\cdot+u)-\varphi(\cdot))\|_{L^p(\mathbb{R}^d,dx)}= 0.
\end{align*}
Next, let us deal with the second term of the sum on the right-hand side of \eqref{eq:cut_two}. By Taylor's formula and Jensen's inequality, for all $u \in \mathbb{R}^d$ such that $\|u\| \leq 1$, 
\begin{align*}
\|f_n(\cdot+u)(\varphi(\cdot+u)-\varphi(\cdot))\|^p_{L^p(\mathbb{R}^d)} & = \int_{\mathbb{R}^d} |f_n(x)|^p |\varphi(x)- \varphi(x-u)|^p dx \\
& =  \int_{\mathbb{R}^d} |f_n(x)|^p \left| \int_0^1 \langle \nabla(\varphi)(x-tu) ; u\rangle dt\right|^p dx \\
& \leq \|u\|^p \int_0^1 \int_{\mathbb{R}^d} |f_n(x)|^p \|\nabla(\varphi)(x-tu)\|^p dt dx. 
\end{align*}
As previously, for all $u \in \mathbb{R}^d$ such that $\|u\|\leq 1$ and all $R>0$, 
\begin{align*}
\int_{\mathbb{R}^d} |f_n(x)|^p \left(\int_0^1\|\nabla(\varphi)(x-tu)\|^p dt\right) dx & = \int_{\|x\| \leq R} |f_n(x)|^p \left(\int_0^1\|\nabla(\varphi)(x-tu)\|^p dt\right) dx \\
& \quad + \int_{\|x\| > R} |f_n(x)|^p \left(\int_0^1\|\nabla(\varphi)(x-tu)\|^p dt\right) dx.
\end{align*}
By Lemma \ref{lem:almost_everywhere_convergence_p}, for all $u \in \mathbb{R}^d$ such that $\|u\|\leq 1$ and all $R>0$,
\begin{align*}
\int_{\|x\| \leq R} |f_n(x)|^p \left(\int_0^1\|\nabla(\varphi)(x-tu)\|^p dt\right) dx \leq \|\nabla(\varphi)\|^p_{\infty} \int_{\|x\| \leq R} |f_n(x)|^p dx \longrightarrow 0,
\end{align*}
as $n$ tends to $+\infty$. Finally, by the H\"older inequality with $r = d/(d-p(\alpha-1))$ and $r' = d/(p(\alpha-1))$,
\begin{align*}
\int_{\|x\| > R} |f_n(x)|^p \left(\int_0^1\|\nabla(\varphi)(x-tu)\|^p dt\right) dx & \leq \left(\int_{\|x\| > R} |f_n(x)|^{rp}dx\right)^{\frac{1}{r}}\\
& \quad\quad \times \left(\int_{\|x\| > R} \left(\int_0^1\|\nabla(\varphi)(x-tu)\|^p dt\right)^{r'}dx\right)^{\frac{1}{r'}} \\
& \leq \sup_{n \geq 1} \|f_n\|^p_{L^{p_\alpha^*}(\mathbb{R}^d,dx)}\\
& \quad\quad \times \left(\int_{\|x\| > R} \left(\int_0^1\|\nabla(\varphi)(x-tu)\|^{pr'} dt\right)dx\right)^{\frac{1}{r'}}.
\end{align*}
But, for all $u \in \mathbb{R}^d$ such that $\|u\|\leq 1$, 
\begin{align*}
\int_{\mathbb{R}^d} \left(\int_0^1\|\nabla(\varphi)(x-tu)\|^{pr'} dt\right)dx = \int_{\mathbb{R}^d}\|\nabla(\varphi)(x)\|^{pr'}dx<+\infty.
\end{align*}
Then, for all $u \in \mathbb{R}^d$ such that $\|u\|\leq 1$,
\begin{align*}
\underset{R \longrightarrow +\infty}{\lim} \sup_{n \geq 1} \int_{\|x\| > R} |f_n(x)|^p \left(\int_0^1\|\nabla(\varphi)(x-tu)\|^p dt\right) dx = 0. 
\end{align*}
As previously, the Lebesgue dominated convergence theorem ensures that 
\begin{align*}
\underset{n \longrightarrow +\infty}{\lim} \int_{\|u\| \leq 1} \|u\|\nu_\alpha(du) \|f_n(\cdot + u)(\varphi(\cdot+u)-\varphi(\cdot))\|_{L^p(\mathbb{R}^d,dx)}= 0.
\end{align*}
Thus, 
\begin{align*}
\underset{n \longrightarrow +\infty}{\lim} \left\|f_n D^{\alpha-1}\left(\varphi\right) + R^{\alpha}\left(\varphi,f_n\right) \right\|_{L^p\left(\mathbb{R}^d,dx\right)} = 0.
\end{align*}
This concludes the proof of the inequality \eqref{ineq:reverse-Holder}. Moreover, by \cite[Lemma 1.2]{Lions_RMI185}, there exist an, at most enumerable, set of indices $I$, nonnegative reals $\tau_i$, $\rho_i$, $i \in I$, and distinct points $x_i \in \mathbb{R}^d$, $i \in I$, such that
\begin{align*}
\tau(dx) = \sum_{i \in I} \tau_i \delta_{x_i}(dx), \quad \rho(dx) \geq \sum_{i \in I} \rho_i \delta_{x_i}(dx).  
\end{align*}
Next, let us prove the inequality \eqref{ineq:rH_atomic} still for the case $f = 0$. Let $i \in I$ and, without loss of generality, assume that $x_i = 0$. Let $\eta$ be a function in $\mathcal{C}^{\infty}_c(\mathbb{R}^d)$ such that $0 \leq \eta \leq 1$, $\eta(0) = 1$ and $\operatorname{supp}\left(\eta\right) = B(0,1)$. For all $\varepsilon \in (0,1)$, let $\eta_\varepsilon$ be defined, for all $x \in \mathbb{R}^d$, by
\begin{align*}
\eta_{\varepsilon}(x) = \eta \left(\frac{x}{\varepsilon}\right). 
\end{align*}
Then, by the reverse H\"older inequality, for all $\varepsilon \in (0,1)$, 
\begin{align*}
S_{p, \alpha,d,\gamma}^{\frac{1}{p}} \|\eta_\varepsilon\|_{L^{p_\alpha^*}(\mathbb{R}^d,\tau(dx))} \leq \|\eta_\varepsilon\|_{L^{p}(\mathbb{R}^d,\rho(dx))}.
\end{align*}
Now, for $\varepsilon \in (0,1)$ small enough, 
\begin{align*}
\int_{\mathbb{R}^d} \left|\eta_\varepsilon(x)\right|^{p_\alpha^*} \tau(dx) = \sum_{j \in I} \tau_j \left|\eta_{\varepsilon}(x_j)\right|^{p_\alpha^*} = \tau_i. 
\end{align*}
Moreover, for all $\varepsilon>0$, 
\begin{align*}
\int_{B(0,\varepsilon)} \left|\eta_\varepsilon\left(x\right)\right|^p \rho(dx) \leq \rho\left(B\left(0,\varepsilon\right)\right). 
\end{align*} 
Thus, for $\varepsilon \in (0,1)$ small enough, 
\begin{align*}
S_{p,\alpha,d,\gamma}^{\frac{1}{p}}\tau_i^{\frac{1}{p_\alpha^*}} \leq \rho\left(B\left(0,\varepsilon\right)\right)^{\frac{1}{p}}. 
\end{align*}
But, $\underset{\varepsilon \rightarrow 0^+}{\lim} \rho \left(B(0, \varepsilon)\right) = \rho \left(\{0\}\right)$, thus, 
\begin{align*}
S_{p,\alpha,d,\gamma}\tau_i^{\frac{p}{p_\alpha^*}} \leq \rho \left(\{0\}\right),
\end{align*}
and the inequality \eqref{ineq:rH_atomic} is proved in this case. Next, let us treat the general case $f \ne 0$. Let $(g_n)_{n \geq 1}$ be the sequence of functions defined, for all $n \geq 1$, by $g_n = f_n - f$. By the Banach-Alaoglu-Bourbaki theorem, in the sense of the weak-$\ast$ convergence of measures, 
\begin{align*}
\left| g_n(x) \right|^{p_\alpha^*} dx \longrightarrow \tilde{\tau}(dx),\quad \left\|D^{\alpha-1}\left(g_n\right)(x)\right\|^p dx \longrightarrow \tilde{\rho}(dx), 
\end{align*} 
as $n$ tends to $+\infty$ (up to a subsequence). Now, by the Br\'ezis-Lieb lemma (see \cite[Theorem 1.9.]{Lieb_Loss_book01}), for all $\varphi \in \mathcal{C}^{\infty}_c(\mathbb{R}^d)$, 
\begin{align*}
& \left| \int_{\mathbb{R}^d} |\varphi(x)|^{p_\alpha^*} |f_n(x)|^{p_\alpha^*}dx - \int_{\mathbb{R}^d} |\varphi(x)|^{p_\alpha^*} |g_n(x)|^{p_\alpha^*}dx - \int_{\mathbb{R}^d} |\varphi(x)|^{p_\alpha^*} |f(x)|^{p_\alpha^*}dx \right| \\& \leq \| \varphi\|^{p_\alpha^*}_{\infty} \int_{\mathbb{R}^d} \left| |f_n(x)|^{p_\alpha^*} - |g_n(x)|^{p_\alpha^*} -  |f(x)|^{p_\alpha^*} \right| dx \longrightarrow 0, 
\end{align*}
as $n$ tends to $+\infty$. Then, for all $\varphi \in \mathcal{C}^{\infty}_c\left(\mathbb{R}^d\right)$, 
\begin{align*}
\int_{\mathbb{R}^d} \left|\varphi(x)\right|^{p_\alpha^*} \tilde{\tau}(dx) = \int_{\mathbb{R}^d} \left| \varphi(x)\right|^{p_\alpha^*} \tau(dx) - \int_{\mathbb{R}^d} |\varphi(x)|^{p_\alpha^*} |f(x)|^{p_\alpha^*}dx.   
\end{align*}
Moreover, the sequence $(g_n)_{n \geq 1}$ satisfies the conditions of the previous case ($f=0$) so that, for all $\varphi \in \mathcal{C}^{\infty}_c(\mathbb{R}^d)$, 
\begin{align*}
S_{p,\alpha,d,\gamma}^{\frac{1}{p}} \|\varphi\|_{L^{p_\alpha^*}\left(\mathbb{R}^d,\tilde{\tau}(dx)\right)} \leq \|\varphi\|_{L^{p}\left(\mathbb{R}^d,\tilde{\rho}(dx)\right)}.
\end{align*}
From the previous inequality, again by \cite[Lemma 1.2]{Lions_RMI185}, 
\begin{align*}
\tilde{\tau}(dx) = \sum_{j \in J} \tilde{\tau}_j \delta_{y_j}(dx),\quad \tilde{\rho}(dx) \geq \sum_{j \in J} \tilde{\rho}_j \delta_{y_j}(dx), 
\end{align*}
for an at most enumerable set of indices $J$, nonnegative reals $\tilde{\tau}_j$, $\tilde{\rho}_j$, $j \in J$, and distinct points $y_j \in \mathbb{R}^d$, $j \in J$ and, for all $j \in J$, 
\begin{align*}
S_{p,\alpha,d,\gamma} \tilde{\tau}_j^{\frac{p}{p_\alpha^*}} \leq \tilde{\rho}_j. 
\end{align*} 
Thus, in the sense of weak-$\ast$ convergence of measures,  
\begin{align*}
|f_n(x)|^{p_\alpha^*}dx \longrightarrow \tau(dx) = |f(x)|^{p_\alpha^*} dx + \sum_{j \in J} \tilde{\tau}_j \delta_{y_j}(dx), \quad n \longrightarrow +\infty.  
\end{align*}
Now, by the fractional Sobolev inequality, for all $n \geq 1$ and all $\varphi \in \mathcal{C}_c^{\infty}\left(\mathbb{R}^d\right)$, 
\begin{align*}
S_{p,\alpha,d,\gamma}^{\frac{1}{p}} \|\varphi f_n\|_{L^{p_\alpha^*}(\mathbb{R}^d,dx)} & \leq \|D^{\alpha-1}(\varphi f_n)\|_{L^p\left(\mathbb{R}^d,dx\right)} \\
& \leq \|\varphi D^{\alpha-1}(f_n)\|_{L^p\left(\mathbb{R}^d,dx\right)} + \|f_n D^{\alpha-1}\left(\varphi\right) + R^{\alpha}\left(\varphi,f_n\right)\|_{L^p\left(\mathbb{R}^d,dx\right)}.
\end{align*}
As previously, for all $\varphi \in \mathcal{C}_c^{\infty}\left(\mathbb{R}^d\right)$, 
\begin{align*}
\underset{n \longrightarrow +\infty}{\lim}\|\varphi f_n\|_{L^{p_\alpha^*}(\mathbb{R}^d,dx)} = \|\varphi\|_{L^{p_\alpha^*}\left(\mathbb{R}^d, \tau(dx)\right)} , \quad \underset{n \longrightarrow +\infty}{\lim}\|\varphi D^{\alpha-1}(f_n)\|_{L^{p}(\mathbb{R}^d,dx)} = \|\varphi\|_{L^{p}\left(\mathbb{R}^d, \rho(dx)\right)}.
\end{align*}
Next, let us prove that 
\begin{align*}
\underset{n\longrightarrow +\infty}{\lim} \|f_n D^{\alpha-1}\left(\varphi\right) + R^{\alpha}\left(\varphi,f_n\right)\|_{L^p\left(\mathbb{R}^d,dx\right)} = \|f D^{\alpha-1}\left(\varphi\right) + R^{\alpha}\left(\varphi,f\right)\|_{L^p\left(\mathbb{R}^d,dx\right)}.
\end{align*}
For all $n \geq 1$ and all $\varphi \in \mathcal{C}_c^{\infty}\left(\mathbb{R}^d\right)$, 
\begin{align*}
&\left| \|f_n D^{\alpha-1}\left(\varphi\right) + R^{\alpha}\left(\varphi,f_n\right)\|_{L^p\left(\mathbb{R}^d,dx\right)} - \|f D^{\alpha-1}\left(\varphi\right) + R^{\alpha}\left(\varphi,f\right)\|_{L^p\left(\mathbb{R}^d,dx\right)} \right| \\
&\quad\quad \leq \|g_n D^{\alpha-1}\left(\varphi\right) + R^{\alpha}\left(\varphi,g_n\right)\|_{L^p\left(\mathbb{R}^d,dx\right)}.  
\end{align*}
Now, from the previous case ($f=0$), 
\begin{align*}
\underset{n \longrightarrow +\infty}{\lim} \|g_n D^{\alpha-1}\left(\varphi\right) + R^{\alpha}\left(\varphi,g_n\right)\|_{L^p\left(\mathbb{R}^d,dx\right)} = 0.
\end{align*}
Then, for all $\varphi \in \mathcal{C}^{\infty}_c(\mathbb{R}^d)$, 
\begin{align*}
S_{p,\alpha,d,\gamma}^{\frac{1}{p}} \|\varphi \|_{L^{p_\alpha^*}(\mathbb{R}^d,\tau(dx))} \leq \|\varphi \|_{L^p\left(\mathbb{R}^d,\rho(dx)\right)} + \| f D^{\alpha-1}\left(\varphi\right) + R^{\alpha}\left(\varphi,f\right)\|_{L^p\left(\mathbb{R}^d,dx\right)}.
\end{align*}
As before, let $\eta \in \mathcal{C}_c^{\infty}\left(\mathbb{R}^d\right)$ be such that $0 \leq \eta \leq 1$, $\eta(0) = 1$ and $\operatorname{supp}(\eta) = B(0,1)$. For all $\varepsilon \in (0,1)$, let $\eta_\varepsilon$ be defined, for all $x \in \mathbb{R}^d$, by 
\begin{align*}
\eta_\varepsilon(x) = \eta\left(\frac{x}{\varepsilon}\right). 
\end{align*}
Let $j$ be fixed in $J$ and, without loss of generality, let us assume that $y_j=0$. Then, for all $\varepsilon \in (0,1)$, 
\begin{align*}
S_{p,\alpha,d,\gamma}^{\frac{1}{p}} \| \eta_\varepsilon \|_{L^{p_\alpha^*}(\mathbb{R}^d,\tau(dx))} \leq \| \eta_\varepsilon \|_{L^p\left(\mathbb{R}^d,\rho(dx)\right)} +  \| f D^{\alpha-1}\left(\eta_\varepsilon\right) + R^{\alpha}\left(\eta_\varepsilon,f\right)\|_{L^p\left(\mathbb{R}^d,dx\right)}.
\end{align*}
Now, 
\begin{align*}
\int_{\mathbb{R}^d} \left|\eta\left(\frac{x}{\varepsilon}\right)\right|^{p_\alpha^*} \tau(dx) = \int_{\mathbb{R}^d} \left|\eta\left(\frac{x}{\varepsilon}\right)\right|^{p_\alpha^*} |f(x)|^{p_\alpha^*} dx + \sum_{k \in J} \tilde{\tau}_k \left|\eta_\varepsilon(y_k)\right|^{p_\alpha^*}. 
\end{align*}
Moreover, for $\varepsilon \in (0,1)$ small enough, 
\begin{align*}
\sum_{k \in J} \tilde{\tau}_k |\eta_\varepsilon(y_k)|^{p_\alpha^*} = \tilde{\tau}_j.
\end{align*}
Thus, for $\varepsilon \in (0,1)$ small enough, 
\begin{align*}
S_{p,\alpha,d,\gamma}^{\frac{1}{p}} \| \eta_\varepsilon \|_{L^{p_\alpha^*}(\mathbb{R}^d,\tau(dx))} \geq S_{p,\alpha,d,\gamma}^{\frac{1}{p}} (\tilde{\tau}_j)^{\frac{1}{p_\alpha^*}}. 
\end{align*}
Furthermore, for all $\varepsilon \in (0,1)$, 
\begin{align*}
\| \eta_\varepsilon \|_{L^p\left(\mathbb{R}^d,\rho(dx)\right)} \leq \rho\left(B(0,\varepsilon)\right)^{\frac{1}{p}}. 
\end{align*}
Then, for $\varepsilon \in (0,1)$ small enough, 
\begin{align}\label{ineq:reverse_Holder_inequality_epsilon_gen}
S_{p,\alpha,d,\gamma}^{\frac{1}{p}} (\tilde{\tau}_j)^{\frac{1}{p_\alpha^*}} \leq \rho\left(B(0,\varepsilon)\right)^{\frac{1}{p}} + \| f D^{\alpha-1}\left(\eta_\varepsilon\right) + R^{\alpha}\left(\eta_\varepsilon,f\right)\|_{L^p\left(\mathbb{R}^d,dx\right)}. 
\end{align}
It remains to prove that 
\begin{align*}
\underset{\varepsilon \rightarrow 0^+}{\lim} \| f D^{\alpha-1}\left(\eta_\varepsilon\right) + R^{\alpha}\left(\eta_\varepsilon,f\right)\|_{L^p\left(\mathbb{R}^d,dx\right)} = 0.
\end{align*}
By Minkowski's integral inequality, for all $\varepsilon \in (0,1)$, 
\begin{align}\label{ineq:cut_two_2}
\| f D^{\alpha-1}\left(\eta_\varepsilon\right) + R^{\alpha}\left(\eta_\varepsilon,f\right)\|_{L^p\left(\mathbb{R}^d,dx\right)} & \leq \varepsilon^{1- \alpha} \int_{\mathbb{R}^d} \|u\| \nu_\alpha(du) \|f(\cdot+\varepsilon u) \Delta_u(\eta)(\cdot/\varepsilon) \|_{L^p(\mathbb{R}^d,dx)} \nonumber \\
& \leq \varepsilon^{1- \alpha} \bigg( \int_{\|u\|>1} \|u\| \nu_\alpha(du) \|f(\cdot+\varepsilon u) \Delta_u(\eta)(\cdot/\varepsilon) \|_{L^p(\mathbb{R}^d,dx)} \nonumber \\
& \quad\quad + \int_{\|u\| \leq 1} \|u\| \nu_\alpha(du) \|f(\cdot+\varepsilon u) \Delta_u(\eta)(\cdot/\varepsilon) \|_{L^p(\mathbb{R}^d,dx)} \bigg).
\end{align}
First, let us deal with the first term of the sum on the right-hand side of the previous inequality. For all $u \in \mathbb{R}^d$ such that $\|u\|>1$, all $\varepsilon \in (0,1)$ and all $R>0$, 
\begin{align*}
\|f(\cdot+\varepsilon u) \Delta_u(\eta)(\cdot/\varepsilon) \|^p_{L^p(\mathbb{R}^d,dx)} & = \int_{B(0,R\varepsilon)} |f(x+\varepsilon u)|^p |\Delta_u(\eta)\left(\frac{x}{\varepsilon}\right)|^p dx \\
&\quad \quad + \int_{B(0,R\varepsilon)^c} |f(x+\varepsilon u)|^p |\Delta_u(\eta)\left(\frac{x}{\varepsilon}\right)|^p dx. 
\end{align*}
Now, by the H\"older inequality with $r = d/(d-p(\alpha-1))$ and $r' = d/(p(\alpha-1))$, 
\begin{align*}
\int_{B(0,R\varepsilon)} |f(x+\varepsilon u)|^p |\Delta_u(\eta)\left(\frac{x}{\varepsilon}\right)|^p dx & \leq \left(\int_{B(0,R\varepsilon)} |f(x+\varepsilon u)|^{p_\alpha^*} dx\right)^{\frac{1}{r}} \left(\int_{B(0,R\varepsilon)} |\Delta_u(\eta)\left(\frac{x}{\varepsilon}\right)|^{pr'} dx\right)^{\frac{1}{r'}} \\
& \leq 2^p \left(\int_{B(0,R\varepsilon)} |f(x+\varepsilon u)|^{p_\alpha^*} dx\right)^{\frac{1}{r}} \left(\int_{\mathbb{R}^d} |\eta\left(\frac{x}{\varepsilon}\right)|^{pr'} dx\right)^{\frac{1}{r'}} \\
& \leq 2^p \varepsilon^{p(\alpha-1)} \left(\int_{B(0,R\varepsilon)} |f(x+\varepsilon u)|^{p_\alpha^*} dx\right)^{\frac{1}{r}} \left(\int_{\mathbb{R}^d} |\eta\left(x\right)|^{pr'} dx\right)^{\frac{1}{r'}}.
\end{align*}
Moreover, for all $u \in \mathbb{R}^d$ fixed and all $R>0$, 
\begin{align*}
\underset{\varepsilon \rightarrow 0^+}{\lim} \int_{B(0,R\varepsilon)} |f(x+\varepsilon u)|^{p_\alpha^*} dx = 0. 
\end{align*}
Now, by the H\"older inequality again, 
\begin{align*}
\int_{B(0,R\varepsilon)^c} |f(x+\varepsilon u)|^p |\Delta_u(\eta)\left(\frac{x}{\varepsilon}\right)|^p dx & \leq \left( \int_{\mathbb{R}^d} |f(x)|^{p_\alpha^*} dx \right)^{\frac{1}{r}} \left(\int_{B(0,R\varepsilon)^c} |\Delta_u(\eta)\left(\frac{x}{\varepsilon}\right)|^{pr'}dx\right)^{\frac{1}{r'}} \\
& \leq \varepsilon^{(\alpha-1)p} \left( \int_{\mathbb{R}^d} |f(x)|^{p_\alpha^*} dx \right)^{\frac{1}{r}} \left(\int_{B(0,R)^c} |\Delta_u(\eta)\left(x\right)|^{pr'}dx\right)^{\frac{1}{r'}}.
\end{align*}
Thus, 
\begin{align*}
\underset{R \longrightarrow +\infty}{\lim} \sup_{\varepsilon \in (0,1)} \varepsilon^{(1-\alpha)p}  \int_{B(0,R\varepsilon)^c} |f(x+\varepsilon u)|^p |\Delta_u(\eta)\left(\frac{x}{\varepsilon}\right)|^p dx = 0.
\end{align*}
Then, for all $u \in \mathbb{R}^d$ such that $\|u\|>1$, 
\begin{align*}
\underset{\varepsilon \rightarrow 0^+}{\lim} \varepsilon^{1-\alpha} \|f(\cdot+\varepsilon u) \Delta_u(\eta)(\cdot/\varepsilon) \|_{L^p(\mathbb{R}^d,dx)} = 0.
\end{align*} 
Moreover, by the H\"older inequality, for all $u \in \mathbb{R}^d$ such that $\|u\|>1$ and all $\varepsilon \in (0,1)$, 
\begin{align*}
\varepsilon^{1-\alpha} \|f(\cdot+\varepsilon u) \Delta_u(\eta)(\cdot/\varepsilon) \|_{L^p(\mathbb{R}^d,dx)} \leq 2 \|f\|_{L^{p_\alpha^*}(\mathbb{R}^d,dx)} \|\eta\|_{L^{pr'}(\mathbb{R}^d,dx)}. 
\end{align*}
The Lebesgue dominated convergence theorem ensures that 
\begin{align*}
\underset{\varepsilon \rightarrow 0^+}{\lim} \varepsilon^{1- \alpha} \int_{\|u\|>1} \|u\| \nu_\alpha(du) \|f(\cdot+\varepsilon u) \Delta_u(\eta)(\cdot/\varepsilon) \|_{L^p(\mathbb{R}^d,dx)} = 0.
\end{align*}
A similar analysis allows to prove that 
\begin{align*}
\underset{\varepsilon \rightarrow 0^+}{\lim} \varepsilon^{1- \alpha} \int_{\|u\| \leq 1} \|u\| \nu_\alpha(du) \|f(\cdot+\varepsilon u) \Delta_u(\eta)(\cdot/\varepsilon) \|_{L^p(\mathbb{R}^d,dx)} = 0.
\end{align*}
Passing to the limit in \eqref{ineq:reverse_Holder_inequality_epsilon_gen}, gives
\begin{align*}
S_{p,\alpha,d,\gamma} (\tilde{\tau}_j)^{\frac{p}{p_\alpha^*}} \leq \rho\left(\{0\}\right). 
\end{align*}
Moreover, by weak convergence, $\rho(dx) \geq \|D^{\alpha-1}\left(f\right)(x)\|^pdx$. Since $\|D^{\alpha-1}\left(f\right)(x)\|^pdx$ and\\
$\sum_{j \in J} \rho_j \delta_{y_j}(dx)$ are mutually singular, the first point of the statement of Theorem \ref{thm:ccp_homo_frac_p} follows.\\
\\
\textbf{Step 2:} In this step, let us prove the last two points of the statement of Theorem \ref{thm:ccp_homo_frac_p} following the strategy of the fifth step of the proof of \cite[Theorem $3.1$]{PQR_jmaa22} (see also \cite[Proof of Theorem 1.1.]{BSS_NoDEA18}).~Let $\Phi$ be an infinitely continuously differentiable function on $\mathbb{R}_+$ such that $\Phi(x) = 0$, for $x \in [0,1]$, $\Phi(x) = 1$, for $x \in [2,+\infty)$ and $0\leq \Phi \leq 1$. Next, for all $R>0$, let $\Phi_R$ be defined, for all $x \in \mathbb{R}^d$, by 
\begin{align*}
\Phi_R(x) = \Phi\left(\frac{\|x\|}{R}\right). 
\end{align*}  
Then, for all $n \geq 1$ an integer and all $R>0$, 
\begin{align*}
\int_{\mathbb{R}^d} \|D^{\alpha-1}\left(f_n\right)(x)\|^pdx  = \int_{\mathbb{R}^d} \|D^{\alpha-1}\left(f_n\right)(x)\|^p \Phi_R(x)^p dx + \int_{\mathbb{R}^d} D^{\alpha-1}\left(f_n\right)(x)\|^p \left(1-\Phi_R(x)^p\right) dx. 
\end{align*}  
Now, by the definition of $\Phi$, for all $R>0$ and all $n \geq 1$ integer, 
\begin{align*}
\int_{\|x\| \geq 2 R} \|D^{\alpha-1}\left(f_n\right)(x)\|^pdx\leq \int_{\mathbb{R}^d} \|D^{\alpha-1}\left(f_n\right)(x)\|^p \Phi_R(x)^p dx \leq\int_{\|x\| \geq R} \|D^{\alpha-1}\left(f_n\right)(x)\|^pdx.  
\end{align*} 
Thus, 
\begin{align*}
\rho_\infty = \underset{R \longrightarrow +\infty}{\lim}\, \underset{n \longrightarrow +\infty}{\limsup} \int_{\mathbb{R}^d} \|D^{\alpha-1}\left(f_n\right)(x)\|^p \Phi_R(x)^p dx.
\end{align*}
Similarly, 
\begin{align*}
\tau_\infty = \underset{R \longrightarrow +\infty}{\lim}\, \underset{n \longrightarrow +\infty}{\limsup} \int_{\mathbb{R}^d} |f_n(x)|^{p_\alpha^*} \Phi_R(x)^{p_\alpha^*} dx.
\end{align*}
Moreover, since $(1-\Phi_R^p)$ is continuous with compact support on $\mathbb{R}^d$, for all $R>0$, 
\begin{align*}
\underset{n \longrightarrow +\infty}{\lim} \int_{\mathbb{R}^d} \|D^{\alpha-1}\left(f_n\right)(x)\|^p (1-\Phi_R(x)^p) dx = \int_{\mathbb{R}^d} (1-\Phi_R(x)^p) \rho(dx).
\end{align*}
Then, by the Lebesgue dominated convergence theorem, 
\begin{align*}
\underset{R \longrightarrow +\infty}{\lim} \int_{\mathbb{R}^d} (1-\Phi_R(x)^p) \rho(dx) = \rho\big(\mathbb{R}^d\big). 
\end{align*}
Thus, 
\begin{align*}
\underset{n \longrightarrow +\infty}\limsup \int_{\mathbb{R}^d} \|D^{\alpha-1}(f_n)(x)\|^p dx = \rho\big(\mathbb{R}^d\big) + \rho_\infty,
\end{align*}
and 
\begin{align*}
\underset{n \longrightarrow +\infty}\limsup \int_{\mathbb{R}^d} |f_n(x)|^{p_\alpha^*} dx = \tau\big(\mathbb{R}^d\big) + \tau_\infty.
\end{align*}
This concludes the proof of the theorem. 
\end{proof}
\noindent
\begin{lem}\label{lem:reverse_Holder_infinity}
Let $\rho_\infty$ and $\tau_\infty$ be as defined in Theorem \ref{thm:ccp_homo_frac_p}. Then, 
\begin{align}\label{ineq:reverse_Holder_inf}
S_{p,\alpha,d,\gamma} \tau_\infty^{\frac{p}{p_\alpha^*}} \leq \rho_\infty.
\end{align}
\end{lem}

\begin{proof}
Let $\Phi$ be an infinitely continuously differentiable function on $\mathbb{R}_+$ such that $\Phi(x) = 0$, for $x \in [0,1]$, $\Phi(x) = 1$, for $x \in [2,+\infty)$ and $0\leq \Phi \leq 1$. Next, for all $R>0$, let $\Phi_R$ be defined, for all $x \in \mathbb{R}^d$, by 
\begin{align*}
\Phi_R(x) = \Phi\left(\frac{\|x\|}{R}\right). 
\end{align*}  
Let $(f_n)_{n \geq 1}$ and $f$ be as before and let $(g_n)_{n \geq 1}$ be defined, for all $n \geq 1$, by $g_n=f_n-f$. By the fractional Sobolev inequality, for all $n \geq 1$ and all $R>0$, 
\begin{align}\label{ineq:pre_reverse_Holder_inf}
S_{p,\alpha,d,\gamma}^{\frac{1}{p}} \|g_n \Phi_R\|_{L^{p_\alpha^*}\left(\mathbb{R}^d,dx\right)} \leq \|D^{\alpha-1}\left(g_n \Phi_R\right)\|_{L^p(\mathbb{R}^d,dx)}.
\end{align}
Moreover, for all $R>0$ and all $n \geq 1$ integer, 
\begin{align*}
\|D^{\alpha-1}\left(g_n \Phi_R\right)\|_{L^p\left(\mathbb{R}^d,dx\right)} \leq \|g_n D^{\alpha-1}\left(\Phi_R\right)+R^\alpha(g_n,\Phi_R)\|_{L^p\left(\mathbb{R}^d,dx\right)} + \|\Phi_R D^{\alpha-1}\left(g_n\right)\|_{L^p\left(\mathbb{R}^d,dx\right)}.
\end{align*}
Now, by the first step in the proof of Theorem \ref{thm:ccp_homo_frac_p},  
\begin{align*}
\underset{R\rightarrow +\infty}{\lim}\, \underset{n \rightarrow +\infty}{\lim} \|g_n D^{\alpha-1}\left(\Phi_R\right)+R^\alpha(g_n,\Phi_R)\|_{L^p\left(\mathbb{R}^d,dx\right)} & = \underset{R\rightarrow +\infty}{\lim}\, \underset{n \rightarrow +\infty}{\lim} \|g_n D^{\alpha-1}\left(1-\Phi_R\right)\\
&\quad\quad +R^\alpha(g_n,1-\Phi_R)\|_{L^p\left(\mathbb{R}^d,dx\right)} = 0. 
\end{align*}
To conclude, let us prove that: 
\begin{align*}
\underset{R \longrightarrow +\infty}{\lim}\, \underset{n \longrightarrow +\infty}{\limsup} \left(\int_{\mathbb{R}^d} |g_n(x)|^{p_\alpha^*} |\Phi_R(x)|^{p_\alpha^*}dx\right)^{\frac{1}{p_\alpha^*}} = \tau_\infty^{\frac{1}{p_\alpha^*}}
\end{align*}
and 
\begin{align*}
\underset{R \longrightarrow +\infty}{\lim}\, \underset{n \longrightarrow +\infty}{\limsup} \left(\int_{\mathbb{R}^d} \|D^{\alpha-1}(g_n)(x)\|^{p} |\Phi_R(x)|^{p}dx\right)^{\frac{1}{p}} = \rho_\infty^{\frac{1}{p}}.
\end{align*}
Recall that, for all $\varepsilon>0$ and all $r>1$, there exists $C_{\varepsilon,r}>0$ such that, for all $a,b \in \mathbb{R}$, 
\begin{align}\label{ineq:error_triangular_inequality}
\left| |a+b|^r - |b|^r \right| \leq \varepsilon |b|^r + C_{\varepsilon,r} |a|^{r}. 
\end{align}
Thus, taking $a = f(x)$ and $b = g_n(x)$, for all $n \geq 1$ and all $x \in \mathbb{R}^d$, 
\begin{align*}
\left| |f_n(x)|^{p_\alpha^*} - |g_n(x)|^{p_\alpha^*} \right| \leq \varepsilon |g_n(x)|^{p_\alpha^*} + C_{\varepsilon,p_{\alpha}^*} |f(x)|^{p_\alpha^*},
\end{align*}
which implies, for all $R>0$, 
\begin{align*}
\int_{\mathbb{R}^d}\left| |f_n(x)|^{p_\alpha^*} - |g_n(x)|^{p_\alpha^*} \right| \Phi_R(x)^{p_\alpha^*}dx \leq \varepsilon \int_{\mathbb{R}^d} |g_n(x)|^{p_\alpha^*} \Phi_R(x)^{p_\alpha^*} dx + C_{\varepsilon,p_{\alpha}^*} \int_{\mathbb{R}^d} |f(x)|^{p_\alpha^*}\Phi_R(x)^{p_\alpha^*} dx. 
\end{align*}
Now, by the uniform boundedness principle, for all $n \geq 1$ and all $R>0$, 
\begin{align*}
\int_{\mathbb{R}^d} |g_n(x)|^{p_\alpha^*} \Phi_R(x)^{p_\alpha^*} dx & \leq \int_{\mathbb{R}^d} |g_n(x)|^{p_\alpha^*} dx \leq 2^{p-1} \left( \sup_{n \geq 1} \left(\int_{\mathbb{R}^d} |f_n(x)|^{p_\alpha^*} dx\right) + \int_{\mathbb{R}^d} |f(x)|^{p_{\alpha}^*} dx \right).
\end{align*}
Thus, for all $R>0$, 
\begin{align*}
\underset{n \rightarrow +\infty}{\limsup}\int_{\mathbb{R}^d}\left| |f_n(x)|^{p_\alpha^*} - |g_n(x)|^{p_\alpha^*} \right| \Phi_R(x)^{p_\alpha^*}dx & \leq 2^{p-1} \varepsilon \left(\sup_{n \geq 1} \|f_n\|^{p_\alpha^*}_{L^{p_\alpha^*}(\mathbb{R}^d,dx)} + \|f\|^{p_\alpha^*}_{L^{p_\alpha^*}(\mathbb{R}^d,dx)} \right) \\ 
&\quad\quad +  C_{\varepsilon,p_{\alpha}^*} \int_{\mathbb{R}^d} |f(x)|^{p_\alpha^*}\Phi_R(x)^{p_\alpha^*} dx.
\end{align*}
By the very definition of $\Phi_R$, for all $R>0$, 
\begin{align*}
\int_{\mathbb{R}^d} |f(x)|^{p_\alpha^*}\Phi_R(x)^{p_\alpha^*} dx \leq \int_{\|x\| \geq R} |f(x)|^{p_\alpha^*} dx. 
\end{align*}
Since $f \in L^{p_\alpha^*}(\mathbb{R}^d,dx)$, 
\begin{align*}
\underset{R \rightarrow +\infty}{\lim}\int_{\|x\| \geq R} |f(x)|^{p_\alpha^*} dx = 0. 
\end{align*}
Thus, for all $\varepsilon>0$,  
\begin{align*}
\underset{R \rightarrow +\infty}{\lim}\underset{n \rightarrow +\infty}{\limsup}\int_{\mathbb{R}^d}\left| |f_n(x)|^{p_\alpha^*} - |g_n(x)|^{p_\alpha^*} \right| \Phi_R(x)^{p_\alpha^*}dx \leq 2^{p-1} \varepsilon \big(\sup_{n \geq 1} \|f_n\|^{p_\alpha^*}_{L^{p_\alpha^*}(\mathbb{R}^d,dx)} + \|f\|^{p_\alpha^*}_{L^{p_\alpha^*}(\mathbb{R}^d,dx)} \big). 
\end{align*}
Letting $\varepsilon$ tend to $0^+$ ensures that 
\begin{align*}
\underset{R \rightarrow +\infty}{\lim}\underset{n \rightarrow +\infty}{\limsup}\int_{\mathbb{R}^d}\left| |f_n(x)|^{p_\alpha^*} - |g_n(x)|^{p_\alpha^*} \right| \Phi_R(x)^{p_\alpha^*}dx = 0
\end{align*}
which implies that
\begin{align*}
\underset{R \rightarrow +\infty}{\lim}\underset{n \rightarrow +\infty}{\limsup}\int_{\mathbb{R}^d} |g_n(x)|^{p_\alpha^*}\Phi_R(x)^{p_\alpha^*}dx = \underset{R \rightarrow +\infty}{\lim}\underset{n \rightarrow +\infty}{\limsup}\int_{\mathbb{R}^d} |f_n(x)|^{p_\alpha^*}\Phi_R(x)^{p_\alpha^*}dx = \tau_\infty. 
\end{align*}
Similarly, using a $d$-dimensional version of the inequality \eqref{ineq:error_triangular_inequality}, 
\begin{align*}
\underset{R \rightarrow +\infty}{\lim}\underset{n \rightarrow +\infty}{\limsup}\int_{\mathbb{R}^d} \|D^{\alpha-1}(g_n)(x)\|^{p}\Phi_R(x)^{p}dx = \underset{R \rightarrow +\infty}{\lim}\underset{n \rightarrow +\infty}{\limsup}\int_{\mathbb{R}^d} \|D^{\alpha-1}(f_n)(x)\|^{p}\Phi_R(x)^{p}dx = \rho_\infty.
\end{align*}
Passing to the limits in \eqref{ineq:pre_reverse_Holder_inf} gives 
\begin{align*}
S_{p, \alpha,d,\gamma}^{\frac{1}{p}} \tau_{\infty}^{\frac{1}{p_\alpha^*}} \leq \rho_\infty^{\frac{1}{p}}. 
\end{align*}
This concludes the proof of the lemma. 
\end{proof}
\noindent
To end this section, let us prove Theorem \ref{thm:existence_non-trivial_optimizer_p}.\\

\noindent
\textit{Proof of Theorem \ref{thm:existence_non-trivial_optimizer_p}.}
Let $(f_n)_{n \geq 1}$ be a minimizing sequence for the optimization problem \eqref{eq:minimization_problem_frac_sobolev-p}, i.e., for all $n \geq 1$, $f_n \in \dot{W}^{\alpha-1,p} \left(\mathbb{R}^d,dx\right)$, $\|f_n\|_{L^{p_\alpha^*}(\mathbb{R}^d,dx)} = 1$ and
\begin{align*}
\underset{n \rightarrow +\infty}{\lim}\|D^{\alpha-1}(f_n)\|^p_{L^p(\mathbb{R}^d,dx)}= S_{p,\alpha,d,\gamma}. 
\end{align*} 
Then, by Lemma \ref{lem:there_is_something_somewhere_non-trivial_p} and Lemma \ref{lem:almost_everywhere_convergence_p}, there exists $f \in \dot{W}^{\alpha-1,p}\left(\mathbb{R}^d,dx\right) \setminus \{0\}$ such that  
\begin{align*}
f_n \longrightarrow f, \quad n \longrightarrow +\infty,
\end{align*}
weakly in $L^{p_\alpha^*}\left(\mathbb{R}^d,dx\right)$ and
\begin{align*}
D^{\alpha-1}(f_n) \longrightarrow D^{\alpha-1}(f), \quad n \longrightarrow +\infty,
\end{align*}
weakly in $L^p\left(\mathbb{R}^d,\mathbb{R}^d,dx\right)$ (up to symmetries and subsequences). Now, by the lower semicontinuity of the norm and weak convergence, 
\begin{align*}
\|D^{\alpha-1}(f)\|^p_{L^p(\mathbb{R}^d,dx)} \leq S_{p, \alpha,d,\gamma},
\end{align*}
and moreover, by the very definition of $S_{p,\alpha,d,\gamma}$, 
\begin{align*}
S_{p,\alpha,d,\gamma} \|f\|^p_{L^{p_\alpha^*}\left(\mathbb{R}^d,dx\right)} \leq \|D^{\alpha-1}(f)\|^p_{L^p(\mathbb{R}^d,dx)}. 
\end{align*}
Next, let us prove that $\|f\|_{L^{p_\alpha^*}(\mathbb{R}^d,dx)} = 1$. Since $(f_n)_{n \geq 1}$ is a minimizing sequence, thanks to the Banach-Alaoglu–Bourbaki theorem and Theorem \ref{thm:ccp_homo_frac_p}, 
\begin{align*}
1 = \dfrac{\rho\left(\mathbb{R}^d\right) + \rho_\infty}{S_{p,\alpha,d,\gamma}}.
\end{align*}
Now, by the fractional Sobolev inequality, 
\begin{align*}
\rho\left(\mathbb{R}^d\right) & \geq \int_{\mathbb{R}^d} \|D^{\alpha-1}(f)(x)\|^p dx +  \sum_{i \in I} \rho_i \\
& \geq S_{p,\alpha,d,\gamma} \left(\|f\|^p_{L^{p_\alpha^*}\left(\mathbb{R}^d\right)} + \sum_{i \in I} \tau_i^{\frac{p}{p_\alpha^*}} \right).
\end{align*}
Moreover, $\rho_\infty \geq S_{p,\alpha,d,\gamma} \tau_\infty^{\frac{p}{p_\alpha^*}}$. Thus, 
\begin{align*}
1 \geq \|f\|^p_{L^{p_\alpha^*}\left(\mathbb{R}^d\right)} + \sum_{i \in I} \tau_i^{\frac{p}{p_\alpha^*}} + \tau_\infty^{\frac{p}{p_\alpha^*}}. 
\end{align*}
By subadditivity, 
\begin{align*}
1 \geq \|f\|^p_{L^{p_\alpha^*}\left(\mathbb{R}^d\right)} + \left(\sum_{i \in I} \tau_i\right)^{\frac{p}{p_\alpha^*}} + \tau_\infty^{\frac{p}{p_\alpha^*}}.
\end{align*}
Finally, by Theorem \ref{thm:ccp_homo_frac_p},
\begin{align*}
1 = \underset{n \longrightarrow +\infty}{\lim} \int_{\mathbb{R}^d} |f_n(x)|^{p_\alpha^*}dx = \int_{\mathbb{R}^d} |f(x)|^{p_\alpha^{*}}dx + \sum_{i \in I} \tau_i + \tau_{\infty}. 
\end{align*}
Since $p \in (1, d/(\alpha-1))$, $p/p_\alpha^*<1$ and so, by subbaditivity again, 
\begin{align*}
1 = \left(\int_{\mathbb{R}^d} |f(x)|^{p_\alpha^{*}}dx + \sum_{i \in I} \tau_i + \tau_{\infty}\right)^{\frac{p}{p_\alpha^*}} \leq \|f\|^p_{L^{p_\alpha^*}(\mathbb{R}^d,dx)} + \left(\sum_{i \in I} \tau_i\right)^{\frac{p}{p_{\alpha}^*}} + \tau_\infty^{\frac{p}{p_\alpha^*}}.
\end{align*}
Thus, 
\begin{align*}
1 = \left(\int_{\mathbb{R}^d} |f(x)|^{p_\alpha^{*}}dx + \sum_{i \in I} \tau_i + \tau_{\infty}\right)^{\frac{p}{p_\alpha^*}} = \|f\|^p_{L^{p_\alpha^*}\left(\mathbb{R}^d,dx\right)} + \left(\sum_{i \in I} \tau_i\right)^{\frac{p}{p_\alpha^*}} + \tau_\infty^{\frac{p}{p_\alpha^*}},
\end{align*}
which ensures that $\|f\|_{L^{p_\alpha^*}\left(\mathbb{R}^d,dx\right)}$, $\tau_\infty$ and $\sum_{i \in I} \tau_i$ all belong to $\{0,1\}$. But, as mentionned before, $f \ne 0$. Therefore, $\|f\|_{L^{p_\alpha^*}\left(\mathbb{R}^d,dx\right)} = 1$ and $\tau_\infty = \sum_{i \in I} \tau_i =0$. This concludes the proof of the theorem.$\qed$

\section{Appendix}
\noindent
\begin{lem}\label{lem:technical_BV}
Let $d \geq 1$ be an integer, let $\alpha \in (1,2)$, let $\sigma_\alpha^*$ be defined by \eqref{def:norm_dual_syn_ndeg} and let $E\in \mathcal{B}(\mathbb{R}^d)$ with $\mathcal{L}_d(E)<+\infty$ such that 
\begin{align*}
\mathcal{P}_\alpha^{\operatorname{var}}(E) : = \sup \left\{ \int_E \operatorname{div}(\varphi) ,\, \varphi \in \mathcal{C}^{\infty}_c(\mathbb{R}^d , \mathbb{R}^d),\, \| \sigma_\alpha^*(\varphi) \|_{\infty} \leq 1\right\} <+\infty.
\end{align*}
Then,
\begin{align}\label{eq:C1_compact_replocal}
\mathcal{P}_\alpha^{\operatorname{var}}(E) = \sup \left\{ \int_E \operatorname{div}(\varphi) ,\, \varphi \in \mathcal{C}^{1}_c(\mathbb{R}^d , \mathbb{R}^d),\, \| \sigma_\alpha^*(\varphi) \|_{\infty} \leq 1\right\}
\end{align}
and
\begin{align}\label{eq:Cinfty_replocal}
\mathcal{P}_\alpha^{\operatorname{var}}(E) & = \sup \left\{ \int_E \operatorname{div}(\varphi) ,\, \varphi \in \mathcal{C}^{\infty}(\mathbb{R}^d , \mathbb{R}^d),\, \| \sigma_\alpha^*(\varphi) \|_{\infty} \leq 1, \, \operatorname{div}(\varphi)\in L^1(\mathbb{R}^d,dx)\right\} \nonumber \\
& = \sup \left\{ \int_E \operatorname{div}(\varphi) ,\, \varphi \in \mathcal{C}^{1}(\mathbb{R}^d , \mathbb{R}^d),\, \| \sigma_\alpha^*(\varphi) \|_{\infty} \leq 1, \, \operatorname{div}(\varphi)\in L^1(\mathbb{R}^d,dx)\right\}.
\end{align}
\end{lem}

\begin{proof}
Let us start with the proof of \eqref{eq:C1_compact_replocal}. First, it is clear that 
\begin{align*}
\mathcal{P}_\alpha^{\operatorname{var}}(E) \leq \sup \left\{ \int_E \operatorname{div}(\varphi) ,\, \varphi \in \mathcal{C}^{1}_c(\mathbb{R}^d , \mathbb{R}^d),\, \| \sigma_\alpha^*(\varphi) \|_{\infty} \leq 1\right\}.
\end{align*}
Next, let $\varphi \in \mathcal{C}_c^1(\mathbb{R}^d, \mathbb{R}^d)$ be such that $\| \sigma_\alpha^*(\varphi)\|_{\infty} \leq 1$. Let $\rho$ be a function in $\mathcal{C}_c^{\infty}(\mathbb{R}^d, \bbr_+)$ radially symmetric such that its compact support is a subset of $B(0,1)$ and
\begin{align*}
\int_{\mathbb{R}^d} \rho(x)dx = 1.  
\end{align*} 
Now, let $(\rho_\varepsilon)_{\varepsilon >0}$ be defined, for all $x \in \mathbb{R}^d$ and all $\varepsilon>0$, by
\begin{align*}
\rho_\varepsilon(x) = \frac{1}{\varepsilon^d} \rho\left(\frac{x}{\varepsilon}\right). 
\end{align*} 
Next, let $\left(\varphi_{\varepsilon}\right)_{\varepsilon>0}$ be defined, for all $\varepsilon>0$ and all $x \in \mathbb{R}^d$, by
\begin{align*}
\varphi_\varepsilon(x) = (\varphi \ast \rho_\varepsilon)(x). 
\end{align*} 
Then, $\varphi_\varepsilon$ is infinitely differentiable on $\mathbb{R}^d$ with compact support such that 
\begin{align*}
\| \sigma_\alpha^*(\varphi_{\varepsilon})\|_{\infty} \leq \int_{\mathbb{R}^d} \| \sigma_\alpha^*(\varphi)\|_{\infty} \rho_{\varepsilon}(x)dx \leq 1, 
\end{align*}
and such that, for all $x \in \mathbb{R}^d$, 
\begin{align*}
\operatorname{div} \left(\varphi_{\varepsilon}\right)(x) = (\operatorname{div}\left(\varphi\right) \ast \rho_{\varepsilon})(x).
\end{align*}
Moreover, for all $\varepsilon \in (0,1)$, 
\begin{align*}
\left| \int_E \operatorname{div} \left(\varphi\right)(x)dx - \int_E \operatorname{div} \left(\varphi_{\varepsilon}\right)(x)dx  \right| & = \left| \int_{\mathbb{R}^d} \rho_{\varepsilon}(y) \left(\int_{E} \operatorname{div} \left(\varphi\right)(x) - \operatorname{div} \left(\varphi\right)(x-y) dx \right) dy \right| \\
& =\left|  \int_{\mathbb{R}^d} \rho(y) \left(\int_{E} \operatorname{div} \left(\varphi\right)(x) - \operatorname{div} \left(\varphi\right)(x-\varepsilon y) dx \right) dy \right| \\
& \leq \int_{B(0,1)} \rho(y) \left(\int_{E} \left| \operatorname{div} \left(\varphi\right)(x) - \operatorname{div} \left(\varphi\right)(x-\varepsilon y)\right| dx \right) dy \\
& \leq \int_{B(0,1)} \rho(y) \bigg(\int_{(\operatorname{Supp}(\varphi) + B(0,1)) \cap E} \bigg| \operatorname{div} \left(\varphi\right)(x) \\
&\quad\quad - \operatorname{div} \left(\varphi\right)(x-\varepsilon y)\bigg| dx \bigg) dy.
\end{align*}
By the uniform continuity of $\operatorname{div}(\varphi)$ on $\operatorname{Supp}(\varphi) + B(0,1)$, for all $\eta >0$ and all $\varepsilon \in (0,1)$ small enough,
\begin{align*}
\left| \int_E \operatorname{div} \left(\varphi\right)(x)dx - \int_E \operatorname{div} \left(\varphi_{\varepsilon}\right)(x)dx  \right| & \leq \eta \mathcal{L}_d\left( (\operatorname{Supp}(\varphi) + B(0,1)) \cap E \right) \int_{B(0,1)} \rho(y) dy \leq \eta  \mathcal{L}_d\left( E\right).
\end{align*}
Then, 
\begin{align*}
\int_{E} \operatorname{div}(\varphi)(x) dx \leq  \eta  \mathcal{L}_d\left( E\right) + \mathcal{P}_\alpha^{\operatorname{var}}(E) . 
\end{align*}
This concludes the proof of \eqref{eq:C1_compact_replocal}. The proof of \eqref{eq:Cinfty_replocal} is rather straightforward and is based on a smooth truncation procedure. First, it is clear that
\begin{align*}
\mathcal{P}_\alpha^{\operatorname{var}}(E) \leq \sup \left\{ \int_E \operatorname{div}(\varphi) ,\, \varphi \in \mathcal{C}^{\infty}(\mathbb{R}^d , \mathbb{R}^d),\, \| \sigma_\alpha^*(\varphi) \|_{\infty} \leq 1, \, \operatorname{div}(\varphi)\in L^1(\mathbb{R}^d,dx)\right\},
\end{align*}
since, for all $\varphi \in \mathcal{C}_c^{\infty}(\mathbb{R}^d, \mathbb{R}^d)$ such that $\| \sigma_\alpha^*(\varphi) \|_{\infty} \leq 1$, $\operatorname{div}(\varphi)\in L^1(\mathbb{R}^d,dx)$. Conversely, take $\varphi \in \mathcal{C}^{\infty}(\mathbb{R}^d , \mathbb{R}^d)$ such that $\| \sigma_\alpha^*(\varphi) \|_{\infty} \leq 1$ and $\operatorname{div}(\varphi)\in L^1(\mathbb{R}^d,dx)$. Let $\chi$ be an infinitely differentiable radially symmetric function with compact support on $\mathbb{R}^d$ such that, $\chi(x) = 0$, for $x \in B(0,2)^c$, $\chi(x) = 1$, for $x \in B(0,1)$ and $\chi(x) \in [0,1]$, for $x \in B(0,2) \setminus B(0,1)$. Then, let $\varphi_R$ be defined, for all $x \in \mathbb{R}^d$ and all $R \geq 1$, by
\begin{align*}
\varphi_R(x) = \varphi(x)\chi\left(\frac{x}{R}\right) =  \varphi(x)\chi_R\left(x\right).
\end{align*}
Note that $\varphi_R$, $R \geq 1$, is a compactly supported $\mathbb{R}^d$-valued function which is infinitely differentiable on $\mathbb{R}^d$ and such that
\begin{align*}
\|  \sigma_\alpha^*(\varphi_R) \|_{\infty} \leq \| \sigma_\alpha^*(\varphi) \|_{\infty} \| \chi\|_{\infty} \leq 1,
\end{align*}
and
\begin{align*}
\operatorname{div}(\varphi_R) = \chi_R \operatorname{div}\left(\varphi\right) + \langle \varphi ; \nabla(\chi_R) \rangle \in L^1(\mathbb{R}^d,dx). 
\end{align*}
Then, 
\begin{align*}
\int_E \operatorname{div}(\varphi) & = \int_E \operatorname{div}(\varphi)-\int_E \operatorname{div}(\varphi_R)+\int_E \operatorname{div}(\varphi_R) \\
& \leq \sup \left\{ \int_E \operatorname{div}(\varphi), \, \varphi \in \mathcal{C}^{\infty}_c(\mathbb{R}^d,\mathbb{R}^d) , \, \| \sigma_\alpha^*(\varphi)\|_{\infty} \leq 1 \right\} + \int_E \operatorname{div}(\varphi)-\int_E \operatorname{div}(\varphi_R) \\
& \leq \mathcal{P}_\alpha^{\operatorname{var}}(E) + \int_E \operatorname{div}(\varphi)(x) \left(1 - \chi_R(x)\right)dx  -\frac{1}{R} \int_E \left\langle \varphi(x) ; \nabla(\chi)\left(\frac{x}{R}\right) \right\rangle dx.
\end{align*}
Now, 
\begin{align*}
\left| \frac{1}{R} \int_E \left\langle \varphi(x) ; \nabla(\chi)\left(\frac{x}{R}\right) \right\rangle dx\right| \leq \frac{1}{R} \|\varphi\|_{\infty} \| \nabla(\chi)\|_{\infty} \mathcal{L}_d(E) \longrightarrow 0,
\end{align*}
as $R$ tends to $+\infty$. Moreover, 
\begin{align*}
\left| \int_E \operatorname{div}(\varphi)(x) \left(1 - \chi_R(x)\right)dx \right| \leq \int_{\mathbb{R}^d} \left| \operatorname{div}(\varphi)(x) \right| \left(1 - \chi_R(x)\right)dx \leq \int_{B(0,R)^c} \left| \operatorname{div}(\varphi)(x) \right| dx \longrightarrow 0,
\end{align*}
as $R$ tends to $+\infty$. Thus, 
\begin{align*}
\int_E \operatorname{div}(\varphi)  \leq \mathcal{P}_\alpha^{\operatorname{var}}(E).
\end{align*}
This concludes the proof of the lemma. 
\end{proof}

\begin{prop}\label{prop:nondeg_symmetric_equality}
Let $\alpha \in (1,2)$ and let $\nu_\alpha$ be a non-degenerate symmetric L\'evy measure on $\mathbb{R}^d$, $d \geq 1$, verifying \eqref{eq:scale}. Let $\sigma_\alpha$ and $\sigma_\alpha^*$ be defined by \eqref{eq:rep_spectral_measure} and \eqref{def:norm_dual_syn_ndeg} respectively.  Let $\mathcal{P}_{\alpha}^{\operatorname{var}}$ and $\mathcal{P}_{\operatorname{cl}}$ be defined by \eqref{def:anisotropic_perimeter_stable} and \eqref{eq:classical} respectively. Then, for all $E$ Borel measurable subset of $\mathbb{R}^d$ with finite Lebesgue measure and finite perimeter, 
\begin{align}
\mathcal{P}_{\alpha}^{\operatorname{var}}(E) = \mathcal{P}_{\operatorname{cl}}(E). 
\end{align}
\end{prop}

\begin{proof}
The proof of this proposition is broken down into two intermediary steps. In the first step, let us observe that, for all $u \in \mathcal{C}^1(\mathbb{R}^d, \bbr)$ such that $\|\nabla(u)\|_{L^1(\mathbb{R}^d,dx)}<+\infty$, 
\begin{align}\label{eq:dual_rep_L1anisnorm_gradient}
\|\sigma_\alpha(\nabla(u))\|_{L^1(\mathbb{R}^d,dx)} = \sup \left\{ \int_{\mathbb{R}^d} u(x) \operatorname{div}(\varphi)(x):\, \varphi \in \mathcal{C}^1_c(\mathbb{R}^d, \mathbb{R}^d),\, \|\sigma_\alpha^*(\varphi)\|_{\infty} \leq 1 \right\}. 
\end{align}
Indeed, this is a consequence of \cite[Theorem $5.1$]{AB_AIHPNL94}.~Now, the L\'evy exponent of the $\alpha$-stable probability measure $\mu_\alpha$ satisfies the Hartman-Wintner condition (see \cite[Equation $(\operatorname{HW}_\infty)$]{KS_13}) since the L\'evy measure $\nu_\alpha$ is non-degenerate.  This implies that, for all multi-index $\beta \in \mathbb{N}_0^d$, 
\begin{align*}
\int_{\mathbb{R}^d} | D^{\beta}(p_\alpha)(x)| dx < +\infty,
\end{align*}
where $D^{\beta} = \partial_{x_1}^{\beta_1} \dots \partial_{x_d}^{\beta_d}$ and where $p_\alpha$ is the Lebesgue density of $\mu_\alpha$, so that, for all $t>0$, $P^\alpha_t(\bbone_E) \in \mathcal{C}^{\infty}(\mathbb{R}^d)$.  To lighten the notations, set, for all $x \in \mathbb{R}^d$ and all $t>0$, 
\begin{align*}
u_t(x) = P^\alpha_t(\bbone_E)(x). 
\end{align*}
Then, for all $t>0$, 
\begin{align*}
\|\sigma_\alpha(\nabla(u_t))\|_{L^1(\mathbb{R}^d,dx)} = \sup \left\{ \int_{\mathbb{R}^d} u_t(x) \operatorname{div}(\varphi)(x):\, \varphi \in \mathcal{C}^1_c(\mathbb{R}^d, \mathbb{R}^d),\, \|\sigma_\alpha^*(\varphi)\|_{\infty} \leq 1 \right\}. 
\end{align*}
Next, since $E$ has finite classical perimeter and since $\sigma_\alpha$ is a norm on $\mathbb{R}^d$,  $\mathcal{P}_{\alpha}^{\operatorname{var}}(E) < +\infty$ and $\mathcal{P}_{\operatorname{cl}}(E) < +\infty$. Moreover, the semigroup $(P_t^\alpha)_{t \geq 0}$ is strongly continuous on $L^1(\mathbb{R}^d,dx)$. Then, $u_t$ tends to $\bbone_E$ in $L^1(\mathbb{R}^d,dx)$ as $t$ tends to $0^+$. Thus, 
\begin{align*}
\mathcal{P}_\alpha^{\operatorname{var}}(E) \leq \underset{t \rightarrow 0^+}{\lim} \| \sigma_\alpha \nabla P^\alpha_t \bbone_E \|_{L^1(\mathbb{R}^d,dx)} =  \mathcal{P}_{\operatorname{cl}}(E).
\end{align*}
For the converse inequality,  by symmetry, for all $\varphi \in \mathcal{C}_c^1(\mathbb{R}^d, \mathbb{R}^d)$ such that $\| \sigma^*_\alpha(\varphi) \|_{\infty} \leq 1$ and all $t>0$, 
\begin{align*}
\int_{\mathbb{R}^d} P^\alpha_t(\bbone_E)(x) \operatorname{div}\left(\varphi\right)(x) dx = \int_{E} \operatorname{div} \left( P^\alpha_t(\varphi)\right)(x)dx. 
\end{align*}
Now, for all $t \geq 0$,  $P_t^\alpha(\varphi) \in \mathcal{C}^{1}(\mathbb{R}^d,\mathbb{R}^d)$, 
\begin{align*}
\| \sigma_\alpha^*\left(P_t^\alpha(\varphi)\right) \|_{\infty} \leq \| P_t^\alpha(\sigma_\alpha^*(\varphi)) \|_{\infty} \leq \| \sigma^*_\alpha(\varphi) \|_{\infty} \leq 1.
\end{align*}
Moreover, 
\begin{align*}
\| \operatorname{div}  \left( P^\alpha_t(\varphi)\right)  \|_{L^1(\mathbb{R}^d,dx)} \leq \| P^\alpha_t\left(\operatorname{div}(\varphi)\right)\|_{L^1(\mathbb{R}^d,dx)} \leq \| \operatorname{div}(\varphi) \|_{L^1(\mathbb{R}^d,dx)} <+\infty. 
\end{align*}
Then, 
\begin{align*}
\int_{\mathbb{R}^d} P^\alpha_t(\bbone_E)(x) \operatorname{div}\left(\varphi\right)(x) dx & \leq \sup \bigg\{ \int_E  \operatorname{div} \left( \varphi\right)(x)dx, \, \varphi \in \mathcal{C}^{1}(\mathbb{R}^d, \mathbb{R}^d),\,  \|\sigma_\alpha^*(\varphi)\|_{\infty} \leq 1, \\ 
&\quad\quad \|\operatorname{div}(\varphi)\|_{L^1(\mathbb{R}^d,dx)}<+\infty \bigg\} \\
& \leq  \mathcal{P}_\alpha^{\operatorname{var}}(E),
\end{align*}
using Lemma \ref{lem:technical_BV} for the last inequality. This concludes the proof of the proposition. 
\end{proof}
\noindent
In the sequel of this Appendix, let us prove similar results for the fractional perimeter $\mathcal{P}_{\operatorname{frac}}$ defined by \eqref{eq:perim_frac}. 

\begin{lem}\label{lem:technical_BVfrac}
Let $d \geq 1$ be an integer, let $\alpha \in (1,2)$, let $\operatorname{div}_\alpha$ be given by \eqref{eq:perim_fractional_divergence} and let $E$ be a Borel measurable subset of $\mathbb{R}^d$ with finite Lebesgue measure such that 
\begin{align}\label{eq:frac_perim_var_rep}
\mathcal{P}^{\operatorname{var}}_{\operatorname{frac}}(E) : = \sup\left\{\int_E \operatorname{div}_\alpha(\varphi)(x) dx: \, \varphi \in \mathcal{C}^{\infty}_c(\mathbb{R}^d, \mathbb{R}^d),\, \|\varphi\|_{\infty} \leq 1 \right\}< +\infty.  
\end{align}
Then, 
\begin{align}\label{eq:frac_perim_var_repC1}
\mathcal{P}^{\operatorname{var}}_{\operatorname{frac}}(E)  = \sup\left\{\int_E \operatorname{div}_\alpha(\varphi)(x) dx: \, \varphi \in \mathcal{C}^{1}_c(\mathbb{R}^d, \mathbb{R}^d),\, \|\varphi\|_{\infty} \leq 1 \right\},
\end{align}
and 
\begin{align}\label{eq:frac_perim_var_rep2}
\mathcal{P}^{\operatorname{var}}_{\operatorname{frac}}(E)  = \sup\left\{\int_E \operatorname{div}_\alpha(\varphi)(x) dx: \, \varphi \in \mathcal{C}^{\infty}(\mathbb{R}^d, \mathbb{R}^d),\, \|\varphi\|_{\infty} \leq 1, \|\nabla(\varphi)\|_{\infty}<+\infty \right\}. 
\end{align}
\end{lem}

\begin{proof}
As for Lemma \ref{lem:technical_BV}, the proof is also based on a regularization procedure. Clearly, for all Borel set $E$ as in the statement of the lemma, 
\begin{align*}
\mathcal{P}_{\operatorname{frac}}^{\operatorname{var}}(E) \leq \sup\left\{\int_E \operatorname{div}_\alpha(\varphi)(x) dx: \, \varphi \in \mathcal{C}^{1}_c(\mathbb{R}^d, \mathbb{R}^d),\, \|\varphi\|_{\infty} \leq 1 \right\}. 
\end{align*}
Next, let $\varphi \in \mathcal{C}^1_c(\mathbb{R}^d, \mathbb{R}^d)$ be such that $\|\varphi\|_{\infty}\leq 1$, let $\rho$ be as in the proof of Lemma \ref{lem:technical_BV} and let $(\rho_\varepsilon)_{\varepsilon>0}$ be the associated sequence of standard mollifiers.~Let $(\varphi_\varepsilon)_{\varepsilon>0}$ be defined, for all $x \in \mathbb{R}^d$ and all $\varepsilon>0$, by 
\begin{align*}
\varphi_{\varepsilon}(x) = (\varphi \ast \rho_\varepsilon)(x).
\end{align*} 
Then, for all $\varepsilon>0$, $\varphi_{\varepsilon} \in \mathcal{C}_c^{\infty}(\mathbb{R}^d, \mathbb{R}^d)$ with $\|\varphi_{\varepsilon}\|_{\infty} \leq 1$. Moreover, for all $\varepsilon >0$ and all $x \in \mathbb{R}^d$, 
\begin{align*}
\operatorname{div}_\alpha\left(\varphi_{\varepsilon}\right)(x) & = \sum_{k = 1}^d D_k^{\alpha-1}\left(\varphi_{\varepsilon,k}\right)(x) = \sum_{k = 1}^d \int_{\mathbb{R}^d} \left(\varphi_{\varepsilon,k}(x+u)-\varphi_{\varepsilon,k}(x)\right) u_k \nu_\alpha(du) \\
& =  \sum_{k = 1}^d \int_{\mathbb{R}^d} ((\varphi_k \ast \rho_{\varepsilon})(x+u) - (\varphi_k \ast \rho_{\varepsilon})(x)) u_k\nu_\alpha(du) \\
& =  \sum_{k=1}^d \int_{\mathbb{R}^d} \rho_\varepsilon(y) \left( \int_{\mathbb{R}^d} \left(\varphi_k(x+u-y) - \varphi_k(x-y)\right) u_k \nu_\alpha(du) \right)dy \\
& = \sum_{k = 1}^d \int_{\mathbb{R}^d} \rho_{\varepsilon}(y) D_k^{\alpha-1}(\varphi_k)(x-y) dy = \sum_{k=1}^d (\rho_\varepsilon \ast D_k^{\alpha-1}(\varphi_k))(x) = (\rho_{\varepsilon} \ast \operatorname{div}_\alpha(\varphi))(x).  
\end{align*}
Thus, for all $\varepsilon \in (0,1)$, 
\begin{align*}
\left| \int_{E} \operatorname{div}_\alpha(\varphi)(x)dx -  \int_{E} \operatorname{div}_\alpha(\varphi_\varepsilon)(x)dx \right| \leq \int_{B(0,1)} \rho(y) \left( \int_{E} \left| \operatorname{div}_\alpha(\varphi)(x) - \operatorname{div}_\alpha(\varphi)(x-\varepsilon y)\right| dx \right) dy. 
\end{align*}
Now, by the Lebesgue dominated convergence theorem, for all $x \in \mathbb{R}^d$ and all $y \in B(0,1)$, 
\begin{align*}
\underset{\varepsilon \rightarrow 0^+}{\lim} \operatorname{div}_\alpha(\varphi)(x-\varepsilon y) =  \operatorname{div}_\alpha(\varphi)(x).  
\end{align*}
Moreover, for all $x \in E$, all $y \in B(0,1)$ and all $\varepsilon \in (0,1)$, 
\begin{align*}
\left|  \operatorname{div}_\alpha(\varphi)(x-\varepsilon y) -  \operatorname{div}_\alpha(\varphi)(x) \right| \leq 2 \|\operatorname{div}_\alpha (\varphi)\|_{\infty} < + \infty. 
\end{align*}
Hence, since $\mathcal{L}_d(E) <+\infty$, 
\begin{align*}
\underset{\varepsilon \rightarrow 0^+}{\lim} \int_{B(0,1)} \rho(y) \left( \int_{E} \left| \operatorname{div}_\alpha(\varphi)(x) - \operatorname{div}_\alpha(\varphi)(x-\varepsilon y)\right| dx \right) dy = 0.
\end{align*}
Finally, let us prove the representation \eqref{eq:frac_perim_var_rep2}.~First of all, it is clear that 
\begin{align*}
\mathcal{P}_{\operatorname{frac}}^{\operatorname{var}}(E) \leq \sup\left\{\int_E \operatorname{div}_\alpha(\varphi)(x) dx: \, \varphi \in \mathcal{C}^{\infty}(\mathbb{R}^d, \mathbb{R}^d),\, \|\varphi\|_{\infty} \leq 1, \|\nabla(\varphi)\|_{\infty}<+\infty \right\}.
\end{align*}
Conversely, let $\varphi \in \mathcal{C}^{\infty}(\mathbb{R}^d, \mathbb{R}^d)$ be such that $\|\varphi\|_{\infty} \leq 1$ and that $\|\nabla(\varphi)\|_{\infty}<+\infty$. Now, let $\chi$ be a radially symmetric infinitely differentiable function with compact support such that $\chi =1$ on $B(0,1)$, $\chi = 0$ on $B(0,2)^c$ and $\chi \in [0,1]$ for $B(0,2) \setminus B(0,1)$. Then, for all $M\geq 1$ and all $x \in \mathbb{R}^d$, 
\begin{align*}
\varphi_M(x) = \chi \left(\frac{x}{M}\right) \varphi(x). 
\end{align*}  
Since the fractional gradient operator $D^{\alpha-1}$ is non-local, it does not follow the classical product rule. Nevertheless from \cite[proof of Lemma $3.2$]{AH23}, for all $x \in \mathbb{R}^d$ and all $M\geq 1$, 
\begin{align*}
\operatorname{div}_{\alpha} \left(\varphi_M\right)(x) & = \sum_{k = 1}^d D_k^{\alpha-1} \left(\varphi_{M,k}\right)(x) = \sum_{k=1}^d D^{\alpha-1}_k\left(\chi_M \varphi_k\right)(x) \\
& = \sum_{k = 1}^d \left(\chi_M(x) D^{\alpha-1}_k \left(\varphi_k\right)(x) + \varphi_k(x) D^{\alpha-1}_k(\chi_M)(x) + R_k^\alpha\left(\chi_M , \varphi_k\right)(x) \right),
\end{align*}
where, for all $k \in \{1, \dots, d\}$, all $M \geq 1$ and all $x \in \mathbb{R}^d$, 
\begin{align*}
R_k^\alpha\left(\chi_M , \varphi_k\right)(x) := \int_{\mathbb{R}^d} \left(\chi_M(x+u) - \chi_M(x)\right)\left(\varphi_k(x+u) - \varphi_k(x)\right) u_k \nu_\alpha(du). 
\end{align*}
Next, for all $M \geq 1$, 
\begin{align*}
\left| \int_E \operatorname{div}_\alpha(\varphi)(x)dx - \int_E \operatorname{div}_\alpha\left(\varphi_M\right)(x)dx \right| &\leq \left| \int_E \left(1- \chi\left(\frac{x}{M}\right)\right) \operatorname{div}_\alpha(\varphi)(x) dx \right| \\
&\quad\quad + \left| \int_E \langle \varphi(x) ; D^{\alpha-1}(\chi_M)(x) \rangle  dx \right| \\
&\quad\quad + \left| \int_E \sum_{k=1}^d R_k^\alpha\left(\chi_M , \varphi_k\right)(x)  dx \right|. 
\end{align*}
First, by the Lebesgue dominated convergence theorem, 
\begin{align*}
\underset{M \longrightarrow +\infty}{\lim} \int_E \left(1 - \chi_M(x)\right) \operatorname{div}_\alpha(\varphi)(x)dx = 0. 
\end{align*}
Moreover, for all $x \in \mathbb{R}^d$, 
\begin{align*}
\underset{M\rightarrow+\infty}{\lim} D^{\alpha-1}(\chi_M)(x) = 0
\end{align*}
and 
\begin{align*}
\underset{M \geq 1}{\sup} \| D^{\alpha-1}(\chi_M) \|_{\infty} <+\infty. 
\end{align*}
Thus, by the Lebesgue dominated convergence theorem again, 
\begin{align*}
\underset{M \rightarrow +\infty}{\lim}  \int_E \langle \varphi(x) ; D^{\alpha-1}(\chi_M)(x) \rangle  dx = 0. 
\end{align*}
An analogous reasoning proves that, for all $k \in \{1, \dots, d\}$, 
\begin{align*}
\underset{M \rightarrow +\infty}{\lim} \int_{E} \left| R_k^\alpha\left(\chi_M , \varphi_k\right)(x) \right| dx = 0,
\end{align*}
which concludes the proof of the lemma. 
\end{proof}

\begin{prop}\label{prop:variational_representation_fractional_perimeter}
Let $d \geq 1$ be an integer, let $\alpha \in (1,2)$, let $\nu_\alpha$ be a non-degenerate symmetric L\'evy measure on $\mathbb{R}^d$ verifying \eqref{eq:scale}. Let $\mathcal{P}_{\operatorname{frac}}^{\operatorname{var}}$ and $\mathcal{P}_{\operatorname{frac}}$ be defined by \eqref{eq:frac_perim_var_rep} and \eqref{eq:perim_frac} respectively. Then, for all $E\in \mathcal{B}(\mathbb{R}^d)$ with $\mathcal{L}_d(E)<+\infty$, 
\begin{align}\label{eq:pfrac_pfracvar}
\mathcal{P}_{\operatorname{frac}}^{\operatorname{var}}(E) = \mathcal{P}_{\operatorname{frac}}(E). 
\end{align}
\end{prop}

\begin{proof}
\textbf{Step 1:} First, let us treat the rotationally invariant case for which partial results are already available in the literature (see \cite{comi_stefani}).~As a first step, let us prove that, for all $u \in \mathcal{C}^1(\mathbb{R}^d) \cap L^1(\mathbb{R}^d,dx)$ such that $\|u\|_{\infty}<+\infty$, $\|\nabla(u)\|_{\infty} <+\infty$ and $\|D^{\alpha-1, \operatorname{rot}}(u)\|_{L^1(\mathbb{R}^d,dx)}<+\infty$, 
\begin{align}\label{eq:representation_integral_fractional_total_variation}
\int_{\mathbb{R}^d} \left\| D^{\alpha-1, \operatorname{rot}}(u)(x) \right\| dx = \sup \left\{ \int_{\mathbb{R}^d} u(x) \operatorname{div}^{\operatorname{rot}}_\alpha(\varphi)(x)dx:\, \varphi \in \mathcal{C}_c^{\infty}(\mathbb{R}^d, \mathbb{R}^d),\, \|\varphi\|_{\infty} \leq 1 \right\},
\end{align}
where $\operatorname{div}_\alpha^{\operatorname{rot}}$ is defined, for all $\varphi \in \mathcal{C}_c^{\infty}(\mathbb{R}^d, \mathbb{R}^d)$ and all $x \in \mathbb{R}^d$, by 
\begin{align*}
\operatorname{div}_\alpha^{\operatorname{rot}}(\varphi)(x) = \sum_{k = 1}^d D_k^{\alpha-1, \operatorname{rot}}(\varphi_k)(x). 
\end{align*}
For this purpose, let $BV^\alpha(\mathbb{R}^d)$ be the fractional bounded variation space defined by
\begin{align*}
BV^\alpha(\mathbb{R}^d) = \{f \in L^1(\mathbb{R}^d,dx):\, |f|_{BV^\alpha(\mathbb{R}^d)}<+\infty\}, 
\end{align*}
where $|f|_{BV^\alpha(\mathbb{R}^d)}$ is given by 
\begin{align*}
|f|_{BV^\alpha(\mathbb{R}^d)} = \sup \left\{  \int_{\mathbb{R}^d} f(x) \operatorname{div}^{\operatorname{rot}}_\alpha(\varphi)(x)dx:\, \varphi \in \mathcal{C}_c^{\infty}(\mathbb{R}^d, \mathbb{R}^d),\, \|\varphi\|_{\infty} \leq 1 \right\}. 
\end{align*}
Note that, up to a normalization constant, the space $BV^\alpha(\mathbb{R}^d)$ is the one introduced in \cite[Definition $3.1$]{comi_stefani}. First, for all $u \in \mathcal{C}^1(\mathbb{R}^d) \cap L^1(\mathbb{R}^d,dx)$ such that $\|u\|_{\infty}<+\infty$, $\|\nabla(u)\|_{\infty} <+\infty$ and $\|D^{\alpha-1, \operatorname{rot}}(u)\|_{L^1(\mathbb{R}^d,dx)}<+\infty$ and all $\varphi \in \mathcal{C}_c^{\infty}(\mathbb{R}^d, \mathbb{R}^d)$ such that $\| \varphi \|_{\infty}\leq 1$,  
\begin{align}\label{eq:ipp_fraction}
\left| \int_{\mathbb{R}^d} u(x) \operatorname{div}^{\operatorname{rot}}_{\alpha}(\varphi)(x)dx \right|  & = \left| \int_{\mathbb{R}^d} \langle D^{\alpha-1, \operatorname{rot}}(u)(x) ; \varphi(x) \rangle dx \right|\nonumber \\
& \leq  \int_{\mathbb{R}^d} \left\| D^{\alpha-1, \operatorname{rot}}(u)(x) \right\| dx <+\infty. 
\end{align}
Then, $u$ belongs to $BV^\alpha(\mathbb{R}^d)$ and 
\begin{align*}
|u|_{BV^\alpha(\mathbb{R}^d)} \leq  \int_{\mathbb{R}^d} \left\| D^{\alpha-1, \operatorname{rot}}(u)(x) \right\| dx. 
\end{align*}
Now, \cite[Theorem $3.2$]{comi_stefani} asserts that there exists a finite vector-valued Radon measure $d_{\alpha-1}^{\operatorname{rot}}(u)$ such that, for all $\varphi \in \mathcal{C}^{\infty}_c(\mathbb{R}^d , \mathbb{R}^d)$ with $\|\varphi\|_{\infty}\leq 1$, 
\begin{align*}
\int_{\mathbb{R}^d} u(x) \operatorname{div}_\alpha(\varphi)(x)dx =  - \int_{\mathbb{R}^d} \langle \varphi(x); d_{\alpha-1}^{\operatorname{rot}}(u)(dx) \rangle \Rightarrow d_{\alpha-1}^{\operatorname{rot}}(u)(dx) = D^{\alpha-1, \operatorname{rot}}(u)(x) \mathcal{L}_d(dx). 
\end{align*}
Moreover, by \cite[Theorem $3.2$, Equation $(3.2)$]{comi_stefani}, 
\begin{align*}
\left| d_{\alpha-1}^{\operatorname{rot}}(u)\right| (\mathbb{R}^d) = \left| u \right|_{BV^\alpha(\mathbb{R}^d)} \Rightarrow \int_{\mathbb{R}^d} \left\| D^{\alpha-1, \operatorname{rot}}(u)(x) \right\| dx =  \left| u \right|_{BV^\alpha(\mathbb{R}^d)},
\end{align*}
and, so equality \eqref{eq:representation_integral_fractional_total_variation} is proved. Next, for all $t>0$ and all $E \in \mathcal{B}(\mathbb{R}^d)$ such that $\mathcal{L}_d(E)<+\infty$, $P^{\alpha, \operatorname{rot}}_t(\bbone_E) \in \mathcal{C}^1(\mathbb{R}^d) \cap L^1(\mathbb{R}^d,dx)$ with $\| P^{\alpha, \operatorname{rot}}_t(\bbone_E) \|_{\infty} < +\infty$, $\| \nabla \left(P^{\alpha, \operatorname{rot}}_t(\bbone_E)\right)\|_{\infty} < +\infty$ and, for all $t>0$,  
\begin{align*}
\int_{\mathbb{R}^d} \left\| D^{\alpha - 1, \operatorname{rot}} P^{\alpha, \operatorname{rot}}_t(\bbone_E)(x) \right\| dx & \leq \frac{1}{t^{1- \frac{1}{\alpha}}} \mathcal{L}_d(E) \left(\int_{\mathbb{R}^d} \left\| D^{\alpha-1, \operatorname{rot}} (p_\alpha) \left( x \right) \right\| dx\right) < +\infty,
\end{align*}
since $\alpha \in (1,2)$. Next, setting $u_t = P^{\alpha, \operatorname{rot}}_t(\bbone_E)$, for all $t>0$, 
\begin{align*}
\int_{\mathbb{R}^d} \left\| D^{\alpha-1, \operatorname{rot}}(u_t)(x) \right\| dx = |u_t|_{BV^\alpha(\mathbb{R}^d)}. 
\end{align*}
Therefore, 
\begin{align*}
\underset{t\rightarrow 0^+}{\liminf}  |u_t|_{BV^\alpha(\mathbb{R}^d)} = \underset{t\rightarrow 0^+}{\liminf}  \int_{\mathbb{R}^d} \left\| D^{\alpha-1, \operatorname{rot}}(u_t)(x) \right\| dx = \mathcal{P}_{\operatorname{frac}}(E). 
\end{align*}
Moreover, $u_t$ converges in $L^1(\mathbb{R}^d,dx)$ to $\bbone_E$ as $t$ tends to $0^+$, and so \cite[Proposition $3.3$]{comi_stefani} implies, 
\begin{align*}
\mathcal{P}_{\operatorname{frac}}^{\operatorname{var}}(E) \leq \mathcal{P}_{\operatorname{frac}}(E). 
\end{align*}
Now, by symmetry and commutation, for all $t>0$ and all $\varphi \in \mathcal{C}_c^{\infty}(\mathbb{R}^d, \mathbb{R}^d)$ with $\| \varphi\|_{\infty} \leq 1$,
\begin{align*}
\int_{\mathbb{R}^d} u_t(x) \operatorname{div}_\alpha^{\operatorname{rot}}(\varphi)(x)dx = \int_{E} P^{\alpha, \operatorname{rot}}_t \left( \operatorname{div}_\alpha^{\operatorname{rot}}(\varphi)\right)(x)dx = \int_E \operatorname{div}_\alpha^{\operatorname{rot}}\left(P^{\alpha, \operatorname{rot}}_t(\varphi)\right)(x)dx.   
\end{align*} 
Moreover, for all $t \geq 0$, $P^{\alpha, \operatorname{rot}}_t(\varphi) \in \mathcal{C}^{\infty}(\mathbb{R}^d, \mathbb{R}^d)$ and, by contractivity,
\begin{align*}
\left\| P^{\alpha, \operatorname{rot}}_t(\varphi)\right\|_{\infty} \leq \left\| \varphi\right\|_{\infty} \leq 1. 
\end{align*}
Finally, for all $t>0$ and all $x \in \mathbb{R}^d$, 
\begin{align*}
\nabla P^{\alpha, \operatorname{rot}}_t\left(\bbone_E\right)(x) = \frac{1}{t^{\frac{1}{\alpha}}} \int_E \nabla(p_\alpha^{\operatorname{rot}}) \left(\frac{x-y}{t^{\frac{1}{\alpha}}}\right) \frac{dy}{t^{\frac{d}{\alpha}}} \Rightarrow \left\| \nabla P^{\alpha, \operatorname{rot}}_t\left(\bbone_E\right) \right\|_{\infty} \leq \frac{1}{t^{\frac{1}{\alpha}}} \int_{\mathbb{R}^d} \|\nabla(p_\alpha^{\operatorname{rot}})(x)\|dx. 
\end{align*}
Thus, for all $t>0$, 
\begin{align*}
\|D^{\alpha-1, \operatorname{rot}}(u_t)\|_{L^1(\mathbb{R}^d,dx)} \leq \sup \left\{ \int_E \operatorname{div}^{\operatorname{rot}}_{\alpha}(\varphi)(x)dx: \, \varphi \in \mathcal{C}^{\infty}(\mathbb{R}^d, \mathbb{R}^d),\, \|\varphi\|_{\infty}\leq 1,\, \|\nabla(\varphi)\|_{\infty}<+\infty \right\}, 
\end{align*}
which implies that 
\begin{align*}
\mathcal{P}_{\operatorname{frac}}(E) \leq \sup \left\{ \int_E \operatorname{div}^{\operatorname{rot}}_{\alpha}(\varphi)(x)dx: \, \varphi \in \mathcal{C}^{\infty}(\mathbb{R}^d, \mathbb{R}^d),\, \|\varphi\|_{\infty}\leq 1,\, \|\nabla(\varphi)\|_{\infty}<+\infty \right\}.
\end{align*}
To conclude the proof of the rotationally invariant case, it remains to prove that 
\begin{align*}
\sup \left\{ \int_E \operatorname{div}^{\operatorname{rot}}_{\alpha}(\varphi)(x)dx: \, \varphi \in \mathcal{C}^{\infty}(\mathbb{R}^d, \mathbb{R}^d),\, \|\varphi\|_{\infty}\leq 1,\, \|\nabla(\varphi)\|_{\infty}<+\infty \right\} \leq \mathcal{P}_{\operatorname{frac}}^{\operatorname{var}}(E). 
\end{align*}
But this is a consequence of \eqref{eq:frac_perim_var_rep2} in Lemma \ref{lem:technical_BVfrac} and so the rotationally invariant case is done.\\
\textbf{Step 2:} Now, let us deal with the general non-degenerate and symmetric case. The proof of the equality \eqref{eq:pfrac_pfracvar} for the non-degenerate and symmetric case is similar to the one of the rotationally invariant case once a structure theorem and the lower semicontinuity of the general fractional variation $| \cdot |_\alpha$ are proved for $BV_\alpha (\mathbb{R}^d)$.\\
\textbf{Step 3:} Let us start with the proof of a structure theorem for $BV_\alpha(\mathbb{R}^d)$ functions.~Namely, $f \in BV_\alpha (\mathbb{R}^d)$ if and only if there exists a finite $\mathbb{R}^d$-valued Radon measure $d_{\alpha-1}(f)$ on $\mathcal{B}(\mathbb{R}^d)$ such that, for all $\varphi \in \mathcal{C}_c^{\infty}(\mathbb{R}^d, \mathbb{R}^d)$, 
\begin{align}\label{eq:frac_ipp_NDS}
\int_{\mathbb{R}^d} f(x) \operatorname{div}_\alpha(\varphi)(x)dx = - \int_{\mathbb{R}^d} \langle \varphi(x) ; d_{\alpha-1}(f)(dx) \rangle. 
\end{align}
First of all, let $f \in L^1(\mathbb{R}^d,dx)$ be such that there exists a finite $\mathbb{R}^d$-valued Radon measure $d_{\alpha-1}(f)$ verifying \eqref{eq:frac_ipp_NDS}.~Then, for all $\varphi \in \mathcal{C}^{\infty}_c(\mathbb{R}^d,\mathbb{R}^d)$ such that $\|\varphi\|_{\infty} \leq 1$, 
 \begin{align*}
 \left| \int_{\mathbb{R}^d} f(x) \operatorname{div}_\alpha(\varphi)(x)dx \right| \leq \left|  \int_{\mathbb{R}^d} \langle \varphi(x) ; d_{\alpha-1}(f)(dx) \rangle \right| \leq \|\varphi\|_{\infty} \left|d_{\alpha-1}(f)\right|(\mathbb{R}^d)<+\infty. 
 \end{align*}
 Thus, $f \in BV_\alpha(\mathbb{R}^d)$. Conversely, let $f$ be a function in $L^1(\mathbb{R}^d,dx)$ such that $|f|_{\alpha}<+\infty$ and let $L_f$ be the linear functional defined, for all $\varphi \in \mathcal{C}_c^{\infty}(\mathbb{R}^d , \mathbb{R}^d)$, by 
\begin{align*}
L_f(\varphi) : = -  \int_{\mathbb{R}^d} f(x) \operatorname{div}_\alpha(\varphi)(x)dx. 
\end{align*}
Since $f \in BV_\alpha(\mathbb{R}^d)$, for all $\varphi \in \mathcal{C}^{\infty}_c(\mathbb{R}^d, \mathbb{R}^d)$, 
\begin{align*}
\left| L_f(\varphi) \right| \leq |f|_{\alpha} \|\varphi\|_{\infty}. 
\end{align*}
Moreover, by a standard approximation argument, $\mathcal{C}_c^{\infty}(\mathbb{R}^d, \mathbb{R}^d)$ is a dense subset of $\mathcal{C}_c(\mathbb{R}^d, \mathbb{R}^d)$ with respect to uniform convergence. Then, $L_f$ extends continuously to a linear functional $\tilde{L}_f$ from $\mathcal{C}_c(\mathbb{R}^d, \mathbb{R}^d)$ to $\bbr$. By a version of the Riesz representation theorem (see, e.g., \cite[Theorem $1.38$]{EG_15}), there exists a finite $\mathbb{R}^d$-valued Radon measure $d_{\alpha-1}(f)$ such that, for all $\varphi \in \mathcal{C}_c(\mathbb{R}^d, \mathbb{R}^d)$,
\begin{align*}
\tilde{L}_f(\varphi) = \int_{\mathbb{R}^d} \langle \varphi(x) ; d_{\alpha-1}(f)(dx) \rangle.  
\end{align*}
In particular, for all $\varphi \in \mathcal{C}_c^{\infty}(\mathbb{R}^d, \mathbb{R}^d)$, 
\begin{align*}
\int_{\mathbb{R}^d} f(x) \operatorname{div}_{\alpha}(\varphi)(x)dx = - \int_{\mathbb{R}^d} \langle \varphi(x) ; d_{\alpha-1}(f)(dx) \rangle.  
\end{align*}
Moreover, 
\begin{align*}
\left| d_{\alpha-1}(f) \right|(\mathbb{R}^d) = \sup \left\{ \tilde{L}_f(\varphi):\, \varphi \in \mathcal{C}_c(\mathbb{R}^d, \mathbb{R}^d),\, \| \varphi\|_{\infty} \leq 1 \right\} = |f|_{\alpha}. 
\end{align*}
This concludes the proof of the structure result for $BV_\alpha(\mathbb{R}^d)$.\\ 
\textbf{Step 4:} In this last step, let us prove the lower semicontinuity of the general fractional variation $|\cdot|_\alpha$ as in \cite[Proposition $3.3$]{comi_stefani}.~Let $(f_n)_{n \geq 1}$ be a sequence of elements in $BV_\alpha(\mathbb{R}^d)$ and let $f \in L^1(\mathbb{R}^d, dx)$ such that $f_n \longrightarrow f$ in $L^1(\mathbb{R}^d,dx)$, as $n$ tends to $+\infty$. Then, for all $\varphi \in \mathcal{C}_c^{\infty}(\mathbb{R}^d , \mathbb{R}^d)$ such that $\|\varphi\|_{\infty} \leq 1$, 
\begin{align*}
- \int_{\mathbb{R}^d} f(x) \operatorname{div}_\alpha(\varphi)(x)dx = - \underset{n \rightarrow +\infty}{\lim} \int_{\mathbb{R}^d} f_n(x) \operatorname{div}_\alpha(\varphi)(x) dx & = \underset{n \rightarrow +\infty}{\lim} \int_{\mathbb{R}^d} \langle \varphi ; d_{\alpha-1}(f_n)(dx)\rangle \\
& \leq \underset{n\longrightarrow +\infty}{\liminf} \left|f_n \right|_{\alpha},
\end{align*}
which implies that
\begin{align*}
\left| f\right|_{\alpha} \leq \underset{n \rightarrow +\infty}{\liminf} \left|f_n \right|_{\alpha} <+\infty.
\end{align*}
This concludes the proof of the proposition. 
\end{proof}
\noindent
The next proposition is an approximation result in the space of fractional bounded variation functions $BV_\alpha(\mathbb{R}^d)$ needed in order to obtain compactness. 

\begin{prop}\label{prop:approximation_general_case_BV}
Let $\alpha \in (1,2)$ and let $f \in BV_\alpha(\mathbb{R}^d)$, $d \geq 1$. Then, there exists $(f_k)_{k \geq 1} \in \mathcal{C}_c^{\infty}(\mathbb{R}^d)$ such that
\begin{align*}
f_k \longrightarrow f\, \operatorname{in} \, L^1(\mathbb{R}^d,dx)
\end{align*}
and
\begin{align*}
|f_k|_{\alpha} \longrightarrow |f|_{\alpha},
\end{align*}
as $k$ tends to $+\infty$. 
\end{prop}

\begin{proof}
First, let $(\rho_\varepsilon)_{\varepsilon >0}$ be a sequence of standard mollifiers as in the proof of Lemma \ref{lem:technical_BV} and let $(f_\varepsilon)_{\varepsilon>0}$ be defined, for all $\varepsilon>0$, by $f_\varepsilon = f \ast \rho_\varepsilon$. Clearly, $f_\varepsilon \rightarrow f$ in $L^1(\mathbb{R}^d,dx)$, as $\varepsilon$ tends to $0^+$. Moreover, for all $\varepsilon>0$, $f_\varepsilon \in \mathcal{C}^{\infty}(\mathbb{R}^d)$. Thus, for all $\varphi \in \mathcal{C}_c^{\infty}(\mathbb{R}^d, \mathbb{R}^d)$,
\begin{align*}
- \int_{\mathbb{R}^d} f_\varepsilon(x) \operatorname{div}_\alpha(\varphi)(x)dx & = - \int_{\mathbb{R}^d} (f\ast \rho_\varepsilon)(x) \operatorname{div}_\alpha(\varphi)(x)dx  \\
& = - \int_{\mathbb{R}^d} f(x) (\rho_\varepsilon \ast \operatorname{div}_\alpha(\varphi))(x) dx \\
& = - \int_{\mathbb{R}^d} f(x) \operatorname{div}_\alpha \left(\rho_\varepsilon \ast \varphi\right)(x) dx \\
& =  \int_{\mathbb{R}^d} \langle (\rho_\varepsilon \ast \varphi)(x) ; d_{\alpha-1}(f)(dx) \rangle \\
& = \int_{\mathbb{R}^d} \langle \varphi(x) ; (\rho_\varepsilon \ast d_{\alpha-1}(f))(x) \rangle dx. 
\end{align*} 
Then, for all $\varepsilon>0$, $f_\varepsilon \in BV_\alpha(\mathbb{R}^d)$ and $d_{\alpha-1}(f_\varepsilon)(dx) = (\rho_\varepsilon \ast d_{\alpha-1}(f))(x)dx$. But, for all $\varepsilon>0$, 
\begin{align*}
\left| d_{\alpha-1}(f_\varepsilon) \right|(\mathbb{R}^d) = \int_{\mathbb{R}^d} \| (\rho_\varepsilon \ast d_{\alpha-1}(f))(x) \| dx \leq |d_{\alpha-1}(f)(\mathbb{R}^d)|<+\infty.
\end{align*}
Finally, by lower semicontinuity, 
\begin{align*}
|d_{\alpha-1}(f)|(\mathbb{R}^d) \leq \underset{\varepsilon \rightarrow 0^+}{\liminf}  |d_{\alpha-1}(f_\varepsilon)|(\mathbb{R}^d). 
\end{align*}
Then, $|d_{\alpha-1}(f_\varepsilon)|(\mathbb{R}^d) \rightarrow |d_{\alpha-1}(f)|(\mathbb{R}^d)$, as $\varepsilon$ tends to $0^+$. To conclude this step, let us proceed using a smooth truncation argument. Assume that $f \in BV_\alpha(\mathbb{R}^d) \cap\, \mathcal{C}^{\infty}(\mathbb{R}^d)$ and let $\chi \in \mathcal{C}^{\infty}_c(\mathbb{R}^d, [0,1])$ be such that $\chi(x) = 1$, for $x \in B(0,1)$, and $\chi(x) = 0$, for $x \in B(0,2)^c$. For all $R \geq 1$, let $f_R$ be defined, for all $x \in \mathbb{R}^d$, by
\begin{align*}
f_R(x) = \chi \left(\frac{x}{R}\right) f(x). 
\end{align*}
For all $R \geq 1$, $f_R \in \mathcal{C}_c^{\infty}(\mathbb{R}^d)$.~By the Lebesgue dominated convergence theorem, $f_R \rightarrow f$ in $L^1(\mathbb{R}^d,dx)$, as $R$ tends to $+\infty$. Then, by lower semicontinuity, 
\begin{align*}
|d_{\alpha-1}(f)|(\mathbb{R}^d) \leq \underset{R \rightarrow +\infty}{\liminf} |d_{\alpha-1}(f_R)|(\mathbb{R}^d).
\end{align*}
Finally, for all $R\geq 1$ and all $\varphi \in \mathcal{C}^{\infty}_c(\mathbb{R}^d, \mathbb{R}^d)$ such that $\| \varphi \|_{\infty} \leq 1$, 
\begin{align*}
- \int_{\mathbb{R}^d} f_R(x) &\operatorname{div}_\alpha(\varphi)(x)dx = - \sum_{k = 1}^d \int_{\mathbb{R}^d} f(x) \chi_R(x) D^{\alpha-1}_k(\varphi_k)(x) dx \\
& = - \sum_{k = 1}^d \int_{\mathbb{R}^d} f(x) \left(D^{\alpha-1}_k \left(\chi_R \varphi_k\right)(x) - \varphi_k(x) D^{\alpha-1}_k(\chi_R)(x) - R^\alpha_k(\chi_R , \varphi_k)(x) \right)dx. 
\end{align*}
Now, for all $R \geq 1$, 
\begin{align*}
\int_{\mathbb{R}^d} f(x) \operatorname{div}_\alpha(\chi_R \varphi)(x)dx \leq \| \chi_R \varphi\|_{\infty} |f|_{\alpha} \leq |f|_{\alpha}.  
\end{align*}
Moreover, for all $R \geq 1$, 
\begin{align*}
\left| \int_{\mathbb{R}^d} f(x) \langle \varphi(x) ; D^{\alpha-1}(\chi_R)(x) \rangle dx \right| \leq \int_{\mathbb{R}^d} |f(x)| \| D^{\alpha-1}(\chi_R)(x) \| dx.
\end{align*}
A direct application of the Lebesgue dominated convergence theorem ensures that 
\begin{align*}
\underset{R \rightarrow +\infty}{\lim} \int_{\mathbb{R}^d} |f(x)| \| D^{\alpha-1}(\chi_R)(x) \| dx = 0. 
\end{align*} 
Finally, for all $R \geq 1$, 
\begin{align*}
\left| \sum_{k = 1}^d \int_{\mathbb{R}^d} f(x) R^\alpha_k(\chi_R , \varphi_k)(x) dx \right| \leq 2 \int_{\mathbb{R}^d} |f(x)| \int_{\mathbb{R}^d} |\chi_R(x+u) - \chi_R(x)|  \|u\| \nu_\alpha(du) dx. 
\end{align*}
Once again, an application of the Lebesgue dominated convergence theorem ensures that
\begin{align*}
\underset{R \rightarrow +\infty}{\lim}\int_{\mathbb{R}^d} |f(x)| \left(\int_{\mathbb{R}^d} |\chi_R(x+u) - \chi_R(x)|  \|u\| \nu_\alpha(du)\right)dx=0. 
\end{align*} 
Thus, 
\begin{align*}
\underset{R \rightarrow +\infty}{\limsup} |f_R|_\alpha \leq |f|_\alpha,
\end{align*}
which concludes the proof of the proposition.
\end{proof}
\noindent
To finish this appendix, let us state and prove:

\begin{lem}\label{lemma_technical_Lorentz_Hardy}
Let $H$ be a norm on $\mathbb{R}^d$, $d \geq 2$, and let $p \in [1,d)$. Let $f$ be a non-negative convex symmetric function with respect to the norm $H$. Then, 
\begin{align}\label{eq:Lorentz_Hardy}
\|f\|_{p^*,p} = \mathcal{L}_d(B_H)^{-\frac{1}{d}} \left(\int_{\mathbb{R}^d} \dfrac{f(x)^p}{H(x)^p}dx\right)^{\frac{1}{p}},
\end{align}
where $p^* = dp /(d-p)$,  $B_H$ is the unit ball with respect to the norm $H$ and $\|.\|_{p^*,p}$ is the Lorentzian quasi-norm $(p^*,p)$.
\end{lem}

\begin{proof}
The proof is similar to the one of \cite[Lemma 4.3]{FS_08} using formula \eqref{eq:polar_coordinate_cone_measure}. Let $g(x) = f(x)^p$, $x \in \mathbb{R}^d$,  and let $G$ be the distribution function of $g$ defined, for all $t \geq 0$, by 
\begin{align}\label{ineq:distribution_function_G}
G(t) = \mathcal{L}_d\left(\{x \in \mathbb{R}^d : \, g(x) > t\}\right). 
\end{align}
Moreover, let $F$ be the distribution function of $f$. Now, by the very definition of the Lorentzian quasi-norm $(p^*,p)$, 
\begin{align}\label{eq:Lorentz_qn}
\|f\|^p_{p^*,p} = p^* \int_0^{+\infty} F(t)^{\frac{p}{p^*}} t^{p-1} dt = \frac{p^*}{p} \int_0^{+\infty} G(t)^{\frac{p}{p^*}} dt.
\end{align}
Next, using the layer cake representation formula for $g$ and Fubini's theorem, 
\begin{align*}
\int_{\mathbb{R}^d} \dfrac{f(x)^p}{H(x)^p}dx = \int_0^{+\infty}\left( \int_{\mathbb{R}^d} \bbone_{g(x)>t} \frac{dx}{H(x)^p} \right) dt.
\end{align*}
By the polar-coordinates formula \eqref{eq:polar_coordinate_cone_measure}, 
\begin{align*}
 \int_{\mathbb{R}^d} \bbone_{g(x)>t} \frac{dx}{H(x)^p} = d \mathcal{L}_d(B_H) \int_0^{+\infty} s^{d-p-1} \left(\int_{\partial B_H} \bbone_{g(x)>t}(s \omega) \mu_{B_H}(d\omega) \right) ds,
\end{align*} 
where $B_H = \{x \in \mathbb{R}^d: \, H(x) \leq 1\}$ and where $\mu_{B_H}$ is the cone measure on $\partial B_H$.~Now, $f$ is convex symmetric with respect to the norm $H$ and so is $g$. In particular, the level sets $\{g(x) >t\}$, for $t >0$, are $H$-balls with radius $R(t)$ given, for all $t>0$, by
\begin{align*}
R(t) : =  \left(  \frac{G(t)}{\mathcal{L}_d(B_H)} \right)^{\frac{1}{d}}. 
\end{align*} 
Thus, for all $t>0$, 
\begin{align*}
 \int_{\mathbb{R}^d} \bbone_{g(x)>t} \frac{dx}{H(x)^p} & = d \mathcal{L}_d(B_H) \int_0^{R(t)} s^{d-p-1} ds = \frac{d}{d-p} \mathcal{L}_d(B_H)^{\frac{p}{d}} G(t)^{\frac{d-p}{d}}. 
\end{align*}
The end of the proof easily follows. 
\end{proof}

\end{document}